\documentclass[a4paper,12pt,amsart,frenchb]{article}

\usepackage{amsmath,amsbsy,amsfonts,amssymb}
\usepackage[french]{babel}
\newcommand{\ass}[2]{\vskip0.3cm\noindent
{\bf {#1}}. { \sl {#2}}\vskip0.3cm\noindent
}
  
\begin{document}

    \title{Stabilisation de la formule des traces tordue VII:   descente globale}
\author{ J.-L. Waldspurger}
\date{2 septembre 2014}
\maketitle

\bigskip

{\bf Introduction}

Nous commen\c{c}ons la preuve des th\'eor\`emes [VI] 5.2 et [VI] 5.4. Rappelons-en les \'enonc\'es, en renvoyant \`a [VI] pour les d\'efinitions. Le triplet $(G,\tilde{G},{\bf a})$ est d\'efini sur un corps de nombres $F$. On fixe un ensemble fini $V$ de places de $F$ contenant l'ensemble $V_{ram}$ des "mauvaises" places.

\ass{Th\'eor\`eme [VI] 5.2 (\`a prouver)}{On suppose $(G,\tilde{G},{\bf a})$ quasi-d\'eploy\'e et \`a torsion int\'erieure. Soit ${\cal O}_{V}$ une classe de conjugaison stable semi-simple dans $\tilde{G}(F_{V})$. Alors $SA^{\tilde{G}}(V,{\cal O}_{V})$ est stable.}

\ass{Th\'eor\`eme [VI] 5.4 (\`a prouver)}{Soit ${\cal O}_{V}$ une classe de conjugaison stable semi-simple dans $\tilde{G}(F_{V})$. Alors on a l'\'egalit\'e $A^{\tilde{G},{\cal E}}(V,{\cal O}_{V},\omega)=A^{\tilde{G}}(V,{\cal O}_{V},\omega)$.}

Dans cet article, nous prouverons le premier th\'eor\`eme. Nous prouverons aussi le second sauf pour certains triplets $(G,\tilde{G},{\bf a})$ particuliers. Pour ceux-ci, nous le prouverons sauf pour un nombre fini de classes ${\cal O}_{V}$ exceptionnelles. On renvoie \`a 3.5 pour des assertions pr\'ecises. Les cas restants de ce th\'eor\`eme seront prouv\'es plus tard en utilisant la formule des traces. 

  On introduit en 1.1  un ensemble de param\`etres ${\bf Stab}(\tilde{G}(F))$. On d\'efinit  en 1.2  une application qui, \`a une classe de conjugaison stable semi-simple dans $\tilde{G}(F)$, associe un \'el\'ement de notre ensemble de param\`etres. Cette application est toujours injective. En g\'en\'eral, elle n'est pas surjective. Elle l'est toutefois si $(G,\tilde{G},{\bf a})$ est quasi-d\'eploy\'e et \`a torsion int\'erieure. La question de savoir si un param\`etre provient bel et bien d'une classe de conjugaison stable semi-simple est d\'elicate et nous ne la r\'esoudrons pas. L'int\'er\^et d'introduire cet ensemble de param\`etres est justement de la contourner. Si ${\bf G}'=(G',{\cal G}',\tilde{s})$ est une donn\'ee endoscopique de $(G,\tilde{G},{\bf a})$, la correspondance entre classes de conjugaison stable semi-simples dans $\tilde{G}'(F)$ et $\tilde{G}(F)$ est compliqu\'ee 
 pr\'ecis\'ement parce qu'il y a des classes dans $\tilde{G}'(F)$ qui ne correspondent \`a rien dans $\tilde{G}(F)$. Mais cette correspondance se traduit par une v\'eritable application ${\bf Stab}(\tilde{G}'(F))\to {\bf Stab}(\tilde{G}(F))$ entre ensembles de param\`etres. 
 
 On \'enonce dans la section 3 de nouveaux th\'eor\`emes qui sont pour l'essentiel des reformulations des th\'eor\`emes ci-dessus. Mais les donn\'ees sont cette fois un \'el\'ement ${\cal X}$ de ${\bf Stab}(\tilde{G}(F))$ et un ensemble fini de places $V$. Ces donn\'ees sont ind\'ependantes l'une de l'autre. Pour ${\cal X}$ fix\'e, on peut faire varier $V$. On prouve dans la section 2 des formules de scindage qui entra\^{\i}nent que, si les th\'eor\`emes sont v\'erifi\'es pour $V$ grand (cette notion d\'ependant de ${\cal X}$), alors ils le sont pour tout $V$. Les r\'esultats de cette section 2 s'appuient sur le lemme fondamental pond\'er\'e tordu d\^u \`a Chaudouard et Laumon (cf. [CL] ; f\^acheusement, ces auteurs n'ont pas encore publi\'e la preuve compl\`ete annonc\'ee dans cette r\'ef\'erence). Ils s'appuient aussi sur une version tordue d'un argument d'annulation d\^u \`a Kottwitz. 
 
 Les sections 4 \`a 8 sont consacr\'ees \`a la preuve des th\'eor\`emes de la section 3 pour un param\`etre ${\cal X}$ et un ensemble de places $V$ assez grand. Cette hypoth\`ese sur $V$ entra\^{\i}ne que toutes les distributions intervenant se calculent par la m\'ethode de descente d'Harish-Chandra. Cela nous ram\`ene \`a deux probl\`emes. D'abord, la m\'ethode de descente appliqu\'ee globalement, c'est-\`a-dire sur $F$, fait intervenir une combinatoire compliqu\'ee contr\^ol\'ee par divers groupes de cohomologie. Dans notre situation tordue, il s'agit de cohomologie de complexes de tores. Fort heureusement, cette combinatoire a \'et\'e enti\`erement \'elucid\'ee par Labesse dans une s\'erie d'articles (cf. [Lab1], [Lab2], [Lab3]).  On est alors ramen\'e \`a des probl\`emes similaires aux probl\`emes initiaux, mais dans une situation non tordue et pour des distributions \`a support unipotent. Puisque la situation n'est plus tordue, on peut utiliser les r\'esultats d'Arthur ([A1]). Cela ne suffit toutefois pas car la descente d'Harish-Chandra appliqu\'ee au cas tordu fait in\'evitablement  intervenir un ph\'enom\`ene qui  ne se produit pas dans le cas non tordu: il appara\^{\i}t des triplets endoscopiques non standard. Pour ceux-ci, on utilise le th\'eor\`eme [VI] 5.6 qui permet d'achever la preuve.  Bien s\^ur, ce dernier th\'eor\`eme n'est pas encore prouv\'e. Toutefois, les hypoth\`eses de r\'ecurrence sophistiqu\'ees que l'on a pos\'ees en [VI] 5.8 permettent de l'appliquer, sauf dans quelques cas particuliers. Cette restriction est pr\'ecis\'ement la raison pour laquelle la d\'emonstration du 
  th\'eor\`eme [VI] 5.4 restera incompl\`ete. 
  
  Dans la section 9, nous prouverons le th\'eor\`eme [VI] 5.6.  La preuve est la m\^eme que ci-dessus, mais invers\'ee. Ci-dessus, on d\'eduit le th\'eor\`eme [VI] 5.4 du th\'eor\`eme [VI] 5.6. Maintenant, on d\'eduit le th\'eor\`eme [VI] 5.6 du th\'eor\`eme [VI] 5.4. 
 Les hypoth\`eses de r\'ecurence assurent cette fois que ce th\'eor\`eme est valide pour les triplets $(G,\tilde{G},{\bf a})$ que nous utiliserons. Il faut toutefois prendre garde que, toutes ces d\'emonstrations se faisant par r\'ecurrence, elles ne deviendront de vraies d\'emonstrations que quand tous les pas de la r\'ecurrence auront \'et\'e trait\'es. Comme on l'expliquera davantage en 3.7, il suffira pour cela d'achever la preuve du th\'eor\`eme [VI] 5.4. 
 
 Comme on le voit, notre preuve suit de tr\`es pr\`es celle d'Arthur dans son deuxi\`eme article sur la stabilisation ([A2]). Il s'agissait seulement  d'y ins\'erer les id\'ees de Labesse afin de  l'adapter \`a la situation tordue. 
 
 J'ai re\c{c}u l'aide de C. Moeglin pour la section 2. Je l'en remercie.
 
 \bigskip

 \section{Coefficients globaux et classes de conjugaison stable}
 
 \bigskip
 
 \subsection{Ensemble de param\`etres}
 Pour tout l'article, sauf mention expresse du contraire, $F$ est un corps de nombres, $G$ est un groupe r\'eductif connexe d\'efini sur $F$, $\tilde{G}$ est un espace tordu sous $G$ d\'efini sur $F$ et ${\bf a}$ est un \'el\'ement de 
 $H^1(W_{F};Z(\hat{G}))/ker^1(W_{F};Z(\hat{G}))$. On utilise les d\'efinitions de [VI] 1.1 et on adjoint au triplet $(G,\tilde{G},{\bf a})$ diverses donn\'ees suppl\'ementaires comme dans cette r\'ef\'erence.
 
 Consid\'erons  $\underline{la}$ paire de Borel \'epingl\'ee ${\cal E}^*=(B^*,T^*,(E^*_{\alpha})_{\alpha\in \Delta})$   de $G$. Elle est munie de l'action galoisienne quasi-d\'eploy\'ee, not\'ee $\sigma\mapsto \sigma_{G^*}$, cf. [I] 1.2. Elle est aussi munie d'un automorphisme $\theta^*$ qui commute \`a l'action galoisienne. On note $\Sigma(T^*)$ l'ensemble des racines de $T^*$ dans $\mathfrak{g}$. Pour $\alpha\in \Sigma(T^*)$, on note $\alpha_{res}$ sa restriction \`a $T^{*,\theta^*,0}$ et on pose $\Sigma(T^*)_{res}=\{\alpha_{res}; \alpha\in \Sigma(T^*)\}$. On note $\Sigma_{+}(T^*)$ et $\Sigma_{res,+}(T^*)$ les sous-ensembles positifs d\'etermin\'es par $B^*$. Pour $\alpha\in \Sigma(T^*)$, on note $N\alpha$ la somme des \'el\'ements de l'orbite de $\alpha$ sous l'action du groupe d'automorphismes engendr\'e par $\theta^*$. Ce caract\`ere de $T^*$ se descend en un caract\`ere de $T^*/(1-\theta^*)(T^*)$.
 
 Soit $\mu\in (T^*/(1-\theta^*)(T^*))\times_{{\cal Z}(G)}{\cal Z}(\tilde{G})$ (cf. [I] 1.2 pour la d\'efinition de ${\cal Z}(\tilde{G})$), fixons $\nu\in T^*$ et $\bar{e}\in {\cal Z}(\tilde{G})$ de sorte que $\mu$ soit l'image de $(\nu,\bar{e})$. On note $\Sigma(\mu)$ l'ensemble des $\alpha_{res}$ pour $\alpha\in \Sigma(T^*)$ v\'erifiant l'une des conditions suivantes
 
 - $\alpha$ est de type $1$ ou $2$ et $(N\alpha)(\nu)=1$;
 
 - $\alpha$ est de type $3$ et $(N\alpha)(\nu)=-1$.
 
Cet ensemble s'interpr\`ete de la fa\c{c}on suivante. Identifions ${\cal E}^*$ \`a une paire de Borel \'epingl\'ee particuli\`ere de $G$ et relevons $\bar{e}$ en un \'el\'ement $e\in Z(\tilde{G},{\cal E}^*)$. Posons $\eta=\nu e\in \tilde{G}$, $\bar{G}=G_{\eta}$, $\bar{B}=\bar{G}\cap B^*$ et $\bar{T}=T^{*,\theta^*,0}$.  Alors $\Sigma(\mu)$ est l'ensemble des racines $\Sigma^{\bar{G}}(\bar{T})$ de $\bar{T}$ dans $\bar{\mathfrak{g}}$. En particulier, $\Sigma(\mu)$ est un honn\^ete syst\`eme de racines et $\Sigma_{+}(\mu)=\Sigma(\mu)\cap \Sigma_{res,+}(T^*)$  est le sous-ensemble positif associ\'e \`a $\bar{B}$. On introduit le groupe de Weyl $W(\mu)$ de $\Sigma(\mu)$, c'est-\`a-dire le sous-groupe de $W^{\theta^*}$ engendr\'e par les sym\'etries relatives aux $\alpha_{res}\in \Sigma(\mu)$. C'est aussi le groupe de Weyl $W^{\bar{G}}$ du groupe $\bar{G}$.

Le groupe $W^{\theta^*}$ agit sur $(T^*/(1-\theta^*)(T^*))\times_{{\cal Z}(G)}{\cal Z}(\tilde{G})$ par l'action naturelle sur le premier facteur et  par l'action triviale sur le second. Tous les objets sont  munis de l'action galoisienne quasi-d\'eploy\'ee. Pour  un cocycle $\omega_{\bar{G}}:\Gamma_{F}\to W^{\theta^*}$, consid\'erons les conditions

(1) $\omega_{\bar{G}}(\sigma)\sigma_{G^*}$ fixe $\mu$ pour tout $\sigma\in \Gamma_{F}$;

(2) $\omega_{\bar{G}}(\sigma)\sigma_{G^*}$ fixe $\mu$ et conserve $\Sigma_{+}(\mu)$ pour tout $\sigma\in \Gamma_{F}$.

\noindent Remarquons que la condition (1) implique en tout cas que  $\omega_{\bar{G}}(\sigma)\sigma_{G^*}$ conserve $\Sigma(\mu) $ pour tout $\sigma\in \Gamma_{F}$. On note $\underline{Stab}(\tilde{G}(F))$ l'ensemble des couples $(\mu,\omega_{\bar{G}})$ tels que  $\omega_{\bar{G}}$ v\'erifie (1) et $Stab(\tilde{G}(F))$ celui des couples $(\mu,\omega_{\bar{G}})$ tels que  $\omega_{\bar{G}}$ v\'erifie (2). Disons que deux \'el\'ements $(\mu,\omega_{\bar{G}})$ et $(\mu',\omega_{\bar{G}'})$ de $\underline{Stab}(\tilde{G}(F))$ sont \'equivalents si et seulement si $\mu=\mu'$ et $\omega_{\bar{G}'}(\sigma)\in W(\mu)\omega_{\bar{G}}(\sigma)$ pour tout $\sigma\in \Gamma_{F}$. On a

(3) $Stab(\tilde{G}(F))$ est un ensemble de repr\'esentants des classes d'\'equivalence dans $\underline{Stab}(\tilde{G}(F))$.  

Preuve. Pour $(\mu,\omega_{\bar{G}})\in \underline{Stab}(\tilde{G}(F))$ et $\sigma\in \Gamma_{F}$, il existe un unique $u(\sigma)\in W(\mu)$ tel que $u(\sigma)\omega_{\bar{G}}(\sigma)\sigma_{G^*}$ conserve $\Sigma_{+}(\mu)$. Posons $\omega'_{\bar{G}}(\sigma)=u(\sigma)\omega_{\bar{G}}(\sigma)$. L'unicit\'e de $u(\sigma)$ entra\^{\i}ne facilement que $\omega'_{\bar{G}}$ est encore un cocycle. Alors $(\mu,\omega'_{\bar{G}})$ est un \'el\'ement de $Stab(\tilde{G}(F))$ qui est \'equivalent \`a $(\mu,\omega_{\bar{G}})$. D'autre part, il est imm\'ediat que deux \'el\'ements de $Stab(\tilde{G}(F))$ sont \'equivalent si et seulement s'ils sont \'egaux. D'o\`u (3). $\square$

Le groupe $W^{\theta^*}$ agit  sur $\underline{Stab}(\tilde{G}(F))$ de la fa\c{c}on suivante. Pour $(\mu,\omega_{\bar{G}})\in \underline{Stab}(\tilde{G}(F))$ et $w\in W^{\theta^*}$, l'image de $(\mu,\omega_{\bar{G}})$ par $w$ est le couple $(\mu',\omega_{\bar{G}'})$ d\'efini par $\mu'=w(\mu)$ et $\omega_{\bar{G}'}(\sigma)=w\omega_{\bar{G}}(\sigma)\sigma_{G^*}(w^{-1})$ pour tout $\sigma$. Cette action  respecte la relation d'\'equivalence introduite ci-dessus. Montrons que

(4) soient $(\mu,\omega_{\bar{G}})$ et $(\mu',\omega_{\bar{G}'})$ deux \'el\'ements de $Stab(\tilde{G}(F))$; alors les trois conditions suivantes sont \'equivalentes:

(i) il existe $w\in W^{\theta^*}$ tel que $(\mu',\omega_{\bar{G}'})=w(\mu,\omega_{\bar{G}})$;

(ii) il existe $w\in W^{\theta^*}$ tel que $(\mu',\omega_{\bar{G}'})=w(\mu,\omega_{\bar{G}})$ et $w(\Sigma_{+}(\mu))=\Sigma_{+}(\mu')$;

(iii) il existe $w\in W^{\theta^*}$ tel que $(\mu',\omega_{\bar{G}'})$ et $w(\mu,\omega_{\bar{G}})$  soient \'equivalents dans $\underline{Stab}(\tilde{G}(F))$.

Preuve. Evidemment, (ii) entra\^{\i}ne (i) et (i) entra\^{\i}ne (iii). Soit $w$ v\'erifiant (iii). On a $w(\mu)=\mu'$ donc $w$ envoie $\Sigma(\mu)$ dans $\Sigma(\mu')$. Il existe un unique $u\in W(\mu')$ tel que $uw$ envoie $\Sigma_{+}(\mu)$ dans $\Sigma_{+}(\mu')$.  Il est clair que $(\mu',\omega_{\bar{G}'})$ et $u(\mu',\omega_{\bar{G}'})$ sont \'equivalents. Puisque  $(\mu',\omega_{\bar{G}'})$ et $w(\mu,\omega_{\bar{G}})$  sont \'equivalents et que l'action de $W^{\theta^*}$ conserve l'\'equivalence, $u(\mu',\omega_{\bar{G}'})$ et $uw(\mu,\omega_{\bar{G}})$  sont \'equivalents. Donc aussi $(\mu',\omega_{\bar{G}'})$ et $uw(\mu,\omega_{\bar{G}})$.
 On peut donc remplacer $w$ par $uw$, la condition (iii) reste v\'erifi\'ee et maintenant, $w$ envoie $\Sigma_{+}(\mu)$ dans $\Sigma_{+}(\mu')$. On voit que cette derni\`ere conditon implique que    $w(\mu,\omega_{\bar{G}})$ v\'erifie (2), donc appartient \`a $Stab(\tilde{G}(F))$. Alors $(\mu',\omega_{\bar{G}'})$ et $w(\mu,\omega_{\bar{G}})$  sont deux \'el\'ements \'equivalents de $Stab(\tilde{G}(F))$. Ils sont donc \'egaux et la conclusion de (ii) est v\'erifi\'ee. $\square$
 
   Pour deux \'el\'ements $(\mu,\omega_{\bar{G}}),(\mu',\omega_{\bar{G}'})\in Stab(\tilde{G}(F))$, on dit qu'ils sont conjugu\'es si et seulement s'ils  v\'erifient les trois conditions \'equivalentes ci-dessus. On note ${\bf Stab}(\tilde{G}(F))$ l'ensemble des classes de conjugaison. Cet ensemble est en bijection avec celui des classes de conjugaison par $W^{\theta^*}$ dans $\underline{Stab}(\tilde{G}(F))$. 

Soit $(\mu,\omega_{\bar{G}})\in Stab(\tilde{G}(F))$.  Comme plus haut, relevons $\mu$ en un \'el\'ement $\eta\in \tilde{G}$  et d\'efinissons le groupe $\bar{G}$. Compl\'etons la paire de Borel $(\bar{B},\bar{T})$ de $\bar{G}$ en une paire de Borel \'epingl\'ee. Alors il existe une unique action de $\Gamma_{F}$ sur $\bar{G}$ qui conserve cette paire de Borel \'epingl\'ee,  de sorte que cette action et l'action $\sigma\mapsto \omega_{\bar{G}}(\sigma)\sigma_{G^*}$ co\"{\i}ncident sur $\bar{T}=T^{*,\theta^*,0}$ et induisent la m\^eme action sur $\Sigma_{+}(\mu)$. On note cette action $\sigma\mapsto \sigma_{\bar{G}}$. Pour celle-ci, $\bar{G}$ est quasi-d\'eploy\'e. Le centre $Z(\bar{G})$ est ind\'ependant du rel\`evement $\eta$: c'est
le sous-groupe des $x\in T^{*,\theta^*,0}$ tels que $\alpha_{res}(x)=1$ pour tout $\alpha_{res}\in \Sigma(\mu)$.    L'action galoisienne sur ce centre est $\sigma\mapsto \omega_{\bar{G}}(\sigma)\sigma_{G^*}$. On dit que $(\mu,\omega_{\bar{G}})$ est elliptique si et seulement si on a l'\'egalit\'e $Z(\bar{G})^{\Gamma_{F},0}=Z(G)^{\Gamma_{F},\theta, 0}$ (on note simplement $\theta$ la restriction de $\theta^*$ \`a $Z(G)$, ce qui est justifi\'e par le fait que c'est aussi la restriction de $ad_{g}\circ \theta^*$ pour tout $g\in G$). Cette propri\'et\'e est conserv\'ee par conjugaison. On note $Stab_{ell}(\tilde{G}(F))$ le sous-ensemble des \'el\'ements elliptiques de $Stab(\tilde{G}(F))$ et ${\bf Stab}_{ell}(\tilde{G}(F))$ l'ensemble des classes de conjugaison dans ce sous-ensemble.

 \bigskip
 
 \subsection{Classes de conjugaison stable semi-simples}
 Rappelons que l'on note $\tilde{G}_{ss}$ l'ensemble des \'el\'ements semi-simples de $\tilde{G}$. Pour $\eta\in \tilde{G}_{ss}$, on pose $I_{\eta}=Z(G)^{\theta}G_{\eta}$. Pour $\eta\in \tilde{G}_{ss}(F)$,  on note ${\cal Y}_{\eta}$ l'ensemble des $y\in G(\bar{F})$ tels que $y\sigma(y)^{-1}\in I_{\eta}=I_{\eta}(\bar{F})$ pour tout $\sigma\in \Gamma_{F}$. On dit que deux \'el\'ements $\eta,\eta'\in \tilde{G}_{ss}(F)$ sont stablement conjugu\'es
  si et seulement s'il existe $y\in {\cal Y}_{\eta}$ tel que $y^{-1}\eta y=\eta'$.
     
 Soit $\eta\in \tilde{G}_{ss}(F)$. Fixons une paire de Borel $(B,T)$ de $G$ conserv\'ee par $ad_{\eta}$. Compl\'etons cette paire en une paire de Borel  \'epingl\'ee ${\cal E}=(B,T,(E_{\alpha})_{\alpha\in \Delta})$.  Fixons $e\in Z(\tilde{G},{\cal E})$ et posons $\theta=ad_{e}$. On introduit une cocha\^{\i}ne $u_{{\cal E}}:\Gamma_{F}\to G_{SC}$ de sorte que $ad_{u_{{\cal E}}(\sigma)}\circ\sigma$ conserve ${\cal E}$ pour tout $\sigma\in \Gamma_{F}$. L'action galoisienne naturelle fixe $\eta$ donc conserve $G_{\eta}$. Le couple $(B\cap G_{\eta},T^{\theta,0})$ est une paire de Borel de $G_{\eta}$. On peut choisir une cocha\^{\i}ne $u_{\eta}:\Gamma_{F}\to G_{SC,\eta}$ ($G_{SC,\eta}$ est le groupe des points fixes de $ad_{\eta}$ dans $G_{SC}$)  de sorte que $ad_{u_{\eta}(\sigma)}\circ\sigma$ conserve cette paire de Borel pour tout $\sigma\in \Gamma_{F}$. Posons $v_{\eta}(\sigma)=u_{\eta}(\sigma)u_{{\cal E}}(\sigma)^{-1}$. On a
 
 (1) pour tout $\sigma\in \Gamma_{F}$, $v_{\eta}(\sigma)$ normalise $T$ et son image dans $W$ est fixe par $\theta$.
 
Preuve. On a $ad_{v_{\eta}(\sigma)}=(ad_{u_{\eta}(\sigma)}\circ\sigma)(ad_{u_{{\cal E}}(\sigma)}\circ\sigma)^{-1}$. Les deux facteurs conservent $T$, donc $ad_{v_{\eta}(\sigma)}$ aussi. Pour la seconde assertion, on doit montrer qu'il existe $t(\sigma)\in T$ de sorte que $\theta^{-1}ad_{v_{\eta}(\sigma)}\theta=ad_{t(\sigma)}ad_{v_{\eta}(\sigma)}$. Il suffit de v\'erifier la m\^eme assertion pour chacun des facteurs $ad_{u_{\eta}(\sigma)}\circ\sigma$ et $ad_{u_{{\cal E}}(\sigma)}\circ\sigma$. Pour le deuxi\`eme, c'est clair: il commute \`a $\theta$. Pour le premier, on \'ecrit $\eta=\nu e$ avec $\nu\in T$. Alors $\theta=ad_{e}=ad_{\nu}^{-1}ad_{\eta}$. Le terme $ad_{\eta}$ commute \`a $ad_{u_{\eta}(\sigma)}\circ\sigma$ car $\eta$ est fixe par $\sigma$ et $u_{\eta}(\sigma)\in G_{\eta}$. Le terme $ad_{\nu}^{-1}$ ne commute pas \`a $ad_{u_{\eta}(\sigma)}\circ\sigma$ mais v\'erifie la relation plus faible que l'on souhaite simplement parce que $ad_{u_{\eta}(\sigma)}\circ\sigma$ normalise $T$. $\square$
 
 Ecrivons comme ci-dessus $\eta=\nu e$, avec $\nu \in T$. Notons $\mu_{\eta}$ l'image de $(\nu,e)$ par les applications
 $$T\times Z(\tilde{G},{\cal E})\to (T/(1-\theta)(T))\times_{{\cal Z}(G)}{\cal Z}(\tilde{G},{\cal E})\simeq (T^*/(1-\theta^*)(T^*))\times_{{\cal Z}(G)}{\cal Z}(\tilde{G}).$$
 Pour $\sigma\in \Gamma_{F}$, notons $\omega_{\eta}(\sigma)$ l'image dans $W^{\theta^*}$ de $v_{\eta}(\sigma)$, modulo l'isomorphisme $W^{\theta}\simeq W^{\theta^*}$.
 
  Le couple $(\mu_{\eta},\omega_{\eta})$ ne d\'epend pas du choix de $e$: on ne peut changer $e$ qu'en le multipliant par un \'el\'ement de $Z(G)$ et on voit que la construction est insensible \`a une telle multiplication. De m\^eme, il ne d\'epend pas des choix de cocha\^{\i}nes $u_{{\cal E}}(\sigma)$ et $u_{\eta}(\sigma)$.  Il ne d\'epend pas de toute la paire de Borel \'epingl\'ee ${\cal E}$ mais seulement de la paire de Borel sous-jacente $(B,T)$. En effet, si on change seulement l'\'epinglage, sans changer la paire de Borel sous-jacente, $\mu_{\eta}$ ne change pas (cf. [I] preuve de 1.10(1)) et les \'el\'ements $u_{{\cal E}}(\sigma)$ et $u_{\eta}(\sigma)$ sont multipli\'es par des \'el\'ements de $T$, ce qui ne change pas $\omega_{\eta}$.
 
 \ass{Proposition}{(i) Pour $\eta\in \tilde{G}_{ss}(F)$, le couple $(\mu_{\eta},\omega_{\eta})$ construit ci-dessus appartient \`a $Stab(\tilde{G}(F))$.
 
 (ii) L'image de ce couple dans ${\bf Stab}(\tilde{G}(F))$ ne d\'epend pas de la paire $(B,T)$ utilis\'ee dans sa construction.
 
 (iii) Soient $\eta,\eta'\in \tilde{G}_{ss}(F)$ et  $(B,T)$, $(B',T')$ des paires de Borel conserv\'ees respectivement par $ad_{\eta}$ et $ad_{\eta'}$. Les couples $(\mu_{\eta},\omega_{\eta})$  et $(\mu_{\eta'},\omega_{\eta'})$ d\'eduits de ces donn\'ees sont \'egaux si et seulement s'il existe $y\in {\cal Y}_{\eta}$ tel que $\eta'=y^{-1}\eta y$ et $(B',T')=ad_{y^{-1}}(B,T)$.  
  
 (iv) Soient $\eta,\eta'\in \tilde{G}_{ss}(F)$ et  $(B,T)$, $(B',T')$ des paires de Borel conserv\'ees respectivement par $ad_{\eta}$ et $ad_{\eta'}$. Les couples $(\mu_{\eta},\omega_{\eta})$  et $(\mu_{\eta'},\omega_{\eta'})$ d\'eduits de ces donn\'ees  ont m\^eme image dans ${\bf Stab}(\tilde{G}(F))$ si et seulement si $\eta$ et $\eta'$ sont stablement conjugu\'es. }
 
 Preuve. Pour d\'emontrer (i), on peut identifier ${\cal E}$ \`a ${\cal E}^*$. On a alors $\sigma_{G^*}=ad_{u_{{\cal E}}(\sigma)}\circ\sigma$ pour tout $\sigma\in \Gamma_{F}$. Soulignons que cette formule d\'efinit une action galoisienne sur $G$ mais pas forc\'ement sur $\tilde{G}$ car l'action par conjugaison du cobord de la cocha\^{\i}ne $u_{{\cal E}}$ peut ne pas \^etre triviale. Pour cette raison, on n'utilise la notation $\sigma_{G^*}$ que pour l'action sur $G$.
  
On doit montrer que $\omega_{\eta}(\sigma)\sigma_{G^*}$ conserve $\mu_{\eta}$ pour tout $\sigma$. L'\'el\'ement $\omega_{\eta}(\sigma)\sigma_{G^*}(\mu_{\eta})$ est l'image dans $(T/(1-\theta)(T))\times_{{\cal Z}(G)}{\cal Z}(\tilde{G},{\cal E})$ de
 $$(2) \qquad (ad_{v_{\eta}(\sigma)}\circ\sigma_{G^*}(\nu),ad_{u_{{\cal E}}(\sigma)}\circ\sigma(e)).$$
 On sait que $\sigma_{G^*}$ conserve $Z(\tilde{G},{\cal E}^*)$, donc il existe $z(\sigma)\in Z(G)$ tel que $ad_{u_{{\cal E}}(\sigma)}\circ\sigma(e)=z(\sigma)^{-1}e$.  Parce que $\omega_{\eta}(\sigma)\in W^{\theta}$, on peut choisir $n(\sigma)\in G_{e}$ qui le rel\`eve et $t(\sigma)\in T$ de sorte que $v_{\eta}(\sigma)=t(\sigma)n(\sigma)$. On a alors $ad_{v_{\eta}(\sigma)}\circ ad_{u_{{\cal E}}(\sigma)}\circ\sigma(e)=t(\sigma)\theta(t(\sigma))^{-1}z(\sigma)^{-1}e$. Parce que $\eta$ est fixe par $\Gamma_{F}$ et que $u_{\eta}(\sigma)\in G_{\eta}$, on a $ad_{u_{\eta}(\sigma)}\circ\sigma(\eta)=\eta$, ou encore $ad_{v_{\eta}(\sigma)}\circ ad_{u_{{\cal E}}(\sigma)}\circ\sigma(\nu e)=\nu e$. Les deux relations pr\'ec\'edentes entra\^{\i}nent 
 $$ad_{v_{\eta}(\sigma)}\circ\sigma_{G^*}(\nu)=t(\sigma)^{-1}\theta(t(\sigma))z(\sigma)\nu.$$
Mais alors le couple (2) a bien pour image $\mu_{\eta}$ dans  $(T/(1-\theta)(T))\times_{{\cal Z}(G)}{\cal Z}(\tilde{G},{\cal E})$.
 
 On doit montrer que $\omega_{\eta}$ est un cocycle. L'application qui, \`a $\sigma\in \Gamma_{F}$, associe l'automorphisme $ad_{u_{\eta}(\sigma)}\circ\sigma$ de $T$, est un homomorphisme: son cobord est donn\'e par des automorphismes int\'erieurs de $G_{\eta}$ qui conservent $(B\cap G_{\eta},T^{\theta,0})$, donc commutent \`a $T$. Mais l'automorphisme $\omega_{\eta}(\sigma)\sigma_{G^*}$ de $T$ est \'egal \`a $ad_{u_{\eta}(\sigma)}\circ\sigma$. Donc $\sigma\mapsto \omega_{\eta}(\sigma)\sigma_{G^*}$ est un homomorphisme de $\Gamma_{F}$ \`a valeurs dans le groupe d'automorphismes de $T$. Cela signifie que  $\omega_{\eta}$ est un cocycle.

 De m\^eme, $\omega_{\eta}(\sigma)\sigma_{G^*}$ conserve $\Sigma_{+}(\mu_{\eta})$ parce que $ad_{u_{\eta}(\sigma)}\circ\sigma$ conserve $B\cap G_{\eta}$. Cela prouve (i).
 
 Prouvons (ii).  Changeons  la paire $(B,T)$. D'apr\`es [I] 1.3(2), on ne peut que la remplacer par une paire $ad_{x}\circ w(B,T)$, o\`u $x\in G_{\eta}$ et $w\in W^{\theta}$. On sait que $G_{\eta}\subset Z(G)G_{SC,\eta}$ et que tout \'el\'ement de $W^{\theta}$ se rel\`eve en un \'el\'ement de $G_{SC,e}$.  On est ramen\'e \`a voir ce qui se passe quand on remplace ${\cal E}$ par $ad_{x}{\cal E}$, avec ou bien $x\in G_{SC,\eta}$, ou bien $x\in G_{SC,e}$ et  normalise $T$. Comme dans la preuve de (i), on peut supposer ${\cal E}={\cal E}^*$. Dans les deux cas, on doit construire les objets relatifs \`a $\underline{{\cal E}}^*=ad_{x}({\cal E}^*)$, notons-les en les soulignant, puis les ramener \`a ${\cal E}^*$ par l'isomorphisme canonique $ad_{x}^{-1}$. Dans le cas o\`u $x\in G_{SC,\eta}$, on a $\eta=ad_{x}(\eta)=ad_{x}(\nu)ad_{x}(e)$ et on peut prendre $\underline{\nu}=ad_{x}(\nu)$, $\underline{e}=ad_{x}(e)$. On peut aussi choisir $u_{\underline{{\cal E}}^*}(\sigma)=xu_{{\cal E}^*}(\sigma)\sigma(x)^{-1}$ et $\underline{u}_{\eta}(\sigma)=xu_{\eta}(\sigma)\sigma(x)^{-1}$. D'o\`u $\underline{v}_{\eta}(\sigma)=xv_{\eta}(\sigma)x^{-1}$. En ramenant les objets $\underline{\nu}$, $\underline{e}$ et $\underline{v}_{\eta}$ par l'isomorphisme $ad_{x}^{-1}$, on retrouve les objets initiaux, rien n'a chang\'e. Dans le cas o\`u $x\in G_{SC,e}$ et $x$ normalise $T$, on a $e=ad_{x}(e)\in Z(\tilde{G},\underline{{\cal E}}^*)$. On peut prendre $\underline{\nu}=\nu$, $\underline{e}=e$. On peut encore prendre $u_{\underline{{\cal E}}^*}(\sigma)=xu_{{\cal E}^*}(\sigma)\sigma(x)^{-1}$. En posant $\underline{B}=xB^*x^{-1}$, le couple $(\underline{B}\cap G_{\eta},T^{*,\theta^*,0})$ est une paire de Borel de $G_{\eta}$ et on peut fixer $y\in G_{SC,\eta}$ tel que cette paire soit \'egale \`a $ad_{y}(B^*\cap G_{\eta},T^{*,\theta^*,0})$. On peut alors prendre $\underline{u}_{\eta}(\sigma)=yu_{\eta}(\sigma)\sigma(y)^{-1}$. Alors
 $$\underline{v}_{\eta}(\sigma)=yu_{\eta}(\sigma)\sigma(y^{-1}x)u_{{\cal E}^*}(\sigma)^{-1}x^{-1}.$$
 Ramenons les \'el\'ements $\underline{\nu}$, $\underline{e}$ et $\underline{v}_{\eta}$ par l'isomorphisme $ad_{x}^{-1}$. On obtient respectivement les \'el\'ements $\nu'=x^{-1}\nu x$, $e'=e$, et
 $$v'_{\eta}(\sigma)=x^{-1}yu_{\eta}(\sigma)\sigma(y^{-1}x)u_{{\cal E}^*}(\sigma)^{-1}$$
 $$=x^{-1}yv_{\eta}(\sigma)ad_{u_{{\cal E}^*}(\sigma)}\circ\sigma(y^{-1}x).$$
Les \'el\'ements $x$ et $y$ normalisent $T^*$. Leurs images dans $W$ sont fixes par $\theta^*$, par hypoth\`ese pour $x$, parce que $y\in G_{\eta}$ pour $y$. Donc $x^{-1}y$ d\'efinit un \'el\'ement $w\in W^{\theta^*}$. En notant $\omega'_{\eta}(\sigma)$ l'image de $v'_{\eta}(\sigma)$ dans $W$, on obtient
$$\omega'_{\eta}(\sigma)=w\omega_{\eta}(\sigma)\sigma_{G^*}(w)^{-1}.$$
Parce que $y$ normalise $T^*$ et appartient \`a $ G_{\eta}$, un calcul d\'ej\`a fait plusieurs fois montre que $y\nu y^{-1}$ a m\^eme image que $\nu $ dans $T^*/(1-\theta^*)(T^*)$. Donc $\nu '=x^{-1}\nu x$ a m\^eme image que $w(\nu )$. Mais alors le couple $(\mu'_{\eta},\omega'_{\eta})$ construit \`a l'aide de $\underline{{\cal E}}$ est $(w(\mu_{\eta}),\omega'_{\eta})$, o\`u $\omega'_{\eta}$ est comme ci-dessus. Ce couple est conjugu\'e \`a $(\mu_{\eta},\omega_{\eta})$, ce qui prouve (ii). 

{\bf Remarque.} En inversant la preuve de (ii), on obtient le r\'esultat suivant. 

(3) Soit $(\mu_{\eta},\omega_{\eta})$ le couple associ\'e \`a $\eta$ et \`a une paire de Borel $(B,T)$.  Soit $ (\mu',\omega_{\bar{G}'})\in Stab(\tilde{G}(F))$ un couple conjugu\'e \`a $(\mu_{\eta},\omega_{\eta})$. Alors il existe une autre paire de Borel   $(\underline{B},\underline{T})$ conserv\'ee par $ad_{\eta}$ de sorte que 

- $(\mu',\omega_{\bar{G}'})$ soit le couple associ\'e \`a $\eta$ et \`a cette paire $ (\underline{B},\underline{T})$;

- $\underline{T}=T$ et $\underline{B}\cap G_{\eta}=B\cap G_{\eta}$.

Preuve.  On introduit ${\cal E}$ et $e$ comme au d\'ebut du paragraphe et on identifie ${\cal E}$  \`a ${\cal E}^*$. D'apr\`es 1.1(4)(ii), on peut fixer $w\in W^{\theta^*}$ de sorte que $(\mu',\omega_{\bar{G}'})=w(\mu_{\eta},\omega_{\eta})$ et $w(\Sigma_{+}(\mu_{\eta}))=\Sigma_{+}(\mu')$. Relevons l'\'el\'ement $w^{-1}$ en un \'el\'ement $x\in G_{SC,e}$ qui normalise $T=T^*$. On pose $\underline{{\cal E}}=ad_{x}({\cal E})$. La preuve de (ii)  montre que cette paire v\'erifie les propri\'et\'es requises, pourvu que l'on montre que $\underline{B}\cap G_{\eta}=B\cap G_{\eta}$ (cette propri\'et\'e entra\^{\i}ne que l'on peut prendre $y=1$ donc le $w$ que l'on vient de choisir est le m\^eme que plus haut). Or l'\'egalit\'e $w(\Sigma_{+}(\mu_{\eta}))=\Sigma_{+}(\mu')$ entra\^{\i}ne que les racines de $T^*$ dans le radical unipotent de $\underline{B}\cap G_{\eta}$ sont positives pour $B=B^*$. Autrement dit $\underline{B}\cap G_{\eta}\subset B$. D'o\`u forc\'ement  $\underline{B}\cap G_{\eta}=B\cap G_{\eta}$. Cela prouve (3).

 \bigskip

Prouvons (iii). Soient $\eta$, $\eta'$, $(B,T)$ et $(B',T')$ comme en (iii). Supposons qu'il existe $y\in {\cal Y}_{\eta}$ tel que $\eta'=y^{-1}\eta y$ et $(B',T')=ad_{y^{-1}}(B,T)$.  On fixe  un tel $y$ et on  le d\'ecompose  en $y=y_{sc}z$, avec $y_{sc}\in G_{SC}$ et $z\in Z(G)$ (on identifie ici $y_{sc}$ \`a son image   dans $G$). On compl\`ete $(B,T)$  en une paire de Borel \'epingl\'ee   que l'on identifie \`a ${\cal E}^*$. On compl\`ete $(B',T')$ en la paire de Borel \'epingl\'ee  ${\cal E}'=ad_{y_{sc}^{-1}}({\cal E}^*)$. On affecte d'un $'$ les objets relatifs \`a $\eta'$ et ${\cal E}'$. On a $\eta'=y^{-1}\eta y=y_{sc}^{-1}\nu y_{sc}z^{-1}\theta(z)y_{sc}^{-1}ey_{sc}$. On peut prendre $\nu'=y_{sc}^{-1}\nu y_{sc}z^{-1}\theta(z)$ et $e'=y_{sc}^{-1}ey_{sc}$. On peut  prendre  $u_{{\cal E}'}(\sigma)=y_{sc}^{-1}u_{{\cal E}^*}(\sigma)\sigma(y_{sc})$. Puisque $ad_{u_{\eta}(\sigma)}\circ\sigma$ conserve $(B^*\cap G_{\eta},T^{*,\theta^*,0})$, l'automorphisme $ad_{y^{-1}u_{\eta}(\sigma)}\circ\sigma \circ ad_{y}$ conserve $(B'\cap G_{\eta'},T'\cap G_{\eta'})$.  On a
$$ad_{y^{-1}u_{\eta}(\sigma)}\circ\sigma \circ ad_{y}=ad_{y^{-1}u_{\eta}(\sigma)\sigma(y)}\circ \sigma.$$
Puisque $\sigma(y)\in I_{\eta}y$, on peut \'ecrire $\sigma(y)=z(\sigma)x(\sigma)y$, avec $z(\sigma)\in Z(G)$ et $x(\sigma)\in G_{SC,\eta}$. Alors
$$ad_{y^{-1}u_{\eta}(\sigma)}\circ\sigma \circ ad_{y}=ad_{y^{-1}u_{\eta}(\sigma)x(\sigma)y}\circ \sigma.$$
L'\'el\'ement $y^{-1}u_{\eta}(\sigma)x(\sigma)y=y_{sc}^{-1}u_{\eta}(\sigma)x(\sigma)y_{sc}$ appartient \`a $G_{SC,\eta'}$ et on peut choisir $u_{\eta'}(\sigma)=y_{sc}^{-1}u_{\eta}(\sigma)x(\sigma)y_{sc}$. D'o\`u
$$v_{\eta'}(\sigma)=y_{sc}^{-1}u_{\eta}(\sigma)x(\sigma)y_{sc}\sigma(y_{sc})^{-1}u_{{\cal E}^*}(\sigma)^{-1}y_{sc}^{-1}.$$
On a
$$x(\sigma)y_{sc}\sigma(y_{sc})^{-1}=z(\sigma)^{-1}\sigma(y)y^{-1}y_{sc}\sigma(y_{sc})^{-1}.$$
Ceci est un \'el\'ement de $Z(G)$, notons-le $\zeta(\sigma)$. Alors 
$$v_{\eta'}(\sigma)=\zeta(\sigma)y_{sc}^{-1}v_{\eta}(\sigma)y_{sc}.$$
Le couple $(\mu_{\eta'},\omega_{\eta'})$ se d\'eduit de $\nu'$, $e'$, $v_{\eta'}$ en ramenant ces objets \`a ${\cal E}^*$ par l'isomorphisme $ad_{y_{sc}}$. Autrement dit, $\mu_{\eta'}$ est l'image naturelle de $(\nu z^{-1}\theta(z),e)$ et $\omega_{\eta'}(\sigma)$ est l'image dans $W$ de $\zeta(\sigma)v_{\eta}(\sigma)$. Ces images ne sont autres que $\mu_{\eta}$ et $\omega_{\eta}$. Cela prouve que  $(\mu_{\eta'},\omega_{\eta'})=(\mu_{\eta},\omega_{\eta})$. 

Inversement, soient $\eta$, $\eta'$, $(B,T)$ et $(B',T')$ comme en (iii) et supposons que $(\mu_{\eta'},\omega_{\eta'})=(\mu_{\eta},\omega_{\eta})$. On compl\`ete $(B,T)$ en une paire de Borel \'epingl\'ee que l'on peut supposer \^etre \'egale \`a ${\cal E}^*$ et on compl\`ete $(B',T')$ en une paire de Borel \'epingl\'ee ${\cal E}'$. On note comme plus haut les termes relatifs \`a $\eta$. 
   On choisit $d_{sc}\in G_{SC}$ tel que   $ad_{d_{sc}}({\cal E}')={\cal E}^*$. L'\'el\'ement $e'=d_{sc}^{-1}ed_{sc}$ appartient \`a $Z(\tilde{G},{\cal E}')$. On peut \'ecrire $\eta'=\nu 'e'$, avec $\nu'\in T'$. L'hypoth\`ese $\mu_{\eta'}=\mu_{\eta}$ signifie que $\nu$ et $d_{sc}\nu'd_{sc}^{-1}$ ont m\^eme image dans $T^*/(1-\theta^*)(T^*)$. On peut choisir $\tau\in T^*$ de sorte que $\tau d_{sc}\nu'd_{sc}^{-1}\theta^*(\tau)^{-1}=\nu$. Posons $y=\tau d_{sc}$. Alors $y^{-1}\eta y=\eta'$ et $ad_{y^{-1}}(B,T)=ad_{y^{-1}}(B^*,T^*)=(B',T')$.  On peut choisir $u_{{\cal E}'}(\sigma)=d_{sc}^{-1}u_{{\cal E}^*}(\sigma)\sigma(d_{sc})$. Choisissons  une cocha\^{\i}ne $u_{\eta'}$ relative \`a $\eta'$ et ${\cal E}'$. Alors
$$v_{\eta'}(\sigma)=u_{\eta'}(\sigma)\sigma(d_{sc})^{-1}u_{{\cal E}^*}(\sigma)^{-1}d_{sc}.$$
L'hypoth\`ese $\omega_{\eta'}=\omega_{\eta}$ signifie que les deux \'el\'ements $d_{sc}v_{\eta'}(\sigma)d_{sc}^{-1}$ et $v_{\eta}(\sigma)$ ont m\^eme image dans $W$. Autrement dit,  on a la relation
$$d_{sc}u_{\eta'}(\sigma)\sigma(d_{sc})^{-1}u_{{\cal E}^*}(\sigma)^{-1}\in T^* u_{\eta}(\sigma)u_{{\cal E}^*}(\sigma)^{-1},$$
ou encore
$$d_{sc}u_{\eta'}(\sigma)\sigma(d_{sc})^{-1}\in T^*u_{\eta}(\sigma).$$
On peut multiplier cette relation \`a gauche par $\tau$ et \`a droite par $\sigma(\tau)^{-1}$. Parce que $u_{\eta}(\sigma)\circ\sigma$ conserve $T^*$, on obtient
$$yu_{\eta'}(\sigma)\sigma(y)^{-1}\in T^*u_{\eta}(\sigma),$$
ou encore
$$yu_{\eta'}(\sigma)y^{-1} y\sigma(y)^{-1}\in T^*u_{\eta}(\sigma).$$
Les \'el\'ements $yu_{\eta'}(\sigma)y^{-1}$ et $u_{\eta}(\sigma)$ appartiennent \`a $G_{\eta}$. Parce que $y^{-1}\eta y=\eta'$, l'\'el\'ement $y\sigma(y)^{-1}$ appartient au commutant $Z_{G}(\eta)$. On peut donc remplacer dans la relation ci-dessus le tore $T^*$ par son intersection avec $Z_{G}(\eta)$, autrement dit par $T^{*,\theta^*}$. Mais ce groupe est \'egal \`a $Z(G)^{\theta}T^{*,\theta^*,0}$, donc contenu dans $I_{\eta}$. La relation pr\'ec\'edente entra\^{\i}ne alors que $y\sigma(y)^{-1}\in I_{\eta}$. Donc  $y\in {\cal Y}_{\eta}$, ce qui ach\`eve la preuve de  (iii).

Prouvons (iv). Soient $\eta$, $\eta'$, $(B,T)$ et $(B',T')$ comme en (iv). Supposons $\eta$ et $\eta'$ stablement conjugu\'es, fixons $y\in {\cal Y}_{\eta}$ tel que $\eta'=y^{-1}\eta y$. La paire $ad_{y^{-1}}(B,T)$ est conserv\'ee par $ad_{\eta'}$. On peut remplacer $(B',T') $ par cette paire $ad_{y^{-1}}(B,T)$ puisque, d'apr\`es (ii),  cela ne change pas l'image de $(\mu_{\eta'},\omega_{\eta'})$ dans ${\bf Stab}(\tilde{G}(F))$. Mais alors, d'apr\`es (iii), on a $(\mu_{\eta'},\omega_{\eta'})=(\mu_{\eta},\omega_{\eta})$. A fortiori, ces termes ont m\^eme image dans ${\bf Stab}(\tilde{G}(F))$. Inversement, supposons que $(\mu_{\eta},\omega_{\eta})$ et $(\mu_{\eta'},\omega_{\eta'})$ ont m\^eme image dans  ${\bf Stab}(\tilde{G}(F))$. D'apr\`es (3), on peut changer la paire $(B',T')$ de sorte que l'on ait l'\'egalit\'e $(\mu_{\eta'},\omega_{\eta'})=(\mu_{\eta},\omega_{\eta})$. Alors (iii) implique que $\eta$ et $\eta'$ sont stablement conjugu\'es. $\square$

Notons $\tilde{G}_{ss}(F)/st-conj$ l'ensemble des classes de conjugaison stable semi-simples dans $\tilde{G}(F)$. La proposition nous fournit une application que l'on note
$$\chi^{\tilde{G}}:\tilde{G}_{ss}(F)/st-conj \to {\bf Stab}(\tilde{G}(F)).$$
Celle-ci est injective d'apr\`es le (iv) de la proposition. 

Soit ${\cal X}\in {\bf Stab}(\tilde{G}(F))$. Modulo divers choix, on lui associe un groupe $\bar{G}$ comme en 1.1. En introduisant les tores "standard" $\hat{T}$ et $\hat{\bar{T}}$ des groupes duaux $\hat{G}$ et $\hat{\bar{G}}$, on peut identifier $\hat{\bar{T}}$ \`a $\hat{T}/(1-\hat{\theta})(\hat{T})$, muni d'une action galoisienne $\sigma\mapsto \omega_{\bar{G}}(\sigma)\circ \sigma_{G^*}$. Un calcul de syst\`eme de racines d\'ej\`a fait plusieurs fois montre que $Z(\hat{G})$ s'envoie dans $Z(\hat{\bar{G}})$. Ainsi, la donn\'ee ${\bf a}$ se pousse en un \'el\'ement de $H^1(W_{F};Z(\hat{\bar{G}}))/ker^1(W_{F};Z(\hat{\bar{G}}))$. Celui-ci d\'etermine un caract\`ere automorphe de $\bar{G}({\mathbb A}_{F})$. On voit que la paire form\'ee du groupe $\bar{G}$ et de ce caract\`ere ne d\'epend des choix qu'\`a isomorphisme pr\`es.

Soient $\eta\in \tilde{G}_{ss}(F)$ et $(B,T)$ une paire de Borel conserv\'ee par $ad_{\eta}$. On peut utiliser  l'\'el\'ement  $\eta$ pour construire le groupe $\bar{G}$ associ\'e comme en 1.1 au couple $(\mu_{\eta},\omega_{\eta})$. Ce groupe n'est autre qu'une forme quasi-d\'eploy\'ee du   groupe $G_{\eta}$.   Le caract\`ere automorphe ci-dessus de $\bar{G}({\mathbb A}_{F})$ se transf\`ere en un caract\`ere automorphe de $G_{\eta}({\mathbb A}_{F})$, qui n'est autre que
 la restriction de $\omega$ \`a ce groupe.  

Rappelons que l'on dit qu'un \'el\'ement $\eta\in \tilde{G}_{ss}(F)$ est elliptique s'il v\'erifie les conditions \'equivalentes

(4) il existe un sous-tore tordu $\tilde{T}$ de $\tilde{G}$ d\'efini sur $F$ et elliptique (c'est-\`a-dire tel que $T^{\Gamma_{F},\theta,0}=Z(G)^{\Gamma_{F},\theta,0}$) tel que $\eta\in \tilde{T}(F)$;

(5) on a l'\'egalit\'e $Z(G_{\eta})^{\Gamma_{F},0}=Z(G)^{\Gamma_{F},\theta,0}$.

On note  $\tilde{G}(F)_{ell}$ l'ensemble des \'el\'ements elliptiques de $\tilde{G}_{ss}(F)$. 
Cet ensemble est invariant par conjugaison stable. On voit qu'une classe de conjugaison stable semi-simple est form\'ee d'\'el\'ements elliptiques si et seulement si son image par $\chi^{\tilde{G}}$ appartient \`a ${\bf Stab}_{ell}(\tilde{G}(F))$.

\bigskip

\subsection{Le cas quasi-d\'eploy\'e \`a torsion int\'erieure}
\ass{Lemme}{Supposons $(G,\tilde{G},{\bf a})$ quasi-d\'eploy\'e et \`a torsion int\'erieure.  

(i) L'application $\chi^{\tilde{G}}$ est bijective.

(ii) Plus pr\'ecis\'ement, pour $(\mu,\omega_{\bar{G}})\in Stab(\tilde{G}(F))$, il existe $\epsilon\in \tilde{G}_{ss}(F)$ et une paire de Borel  $(B_{\epsilon},T_{\epsilon} )$ tels que

- $ad_{\epsilon}$ conserve $(B_{\epsilon},T_{\epsilon})$;

- en utilisant la paire  $(B_{\epsilon},T_{\epsilon})$ dans la construction de 1.2, on a l'\'egalit\'e  $(\mu_{\epsilon},\omega_{\epsilon})=(\mu,\omega_{\bar{G}})$;

- $G_{\epsilon}$ est quasi-d\'eploy\'e et la paire de Borel $(B_{\flat}^*,T_{\epsilon})=(G_{\epsilon}\cap B_{\epsilon},T_{\epsilon})$ de ce groupe est d\'efinie sur $F$.}

Preuve. On peut identifier ${\cal E}^*$ \`a une paire de Borel \'epingl\'ee de $G$ d\'efinie sur $F$. L'action galoisienne quasi-d\'eploy\'ee n'est autre que l'action naturelle. Soit $(\mu,\omega_{\bar{G}})\in Stab(\tilde{G}(F))$. Le terme $\omega_{\bar{G}}$ est un cocycle de $\Gamma_{F}$ dans $W$. D'apr\`es [K1] corollaire 1.2, on peut fixer $g_{sc}\in G_{SC}$ tel que, pour tout $\sigma\in \Gamma_{F}$, $g_{sc}\sigma(g_{sc})^{-1}$ normalise $T^*$ et ait $\omega_{\bar{G}}(\sigma)$ pour image dans $W$. Le terme $\mu$ est un \'el\'ement de $T^*\times_{Z(G)}Z(\tilde{G},{\cal E}^*)$ et cet ensemble n'est autre que le centralisateur $\tilde{T}^*$ de $T^*$ dans $\tilde{G}$. Posons $\eta=g_{sc}\mu g_{sc}^{-1}$. Parce que $\omega_{\bar{G}}(\sigma)\sigma$ fixe $\mu$ pour tout $\sigma$, on voit que $\eta\in \tilde{G}(F)$. C'est un \'el\'ement semi-simple. On pose $(B_{\epsilon},T_{\epsilon})=ad_{g_{sc}}(B^*,T^*)$. On v\'erifie imm\'ediatement qu'en utilisant cette paire, la construction de 1.2 envoie $\epsilon$ sur $(\mu,\omega_{\bar{G}})$. En utilisant les notations de l'\'enonc\'e, on a $B_{\flat}^*=ad_{g_{sc}}(B^*\cap G_{\mu})$ et $T_{\epsilon}=ad_{g_{sc}}(T^*)$. Parce que $\omega_{\bar{G}}(\sigma)\sigma$ conserve $\Sigma_{+}(\mu)$ pour tout $\sigma$, l'automorphisme $ad_{g_{sc}\sigma(g_{sc})^{-1}}\circ \sigma$ conserve $(B^*\cap G_{\mu},T^*)$. Donc $\sigma$ conserve $(B_{\flat}^*,T_{\epsilon})$.   $\square$

\bigskip

\subsection{Le cas local}
Dans les trois paragraphes pr\'ec\'edents, le corps de base $F$ \'etait notre corps de nombres. En fait, on peut le remplacer par un corps local de caract\'eristique nulle, tout reste vrai \`a l'exception suivante pr\`es.

{\bf Exception.} C'est la derni\`ere assertion de 1.2 dans le cas o\`u $F$ est archim\'edien. Dans ce cas, les conditions (4) et (5) de 1.2 ne sont plus \'equivalentes. 
Dans les articles pr\'ec\'edents,  on a choisi la  condition (4) pour d\'efinir l'ellipticit\'e. Mais la derni\`ere assertion de 1.2 n'est vraie que  si on utilise la condition (5).
\bigskip

 En particulier, pour une place $v$ de notre corps de nombres $F$, on d\'efinit les ensembles $Stab(\tilde{G}(F_{v}))$ et ${\bf Stab}(\tilde{G}(F_{v}))$. Il y a une application naturelle de localisation $ Stab(\tilde{G}(F))\to  Stab(\tilde{G}(F_{v}))$: \`a $(\mu,\omega_{\bar{G}})$, on associe $(\mu,\omega_{\bar{G}_{v}})$, o\`u $\omega_{\bar{G}_{v}}$ est la restriction de $\omega_{\bar{G}}$ \`a $\Gamma_{F_{v}}$. Cette application se quotiente en une application ${\bf Stab}(\tilde{G}(F))\to {\bf Stab}(\tilde{G}(F_{v}))$. Le diagramme suivant est commutatif
$$\begin{array}{ccc}\tilde{G}_{ss}(F)/st-conj&\to&\tilde{G}_{ss}(F_{v})/st-conj\\ \chi^{\tilde{G}}\downarrow&&\downarrow \chi^{\tilde{G}_{v}}\\ {\bf Stab}(\tilde{G}(F))&\to&{\bf Stab}(\tilde{G}(F_{v}))\\ \end{array}$$
On peut remplacer la place $v$ par un ensemble fini $V$ de places de $F$: on d\'efinit ${\bf Stab}(\tilde{G}(F_{V}))=\prod_{v\in V}{\bf Stab}(\tilde{G}(F_{v}))$ et on a des propri\'et\'es analogues.

\bigskip

\subsection{Rappels sur le cas local non ramifi\'e}
Fixons une place finie $v\not\in V_{ram}$. Nous allons d'abord fixer les notations qui seront utilis\'ees dans toute la suite de l'article. On note $\mathfrak{o}_{v}$ l'anneau des entiers de $F_{v}$, $\mathfrak{o}_{v}^{\times}$ le groupe d'unit\'es et ${\mathbb F}_{v}$ le corps r\'esiduel.  On note $\bar{\mathfrak{o}}_{v}$, $\bar{\mathfrak{o}}_{v}^{\times}$ et $\bar{{\mathbb F}}_{v}$, resp. $\mathfrak{o}_{v}^{nr}$, $\mathfrak{o}_{v}^{nr ,\times}$, $\bar{{\mathbb F}}_{v}$, les objets analogues pour la cl\^oture alg\'ebrique $\bar{F}_{v}$, resp. pour la plus grande extension  non-ramifi\'ee $F^{nr}_{v}$ contenue dans $\bar{F}_{v}$.  On pose $\Gamma_{v}^{nr}=Gal(F_{v}^{nr}/F_{v})$ et on note $I_{v}\subset \Gamma_{F_{v}}$ le groupe d'inertie. On a aussi $I_{v}\subset W_{F_{v}}$ et on pose $W_{v}^{nr}=W_{F_{v}}/I_{v}$. On note $p$ la caract\'eristique de ${\mathbb F}_{v}$.  En [VI] 1.1, on a fix\'e un sous-groupe compact  hypersp\'ecial $K_{v}$ de $G(F_{v})$ et un sous-espace hypersp\'ecial $\tilde{K}_{v}$ de $\tilde{G}(F_{v})$.  Au groupe $K_{v}$ est associ\'e un sch\'ema en groupes lisse  ${\cal K}_{v}$ sur $\mathfrak{o}_{v}$. On note $K_{v}^{nr}={\cal K}_{v}(\mathfrak{o}_{v}^{nr})$. Si $E$ est une extension finie de $F$ non ramifi\'ee en $v$ et si $w$ est une place de $E$ au-dessus de $v$, on utilise les notations $\mathfrak{o}_{w}$ etc... et $K_{w}={\cal K}_{v}(\mathfrak{o}_{w})$. Si $E$ est une extension non ramifi\'ee de $F_{v}$, on utilisera plut\^ot les notations $\mathfrak{o}_{E}$ etc... et $K_{v}^{E}={\cal K}_{v}(\mathfrak{o}_{E})$. 
Le groupe $K_{v}$ d\'etermine des sous-groupes compacts hypersp\'eciaux des groupes $G_{SC}$, $G_{AD}$ et $G_{\sharp}=G/Z(G)^{\theta}$. On les note $K_{sc,v}$, $K_{ad,v}$, $K_{\sharp,v}$.

On se rappelle que le groupe $K_{v}$ est issu d'une paire de Borel \'epingl\'ee de $G$ d\'efinie sur $F_{v}$. On fixe une telle paire ${\cal E}_{0}=(B_{0},T_{0},(E_{0,\alpha})_{\alpha\in \Delta})$. Le tore $T_{0}$ est non ramifi\'e et a donc une structure naturelle sur $\mathfrak{o}_{v}$. On a $T_{0}(\mathfrak{o}_{v})=T_{0}(F_{v})\cap K_{v}$. D'apr\`es les th\'eor\`emes de structure de Bruhat et Tits, l'application
$$(1)\qquad  \begin{array}{ccc}T_{0}(\mathfrak{o}_{v})\times K_{sc,v}&\to& K_{v}\\(t,x)&\mapsto &t\pi(x)\\ \end{array}$$
est surjective, o\`u $\pi:G_{SC}\to G$ est l'homomorphisme naturel. 

On se rappelle le groupe $G_{ab}(F_{v})$ de [I] 1.12. On a $G_{ab}(F_{v})=G(F_{v})/\pi(G_{SC}(F_{v}))$. Soit $S$ un sous-tore maximal  de $G$ d\'efini sur $F_{v}$ et non ramifi\'e. On a

(2) $S(\mathfrak{o}_{v})$ et $K_{v}$ ont m\^eme image dans $G_{ab}(F_{v})$. 

Preuve. C'est clair d'apr\`es (1) si $S=T_{0}$. Il suffit donc de prouver que, si $S_{1}$ et $S_{2}$ sont deux sous-tores maximaux de $G$ d\'efinis sur $F_{v}$ et non ramifi\'es, on a

(3) $S_{1}(\mathfrak{o}_{v})$ et $S_{2}(\mathfrak{o}_{v})$ ont m\^eme image dans $G_{ab}(F_{v})$.

 Puisque $S_{1}$ et $S_{2}$ sont d\'eploy\'es sur $F_{v}^{nr}$, on peut fixer $x\in G(F_{v}^{nr})$ de sorte que $S_{2}=ad_{x}(S_{1})$. Soit $s_{1}\in S_{1}(\mathfrak{o}_{v})$, posons $s_{2}=xs_{1}x^{-1}$. Alors $s_{2}\in S_{2}(\mathfrak{o}_{v}^{nr})$. On sait que tout commutateur se rel\`eve canoniquement dans $G_{SC}$. Cela entra\^{\i}ne qu'il existe $y\in G_{SC}(F_{v}^{nr})$ tel que $xs_{1}x^{-1}s_{1}^{-1}=\pi(y)$. On a $s_{2}=\pi(y)s_{1}$. Soit $\sigma\in \Gamma_{v}^{nr}$. On a $\sigma(s_{2})=\pi(\sigma(y))s_{1}$ donc $s_{2}\sigma(s_{2})^{-1}=\pi(y\sigma(y)^{-1})$. Mais $s_{2}\sigma(s_{2})^{-1}$ appartient \`a $S_{2}(\mathfrak{o}_{v}^{nr})$ et l'image r\'eciproque de ce groupe dans $G_{SC}(F_{v}^{nr})$ est $S_{2,sc} (\mathfrak{o}_{v}^{nr})$. Donc $y\sigma(y)^{-1}\in S_{2,sc}(\mathfrak{o}_{v}^{nr})$. L'application $\sigma\mapsto y\sigma(y)^{-1}$ est un cocycle et $H^1(\Gamma_{v}^{nr};S_{2,sc}(\mathfrak{o}_{v}^{nr}))=\{1\}$. On peut donc fixer $u\in S_{2,sc}(\mathfrak{o}_{v}^{nr})$ tel que $y\sigma(y)^{-1}=u^{-1}\sigma(u)$ pour tout $\sigma$. Posons $s'_{2}=\pi(u)s_{2}$ et $y'=uy$. Alors $y'\in G_{SC}(F_{v})$, $s'_{2}\in S_{2}(\mathfrak{o}_{v}^{nr})$ et $s'_{2}=\pi(y')s_{1}$. Ces relations entra\^{\i}nent que $s'_{2}\in S_{2}(\mathfrak{o}_{v})$. De plus, puisque  $\pi(G_{SC}(F_{v}))$ est le noyau de la projection $G(F_{v})\to G_{ab}(F_{v})$, $s'_{2}$ et $s_{1}$ ont m\^eme image dans $G_{ab}(F_{v})$. Cela d\'emontre que l'image dans ce groupe de $S_{1}(\mathfrak{o}_{v})$ est contenue dans celle de $S_{2}(\mathfrak{o}_{v})$. En \'echangeant les r\^oles de $S_{1}$ et $S_{2}$, on obtient l'\'egalit\'e de ces images. $\square$
 
 Rappelons que $H^1(W_{F_{v}};Z(\hat{G}))$ est le groupe dual de $G_{ab}(F_{v})$. Notons  $Res_{I_{v}}:H^1(W_{F_{v}};Z(\hat{G}))\to H^1(I_{v};Z(\hat{G}))$ l'homomorphisme de restriction. On a

    (4) le noyau de $Res_{I_{v}}$ est l'annulateur dans $H^1(W_{F_{v}};Z(\hat{G}))$ de l'image de $K_{v}$ dans $G_{ab}(F_{v})$. 
    
    Preuve. Supposons d'abord que $G$ soit un tore, notons-le plut\^ot $T_{0}$.  On a alors $K_{v}=T_{0}(\mathfrak{o}_{v})$.   On a le diagramme de suites exactes
    $$\begin{array}{ccccccc}1&\to&H^1(W_{v}^{nr};\hat{T}_{0})&\to&H^1(W_{F_{v}};\hat{T}_{0})&\stackrel{Res_{I_{v}}}{\to}&H^1(I_{v};\hat{T}_{0})\\0&\leftarrow& X_{*}(T_{0})^{\Gamma_{v}^{nr}}&\leftarrow&T_{0}(F_{v})&\leftarrow&T_{0}(\mathfrak{o}_{v})\\ \end{array}$$
    Le groupe $H^1(W_{v}^{nr};\hat{T}_{0})$ s'identifie au quotient des coinvariants $\hat{T}_{0,\Gamma_{v}^{nr}}$, qui s'identifie lui-m\^eme \`a $Hom(X_{*}(T_{0})^{\Gamma_{v}^{nr}};{\mathbb C}^{\times})$. 
    Les fl\`eches de gauche du diagramme ci-dessus sont compatibles \`a cette dualit\'e et \`a celle entre $H^1(W_{F_{v}};\hat{T}_{0})$ et $T_{0}(F_{v})$. Un \'el\'ement de $H^1(W_{F_{v}};\hat{T}_{0})$ appartient au noyau de $Res_{I_{v}}$ si et seulement s'il provient d'un \'el\'ement de $H^1(W_{v}^{nr};\hat{T}_{0})$, ou encore si et seulement si le caract\`ere de $T_{0}(F_{v})$ qu'il d\'efinit se quotiente en un caract\`ere de $  X_{*}(T_{0})^{\Gamma_{v}^{nr}}$, ou encore si et seulement si ce caract\`ere de $T_{0}(F_{v})$  annule $T_{0}(\mathfrak{o}_{v})$. Cela prouve (4) pour un tore. 
    
    Passons au cas g\'en\'eral. Avec les notations introduites plus haut, on a le diagramme commutatif
    $$\begin{array}{ccc}H^1(W_{F_{v}};Z(\hat{G}))&\to&H^1(W_{F_{v}};\hat{T}_{0})\\ Res_{I_{v}}\downarrow\,\,&&Res_{I_{v}}\downarrow\,\,\\ H^1(I_{v};Z(\hat{G}))&\to&H^1(I_{v};\hat{T}_{0})\\ \end{array}$$
    Les fl\`eches horizontales sont injectives: par des suites exactes de cohomologie, cela r\'esulte de la connexit\'e de $\hat{T}_{0,ad}^{\Gamma_{F_{v}}}$ et $\hat{T}_{0,ad}^{I_{v}}$. Un \'el\'ement  $\chi\in H^1(W_{F_{v}};Z(\hat{G}))$ appartient donc au noyau de $Res_{I_{v}}$ si et seulement si son image dans $H^1(W_{F_{v}};\hat{T}_{0})$ appartient au noyau de l'application similaire. D'apr\`es ce que l'on a d\'ej\`a prouv\'e, cela \'equivaut \`a ce que $\chi$ annule l'image de $T_{0}(\mathfrak{o}_{v})$ dans $G_{ab}(F_{v})$.  D'apr\`es (2), cette image est aussi celle de $K_{v}$.     $\square$

 On note pr\'ecis\'ement $W$ le groupe de Weyl de $G$ relatif \`a $T_{0}$. Soit $E$ une extension finie  non ramifi\'ee de $F_{v}$ telle que $G$ soit d\'eploy\'e sur $E$.   Montrons que

(5) soit $u:Gal(E/F_{v})\to W$ un cocycle; alors il existe $x\in K_{v}^{E}$ tel que, pour tout $\sigma\in Gal(E/F_{v})$, $x\sigma(x)^{-1}$ normalise $T_{0}$ et ait $u(\sigma)$ pour  image dans $W$.

Fixons un Frobenius $\phi\in \Gamma_{F_{v}}$. On peut relever $u(\phi)$ en un \'el\'ement de $K_{v}^{nr}$ qui normalise $T_{0}$. On peut m\^eme supposer que cet \'el\'ement appartient \`a un sous-groupe invariant par $ \Gamma_{v}^{nr}$ dont tous les \'el\'ements sont d'ordre fini (le groupe engendr\'e par l'image d'une section de Springer et tous les \'el\'ements d'ordre au plus $2$ de $T_{0}(\mathfrak{o}_{v}^{nr})$ convient). Appliquant [W1] 4.2(2), il existe $y\in K_{v}^{nr}$ tel que $y\phi(y)^{-1}$ soit un rel\`evement de $u(\phi)$ dans le normalisateur de $T_{0}$. Notons $N=[E:F_{v}]$. Alors $y\phi^N(y)^{-1}$ rel\`eve $u(\phi^N)=1$ donc appartient \`a $T_{0}\cap K_{v}^{nr}=
 T_{0}(\mathfrak{o}_{v}^{nr})$. L'application $\sigma\mapsto y\sigma(y)^{-1}$ est un cocycle  de $ Gal(F_{v}^{nr}/E)$ dans $T_{0}(\mathfrak{o}_{v}^{nr})$. Un tel cocycle est un cobord. Il existe donc $t\in T_{0}(\mathfrak{o}_{v}^{nr})$ tel que $y\sigma(y)^{-1}=t\sigma(t)^{-1}$ pour tout $\sigma\in Gal(F_{v}^{nr}/E)$. On pose $x=t^{-1}y$.  Alors $x\in K_{v}^{nr}$ et $x\phi^N(x)^{-1}=1$, donc $x\in  K_{v}^{E}$. De plus $x\phi(x)^{-1}$ rel\`eve $u(\phi)$. Par la relation de cocycle, $x\sigma(x)^{-1}$ rel\`eve donc $u(\sigma)$ pour tout $\sigma\in Gal(E/F_{v})$. Cela prouve (5). $\square$

  On a vu en [I] 6.2(2) qu'il existait un couple $(\nu_{0},e_{0})$ tel que $\nu_{0}\in T_{0}(\mathfrak{o}_{v}^{nr})$, $e_{0}\in Z(\tilde{G},{\cal E}_{0})(F_{v}^{nr})$ et $\nu_{0}e_{0}\in \tilde{K}_{v}$. L'hypoth\`ese $v\not\in V_{ram}$ implique que $p$ est "grand", donc que l'image naturelle de $X_{*}(T_{0,sc})$ dans $X_{*}(T_{0,ad})$ est d'indice premier \`a $p$. Puisqu'extraire des racines d'ordre premier \`a $p$ d'unit\'es ne cr\'ee que des extensions non-ramifi\'ees, l'application
$$T_{0,sc}(\mathfrak{o}_{v}^{nr})\to T_{0,ad}(\mathfrak{o}_{v}^{nr})$$
est surjective. Quitte \`a multiplier $\nu_{0}$ par un \'el\'ement de $Z(G;F_{v}^{nr})\cap T_{0}(\mathfrak{o}_{v}^{nr})$ et $e_{0}$ par l'inverse de cet \'el\'ement, on peut donc supposer qu'il existe $\nu_{0,sc}\in T_{0,sc}(\mathfrak{o}_{v}^{nr})$ tel que $\nu_{0}=\pi(\nu_{0,sc})$. La condition $\nu_{0}e_{0}\in \tilde{G}(F_{v})$ implique alors qu'il existe un cocycle non ramifi\'e $\sigma\mapsto z(\sigma)$ de $\Gamma_{F_{v}}$ dans $Z(G_{SC})$  ($=Z(G_{SC};F_{v}^{nr})$ d'apr\`es l'hypoth\`ese $v\not\in V_{ram}$) tel que $\sigma(\nu_{0,sc})=z(\sigma)\nu_{0,sc}$ et $\sigma(e_{0})=z(\sigma)^{-1}e_{0}$ pour tous $\sigma$. On suppose d\'esormais que $(\nu_{0},e_{0})$ v\'erifie cette hypoth\`ese.

On peut identifier ${\cal E}_{0}$ \`a $\underline{la}$ paire de Borel \'epingl\'ee de $G$, munie de son action galoisienne quasi-d\'eploy\'ee.   
  Notons $\mu_{0}$ l'image du couple $(\nu_{0},e_{0})$ dans $(T^*/(1-\theta^*)(T^*))\times_{{\cal Z}(G)}{\cal Z}(\tilde{G})$. Ce terme d\'epend des choix effectu\'es. Mais les groupes $T^*(\mathfrak{o}_{v}^{nr})$ et $T^*(\bar{\mathfrak{o}}_{v})$ agissent sur $T^*$ et ces actions se descendent en des actions sur $(T^*/(1-\theta^*)(T^*))\times_{{\cal Z}(G)}{\cal Z}(\tilde{G})$. On v\'erifie que la classe $T^*(\mathfrak{o}^{nr}_{v})\mu_{0}$ ne d\'epend pas des choix, a fortiori la classe $T^*(\bar{\mathfrak{o}}_{v})\mu_{0}$ n'en d\'epend pas non plus.   Remarquons que l'on obtiendrait les m\^emes classes en rempla\c{c}ant $\mu_{0}$ par l'image du couple $(1,e_{0})$. On a

(6) ces classes sont invariantes par l'action de $W^{\theta^*}$.

Preuve.  Le couple $(1,e_{0})$ \'etant invariant par $W^{\theta^*}$, la seule chose \`a prouver est que les sous-groupes $T^*(\mathfrak{o}_{v}^{nr})$ et $T^*(\bar{\mathfrak{o}}_{v})$ le sont aussi. C'est imm\'ediat puisque, le groupe $G$ \'etant non ramifi\'e en $v$, tout \'el\'ement de $W$ d\'efinit un  automorphisme de $T^*$ d\'efini sur $F_{v}^{nr}$. $\square$

On pose $\mu(\tilde{K}_{v})=T^*(\bar{\mathfrak{o}}_{v})\mu_{0}$.

\bigskip

\subsection{Param\`etres dans le cas local non ramifi\'e}

On fixe $v\not\in V_{ram}$. 
Soit ${\cal X}\in {\bf Stab}(\tilde{G}(F_{v}))$. Fixons $(\mu,\omega_{\bar{G}})\in Stab(\tilde{G}(F_{v}))$ d'image ${\cal X}$ et relevons $\mu$ en un couple $(\nu,\bar{e})$, avec $\nu\in T^*$ et $\bar{e}\in {\cal Z}(\tilde{G})$.  Consid\'erons les conditions suivantes:

(nr1) pour tout $\alpha\in \Sigma(T^*)$, $(N\alpha)(\nu)\in \bar{\mathfrak{o}}_{v}^{\times}$;

(nr2)(type 1) pour $\alpha\in \Sigma(T^*)$  de type $1$, la relation $(N\alpha)(\nu)\not=1$ entra\^{\i}ne que la r\'eduction dans $\bar{{\mathbb F}}_{v}$ de $(N\alpha)(\nu)$ est diff\'erente de $1$;

(nr2)(types 2 et 3) pour $\alpha\in \Sigma(T^*)$ de type $2$ ou $3$, et pour $\epsilon=\pm 1$, la relation $(N\alpha)(\nu)+ \epsilon\not=0$ entra\^{\i}ne que la r\'eduction dans $\bar{{\mathbb F}}_{v}$ de  $(N\alpha)(\nu)+\epsilon$ est non nulle;

 (nr3) $\mu\in \mu(\tilde{K}_{v})$;

(nr4) le cocycle $\omega_{\bar{G}}$ est non ramifi\'e.

Ces conditions ne d\'ependent pas des rel\`evements choisis.

\ass{Lemme}{(i) Supposons v\'erifi\'ees (nr3) et (nr4). Alors il existe  une classe de conjugaison stable ${\cal O}\in \tilde{G}_{ss}(F_{v})/st-conj$ et un \'el\'ement $\eta\in {\cal O}$ tels que $\chi^{\tilde{G}}({\cal O})={\cal X}$,  que $\eta\in \tilde{K}_{v}$ et  que $G_{\eta}$ soit quasi-d\'eploy\'e.

(ii) Supposons v\'erifi\'ees (nr1) et (nr2) et supposons qu'il existe une classe de conjugaison stable ${\cal O}\in \tilde{G}_{ss}(F_{v})/st-conj$ telle que $\chi^{\tilde{G}}({\cal O})={\cal X}$ et que ${\cal O} $ coupe $\tilde{K}_{v}$. Alors  (nr3) et (nr4) sont v\'erifi\'ees. Pour $\eta\in {\cal O}\cap \tilde{K}_{v}$, le groupe $G_{\eta}$ est non ramifi\'e et $K_{v}\cap G_{\eta}(F_{v})$ en est un sous-groupe compact hypersp\'ecial. Plus pr\'ecis\'ement, en notant ${\cal K}_{\eta}$ le sch\'ema en groupes associ\'e \`a ce groupe hypersp\'ecial, on a ${\cal K}_{\eta}(\mathfrak{o}_{E})={\cal K}_{v}(\mathfrak{o}_{E})\cap G_{\eta}(E)$ pour toute extension non ramifi\'ee $E$ de $F_{v}$. }

Preuve de (i).   Avec les notations de 1.5, on peut identifier ${\cal E}^*$ \`a ${\cal E}_{0}$ et supposer  que $\bar{e}$ est l'image de $e_{0}$. Notons $\bar{\nu}$ l'image de $\nu$ dans $T_{0}/(1-\theta)(T_{0})$, o\`u $\theta=ad_{e_{0}}$. Puisque $\mu$ est fixe par l'action $\sigma\mapsto \omega_{\bar{G}}(\sigma)\sigma$ et puisque $\sigma(e_{0})=z(\sigma)^{-1}e_{0}$, on a l'\'egalit\'e
$$(1) \qquad \omega_{\bar{G}}(\sigma)\sigma(\bar{\nu})=z(\sigma)\bar{\nu}.$$
Puisque $\omega_{\bar{G}}$ et $z$ sont non ramifi\'es, cette relation implique que $\bar{\nu}\in (T_{0}/(1-\theta)(T_{0}))(F_{v}^{nr})$. D'apr\`es (nr3), on a aussi $\bar{\nu}\in (T_{0}/(1-\theta)(T_{0}))(\bar{\mathfrak{o}}_{v})$. Donc $\bar{\nu}$ appartient \`a l'intersection de ces deux groupes, qui n'est autre que $(T_{0}/(1-\theta)(T_{0}))(\mathfrak{o}_{v}^{nr})$. De la suite exacte de tores non ramifi\'es
$$1\to (1-\theta)(T_{0})\to T_{0}\to T_{0}/(1-\theta)(T_{0})\to 1$$
se d\'eduit une  suite exacte
$$1\to ((1-\theta)(T_{0}))(\mathfrak{o}_{v}^{nr})\to T_{0}(\mathfrak{o}_{v}^{nr})\to (T_{0}/(1-\theta)(T_{0}))(\mathfrak{o}_{v}^{nr})\to 1.$$
Quitte \`a changer $\nu$, on peut donc supposer $\nu\in T_{0}(\mathfrak{o}_{v}^{nr})$. La relation (1) implique l'existence d'une cocha\^{\i}ne non ramifi\'ee $t:\Gamma_{F_{v}}\to ((1-\theta)(T_{0}))(\mathfrak{o}_{v}^{nr})$ telle que
$$\omega_{\bar{G}}(\sigma)\sigma(\nu)=z(\sigma)t(\sigma)\nu.$$
Puisque $z$ est un cocycle \`a valeurs centrales, cette \'egalit\'e implique que $t$ est un cocycle si l'on munit $((1-\theta_{0})(T_{0}))(\mathfrak{o}_{v}^{nr})$ de la structure galoisienne $\sigma\mapsto \omega_{\bar{G}}(\sigma)\sigma_{G}$.  Le th\'eor\`eme de Lang implique qu'un tel cocycle est un cobord. Donc, quitte \`a changer encore $\nu$, on peut supposer
$$\omega_{\bar{G}}(\sigma)\sigma_{G}(\nu)=z(\sigma)\nu$$
pour tout $\sigma$.  Introduisons le groupe $G_{SC}^{\theta}$ des points fixes de $\theta =ad_{e_{0}}$ dans $G_{SC}$. De la paire de Borel \'epingl\'ee ${\cal E}_{0}$ se d\'eduit une telle paire pour $G_{SC}^{\theta}$, puis un sous-groupe compact hypersp\'ecial de $G_{SC}^{\theta}(F_{v})$, notons-le $K^1_{v}$. Comme pr\'ec\'edemment, il d\'etermine un sous-groupe $K^{1,nr}_{v}$ de $G_{SC}^{\theta}(F_{v}^{nr})$. En appliquant 1.5(5), on obtient un \'el\'ement $k\in K_{v}^{1,nr}$ tel que, pour tout $\sigma\in  \Gamma_{v}^{nr}$, $k^{-1}\sigma(k)$ normalise $T_{0}$ et ait $\omega_{\bar{G}}(\sigma)$ pour image dans $W^{\theta}$. 
  Posons  $\eta=k\nu e_{0}k^{-1}$. On a $\eta\in \tilde{G}(F_{v}^{nr})$. De plus
$$\sigma(\eta)=\sigma(k)\sigma(\nu)\sigma(k^{-1})\sigma(e_{0})$$
pour tout $\sigma\in \Gamma_{v}^{nr}$ parce  $k\in G_{SC}^{\theta}$. Puis
$$\sigma(\eta)=\sigma(k)\omega_{\bar{G}}(\sigma)^{-1}(z(\sigma)\nu)\sigma(k^{-1})z(\sigma)^{-1}e_{0}.$$
En utilisant l'\'egalit\'e $\omega_{\bar{G}}(\sigma)^{-1}=ad_{\sigma(k)^{-1}k}$, on obtient $\sigma(\eta)=\eta$, donc  $\eta$ appartient \`a $\tilde{G}(F_{v})$. Par un calcul analogue, la relation 1.1(2) entra\^{\i}ne que la paire de Borel $(k(B_{0}\cap G_{\nu e_{0}})k^{-1},kT^{\theta,0}_{0}k^{-1})$ de $G_{\eta}$ est d\'efinie sur $F_{v}$. Donc $G_{\eta}$ est quasi-d\'eploy\'e.
On a aussi $\eta=k\nu k^{-1}e_{0} =k\nu k^{-1}\nu_{0}^{-1}\nu_{0}e_{0}$. On sait que $\nu_{0}e_{0}\in \tilde{G}(F_{v})$, donc $k\nu k^{-1}\nu_{0}^{-1}\in G(F_{v})$. Or c'est un \'el\'ement de $K_{v}^{nr}$. Donc il appartient \`a $K_{v}$. Puisque $\nu_{0}e_{0}\in \tilde{K}_{v}$, cela entra\^{\i}ne que $\eta\in \tilde{K}_{v}$. Il est clair que la classe de conjugaison stable de $\eta$ a pour image ${\cal X}$ dans ${\bf Stab}(\tilde{G}(F_{v}))$. La conclusion de (i) est v\'erifi\'ee.

Preuve de (ii). Soit $\eta\in \tilde{K}_{v}$ dont la classe de conjugaison stable s'envoie sur ${\cal X}$. On peut fixer un entier $N$ premier \`a $p$ de sorte que

- $\theta^{N}=1$;

- le nombre d'\'el\'ements de $Z(G_{SC})$ divise $N$.

On introduit le groupe non connexe $G^+=G\rtimes \{1,\theta,...,\theta^{N-1}\}$, muni de l'action de  $\Gamma_{F_{v}}$ d\'efinie  par $\sigma(g,\theta^{j})=(z(\sigma)^{-j}\sigma(g),\theta^{j})$. Alors $\tilde{G}$ s'identifie \`a la composante $G\theta$, $ge_{0}$ s'identifiant \`a $g\theta$ pour $g\in G$. Dans cette situation, on a d\'efini en [W1] 5.2 la notion d'\'el\'ement compact de $\tilde{G}(F_{v})$: un \'el\'ement est compact si et seulement si le sous-groupe qu'il engendre dans $G^+(F_{v})$ est d'adh\'erence compacte. La condition $\eta\in \tilde{K}_{v}$ entra\^{\i}ne que le sous-groupe engendr\'e par $\eta$ est inclus dans $G^+(F_{v})\cap (K_{v}^{nr}\times \{1,\theta,...,\theta^{N-1}\})$, donc $\eta$ est compact. D'apr\`es [W1] 5.2, on peut d\'ecomposer $\eta$ en $u\eta_{p'}$, o\`u $\eta_{p'}$ est d'ordre fini premier \`a $p$ et $u\in G(F_{v})$ est topologiquement unipotent. Ces \'el\'ements appartiennent \`a l'adh\'erence du groupe engendr\'e par $\eta$. En particulier, ils commutent entre eux et le commutant de $\eta$ dans $G$ est l'intersection des commutants de $u$ et $\eta_{p'}$. La description ci-dessus de l'adh\'erence du sous-groupe engendr\'e par $\eta$ montre que $u\in K_{v}$ et $\eta_{p'}\in \tilde{K}_{v}$. Le lemme [W1] 5.4 nous dit qu'il existe $k\in K_{v}^{nr}$ et un \'el\'ement $\nu_{p'}\in T_{0}(\mathfrak{o}_{v}^{nr})$ d'ordre fini premier \`a $p$ de sorte que $k\eta_{p'}k^{-1}=\nu_{p'} e_{0}$ et que le tore $T=ad_{k^{-1}}(T_{0})$ soit d\'efini sur $F_{v}$. L'\'el\'ement $u$ appartient \`a $Z_{G}(\eta_{p'};F_{v})$. La condition $v\not\in V_{ram}$ implique que l'indice de $G_{\eta_{p'}}(F_{v})$ dans ce groupe est premier \`a $p$. Puisque $u$ est topologiquement unipotent, il appartient \`a $G_{\eta_{p'}}(F_{v})$. On peut fixer $x\in G_{\eta_{p'}}$ tel que $xux^{-1}\in T$. Posons $u'=kxux^{-1}k^{-1}$. C'est un \'el\'ement de $T_{0}$ qui est topologiquement unipotent. Posons $\nu=u'\nu_{p'}$. On a $kx\eta (kx)^{-1}=\nu e_{0}$. Par construction, les hypoth\`eses (nr1) et (nr2) s'appliquent \`a cet \'el\'ement $\nu$. Pour $\alpha\in \Sigma(T_{0})$,  $(N\alpha)(\nu_{p'})$  est  un \'el\'ement d'ordre premier \`a $p$ de $\mathfrak{o}_{v}^{nr,\times}$, tandis que $(N\alpha)(u')$ est un \'el\'ement topologiquement unipotent de $\bar{F}_{v}^{\times}$, donc une unit\'e de r\'eduction $1$ dans le corps r\'esiduel. La condition (nr2) nous dit donc  que, pour $\epsilon=1$ dans le cas o\`u $\alpha$ est de type $1$ et pour $\epsilon=\pm 1$ dans le cas o\`u $\alpha$ est de type $2$ ou $3$, la condition $(N\alpha)(\nu)=\epsilon$ est \'equivalente \`a $(N\alpha)(\nu_{p'})=\epsilon$. D'apr\`es la description des commutants de $\nu e_{0}$ et $\nu_{p'}e_{0}$, cela entra\^{\i}ne que ces deux commutants ont m\^eme syst\`eme de racines. On sait d\'ej\`a que $G_{\eta}\subset G_{\eta_{p'}}$, donc $G_{\nu e_{0}}\subset G_{\nu_{p'}e_{0}}$. Ces deux groupes sont donc \'egaux et aussi $G_{\eta}=G_{\eta_{p'}}$.  Les relations  $u\in G_{\eta_{p'}}=G_{\eta}\subset G_{u}$ entra\^{\i}nent que $u$ appartient au centre de $G_{\eta_{p'}}$. Mais alors l'\'el\'ement $x$ de la construction ci-dessus ne sert \`a rien: on a $xux^{-1}=u$. On reprend la construction avec $x=1$. L'\'el\'ement $u'$ appartient maintenant \`a $T_{0}(\mathfrak{o}_{v}^{nr})$ et $\nu$ aussi. On a $k\eta k^{-1}=\nu e_{0}$. En reprenant les d\'efinitions, on voit que la relation $\nu\in T_{0}(\mathfrak{o}_{v}^{nr})$ entra\^{\i}ne la condition (nr3) tandis que la relation $k\in K_{v}^{nr}$ entra\^{\i}ne la condition (nr4). Enfin, les derni\`eres assertions de (ii) r\'esultent de l'\'egalit\'e $G_{\eta}=G_{\eta_{p'}}$ et du lemme [W1] 5.4(ii)  ou plus exactement de sa preuve, qui montre que ces assertions sont v\'erifi\'ees par le groupe $G_{\eta_{p'}}$. $\square$

\bigskip

\subsection{Param\`etres et endoscopie}
Soit ${\bf G}'=(G',{\cal G}',\tilde{s})$ une donn\'ee endoscopique de $(G,\tilde{G},{\bf a})$. On fixe une paire de Borel \'epingl\'ee $\hat{{\cal E}}=(\hat{B},\hat{T},(\hat{E}_{\alpha})_{\alpha\in \Delta})$ de $\hat{G}$ de sorte que $ad_{\tilde{s}}$ conserve $(\hat{B},\hat{T})$. On en d\'eduit un automorphisme $\hat{\theta}$ de $\hat{G}$ et une action galoisienne modifi\'ee qui conserve $\hat{{\cal E}}$, cf. [I] 1.4. On \'ecrit $\tilde{s}=s\hat{\theta}$. On choisit une paire de Borel \'epingl\'ee $\hat{{\cal E}}'$ de $\hat{G}'$ dont la paire de Borel sous-jacente soit $(\hat{B}\cap \hat{G}', \hat{T}^{\hat{\theta},0})$. On note ${\cal E}^{_{'}*}=(B^{_{'}*},T^{_{'}*},(E'_{\alpha})_{\alpha\in \Delta'})$ $\underline{la}$ paire de Borel \'epingl\'ee de $G'$. De l'injection naturelle $\hat{T}^{\hat{\theta},0}\subset \hat{T}$ se d\'eduit un homomorphisme
$$\xi:T^*\to T^*/(1-\theta^*)(T^*)\simeq T^{_{'}*}.$$
En munissant ces objets des actions quasi-d\'eploy\'ees, il y a un cocycle $\omega_{G'}:\Gamma_{F}\to W^{\theta^*}$ de sorte que $\sigma_{G^{_{'}*}}\circ\xi=\xi\circ \omega_{G'}(\sigma)\circ\sigma_{G^*}$ pour tout $\sigma\in \Gamma_{F}$. Le groupe $W^{G'}$ s'identifie \`a un sous-groupe de $W^{\theta^*}$ en identifiant $w'\in W^{G'}$ \`a l'unique \'el\'ement $w\in W^{\theta^*}$ tel que $\xi \circ w=w'\circ \xi$.  On a d\'ecrit en [I] 1.6 l'ensemble de racines $\Sigma(T^{_{'}*})$. 

Il y a aussi un homomorphisme naturel ${\cal Z}(\tilde{G})\to {\cal Z}(\tilde{G}')$ qui est \'equivariant pour les actions galoisiennes. On en d\'eduit un isomorphisme
$$\tilde{\xi}:( T^*/(1-\theta^*)(T^*))\times_{{\cal Z}(G)}{\cal Z}(\tilde{G})\simeq T^{_{'}*}\times_{Z(G')}{\cal Z}(\tilde{G}').$$

Soit $(\mu',\omega_{\bar{G}'})\in Stab(\tilde{G}'(F))$. Par l'inverse de l'isomorphisme pr\'ec\'edent, $\mu'$ s'identifie \`a un \'el\'ement $\mu\in ( T^*/(1-\theta^*)(T^*))\times_{{\cal Z}(G)}{\cal Z}(\tilde{G})$. L'ensemble de racines $\Sigma(\mu')$   ne s'identifie pas \`a un sous-ensemble de $\Sigma(\mu)$ car le premier ensemble est form\'e d'\'el\'ements $N\alpha$ ou $2N\alpha$ pour $\alpha\in \Sigma(T^*)$ tandis que le second est form\'e d'\'el\'ements $\alpha_{res}$. Mais, pour tout $\beta'\in \Sigma(\mu')$, il existe un unique $\beta\in \Sigma(\mu)$ de sorte que la restriction de $\beta'$ \`a $T^{*,\theta^*,0}$ soit de la forme $b\beta$, avec $b\in {\mathbb Q}$, $b>0$. Pour une raison qui appara\^{\i}tra plus tard, on note $\Sigma^{\bar{H}}$  l'ensemble de ces racines $\beta$.
Rappelons plus pr\'ecis\'ement cette correspondance. On fixe $\nu\in T^*$ et $\bar{e}\in {\cal Z}(\tilde{G})$ tel que $\mu$ soit l'image de $(\nu,\bar{e})$. Alors
$$\Sigma(\mu')=\{N\alpha; \alpha\in \Sigma(T^*),\alpha \text{ de type }1,(N\hat{\alpha})(s)=1,(N\alpha)(\nu)=1\}$$
$$\cup \{2N\alpha; \alpha\in \Sigma(T^*),\alpha \text{ de type }2,(N\hat{\alpha})(s)=1,(N\alpha)(\nu)= 1\}$$
$$\cup \{2N\alpha; \alpha\in \Sigma(T^*),\alpha \text{ de type }2,(N\hat{\alpha})(s)=1,(N\alpha)(\nu)= -1\}$$
$$\cup \{N\alpha; \alpha\in \Sigma(T^*),\alpha \text{ de type }3,(N\hat{\alpha})(s)=-1,(N\alpha)(\nu)=1\}.$$
On envoie un \'el\'ement $N\alpha$ du premier ensemble sur $\alpha_{res}$, un \'el\'ement $2N\alpha$ du deuxi\`eme sur $\alpha_{res}$, un \'el\'ement $2N\alpha$ du troisi\`eme sur $2\alpha_{res}$, un \'el\'ement $N\alpha$ du quatri\`eme sur $\alpha_{res}/2$.  On v\'erifie que $\Sigma^{\bar{H}}$ est un sous-syst\`eme de racines de $\Sigma(\mu)$, dont le groupe de Weyl n'est autre que $W^{G'}(\mu')$. Cela entra\^{\i}ne que $W^{G'}(\mu')$ est un sous-groupe de $W(\mu)$.

L'application $\sigma\mapsto \omega_{\bar{G}'}(\sigma)\omega_{G'}(\sigma)$ est un cocycle de $\Gamma_{F}$ dans $W^{\theta^*}$. Pour tout $\sigma\in \Gamma_{F}$, $\omega_{\bar{G}'}(\sigma)\omega_{G'}(\sigma)\sigma_{G^*}=\omega_{\bar{G}'}(\sigma)\circ \sigma_{G^{_{'}*}}$ fixe $\mu$. 
Le couple form\'e de $\mu$ et de ce cocycle appartient donc \`a $\underline{Stab}(\tilde{G}(F))$. On note $(\mu,\omega_{\bar{G}})$ l'unique \'el\'ement de $Stab(\tilde{G}(F))$ qui lui est \'equivalent. On note $\omega_{\tilde{H}}$ l'unique cocycle de $\Gamma_{F}$ dans $W(\mu)$ tel que $\omega_{\bar{H}}(\sigma)\circ\omega_{\bar{G}}(\sigma)=\omega_{\bar{G}'}(\sigma)\omega_{G'}(\sigma)$ pour tout $\sigma$. 
 On a ainsi d\'efini une application
$$\begin{array}{ccc}Stab(\tilde{G}'(F))&\to&Stab(\tilde{G}(F))\\ (\mu',\omega_{\bar{G}'})&\mapsto &(\mu,\omega_{\bar{G}}).\\ \end{array}$$
On v\'erifie qu'elle se quotiente en une application
$$(1) \qquad {\bf Stab}(\tilde{G}'(F))\to {\bf Stab}(\tilde{G}(F)).$$
Celle-ci ne d\'epend pas du choix de la paire de Borel \'epingl\'ee $\hat{{\cal E}}$. Ses fibres sont \'evidemment finies. 

On a

(2) supposons ${\bf G}'$ elliptique; soit ${\cal X}'\in {\bf Stab}(\tilde{G}'(F))$, notons ${\cal X}$ son image par l'application ci-dessus; si ${\cal X}'$ est elliptique, alors ${\cal X}$ l'est.

Preuve. On reprend les constructions pr\'ec\'edentes en supposant que $(\mu',\omega_{\bar{G}'})$ est elliptique. Avec les notations de 1.1, on  a l'inclusion  $ Z(G)^{\Gamma_{F},\theta,0}\subset Z(\bar{G})^{\Gamma_{F},0}$ et on doit prouver que c'est une \'egalit\'e. Puisqu'il s'agit de groupes connexes et que l'homomorphisme naturel $T^{*,\theta^*,0}\to T^*/(1-\theta^*)(T^*)$ a un noyau fini, il revient au m\^eme de prouver que les images de ces groupes dans $T^*/(1-\theta^*)(T^*)$  sont \'egales. L'action galoisienne sur $Z(\bar{G})$ est $\sigma\mapsto \omega_{\bar{G}}(\sigma)\sigma_{G^*}$. Puisque $\omega_{\bar{H}}(\sigma)$ appartient au groupe $W(\mu)=W^{\bar{G}}$, son action sur $Z(\bar{G})$ est triviale et on peut aussi bien remplacer l'action pr\'ec\'edente par $\sigma\mapsto \omega_{\bar{H}}(\sigma)\omega_{\bar{G}}(\sigma)\sigma_{G^*}=\omega_{\bar{G}'}(\sigma)\sigma_{G^{_{'}*}}$. Un \'el\'ement $x_{*}\in X_{*}(Z(\bar{G})^0)$ annule l'ensemble $\Sigma(\mu)$, donc aussi $\Sigma^{\bar{H}}$. Son image dans $X_{*}(T^*/(1-\theta^*)(T^*))$ annule donc $\Sigma^{G'}(\mu')$. Cela montre que $Z(\bar{G})^0$ s'envoie dans $Z(\bar{G}')^0$. Donc $Z(\bar{G})^{\Gamma_{F},0}$ s'envoie dans $Z(\bar{G}')^{\Gamma_{F},0}$, o\`u  l'action de $\Gamma_{F}$ est  l'action ci-dessus. Cette action sur $
Z(\bar{G}')$ est la m\^eme qu'en 1.1. L'ellipticit\'e de ${\cal X}'$ nous dit que 
$Z(\bar{G}')^{\Gamma_{F},0}=Z(G')^{\Gamma_{F},0}$. Mais l'ellipticit\'e de ${\bf G}'$ signifie que ce groupe est pr\'ecis\'ement l'image de $Z(G)^{\Gamma_{F},\theta,0}$. D'o\`u la conclusion. $\square$

Si l'on remplace le corps de base $F$ par un compl\'et\'e $F_{v}$, on a une application similaire \`a (1). Soit $v$ une place de $F$ telle que $v\not\in V_{ram}({\bf G}')$. Soit ${\cal X}'\in {\bf Stab}(\tilde{G}'(F_{v}))$, notons ${\cal X}\in {\bf Stab}(\tilde{G}(F_{v}))$ son image par cette application. On a

(3) si ${\cal X}$ v\'erifie la condition (nr1) de 1.6, resp. (nr2), alors ${\cal X}'$ v\'erifie la m\^eme condition; ${\cal X}$ v\'erifie la condition (nr3) si et seulement si ${\cal X}'$ v\'erifie la m\^eme condition.

Quand on passe de ${\cal X}$ \`a ${\cal X}'$, les termes $(N\alpha)(\nu)$ de ces conditions sont remplac\'es par les m\^emes termes ou \'eventuellement par $(2N\alpha)(\nu)$ pour $\alpha$ de type $2$ ou $3$. De plus, pour $G'$, toutes les racines sont de type $1$. La premi\`ere assertion  s'ensuit. D'autre part, l'espace $\tilde{K}'_{v}$ est d\'efini de telle sorte que, par l'isomorphisme $\tilde{\xi}$ d\'efini plus haut, le terme $\mu(\tilde{K}_{v})$ s'identifie \`a $\mu(\tilde{K}'_{v})$. D'o\`u la seconde assertion.

On a aussi

(4) si ${\cal X}'$ v\'erifie la condition (nr4), alors ${\cal X}$ la v\'erifie aussi.

Soit $\sigma$ dans le groupe d'inertie $I_{v}$. Avec les notations de la construction ci-dessus, $\sigma_{G^*}$ agit trivialement sur $\Sigma(T^*)$  puisque $v\not\in V_{ram}$, $\omega_{G'}(\sigma)=1$ puisque $v\not\in V_{ram}({\bf G}')$ et $\omega_{\bar{G}'}(\sigma)=1$ puisque ${\cal X}'$ v\'erifie (nr4). Donc $\omega_{\bar{G}}(\sigma)=\omega_{\bar{H}}(\sigma)^{-1}$. Cet \'el\'ement appartient \`a $W(\mu)$ et conserve $\Sigma_{+}(\mu)$. C'est donc l'identit\'e. $\square$

\bigskip

\subsection{Retour sur la correspondance entre classes de conjugaison stable}
 
Soit ${\bf G}'=(G',{\cal G}',\tilde{s})$ une donn\'ee endoscopique de $(G,\tilde{G},{\bf a})$ et soit $V$ un ensemble fini de places de $F$. On a un diagramme
$$\begin{array}{ccc}\tilde{G}'_{ss}(F_{V})/st-conj&...&\tilde{G}_{ss}(F_{V})/st-conj\\ \chi^{\tilde{G}'_{V}}\downarrow&&\downarrow\chi^{\tilde{G}_{V}}\\ {\bf Stab}(\tilde{G}'(F_{V}))&\to&{\bf Stab}(\tilde{G}(F_{V}))\\ \end{array}$$
La ligne horizontale du haut n'est qu'une correspondance, pr\'ecis\'ement une fonction d\'efinie sur un sous-ensemble de $\tilde{G}'_{ss}(F_{V})/st-conj$, \`a valeurs dans $\tilde{G}_{ss}(F_{V})/st-conj$. L'application $\chi^{\tilde{G}'_{V}}$ est bijective (lemme 1.3). On consid\`ere maintenant le cas o\`u $V$ est r\'eduit \`a une seule place $v$. 

\ass{Lemme}{Soit ${\cal O}'\in \tilde{G}'_{ss}(F_{v})/st-conj$. Notons ${\cal X}'$ son image par $\chi^{\tilde{G}'_{v}}$ et ${\cal X}$ l'image de ${\cal X}'$ dans ${\bf Stab}(\tilde{G}(F_{v}))$.

(i)  Supposons que ${\cal O}'$ corresponde \`a une classe ${\cal O}\in \tilde{G}_{ss}(F_{v})/st-conj$  par
 la correspondance sup\'erieure du diagramme. Alors $\chi^{\tilde{G}_{v}}({\cal O})={\cal X}$. 
 
 (ii) Supposons que $v\not\in V_{ram}$ et que ${\cal X}$ v\'erifie les conditions (nr3) et (nr4) de 1.6. Alors ${\cal O}'$ correspond \`a une classe ${\cal O}\in \tilde{G}_{ss}(F_{v})/st-conj$ et on a $\chi^{\tilde{G}_{v}}({\cal O})={\cal X}$.}

Preuve.  Prouvons (i).  On fixe des paires de Borel \'epingl\'ees dans les groupes duaux comme en  1.7. Puisque ${\cal O}$ correspond \`a ${\cal O}'$,  on peut fixer un diagramme $(\epsilon,B',T',B,T,\eta)$ (sur le corps de base $F_{v}$) avec $\epsilon\in {\cal O}'$ et $\eta\in {\cal O}$. On compl\`ete $(B,T)$ et $(B',T')$ en des paires de Borel \'epingl\'ees ${\cal E}$ et ${\cal E}'$, que l'on peut  identifier \`a ${\cal E}^*$ et ${\cal E}^{_{'}*}$.  On \'ecrit $\eta=\nu e$, avec $\nu\in T$ et $e\in Z(\tilde{G},{\cal E})$. On a alors $\epsilon=\xi(\nu)e'$ o\`u $e'$ est l'image de $e$ dans ${\cal Z}(\tilde{G}')\simeq Z(\tilde{G}',{\cal E}')$. Donc les termes $\mu$ et $\mu'$ co\"{\i}ncident dans $(T/(1-\theta)(T))\times_{{\cal Z}(G)}{\cal Z}(\tilde{G})$. On introduit les cocha\^{\i}nes $u_{{\cal E}}$, $u_{\eta}$, $u_{{\cal E}'}$ et $u_{\epsilon}$ comme en 1.2. Soit $\sigma\in \Gamma_{F}$.  Parce que les \'el\'ements $\eta$ et $\epsilon$ et les tores $T$ et $T'$ sont d\'efinis sur $F_{v}$,    les \'el\'ements $u_{{\cal E}}(\sigma)$ et  $u_{\eta}(\sigma)$ normalisent $T$ et  les \'el\'ements $u_{{\cal E}'}(\sigma)$ et $u_{\epsilon}(\sigma)$ normalisent $T'$. Leurs images dans $W$ sont invariantes par $\theta$  (cf. [I] 1.3(3)). On les note par des lettres grasses: ${\bf u}_{{\cal E}}(\sigma)$ est l'image de $u_{{\cal E}}(\sigma)$ dans $W^{\theta}$. Les d\'efinitions entra\^{\i}nent $\omega_{\eta}(\sigma)={\bf u}_{\eta}(\sigma){\bf u}_{{\cal E}}(\sigma)^{-1}$, $\omega_{\epsilon}(\sigma)={\bf u}_{\epsilon}(\sigma){\bf u}_{{\cal E}'}(\sigma)^{-1}$. L'homomorphisme $\xi$ est \'equivariant pour les actions naturelles. Donc
$$\xi \circ{\bf u}_{{\cal E}}(\sigma)^{-1}\circ\sigma_{G^*}={\bf u}_{{\cal E}'}(\sigma)^{-1}\circ\sigma_{G^{_{'}*}}\circ\xi.$$
Il en r\'esulte que
$$\omega_{G'}(\sigma)={\bf u}_{{\cal E}'}(\sigma){\bf u}_{{\cal E}}(\sigma)^{-1}.$$
On v\'erifie que $\omega_{\epsilon}(\sigma)\omega_{G'}(\sigma)\sigma_{G^*}$ conserve $\mu$. 
On peut donc introduire la cocha\^{\i}ne $\omega_{\bar{H}}$ \`a valeurs dans $W(\mu)$ telle que $\omega_{\bar{H}}(\sigma)^{-1}\omega_{\epsilon}(\sigma)\omega_{G'}(\sigma)\sigma_{G^*}$ conserve $\Sigma_{+}(\mu)$. Pour prouver que $\chi^{\tilde{G}_{v}}({\cal O})={\cal X}$, il suffit de prouver que 
$$\omega_{\bar{H}}(\sigma)^{-1}\omega_{\epsilon}(\sigma)\omega_{G'}(\sigma)=\omega_{\eta}(\sigma)$$
pour tout $\sigma\in \Gamma_{F_{v}}$. Puisque, compos\'es avec $\sigma_{G^*}$, ces deux \'el\'ements conservent $\Sigma_{+}(\mu)$, il suffit de prouver que leurs images dans $W(\mu)\backslash W^{\theta^*}$ sont \'egales. On a
$$\omega_{\bar{H}}(\sigma)^{-1}\omega_{\epsilon}(\sigma)\omega_{G'}(\sigma)=\omega_{\bar{H}}(\sigma)^{-1}{\bf u}_{\epsilon}(\sigma){\bf u}_{{\cal E}'}(\sigma)^{-1}\omega_{G'}(\sigma)=\omega_{\bar{H}}(\sigma)^{-1}{\bf u}_{\epsilon}(\sigma){\bf u}_{{\cal E}}(\sigma)^{-1},$$
$$\omega_{\eta}(\sigma)={\bf u}_{\eta}(\sigma){\bf u}_{{\cal E}}(\sigma)^{-1}.$$
Dans le membre de droite de la premi\`ere \'egalit\'e, les deux premiers termes appartiennent \`a $W^{G'}(\mu')\subset W(\mu)$. Dans le membre de droite de la seconde, le premier terme appartient \`a $W(\mu)$.  Les deux termes appartiennent donc \`a $W(\mu){\bf u}_{{\cal E}}(\sigma)^{-1}$, ce qui prouve l'assertion cherch\'ee. 

Prouvons (ii). On peut choisir $\epsilon\in {\cal O}'$ tel que $G'_{\epsilon}$ soit quasi-d\'eploy\'e. On fixe une paire de Borel $(B',T')$ de $G'$ conserv\'ee par $ad_{\epsilon}$ et telle que $(B'\cap G_{\epsilon},T')$ soit d\'efinie sur $F_{v}$. D'apr\`es le lemme 1.6(i), ${\cal X}$ est l'image par $\chi^{\tilde{G}_{v}}$ d'une classe ${\cal O}\in \tilde{G}_{ss}(F_{v})/st-conj$ et on peut choisir un \'el\'ement $\eta\in {\cal O}$ tel que $G_{\eta}$ soit quasi-d\'eploy\'e. On choisit une paire de Borel $(B_{\star},T_{\star})$ de $G$ conserv\'ee par $ad_{\eta}$ et telle que $(B_{\star}\cap G_{\eta},T_{\star}\cap G_{\eta})$ soit d\'efinie sur $F_{v}$.   A l'aide des paires $(B',T')$ et $(B_{\star},T_{\star})$, on construit les couples $(\mu_{\epsilon},\omega_{\epsilon})$ et $(\mu_{\eta},\omega_{\eta})$ comme en 1.2. De $(\mu_{\epsilon},\omega_{\epsilon})$ se d\'eduit comme en 1.7 un \'el\'ement de $Stab(\tilde{G}(F_{v}))$. L'hypoth\`ese que ${\cal X}$ est l'image de ${\cal X}'$ signifie que cet \'el\'ement est conjugu\'e \`a $(\mu_{\eta},\omega_{\eta})$. En appliquant 1.2(3), on voit que l'on peut remplacer la paire $(B_{\star},T_{\star})$ par une autre qui poss\`ede les m\^emes propri\'et\'es que ci-dessus, de sorte qu'apr\`es ce remplacement, le couple $(\mu_{\eta},\omega_{\eta})$ soit \'egal \`a celui d\'eduit de $(\mu_{\epsilon},\omega_{\epsilon})$. On compl\`ete les paires $(B',T')$ et $(B_{\star},T_{\star})$ en des paires de Borel \'epingl\'ees ${\cal E}'$ et ${\cal E}_{\star}$.  Ecrivons $\eta=\nu_{\star} e_{\star}$ avec $\nu_{\star}\in T_{\star}$ et $e_{\star}\in Z(\tilde{G},{\cal E}_{\star})$ et $\epsilon=\nu' e'$, o\`u $\nu'\in T'$ et  $e'\in Z(\tilde{G}',{\cal E}')$ est l'image de $e_{\star}$. On a l'homomorphisme usuel $\xi_{T_{\star},T'}:T_{\star}\to T'$. On voit que l'\'egalit\'e des couples ci-dessus signifie que $\xi_{T_{\star},T'}(\nu_{\star})=\nu'$ et qu'il existe un cocycle $\sigma\mapsto \omega_{\bar{H}}(\sigma)$ de $\Gamma_{F_{v}}$ dans $W^{G_{\eta}}$ de sorte que l'on ait la relation $\sigma_{G'}\circ \xi_{T_{\star},T'}=\xi_{T_{\star},T'}\circ \omega_{\bar{H}}(\sigma)\circ \sigma_{G}$. Parce que $G_{\eta}$ est quasi-d\'eploy\'e, on peut appliquer [K1] corollaire 2.2: on peut fixer $g\in G_{\eta,SC}$ de sorte que le tore $T=ad_{g^{-1}}(T_{\star})$ soit d\'efini sur $F_{v}$ et que, pour tout $\sigma\in \Gamma_{F_{v}}$, l'\'el\'ement $g\sigma(g)^{-1}$ normalise $T_{\star}$ et ait $\omega_{\bar{H}}(\sigma)$ pour image dans $W^{G_{\eta}}$. Posons ${\cal E}=ad_{g^{-1}}({\cal E}_{\star})$. Son tore sous-jacent est $T$ et on note $B$ le sous-groupe de Borel sous-jacent. Parce que $ad_{g}$ fixe $\eta$, on a l'\'egalit\'e $\eta=\nu e$, o\`u $\nu=ad_{g^{-1}}(\nu_{\star})$ et $e=ad_{g^{-1}}(e_{\star})$. On a aussi $\xi_{T,T'}=\xi_{T_{\star},T'}\circ ad_{g}$. On v\'erifie alors que $\xi_{T,T'}(\nu)=\nu'$ et que l'on a l'\'egalit\'e $\sigma_{G'}\circ \xi_{T,T'}=\xi_{T,T'}\circ \sigma$. Mais cela signifie que $(\epsilon,B',T',B,T,\eta)$ est un diagramme. Donc ${\cal O}$ correspond \`a ${\cal O}'$. Cela prouve la premi\`ere assertion de (ii). La seconde r\'esulte de (i). 
$\square$

\bigskip

\subsection{Distributions associ\'ees \`a un param\`etre}
Soit ${\cal X}\in {\bf Stab}(\tilde{G}(F))$.  On peut  repr\'esenter ${\cal X}$ par un \'el\'ement $(\mu,\omega_{\bar{G}})\in Stab(\tilde{G}(F))$ et relever $\mu$ en $(\nu,\bar{e})$, avec $\nu\in T^*$ et $\bar{e}\in {\cal Z}(\tilde{G})$.   On consid\`ere les conditions  (nr1),...,(nr4) de 1.6 pour une place $v\in Val(F)-V_{ram}$. 
Il est clair que chacune d'elles est v\'erifi\'ee sauf pour un ensemble fini de places. On note $S({\cal X})$ le plus petit ensemble de places contenant $V_{ram}$ tel que (nr1) et (nr2)  soient v\'erifi\'ees hors de $S({\cal X})$. On note $S({\cal X},\tilde{K})$ le plus petit ensemble de places contenant  $V_{ram}$ tel que les conditions (nr1), (nr2), (nr3) et (nr4)  soient v\'erifi\'ees hors de $S({\cal X},\tilde{K})$. On a \'evidemment $S({\cal X})\subset S({\cal X},\tilde{K})$.

Si ${\cal X}$ est l'image  par $\chi^{\tilde{G}}$ d'une classe de conjugaison stable dans $\tilde{G}_{ss}(F)$, l'ensemble $S({\cal X})$ co\"{\i}ncide avec l'ensemble $S({\cal C})$ d\'efini en [VI] 2.3 pour toute classe de conjugaison ${\cal C}$ contenue dans cette classe de conjugaison stable. 

Soient ${\cal X}\in {\bf Stab}(\tilde{G}(F))$ et $V$ un ensemble fini de places de $F$ contenant $V_{ram}$. On d\'efinit une distribution $A^{\tilde{G}}(V,{\cal X},\omega)\in D_{g\acute{e}om}(\tilde{G}(F_{V}))\otimes Mes(G(F_{V}))^*$. Elle est nulle si ${\cal X}$ n'appartient pas \`a l'image de l'application $\chi^{\tilde{G}}$. Supposons que ${\cal X}=\chi^{\tilde{G}}({\cal O})$, o\`u ${\cal O}$ est une classe de conjugaison stable semi-simple. Pour toute classe de conjugaison ordinaire ${\cal C}\in \tilde{G}_{ss}(F)/conj$, on a d\'efini la distribution $A^{\tilde{G}}(V,{\cal C},\omega)$ en [VI] 2.3. On pose
$$A^{\tilde{G}}(V,{\cal X},\omega)=\sum_{{\cal C}\in \tilde{G}_{ss}(F)/conj, {\cal C}\subset {\cal O}}A^{\tilde{G}}(V,{\cal C},\omega).$$
Les termes de cette somme sont presque tous nuls. En effet, pour tout $v\in V$, la classe de conjugaison stable dans $\tilde{G}(F_{v})$ engendr\'ee par ${\cal O}$ se d\'ecompose en un nombre fini de classes de conjugaison par $G(F_{v})$. Il existe donc un sous-ensemble fini $\tilde{U}_{V}$ de $\tilde{G}(F_{V})$ tel que, pour tout $\eta\in {\cal O}$, la classe de conjugaison par $G(F_{V})$ de $\eta$ coupe $\tilde{U}_{V}$.  De plus, un terme $A^{\tilde{G}}(V,{\cal C},\omega)$ n'est non nul que si, pour tout $v\not\in V$, la classe de conjugaison par $G(F_{v})$ engendr\'ee par ${\cal C}$ coupe $\tilde{K}_{v}$, cf. [VI] 2.3(6).  Le lemme [VI] 2.1 entra\^{\i}ne la finitude affirm\'ee.

On a

(1) si $S({\cal X},\tilde{K})-S({\cal X})\not\subset V$, alors  $A^{\tilde{G}}(V,{\cal X},\omega)=0$. 

C'est clair si ${\cal X}$ n'est pas dans l'image de $\chi^{\tilde{G}}$. Si ${\cal X}=\chi^{\tilde{G}}({\cal O})$,  consid\'erons une place $v\in S({\cal X},\tilde{K})-S({\cal X})$ qui n'appartient pas \`a $  V$. En $v$, les conditions (nr1) et (nr2) sont v\'erifi\'ees, mais au moins une des conditions (nr3) ou (nr4) ne l'est pas. Le lemme 1.6(ii) nous dit qu'il n'existe aucune classe ${\cal C}\subset {\cal O}$ telle que la classe engendr\'ee par ${\cal C}$ dans $\tilde{G}(F_{v})$ coupe $\tilde{K}_{v}$. Toutes les distributions $A^{\tilde{G}}(V,{\cal C},\omega)$ intervenant sont donc nulles. $\square$

On a

(2) si $S({\cal X})\subset V$ et  si ${\cal X}$ n'est pas elliptique, alors $A^{\tilde{G}}(V,{\cal X},\omega)=0$.

  Pour toute classe ${\cal C}\subset {\cal O}$, on a $S({\cal X})=S({\cal C})$ et  l'hypoth\`ese que ${\cal X}$ n'est pas elliptique entra\^{\i}ne que ${\cal C}$ ne l'est pas non plus. D'apr\`es la d\'efinition de [VI] 2.3, on a alors $A^{\tilde{G}}(V,{\cal C},\omega)=0$. $\square$
  
  En 1.1 et 1.2 on a associ\'e \`a ${\cal X}$ un groupe $\bar{G}$ et un caract\`ere automorphe de $\bar{G}({\mathbb A}_{F})$. Le couple form\'e de ce groupe et de ce caract\`ere n'est d\'efini qu'\`a isomorphisme pr\`es mais la condition que la restriction du caract\`ere \`a $Z(\bar{G};{\mathbb A}_{F})$ est triviale est insensible \`a un tel isomorphisme. Par abus de langage, nous la formulerons: la restriction de $\omega$ \`a  $Z(\bar{G};{\mathbb A}_{F})$ est triviale. On a
  
  (3) si les restrictions de $\omega$ \`a $Z(G;{\mathbb A}_{F})^{\theta}$ et \`a $Z(\bar{G};{\mathbb A}_{F})$ ne sont pas toutes deux triviales, alors $A^{\tilde{G}}(V,{\cal X},\omega)=0$.
  
  En effet, soient ${\cal C}\subset {\cal O}$ et $\dot{\gamma}\in {\cal C}$. Alors 
  la restriction de $\omega$ \`a  $Z(\bar{G};{\mathbb A}_{F})$ est triviale si et seulement si  la restriction de $\omega$ \`a $Z(G_{\dot{\gamma}};{\mathbb A}_{F})$ est triviale. L'assertion r\'esulte alors de [VI] proposition 2.3(iv).

\bigskip

\subsection{Distributions stables et endoscopiques associ\'ees \`a un param\`etre}
On fixe un ensemble fini $V$ de places de $F$ contenant $V_{ram}$.

Supposons $(G,\tilde{G},{\bf a})$ quasi-d\'eploy\'e et \`a torsion int\'erieure. Soit ${\cal X}\in {\bf Stab}(\tilde{G}(F))$. On va d\'efinir une distribution $SA^{\tilde{G}}(V,{\cal X})\in D_{g\acute{e}om}(\tilde{G}(F_{V}))\otimes Mes(G(F_{V}))^*$. Comme toujours, on a besoin de supposer par r\'ecurrence que cette distribution v\'erifie certaines propri\'et\'es. Il y a des propri\'et\'es formelles qui permettent de "recoller" ces distributions dans la situation habituelle, cf. [VI] 1.15. Ce sont les m\^emes que dans cette r\'ef\'erence et on les abandonne au lecteur. On  donnera toutefois  dans le paragraphe suivant des formules plus explicites  dans la situation "avec caract\`ere central".  Comme en [VI] 5.2(1), il y a une condition concernant les espaces hypersp\'eciaux $\tilde{K}_{v}$ pour $v\not\in V$. La d\'efinition fournit une distribution qui d\'epend de ces espaces. On doit savoir que

(1) elle ne d\'epend que des classes de conjugaison par $G_{AD}(F_{v})$ des $\tilde{K}_{v}$ pour $v\not\in V$.

Surtout, on doit supposer par r\'ecurrence que cette distribution $SA^{\tilde{G}}(V,{\cal X})$ est stable. Modulo ces hypoth\`eses, soit ${\bf G}'=(G',{\cal G}',s)$ une donn\'ee endoscopique de $(G,\tilde{G},{\bf a})$ relevante et non ramifi\'ee hors de $V$,  avec $dim(G'_{SC})<dim(G_{SC})$. Soit  ${\cal X}'\in {\bf Stab}(\tilde{G}'(F))$. Alors on peut d\'efinir $SA^{{\bf G}'}(V,{\cal X}')\in D_{g\acute{e}om}^{st}({\bf G}'_{V})\otimes Mes(G'(F_{V}))^*$.
On pose (avec les notations de [VI] 5.1):
$$(2) \qquad SA^{\tilde{G}}(V,{\cal X})=A^{\tilde{G}}(V,{\cal X})-\sum_{{\bf G}'\in {\cal E}(\tilde{G},V), G'\not=G}$$
$$\sum_{{\cal X}'\in {\bf Stab}(\tilde{G}'(F)), {\cal X}'\mapsto {\cal X}}i(\tilde{G},\tilde{G}')transfert(SA^{{\bf G}'}(V,{\cal X}')),$$
o\`u on a not\'e ${\cal X}'\mapsto {\cal X}$ l'application (1)  de 1.7.

On revient au cas o\`u $(G,\tilde{G},{\bf a})$ est quelconque.   Pour ${\cal X}\in {\bf Stab}(\tilde{G}(F))$, on pose  
$$(3) \qquad A^{\tilde{G},{\cal E}}(V,{\cal X},\omega)=\sum_{{\bf G}'\in {\cal E}(\tilde{G},{\bf a},V)}\sum_{{\cal X}'\in {\bf Stab}(\tilde{G}'(F)), {\cal X}'\mapsto {\cal X}}i(\tilde{G},\tilde{G}')transfert(SA^{{\bf G}'}(V,{\cal X}')).$$

{\bf Remarque.} Comme souvent, le cas o\`u $(G,\tilde{G},{\bf a})$ est quasi-d\'eploy\'e et \`a torsion int\'erieure est un peu particulier. Dans ce cas, les hypoth\`eses de r\'ecurrence ne s'appliquent pas \`a la donn\'ee endoscopique principale ${\bf G}$. Il convient de remplacer le terme $transfert(SA^{{\bf G}}(V,{\cal X}))$ intervenant dans la somme par $SA^{\tilde{G}}(V,{\cal X})$. On a alors $A^{\tilde{G},{\cal E}}(V,{\cal X})=A^{\tilde{G}}(V,{\cal X})$ par d\'efinition de ce terme $SA^{\tilde{G}}(V,{\cal X})$.

 \ass{Th\'eor\`eme (i) (\`a prouver)}{Pour tout  ${\cal X}\in {\bf Stab}(\tilde{G}(F))$, on a l'\'egalit\'e 
 $A^{\tilde{G},{\cal E}}(V,{\cal X},\omega)=A^{\tilde{G}}(V,{\cal X},\omega)$.}

\ass{Th\'eor\`eme (ii)}{Supposons $(G,\tilde{G},{\bf a})$  quasi-d\'eploy\'e et \`a torsion int\'erieure. Alors, pour tout 
  ${\cal X}\in {\bf Stab}(\tilde{G}(F))$, 
 $SA^{\tilde{G}}(V,{\cal X})$ est stable et v\'erifie (1).}

 Le th\'eor\`eme (ii) sera prouv\'e dans cet article, cf. 3.4. Le th\'eor\`eme (i) ne sera enti\`erement prouv\'e que plus tard. Toutefois, dans cet article, nous prouverons ce th\'eor\`eme sauf pour des triplets $(G,\tilde{G},{\bf a})$ particuliers. Pour ceux-ci, le th\'eor\`eme sera prouv\'e sauf pour des ${\cal X}$ particuliers, qui sont en nombre fini. On renvoie \`a 3.5 pour des assertions pr\'ecises. 

 \bigskip

\subsection{Formules dans la situation avec caract\`ere central}
On suppose $(G,\tilde{G},{\bf a})$ quasi-d\'eploy\'e et \`a torsion int\'erieure. On suppose donn\'es  une extension 
$$1\to C_{1}\to G_{1}\to G$$
o\`u $C_{1}$ est un tore central induit, une extension compatible $\tilde{G}_{1}\to \tilde{G}$ o\`u $\tilde{G}_{1}$ est encore \`a torsion int\'erieure et un caract\`ere automorphe $\lambda_{1}$ de $C_{1}({\mathbb A}_{F})$. On a d\'efini l'ensemble de places $V_{1,ram}$ en [VI] 1.15. Pour $v\not\in V_{1,ram}$, on fixe un espace hypersp\'ecial $\tilde{K}_{1,v}\subset \tilde{G}_{1}(F_{v})$ au-dessus de $\tilde{K}_{v}$. On impose la condition de compatibilit\'e globale habituelle: pour $\gamma_{1}\in \tilde{G}_{1}(F)$, on a $\gamma_{1}\in \tilde{K}_{1,v}$ pour presque tout $v$. On a une suite exacte
$$0\to \mathfrak{A}_{C_{1}}\to \mathfrak{A}_{G_{1}}\to \mathfrak{A}_{G}\to 0$$
On a fix\'e en [VI] 1.3 une mesure de Haar sur $\mathfrak{A}_{G}$. On en fixe sur les deux autres groupes de sorte que la suite soit compatible aux mesures. 

 Introduisons $\underline{les}$ paires de Borel \'epingl\'ees   de $G$ et $G_{1}$, dont on  note les tores $T^*$ et $T_{1}^*$. On a des applications naturelles $T_{1}^*\to T^*$ et ${\cal Z}(\tilde{G}_{1})\to {\cal Z}(\tilde{G})$. On en d\'eduit une application
 $$\begin{array}{ccc}T_{1}^*\times_{Z(G_{1})}{\cal Z}(\tilde{G}_{1}) &\to& T^*\times_{Z(G)}{\cal Z}(\tilde{G})\\ \mu_{1}&\mapsto& \mu\\ \end{array}$$
dont les fibres sont isomorphes \`a $C_{1}$. D'o\`u une application $Stab(\tilde{G}_{1}(F))\to Stab(\tilde{G}(F))$, qui, \`a $(\mu_{1},\omega_{\bar{G}})$, associe $(\mu,\omega_{\bar{G}})$. Elle se quotiente en une application 
$$(1) \qquad {\bf Stab}(\tilde{G}_{1}(F))\to {\bf Stab}(\tilde{G}(F)).$$
Celle-ci traduit simplement l'application de projection $\tilde{G}_{1,ss}(F)/st-conj\to \tilde{G}_{ss}(F)/st-conj$. Remarquons que les fibres de cette application ne sont pas isomorphes \`a $C_{1}(F)$ en g\'en\'eral, deux \'el\'ements de $Stab(\tilde{G}_{1}(F))$ de la forme $(\mu_{1},\omega_{\bar{G}})$ et $(c\mu_{1},\omega_{\bar{G}})$, avec $c\in C_{1}(F)$, $c\not=1$, pouvant avoir la m\^eme image dans ${\bf Stab}(\tilde{G}_{1}(F))$.

Soit $V$ un ensemble fini de places contenant $V_{1,ram}$. Fixons des mesures $dg$ sur $G({\mathbb A}_{F})$ et $dc$ sur $C_{1}({\mathbb A}_{F})$, dont on d\'eduit une mesure $dg_{1}$ sur $G_{1}({\mathbb A}_{F})$. On identifie ces mesures \`a des mesures sur $G(F_{V})$, $C_{1}(F_{V})$ et $G_{1}(F_{V})$, cf. [VI] 1.1.  On rappelle que les distributions d\'efinies en 1.9 et 1.10 d\'ependent de l'espace $\tilde{K}^V=\prod_{v\not\in V}\tilde{K}_{v}$, bien que l'on n'ait pas fait figurer cet espace dans la notation. Dans les formules qui suivent,  on ins\`ere si besoin est cet espace dans la notation, de fa\c{c}on que l'on esp\`ere compr\'ehensible.

Soit ${\cal X}\in {\bf Stab}(\tilde{G}(F))$.
Soit $f\in C_{c,\lambda_{1}}^{\infty}(\tilde{G}_{1}(F_{V}))$. On fixe une fonction $\phi\in C_{c}^{\infty}(\tilde{G}_{1}(F_{V}))$ de sorte que
$$f=\int_{C_{1}(F_{V})}\phi^c\lambda_{1}(c)\,dc.$$
La formule [VI] 2.5(14) donne imm\'ediatement la variante $A_{\lambda_{1}}^{\tilde{G}_{1}}(V,{\cal X})$ de la distribution d\'efinie en 1.9 sous la forme
$$(2) \qquad I^{\tilde{G}_{1}}_{\lambda_{1}}(A^{\tilde{G}_{1}}_{\lambda_{1}}(V,{\cal X}),f\otimes dg)=mes(\mathfrak{A}_{C_{1}}C_{1}(F)\backslash C_{1}({\mathbb A}_{F}))^{-1}\int_{C_{1}(F)\backslash C_{1}({\mathbb A}_{F})}$$
$$\sum_{{\cal X}_{1}\in Fib({\cal X})}I^{\tilde{G}_{1}}(A^{\tilde{G}_{1}}(V,{\cal X}_{1},c^V\tilde{K}_{1}^V),\phi^{c_{V}}\otimes dg_{1})\lambda_{1}(c)\,dc,$$
o\`u $Fib({\cal X})$ est la fibre de l'application (1) au-dessus de ${\cal X}$. 

On obtient par r\'ecurrence une formule analogue pour la variante $SA_{\lambda_{1}}^{\tilde{G}_{1}}(V,{\cal X})$ de la distribution d\'efinie en 1.10:
$$(3) \qquad I_{\lambda_{1}}^{\tilde{G}_{1}}(SA^{\tilde{G}_{1}}_{\lambda_{1}}(V,{\cal X}),f\otimes dg)=mes(\mathfrak{A}_{C_{1}}C_{1}(F)\backslash C_{1}({\mathbb A}_{F}))^{-1}\int_{C_{1}(F)\backslash C_{1}({\mathbb A}_{F})}$$
$$\sum_{{\cal X}_{1}\in Fib({\cal X})}I^{\tilde{G}_{1}}(SA^{\tilde{G}_{1}}(V,{\cal X}_{1},c^V\tilde{K}_{1}^V),\phi^{c_{V}}\otimes dg_{1})\lambda_{1}(c)\,dc.$$

Consid\'erons le cas particulier o\`u $\tilde{G}=G$, $\tilde{G}_{1}=G_{1}$, $\tilde{K}_{v}=K_{v}$ et $\tilde{K}_{1,v}=K_{1,v}$ pour tout $v\not\in V$ et o\`u ${\cal X}$ correspond \`a la classe de conjugaison stable de l'\'el\'ement neutre. Comme toujours, on remplace dans la notation la lettre ${\cal X}$ par un indice $unip$: $A^{\tilde{G}}_{unip}(V)$ au lieu de $A^{\tilde{G}}(V,{\cal X})$ etc... La fibre $Fib({\cal X})$ est alors l'ensemble $\{\xi {\cal X}_{1};\xi\in C_{1}(F)\}$, o\`u ${\cal X}_{1}$ correspond \`a la classe de conjugaison stable de l'unit\'e dans $G_{1}(F)$. D'apr\`es [VI] 2.4(7), on a l'\'egalit\'e:
$$I^{\tilde{G}_{1}}(A^{\tilde{G}_{1}}(V,\xi{\cal X}_{1},\xi^Vc^VK_{1}^V),\phi^{\xi_{V}c_{V}}\otimes dg_{1})=I^{\tilde{G}_{1}}(A^{\tilde{G}_{1}}(V,{\cal X}_{1},c^VK_{1}^V),\phi^{c_{V}}\otimes dg_{1}).$$
Puisque de plus
$$I^{\tilde{G}_{1}}(A^{\tilde{G}_{1}}(V,{\cal X}_{1},c^VK_{1}^V),\phi^{c_{V}}\otimes dg_{1})=\left\lbrace\begin{array}{cc}I^{\tilde{G}_{1}}(A^{\tilde{G}_{1}}(V,{\cal X}_{1},K_{1}^V),\phi^{c_{V}}\otimes dg_{1}),&\text{ si }c^V\in K_{C_{1}}^V,\\ 0,&\text{ sinon,}\\ \end{array}\right.$$
la formule (2) se simplifie en
$$(4) \qquad I^{\tilde{G}_{1}}_{\lambda_{1}}(A^{\tilde{G}_{1}}_{unip,\lambda_{1}}(V),f\otimes dg)=mes(\mathfrak{A}_{C_{1}}C_{1}(F)\backslash C_{1}({\mathbb A}_{F}))^{-1}$$
$$\int_{C_{1}(F_{V})}I^{\tilde{G}_{1}}(A_{unip}^{\tilde{G}_{1}}(V),\phi^{c_{V}}\otimes dg_{1})\lambda_{1}(c)\,dc.$$
La formule (5) se simplifie de m\^eme en 
$$(5) \qquad I^{\tilde{G}}_{\lambda_{1}}(SA^{\tilde{G}_{1}}_{unip,\lambda_{1}}(V),f\otimes dg)=mes(\mathfrak{A}_{C_{1}}C_{1}(F)\backslash C_{1}({\mathbb A}_{F}))^{-1}$$
$$\int_{C_{1}(F_{V})}I^{\tilde{G}_{1}}(SA_{unip}^{\tilde{G}_{1}}(V),\phi^{c_{V}}\otimes dg_{1})\lambda_{1}(c)\,dc.$$
Autrement dit,  la distribution $A^{\tilde{G}_{1}}_{unip,\lambda_{1}}(V) dg$, qui appartient \`a $ D_{unip,\lambda_{1}}(\tilde{G}'_{1}(F_{V}))$ est l'image de l'\'el\'ement $mes(\mathfrak{A}_{C_{1}}C_{1}(F)\backslash C_{1}({\mathbb A}_{F}))^{-1}A^{\tilde{G}_{1}}_{unip}(V) dg_{1}$ de $D_{unip}(\tilde{G}'_{1}(F_{V}))$  par l'application d\'efinie en [II] 1.10(3).  De m\^eme, $SA^{\tilde{G}_{1}}_{unip,\lambda_{1}}(V) dg$ est l'image de 
$$mes(\mathfrak{A}_{C_{1}}C_{1}(F)\backslash C_{1}({\mathbb A}_{F}))^{-1}SA^{\tilde{G}_{1}}_{unip}(V) dg_{1}.$$

 \bigskip

\subsection{Relation avec les distributions associ\'ees aux classes de conjugaison stable locales}

Soit $V$ un ensemble fini de places de $F$ contenant $V_{ram}$. Soit ${\cal O}_{V}\in \tilde{G}_{ss}(F_{V})/st-conj$. Posons ${\cal X}_{V}=\chi^{\tilde{G}_{V}}({\cal O}_{V})$. Pour ${\cal X}\in {\bf Stab}(\tilde{G}(F))$, on note simplement ${\cal X}\mapsto {\cal X}_{V}$ la relation: ${\cal X}_{V}$ est l'image de ${\cal X}$ par localisation. Les distributions des membres de gauche des \'egalit\'es de l'\'enonc\'e ci-dessous ont \'et\'e d\'efinies en [VI] 2.3, 5.2 et 5.4. 

\ass{Proposition}{(i) On a l'\'egalit\'e 
$$A^{\tilde{G}}(V,{\cal O}_{V},\omega) =\sum_{{\cal X}\in {\bf Stab}(\tilde{G}(F)), {\cal X}\mapsto {\cal X}_{V}}A^{\tilde{G}}(V,{\cal X},\omega).$$

(ii) On a l'\'egalit\'e
$$A^{\tilde{G},{\cal E}}(V,{\cal O}_{V},\omega) =\sum_{{\cal X}\in {\bf Stab}(\tilde{G}(F)), {\cal X}\mapsto {\cal X}_{V}}A^{\tilde{G},{\cal E}}(V,{\cal X},\omega).$$

(iii) Supposons $(G,\tilde{G},{\bf a})$ quasi-d\'eploy\'e et \`a torsion int\'erieure. Alors on a l'\'egalit\'e
$$SA^{\tilde{G}}(V,{\cal O}_{V}) =\sum_{{\cal X}\in {\bf Stab}(\tilde{G}(F)), {\cal X}\mapsto {\cal X}_{V}}SA^{\tilde{G}}(V,{\cal X}).$$}

Preuve. Les deux c\^ot\'es de l'\'egalit\'e (i) sont des sommes de $A^{\tilde{G}}(V,{\cal C},\omega)$, o\`u ${\cal C}\in \tilde{G}_{ss}(F)/conj$. Du c\^ot\'e gauche, on somme sur les ${\cal C}$ dont la classe localis\'ee ${\cal C}_{V}$ est contenue dans ${\cal O}_{V}$. Du c\^ot\'e droit, on somme sur les ${\cal C}$ contenus dans une classe de conjugaison stable ${\cal O}$ dont le param\`etre ${\cal X}$ s'envoie par localisation sur ${\cal X}_{V}$. La commutativit\'e du diagramme de 1.4 entra\^{\i}ne que ces ensembles de sommation sont les m\^emes. D'o\`u (i).

Si  $(G,\tilde{G},{\bf a})$ est quasi-d\'eploy\'e et \`a torsion int\'erieure, les termes intervenant dans (ii) sont identiques par d\'efinition aux m\^emes termes o\`u l'on supprime l'exposant ${\cal E}$.  L'assertion (ii) n'est alors autre que (i). Supposons maintenant que $(G,\tilde{G},{\bf a})$ n'est pas quasi-d\'eploy\'e et \`a torsion int\'erieure.
Par d\'efinition
$$A^{\tilde{G},{\cal E}}(V,{\cal O}_{V},\omega) =\sum_{{\bf G}'\in {\cal E}(\tilde{G},{\bf a},V)}i(\tilde{G},\tilde{G}')transfert(SA^{{\bf G}'}(V,{\cal O}_{V,\tilde{G}'})).$$
Fixons ${\bf G}'$. Rappelons que ${\cal O}_{V,\tilde{G}'}$ est la r\'eunion des ${\cal O}_{V}'\in \tilde{G}_{ss}'(F_{V})/st-conj$ qui correspondent \`a ${\cal O}_{V}$. Le lemme 1.8 nous dit que cet ensemble de classes est \'egal \`a celui des classes ${\cal O}_{V}'$ qui v\'erifient les deux conditions suivantes:

- elles correspondent \`a une classe dans $\tilde{G}_{ss}(F_{V})/st-conj$;

- leur param\`etre ${\cal X}'_{V}$ s'envoie sur ${\cal X}_{V}$ par la version locale de l'application  1.7(1) (ce que l'on note ${\cal X}'_{V}\mapsto {\cal X}$).

On se rappelle que $SA^{{\bf G}'}(V,{\cal O}_{V}')$ est \`a support dans l'ensemble des \'el\'ements dont la partie semi-simple appartient \`a ${\cal O}_{V}'$, cf. [VI] 5.2. Si ${\cal O}_{V}'$ ne correspond \`a aucune classe dans $\tilde{G}_{ss}(F_{V})$, le transfert de $SA^{{\bf G}'}(V,{\cal O}_{V}')$ est donc nul. On peut donc aussi bien supprimer la premi\`ere condition ci-dessus:
$$transfert(SA^{{\bf G}'}(V,{\cal O}_{V,\tilde{G}'}))=\sum_{{\cal X}'_{V}\mapsto {\cal X}_{V}}transfert(SA^{{\bf G}'}(V,{\cal O}_{V}')),$$
o\`u ${\cal O}'_{V}$ est l'unique \'el\'ement param\'etr\'e par ${\cal X}'_{V}$. Modulo les formalit\'es habituelles, on peut appliquer (iii) aux termes du membre de droite. On obtient
$$transfert(SA^{{\bf G}'}(V,{\cal O}_{V,\tilde{G}'}))=\sum_{{\cal X}'_{V}\mapsto {\cal X}_{V}}transfert\left(\sum_{{\cal X}'\in {\bf Stab}(\tilde{G}'(F)), {\cal X}'\mapsto {\cal X}'_{V}}SA^{{\bf G}'}(V,{\cal X}')\right).$$
La localisation commute \`a l'application 1.7(1). Sommer en ${\cal X}'_{V}\mapsto {\cal X}_{V}$ puis ${\cal X}'\mapsto {\cal X}'_{V}$ revient \`a sommer sur les ${\cal X}\in {\bf Stab}(\tilde{G}(F))$ tels que ${\cal X}\mapsto {\cal X}_{V}$ puis sur les ${\cal X}'$ tels que ${\cal X}'\mapsto {\cal X}$. Donc
$$transfert(SA^{{\bf G}'}(V,{\cal O}_{V,\tilde{G}'}))=\sum_{{\cal X}\in {\bf Stab}(\tilde{G}(F)), {\cal X}\mapsto {\cal X}_{V}}\sum_{{\cal X}'\in {\bf Stab}(\tilde{G}'(F)), {\cal X}'\mapsto {\cal X}}transfert (SA^{{\bf G}'}(V,{\cal X}')).$$
Puis
$$A^{\tilde{G},{\cal E}}(V,{\cal O}_{V},\omega) =\sum_{{\cal X}\in {\bf Stab}(\tilde{G}(F)), {\cal X}\mapsto {\cal X}_{V}}\sum_{{\bf G}'\in {\cal E}(\tilde{G},{\bf a},V)}i(\tilde{G},\tilde{G}')$$
$$\sum_{{\cal X}'\in {\bf Stab}(\tilde{G}'(F)), {\cal X}'\mapsto {\cal X}}transfert (SA^{{\bf G}'}(V,{\cal X}')).$$
La double somme int\'erieure est par d\'efinition \'egale \`a $A^{\tilde{G},{\cal E}}(V,{\cal X},\omega)$. On obtient (ii).

La preuve de (iii) est analogue. Par d\'efinition, on a cette fois
$$SA^{\tilde{G}}(V,{\cal O}_{V}) =A^{\tilde{G}}(V,{\cal O}_{V})-\sum_{{\bf G}'\in {\cal E}(\tilde{G},V), G'\not=G}i(\tilde{G},\tilde{G}')transfert(SA^{{\bf G}'}(V,{\cal O}_{V,\tilde{G}'})).$$
On conna\^{\i}t (i) pour le premier terme de droite et  (iii) par r\'ecurrence pour les autres termes. Par le m\^eme calcul, on en d\'eduit (iii) pour le terme de gauche. $\square$

 \bigskip

\section{Formules de scindage}

\bigskip

\subsection{ Compl\'ement sur le lemme fondamental pond\'er\'e}

{\bf Par exception, dans ce paragraphe, $F$ est un corps local non-archim\'edien de caract\'eristique nulle}. On note $p$ sa caract\'eristique r\'esiduelle. On consid\`ere un triplet $(G,\tilde{G},{\bf a})$ d\'efini sur $F$ qui est "non ramifi\'e". Pr\'ecis\'ement, comme en [VI] 1.1, on suppose que $G$ et ${\bf a}$ sont non ramifi\'es, que $\tilde{G}(F)$ poss\`ede un sous-espace hypersp\'ecial et que, en posant $e=[F:{\mathbb Q}_{p}]$, on a $p>5$ et $p> N(G)e+1$, o\`u $N(G)$ est d\'efini en [W1] 4.3. On fixe un espace hypersp\'ecial $\tilde{K}$, de groupe associ\'e $K$. Soit $\tilde{M}$ un espace de Levi de $\tilde{G}$ tel que le Levi $M$ associ\'e soit en bonne position relativement \`a $K$. On munit $G(F)$ de la mesure canonique pour laquelle $K$ est de masse totale $1$.  On a d\'efini en [II] 4.1 une forme lin\'eaire $r_{\tilde{M}}^{\tilde{G}}(.,\tilde{K})$ sur $D_{g\acute{e}om}(\tilde{M}(F),\omega)$: on a
$$r_{\tilde{M}}^{\tilde{G}}(\boldsymbol{\gamma},\tilde{K})=J_{\tilde{M}}^{\tilde{G}}(\boldsymbol{\gamma},{\bf 1}_{\tilde{K}}),$$
o\`u ${\bf 1}_{\tilde{K}}$ est la fonction caract\'eristique de $\tilde{K}$. Dans le cas o\`u $(G,\tilde{G},{\bf a})$ est quasi-d\'eploy\'e et \`a torsion int\'erieure, on a d\'efini en 4.2 un avatar stable de cette forme lin\'eaire. Ici, il n'est plus besoin de supposer que $\tilde{M}$ est en bonne position relativement \`a $\tilde{K}$. L'avatar stable est une forme lin\'eaire $s_{\tilde{M}}^{\tilde{G}}(.,\tilde{K})$ sur $D_{g\acute{e}om}^{st}(\tilde{M}(F))$. Elle v\'erifie

(1) $s_{\tilde{M}}^{\tilde{G}}(.,\tilde{K})$ ne d\'epend que de la classe de conjugaison de $\tilde{K}$ par $G_{AD}(F)$. 

 Pour  $\boldsymbol{\delta}\in D_{g\acute{e}om}^{st}(\tilde{M}(F))$, on a la formule famili\`ere
$$s_{\tilde{M}}^{\tilde{G}}(\boldsymbol{\delta},\tilde{K})=r_{\tilde{M}}^{\tilde{G}}(\boldsymbol{\delta},\tilde{K})-\sum_{s\in Z(\hat{M})^{\Gamma_{F}}/Z(\hat{G})^{\Gamma_{F}}, s\not=1}i_{\tilde{M}}(\tilde{G},\tilde{G}'(s))s_{{\bf M}}^{{\bf G}'(s)}(\boldsymbol{\delta},\tilde{K}).$$
Expliquons la signification du dernier terme. Comme on l'a dit  en [II] 4.2, pour $s$ intervenant dans cette somme, le choix d'un sous-espace hypersp\'ecial $\tilde{K}'(s)$ de $\tilde{G}'(s;F)$ permet de d\'efinir par r\'ecurrence un terme $s_{{\bf M}}^{{\bf G}'(s)}(\boldsymbol{\delta},\tilde{K}'(s))$. Mais l'espace $\tilde{K}$ d\'etermine un sous-espace $\tilde{K}'(s)$ de $\tilde{G}'(s;F)$ bien d\'efini \`a conjugaison pr\`es par $G'(s)_{AD}(F)$. La propri\'et\'e (1) permet de noter $s_{{\bf M}}^{{\bf G}'(s)}(\boldsymbol{\delta},\tilde{K})$ le terme  $s_{{\bf M}}^{{\bf G}'(s)}(\boldsymbol{\delta},\tilde{K}'(s))$ pour un tel $\tilde{K}'(s)$. 

Revenons \`a un triplet $(G,\tilde{G},{\bf a})$ quelconque. Soit ${\bf M}'=(M',{\cal M}',\tilde{\zeta})$ une donn\'ee endoscopique de $(M,\tilde{M},{\bf a}_{M})$. Si ${\bf M}'$ est elliptique et non ramifi\'ee (donc relevante d'apr\`es le lemme [I] 6.2), on a d\'efini en [II] 4.3 une forme lin\'eaire $r_{\tilde{M}}^{\tilde{G},{\cal E}}({\bf M}',.,\tilde{K})$ sur $D_{g\acute{e}om}^{st}({\bf M}')$ par l'\'egalit\'e
 $$r_{\tilde{M}}^{\tilde{G},{\cal E}}({\bf M}',\boldsymbol{\delta},\tilde{K})=\sum_{\tilde{s}\in \tilde{\zeta}Z(\hat{M})^{\Gamma_{F},\hat{\theta}}/Z(\hat{G})^{\Gamma_{F},\hat{\theta}}}i_{\tilde{M}'}(\tilde{G},\tilde{G}'(\tilde{s}))s_{{\bf M}'}^{{\bf G}'(\tilde{s})}(\boldsymbol{\delta},\tilde{K}).$$
 Il convient de g\'en\'eraliser la d\'efinition au cas o\`u ${\bf M}'$ est non ramifi\'ee mais pas elliptique. Dans ce cas, il existe un espace de Levi $\tilde{R}$ de $\tilde{M}$ tel que ${\bf M}'$ apparaisse comme une donn\'ee endoscopique elliptique et non ramifi\'ee de $(R,\tilde{R},{\bf a}_{R})$. On peut supposer $R$ en bonne position relativement \`a $K$. Comme en [VI] 4.5, on pose alors
 $$(2) \qquad r_{\tilde{M}}^{\tilde{G},{\cal E}}({\bf M}',\boldsymbol{\delta},\tilde{K})=\sum_{\tilde{L}\in {\cal L}(\tilde{R})}d_{\tilde{R}}^{\tilde{G}}(\tilde{M},\tilde{L})r_{\tilde{R}}^{\tilde{L},{\cal E}}({\bf M}',\boldsymbol{\delta},\tilde{K}^{\tilde{L}}),$$
 o\`u $\tilde{K}^{\tilde{L}}=\tilde{K}\cap \tilde{L}(F)$. 
 
 \ass{Proposition}{Soit ${\bf M}'=(M',{\cal M}',\tilde{\zeta})$ une donn\'ee endoscopique relevante de $(M,\tilde{M},{\bf a}_{M})$ et soit $\boldsymbol{\delta}\in D_{g\acute{e}om}^{st}({\bf M}')$. Alors
 
 (i) si ${\bf M}'$ est non ramifi\'ee, on a l'\'egalit\'e $r_{\tilde{M}}^{\tilde{G}}(transfert(\boldsymbol{\delta}),\tilde{K})=r_{\tilde{M}}^{\tilde{G},{\cal E}}({\bf M}',\boldsymbol{\delta},\tilde{K})$;
 
 (ii) si ${\bf M}'$ n'est pas non ramifi\'ee, on a l'\'egalit\'e $r_{\tilde{M}}^{\tilde{G}}(transfert(\boldsymbol{\delta}),\tilde{K})=0$.}
 
 Preuve. Supposons que ${\bf M}'$ ne soit pas elliptique. Comme ci-dessus, on introduit $\tilde{R}\subset \tilde{M}$ de sorte que ${\bf M}'$ soit une donn\'ee elliptique pour $(R,\tilde{R},{\bf a}_{R})$. Remarquons que ${\bf M}'$ est non ramifi\'ee pour $(M,\tilde{M},{\bf a}_{M})$ si et seulement si elle l'est pour $(R,\tilde{R},{\bf a}_{R})$. Notons $\boldsymbol{\gamma}$ le transfert de $\boldsymbol{\delta}$ \`a $\tilde{R}(F)$. Alors le transfert de $\boldsymbol{\delta}$ \`a $\tilde{M}(F)$ est l'induite $\boldsymbol{\gamma}^{\tilde{M}}$. On a la formule
 $$r_{\tilde{M}}^{\tilde{G}}(\boldsymbol{\gamma}^{\tilde{M}},\tilde{K})=\sum_{\tilde{L}\in {\cal L}(\tilde{R})}d_{\tilde{R}}^{\tilde{G}}(\tilde{M},\tilde{L})r_{\tilde{R}}^{\tilde{L}}(\boldsymbol{\gamma},\tilde{K}^{\tilde{L}}),$$
 cf. [II] 4.1(1). Celle-ci et la formule parall\`ele (2) nous ram\`ene \`a d\'emontrer les assertions de la proposition quand on remplace le couple $(\tilde{G},\tilde{M})$ par un couple $(\tilde{L},\tilde{R})$. En oubliant cette construction, on est ramen\'e au cas ${\bf M}'$ est une donn\'ee endoscopique elliptique de $(M,\tilde{M},{\bf a}_{M})$. L'assertion (i) est le lemme fondamental pond\'er\'e, cf. [II] 4.4. 
 
 Il reste \`a prouver (ii). On suppose donc que ${\bf M}'$ n'est pas non ramifi\'ee. On utilise la m\'ethode d'Arthur qui se base sur un lemme de Kottwitz ([K2] proposition 7.5) que l'on g\'en\'eralise \`a notre cas. Fixons une paire de Borel \'epingl\'ee ${\cal E}=(B,T,(E_{\alpha})_{\alpha\in \Delta})$ de $G$, d\'efinie sur $F$,  telle que $M$ soit standard pour ${\cal E}$ et que le groupe $K$ soit celui issu de ${\cal E}$.  On pose $M_{\sharp}=M/Z(M)^{\theta}$ et $T^M_{\sharp}=T/Z(M)^{\theta}$. On rappelle (cf. [I] 2.7) que le groupe $M_{\sharp}(F)$ op\`ere par conjugaison sur $\tilde{M}(F)$. Si on fixe des donn\'ees auxiliaires $M'_{1}$,...,$\Delta_{1}$ pour ${\bf M}'$, le facteur de transfert $\Delta_{1}(\delta_{1},\gamma)$ se transforme, quand on  conjugue $\gamma$ par un \'el\'ement de $M_{\sharp}(F)$, par un caract\`ere $\omega_{\sharp}$ de $M_{\sharp}(F)$ qui prolonge le caract\`ere $\omega$ de $M(F)$. Notons $\mathfrak{o}$ l'anneau des entiers de $F$. Montrons que
 
 (3) le groupe d'inertie $I_{F}$ est inclus dans ${\cal M}'$ si et seulement si $\omega_{\sharp}$ est trivial sur $T^M_{\sharp}(\mathfrak{o})$, ce qui est encore \'equivalent \`a ce que le cocycle d\'efinissant $\omega_{\sharp}$ soit trivial sur $I_{F}$.
 
 On a une action galoisienne sur $\hat{M}'$. Pour $w\in I_{F}$, fixons $m_{w}=(m(w),w)\in {\cal M}'$ qui agisse sur $\hat{M}'$ comme $w$. Le groupe $I_{F}$ est contenu dans ${\cal M}'$ exactement quand $m(w)\in \hat{M}'$ pour tout $w\in I_{F}$, c'est-\`a-dire quand $m(w)$ est dans la composante neutre du centralisateur de $\tilde{\zeta}$ dans $\hat{M}$.
 
 La deuxi\`eme condition de (3) est \'equivalente \`a ce que la restriction de $\omega_{\sharp}$ \`a $T^M_{\sharp}(F)$ soit non ramifi\'ee c'est-\`a-dire que tout cocycle d\'efinissant ce caract\`ere soit trivial sur $I_{F}$. Un tel cocycle est \`a valeurs dans $Z(\hat{M}_{\sharp})$. Rappelons la construction de la restriction \`a $I_{F}$ d'un tel cocycle. On suppose $\tilde{\zeta}=\zeta\hat{\theta}$ avec $\zeta\in \hat{T}$, avec les notations habituelles. Soient $w\in I_{F}$ et $m_{w}=(m(w),w)$ comme ci-dessus.  Comme pour tout \'el\'ement de ${\cal M}'$, on a la propri\'et\'e de commutation suivante vis-\`a-vis de $\tilde{\zeta}$:
 $$\zeta\hat{\theta}(m(w))\zeta^{-1}=m(w)$$
 puisque ${\bf G}$ et ${\bf a}$ sont non ramifi\'es. Le cocycle associ\'e \`a $\omega_{\sharp}$ est d\'efini en relevant les \'el\'ements intervenant ci-dessus dans le groupe dual de $M_{\sharp}$ et en fait, parce que le groupe $\hat{M}_{\sharp}$ est plus difficile \`a analyser, on rel\`eve m\^eme dans un rev\^etement simplement connexe du groupe d\'eriv\'e de $\hat{M}$ not\'e $\hat{M}_{SC}$. On fixe un \'el\'ement $z(w)\in Z(\hat{M})$ tel que $m(w)z(w)$ soit l'image d'un \'el\'ement $m_{sc}(w)\in \hat{M}_{SC}$. Quitte \`a changer la donn\'ee endoscopique en une donn\'ee \'equivalente, on peut aussi  relever $\zeta$ en un \'el\'ement $\zeta_{sc}$. Et on a une relation
 $$(4) \qquad \zeta_{sc}\hat{\theta}(m_{sc}(w))\zeta_{sc}^{-1}=a_{sc}(w)m_{sc}(w),$$
 avec un \'el\'ement $a_{sc}(w)\in \hat{M}_{SC}$ dont on v\'erifie ais\'ement qu'il est dans le centre de ce groupe. On note $a_{\sharp}(w)$ l'image de $a_{sc}(w)$ dans $\hat{M}_{\sharp}$ (par projection naturelle) et on pose $z_{\sharp}(w)=a_{\sharp}(w)z(w)\hat{\theta}(z(w))^{-1}$. Par d\'efinition (cf. [I] 2.7), $z_{\sharp}$ est la restriction \`a $I_{F}$ d'un cocycle d\'efinissant $\omega_{\sharp}$. En [I] 2.7 (suite exacte suivant la suite (2)), il est montr\'e que $z_{\sharp}(w)$ est l'\'el\'ement neutre de $Z(\hat{M}_{\sharp})$ si et seulement s'il existe $b(w)\in Z(\hat{M}_{SC})$ tel que $z(w)\in b(w)\hat{T}^{\hat{\theta},0}$ et $a_{sc}(w)=\hat{\theta}(b(w))b(w)^{-1}$. 
 
 Supposons que $z_{\sharp}(w)=1$. On fixe les choix pr\'ec\'edents d'o\`u en particulier $b(w)\in Z(\hat{M}_{SC})$. On modifie $z(w)$ en le rempla\c{c}ant par $z(w)b(w)^{-1}$; ainsi $z(w)\in \hat{T}^{\hat{\theta},0}$ et $a_{sc}(w)=1$.  Ainsi, avec (4), $m_{sc}(w)$ est dans le commutant de $\zeta_{sc}\hat{\theta}$ qui est un groupe connexe et son image $m(w)z(w)$ est dans la composante connexe du commutant de $\tilde{\zeta}$ dans $\hat{M}$, c'est-\`a-dire $\hat{M}'$. Mais $\hat{M}'$ contient $\hat{T}^{\hat{\theta},0}$ donc $m(w)$ lui-m\^eme est dans $\hat{M}'$. Puisque $(m(w),w)\in {\cal M}'$, on a $(1,w)\in {\cal M}'$ et $ {\cal M}'$ contient $I_{F}$.
 
 R\'eciproquement, supposons que $(1,w)\in {\cal M}'$ pour tout $w\in I_{F}$. Alors on peut prendre ci-dessus $m(w)=1$. Il est alors clair que $z_{\sharp}(w)=1$ en reprenant la construction rappel\'ee ci-dessus de cet \'el\'ement. Cela d\'emontre (3).

 On termine la preuve du (ii) de la proposition. Le groupe $T^M_{\sharp}(\mathfrak{o})$ agit par conjugaison sur $C_{c}^{\infty}(\tilde{M}(F))$. Il s'en d\'eduit une action sur $I(\tilde{M}(F),\omega)$ et, par dualit\'e, sur $D_{g\acute{e}om}(\tilde{M}(F),\omega)$. La propri\'et\'e de transformation des facteurs de transfert  rappel\'ee ci-dessus  entra\^{\i}ne que, pour $f\in C_{c}^{\infty}(\tilde{M}(F))$, les transferts \`a ${\bf M}'$ de $f$ et de $\omega_{\sharp}(x)ad_{x}(f)$ sont \'egaux. Par dualit\'e, il en r\'esulte que $\omega_{\sharp}(x)ad_{x}(transfert(\boldsymbol{\delta}))=transfert(\boldsymbol{\delta})$. On va faire agir $x$ par conjugaison mais pour que $x$ agisse sur $\tilde{G}(F)$, il faut commencer par relever $x$ en un \'el\'ement de $G_{\sharp}(F)$, o\`u $G_{\sharp}=G/Z(G)^{\theta}$. Pour faire cela, posons $T_{\sharp}=T/Z(G)^{\theta}$. On commence par d\'emontrer que l'application $T_{\sharp}(\mathfrak{o})\to T^M_{\sharp}(\mathfrak{o})$ est surjective: le conoyau s'injecte dans $H^1(Gal(F^{nr}/F);(Z(M)^{\theta}/Z(G)^{\theta})(\mathfrak{o}^{nr}))$, o\`u $F^{nr}$ est l'extension non ramifi\'ee maximale de $F$ et $\mathfrak{o}^{nr}$ est son anneau d'entiers. Or $Z(M)^{\theta}/Z(G)^{\theta}$ est un tore (donc connexe) et ce groupe de cohomologie est nul par le th\'eor\`eme de Lang. On rel\`eve donc $x$ en un \'el\'ement de $T_{\sharp}(\mathfrak{o}_{v})$. Par simple transport de structure, on a l'\'egalit\'e
 $$r_{\tilde{M}}^{\tilde{G}}(transfert(\boldsymbol{\delta}),\tilde{K})=r_{\tilde{M}}^{\tilde{G}}(ad_{x}(transfert(\boldsymbol{\delta})),ad_{x}(\tilde{K})).$$
 Puisque $transfert(\boldsymbol{\delta})$ se transforme par un caract\`ere non trivial de $T^M_{\sharp}(\mathfrak{o})$, il ne reste plus qu'\`a d\'emontrer que $\tilde{K}=ad_{x}(\tilde{K})$. Certainement, $ad_{x}$ conserve $K$ puisque $K$ est associ\'e \`a ${\cal E}$. Il suffit donc de prouver qu'il existe $\gamma\in \tilde{K}$ tel que $ad_{x}(\gamma)\in \tilde{K}$. On a  rappel\'e en 1.5 que l'on pouvait choisir $t\in T(\mathfrak{o}^{nr})$ et $e\in Z(\tilde{G},{\cal E};F^{nr})$ de sorte que $\gamma=te$ appartienne \`a $\tilde{K}$. On rel\`eve $x$ en un \'el\'ement $y\in T(\mathfrak{o}^{nr})$ et on \'ecrit
 $$x\gamma x^{-1}=yt\,ad_{e}(y^{-1})e.$$
 Les \'el\'ements $yt$ et $ad_{e}(y^{-1})$ sont dans $T(\mathfrak{o}^{nr})$. Ainsi $yt\,ad_{e}(y^{-1})\in T(\mathfrak{o}^{nr})$ et aussi $u=yt\,ad_{e}(y^{-1})t^{-1}$. Ainsi
 $$x\gamma x^{-1}\in  ute= u\gamma,$$
 avec $u\in T(\mathfrak{o}^{nr})$. 
 Mais $x\gamma x^{-1}$ et $\gamma$ sont dans $\tilde{G}(F)$ donc $u\in T(\mathfrak{o}^{nr})\cap G(F)=T(\mathfrak{o})\subset K$. Ainsi $x\gamma x^{-1}\in \tilde{K}$, ce qui est l'assertion cherch\'ee. $\square$
 
 \bigskip
 
 \subsection{Version globale du  lemme fondamental pond\'er\'e}
 Le corps de base est de nouveau notre corps de nombres. On consid\`ere un triplet $(G,\tilde{G},{\bf a})$ d\'efini sur $F$.  Soit $U$ un ensemble fini de places tel que, contrairement \`a l'habitude, $U\cap V_{ram}=\emptyset$ et soit $\tilde{M}\in {\cal L}(\tilde{M}_{0})$.  On a d\'efini en [VI] 1.13 une forme lin\'eaire $r_{\tilde{M}}^{\tilde{G}}(.,\tilde{K}_{U})$ sur $D_{g\acute{e}om}(\tilde{M}(F_{U}),\omega)$ par
 $$r_{\tilde{M}}^{\tilde{G}}(\boldsymbol{\gamma},\tilde{K}_{U})=J_{\tilde{M}}^{\tilde{G}}(\boldsymbol{\gamma},{\bf 1}_{\tilde{K}_{U}}).$$

Supposons $(G,\tilde{G},{\bf a})$ quasi-d\'eploy\'e et \`a torsion int\'erieure.   Nous allons d\'efinir une forme lin\'eaire $\boldsymbol{\delta}\mapsto s_{\tilde{M}}^{\tilde{G}}(\boldsymbol{\delta},\tilde{K}_{U})$ sur $D^{st}_{g\acute{e}om}(\tilde{M}(F_{U}))$. Elle doit v\'erifier les propri\'et\'es formelles  habituelles. Elle doit aussi  poss\'eder une propri\'et\'e d'invariance relativement \`a l'action de $G_{AD}(F_{U})$ en le sens suivant. On oublie pour un temps que l'on a fix\'e en 1.1 les espaces hypersp\'eciaux $\tilde{K}_{v}$. Soient  $\tilde{K}_{U}$ et $\tilde{K}'_{U}$ deux sous-espaces hypersp\'eciaux de $\tilde{G}(F_{U})$ dont les groupes sous-jacents $K_{U}$ et $K'_{U}$ sont en bonne position relativement \`a $M_{0}$. Supposons que $\tilde{K}_{U}$ et $\tilde{K}'_{U}$ soient conjugu\'es par un \'el\'ement de $G_{AD}(F_{U})$. On doit alors avoir l'\'egalit\'e

(1) $s_{\tilde{M}}^{\tilde{G}}(\boldsymbol{\delta},\tilde{K}'_{U})=s_{\tilde{M}}^{\tilde{G}}(\boldsymbol{\delta},\tilde{K}_{U})$

\noindent pour tout $\boldsymbol{\delta}$. Comme en [II] 4.2, cette condition permet de g\'en\'eraliser la d\'efinition au cas o\`u le groupe $K_{U}$ n'est pas suppos\'e en bonne position relativement \`a $M_{0}$. 

Soit $s\in Z(\hat{M})^{\Gamma_{F}}$, introduisons la donn\'ee endoscopique ${\bf G}'={\bf G}'(s)$  de $(G,\tilde{G},{\bf a})$. Introduisons des donn\'ees auxiliaires $G'_{1},\tilde{G}'_{1},C_{1},\hat{\xi}_{1},\tilde{K}'_{1,U}$ non ramifi\'ees dans $U$.   Par une extension formelle des d\'efinitions, on d\'efinit $s_{\tilde{M}'_{1},\lambda_{1}}^{\tilde{G}'_{1}}(\boldsymbol{\delta}_{1},\tilde{K}'_{1,U})$ pour  $\boldsymbol{\delta}_{1}\in D^{st}_{g\acute{e}om,\lambda_{1}}(\tilde{M}'_{1}(F_{U}))$. Si on remplace les donn\'ees auxiliaires par d'autres donn\'ees $G'_{2},...,\tilde{K}'_{2,U}$, ces termes se recollent pourvu que la fonction de recollement $\tilde{\lambda}_{12,U}$ v\'erifie l'\'egalit\'e $\tilde{\lambda}_{12,U}(\gamma_{1},\gamma_{2})=1$ pour $\gamma_{1}\in \tilde{K}_{1,U}$ et $\gamma_{2}\in \tilde{K}_{2,U}$. Mais  l'espace $D^{st}_{g\acute{e}om}({\bf M}_{U})$ a \'et\'e d\'efini en consid\'erant des donn\'ees auxiliaires  $G'_{1},\tilde{G}'_{1},C_{1},\hat{\xi}_{1},\Delta_{1,U}$ et des fonctions de recollement identifiant les facteurs de transfert. L'identification entre les deux types de donn\'ees auxiliaires se fait bien s\^ur en utilisant les facteurs de transfert "non ramifi\'es": on d\'eduit des donn\'ees $G'_{1},\tilde{G}'_{1},C_{1},\hat{\xi}_{1},\tilde{K}'_{1,U}$  les donn\'ees $G'_{1},\tilde{G}'_{1},C_{1},\hat{\xi}_{1},\Delta_{1,U}$, o\`u $\Delta_{1,U}$ est le facteur de transfert d\'etermin\'e par le couple $(\tilde{K}_{U},\tilde{K}'_{1,U})$. Les deux notions de recollement co\"{\i}ncident alors. Les formes lin\'eaires pr\'ec\'edemment d\'efinies se recollent en une forme lin\'eaire sur $D^{st}_{g\acute{e}om}({\bf M}_{U})$. La propri\'et\'e (1) montre qu'elle ne d\'epend  que de la classe de conjugaison par $G'_{AD}(F_{U})$ de l'espace hypersp\'ecial $\tilde{K}'_{U}$ de $\tilde{G}'(F_{U})$, laquelle ne d\'epend que de $\tilde{K}_{U}$ lui-m\^eme. Cela justifie de noter cette forme lin\'eaire
$$\boldsymbol{\delta}\mapsto s_{{\bf M}}^{{\bf G}'}(\boldsymbol{\delta},\tilde{K}_{U}).$$

 On peut alors poser la d\'efinition, pour $\boldsymbol{\delta}\in D^{st}_{g\acute{e}om}(\tilde{M}(F_{U}))$:
$$s_{\tilde{M}}^{\tilde{G}}(\boldsymbol{\delta},\tilde{K}_{U})=r_{\tilde{M}}^{\tilde{G}}(\boldsymbol{\delta},\tilde{K}_{U})-\sum_{s\in Z(\hat{M})^{\Gamma_{F}}/Z(\hat{G})^{\Gamma_{F}}, s\not=1}i_{\tilde{M}}(\tilde{G},\tilde{G}'(s))s_{{\bf M}}^{{\bf G}'(s)}(\boldsymbol{\delta},\tilde{K}_{U}).$$

Une preuve similaire \`a celle de la proposition [VI] 4.2 montre que, pour $\boldsymbol{\delta}=\otimes_{v\in U} \boldsymbol{\delta}_{v}$, on a l'\'egalit\'e
$$s_{\tilde{M}}^{\tilde{G}}(\boldsymbol{\delta},\tilde{K}_{U})=\sum_{\tilde{L}^{U}\in {\cal L}(\tilde{M}_{U})}e_{\tilde{M}_{U}}^{\tilde{G}}(\tilde{M},\tilde{L}^{U})\prod_{v\in U}s_{\tilde{M}_{v}}^{\tilde{L}^v}(\boldsymbol{\delta}_{v},\tilde{K}_{v}^{\tilde{L}^v}),$$
o\`u les derniers termes sont les formes lin\'eaires locales d\'efinies en [II] 4.2. Gr\^ace \`a cette \'egalit\'e, les propri\'et\'es formelles requises ainsi que la propri\'et\'e (1) r\'esultent des m\^emes propri\'et\'es de ces formes lin\'eaires locales prouv\'ees en  [II] 4.2.

Revenons au cas o\`u $(G,\tilde{G},{\bf a})$ est quelconque.  Soit ${\bf M}'=(M',{\cal M}',\tilde{\zeta})$ une donn\'ee endoscopique de $(M,\tilde{M},{\bf a}_{M})$ non ramifi\'ee dans $U$. De nouveau, pour $\tilde{s}\in \tilde{\zeta}Z(\hat{M})^{\Gamma_{F},\hat{\theta}}$, on d\'efinit une forme lin\'eaire $\boldsymbol{\delta}\mapsto s_{{\bf M}'}^{{\bf G}'(\tilde{s})}(\boldsymbol{\delta},\tilde{K}_{U})$ sur $D^{st}_{g\acute{e}om}({\bf M}'_{U})$. Pour un tel $\boldsymbol{\delta}$, on pose
$$r_{\tilde{M}}^{\tilde{G},{\cal E}}({\bf M}',\boldsymbol{\delta},\tilde{K}_{U})=\sum_{\tilde{s}\in \tilde{\zeta}Z(\hat{M})^{\Gamma_{F},\hat{\theta}}/Z(\hat{G})^{\Gamma_{F},\hat{\theta}}}i_{\tilde{M}'}(\tilde{G},\tilde{G}'(\tilde{s}))s_{{\bf M}'}^{{\bf G}'(\tilde{s})}(\boldsymbol{\delta},\tilde{K}_{U}).$$
Une preuve similaire \`a celle de la proposition [VI] 4.5 montre que, pour $\boldsymbol{\delta}=\otimes_{v\in U} \boldsymbol{\delta}_{v}$, on a l'\'egalit\'e
$$r_{\tilde{M}}^{\tilde{G},{\cal E}}({\bf M}',\boldsymbol{\delta},\tilde{K}_{U})=\sum_{\tilde{L}^U\in {\cal L}(\tilde{M}_{U})}d_{\tilde{M}_{U}}^{\tilde{G}}(\tilde{M},\tilde{L}^U)\prod_{v\in U}r_{\tilde{M}_{v}}^{\tilde{L}^v,{\cal E}}({\bf M}',\boldsymbol{\delta}_{v},\tilde{K}_{v}^{\tilde{L}^v}),$$
o\`u les derniers termes sont les formes lin\'eaires locales d\'efinies en 2.1. Gr\^ace \`a cette \'egalit\'e et \`a  l'\'egalit\'e parall\`ele concernant les termes $r_{\tilde{M}}^{\tilde{G}}(\boldsymbol{\gamma},\tilde{K}_{U})$ (cf. [VI] 1.13) , la proposition suivante r\'esulte de  celle du paragraphe pr\'ec\'edent.

\ass{Proposition}{Soit ${\bf M}'$ une donn\'ee endoscopique elliptique et relevante  de $(M,\tilde{M},{\bf a}_{M})$ et soit $\boldsymbol{\delta}\in D^{st}_{g\acute{e}om}({\bf M}'_{V})$. Alors

(i) si ${\bf M}'$ est  non ramifi\'ee dans $U$, on a l'\'egalit\'e
$$r_{\tilde{M}}^{\tilde{G}}(transfert(\boldsymbol{\delta}),\tilde{K}_{U})=r_{\tilde{M}}^{\tilde{G},{\cal E}}({\bf M}',\boldsymbol{\delta},\tilde{K}_{U});$$

(ii) si ${\bf M}'$ n'est pas non ramifi\'ee dans $U$, on a l'\'egalit\'e
$$r_{\tilde{M}}^{\tilde{G}}(transfert(\boldsymbol{\delta}),\tilde{K}_{U})=0.$$}

\bigskip

\subsection{Enonc\'e des formules de scindage}
On consid\`ere un triplet quelconque $(G,\tilde{G},{\bf a})$  et deux ensembles finis $V$ et $S$ de places de $F$ tels que $V_{ram}\subset V\subset S$. 

Soit  $\tilde{M}\in {\cal L}(\tilde{M}_{0})$. Choisissons une paire de Borel \'epingl\'ee ${\cal E}=(B,T,(E_{\alpha})_{\alpha\in \Delta})$ de $G$  telle que $M$ soit standard. Alors $M$ d\'etermine un sous-ensemble $\Delta^M\subset \Delta$ et $(B,T,(E_{\alpha})_{a\in \Delta^M})$ est une paire de Borel \'epingl\'ee de $M$. En identifiant ces paires $\underline{aux}$ paires de Borel \'epingl\'ees de $G$ et $M$, on a une inclusion naturelle $Stab(\tilde{M}(F))\subset Stab(\tilde{G}(F))$, qui d\'epend du choix de ${\cal E}$. Il s'en d\'eduit une application ${\bf Stab}(\tilde{M}(F))\to {\bf Stab}(\tilde{G}(F))$. Celle-ci ne d\'epend plus du choix de ${\cal E}$.  On la note simplement ${\cal X}_{M}\mapsto {\cal X}$. Pour ${\cal X}_{M}\in {\bf Stab}(\tilde{M}(F))$, on a d\'efini en 1.9 une distribution $A^{\tilde{M}}(S,{\cal X}_{M},\omega)\in D_{g\acute{e}om}(\tilde{M}(F_{S}),\omega)\otimes Mes(M(F_{S}))^*$. Cet espace est le produit tensoriel de $  D_{g\acute{e}om}(\tilde{M}(F_{S}^V),\omega)\otimes Mes(M(F_{S}^V))^*$  et de $  D_{g\acute{e}om}(\tilde{M}(F_{V}),\omega)\otimes Mes(M(F_{V}))^*$. De plus, le choix des mesures canoniques identifie  $Mes(M(F_{S}^V))^*$ \`a ${\mathbb C}$. On peut donc \'ecrire
$$(1) \qquad A^{\tilde{M}}(S,{\cal X}_{M},\omega)=\sum_{i=1,...,n({\cal X}_{M})}k_{i}^{\tilde{M}}({\cal X}_{M},\omega)_{S}^V\otimes A_{i}^{\tilde{M}}({\cal X}_{M},\omega)_{V}$$
avec des $k_{i}^{\tilde{M}}({\cal X}_{M},\omega)_{S}^V\in  D_{g\acute{e}om}(\tilde{M}(F_{S}^V),\omega)$ et des $A_{i}^{\tilde{M}}({\cal X}_{M},\omega)_{V}\in  D_{g\acute{e}om}(\tilde{M}(F_{V}),\omega)\otimes Mes(M(F_{V}))^*$.  On note $A_{i}^{\tilde{G}}({\cal X}_{M},\omega)_{V}\in  D_{g\acute{e}om}(\tilde{G}(F_{V}),\omega)\otimes Mes(G(F_{V}))^*$ l'induite de $A_{i}^{\tilde{M}}({\cal X}_{M},\omega)_{V}$. On a rappel\'e dans le paragraphe pr\'ec\'edent la forme lin\'eaire $r_{\tilde{M}}^{\tilde{G}}(.,\tilde{K}_{S-V})$ et on en a d\'efini divers avatars. Nous modifions l\'eg\`erement leur notation en rempla\c{c}ant $\tilde{K}_{S-V}$ par $\tilde{K}_{S}^V$.

Soit ${\cal X} \in {\bf Stab}(\tilde{G}(F))$.  Il r\'esulte de la d\'efinition de 1.9 et de [VI] 2.3(9) que l'on a l'\'egalit\'e
$$(2) \qquad A^{\tilde{G}}(V,{\cal X},\omega)=\sum_{\tilde{M}\in {\cal L}(\tilde{M}_{0})}\vert W^{\tilde{M}}\vert \vert W^{\tilde{G}}\vert ^{-1}\sum_{{\cal X}_{M}\in {\bf Stab}(\tilde{M}(F)), {\cal X}_{M}\mapsto {\cal X}}$$
$$\sum_{i=1,...,n({\cal X}_{M})}r_{\tilde{M}}^{\tilde{G}}(k_{i}^{\tilde{M}}({\cal X}_{M},\omega)_{S}^V,\tilde{K}_{S}^V)A_{i}^{\tilde{G}}({\cal X}_{M},\omega)_{V}.$$

 Supposons $(G,\tilde{G},{\bf a})$ quasi-d\'eploy\'e et \`a torsion int\'erieure. Pour $\tilde{M}\in {\cal L}(\tilde{M}_{0})$ et ${\cal X}_{M}\in {\bf Stab}(\tilde{M}(F))$, on a  d\'efini la distribution  $SA^{\tilde{M}}(S,{\cal X}_{M})$ en 1.10. Si $\tilde{M}\not=\tilde{G}$, elle est stable d'apr\`es le th\'eor\`eme 1.10(ii) et nos hypoth\`eses de r\'ecurrence. Si $\tilde{M}=\tilde{G}$, supposons qu'elle est stable. 
  On peut alors la d\'ecomposer comme en (1) en 
$$(3)\qquad S A^{\tilde{M}}(S,{\cal X}_{M})=\sum_{i=1,...,n({\cal X}_{M})}Sk_{i}^{\tilde{M}}({\cal X}_{M})_{S}^V\otimes SA_{i}^{\tilde{M}}({\cal X}_{M})_{V}$$
avec des $Sk_{i}^{\tilde{M}}({\cal X}_{M})_{S}^V\in  D^{st}_{g\acute{e}om}(\tilde{M}(F_{S}^V))$ et des $SA_{i}^{\tilde{M}}({\cal X}_{M})_{V}\in  D^{st}_{g\acute{e}om}(\tilde{M}(F_{V}))\otimes Mes(M(F_{V}))^*$.  On note $SA_{i}^{\tilde{G}}({\cal X}_{M})_{V}\in  D^{st}_{g\acute{e}om}(\tilde{G}(F_{V}))\otimes Mes(G(F_{V}))^*$ l'induite de $SA_{i}^{\tilde{M}}({\cal X}_{M})_{V}$.

Revenons au cas o\`u $(G,\tilde{G},{\bf a})$ est quelconque. Pour ${\cal X}\in {\bf Stab}(\tilde{G}(F))$, 
 on a d\'efini une distribution $A^{\tilde{G},{\cal E}}(S,{\cal X},\omega)$ en 1.10. On peut la d\'ecomposer de la m\^eme fa\c{c}on qu'en (1):
$$(4)  \qquad A^{\tilde{G},{\cal E}}(S,{\cal X},\omega)=\sum_{i=1,...,n({\cal X})}k_{i}^{\tilde{G},{\cal E}}({\cal X},\omega)_{S}^V\otimes A_{i}^{\tilde{G},{\cal E}}({\cal X},\omega)_{V}$$
avec des $k_{i}^{\tilde{G},{\cal E}}({\cal X}_{M},\omega)_{S}^V\in  D_{g\acute{e}om}(\tilde{G}(F_{S}^V),\omega)$ et des $A_{i}^{\tilde{G},{\cal E}}({\cal X},\omega)_{V}\in  D_{g\acute{e}om}(\tilde{G}(F_{V}),\omega)\otimes Mes(G(F_{V}))^*$.   

\ass{Proposition}{Soit ${\cal X}\in {\bf Stab}(\tilde{G}(F))$.

(i) On a l'\'egalit\'e
$$A^{\tilde{G},{\cal E}}(V,{\cal X},\omega)=\sum_{i=1,...,n({\cal X})}r^{\tilde{G}}(k_{i}^{\tilde{G},{\cal E}}({\cal X},\omega)_{S}^V,\tilde{K}_{S}^V)A_{i}^{\tilde{G},{\cal E}}({\cal X},\omega)_{V}$$
$$+\sum_{\tilde{M}\in {\cal L}(\tilde{M}_{0}),\tilde{M}\not=\tilde{G}}\vert W^{\tilde{M}}\vert \vert W^{\tilde{G}}\vert ^{-1}\sum_{{\cal X}_{M}\in {\bf Stab}(\tilde{M}(F)), {\cal X}_{M}\mapsto {\cal X}}\sum_{i=1,...,n({\cal X}_{M})}r_{\tilde{M}}^{\tilde{G}}(k_{i}^{\tilde{M}}({\cal X}_{M},\omega)_{S}^V,\tilde{K}_{S}^V)A_{i}^{\tilde{G}}({\cal X}_{M},\omega)_{V}.$$

(ii) Supposons que   $(G,\tilde{G},{\bf a})$ soit  quasi-d\'eploy\'e et \`a torsion int\'erieure et que $SA^{\tilde{G}}(S,{\cal X})$ soit stable. Alors on a l'\'egalit\'e
$$SA^{\tilde{G}}(V,{\cal X})=\sum_{\tilde{M}\in {\cal L}(\tilde{M}_{0})}\vert W^{\tilde{M}}\vert \vert W^{\tilde{G}}\vert ^{-1}\sum_{{\cal X}_{M}\in {\bf Stab}(\tilde{M}(F)), {\cal X}_{M}\mapsto {\cal X}}$$
$$\sum_{i=1,...,n({\cal X}_{M})}s_{\tilde{M}}^{\tilde{G}}(Sk_{i}^{\tilde{M}}({\cal X}_{M})_{S}^V,\tilde{K}_{S}^V)SA_{i}^{\tilde{G}}({\cal X}_{M})_{V}.$$}

 \bigskip

\subsection{Preuve de la proposition 2.3}
Prouvons le (i) de cette proposition. Si $(G,\tilde{G},{\bf a})$ est quasi-d\'eploy\'e et \`a torsion int\'erieure, on a par d\'efinition  les \'egalit\'es $A^{\tilde{G},{\cal E}}(V,{\cal X},\omega)=A^{\tilde{G}}(V,{\cal X},\omega)$ et $A^{\tilde{G},{\cal E}}(S,{\cal X},\omega)=A^{\tilde{G}}(S,{\cal X},\omega)$. La formule \`a prouver n'est autre que 2.3(2). On suppose maintenant que $(G,\tilde{G},{\bf a})$ n'est pas  quasi-d\'eploy\'e et \`a torsion int\'erieure.

Soient ${\bf G}'\in {\cal E}(\tilde{G},{\bf a},S)$ et $\tilde{M}'\in {\cal L}(\tilde{M}'_{0})$ (o\`u $\tilde{M}'_{0}$ est le Levi minimal fix\'e dans $\tilde{G}'$). Soit ${\cal X}_{M'}\in {\bf Stab}(\tilde{M}'(F))$. Supposons que $\tilde{M}'$ soit relevant. Alors il existe un espace de Levi $\tilde{M}$ de $\tilde{G}$ tel que $\tilde{M}'$ apparaisse comme l'espace associ\'e \`a une donn\'ee endoscopique elliptique ${\bf M}'$ de $(M,\tilde{M},{\bf a}_{M})$. Comme on l'a dit en 1.10, des formalit\'es permettent de d\'efinir des distributions $SA^{{\bf M}'}(V,{\cal X}_{M'})$ et $SA^{{\bf M}'}(S,{\cal X}_{M'})$. On peut d\'ecomposer cette derni\`ere par une formule similaire \`a 2.3(3):
$$(1) \qquad SA^{{\bf M}'}(S,{\cal X}_{M'})=\sum_{i=1,...,n({\cal X}_{M'})}Sk_{i}^{{\bf M}'}({\cal X}_{M'})_{S}^V\otimes SA_{i}^{{\bf M}'}({\cal X}_{M'})_{V}.$$
On note  $SA_{i}^{{\bf G}'}({\cal X}_{M'})_{V}$ l'induite \`a ${\bf G}'$ de $ SA_{i}^{{\bf M}'}({\cal X}_{M'})_{V}$. 

Soient ${\bf G}'\in {\cal E}(\tilde{G},{\bf a},V)$ et  ${\cal X}_{G'}\in {\bf Stab}(\tilde{G}(F))$. Montrons que l'on a l'\'egalit\'e
$$(2) \qquad transfert(SA^{{\bf G}'}(V,{\cal X}_{G'}))=\sum_{\tilde{M}'\in {\cal L}(\tilde{M}'_{0}),\tilde{M}'\text{ relevant }}\vert W^{M'}\vert \vert W^{G'}\vert ^{-1}$$
$$\sum_{{\cal X}_{M'}\in {\bf Stab}(\tilde{M}'(F)), {\cal X}_{M'}\mapsto {\cal X}_{G'}}\sum_{i=1,...,n({\cal X}_{M'})}s_{{\bf M}'}^{{\bf G}'}(Sk_{i}^{{\bf M}'}({\cal X}_{M'})_{S}^V,\tilde{K}_{S}^V)\,transfert( SA_{i}^{{\bf G}'}({\cal X}_{M'})_{V}).$$
Preuve. On fixe des donn\'ees auxiliaires $G'_{1}$ etc... pour ${\bf G}'$. Alors $SA^{{\bf G}'}(V,{\cal X}_{G'})$ s'identifie \`a une distribution que l'on peut noter $SA^{\tilde{G}'_{1}}_{\lambda_{1}}(V,{\cal X}_{G'})$. Gr\^ace \`a nos hypoth\`eses de r\'ecurrence et \`a quelques formalit\'es, on peut lui appliquer le (ii) de la proposition 2.3. On applique ensuite l'application de transfert \`a la formule obtenue. On obtient une formule similaire \`a celle ci-dessus. Plus pr\'ecis\'ement, on obtient une somme sur $\tilde{M}'\in {\cal L}(\tilde{M}'_{0})$ de certains termes, notons-les $X(\tilde{M}')$. On voit facilement que, si $\tilde{M}'$ est relevant,  $X(\tilde{M}')$ est \'egal au terme index\'e par $\tilde{M}'$ intervenant dans (2). Il faut montrer que, si $\tilde{M}'$ n'est pas relevant,  $X(\tilde{M}')$ est nul. En tout cas, $X(\tilde{M}')$ est un transfert d'un \'el\'ement de $D_{g\acute{e}om,\lambda_{1}}^{st}(\tilde{G}'_{1}(F_{V}))\otimes Mes(G'(F_{V}))^*$ qui est induit \`a partir d'un \'el\'ement de 
$D_{g\acute{e}om,\lambda_{1}}^{st}(\tilde{M}'_{1}(F_{V}))\otimes Mes(M'(F_{V}))^*$. Il s'ensuit que, pour que $X(\tilde{M}')$ soit nul, il suffit qu'il existe une place $v\in V$ telle que l'espace de Levi localis\'e $\tilde{M}'_{v}$ ne soit pas relevant. Par hypoth\`ese,  $\tilde{M}'$ n'est pas relevant. Par d\'efinition, cela signifie  soit que $\tilde{M}'(F)=\emptyset$, soit qu'il existe $v\in Val(F)$ tel que $\tilde{M}'_{v}$ n'est pas relevant. La premi\`ere possibilit\'e est exclue: ${\bf G}'$ est relevant, donc $\tilde{G}'(F)\not=\emptyset$ et alors $\tilde{M}'(F)\not=\emptyset$ puisque $\tilde{M}'$ est un espace de Levi de $\tilde{G}'$. Donc  il existe une place  $v\in Val(F)$ telle que $\tilde{M}'_{v}$ n'est pas relevant. Il reste \`a montrer qu'une telle place appartient \`a $V$. Pour $v\not\in V$, $G$  est quasi-d\'eploy\'e sur $F_{v}$ donc  il existe un espace de Levi $\tilde{M}_{v}$ de $\tilde{G}_{v}$ tel que $\tilde{M}_{v}'$ apparaisse comme l'espace associ\'e \`a une donn\'ee endoscopique elliptique ${\bf M}_{v}'$ de $(M_{v},\tilde{M}_{v},{\bf a}_{M_{v}})$. Parce que ${\bf G}'$ est non ramifi\'ee en $v$, ${\bf M}_{v}'$ l'est aussi. Alors ${\bf M}'_{v}$ est relevante d'apr\`es le lemme [I] 6.2. $\square$

On applique la formule de d\'efinition 1.10(3) et la formule (2) ci-dessus. On obtient
$$A^{\tilde{G},{\cal E}}(V,{\cal X},\omega)=\sum_{{\bf G}'\in {\cal E}(\tilde{G},{\bf a},V)}\sum_{{\cal X}_{G'}\in Stab(\tilde{G}'(F)),{\cal X}_{G'}\mapsto {\cal X}}i(\tilde{G},\tilde{G}')$$
$$\sum_{\tilde{M}'\in {\cal L}(\tilde{M}'_{0}),\tilde{M}'\text{ relevant }}\vert W^{M'}\vert \vert W^{G'}\vert ^{-1}\sum_{{\cal X}_{M'}\in {\bf Stab}(\tilde{M}'(F)), {\cal X}_{M'}\mapsto {\cal X}_{G'}}\sum_{i=1,...,n({\cal X}_{M'})}$$
$$s_{{\bf M}'}^{{\bf G}'}(Sk_{i}^{{\bf M}'}({\cal X}_{M'})_{S}^V,\tilde{K}_{S}^V)\,transfert( SA_{i}^{{\bf G}'}({\cal X}_{M'})_{V}).$$
Pour un \'el\'ement ${\bf G}'\in {\cal E}(\tilde{G},{\bf a})$ et un espace de Levi $\tilde{M}'$ de $\tilde{G}'$, que l'on identifiera dans la notation \`a une "donn\'ee de Levi" ${\bf M}'$, d\'efinissons un terme $S({\bf G}',{\bf M}')$ de la fa\c{c}on suivante. Si ${\bf G}'$ n'est pas non ramifi\'ee hors de $V$ ou si ${\bf M}'$ n'est pas relevant, on pose $S({\bf G}',{\bf M}')=0$. Si ${\bf G}'$ est non ramifi\'ee hors de $V$ et si ${\bf M}'$ est relevant, on pose
$$S({\bf G}',{\bf M}')=\sum_{{\cal X}_{G'}\in Stab(\tilde{G}'(F)),{\cal X}_{G'}\mapsto {\cal X}}\sum_{{\cal X}_{M'}\in {\bf Stab}(\tilde{M}'(F)), {\cal X}_{M'}\mapsto {\cal X}_{G'}}\sum_{i=1,...,n({\cal X}_{M'})}s_{{\bf M}'}^{{\bf G}'}(Sk_{i}^{{\bf M}'}({\cal X}_{M'})_{S}^V,\tilde{K}_{S}^V)$$
$$transfert( SA_{i}^{{\bf G}'}({\cal X}_{M'})_{V}).$$
La formule ci-dessus se r\'ecrit
$$A^{\tilde{G},{\cal E}}(V,{\cal X},\omega)=\sum_{{\bf G}'\in {\cal E}(\tilde{G},{\bf a},V)}i(\tilde{G},\tilde{G}')\sum_{\tilde{M}'\in {\cal L}(\tilde{M}'_{0})}\vert W^{M'}\vert \vert W^{G'}\vert ^{-1}S({\bf G}',{\bf M}').$$
On peut appliquer la proposition [VI] 6.5. La formule de cette proposition fait appara\^{\i}tre des Levi $\hat{M}$ du groupe dual $\hat{G}$. Cela parce que l'on consid\'erait une situation g\'en\'erale o\`u le terme $S({\bf G}',{\bf M}')$ pouvait \^etre non nul m\^eme si ${\bf M}'$ n'\'etait pas relevant. Ici, seuls peuvent appara\^{\i}tre des  ${\bf M}'$ qui sont relevants et des $\hat{M}$ correspondant \`a des espaces de Levi $\tilde{M} $ de $\tilde{G}$.   On peut r\'ecrire cette proposition en sommant sur de tels espaces $\tilde{M}$ et de telles donn\'ees endoscopiques de $(M,\tilde{M},{\bf a}_{M})$. On obtient
$$(3) \qquad A^{\tilde{G},{\cal E}}(V,{\cal X},\omega)=\sum_{\tilde{M}\in {\cal L}(\tilde{M}_{0})}\vert W^{\tilde{M}}\vert \vert W^{\tilde{G}}\vert ^{-1}\sum_{{\bf M}'\in {\cal E}(\tilde{M},{\bf a}_{M})}i(\tilde{M},\tilde{M}')$$
$$\sum_{\tilde{s}\in \tilde{\zeta}Z(\hat{M})^{\Gamma_{F},\hat{\theta}}/Z(\hat{G})^{\Gamma_{F},\hat{\theta}}}i_{\tilde{M}'}(\tilde{G},\tilde{G}'(\tilde{s}))S({\bf G}'(\tilde{s}),{\bf M}').$$
On a not\'e simplement $\tilde{\zeta}$ le terme tel que ${\bf M}'=(M',{\cal M}',\tilde{\zeta})$. On peut limiter la somme en ${\bf M}'$ aux \'el\'ements de ${\cal E}(\tilde{M},{\bf a}_{M},V)$. En effet, si ${\bf M}'$ n'est pas ramifi\'e hors de $V$, les donn\'ees ${\bf G}'(\tilde{s})$ apparaissant ne sont pas non plus non ramifi\'ees hors de $V$ et les termes $S({\bf G}'(\tilde{s}),{\bf M}')$ sont nuls. D'autre part, fixons $\tilde{M}$, ${\bf M}'$ et $\tilde{s}$ intervenant ci-dessus. On a un diagramme commutatif
$$\begin{array}{ccc}{\bf Stab}(\tilde{M}'(F))& \to&{\bf Stab}(\tilde{G}'(\tilde{s};F))\\ \downarrow&&\downarrow\\ {\bf Stab}(\tilde{M}(F))&\to&{\bf Stab}(\tilde{G}(F))\\ \end{array}$$
Dans la d\'efinition de $S({\bf G}'(\tilde{s}),{\bf M}')$, on peut donc remplacer la double somme en ${\cal X}_{G'(\tilde{s})}$ tel que ${\cal X}_{G'(\tilde{s})}\mapsto {\cal X}$ et en ${\cal X}_{M'}$ tel que ${\cal X}_{M'}\mapsto {\cal X}_{G'(\tilde{s})}$ par une double somme sur ${\cal X}_{M}\in {\bf Stab}(\tilde{M}(F))$ tel que ${\cal X}_{M}\mapsto {\cal X}$ et sur ${\cal X}_{M'}\in {\bf Stab}(\tilde{M}'(F))$ tel que ${\cal X}_{M'}\mapsto {\cal X}_{M}$. On a aussi l'\'egalit\'e $transfert(SA_{i}^{{\bf G}'(\tilde{s})}({\cal X}_{M'})_{V})=(transfert(SA_{i}^{{\bf M}'}({\cal X}_{M'})_{V}))^{\tilde{G}}$.  On obtient
$$S({\bf G}'(\tilde{s}),{\bf M}')=\sum_{{\cal X}_{M}\in {\bf Stab}(\tilde{M}(F)), {\cal X}_{M}\mapsto {\cal X}}\sum_{{\cal X}_{M'}\in {\bf Stab}(\tilde{M}'(F)),{\cal X}_{M'}\mapsto {\cal X}_{M}}$$
$$\sum_{i=1,...,n({\cal X}_{M'})}s_{{\bf M}'}^{{\bf G}'(\tilde{s})}(Sk_{i}^{{\bf M}'}({\cal X}_{M'})_{S}^V,\tilde{K}_{S}^V)(transfert(SA_{i}^{{\bf M}'}({\cal X}_{M'})_{V}))^{\tilde{G}}.$$
Pour ${\cal X}_{M'}\in {\bf Stab}(\tilde{M}'(F))$, posons
$$(4) \qquad b({\bf M}',{\cal X}_{M'})=\sum_{i=1,...,n({\cal X}_{M'})}(transfert(SA_{i}^{{\bf M}'}({\cal X}_{M'})_{V}))^{\tilde{G}}$$
$$\sum_{\tilde{s}\in \tilde{\zeta}Z(\hat{M})^{\Gamma_{F},\hat{\theta}}/Z(\hat{G})^{\Gamma_{F},\hat{\theta}}}i_{\tilde{M}'}(\tilde{G},\tilde{G}'(\tilde{s}))s_{{\bf M}'}^{{\bf G}'(\tilde{s})}(Sk_{i}^{{\bf M}'}({\cal X}_{M'})_{S}^V,\tilde{K}_{S}^V).$$
Pour ${\cal X}_{M}\in {\bf Stab}(\tilde{M}(F))$, posons
$$(5) \qquad B(\tilde{M},{\cal X}_{M})=\sum_{{\bf M}'\in {\cal E}(\tilde{M},{\bf a}_{M},V)}i(\tilde{M},\tilde{M}')\sum_{{\cal X}_{M'}\in {\bf Stab}(\tilde{M}'(F)),{\cal X}_{M'}\mapsto {\cal X}_{M}}b({\bf M}',{\cal X}_{M'}).$$
Les consid\'erations ci-dessus permettent de r\'ecrire l'\'egalit\'e (3) sous la forme
$$(6) \qquad   A^{\tilde{G},{\cal E}}(V,{\cal X},\omega)=\sum_{\tilde{M}\in {\cal L}(\tilde{M}_{0})}\vert W^{\tilde{M}}\vert \vert W^{\tilde{G}}\vert ^{-1}\sum_{{\cal X}_{M}\in {\bf Stab}(\tilde{M}(F)), {\cal X}_{M}\mapsto {\cal X}}B(\tilde{M},{\cal X}_{M}).$$
Dans la formule (4), on reconna\^{\i}t la somme en $\tilde{s}$: elle est \'egale \`a $r_{\tilde{M}}^{\tilde{G},{\cal E}}({\bf M}',Sk_{i}^{{\bf M}'}({\cal X}_{M'})_{S}^V,\tilde{K}_{S}^V)$. On applique la proposition 2.2(i) et on obtient
$$(7)\qquad  b({\bf M}',{\cal X}_{M'})=\sum_{i=1,...,n({\cal X}_{M'})}r_{\tilde{M}}^{\tilde{G}}(transfert(Sk_{i}^{{\bf M}'}({\cal X}_{M'})_{S}^V),\tilde{K}_{S}^V)$$
$$(transfert(SA_{i}^{{\bf M}'}({\cal X}_{M'})_{V}))^{\tilde{G}}.$$
On a suppos\'e ici ${\bf M}'$ non ramifi\'e hors de $V$. Mais le membre de droite ci-dessus conserve un sens si ${\bf M}'$ est seulement non ramifi\'e hors de $S$. Pour un tel ${\bf M}'$, on d\'efinit $b({\bf M}',{\cal X}_{M'})$ par l'\'egalit\'e (7). Si ${\bf M}'$ est non ramifi\'e hors de $S$ mais pas hors de $V$, on a $b({\bf M}',{\cal X}_{M'})=0$: cela r\'esulte de la proposition 2.2(ii). Dans la d\'efinition (5), on  peut donc \'etendre la somme en ${\bf M}'\in {\cal E}(\tilde{M},{\bf a}_{M},V)$ en une somme en ${\bf M}'\in {\cal E}(\tilde{M},{\bf a}_{M},S)$. C'est-\`a-dire
$$(8) \qquad B(\tilde{M},{\cal X}_{M})=\sum_{{\bf M}'\in {\cal E}(\tilde{M},{\bf a}_{M},S)}i(\tilde{M},\tilde{M}')\sum_{{\cal X}_{M'}\in {\bf Stab}(\tilde{M}'(F)),{\cal X}_{M'}\mapsto {\cal X}_{M}}b({\bf M}',{\cal X}_{M'}).$$

On note $\underline{A}^{\tilde{G},{\cal E}}(V,{\cal X},\omega)$ le membre de droite de l'\'egalit\'e du (i) de la proposition 2.3. Pour $\tilde{M}\in {\cal L}(\tilde{M}_{0})$ et ${\cal X}_{M}\in {\bf Stab}(\tilde{M}(F))$, posons

- si $\tilde{M}\not=\tilde{G}$, 
$$\underline{B}(\tilde{M},{\cal X}_{M})=\sum_{i=1,...,n({\cal X}_{M})}r_{\tilde{M}}^{\tilde{G}}(k_{i}^{\tilde{M}}({\cal X}_{M},\omega)_{S}^V,\tilde{K}_{S}^V)A_{i}^{\tilde{G}}({\cal X}_{M},\omega)_{V};$$

- si $\tilde{M}=\tilde{G}$, 
$$\underline{B}(\tilde{G},{\cal X}_{G})=\sum_{i=1,...,n({\cal X}_{G})}r^{\tilde{G},{\cal E}}(k_{i}^{\tilde{G},{\cal E}}({\cal X}_{G},\omega)_{S}^V,\tilde{K}_{S}^V)A_{i}^{\tilde{G},{\cal E}}({\cal X}_{G},\omega)_{V}.$$

On a alors
$$(9) \qquad \underline{A}^{\tilde{G},{\cal E}}(V,{\cal X},\omega)=\sum_{\tilde{M}\in {\cal L}(\tilde{M}_{0})}\vert W^{\tilde{M}}\vert \vert W^{\tilde{G}}\vert ^{-1}\underline{B}(\tilde{M},{\cal X}_{M}).$$
 Fixons $\tilde{M}\in {\cal L}(\tilde{M}_{0})$ et ${\cal X}_{M}\in {\bf Stab}(\tilde{M}(F))$. Par d\'efinition, on a
 $$A^{\tilde{M},{\cal E}}(S,{\cal X}_{M},\omega)=\sum_{{\bf M}'\in {\cal E}(\tilde{M},{\bf a}_{M},V)}i(\tilde{M},\tilde{M}')\sum_{{\cal X}_{M'}\in {\bf Stab}(\tilde{M}'(F)),{\cal X}_{M'}\mapsto {\cal X}_{M}}transfert(SA^{{\bf M}'}(S,{\cal X}_{M'})).$$
 En utilisant (1), on obtient
 $$A^{\tilde{M},{\cal E}}(S,{\cal X}_{M},\omega)=\sum_{{\bf M}'\in {\cal E}(\tilde{M},{\bf a}_{M},S)}i(\tilde{M},\tilde{M}')\sum_{{\cal X}_{M'}\in {\bf Stab}(\tilde{M}'(F)),{\cal X}_{M'}\mapsto {\cal X}_{M}}$$
 $$\sum_{i=1,...,n({\cal X}_{M'})}transfert(Sk_{i}^{{\bf M}'}({\cal X}_{M'})_{S}^V)transfert(SA_{i}^{{\bf M}'}({\cal X}_{M'})_{V}).$$
 Si $\tilde{M}\not=\tilde{G}$, on a $A^{\tilde{M},{\cal E}}(S,{\cal X}_{M},\omega)=A^{\tilde{M}}(S,{\cal X}_{M},\omega)$ d'apr\`es les hypoth\`eses de r\'ecurrence. Alors la d\'ecomposition ci-dessus est de la forme 2.3(1): l'ensemble d'indices $\{1,...,n({\cal X}_{M})\}$ est la r\'eunion disjointe des $\{1,...,n({\cal X}_{M'})\}$ sur les ${\bf M}'$ et ${\cal X}_{M'}$ intervenant ci-dessus. Si $\tilde{M}=\tilde{G}$, cette d\'ecomposition est de m\^eme de la forme 2.3(4). Il en r\'esulte par d\'efinition que
 $$\underline{B}(\tilde{M},{\cal X}_{M})= \sum_{{\bf M}'\in {\cal E}(\tilde{M},{\bf a}_{M},S)}i(\tilde{M},\tilde{M}')\sum_{{\cal X}_{M'}\in {\bf Stab}(\tilde{M}'(F)),{\cal X}_{M'}\mapsto {\cal X}_{M}}$$
 $$\sum_{i=1,...,n({\cal X}_{M'})}r_{\tilde{M}}^{\tilde{G}}(transfert(Sk_{i}^{{\bf M}'}({\cal X}_{M'})_{S}^V),\tilde{K}_{S}^V)(transfert(SA_{i}^{{\bf M}'}({\cal X}_{M'})_{V}))^{\tilde{G}}.$$
 En utilisant (7) et (8), on voit que
 $$\underline{B}(\tilde{M},{\cal X}_{M})=B(\tilde{M}, {\cal X}_{M}).$$
 Alors les membres de droite de (6) et (9) co\"{\i}ncident. Cela prouve l'\'egalit\'e
$$(10) \qquad A^{\tilde{G},{\cal E}}(V,{\cal X},\omega)=\underline{A}^{\tilde{G},{\cal E}}(V,{\cal X},\omega),$$
ce qui est le (i) de la proposition 2.3.  

Prouvons maintenant le (ii) de cette proposition. On suppose $(G,\tilde{G},{\bf a})$ quasi-d\'eploy\'e et \`a torsion int\'erieure. Soit ${\cal X}\in {\bf Stab}(\tilde{G}(F))$. La formule 1.10(3) se modifie en
$$A^{\tilde{G},{\cal E}}(V,{\cal X})=SA^{\tilde{G}}(V,{\cal X})+\sum_{{\bf G}'\in {\cal E}(\tilde{G},{\bf a},V), {\bf G}'\not={\bf G}}$$
$$\sum_{{\cal X}_{G'}\in {\bf Stab}(\tilde{G}'(F)),{\cal X}_{G'}\mapsto {\cal X}}i(\tilde{G},\tilde{G}')transfert(SA^{{\bf G}'}(V,{\cal X}_{G'})).$$
Pour ${\bf G}'\not={\bf G}$, on a encore l'\'egalit\'e (2) et on peut remplacer le terme $transfert(SA^{{\bf G}'}(V,{\cal X}_{G'}))$ ci-dessus par le membre de droite de cette \'egalit\'e. Pour ${\bf G}'={\bf G}$, le membre de droite de  (2)  est \'egal au membre de droite de l'\'egalit\'e du (ii) de la proposition 2.3. On ne sait pas qu'il est \'egal \`a $SA^{\tilde{G}}(V,{\cal X})$, c'est ce qu'on veut prouver. Mais on peut remplacer $SA^{\tilde{G}}(V,{\cal X})$ par le membre de droite de (2), plus un nombre $C$ dont on veut prouver qu'il est nul. Le calcul se poursuit comme pr\'ec\'edemment et on obtient l'\'egalit\'e
$$A^{\tilde{G},{\cal E}}(V,{\cal X})=C+ \underline{A}^{\tilde{G},{\cal E}}(V,{\cal X}).$$
Mais, comme on l'a dit au d\'ebut de la preuve, dans notre situation quasi-d\'eploy\'ee \`a torsion int\'erieure, le (i) de la proposition 2.3, c'est-\`a-dire l'\'egalit\'e (10),
est tautologique. On conclut $C=0$, ce qui prouve le (ii) de la proposition 2.3. 
\bigskip
\bigskip

\subsection{Extension de l'ensemble fini de places}
\ass{Corollaire}{(i) Soit $V$ un ensemble fini de places de $F$ contenant $V_{ram}$ et soit ${\cal X}\in {\bf Stab}(\tilde{G}(F))$. Supposons qu'il existe un ensemble fini $S$ de places de $F$ contenant $V$ et tel que l'assertion du th\'eor\`eme 1.10(i)  soit v\'erifi\'ee pour le couple $(S,{\cal X})$. Alors cette assertion   est v\'erifi\'ee pour le couple $(V,{\cal X})$.

(ii) Supposons $(G,\tilde{G},{\bf a})$ quasi-d\'eploy\'e et \`a torsion int\'erieure. Alors la m\^eme propri\'et\'e vaut pour le th\'eor\`eme 1.10(ii). }

Preuve. Si $A^{\tilde{G},{\cal E}}(S,{\cal X},\omega)=A^{\tilde{G}}(S,{\cal X},\omega)$, on peut supposer que les d\'ecompositions 2.3 (1) et 2.3 (4)   de cette distribution co\"{\i}ncident. L'\'egalit\'e 2.3(2) et celle du (i) de la proposition 2.3 entra\^{\i}nent alors l'\'egalit\'e $A^{\tilde{G},{\cal E}}(V,{\cal X},\omega)=A^{\tilde{G}}(V,{\cal X},\omega)$. Supposons $(G,\tilde{G},{\bf a})$ quasi-d\'eploy\'e et \`a torsion int\'erieure. L'hypoth\`ese du (ii) de la proposition 2.3 est v\'erifi\'ee. La formule de cette proposition exprime $SA^{\tilde{G}}(V,{\cal X})$ comme combinaison lin\'eaire de distributions stables. Donc $SA^{\tilde{G}}(V,{\cal X})$ est stable. La distribution $SA^{\tilde{G}}(S,{\cal X})$ ne d\'epend que des classes de conjugaison par $G_{AD}(F_{v})$ des $\tilde{K}_{v}$ pour $v\not\in S$. Elle ne d\'epend pas des $\tilde{K}_{v}$ pour $v\in S-V
$. Pour un espace de Levi $\tilde{M}\in {\cal L}(\tilde{M}_{0})$ et pour ${\cal X}_{M}\in {\bf Stab}(\tilde{M}(F))$, la distribution $SA^{\tilde{M}}(S,{\cal X}_{\tilde{M}})$ v\'erifie les m\^emes propri\'et\'es. En effet, elle ne d\'epend que des classes de conjusaison par $M_{AD}(F_{v})$ des  $\tilde{K}_{v}\cap \tilde{M}(F_{v})$ pour $v\not\in S$. On a vu dans la preuve de [II] 4.2(3) que, si $\tilde{K}_{v}$ et $\tilde{K}'_{v}$ sont conjugu\'es par un \'el\'ement de $G_{AD}(F_{v})$ et sont tous deux en bonne position relativement \`a $M$, alors  les espaces $\tilde{K}_{v}\cap \tilde{M}(F_{v})$ et $\tilde{K}'_{v}\cap \tilde{M}(F_{v})$ sont conjugu\'es par un \'el\'ement de $M_{ad}(F_{v})$. L'assertion s'ensuit. Alors, pour $v\not\in S$, la formule du (ii) de la proposition 2.3 ne d\'epend que de la   classe de conjugaison par $G_{AD}(F_{v})$ de $\tilde{K}_{v}$. Pour $v\in S-V$, la formule ne d\'epend des $\tilde{K}_{v}$ que par les formes lin\'eaires $s_{\tilde{M}}^{\tilde{G}}(.,\tilde{K}_{S}^V)$. Or, d'apr\`es 2.2(1), celles-ci ne d\'ependent que  des classes de conjugaison par $G_{AD}(F_{v})$ des $\tilde{K}_{v}$ pour $v\in S-V$. Finalement, $SA^{\tilde{G}}(V,{\cal X})$ ne d\'epend que des classes de conjugaison par $G_{AD}(F_{v})$ des $\tilde{K}_{v}$ pour $v\not\in V$. $\square$

\bigskip

  \section{Enonc\'es de nouveaux  th\'eor\`emes}

\bigskip

\subsection{Le th\'eor\`eme d'Arthur}
Supposons ici $G=\tilde{G}$, ${\bf a}=1$, $\tilde{K}_{v}=K_{v}$ pour tout $v\not\in V$.

\ass{Th\'eor\`eme}{Sous ces hypoth\`eses, les th\'eor\`emes 1.10(i) et (ii) sont v\'erifi\'es.}

C'est l'un des principaux r\'esultats de l'article d'Arthur ([A1] global theorem 1').  La preuve que nous donnerons  des th\'eor\`emes 1.10(i) et (ii) \'etant directement inspir\'ee de celle d'Arthur, nous pourrions aussi bien les red\'emontrer enti\`erement. Mais cela n'aurait aucun int\'er\^et. Nous pr\'ef\'erons simplifier un peu la n\^otre en utilisant le r\'esultat d'Arthur. La propri\'et\'e  1.10(1) de $SA^{G}(V,{\cal X})$ n'est pas clairement \'enonc\'ee par Arthur, mais est incluse dans sa d\'emonstration. On la retrouve en tout cas de la fa\c{c}on suivante. Consid\'erons d'autres sous-groupes hypersp\'eciaux $K'_{v}$ pour $v\not\in V$, soumis aux conditions de [VI] 1.1. Ces conditions impliquent qu'il existe un ensemble fini de places $S$ contenant $V$ tel que $K'_{v}=K_{v}$ pour $v\not\in S$.  La distribution $SA^G(S,{\cal X})$ ne change donc pas quand on remplace les $K_{v}$ par les $K'_{v}$. On sait qu'elle est stable d'apr\`es le th\'eor\`eme d'Arthur. De plus, pour $v\in S-V$, les compacts $K_{v}$ et $K'_{v}$ sont conjugu\'es par un \'el\'ement de $G_{AD}(F_{v})$: il en est ainsi pour tout couple de sous-groupes compacts hypersp\'eciaux. La preuve du corollaire 2.5 montre que la distribution $SA^G(V,{\cal X})$ ne change  pas non plus quand on remplace les $K_{v}$ par les $K'_{v}$.

En fait, nous n'utiliserons le th\'eor\`eme d'Arthur que  pour l'\'el\'ement ${\cal X}$ correspondant \`a la classe de conjugaison stable r\'eduite \`a $\{1\}$. Dans ce cas, on note plut\^ot nos distributions $A^{\tilde{G},{\cal E}}_{unip}(V)$ et $SA^{\tilde{G}}_{unip}(V)$.

  \bigskip
  
  \subsection{D\'efinition d'une autre distribution stable}
   
  On suppose que $(G,\tilde{G},{\bf a})$ est quasi-d\'eploy\'e et \`a torsion int\'erieure. Soient ${\cal X}\in {\bf Stab}(\tilde{G}(F))$ et $V$ un ensemble fini de places contenant $V_{ram}$.   On a vu que ${\cal X}$ correspondait \`a une classe de conjugaison stable dans $\tilde{G}_{ss}(F)$. On sait qu'il existe un \'el\'ement $\epsilon$ de cette classe telle que $G_{\epsilon}$ soit quasi-d\'eploy\'e. On fixe un tel $\epsilon$. Soit $U_{\epsilon}$ un voisinage ouvert de l'unit\'e dans $G_{\epsilon}(F_{V})$ qui v\'erifie les conditions suivantes:
  
  $x\in U_{\epsilon}$ si et seulement si sa partie semi-simple $x_{ss}$ appartient \`a $U_{\epsilon}$;
  
si $x\in U_{\epsilon}$ et $y\in G_{\epsilon}(F_{V})$ sont conjugu\'es par un \'el\'ement de $Z_{G}(\epsilon;\bar{F}_{V})$ (o\`u $\bar{F}_{V}=\prod_{v\in V}\bar{F}_{v}$), alors $y\in U_{\epsilon}$.
 
 On note $\tilde{U}$ l'ensemble des \'el\'ements de $\tilde{G}(F_{V})$ dont la partie semi-simple est stablement conjugu\'ee \`a un \'el\'ement de $U_{\epsilon}\epsilon$. On note $SI(\tilde{U})$, resp. $SI(U_{\epsilon})$, le sous-espace des \'el\'ements de $SI(\tilde{G}(F_{V}))$, resp. $SI(G_{\epsilon}(F_{V}))$, \`a support dans $\tilde{U}$, resp. $U_{\epsilon}$. 
 On pose $\Xi_{\epsilon}=Z_{G}(\epsilon)/G_{\epsilon}$. Ce groupe est naturellement muni d'une action galoisienne. On a \'etabli en [I] lemme 4.8 un isomorphisme de descente $desc_{\epsilon}^{st}:SI(\tilde{U})\otimes Mes(G(F_{V}))\simeq SI(U_{\epsilon})^{\Xi_{\epsilon}^{\Gamma_{F_{V}}}} \otimes Mes(G_{\epsilon}(F_{V}))$ pourvu que $U_{\epsilon}$ soit assez petit. Dans cette r\'ef\'erence, on avait fix\'e les mesures mais l'isomorphisme devient canonique quand on l'\'ecrit sous la forme ci-dessus.
  Rappelons la caract\'erisation de l'isomorphisme. Soit $x\in U_{\epsilon}$ tel que $x\epsilon\in \tilde{G}_{reg}(F_{V})$. Fixons une mesure de Haar sur $(G_{\epsilon})_{x}(F_{V})=G_{x\epsilon}(F_{V})$.  Rappelons que la donn\'ee de $x$  et de la mesure d\'efinit un \'el\'ement de $D_{orb}^{st}(U_{\epsilon})  \otimes Mes(G_{\epsilon}(F_{V}))^*$, \`a savoir l'int\'egrale orbitale stable $S^{G_{\epsilon}}(x,.)$. De m\^eme, la donn\'ee de $x\epsilon$ et de la mesure d\'efinit une int\'egrale orbitale stable $S^{\tilde{G}}(x\epsilon,.)$. Soient alors 
 ${\bf f}\in SI(\tilde{U})\otimes Mes(G(F_{V}))$ et $ {\bf f}_{\epsilon}=desc_{\epsilon}^{st}({\bf f})$.  
 On a l'\'egalit\'e
$$S^{G_{\epsilon}}(x,{\bf f}_{\epsilon})=S^{\tilde{G}}(x\epsilon,{\bf f}).$$

Les distributions de 1.10 d\'ependent d'une mesure sur $\mathfrak{A}_{G}$ fix\'ee en [VI] 1.3. On doit aussi fixer 
 une mesure sur $\mathfrak{A}_{G_{\epsilon}}$. Dans le cas g\'en\'eral, le choix est arbitraire. Mais, si on suppose $\epsilon$ elliptique, on a $\mathfrak{A}_{G_{\epsilon}}=\mathfrak{A}_{G}$ et on choisit la mesure d\'ej\`a fix\'ee sur ce dernier espace.
 
 Rappelons que, pour tout groupe r\'eductif connexe $H$ d\'efini sur $F$, on pose
 $$\tau(H)=\vert \pi_{0}(Z(\hat{H})^{\Gamma_{F}})\vert \vert ker^1(F,Z(\hat{G}))\vert ^{-1},$$
 cf. [VI] 5.1. 
Supposons

(1) $S({\cal X})\subset V$.

On d\'efinit une distribution $\underline{SA}^{\tilde{G}}(V,{\cal X})\in D_{g\acute{e}om}(\tilde{G}(F_{V}))\otimes Mes(G(F_{V}))^*$ de la fa\c{c}on suivante. Si ${\cal X}$ n'est pas elliptique ou si $V$ ne contient pas $S({\cal X},\tilde{K})$, on pose $\underline{SA}^{\tilde{G}}(V,{\cal X})=0$. Supposons

(2) ${\cal X}$ est elliptique et $S({\cal X},\tilde{K})\subset V$.

Soit $ {\bf f}\in C_{c}^{\infty}(\tilde{G}(F_{V}))\otimes Mes(G(F_{V}))$. On restreint $f$ \`a $\tilde{U}$. On consid\`ere l'image dans $SI(\tilde{U})\otimes Mes(G(F_{V}))$ de cette restriction. On note $ {\bf f}_{\epsilon}\in SI(U_{\epsilon})^{\Xi_{\epsilon}^{\Gamma_{F_{V}}}}\otimes Mes(G_{\epsilon}(F_{V}))$  l'image de l'\'el\'ement obtenu par l'isomorphisme  $desc_{\epsilon}^{st}$. On pose
$$I^{\tilde{G}}(\underline{SA}^{\tilde{G}}(V,{\cal X}),{\bf f})= \vert \Xi_{\epsilon}^{\Gamma_{F}}\vert \tau(G)\tau(G_{\epsilon})^{-1} S^{G_{\epsilon}}(SA^{G_{\epsilon}}_{unip}(V),{\bf f}_{\epsilon}).$$

Rappelons que l'on sait que $SA^{G_{\epsilon}}_{unip}(V)$ est stable et ind\'ependante de tout choix de sous-groupes compacts hypersp\'eciaux, cf. 3.1. Cela  donne un sens \`a cette d\'efinition. Puisque $SA^{G_{\epsilon}}_{unip}(V)$ est \`a support unipotent, la d\'efinition ne d\'epend pas du choix de $U_{\epsilon}$. Elle ne d\'epend pas non plus du choix de $\epsilon$. En effet, rempla\c{c}ons $\epsilon$ par $\epsilon'$ v\'erifiant les m\^emes propri\'et\'es. On peut fixer $y\in G$ tel que $y^{-1}\epsilon y=\epsilon'$ et $y\sigma(y)^{-1}\in G_{\epsilon}$ pour tout $\sigma\in \Gamma_{F}$. Parce que l'on suppose $G_{\epsilon}$ et $G_{\epsilon'}$ quasi-d\'eploy\'es, on peut supposer que $ad_{y^{-1}}$ envoie une paire de Borel \'epingl\'ee d\'efinie sur $F$ de $G_{\epsilon}$ sur une telle paire de de $G_{\epsilon'}$. Cela entra\^{\i}ne que $y\sigma(y)^{-1}$ appartient au centre de $G_{\epsilon}$. Alors l'isomorphisme $ad_{y^{-1}}:G_{\epsilon}\to G_{\epsilon'}$ est d\'efini sur $F$. Cet isomorphisme identifie $SA^{G_{\epsilon}}_{unip}(V)$ \`a $SA^{G_{\epsilon'}}_{unip}(V)$ et ${\bf f}_{\epsilon}$ \`a ${\bf f}_{\epsilon'}$. L'ind\'ependance affirm\'ee s'ensuit.

Remarquons que

(3) $\underline{SA}^{\tilde{G}}(V,{\cal X})$ est stable.  

En effet, si l'image de ${\bf f}$ dans $SI(\tilde{G}(F_{V}))\otimes Mes(G(F_{V}))$ est nulle, alors ${\bf f}_{\epsilon}=0$. L'assertion s'ensuit. 

La distribution $\underline{SA}^{\tilde{G}}(V,{\cal X})$ d\'epend \'evidemment de diverses donn\'ees. Mais, quant \`a sa d\'ependance des espaces hypersp\'eciaux $\tilde{K}_{v}$, on a

(4) $\underline{SA}^{\tilde{G}}(V,{\cal X})$ ne d\'epend que des classes de conjugaison par $G_{AD}(F_{v})$ des $\tilde{K}_{v}$ pour $v\not\in V$.

En effet, la d\'efinition ci-dessus ne d\'epend des $\tilde{K}_{v}$ pour $v\not\in V$ que par la condition $S({\cal X},\tilde{K})\subset V$. Or, d'apr\`es la d\'efinition de l'ensemble $S({\cal X},\tilde{K})$, cette condition  ne d\'epend que des classes de conjugaison par $G_{AD}(F_{v})$ des $\tilde{K}_{v}$ pour $v\not\in V$. 

\ass{Th\'eor\`eme  }{Pour tout ${\cal X}\in {\bf Stab}(\tilde{G}(F))$ et tout $V$ contenant $S({\cal X})$, on a l'\'egalit\'e
$$\underline{SA}^{\tilde{G}}(V,{\cal X})=SA^{\tilde{G}}(V,{\cal X}).$$}

 Cela sera d\'emontr\'e en 3.4.  

\bigskip

\subsection{Enonc\'e du th\'eor\`eme principal}
En [III] 6.3, on a d\'efini certains triplets $(G,\tilde{G},{\bf a})$ particuliers. Dans cette r\'ef\'erence, le corps de base \'etait local non-archim\'edien, mais les d\'efinitions et r\'esultats valent aussi bien sur  notre corps de nombres. Consid\'erons un tel triplet. Rappelons que $G$ est quasi-d\'eploy\'e sur $F$ et simplement connexe. On a ${\bf a}=1$. Notons $\Theta_{F}$ l'ensemble des $\eta\in \tilde{G}(F)$ tels qu'il existe une paire de Borel \'epingl\'ee ${\cal E}$ de $G$ d\'efinie sur $F$ de sorte que $ad_{\eta}$ conserve ${\cal E}$. Cet ensemble n'est pas vide. L'ensemble des classes de conjugaison stable contenant un \'el\'ement de $\Theta_{F}$, que l'on note $\Theta_{F}/st-conj$, est en bijection avec ${\cal Z}(\tilde{G})^{\Gamma_{F}}$. De plus, l'application naturelle
 $${\cal Z}(\tilde{G})\to (T^*/(1-\theta^*)(T^*))\times_{{\cal Z}(G)}{\cal Z}(\tilde{G})$$
 est injective, cf. preuve de [III] 6.3(3). On voit alors que l'ensemble $\Theta_{F}/st-conj$ est param\'etr\'e par le sous-ensemble fini ${\bf Stab}_{excep}(\tilde{G}(F))$ de ${\bf Stab}(\tilde{G}(F))$  d\'efini de la fa\c{c}on suivante. C'est l'ensemble des $(\mu,\omega_{\bar{G}})$ tels que $\mu$ appartienne \`a l'image de ${\cal Z}(\tilde{G})^{\Gamma_{F}}$ par l'application pr\'ec\'edente. On a alors $W(\mu)=W^{\theta^*}$ donc $\omega_{\bar{G}}$ est forc\'ement trivial. L'indice $excep$ signifie exceptionnel. Si $(G,\tilde{G},{\bf a})$ n'est pas l'un des triplets d\'efinis en [III] 6.3, on pose ${\bf Stab}_{excep}(\tilde{G}(F))=\emptyset$.

Consid\'erons un triplet $(G,\tilde{G},{\bf a})$ quelconque. 

\ass{Th\'eor\`eme }{Soient ${\cal X}\in {\bf Stab}(\tilde{G}(F))$ et $V$ un ensemble fini de places contenant $S({\cal X})$. On suppose ${\cal X}\not\in {\bf Stab}_{excep}(\tilde{G}(F))$. Alors on a l'\'egalit\'e
$$A^{\tilde{G}}(V,{\cal X},\omega)=\sum_{{\bf G}'\in {\cal E}(\tilde{G},{\bf a},V)}\sum_{{\cal X}'\in {\bf Stab}(G'(F)); {\cal X}'\mapsto {\cal X} }i(\tilde{G},\tilde{G}') transfert(\underline{SA}^{{\bf G}'}(V,{\cal X}')).$$}

D'apr\`es 1.7(3), on a $S({\cal X}')\subset S({\cal X})\subset V$ pour tout ${\cal X}'$ intervenant ci-dessus.  Alors les termes  $\underline{SA}^{{\bf G}'}(V,{\cal X}')$ sont d\'eduits de ceux d\'efinis au paragraphe pr\'ec\'edent par les constructions formelles habituelles. La d\'emonstration du th\'eor\`eme  occupe les sections 5 \`a 8.

\bigskip

\subsection{Le th\'eor\`eme 3.3 implique les  th\'eor\`emes  3.2, 1.10(ii) et [VI] 5.2 }
 On suppose que 
  $(G,\tilde{G},{\bf a})$ est quasi-d\'eploy\'e et \`a torsion int\'erieure. En particulier, ce n'est pas l'un des triplets d\'efinis en [III] 6.3 et l'ensemble $ {\bf Stab}_{excep}(\tilde{G}(F))$ est vide. Soient ${\cal X}\in {\bf Stab}(\tilde{G}(F))$ et $V$ un ensemble fini de places contenant $S({\cal X})$. Soient ${\bf G}'\in {\cal E}(\tilde{G},V)$ tel que ${\bf G}'\not={\bf G}$ et  ${\cal X}'\in {\bf Stab}(\tilde{G}'(F))$ tel que ${\cal X}'\mapsto {\cal X}$. Comme on l'a dit, on a $S({\cal X}')\subset S({\cal X})\subset V$ et on peut appliquer le th\'eor\`eme 3.2 par r\'ecurrence: $\underline{SA}^{{\bf G}'}(V,{\cal X}')=SA^{{\bf G}'}(V,{\cal X}')$.
L'\'egalit\'e du th\'eor\`eme 3.3 se r\'ecrit
$$A^{\tilde{G}}(V,{\cal X})=\underline{SA}^{\tilde{G}}(V,{\cal X})+\sum_{{\bf G}'\in {\cal E}(\tilde{G},V),{\bf G}'\not={\bf G}}\sum_{{\cal X}'\in {\bf Stab}(G'(F)); {\cal X}'\mapsto {\cal X} }i(\tilde{G},\tilde{G}') transfert(SA^{{\bf G}'}(V,{\cal X}')).$$
En utilisant la d\'efinition 1.10(2), cela entra\^{\i}ne
$$\underline{SA}^{\tilde{G}}(V,{\cal X})=SA^{\tilde{G}}(V,{\cal X}),$$
ce qui prouve le th\'eor\`eme 3.2. 

  Gr\^ace au th\'eor\`eme 3.2, les propri\'et\'es 3.2(3) et 3.2(4) de la distribution  $\underline{SA}^{\tilde{G}}(V,{\cal X})$  impliquent l'assertion du th\'eor\`eme 1.10(ii) sous la restriction $S({\cal X})\subset V$. Le corollaire 2.5 permet de supprimer cette restriction, d'o\`u le th\'eor\`eme 1.10(ii).

Soit $V$ un ensemble fini de places de $F$ contenant $V_{ram}$. En utilisant le th\'eor\`eme 1.10(ii), la proposition 1.12(iii)  entra\^{\i}ne que $SA^{\tilde{G}}(V,{\cal O}_{V})$ est stable pour tout ${\cal O}_{V}\in \tilde{G}_{ss}(F_{V})/st-conj$ et qu'elle ne d\'epend que des classes de conjugaison par $G_{AD}(F_{v})$ des $\tilde{K}_{v}$ pour $v\not\in V$.  Ce sont les assertions du th\'eor\`eme [VI] 5.2.

 A ce point, on peut remarquer que le th\'eor\`eme 3.2 permet de retrouver la formule habituelle pour une distribution  $SA^{\tilde{G}}(V,{\cal O}_{V})$ associ\'ee \`a une classe de conjugaison stable elliptique et fortement r\'eguli\`ere. Pr\'ecis\'ement, consid\'erons un \'el\'ement $\delta\in \tilde{G}(F)$ elliptique et fortement r\'egulier. Notons ${\cal X}$ le param\`etre de sa classe de conjugaison stable. 
 Fixons un ensemble fini $V$ de places de $F$ contenant $S({\cal X}, \tilde{K})$. Notons ${\cal O}_{V}$ la classe de conjugaison stable de $\delta$ dans $\tilde{G}(F_{V})$, fixons un ensemble de repr\'esentants $\dot{{\cal Y}}_{\delta}$ des classes de conjugaison par $G(F_{V})$ contenues dans ${\cal O}_{V}$. Fixons une mesure de Haar $dx$ sur $G_{\delta}(F_{V})$. Pour $\delta'\in \dot{{\cal Y}}_{\delta}$, $G_{\delta'}(F_{V})$ est isomorphe \`a $G_{\delta}(F_{V})$ et on munit le premier groupe de la mesure correspondant \`a  $dx$. 
 Soient $f\in C_{c}^{\infty}(\tilde{G}(F_{V}))$ et $dg$ une mesure de Haar sur $G(F_{V})$. Alors on a l'\'egalit\'e
 $$(1) \qquad S^{\tilde{G}}(SA^{\tilde{G}}(V,{\cal O}_{V}),f\otimes dg)=\tau(G)\tau(G_{\delta})^{-1}mes(\mathfrak{A}_{G}G_{\delta}(F)\backslash G_{\delta}({\mathbb A}_{F}))$$
 $$ \sum_{\delta'\in \dot{{\cal Y}}_{\delta}}\int_{G_{\delta'}( F_{V})\backslash G(F_{V})}f(g^{-1}\delta'g)\,dg.$$
 
 Preuve. Notons ${\cal X}_{V}$ l'image naturelle de ${\cal X}$ dans ${\bf Stab}(\tilde{G}(F_{V}))$. La condition de forte r\'egularit\'e impos\'ee \`a $\delta$ implique que l'ensemble de sommation de la proposition 1.12(iii) est r\'eduit \`a $\{{\cal X}\}$. Donc $SA^{\tilde{G}}(V,{\cal O}_{V})=SA^{\tilde{G}}(V,{\cal X})$. Gr\^ace au th\'eor\`eme 3.2, ceci est \'egal \`a $\underline{SA}^{\tilde{G}}(V,{\cal X})$. On utilise la d\'efinition de ce terme. Par forte r\'egularit\'e, le groupe $\Xi_{\delta}$ est r\'eduit \`a $\{1\}$. On obtient
 $$S^{\tilde{G}}(SA^{\tilde{G}}(V,{\cal O}_{V}),f\otimes dg)=\tau(G)\tau(G_{\delta})^{-1}S^{G_{\delta}}(SA_{unip}^{G_{\delta}}(V),(f\otimes dg)_{\delta}).$$
 Puisque $G_{\delta}$ est un tore, on a $SA_{unip}^{G_{\delta}}(V)=A_{unip}^{G_{\delta}}(V)$.  Ecrivons $(f\otimes dg)_{\delta}=\varphi\otimes dx$. D'apr\`es [VI] 2.2, on a 
 $$I^{G_{\delta}}(A_{unip}^{G_{\delta}}(V),\varphi\otimes dx)=mes(\mathfrak{A}_{G}G_{\delta}(F)\backslash G_{\gamma}({\mathbb A}_{F}))\varphi(1).$$
 Il r\'esulte de la d\'efinition de l'application de descente que 
 $$\varphi(1)=\sum_{\delta'\in \dot{{\cal Y}}_{\delta}}\int_{G_{\delta'}( F_{V})\backslash G(F_{V})}f(g^{-1}\delta'g)\,dg.$$
 En mettant ces calculs bout \`a bout, on obtient (1). $\square$
 
  Comme on l'expliquera en 4.1, si la mesure $dx$ et  la mesure sur $\mathfrak{A}_{G}$ sont les mesures de Tamagawa, on a l'\'egalit\'e $mes(\mathfrak{A}_{G}G_{\delta}(F)\backslash G_{\delta}({\mathbb A}_{F})) =\tau(G_{\delta})$. Alors la formule (1) se simplifie en
 $$S^{\tilde{G}}(SA^{\tilde{G}}(V,{\cal O}_{V}),f\otimes dg)=\tau(G) \sum_{\delta'\in \dot{{\cal Y}}_{\delta}}\int_{G_{\delta'}( F_{V})\backslash G(F_{V})}f(g^{-1}\delta'g)\,dg.$$
 On retrouve ainsi les formules de [KS] et [Lab3].

\bigskip

\subsection{Le th\'eor\`eme 3.3 implique presque les th\'eor\`emes 1.10(i) et [VI] 5.4}
Soit $V$ un ensemble fini de places de $F$ contenant $V_{ram}$. On d\'efinit l'ensemble 

\noindent ${\bf Stab}_{excep}(\tilde{G}(F_{V}))$ de la m\^eme fa\c{c}on qu'en 3.3. Il est vide si $(G,\tilde{G},{\bf a})$ n'est pas l'un des triplets d\'efinis en [III] 6.3. Si $(G,\tilde{G},{\bf a})$ est l'un de ces triplets, ${\bf Stab}_{excep}(\tilde{G}(F_{V}))$ param\`etre les classes de conjugaison stable dans $\tilde{G}(F_{V})$ d'\'el\'ements $\eta_{V}=(\eta_{v})_{v\in V}\in \tilde{G}_{ss}(F_{V})$ tels que, pour tout $v\in V$, il existe une paire de Borel \'epingl\'ee ${\cal E}_{v}$ de $G$ d\'efinie sur $F_{v}$ de sorte que $ad_{\eta_{v}}$ conserve ${\cal E}_{v}$. L'ensemble ${\bf Stab}_{excep}(\tilde{G}(F_{V}))$ est en bijection avec $\prod_{v\in V}{\cal Z}(\tilde{G})^{\Gamma_{F_{v}}}$. C'est un ensemble fini. 

\ass{Proposition}{(i) Le th\'eor\`eme 3.3 implique l'assertion du th\'eor\`eme 1.10(i) pour ${\cal X}\not\in {\bf Stab}_{excep}(\tilde{G}(F))$. 

(ii) Le th\'eor\`eme 3.3 implique l'assertion du th\'eor\`eme [VI] 5.4 pour toute classe ${\cal O}_{V}\in \tilde{G}_{ss}(F_{V})/st-conj$ qui est param\'etr\'ee par un \'el\'ement de ${\bf Stab}(\tilde{G}(F_{V}))$ qui n'appartient pas \`a ${\bf Stab}_{excep}(\tilde{G}(F_{V}))$.}

Preuve. On peut supposer que $(G,\tilde{G},{\bf a})$ n'est pas quasi-d\'eploy\'e et \`a torsion int\'erieure, sinon les th\'eor\`emes 1.10(i) et [VI] 5.4 sont tautologiques. Soit ${\cal X}\in {\bf Stab}(\tilde{G}(F))-{\bf Stab}_{excep}(\tilde{G}(F))$, supposons d'abord $S({\cal X})\subset V$. Comme dans le paragraphe pr\'ec\'edent, les hypoth\`eses de r\'ecurrence permettent d'appliquer le th\'eor\`eme 3.2, cette fois pour tout ${\bf G}'$. Alors l'\'egalit\'e du th\'eor\`eme 3.3 devient
$$A^{\tilde{G}}(V,{\cal X},\omega)=\sum_{{\bf G}'\in {\cal E}(\tilde{G},{\bf a},V)}\sum_{{\cal X}'\in {\bf Stab}(G'(F)); {\cal X}'\mapsto {\cal X} }i(\tilde{G},\tilde{G}') transfert(SA^{{\bf G}'}(V,{\cal X}')).$$
En comparant avec 1.10(3), on obtient
$$ A^{\tilde{G}}(V,{\cal X},\omega)=A^{\tilde{G},{\cal E}}(V,{\cal X},\omega).$$
Cela d\'emontre le th\'eor\`eme 1.10(i) sous la restriction $S({\cal X})\subset V$. Celle-ci dispara\^{\i}t gr\^ace au corollaire 2.5. Cela prouve (i).

L'application ${\bf Stab}(\tilde{G}(F))\to {\bf Stab}(\tilde{G}(F_{V}))$ envoie ${\bf Stab}_{excep}(\tilde{G}(F))$ dans ${\bf Stab}_{excep}(\tilde{G}(F_{V}))$. Pour une classe ${\cal O}_{V}$ qui est param\'etr\'ee par un \'el\'ement de ${\bf Stab}(\tilde{G}(F_{V}))$ qui n'est pas exceptionnel, les ${\cal X}$ qui interviennent dans les formules (i) et (ii) de la proposition 1.12 ne sont donc pas exceptionnels. Le th\'eor\`eme 1.10(i) d\'ej\`a d\'emontr\'e pour les \'el\'ements non exceptionnels implique que les membres de droite de ces deux formules sont \'egaux. D'o\`u l'\'egalit\'e des membres de gauche, ce qui est l'assertion du th\'eor\`eme [VI] 5.4. $\square$

\bigskip

\subsection{Le th\'eor\`eme [VI] 5.4 implique le th\'eor\`eme 1.10(i) et \'etend le th\'eor\`eme 3.3}
A la fin du pr\'esent article, nous aurons d\'emontr\'e le th\'eor\`eme 3.3, avec ses cons\'equences d\'ecrites dans les deux paragraphes pr\'ec\'edents. Nous compl\'eterons ult\'erieurement la preuve du th\'eor\`eme [VI] 5.4, c'est-\`a-dire nous le d\'emontrerons pour les classes ${\cal O}_{V}$ param\'etr\'ees par des \'el\'ements de ${\bf Stab}_{excep}(\tilde{G}(F_{V}))$. Montrons que cela suffira pour compl\'eter la preuve du th\'eor\`eme 1.10(i). Soit ${\cal X}\in {\bf Stab}(\tilde{G}(F))$. Si ${\cal X}$ n'est pas exceptionnel, on a vu ci-dessus que l'assertion du th\'eor\`eme 1.10(i) r\'esultait du th\'eor\`eme 3.3. Supposons ${\cal X}\in {\bf Stab}_{excep}(\tilde{G}(F))$. Cet \'el\'ement param\`etre une classe de conjugaison stable ${\cal O}$. Il s'envoie sur un \'el\'ement ${\cal X}_{V}\in {\bf Stab}_{excep}(\tilde{G}(F_{V}))$, qui param\`etre la classe ${\cal O}_{V}$. En g\'en\'eral, l'application ${\bf Stab}(\tilde{G}(F))\to {\bf Stab}(\tilde{G}(F_{V}))$ n'est pas injective. Mais ici, la fibre de cette application au-dessus de ${\cal X}_{V}$ est r\'eduite \`a $\{{\cal X}\}$.   En effet, le d\'efaut d'injectivit\'e est d\^u au fait qu'un cocycle $\omega_{\bar{G}}$ n'est pas toujours d\'etermin\'e par ses restrictions $\omega_{\bar{G},v}$ pour $v\in V$. Ici, ${\cal X}$ est l'image d'un couple $(\mu,\omega_{\bar{G}})$ tel que $\mu\in {\cal Z}(\tilde{G})^{\Gamma_{F}}$. On a $W(\mu)=W^{\theta^*}$ tout entier et tout cocycle $\omega'_{\bar{G}}$ compl\'etant $\mu$ en un \'el\'ement $(\mu,\omega'_{\bar{G}})\in Stab(\tilde{G}(F))$ est forc\'ement trivial. Les assertions (i) et (ii) de la proposition 1.12 se r\'eduisent aux \'egalites
$$A^{\tilde{G}}(V,{\cal O}_{V},\omega)=A^{\tilde{G}}(V,{\cal X},\omega),$$
$$A^{\tilde{G},{\cal E}}(V,{\cal O}_{V},\omega)=A^{\tilde{G},{\cal E}}(V,{\cal X},\omega).$$
Le th\'eor\`eme [VI] 5.4 affirme l'\'egalit\'e des deux membres de gauche. D'o\`u l'\'egalit\'e des membres de droite, ce qui est l'assertion du th\'eor\`eme 1.10(i). 

Sous la m\^eme hypoth\`ese, l'\'egalit\'e du th\'eor\`eme 3.3 est valable pour tout ${\cal X}\in {\bf Stab}(\tilde{G}(F))$. En effet, comme dans le paragraphe pr\'ec\'edent, le membre de droite de cette \'egalit\'e est \'egal \`a $A^{\tilde{G},{\cal E}}(V,{\cal X},\omega)$. Le th\'eor\`eme 1.10(i) affirme que ce terme est \'egal \`a $A^{\tilde{G}}(V,{\cal X},\omega)$.

 \bigskip

\subsection{Quelques cas faciles}

Soient ${\cal X}\in {\bf Stab}(\tilde{G}(F))$ et $V$ contenant $S({\cal X})$. On consid\`ere les hypoth\`eses suivantes:

(1) $V$ contient $S({\cal X},\tilde{K})$;

(2) ${\cal X}$ est elliptique;

(3) pour toute place $v\in Val(F)$, l'image de ${\cal X}$ dans ${\bf Stab}(\tilde{G}(F_{v}))$ appartient \`a l'image de l'application $\chi^{\tilde{G}_{v}}$;

(4) les restrictions de $\omega$ \`a $Z(G;{\mathbb A}_{F})^{\theta}$ et \`a $Z(\bar{G};{\mathbb A}_{F})$ sont triviales. 

\ass{Lemme}{Si l'une de ces hypoth\`eses n'est pas satisfaite, l'\'egalit\'e  du th\'eor\`eme 3.3 est v\'erifi\'ee, les deux membres \'etant nuls.}

 Preuve. Les assertions 1.9 (1), (2)  et (3) nous disent que $A^{\tilde{G}}(V,{\cal X},\omega)=0$  si (1), resp (2), (4),  n'est pas v\'erifi\'ee. On a aussi  $A^{\tilde{G}}(V,{\cal X},\omega)=0$ par d\'efinition si (3) n'est pas v\'erifi\'ee car a fortiori ${\cal X}$ n'est pas dans l'image de $\chi^{\tilde{G}}$.  
 
 Consid\'erons maintenant le membre de droite de l'\'egalit\'e . Soient ${\bf G}'$ et ${\cal X}'$ y intervenant. Supposons $transfert(SA^{{\bf G}'}(V,{\cal X}'))\not=0$. On va prouver que les conditions (1) \`a (4) sont alors v\'erifi\'ees. 
  Les relations 1.7(3) et (4) nous disent que $S({\cal X}')\subset S({\cal X})\subset V$ et $S({\cal X},\tilde{K})-S({\cal X})\subset S({\cal X}',\tilde{K}')-S({\cal X}')$. La non-nullit\'e de $SA^{{\bf G}'}(V,{\cal X}')$ entra\^{\i}ne $S({\cal X}',\tilde{K}')\subset V $ d'apr\`es 3.2(2). Ces inclusions entra\^{\i}nent (1). De m\^eme, 3.2(2) implique que ${\cal X}'$ est elliptique, d'o\`u (2) 
  d'apr\`es 1.7(2). Notons ${\cal O}'$ la classe de conjugaison stable dans $\tilde{G}'(F)$ associ\'ee \`a ${\cal X}'$. Puisque $\underline{SA}^{{\bf G}'}(V,{\cal X}')$ est \`a support dans la classe de conjugaison stable ${\cal O}'_{V}$ dans $\tilde{G}'(F_{V})$  engendr\'ee par ${\cal O}'$, la non-nullit\'e de  $transfert(\underline{SA}^{{\bf G}'}(V,{\cal X}'))$ entra\^{\i}ne que cette classe correspond \`a une classe de conjugaison stable dans $\tilde{G}(F_{V})$. D'apr\`es le lemme 1.8, il existe pour tout $v\in V$ une classe ${\cal O}_{v}\in \tilde{G}_{ss}(F_{v})/st-conj$ qui corresponde \`a ${\cal O}'_{v}$ et dont l'image par  $\chi^{\tilde{G}_{v}} $ soit l'image de ${\cal X}$ dans ${\bf Stab}(\tilde{G}(F_{v}))$. Puisqu'on a d\'ej\`a prouv\'e que (1) \'etait v\'erifi\'ee, le lemme 1.8(2) entra\^{\i}ne la m\^eme propri\'et\'e pour $v\not\in V$. Cela entra\^{\i}ne (3). Cela entra\^{\i}ne aussi que, pour toute place $v$, on peut fixer un diagramme $(\epsilon,B',T',B,T,\eta)$ d\'efini sur $F_{v}$, avec $\epsilon\in \tilde{G}'_{ss}(F_{v})$ et $\eta\in {\cal O}_{v}$. D'apr\`es [KS] lemme 4.4.C, $\omega$ est trivial sur $T^{\theta,0}(F_{v})$. A fortiori $\omega$ est trivial sur le sous-groupe $Z(G_{\eta};F_{v})\subset T^{\theta,0}(F_{v})$. Cela \'equivaut \`a ce que la restriction de $\omega$ \`a $Z(\bar{G};F_{v})$ soit triviale. Enfin,  l'existence de ${\bf G}'$ entra\^{\i}ne que $\omega$ est trivial sur $Z(G;{\mathbb A}_{F})^{\theta}$ (cf. par exemple [I] lemme 2.7). Cela v\'erifie la condition (4). $\square$

 \bigskip

\section{Distributions \`a support unipotent}

\bigskip

\subsection{ Mesures de Tamagawa}
Dans ce paragraphe, $G$ est un groupe r\'eductif connexe d\'efini sur $F$. Pour toute place finie $v$ de $F$, on note $\mathfrak{o}_{v}$ l'anneau des entiers de $F_{v}$, $\mathfrak{p}_{v}$ son id\'eal maximal et $ {\mathbb F}_{v}$ le corps r\'esiduel. Pour toute place $v\in Val(F)$, on munit $F_{v}$ d'une mesure de Haar de sorte que $mes(\mathfrak{o}_{v})=1$ pour presque toute place finie $v$. Le produit de ces mesures est une mesure sur ${\mathbb A}_{F}$. On suppose que $mes({\mathbb A}_{F}/F)=1$. Fixons une forme diff\'erentielle de degr\'e maximal sur $\mathfrak{g}$, d\'efinie sur $F$ et non nulle.  Pour toute  place $v$, on d\'eduit de cette forme diff\'erentielle et de la mesure sur $F_{v}$ une mesure $dX_{v}$ sur $\mathfrak{g}(F_{v})$. Rappelons que l'ensemble des mesures de Haar sur $\mathfrak{g}(F_{v})$ s'identifie \`a celui des mesures de Haar sur $G(F_{v})$: deux mesures se correspondent si et seulement si le jacobien  de l'application exponentielle, calcul\'e pour ces mesures,  vaut $1$ au point $0\in \mathfrak{g}(F_{v})$. On a donc aussi une mesure $dg_{v}$ sur $G(F_{v})$. Fixons un ensemble fini $V$ de places de $F$, contenant les places archim\'ediennes, de sorte que $G$ soit non ramifi\'e hors de $V$. Notons $\rho_{G}$ la repr\'esentation de $\Gamma_{F}$ dans $X^*(G)\otimes_{{\mathbb Z}}{\mathbb C}$. On note $L_{v}(\rho_{G},s)$ sa fonction $L$ en la place $v$ et $L^V(\rho_{G},s)$ sa fonction $L$ partielle hors de $V$.  Notons $r$ l'ordre du p\^ole en $1$ de la fonction $L^V(\rho_{G},s)$ et posons 
$$\ell^V_{G}=lim_{s\to 1}(s-1)^{r}L^V(\rho_{G},s).$$
  La mesure de Tamagawa sur $G({\mathbb A}_{F})$ est  par d\'efinition  \'egale \`a
$$(\ell_{G}^V)^{-1}(\prod_{v\not\in V}L_{v}(\rho_{G},1)dg_{v})\otimes (\prod_{v\in V}dg_{v}),$$
ce produit \'etant convergent.  

 Conform\'ement \`a nos d\'efinitions de [VI] 1.1, il se d\'eduit de la mesure de Tamagawa  sur $G({\mathbb A}_{F})$ une mesure  $dg_{V}^{Tam}$ sur $G(F_{V})$. En effet, fixons pour tout $v\not\in V$ un sous-groupe compact hypersp\'ecial $K_{v}$ de $G(F_{v})$. On a une mesure canonique $dg^{can}_{v}$ sur $G(F_{v})$ telle que la mesure de $K_{v}$ soit $1$. Alors $dg_{V}^{Tam}$ est la mesure telle que $dg_{V}^{Tam}\otimes \otimes_{v\not\in V}dg_{v}^{can}$ soit la mesure de Tamagawa sur $G({\mathbb A}_{F})$.  D'une fa\c{c}on g\'en\'erale, si $X$ est un ensemble muni d'une mesure $dx$ et si $Y$ est un sous-ensemble mesurable de $X$, notons $mes(Y,dx)$ la mesure de $Y$. La d\'efinition entra\^{\i}ne que $dg_{V}^{Tam}$ est \'egale \`a
 $$(1) \qquad (\ell_{G}^V)^{-1}(\prod_{v\not\in V}L_{v}(\rho_{G},1)mes(K_{v},dg_{v}))\prod_{v\in V}dg_{v}.$$
 Soit $v\not\in V$. On sait qu'\`a $K_{v}$ est associ\'e un sch\'ema en groupes ${\cal K}_{v}$ sur $\mathfrak{o}_{v}$. On note $\mathfrak{k}_{v}$ le groupe des points sur $\mathfrak{o}_{v}$ de son alg\`ebre de Lie. C'est une  sous-$\mathfrak{o}_{v}$-alg\`ebre   de $\mathfrak{g}(F_{v})$. On note ${\mathbb K}_{v}$ la fibre r\'esiduelle de ${\cal K}_{v}$. On a un homomorphisme surjectif
 $ K_{v}\to {\mathbb K}_{v}({\mathbb F}_{v})$
 dont on note $K_{v}^1$ le noyau. Alors 
 l'exponentielle se restreint en un isomorphisme $\mathfrak{p}_{v}\mathfrak{k}_{v}\to K_{v}^1$ qui pr\'eserve les mesures.
  La formule pr\'ec\'edente se r\'ecrit sous la forme suivante: notre mesure $dg_{V}^{Tam}$ sur $G(F_{V})$ est \'egale \`a
$$(2) \qquad (\ell_{G}^V)^{-1}(\prod_{v\not\in V}L_{v}(\rho_{G},1)\vert {\mathbb K}_{v}({\mathbb F}_{v})\vert mes(\mathfrak{p}_{v}\mathfrak{k}_{v},dX_{v}))\prod_{v\in V}dg_{v}.$$

{\bf Jusqu'\`a la fin de l'article, pour tout groupe $G$ et tout ensemble $V$ de places comme ci-dessus, on munit $G({\mathbb A}_{F})$ et $G(F_{V})$ des mesures de Tamagawa.} Cela nous d\'ebarrasse des espaces de mesures. 

On   a aussi besoin d'une mesure sur $\mathfrak{A}_{G}$. Rappelons la normalisation habituelle dans le cadre des mesures de Tamagawa. Identifions $\mathfrak{A}_{G}$ \`a $Hom(X^*(G)^{\Gamma_{F}},{\mathbb R})$, o\`u $X^*(G)$ est  le groupe des caract\`eres alg\'ebriques de $G$. On d\'efinit le r\'eseau $\mathfrak{A}_{G,{\mathbb Z}}=Hom(X^*(G)^{\Gamma_{F}},{\mathbb Z})$. La mesure "de Tamagawa" sur $\mathfrak{A}_{G}$ est celle pour laquelle  ce r\'eseau est de covolume $1$. Si on munit $G({\mathbb A}_{F})$ de la mesure de Tamagawa et que l'on choisit cette mesure  sur $\mathfrak{A}_{G}$, on sait que la mesure de
$$\mathfrak{A}_{G}G(F)\backslash G({\mathbb A}_{F})$$
 est \'egale au terme $\tau(G)$ d\'efini en 3.2.
Mais cette  normalisation est peu commode. Par exemple, elle  n'est pas compatible avec la situation de 4.2(2) ci-dessous. On fixe donc la mesure sur $\mathfrak{A}_{G}$ sans supposer qu'il s'agit de la mesure de Tamagawa. On note $covol(\mathfrak{A}_{G,{\mathbb Z}})$ le covolume de ce r\'eseau et $\tau'(G)$ la mesure de $\mathfrak{A}_{G}G(F)\backslash G({\mathbb A}_{F})$ calcul\'ee \`a l'aide de notre mesure sur $\mathfrak{A}_{G}$. On a alors l'\'egalit\'e
$$(3) \qquad  \tau'(G)=covol(\mathfrak{A}_{G,{\mathbb Z}})^{-1}\tau(G).$$

\bigskip

\subsection{Compatibilit\'e des mesures}
Consid\'erons les deux situations suivantes.

(1) On se donne une suite exacte
$$1\to C_{1}\to G_{1}\to G\to 1$$
de groupes r\'eductifs connexes d\'efinis sur $F$, o\`u $C_{1}$ est un tore central induit. On a une suite exacte
$$1\to \mathfrak{A}_{C_{1}}\to \mathfrak{A}_{G_{1}}\to \mathfrak{A}_{G}\to 1.$$
On suppose qu'elle est compatible aux mesures. On fixe un ensemble fini $V$ de places de $F$ tel que les trois groupes soient non ramifi\'es hors de $V$.

(2) On se donne une suite exacte
$$1\to \Xi\to G_{1}\times G_{2}\to G\to 1$$
o\`u $G_{1}$, $G_{2}$ et $G$ sont des groupes r\'eductifs connexes d\'efinis sur $F$ et $\Xi$ est un sous-groupe fini central de $G_{1}\times G_{2}$. On fixe un ensemble fini $V$ de places de $F$ tel que $G_{1}$, $G_{2}$ et $G$ soient non ramifi\'es hors de $V$ et tel que le nombre d'\'el\'ements de $X$ soit premier \`a $p$ pour tout nombre premier $p$ divisant une place hors de $V$.

\ass{Lemme}{(i) Dans la situation (1), la suite exacte
$$1\to C_{1}(F_{V})\to G_{1}(F_{V})\to G(F_{V})\to 1$$
est compatible aux mesures. On a l'\'egalit\'e
$\tau'(G_{1})=\tau'(C_{1})\tau'(G)$. 

(ii) Dans la situation (2), le rev\^etement
$$ G_{1}(F_{V})\times G_{2}(F_{V})\to G(F_{V})$$
pr\'eserve localement les mesures.}

Preuve. On effectue les constructions du paragraphe pr\'ec\'edent en  adaptant les notations de fa\c{c}on \'evidente. 

 Consid\'erons la situation (1). Fixons un isomorphisme d'alg\`ebres de Lie
$$\mathfrak{g}_{1}=\mathfrak{c}_{1}\oplus \mathfrak{g}.$$
On fixe des formes diff\'erentielles de degr\'e maximal sur $\mathfrak{c}_{1}$ et $\mathfrak{g}$, d\'efinies sur $F$ et non nulles. Par produit tensoriel, on en d\'eduit une telle forme sur $\mathfrak{g}_{1}$. Pour chaque place $v$, on associe \`a ces formes des mesures de Haar sur $C_{1}(F_{v})$, $G_{1}(F_{v})$ et $G(F_{v})$ comme dans le paragraphe pr\'ec\'edent. La suite
$$1\to C_{1}(F_{v})\to G_{1}(F_{v})\to G(F_{v})\to 1$$
est compatible \`a ces mesures. On a aussi $\rho_{G_{1}}=\rho_{C_{1}}\oplus \rho_{G}$. On en d\'eduit que la suite
$$1\to C_{1}({\mathbb A}_{F})\to G_{1}({\mathbb A}_{F})\to G({\mathbb A}_{F})\to 1$$
est compatible aux mesures de Tamagawa. La derni\`ere assertion du (i) r\'esulte alors de [Oe] th\'eor\`eme 5.3 (dont l'\'enonc\'e se simplifie gr\^ace \`a l'hypoth\`ese que $C_{1}$ est induit). Pour $v\not\in V$, le sous-groupe compact hypersp\'ecial $K_{C_{1},v}$ de $C_{1}(F_{v})$ est unique et on peut choisir des sous-groupes compacts hypersp\'eciaux dans $G_{1}(F_{v})$ et $G(F_{v})$ qui fixent le m\^eme point hypersp\'ecial de l'immeuble commun du groupe $G_{1,AD}=G_{AD}$. 
Alors la suite
$$1\to {\mathbb K}_{C_{1},v}\to {\mathbb K}_{1,v}\to  {\mathbb K}_{v}\to 1$$
est exacte. Le th\'eor\`eme de Lang entra\^{\i}ne que la suite d\'eduite 
 $$1\to {\mathbb K}_{C_{1},v}({\mathbb F}_{v})\to {\mathbb K}_{1,v}({\mathbb F}_{v})\to {\mathbb K}_{v}({\mathbb F}_{v})\to 1$$
 est exacte. La suite
 $$0\to \mathfrak{k}_{C_{1},v}\to \mathfrak{k}_{1,v}\to \mathfrak{k}_{v}\to 0$$
 est aussi exacte. De ces deux faits, on d\'eduit par un argument familier 
  que la suite
$$1\to K_{C_{1},v}\to K_{1,v}\to K_{v}\to 1$$
est aussi exacte.  D'apr\`es la compatibilit\'e des mesures, on a $mes(K_{C_{1},v},dc_{1,v})mes(K,dg)=mes(K_{1},dg_{1})$. La premi\`ere assertion de (i) r\'esulte alors de la formule (1) du paragraphe pr\'ec\'edent. 

Consid\'erons la situation (2). On a un isomorphisme
$$\mathfrak{g}_{1}\oplus \mathfrak{g_{2}}\simeq \mathfrak{g}.$$
On peut supposer que la forme diff\'erentielle fix\'ee sur $\mathfrak{g}$ est le produit des formes diff\'erentielles fix\'ees sur $\mathfrak{g}_{1}$ et sur $\mathfrak{g}_{2}$. Pour toute place $v$, l'isomorphisme
$$\mathfrak{g}_{1}(F_{v})\oplus \mathfrak{g}_{2}(F_{v})\simeq \mathfrak{g}(F_{v})$$
est alors compatible aux mesures d\'eduites de ces formes. Donc le rev\^etement
$$G_{1}(F_{v})\times G_{2}(F_{v})\to G(F_{v})$$
pr\'eserve localement ces mesures. Pour prouver le (ii) du lemme, il suffit de prouver que la constante figurant dans la formule (2)  du paragraphe pr\'ec\'edent pour $G$ est le produit des constantes pour $G_{1}$ et pour $G_{2}$.   On a l'\'egalit\'e $\rho_{G}=\rho_{G_{1}}\oplus \rho_{G_{2}}$,  les constantes provenant des fonctions $L$ sont donc compatibles.  Soit $v\not\in V$.  On peut de nouveau supposer que les sous-groupes hypersp\'eciaux de $G_{1}(F_{v})\times G_{2}(F_{v})$ et de $G(F_{v})$ fixent le m\^eme point hypersp\'ecial de l'immeuble commun du groupe $G_{AD}$. L'hypoth\`ese que $\vert \Xi\vert $ est premier \`a la caract\'eristique r\'esiduelle de $F_{v}$ entra\^{\i}ne que l'on a l'\'egalit\'e
$$\mathfrak{k}_{v}=\mathfrak{k}_{1,v}\oplus \mathfrak{k}_{2,v}$$ et que la suite
$$1\to \Xi\to {\mathbb K}_{1,v}\times {\mathbb K}_{2,v}\to {\mathbb K}_{v}\to 1$$
est exacte. L'\'egalit\'e ci-dessus entra\^{\i}ne
$$(1) \qquad mes(\mathfrak{p}_{v}\mathfrak{k}_{v})=mes(\mathfrak{p}_{v}\mathfrak{k}_{1,v})mes(\mathfrak{p}_{v}\mathfrak{k}_{2,v}).$$
La suite exacte ci-dessus, jointe au th\'eor\`eme de Lang, entra\^{\i}ne l'exactitude de la suite
$$1\to \Xi^{\Gamma_{v}^{nr}}\to {\mathbb K}_{1,v}({\mathbb F}_{v})\times {\mathbb K}_{2,v}({\mathbb F}_{v})\to {\mathbb K}({\mathbb F}_{v})\to H^1(\Gamma_{v}^{nr};\Xi)\to 1.$$
On v\'erifie facilement l'\'egalit\'e
$$\vert  \Xi^{\Gamma_{v}^{nr}}\vert =\vert H^1(\Gamma_{v}^{nr};\Xi)\vert $$
qui est valable pour tout $\Gamma_{v}^{nr}$-module ab\'elien fini $\Xi$. 
D'o\`u
$$(2) \qquad \vert {\mathbb K}_{1,v}({\mathbb F}_{v})\times {\mathbb K}_{2,v}({\mathbb F}_{v})\vert =\vert  {\mathbb K}({\mathbb F}_{v})\vert .$$
Les \'egalit\'es (1) et (2) entra\^{\i}nent l'\'egalit\'e requise des constantes. $\square$

 \bigskip
 
 \subsection{Coefficients et rev\^etement}
 Dans la suite de cette section, on consid\`ere un triplet $(G,\tilde{G},{\bf a})$ tel que $\tilde{G}=G$. Mais on n'impose pas que ${\bf a}=1$. On impose toutefois que $\omega$ est trivial sur $Z(G;{\mathbb A}_{F})$. On suppose $\tilde{K}_{v}=K_{v}$ pour tout $v\not\in V_{ram}$.   On consid\`ere un sous-tore $Z\subset Z(G)$ et un groupe r\'eductif connexe $G_{\sharp}$. On suppose donn\'e un homomorphisme $q:G_{\sharp}\to G$. Ces trois donn\'ees sont d\'efinies sur $F$. On pose $G_{\flat}=Z\times G_{\sharp}$ et on  prolonge $q$ par l'identit\'e sur $Z$. On obtient ainsi un homomorphisme   not\'e $q_{\flat}:G_{\flat}\to G$. On suppose qu'il s'inscrit dans une suite exacte
 $$1\to \Xi_{\flat}\to G_{\flat}\stackrel{q_{\flat}}{\to} G\to 1$$
 o\`u $\Xi_{\flat}$ est un sous-groupe fini central. On suppose que $\omega$ est trivial sur $q(G_{\flat}({\mathbb A}_{F}))$. On fixe un ensemble fini $V$ de places de $F$ tel que $G_{\flat}$ et $G$ soient non ramifi\'es hors de $V$ et tel que le nombre d'\'el\'ements de $\Xi_{\flat}$ soit premier \`a tout nombre premier divisant une place $v\not\in V$.
 
 {\bf Exemple.} On peut prendre $Z=Z(G)^0$, $G_{\sharp}=G_{SC}$ et $V\supset V_{ram}$. Ces donn\'ees v\'erifient les conditions ci-dessus.

\bigskip 
On note $\Xi$ la projection de $\Xi_{\flat}$ dans $G_{\sharp}$.  Pour toute place $v\in V$, fixons un voisinage ouvert $\Omega_{\sharp,v}$ de $1$ dans $G_{\sharp}(F_{v})$. On suppose que $x$ appartient \`a $\Omega_{\sharp,v}$ si et seulement si la partie semi-simple de $x$ appartient \`a $\Omega_{\sharp,v}$. On suppose que, si $x$ et $x'$ sont deux \'el\'ements de $G_{\sharp}(F_{v})$ qui sont conjugu\'es par un \'el\'ement de $G_{\sharp}(\bar{F}_{v})$, alors $x\in \Omega_{\sharp,v}$ si et seulement si $x'\in \Omega_{\sharp,v}$. On suppose enfin que, si $\xi\in \Xi(F_{v})$ est diff\'erent de $1$, alors $\Omega_{\sharp,v}\cap \xi\Omega_{\sharp,v}=\emptyset$. On pose $\Omega_{v}=q_{\flat}(Z(F_{v})\times \Omega_{\sharp,v})$.  On pose $\Omega_{V}=\prod_{v\in V}\Omega_{v}$, $\Omega_{\sharp,V}=\prod_{v\in V}\Omega_{\sharp,v}$. 
 Fixons un ensemble de repr\'esentants ${\cal U}_{V}$ du quotient fini
 $$ q_{\flat}(G_{\flat}(F_{V}))\backslash G(F_{V}).$$
 Pour $f\in C_{c}^{\infty}( \Omega_{V})$ et $u\in {\cal U}_{V}$, on d\'efinit la fonction $(^{u}f)_{G_{\sharp}}$ sur $ \Omega_{\sharp,V}$ par $(^{u}f)_{G_{\sharp}}( x)=f(u^{-1}q(x)u)$ pour tout $x\in \Omega_{\sharp,V}$.  On pose
 $$(1) \qquad \iota_{G_{\sharp},G}(f)=\vert {\cal U}_{V}\vert ^{-1}\sum_{u\in {\cal U}_{V}}\omega(u)(^{u}f)_{G_{\sharp}} .$$
 Cette application est le produit sur les places $v\in V$ de celles que l'on a d\'efinies et \'etudi\'ees en [III] 3.1. Elle d\'epend du choix de ${\cal U}_{V}$, mais il s'en d\'eduit une application
 $$\iota_{G_{\sharp},G}:I(\Omega_{V},\omega)\to I(\Omega_{\sharp,V})$$
 qui n'en d\'epend pas (les notations sont celles de [III] 3.1, adapt\'ees \`a un ensemble fini de places). Elles sont aussi compatibles, en un sens facile \`a pr\'eciser, \`a un changement de voisinage $\Omega_{\sharp,V}$. Il s'en d\'eduit dualement un homomorphisme
 $$D_{g\acute{e}om}(\Omega_{\sharp,V})\to D_{g\acute{e}om}(\Omega_{V},\omega).$$
 Celui-ci se restreint en un homomorphisme
 $$\iota_{G_{\sharp},G}^*:D_{unip}(G_{\sharp}(F_{V}))\to D_{unip}(G(F_{V}),\omega).$$

 On a fix\'e les mesures en 4.1, la distribution $A^G_{unip}(V,\omega)$ est donc un \'el\'ement de $D_{unip}(G(F_{V}),\omega)$. Elle d\'epend de  la mesure fix\'ee sur $\mathfrak{A}_{G}$, mais $\tau'(G)^{-1}A^G_{unip}(V,\omega)$ n'en d\'epend pas. Elle d\'epend du groupe $K^V=\prod_{v\not\in V}K_{v}$. Si  n\'ecessaire, on fait figurer ce groupe dans la notation. Pour le groupe $G_{\sharp}$ et pour $v\not\in V$, on choisit pour sous-groupe hypersp\'ecial $K_{\sharp,v}$ de $G_{\sharp}(F_{v})$  le groupe tel que $K_{\sharp,v}$ et $K_{v}$ fixent le m\^eme point hypersp\'ecial de l'immeuble commun du groupe $G_{AD}$. Autrement dit $K_{\sharp,v}=q^{-1}(K_{v})$. 
  Soit $S$ un sous-ensemble fini de places de $F$ contenant $V$.   Choisissons un ensemble de repr\'esentants ${\cal U}_{S}^V$ du quotient
 $$ q_{\flat}(G_{\flat}(F_{S}^V)\backslash G(F_{S}^V).$$
 Pour $u=(u_{v})_{v\in S-V}\in {\cal U}_{S}^V$, posons 
 $$^{u}K_{\sharp}^V=(\prod_{v\not\in S}K_{\sharp,v})(\prod_{v\in S-V}u_{v}K_{\sharp,v}u_{v}^{-1}).$$
 Posons
 $$A_{unip;G,\omega,S}^{G_{\sharp}}(V)=\vert {\cal U}_{S}^V\vert ^{-1}\sum_{u\in {\cal U}_{S}^V}\omega(u)A^{G_{\sharp}}_{unip}(V,{^uK}_{\sharp}^V).$$
 Cela ne d\'epend pas du choix de l'ensemble de repr\'esentants.
 
 \ass{Proposition}{  Si l'ensemble $S$ est assez grand, la distribution
  $ \tau'(G)^{-1} A^G_{unip}(V,\omega)$
  est l'image par l'homomorphisme $\iota_{G_{\sharp},G}^*$  de
  $ \tau'(G_{\sharp})^{-1}A_{unip;G,\omega,S}^{G_{\sharp}}(V)$.}

  \bigskip
  
  \subsection{Preuve de la proposition 4.3}
 On doit commencer par quelques pr\'eliminaires. De tout Levi $M$ de $G$ se d\'eduisent des Levi $M_{\flat}=q_{\flat}^{-1}(M)$ et $M_{\sharp}=q^{-1}(M)$ de $G_{\flat}$ et $G_{\sharp}$. Pour $v\in Val(F)$, on pose $K_{\sharp,v}=q^{-1}(K_{v})$ et  $K_{\flat,v}=q_{\flat}^{-1}(K_{v})$. On note $K_{v}^M=M(F_{v})\cap K_{v}$. On rappelle que l'on a fix\'e un Levi minimal $M_{0}$.  On a l'\'egalit\'e
 $$q_{\flat}(M_{0,\flat}(F_{v}))\backslash M_{0}(F_{v})=q_{\flat}(G_{\flat}(F_{v}))\backslash G(F_{v}).$$
 On fixe un ensemble de repr\'esentants ${\cal U}_{v}$ de ce quotient, contenu dans $M_{0}(F_{v})$.     Consid\'erons un ensemble fini $S$ de places de $F$, contenant $V$ et tel que:
 
 (1) pour tout $M\in {\cal L}(M_{0})$, on a les \'egalit\'es
 $$M({\mathbb A}_{F})=M(F)(M(F_{S})\times K^{M,S}),\,\, M_{\flat}({\mathbb A}_{F})=M_{\flat}(F)(M_{\flat}(F_{S})\times K_{\flat}^{M_{\flat},S}).$$
 
 Cette condition est v\'erifi\'ee si $S$ est assez grand. On va prouver l'assertion de la proposition pour un tel $S$.
On pose ${\cal U}_{S}=\prod_{v\in S}{\cal U}_{v}$. 
 
 Soit $\varphi$ une fonction int\'egrable sur $G(F)\mathfrak{A}_{G}\backslash G({\mathbb A}_{F})$. On suppose qu'elle est invariante \`a droite par $Z(G)^0({\mathbb A}_{F})K^S$. Montrons qu'on a l'\'egalit\'e
 $$(2)\qquad \tau'(G)^{-1}\int_{\mathfrak{A}_{G}G(F)\backslash G({\mathbb A}_{F})}\varphi(g)\,dg\,=\tau'(G_{\sharp})^{-1}\vert {\cal U}_{S}\vert ^{-1}\sum_{u\in {\cal U}_{S}}\int_{\mathfrak{A}_{G_{\sharp}}G_{\sharp}(F)\backslash G_{\sharp}({\mathbb A}_{F})}\varphi(q(x)u)\,dx.$$
 
 Preuve. On commence par faire un calcul \`a des constantes multiplicatives pr\`es, \'etant entendu que ces constantes ne d\'ependent pas de $\varphi$.  Notons $\Delta$ la projection dans $G(F_{S})$ de $G(F)\cap (G(F_{S})\times K^S)$ et  d\'efinissons de m\^eme $\Delta_{\flat}$.  En vertu de (1) et de l'invariance de $\varphi$ par $K^S$,   on a l'\'egalit\'e
   $$\int_{\mathfrak{A}_{G}G(F)\backslash G({\mathbb A}_{F})}\varphi(g)\,dg =\int_{\mathfrak{A}_{G}\Delta\backslash G(F_{S})}\varphi(g)\,dg.$$
   On montrera ci-dessous que
   
   (3) $q_{\flat}(\Delta_{\flat})$ est d'indice fini dans $\Delta$.

Il existe donc $c_{1}>0$ tel que l'int\'egrale pr\'ec\'edente soit \'egale \`a
$$c_{1}\int_{\mathfrak{A}_{G}q_{\flat}(\Delta_{\flat})\backslash G(F_{S})}\varphi(g)\,dg.$$ 
Les hypoth\`eses entra\^{\i}nent   que $\mathfrak{A}_{G}=q_{\flat}(\mathfrak{A}_{G_{\flat}})$. Puisque $G(F_{S})=\sqcup_{u\in {\cal U}_{S}}q_{\flat}(G_{\flat}(F_{S}))u$, on peut d\'ecomposer l'int\'egrale pr\'ec\'edente en
$$c_{1}\sum_{u\in {\cal U}_{S}}\int_{q_{\flat}(\mathfrak{A}_{G_{\flat}}\Delta_{\flat})\backslash q_{\flat}(G_{\flat}(F_{S}))}\varphi(gu)\,dg.$$
Puisque $G_{\flat}(F_{S})\to q_{\flat}(G_{\flat}(F_{S}))$ est un honn\^ete rev\^etement, il existe une constante $c_{2}>0$ tel que l'expression pr\'ec\'edente soit \'egale \`a
$$c_{2}\sum_{u\in {\cal U}_{S}}\int_{\mathfrak{A}_{G_{\flat}}\Delta_{\flat}\backslash G_{\flat}(F_{S})}\varphi(q_{\flat}(y)u)\,dy.$$
En utilisant (1) et l'invariance de la fonction \`a int\'egrer par $K_{\flat}^S$, on peut reconstituer cette expression comme
$$c_{2}\sum_{u\in {\cal U}_{S}}\int_{\mathfrak{A}_{G_{\flat}}G_{\flat}(F)\backslash G_{\flat}({\mathbb A}_{F})}\varphi(q_{\flat}(y)u)\,dy.$$
Les int\'egrales se d\'ecomposent en produit d'int\'egrales sur $z\in \mathfrak{A}_{Z}Z(F)\backslash Z({\mathbb A}_{F})$ et sur $x\in \mathfrak{A}_{G_{\sharp}}G_{\sharp}(F)\backslash G_{\sharp}({\mathbb A}_{F})$. D'apr\`es l'hypoth\`ese d'invariance de $\varphi$, la fonction \`a int\'egrer est constante en $z$. Puisque le volume de $\mathfrak{A}_{Z}Z(F)\backslash Z({\mathbb A}_{F})$ est fini, il existe $c_{3}>0$ tel que l'expression pr\'ec\'edente soit \'egale \`a
 $$c_{3}\sum_{u\in {\cal U}_{S}}\int_{\mathfrak{A}_{G_{\sharp}}G_{\sharp}(F)\backslash G_{\sharp}({\mathbb A}_{F})}\varphi(q(x)u)\,dx.$$
 Cela d\'emontre l'\'egalit\'e (2), au calcul pr\`es de la constante $c_{3}$. Pour calculer celle-ci, il suffit d'appliquer la relation obtenue \`a la fonction $\varphi$ constante de valeur $1$. $\square$
 
 Preuve de (3).  On peut d\'efinir
 $$vol(\mathfrak{A}_{G}q_{\flat}(\Delta_{\flat})\backslash G(F_{S}))=\int_{\mathfrak{A}_{G}q_{\flat}(\Delta_{\flat})\backslash G(F_{S})}\,dg$$
 ce terme pouvant valoir $+\infty$. L'assertion (3) \'equivaut \`a dire que ce volume est fini. 
 On reprend le calcul ci-dessus en l'appliquant \`a la fonction $\varphi$ constante de valeur $1$. Il montre que $vol(\mathfrak{A}_{G}q_{\flat}(\Delta_{\flat})\backslash G(F_{S}))$ est le produit d'une constante finie et de $vol(\mathfrak{A}_{G_{\sharp}}G_{\sharp}(F)\backslash G_{\sharp}({\mathbb A}_{F}))$. On sait bien que ce dernier volume est fini. Donc $vol(\mathfrak{A}_{G}q_{\flat}(\Delta_{\flat})\backslash G(F_{S}))$ l'est aussi et (3) est v\'erifi\'ee. 
 $\square$
 
 Soit $f\in C_{c}^{\infty}(G({\mathbb A}_{F}))$. Rappelons qu'en [VI] 2.1, pour un param\`etre $T\in {\cal A}_{M_{0}}$ dans un certain c\^one, on a d\'efini une fonction $g\mapsto k^T_{unip}(f,g)$ sur $Z(G;{\mathbb A}_{F})\backslash G({\mathbb A}_{F})$, puis l'int\'egrale
 $$J^T_{unip}(f,\omega)=\int_{\mathfrak{A}_{G}G(F)\backslash G({\mathbb A}_{F})}k^T_{unip}(f,g)\omega(g)\,dg.$$
 Cette expression est asymptote \`a un polyn\^ome en $T$, dont on note 
   $J_{unip}^G(f,\omega) $ la valeur en un certain point $T_{0}$. Tous ces objets, y compris le point $T_{0}$, d\'ependent de $K=\prod_{v\in Val(F)}K_{v}$. Si n\'ecessaire, on fait figurer ce groupe dans la notation.   En appliquant (2), on obtient
   $$\tau'(G)^{-1}J^T_{unip}(f,\omega)=\tau'(G_{\sharp})^{-1}\vert {\cal U}_{S}\vert ^{-1}\sum_{u\in {\cal U}_{S}}\int_{\mathfrak{A}_{G_{\sharp}}G_{\sharp}(F)\backslash G_{\sharp}({\mathbb A}_{F})}k^T_{unip}(f,q(x)u) \omega(u)\,dx$$
   puisque $\omega$ est trivial sur $q_{\flat}(G_{\flat}({\mathbb A}_{F}))$. Pour tout $u\in {\cal U}_{S}$, on pose $^{u}K_{\sharp}=uK_{\sharp}u^{-1}$ et on d\'efinit une fonction $(^{u}f)_{G_{\sharp}}$ sur $G_{\sharp}({\mathbb A}_{F})$ par $(^{u}f)_{G_{\sharp}}(x)=f(u^{-1}q(x)u)$. En reprenant les d\'efinitions de [LW], on voit qu'on a l'\'egalit\'e
   $$k^T_{unip}(f,q(x)u)=k_{unip}^{G_{\sharp},T-T_{0}+T_{0}(^{u}K_{\sharp})}((^{u}f)_{G_{\sharp}},x,{^{u}K}_{\sharp}).$$
   On en d\'eduit l'\'egalit\'e
   $$(4) \qquad \tau'(G)^{-1}J^G_{unip}(f,\omega)=\tau'(G_{\sharp})^{-1}\vert {\cal U}_{S}\vert ^{-1}\sum_{u\in {\cal U}_{S}}\omega(u)J^{G_{\sharp}}_{unip}((^{u}f)_{G_{\sharp}},{^{u}K}_{\sharp}).$$

Supposons maintenant $f=f_{V}\otimes {\bf 1}_{K^V}$, pour $f_{V}\in C_{c}^{\infty}(\Omega_{V})$.   En posant ${\cal U}_{V}=\prod_{v\in V}{\cal U}_{v}$ et ${\cal U}_{S}^V=\prod_{v\in S-V}{\cal U}_{v}$, on a ${\cal U}_{S}={\cal U}_{V}\times {\cal U}_{S}^V$. 
Pour $u\in {\cal U}_{S}$, que l'on \'ecrit $u=u'u''$, avec $u'\in {\cal U}_{V}$ et $u''\in {\cal U}_{S}^V$, on a $(^{u}f)_{G_{\sharp}}=(^{u'}f_{V})_{G_{\sharp}}\otimes {\bf 1}_{^{u''}K^V_{\sharp}}$, avec des notations naturelles. Pour le membre de gauche de (4), on peut utiliser le d\'eveloppement [VI] 2.2(1) relatif \`a $K^V$. Pour le terme index\'e par $u$ du membre de droite, on utilise le m\^eme d\'eveloppement relatif \`a $^{u}K_{\sharp}^V$. On obtient  l'\'egalit\'e
$$(5) \qquad \sum_{M\in {\cal L}(M_{0})}\vert W^M\vert \vert W^G\vert ^{-1}X_{M}=0,$$
o\`u
$$X_{M}=\tau'(G)^{-1}J_{M}^G(A^M_{unip}(V,\omega),f_{V})-\tau'(G_{\sharp})^{-1}\vert {\cal U}_{S}\vert ^{-1}$$
$$\sum_{u\in {\cal U}_{S}}\omega(u)J_{M_{\sharp}}^{G_{\sharp}}(A^{M_{\sharp}}_{unip}(V,{^{u''}K}_{\sharp}^{M_{\sharp},V}),(^{u'}f_{V})_{G_{\sharp}},{^{u'}K}_{\sharp,V}).$$
  Avec la d\'efinition  de 4.3, on r\'ecrit
$$X_{M}=\tau'(G)^{-1}J_{M}^G(A^M_{unip}(V,\omega),f_{V})-\tau'(G_{\sharp})^{-1}\vert {\cal U}_{V}\vert ^{-1}$$
$$\sum_{u\in {\cal U}_{V}}\omega(u)J_{M_{\sharp}}^{G_{\sharp}}(A^{M_{\sharp}}_{unip;M,\omega,S}(V),(^{u}f_{V})_{G_{\sharp}},{^{u}K}_{\sharp,V}).$$
Soit $M\in {\cal L}(M_{0})$. Les distributions et int\'egrales orbitales pond\'er\'ees intervenant d\'ependent de mesures sur $\mathfrak{A}_{G}$, $\mathfrak{A}_{M}$, $\mathfrak{A}_{G_{\sharp}}$ et $\mathfrak{A}_{M_{\sharp}}$. L'homomorphisme $q$ d\'efinit un isomorphisme de $\mathfrak{A}_{M}/\mathfrak{A}_{G}$ sur $\mathfrak{A}_{M_{\sharp}}/\mathfrak{A}_{G_{\sharp}}$. On peut supposer que cet isomorphisme pr\'eserve les mesures. Pour tout $\boldsymbol{\gamma}\in D_{unip}(M_{\sharp}(F_{V}))$, on a alors l'\'egalit\'e
$$(6)\qquad J_{M}^G(\iota_{M_{\sharp},M}^*(\boldsymbol{\gamma}),\omega,f_{V})=\vert {\cal U}_{V}\vert ^{-1}\sum_{u\in {\cal U}_{V}}\omega(u)J_{M_{\sharp}}^{G_{\sharp}}(\boldsymbol{\gamma}, (^{u}f_{V})_{G_{\sharp}},{^{u}K}_{\sharp,V}).$$
En effet, on a vu l'\'egalit\'e analogue dans le cas local en [III] 3.3. Pour obtenir (5), on peut soit reprendre la preuve de ce cas local, soit utiliser les formules de scindage habituelles. On laisse les d\'etails au lecteur.

Pour $M\not=G$, on peut utiliser la proposition 4.3 par r\'ecurrence: pour 

\noindent $\boldsymbol{\gamma}=\tau'(M_{\sharp})^{-1}A^{M_{\sharp}}_{unip;M,\omega,S}(V)$, on a $\iota_{M_{\sharp},M}^*(\boldsymbol{\gamma})=\tau'(M)^{-1}A^M_{unip}(V,\omega)$ (notons que la condition (1) impos\'ee \`a $S$ implique la m\^eme condition quand on remplace $G$ par $M$). En utilisant (6), on transforme $X_{M}$ en
$$(\tau'(G)^{-1}-\tau'(G_{\sharp})^{-1}\tau'(M)^{-1}\tau'(M_{\sharp}))J_{M}^G(A^M_{unip}(V,\omega),f_{V}).$$
On montrera ci-dessous que
$$(7) \qquad \tau'(G)^{-1}\tau'(M)=\tau'(G_{\sharp})^{-1}\tau'(M_{\sharp}).$$
On obtient alors $X_{M}=0$ pour tout Levi $M\not=G$. L'\'egalit\'e (5) entra\^{\i}ne alors $X_{G}=0$. D'apr\`es les d\'efinitions, cela signifie que
$$\tau'(G)^{-1}I^G(A^G_{unip}(V,\omega),f_{V})=\tau'(G_{\sharp})^{-1}I^{G_{\sharp}}(A^{G_{\sharp}}_{unip;G,\omega,S}(V), \iota_{G_{\sharp},G}(f_{V})),$$
ou encore
$$\tau'(G)^{-1}I^G(A^G_{unip}(V,\omega),f_{V})=\tau'(G_{\sharp})^{-1}I^G(\iota_{G_{\sharp},G}^*(A^{G_{\sharp}}_{unip;G,\omega,S}(V)),f_{V}).$$
Cela \'etant vrai pour tout $f_{V}\in C_{c}^{\infty}(\Omega_{V})$, on en d\'eduit 
$$\tau'(G)^{-1}A^G_{unip}(V,\omega)=\tau'(G_{\sharp})^{-1}\iota_{G_{\sharp},G}^*(A^{G_{\sharp}}_{unip;G,\omega,S}(V)),$$
ce qui prouve la proposition 4.3.

Il reste \`a prouver (7). Puisque $G_{\flat}$, resp. $M_{\flat}$, est le produit de $Z$ et de $G_{\sharp}$, resp. $M_{\sharp}$, on a 
$$  \tau'(G_{\sharp})^{-1}\tau'(M_{\sharp})=\tau'(G_{\flat})^{-1}\tau'(M_{\flat}).$$
On peut aussi bien d\'emontrer l'\'egalit\'e
$$(8)\qquad \tau'(G)^{-1}\tau'(M)=\tau'(G_{\flat})^{-1}\tau'(M_{\flat}).$$
Rappelons que l'on a identifi\'e $\mathfrak{A}_{G}$ \`a $Hom(X^*(G)^{\Gamma_{F}},{\mathbb R})$ et que l'on a d\'efini le r\'eseau $\mathfrak{A}_{G,{\mathbb Z}}=Hom(X^*(G)^{\Gamma_{F}},{\mathbb Z})$. On a une injection $X^*(G)^{\Gamma_{F}}\to X^*(M)^{\Gamma_{F}}$. Son conoyau est sans torsion. En effet, si on introduit un tore maximal $T\subset M$, ce conoyau est l'image de l'homomorphisme naturel $X^*(M)^{\Gamma_{F}}\to X^*(T_{sc})^{\Gamma_{F}}$, o\`u, comme toujours, $T_{sc}$ est l'image r\'eciproque de $T$ dans $G_{SC}$. De l'injection pr\'ec\'edente se d\'eduisent des homomorphismes $\mathfrak{A}_{M}\to \mathfrak{A}_{G}$ et $\mathfrak{A}_{M,{\mathbb Z}}\to \mathfrak{A}_{G,{\mathbb Z}}$.  Le premier est trivialement surjectif. Le second est lui aussi surjectif, parce que le conoyau de l'injection $X^*(G)^{\Gamma_{F}}\to X^*(M)^{\Gamma_{F}}$ est sans torsion. En notant $\mathfrak{A}_{M}^G$ et $\mathfrak{A}_{M,{\mathbb Z}}^G$ les noyaux de ces homomorphismes, on obtient une suite exacte
$$0\to \mathfrak{A}_{M}^G/\mathfrak{A}_{M,{\mathbb Z}}^G\to \mathfrak{A}_{M}/\mathfrak{A}_{M,{\mathbb Z}}\to \mathfrak{A}_{G}/\mathfrak{A}_{G,{\mathbb Z}}\to 0$$
Donc
$$covol(\mathfrak{A}_{M,{\mathbb Z}})^{-1}covol(\mathfrak{A}_{G,{\mathbb Z}})=covol(\mathfrak{A}_{M,{\mathbb Z}}^G)^{-1}.$$
En se rappelant la d\'efinition de 4.1, on obtient
$$  \tau'(G)^{-1}\tau'(M)=\tau(G)^{-1}\tau(M)covol(\mathfrak{A}_{M,{\mathbb Z}}^G)^{-1}.$$
On a d\'emontr\'e en [VI] 6.1 l'\'egalit\'e $ker^1(F,Z(\hat{G}))=ker^1(F,Z(\hat{M}))$. L'\'egalit\'e pr\'ec\'edente se r\'ecrit donc
$$(9)\qquad  \tau'(G)^{-1}\tau'(M)=\vert \pi_{0}(Z(\hat{G})^{\Gamma_{F}})\vert ^{-1}\vert \pi_{0}(Z(\hat{M})^{\Gamma_{F}})\vert covol(\mathfrak{A}_{M,{\mathbb Z}}^G)^{-1}.$$
On a bien s\^ur une relation analogue pour $G_{\flat}$ et $M_{\flat}$.  
On a un diagramme commutatif
$$(10) \qquad \begin{array}{ccc}X^*(G)^{\Gamma_{F}}&\to&X^*(M)^{\Gamma_{F}}\\ \downarrow&&\downarrow\\ X^*(G_{\flat})^{\Gamma_{F}}&\to&X^*(M_{\flat})^{\Gamma_{F}}\\ \end{array}$$
dont les fl\`eches sont injectives. Les fl\`eches verticales sont de conoyaux finis. On en d\'eduit un diagramme d'isomorphismes
$$\begin{array}{ccccccccc}1&\to&\mathfrak{A}_{M}^G&\to&\mathfrak{A}_{M}&\to&\mathfrak{A}_{G}&\to&1\\ &&\parallel &&\parallel&&\parallel&&\\1&\to&\mathfrak{A}_{M_{\flat}}^{G_{\flat}}&\to&\mathfrak{A}_{M_{\flat}}&\to&\mathfrak{A}_{G_{\flat}}&\to&1\\ \end{array}$$
et un diagramme commutatif
$$\begin{array}{ccccccccc}1&\to&\mathfrak{A}_{M,{\mathbb Z}}^G&\to&\mathfrak{A}_{M,{\mathbb Z}}&\to&\mathfrak{A}_{G,{\mathbb Z}}&\to&1\\ &&\uparrow &&\uparrow&&\uparrow&&\\1&\to&\mathfrak{A}_{M_{\flat},{\mathbb Z}}^{G_{\flat}}&\to&\mathfrak{A}_{M_{\flat},{\mathbb Z}}&\to&\mathfrak{A}_{G_{\flat},{\mathbb Z}}&\to&1\\ &&\uparrow&&\uparrow&&\uparrow&&\\ &&1&&1&&1&&\\ \end{array}$$
Les suites ci-dessus sont exactes. On les compl\`ete en notant $B_{M}^G$, $B_{M}$ et $B_{G}$ les conoyaux des suites verticales. Ils sont finis et  on a la suite exacte
$$1\to B_{M}^G\to B_{M}\to B_{G}\to 1.$$
   On a normalis\'e les mesures de sorte que l'isomorphisme  $\mathfrak{ A}_{M_{\flat}}^{G_{\flat}}\to \mathfrak{ A}_{M}^G$   les conserve. On en d\'eduit
$$ covol(\mathfrak{A}_{M_{\flat},{\mathbb Z}}^{G_{\flat}})=covol(\mathfrak{A}_{M,{\mathbb Z}}^G)\vert B_{M}^G\vert ,$$
d'o\`u aussi
$$(11) \qquad covol(\mathfrak{A}_{M_{\flat},{\mathbb Z}}^{G_{\flat}})=covol(\mathfrak{A}_{M,{\mathbb Z}}^G)\vert  B_{G}\vert ^{-1}\vert B_{M}\vert .$$
 On sait que $X^*(G)\simeq X_{*}(Z(\hat{G})^0)$. Dualement au diagramme (10), on a un diagramme commutatif
$$\begin{array}{ccc}Z(\hat{G})^{\Gamma_{F},0}&\to &Z(\hat{M})^{\Gamma_{F},0}\\ \downarrow&&\downarrow\\ Z(\hat{G}_{\flat})^{\Gamma_{F},0}&\to&Z(\hat{M}_{\flat})^{\Gamma_{F},0}\\ \end{array}$$
Les fl\`eches horizontales sont injectives. Les fl\`eches verticales sont surjectives. On note leurs noyaux $\hat{B}_{G}$ et $\hat{B}_{M}$. Remarquons que, puisque $\hat{G}$ et $\hat{G}_{\flat}$ ont m\^eme groupe adjoint,  l'image r\'eciproque par la deuxi\`eme fl\`eche verticale de $Z(\hat{M}_{\flat})^{\Gamma_{F},0}\cap Z(\hat{G}_{\flat})^{\Gamma_{F}}$ n'est autre que $Z(\hat{M})^{\Gamma_{F},0}\cap Z(\hat{G})^{\Gamma_{F}}$. Remarquons aussi que le quotient
$$Z(\hat{G})^{\Gamma_{F},0}\backslash (Z(\hat{M})^{\Gamma_{F},0}\cap Z(\hat{G})^{\Gamma_{F}})$$
n'est autre que le groupe des composantes connexes $\pi_{0}(Z(\hat{M})^{\Gamma_{F},0}\cap Z(\hat{G})^{\Gamma_{F}})$.
On obtient alors un diagramme commutatif
$$\begin{array}{ccccccccc}&&1&&1&&1&&\\ &&\downarrow&&\downarrow&&\downarrow&&\\ 1&\to &\hat{B}_{G}&\to&\hat{B}_{M}&\to &\hat{B}_{M}^G&\to&1\\ &&\downarrow&&\downarrow&&\downarrow&&\\
1&\to&Z(\hat{G})^{\Gamma_{F},0}&\to&Z(\hat{M})^{\Gamma_{F},0}\cap Z(\hat{G})^{\Gamma_{F}}&\to&\pi_{0}(Z(\hat{M})^{\Gamma_{F},0}\cap Z(\hat{G})^{\Gamma_{F}})&\to&1\\ &&\downarrow&&\downarrow&&\downarrow&&\\ 1&\to& Z(\hat{G}_{\flat})^{\Gamma_{F},0}&\to&Z(\hat{M}_{\flat})^{\Gamma_{F},0}\cap Z(\hat{G}_{\flat})^{\Gamma_{F}}&\to&\pi_{0}(Z(\hat{M}_{\flat})^{\Gamma_{F},0}\cap Z(\hat{G}_{\flat})^{\Gamma_{F}})&\to&1\\ &&\downarrow&&\downarrow&&\downarrow&&\\  &&1&&1&&1&&\\ \end{array}$$
o\`u $\hat{B}_{M}^G$ est d\'efini comme le noyau de la derni\`ere suite verticale. Les lignes verticales sont exactes. Les deux derni\`eres lignes horizontales aussi. Donc la premi\`ere ligne horizontale aussi. 

Revenons \`a la d\'efinition de $B_{G}$, qui est le conoyau de l'homomrophisme $\mathfrak{A}_{G_{\flat},{\mathbb Z}}\to \mathfrak{A}_{G,{\mathbb Z}}$.
Par dualit\'e, $\vert B_{G}\vert $ est aussi le nombre d'\'el\'ements du conoyau de l'homomorphisme $X^*(G)^{\Gamma_{F}}\to X^*(G_{\flat})^{\Gamma_{F}}$.  En vertu de l'isomorphisme $X^*(G)\simeq X_{*}(Z(\hat{G})^0)$ et de l'isomorphisme analogue pour $G_{\flat}$, on voit que $B_{G}$ a m\^eme nombre d'\'el\'ements que $\hat{B}_{G}$.  De m\^eme, $B_{M}$ a m\^eme nombre d'\'el\'ements que $\hat{B}_{M}$. D'apr\`es le diagramme ci-dessus, on obtient
$$\vert  B_{G}\vert^{-1} \vert B_{M}\vert =\vert \hat{B}_{M}^G\vert =\vert \pi_{0}(Z(\hat{M}_{\flat})^{\Gamma_{F},0}\cap Z(\hat{G}_{\flat})^{\Gamma_{F}})\vert^{-1} \vert \pi_{0}(Z(\hat{M})^{\Gamma_{F},0}\cap Z(\hat{G})^{\Gamma_{F}})\vert .$$
On a la suite exacte
$$1\to \pi_{0}(Z(\hat{M})^{\Gamma_{F},0}\cap Z(\hat{G})^{\Gamma_{F}})\to \pi_{0}(Z(\hat{M})^{\Gamma_{F}})\to \pi_{0}(Z(\hat{G})^{\Gamma_{F}})\to 1.$$
D'o\`u l'\'egalit\'e
$$\vert  B_{G}\vert \vert B_{M}\vert ^{-1}=\vert \pi_{0}(Z(\hat{M}_{\flat})^{\Gamma_{F}})\vert ^{-1}\vert \pi_{0}(Z(\hat{G}_{\flat})^{\Gamma_{F}})\vert \vert \pi_{0}(Z(\hat{M})^{\Gamma_{F}})\vert \vert \pi_{0}(Z(\hat{G})^{\Gamma_{F}})\vert^{-1} .$$
En ins\'erant cette \'egalit\'e dans (11) et en utilisant (9) ainsi que la relation analogue pour $G_{\flat}$, on obtient (8). Cela ach\`eve la d\'emonstration. $\square$

\bigskip

\subsection{Donn\'ees endoscopiques et rev\^etement}
On a \'etudi\'e cette question dans le cas local en [III] 3.6 et 3.7 et en [V] 3.3. On se contente ici de reprendre bri\`evement les constructions dans notre cadre global. 

On suppose $G$ quasi-d\'eploy\'e et ${\bf a}=1$. Comme on l'a vu en [III] 3.5, l'homomorphisme $\iota_{G_{\sharp},G}:C_{c}^{\infty}(\Omega_{V})\to C_{c}^{\infty}(\Omega_{\sharp,V})$  de 4.3 se quotiente en un homomorphisme $\iota_{G_{\sharp},G}:SI(\Omega_{V})\to SI(\Omega_{\sharp,V})$. Il est plus commode de noter cet homomorphisme
$\iota_{G_{\sharp},G}:SI(G(F_{V}))\to SI(G_{\sharp}(F_{V}))$, \'etant entendu qu'il n'est d\'efini que sur les fonctions \`a support assez voisin de l'origine. 
Dualement, on a un homomorphisme $D_{g\acute{e}om}^{st}( G_{\sharp}(F_{V}))\to D_{g\acute{e}om}^{st}(G(F_{V}))$, d\'efini sur les distributions \`a support voisin de l'origine. Il se restreint en un homomorphisme
$$\iota^*_{G_{\sharp},G}:D_{unip}^{st}(G_{\sharp}(F_{V}))\to D_{unip}^{st}(G(F_{V})).$$
C'est aussi la restriction \`a $D_{unip}^{st}(G_{\sharp}(F_{V}))$ de l'homomorphisme $\iota^*_{G_{\sharp},G}$ d\'efini en 4.3.

Dualement \`a la suite exacte
$$1\to \Xi_{\flat}\to Z\times G_{\sharp}\to G\to 1,$$
on a une suite exacte
$$1\to \hat{\Xi}_{\flat}\to \hat{G}\to \hat{Z}\times \hat{G}_{\sharp}\to 1,$$
o\`u $\hat{\Xi}_{\flat}$ est un sous-groupe fini central de $\hat{G}$.
Soit ${\bf G}'=(G',{\cal G}',s)$ une donn\'ee endoscopique de $G$. L'\'el\'ement $s\in \hat{G}$ s'envoie sur un \'el\'ement $(z,s_{\sharp})\in \hat{Z}\times \hat{G}_{\sharp}$. En notant $\hat{G}'_{\sharp}$ la composante neutre de $Z_{G_{\sharp}}(s_{\sharp})$, on a la suite exacte
$$(1) \qquad 1\to \hat{\Xi}_{\flat}\to\hat{G}'\to \hat{Z}\times \hat{G}'_{\sharp}\to 1.$$
Le groupe ${\cal G}'/\hat{\Xi}_{\flat}$ contient $\hat{Z}$ donc est de la forme $\hat{Z}\times {\cal G}'_{\sharp}$, o\`u ${\cal G}'_{\sharp}$ est un sous-groupe de $^LG_{\sharp}$. Ce groupe d\'efinit une action galoisienne sur $\hat{G}'_{\sharp}$. On introduit un groupe quasi-d\'eploy\'e $G'_{\sharp}$ sur $F$ dont  $\hat{G}'_{\sharp}$ soit le groupe dual. Alors ${\bf G}'_{\sharp}=(G'_{\sharp},{\cal G}'_{\sharp},s_{\sharp})$ est une donn\'ee endoscopique de $G_{\sharp}$. En [I] 2.7, on a associ\'e \`a une telle donn\'ee un caract\`ere de $G_{\sharp,AD}({\mathbb A}_{F})$ not\'e $\omega_{\sharp}$.   

{\bf Remarque.} Dans cette r\'ef\'erence, le corps de base \'etait local mais la construction vaut aussi bien sur notre corps de nombres. D'autre part, l'indice $\sharp$ n'avait pas de rapport avec le pr\'esent indice mais on peut aussi bien conserver ici cette notation.

\bigskip
Rappelons que $G_{\sharp,AD}=G_{AD}$. Pour la donn\'ee que l'on vient de construire, on a

(2) la restriction de $\omega_{\sharp}$ \`a l'image de $G({\mathbb A}_{F})$ dans $G_{AD}({\mathbb A}_{F})$ est triviale.

Inversement, soit  ${\bf G}'_{\sharp}=(G'_{\sharp},{\cal G}'_{\sharp},s_{\sharp})$  une donn\'ee endoscopique de $G_{\sharp}$. On fixe une image r\'eciproque $s\in \hat{G}$ de $(1,s_{\sharp})\in \hat{Z}\times \hat{G}_{\sharp}$. On note ${\cal G}'$ l'image r\'eciproque de $\hat{Z}\times {\cal G}'_{\sharp}$ dans $^LG$. Ce groupe agit sur $\hat{G}_{s}$ et munit ce groupe d'une action galoisienne. On introduit un groupe $G'$ quasi-d\'eploy\'e sur $F$ dont $\hat{G}_{s}$ est le groupe dual. Le triplet $(G',{\cal G}',s)$ est une donn\'ee endoscopique pour $G$ muni d'un certain cocycle ${\bf a}$. C'est une donn\'ee endoscopique pour $G$, c'est-\`a-dire ce cocycle est trivial, si et seulement si la condition (2) est v\'erifi\'ee.

Ces constructions d\'efinissent des bijections inverses l'une de l'autre entre les classes d'\'equivalence de donn\'ees endoscopiques pour $G$ et les classes d'\'equivalence de donn\'ees endoscopiques pour $G_{\sharp}$ v\'erifiant (2). Ces bijections pr\'eservent l'ellipticit\'e et la non-ramification hors de $V$. On a fix\'e des ensembles de repr\'esentants ${\cal E}(G,V)$ et ${\cal E}(G_{\sharp},V)$ des classes d'\'equivalence de donn\'ees endoscopiques pour $G$ et $G_{\sharp}$ qui sont elliptiques et non ramifi\'ees hors de $V$. On note ${\cal E}_{G}(G_{\sharp},V)$ le sous-ensemble des \'el\'ements de ${\cal E}(G_{\sharp},V)$ qui v\'erifient la condition (2). Les ensembles ${\cal E}(G,V)$ et ${\cal E}_{G}(G_{\sharp},V)$ sont en bijection. 

Consid\'erons une donn\'ee ${\bf G}'=(G',{\cal G}',s)$ et la donn\'ee ${\bf G}'_{\sharp}=(G'_{\sharp},{\cal G}'_{\sharp},s_{\sharp})$ construite ci-dessus. Dualement \`a la suite (1), on a la suite exacte
$$1\to \Xi_{\flat}\to Z\times G'_{\sharp}\stackrel{q_{\sharp}}{\to} G'\to 1.$$
Les groupes $G'$ et $G'_{\sharp}$ sont donc dans la m\^eme situation que les groupes $G$ et $G_{\sharp}$. Fixons des donn\'ees auxiliaires $G'_{1}$, $C_{1}$, $\hat{\xi}_{1}$ pour ${\bf G}'$. Notons $G'_{\sharp,1}$ le produit fibr\'e de $G'_{1}$ et $G'_{\sharp}$ au-dessus de $G'$. Le groupe dual $\hat{G}'_{\sharp,1}$ s'identifie \`a $\hat{G}'_{1}/\hat{\xi}_{1}(\hat{Z}_{\flat})$, o\`u $\hat{Z}_{\flat}$ est l'image r\'eciproque de $\hat{Z}$ dans $\hat{G}$. Puisque ${\cal G}'_{\sharp}$ est aussi isomorphe \`a ${\cal G}'/\hat{Z}_{\flat}$, le plongement $\hat{\xi}_{1}$ se quotiente en un plongement $\hat{\xi}_{\sharp,1}:{\cal G}'_{\sharp}\to {^LG}'_{\sharp,1}$. Les donn\'ees $G'_{\sharp,1}$, $C_{1}$, $\hat{\xi}_{\sharp,1}$ sont des donn\'ees auxiliaires pour ${\bf G}'_{\sharp}$. Supposons les donn\'ees endoscopiques non ramifi\'ees hors de $V$, ainsi que les donn\'ees auxiliaires pour ${\bf G}'$. Alors les donn\'ees auxiliaires pour ${\bf G}'_{\sharp}$ sont elles-aussi non ramifi\'ees hors de $V$. On a vu en [VI] 3.6 que le choix des groupes $K_{v}$ pour $v\not\in V$ permettait de d\'efinir un  facteur de transfert canonique $\Delta_{1,V}$ sur $G'_{1}(F_{V})\times G(F_{V})$ (notons que, dans le cas non tordu, l'hypoth\`ese {\bf Hyp} de cette r\'ef\'erence est toujours v\'erifi\'ee). On a relev\'e les $K_{v}$ en des sous-groupes $K_{\sharp,v}$. Ils d\'eterminent de m\^eme un facteur de transfert canonique $\Delta_{\sharp, 1,V}$ sur $G'_{\sharp,1}(F_{V})\times G_{\sharp}(F_{V})$. On v\'erifie que si $(\delta_{\sharp,1},\gamma_{\sharp})\in  G'_{\sharp,1}(F_{V})\times G_{\sharp}(F_{V})$ est un couple d'\'el\'ements semi-simples $G_{\sharp}$-fortement r\'eguliers qui se correspondent, on a l'\'egalit\'e
$$\Delta_{\sharp,1,V}(\delta_{\sharp,1},\gamma_{\sharp})=\Delta_{1,V}(\delta_{1},\gamma),$$
o\`u $\delta_{1}$ et $\gamma$ sont les projections naturelles de $\delta_{\sharp,1}$ et $\gamma_{\sharp}$. 

{\bf Remarque.} Il se peut que les \'el\'ements $\delta_{1}$ et $\gamma$ se correspondent alors que $\delta_{\sharp,1}$ et $\gamma_{\sharp}$ ne se correspondent pas. L'\'egalit\'e ci-dessus devient fausse, le membre de gauche \'etant nul alors que celui de droite ne l'est pas.
\bigskip

On peut adapter les constructions de 4.3 et d\'efinir un homomorphisme
$$\iota_{G'_{\sharp,1},G'_{1}}:SI_{\lambda_{1}}(G'_{1}(F_{V}))\to SI_{\lambda_{1}}(G'_{\sharp,1}(F_{V}))$$
bien d\'efini sur les fonctions dont le support dans $G'(F_{V})$ est voisin de l'origine. Le diagramme suivant est commutatif
$$(3) \qquad \begin{array}{ccc}I(G(F_{V}))&\stackrel{transfert}{\to}&SI_{\lambda_{1}}(G'_{1}(F_{V}))\\ \iota_{G_{\sharp},G}\downarrow\,\,\,&&\,\,\,\downarrow \iota_{G'_{\sharp,1},G'_{1}}\\ I(G_{\sharp}(F_{V}))&\stackrel{transfert}{\to}&SI_{\lambda_{1}}(G'_{\sharp,1}(F_{V}))\\ \end{array}$$
Supposons que nos donn\'ees endoscopiques soient elliptiques. Les espaces $\mathfrak{A}_{G}$ et $\mathfrak{A}_{G'}$ sont isomorphes et, conform\'ement aux conventions de [VI], on suppose que l'isomorphisme pr\'eserve les mesures. De m\^eme pour les espaces $\mathfrak{A}_{G_{\sharp}}$ et $\mathfrak{A}_{G'_{\sharp}}$. On fixe une mesure sur $\mathfrak{A}_{C_{1}}$. En utilisant les suites exactes habituelles, des mesures d\'ej\`a fix\'ees se d\'eduisent des mesures sur $\mathfrak{A}_{G'_{1}}$ et $\mathfrak{A}_{G'_{\sharp,1}}$. Avec ces normalisations,
montrons que l'on a l'\'egalit\'e
$$(4)\qquad i(G,G')\tau'(G)^{-1}\tau'(G'_{1})=i(G_{\sharp},G'_{\sharp})\tau'(G_{\sharp})^{-1}\tau'(G'_{\sharp,1}).$$
Preuve.  Rappelons (cf. [VI] 5.1) que, dans notre situation quasi-d\'eploy\'ee et sans torsion, on a simplement 
$$i(G,G')=\vert Out({\bf G}')\vert ^{-1} \tau(G)\tau(G')^{-1}.$$
D'autre part, le lemme 4.2 nous dit que $\tau'(G'_{1})=\tau'(G')\tau'(C_{1})$. 
Le membre de gauche de (4) est donc \'egal \`a
$$\vert Out({\bf G}')\vert ^{-1}covol(\mathfrak{A}_{G,{\mathbb Z}})covol(\mathfrak{A}_{G',{\mathbb Z}})^{-1} \tau'(C_{1}).$$
 Parce que ${\bf G}'$ est elliptique, on a les isomorphismes
$$X^*(G)^{\Gamma_{F}}\simeq X_{*}(Z(\hat{G})^{\Gamma_{F},0})=X_{*}(Z(\hat{G}')^{\Gamma_{F},0})\simeq X^*(G')^{\Gamma_{F}}.$$
Il en r\'esulte que $covol(\mathfrak{A}_{G,{\mathbb Z}})=covol(\mathfrak{A}_{G',{\mathbb Z}})$. Donc le membre de gauche de (4) est \'egal \`a $\vert Out({\bf G}')\vert ^{-1}\tau'(C_{1})$. On a une formule analogue pour le membre de droite. 
On v\'erifie imm\'ediatement que les groupes $Out({\bf G}')$ et $Out({\bf G}'_{\sharp})$ sont isomorphes. L'\'egalit\'e (4) s'ensuit. $\square$

\bigskip

\subsection{ Coefficients stables et rev\^etement}
On suppose encore $G$ quasi-d\'eploy\'e et ${\bf a}=1$.

 \ass{Proposition}{ La distribution
 $ \tau'(G)^{-1}SA^G_{unip}(V)$
 est l'image par l'homomorphisme $\iota^*_{G_{\sharp},G}$ de
 $\tau'(G_{\sharp})^{-1} SA^{G_{\sharp}}_{unip}(V)$.}

 Preuve.   Soit $f\in C_{c}^{\infty}(G(F_{V}))$. On a par d\'efinition
 $$(1) \qquad \tau'(G)^{-1}I^{G}(SA^G_{unip}(V),f)=\tau'(G)^{-1}I^{G}(A^G_{unip}(V),f)$$
 $$-\sum_{{\bf G}'\in {\cal E}(G,V),G'\not=G}\tau'(G)^{-1}i(G,G')S^{{\bf G}'}( SA^{{\bf G}'}_{unip}(V),f^{{\bf G}'}).$$
Fixons ${\bf G}'=(G',{\cal G}',s)\in {\cal E}(G,V)$ avec $G'\not=G$.  Introduisons des donn\'ees auxiliaires pour cette donn\'ee, non ramifi\'ees hors de $V$. Introduisons aussi la donn\'ee ${\bf G}'_{\sharp}$ comme dans le paragraphe pr\'ec\'edent. On utilise les notations de ce paragraphe. On identifie $f^{{\bf G}'}$ \`a un \'el\'ement $f_{1}\in SI_{\lambda_{1}}(G'_{1}(F_{V}))$. On a
$$S^{{\bf G}'}( SA^{{\bf G}'}_{unip}(V),f^{{\bf G}'})=S^{G'_{1}}_{\lambda_{1}}(SA^{G'_{1}}_{unip,\lambda_{1}}(V),f_{1}).$$
On peut appliquer la proposition par r\'ecurrence \`a $G'$ puisque $G'\not=G$. On travaille ici avec des distributions qui se transforment selon le caract\`ere $\lambda_{1}$ de $C_{1}(F_{V})$ mais la derni\`ere \'egalit\'e du paragraphe 1.11 permet d'adapter la proposition \`a ce cas. On a donc 
$$SA^{G'_{1}}_{unip,\lambda_{1}}(V)=\tau'(G'_{1})\tau'(G'_{1,\sharp})^{-1}\iota_{G'_{\sharp,1},G'_{1}}^*(SA^{G'_{\sharp,1}}_{unip,\lambda_{1}}(V)).$$
D'o\`u
$$S^{G'_{1}}_{\lambda_{1}}(SA^{G'_{1}}_{unip,\lambda_{1}}(V),f_{1})=\tau'(G'_{1})\tau'(G'_{1,\sharp})^{-1}S^{G'_{\sharp,1}}_{\lambda_{1}}(SA^{G'_{\sharp,1}}_{unip,\lambda_{1}}(V),f_{\sharp,1}),$$
o\`u  $f_{\sharp,1}=\iota_{G'_{\sharp,1},G'_{1}}(f_{1})$.
En utilisant ces formules et 4.5(4), on transforme la somme en ${\bf G}'$ de l'expression (1) en
$$(2) \qquad \sum_{{\bf G}'_{\sharp}\in {\cal E}_{G}(G_{\sharp},V)}\tau'(G_{\sharp})^{-1}i(G_{\sharp},G'_{\sharp})S^{G'_{\sharp,1}}_{\lambda_{1}}(SA^{G'_{\sharp,1}}_{unip,\lambda_{1}}(V),f_{\sharp,1}).$$

Posons $f_{\sharp}=\iota_{G_{\sharp},G}(f)$. Remarquons que, d'apr\`es 4.5(3), les fonctions $f_{\sharp,1}$ intervenant peuvent aussi se d\'efinir comme le transfert de $f_{\sharp}$ \`a $G'_{\sharp,1}(F_{V})$. Le terme que l'on somme ne fait alors plus r\'ef\'erence \`a $G$. En particulier, on peut le d\'efinir pour toute donn\'ee ${\bf G}'_{\sharp}\in {\cal E}(G_{\sharp},V)$ et pas seulement pour celles qui proviennent d'une donn\'ee endoscopique de $G$. D'autre part, les facteurs de transfert utilis\'es pour d\'efinir ces fonctions $f_{\sharp,1}$ d\'ependent des compacts $K_{v}$ pour $v\not\in V$. Notons-les plus pr\'ecis\'ement $f_{\sharp,1}[K^V]$. On sait par contre que les distributions $SA^{G'_{\sharp,1}}_{unip,\lambda_{1}}(V)$ ne d\'ependent pas des compacts. Soit $S$ un ensemble fini de places de $F$ contenant $V$. Soit ${\cal U}_{S}^V$ un ensemble de repr\'esentants du quotient $q_{\flat}(G_{\flat}(F_{S}^V))\backslash G(F_{S}^V)$. Pour $u\in {\cal U}_{S}^V$, on d\'efinit pour $v\in S-V$ le groupe $^{u}K_{v}=uK_{v}u^{-1}$, qui se rel\`eve en le compact $^{u}K_{\sharp,v}=uK_{\sharp,v}u^{-1}$ de $G_{\sharp}(F_{v})$. On pose $^{u}K^V=(\prod_{v\in S-V}{^{u}K}_{v})(\prod_{v\not\in S}K_{v})$. Soit ${\bf G}'_{\sharp}\in {\cal E}(G_{\sharp},V)$. Montrons que

(3) si $S$ est assez grand, on a l'\'egalit\'e
$$\vert {\cal U}_{S}^V\vert ^{-1}\sum_{u\in {\cal U}_{S}^V}S^{G'_{\sharp,1}}_{\lambda_{1}}(SA^{G'_{\sharp,1}}_{unip,\lambda_{1}}(V),f_{\sharp,1}[^{u}K^V])=$$
$$\left\lbrace\begin{array}{cc}S^{G'_{\sharp,1}}_{\lambda_{1}}(SA^{G'_{\sharp,1}}_{unip,\lambda_{1}}(V),f_{\sharp,1}[K^{V}]),&\text{ si }{\bf G}'_{\sharp}\in {\cal E}_{G}(G_{\sharp},V),\\ 0,&\text{ sinon.}\\ \end{array}\right.$$

 En 4.5, on a rappel\'e l'existence d'un caract\`ere $\omega_{\sharp}$ de $G_{AD}({\mathbb A}_{F})$ associ\'e \`a ${\bf G}'_{\sharp}$. La fonction $f_{\sharp}$ est par construction invariante par l'action de $G(F_{V})$. Il r\'esulte de cela et de la d\'efinition du caract\`ere $\omega_{\sharp}$ que tous les transferts $f_{\sharp,1}[^{u}K^V]$ sont nuls sauf si $\omega_{\sharp}$ est trivial sur $G(F_{V})$. Si cette condition n'est pas v\'erifi\'ee, la relation (2) est donc \'evidente. Supposons que $\omega_{\sharp}$ est trivial sur $G(F_{V})$. 
Pour $v\in S-V$, notons $^{u}\Delta_{\sharp,1,v}$ le facteur de transfert local associ\'e \`a $^{u}K_{\sharp,v}$. Ce n'est autre que $(\delta_{\sharp,1},\gamma_{\sharp})\mapsto \Delta_{\sharp,v}(\delta_{\sharp,1},u^{-1}\gamma_{\sharp}u)$, o\`u $\Delta_{\sharp,v}$ est associ\'e \`a $K_{\sharp,v}$. On a simplement $\Delta_{\sharp,v}(\delta_{\sharp,1},u^{-1}\gamma_{\sharp}u)=\omega_{\sharp,v}(u)\Delta_{\sharp,v}(\delta_{\sharp,1},\gamma_{\sharp})$. Donc $^{u}\Delta_{\sharp,1,v}=\omega_{\sharp,v}(u)\Delta_{\sharp,1,v}$. Ces facteurs locaux en $v\in S-V$ sont les seules donn\'ees qui changent quand on change de compacts. En se rappelant la d\'efinition des facteurs de transfert canoniques, on voit alors que $f_{\sharp,1}[^{u}K^V]=\omega_{\sharp}(u)^{-1}f_{\sharp,1}[K^V]$. La somme de (3) vaut donc
$$S^{G'_{\sharp,1}}_{\lambda_{1}}(SA^{G'_{\sharp,1}}_{unip,\lambda_{1}}(V),f_{\sharp,1}[K^{V}])
\vert {\cal U}_{S}^V\vert ^{-1}\sum_{u\in {\cal U}_{S}^V}\omega_{\sharp}(u)^{-1}.$$
C'est-\`a-dire qu'elle vaut $S^{G'_{\sharp,1}}_{\lambda_{1}}(SA^{G'_{\sharp,1}}_{unip,\lambda_{1}}(V),f_{\sharp,1}[K^{V}])$ si $\omega_{\sharp}$ est trivial sur $G(F_{S}^V)$ et $0$ sinon.   Le caract\`ere $\omega_{\sharp}$ est   automorphe et, parce que ${\bf G}'_{\sharp}$ est non ramifi\'e hors de $V$,  il est trivial sur $K^V$. Choisissons $S$ tel que $G({\mathbb A}_{F})=G(F)(G(F_{S})\times K^S)$. Parce que l'on a suppos\'e  que $\omega_{\sharp}$ \'etait trivial sur $G(F_{V})$, la condition que ce caract\`ere soit  trivial sur $G(F_{S}^V)$ \'equivaut alors \`a ce qu'il soit trivial sur tout $G({\mathbb A}_{F})$. C'est la condition 4.5(2), dont on a vu qu'elle \'equivalait \`a l'appartenance de ${\bf G}'_{\sharp}$ \`a ${\cal E}_{G}(G_{\sharp},V)$. Cela d\'emontre (3).

On fixe $S$ assez grand. Gr\^ace \`a (3), l'expression (2) se r\'ecrit 
 $$\tau'(G_{\sharp})^{-1}\vert {\cal U}_{S}^V\vert ^{-1}\sum_{u\in {\cal U}_{S}^V}  \sum_{{\bf G}'_{\sharp}\in {\cal E}(G_{\sharp},V)} i(G_{\sharp},G'_{\sharp})S^{G'_{\sharp,1}}_{\lambda_{1}}(SA^{G'_{\sharp,1}}_{unip,\lambda_{1}}(V),f_{\sharp,1}[^{u}K^V]).$$
 La proposition 4.3 nous dit que le premier terme du membre de droite de (1) est \'egal \`a
$$\tau'(G_{\sharp})^{-1}\vert {\cal U}_{S}^V\vert ^{-1}\sum_{u\in {\cal U}_{S}^V}I^{G_{\sharp}}(A_{unip}^{G_{\sharp}}(V,{^{u}K}_{\sharp}^V),f_{\sharp}). $$
Ainsi le membre de droite de la formule (1) est le produit de $\tau'(G_{\sharp})^{-1}\vert {\cal U}_{S}^V\vert ^{-1}$ et de la somme en $u\in {\cal U}_{S}^V$ de l'expression
$$I^{G_{\sharp}}(A_{unip}^{G_{\sharp}}(V,{^{u}K}_{\sharp}^V),f_{\sharp})-\sum_{{\bf G}'_{\sharp}\in {\cal E}(G_{\sharp},V)} i(G_{\sharp},G'_{\sharp})S^{G'_{\sharp,1}}_{\lambda_{1}}(SA^{G'_{\sharp,1}}_{unip,\lambda_{1}}(V),f_{\sharp,1}[^{u}K^V]).$$
Celle-ci est par d\'efinition \'egale \`a $I^{G_{\sharp}}(SA^{G_{\sharp}}_{unip}(V),f_{\sharp})$. La r\'ef\'erence aux compacts dispara\^{\i}t puisqu'on sait que la distribution $SA^{G_{\sharp}}_{unip}(V)$ n'en d\'epend pas. Ce terme \'etant ind\'ependant de $u$, l'\'egalit\'e (1) devient simplement
$$\tau'(G)I^G(SA^G_{unip}(V),f)=\tau'(G_{\sharp})^{-1}I^{G_{\sharp}}(SA^{G_{\sharp}}_{unip}(V),f_{\sharp}).$$
Dire que cette \'egalit\'e est v\'erifi\'ee pour tout $f$ \'equivaut \`a l'assertion de l'\'enonc\'e. $\square$

    \bigskip
  
  \section{Descente}

\bigskip

\subsection{Une premi\`ere transformation}
On commence la preuve du th\'eor\`eme 3.3. Le triplet $(G,\tilde{G},{\bf a})$ est quelconque. On fixe ${\cal X}\in {\bf Stab}(\tilde{G}(F))$ qui n'appartient pas \`a ${\bf Stab}_{excep}(\tilde{G}(F))$. On fixe un ensemble fini $V$ de places de $F$ contenant $S({\cal X})$.  En vertu du lemme 3.7, on impose de plus

(1) $V$ contient $S({\cal X},\tilde{K})$;

(2) ${\cal X}$ est elliptique;

(3) pour toute place $v\in Val(F)$, l'image de ${\cal X}$ dans ${\bf Stab}(\tilde{G}(F_{v}))$ appartient \`a l'image de l'application $\chi^{\tilde{G}_{v}}$;

(4) les restrictions de $\omega$ \`a $Z(G;{\mathbb A}_{F})^{\theta}$ et $Z(\bar{G};{\mathbb A}_{F})$ sont triviales. 

On pose
$$\underline{A}^{G,{\cal E}}(V,{\cal X},\omega)=\sum_{{\bf G}'\in {\cal E}(\tilde{G},{\bf a},V)}\sum_{{\cal X}'\in {\bf Stab}(\tilde{G}'(F)); {\cal X}'\mapsto {\cal X}}i(\tilde{G},\tilde{G}')transfert(\underline{SA}^{{\bf G}'}(V,{\cal X}')).$$
 Il s'agit de d\'emontrer l'\'egalit\'e
$$A^{G}(V,{\cal X},\omega)=\underline{A}^{G,{\cal E}}(V,{\cal X},\omega).$$

  On fixe une paire de Borel \'epingl\'ee $\hat{\cal E}=(\hat{B},\hat{T},(\hat{E}_{\alpha})_{\alpha\in \Delta})$ de $\hat{G}$ et on d\'efinit l'automorphisme $\hat{\theta}$ ainsi que l'action galoisienne relativement \`a cette paire, cf. [I] 1.2. Notons $E_{\hat{T}}(\tilde{G},{\bf a},V)$ l'ensemble des donn\'ees endoscopiques de $(G,\tilde{G},{\bf a})$ qui sont de la forme ${\bf G}'=(G',{\cal G}',s\hat{\theta})$ avec $s\in \hat{T}$ et qui sont elliptiques, relevantes  et non ramifi\'ees hors de $V$.  Pour deux donn\'ees  ${\bf G}'_{1}=(G'_{1},{\cal G}'_{1},s_{1}\hat{\theta})$ et ${\bf G}'_{2}=(G'_{2},{\cal G}'_{2},s_{2}\hat{\theta})$ dans cet ensemble, 
  disons qu'elles sont $\hat{T}$-\'equivalentes s'il existe $x\in \hat{T}$ tel que $x{\cal G}'_{1}x^{-1}={\cal G}'_{2}$ et $xs_{1}\hat{\theta}(x)^{-1}\in s_{2}Z(\hat{G})$. Fixons un ensemble de repr\'esentants  ${\cal E}_{\hat{T}}(\tilde{G}, {\bf a},V)$ des     classes de $\hat{T}$-\'equivalence dans $E_{\hat{T}}(\tilde{G},{\bf a},V)$.     Toute donn\'ee endoscopique elliptique, relevante et non ramifi\'ee hors de $V$  est \'equivalente \`a une donn\'ee  appartenant \`a $E_{\hat{T}}(\tilde{G},{\bf a},V)$. Il y a donc une application surjective
$${\cal E}_{\hat{T}}(\tilde{G}, {\bf a},V)\to {\cal E}(\tilde{G},{\bf a},V).$$
 Soit ${\bf G}'=(G',{\cal G}',s\hat{\theta})\in {\cal E}_{\hat{T}}(\tilde{G},{\bf a},V)$. La fibre de cette application au-dessus de l'image de ${\bf G}'$ est form\'e des $(G'_{1}, {\cal G}'_{1},s_{1}\hat{\theta})$, \`a $\hat{T}$-\'equivalence pr\`es, pour lesquels il existe $x\in \hat{G}$ de sorte que $x{\cal G}'x^{-1}={\cal G}'_{1}$ et $xs\hat{\theta}(x)^{-1}\in s_{1}Z(\hat{G})$. Puisque $\hat{G'}$ et $\hat{G}'_{1}$ contiennent  $\hat{T}^{\hat{\theta},0}$, on  peut supposer que $x$ normalise $\hat{T}^{\hat{\theta},0}$, donc aussi $\hat{T}$. La condition  $xs\hat{\theta}(x)^{-1}\in s_{1}Z(\hat{G})$ implique alors que l'image de $x$ dans $W$ est fixe par $\hat{\theta}$. Inversement, un \'el\'ement $x\in Norm_{\hat{G}}(\hat{T})$ dont l'image dans $W$ est fixe par $\hat{\theta}$ donne naissance \`a une donn\'ee $(G'_{1}, {\cal G}'_{1},s_{1}\hat{\theta})$. Cette donn\'ee est $\hat{T}$-\'equivalente \`a ${\bf G}'$ si et seulement si $x\in \hat{T}(Aut({\bf G}')\cap Norm_{\hat{G}}(\hat{T}))$.  On a une suite exacte
 $$1\to W^{G'}\to \hat{T}\backslash \hat{T}(Aut({\bf G}')\cap Norm_{\hat{G}}(\hat{T}))\to Out({\bf G}')\to 1$$
  Ainsi la fibre de l'application pr\'ec\'edente au-dessus de l'image de ${\bf G}'$  a pour nombre d'\'el\'ements $\vert W^{\theta}\vert \vert Out({\bf G}')\vert ^{-1}\vert W^{G'}\vert ^{-1}$. On peut donc r\'ecrire
 $$\underline{A}^{G,{\cal E}}(V,{\cal X},\omega)=\vert W^{\theta}\vert ^{-1}\sum_{{\bf G}'\in {\cal E}_{\hat{T}}(G,{\bf a},V)} \vert Out({\bf G}')\vert\vert W^{G'}\vert $$
 $$i(\tilde{G},\tilde{G}')\sum_{{\cal X}'\in {\bf Stab}(G'(F)); {\cal X}'\mapsto {\cal X} }  transfert(\underline{SA}^{{\bf G}'}(V,{\cal X}')).$$
 Soit ${\bf G}'=(G',{\cal G}',s\hat{\theta})\in {\cal E}_{\hat{T}}(G,{\bf a},V)$. On fixe une paire de Borel \'epingl\'ee de $\hat{G}'$ dont la paire de Borel sous-jacente $(\hat{B}',\hat{T}')$ est $(\hat{B}\cap \hat{G}',\hat{T}^{\hat{\theta},0})$. En utilisant les notations de 1.7, cela d\'efinit l'application $\xi:T^*\to T^{_{'}*}$ et une application $Stab(\tilde{G}'(F))\to Stab(\tilde{G}(F))$ que l'on note simplement $(\mu',\omega_{\bar{G}'})\mapsto (\mu,\omega_{\bar{G}})$. Notons ici $p_{\tilde{G}}:Stab(\tilde{G}(F))\to {\bf Stab}(\tilde{G}(F))$ l'application naturelle. 
 
 {\bf Remarque.} Les hypoth\`eses que ${\bf G}'$ est relevante et non ramifi\'ee hors de $V$ n'interviennent pas ici. Les m\^emes notations seront utilis\'ees plus loin pour des donn\'ees ne v\'erifiant pas ces hypoth\`eses.
 
 \bigskip
 
 Pour ${\cal Y}\in{\bf  Stab}(\tilde{G}(F))$, notons $Fib( {\cal Y})$ la fibre de $p_{\tilde{G}}$ au-dessus de $ {\cal Y}$. On adopte de m\^emes notations pour $\tilde{G}'$. Sommer sur les ${\cal X}'\in {\bf Stab}(\tilde{G}'(F))$ tels que ${\cal X}'\mapsto {\cal X}$ revient \`a sommer sur les $(\mu,\omega_{\bar{G}})\in Stab(\tilde{G}(F))$ tels que $p_{\tilde{G}}(\mu,\omega_{\bar{G}})={\cal X}$ et sur les $(\mu',\omega_{\bar{G}'})$ tels que $(\mu',\omega_{\bar{G}'})\mapsto (\mu,\omega_{\bar{G}})$, \`a condition d'affecter les termes d'un coefficient $\vert Fib(p_{\tilde{G}'}(\mu',\omega_{\bar{G}'}))\vert ^{-1}$. Ce nombre d'\'el\'ements se calcule ais\'ement \`a l'aide de 1.1(4). Notons $Fix^{G'}(\mu',\omega_{\bar{G}'})$ le groupe des $w\in W^{G'}$ tels que $w\mu'=\mu'$, $w(\Sigma_{+}(\mu'))=\Sigma_{+}(\mu')$ et $w\omega_{\bar{G}'}(\sigma)\sigma_{G^{_{'}*}}(w)^{-1}=\omega_{\bar{G}'}(\sigma)$ pour tout $\sigma\in \Gamma_{F}$. Alors
 $$(4) \qquad \vert Fib(p_{\tilde{G}'}(\mu',\omega_{\bar{G}'}))\vert=\vert W^{G'}\vert\vert W^{G'}(\mu')\vert ^{-1}\vert Fix^{G'}(\mu',\omega_{\bar{G}'})\vert ^{-1}.$$
On obtient
  $$\sum_{{\cal X}'\in {\bf Stab}(G'(F)); {\cal X}'\mapsto {\cal X} } transfert(\underline{SA}^{{\bf G}'}(V,{\cal X}'))=\sum_{(\mu,\omega_{\bar{G}})\in Stab(\tilde{G}(F)), p_{\tilde{G}}(\mu,\omega_{\bar{G}})={\cal X}}$$
  $$ \sum_{(\mu',\omega_{\bar{G}'})\in Stab(\tilde{G}'(F)), (\mu',\omega_{\bar{G}'})\mapsto (\mu,\omega_{\bar{G}})} \vert W^{G'}\vert^{-1}\vert W^{G'}(\mu')\vert \vert Fix^{G'}(\mu',\omega_{\bar{G}'})\vert $$
  $$transfert(\underline{SA}^{{\bf G}'}(V,p_{\tilde{G}'}(\mu',\omega_{\bar{G}'}))).$$
  Cela conduit \`a l'\'egalit\'e
$$\underline{A}^{\tilde{G},{\cal E}}(V,{\cal X},\omega)=\vert W^{\theta}\vert ^{-1}  \sum_{{\bf G}'\in {\cal E}_{\hat{T}}(G,{\bf a},V)}\sum_{(\mu,\omega_{\bar{G}})\in Stab(\tilde{G}(F)), p_{\tilde{G}}(\mu,\omega_{\bar{G}})={\cal X}}$$
$$ \sum_{(\mu',\omega_{\bar{G}'})\in Stab(\tilde{G}'(F)), (\mu',\omega_{\bar{G}'})\mapsto (\mu,\omega_{\bar{G}})} i(\tilde{G},\tilde{G}',\mu',\omega_{\bar{G}'}) transfert(\underline{SA}^{{\bf G}'}(V,p_{\tilde{G}'}(\mu',\omega_{\bar{G}'}))),$$
o\`u
$$(5) \qquad  i(\tilde{G},\tilde{G}',\mu',\omega_{\bar{G}'})=i(\tilde{G},\tilde{G}')\vert Out({\bf G}')\vert \vert W^{G'}(\mu')\vert \vert Fix^{G'}(\mu',\omega_{\bar{G}'})\vert.$$
  Pour tout $(\mu,\omega_{\bar{G}})\in Fib({\cal X})$, on a la formule parall\`ele \`a (4):
$$\vert Fib({\cal X})\vert =\vert W^{\theta}\vert \vert W^{G}(\mu)\vert ^{-1}\vert Fix^G(\mu,\omega_{\bar{G}})\vert ^{-1}.$$
Posons
 $$(6) \qquad \underline{A}^{\tilde{G},{\cal E}}(V,\mu,\omega_{\bar{G}},\omega)=\vert W^{G}(\mu)\vert^{-1} \vert Fix^G(\mu,\omega_{\bar{G}})\vert^{-1} \sum_{{\bf G}'\in {\cal E}_{\hat{T}}(G,{\bf a},V)}$$
 $$ \sum_{(\mu',\omega_{\bar{G}'})\in Stab(\tilde{G}'(F)), (\mu',\omega_{\bar{G}'})\mapsto (\mu,\omega_{\bar{G}})} i(\tilde{G},\tilde{G}',\mu',\omega_{\bar{G}'}) transfert(\underline{SA}^{{\bf G}'}(V,p_{\tilde{G}'}(\mu',\omega_{\bar{G}'}))).$$
 On obtient la formule
 $$  \underline{A}^{G,{\cal E}}(V,{\cal X},\omega)=\vert Fib({\cal X})\vert ^{-1}\sum_{(\mu,\omega_{\bar{G}})\in Fib({\cal X})}\underline{A}^{\tilde{G},{\cal E}}(V,\mu,\omega_{\bar{G}},\omega).$$
 Enon\c{c}ons une forme plus pr\'ecise du th\'eor\`eme 3.3, qui implique celui-ci d'apr\`es la formule pr\'ec\'edente.
 
 \ass{Proposition}{Soient ${\cal X}$ et $V$ comme plus haut. Pour tout $(\mu,\omega_{\bar{G}})\in Fib({\cal X})$, on a l'\'egalit\'e
 $$  \underline{A}^{\tilde{G},{\cal E}}(V,\mu,\omega_{\bar{G}},\omega)=A^{\tilde{G}}(V,{\cal X},\omega).$$}

 C'est cette assertion que nous allons prouver. On fixe jusqu'\`a la fin de la preuve un \'el\'ement 
 $(\mu,\omega_{\bar{G}})\in Fib({\cal X})$.
 
 \bigskip
 
 \subsection{Descente des donn\'ees endoscopiques}
 On identifie $\underline{la}$ paire de Borel \'epingl\'ee ${\cal E}^*=(B^*,T^*,(E_{\alpha})_{\alpha\in \Delta})$ de $G$ \`a une paire particuli\`ere. 
  On rel\`eve $\mu$ en $(\nu,e)$ avec $\nu\in T^*$ et $e\in Z(\tilde{G},{\cal E}^*)$. On pose $\eta=\nu e$ et $\bar{G}=G_{\eta}$. On note $\bar{T}$ le tore $T^{*,\theta^*,0}$ vu comme un sous-tore maximal de $\bar{G}$. On rappelle que $\Sigma(\mu)$ s'identifie \`a l'ensemble de racines $\Sigma^{\bar{G}}(\bar{T})$ de $\bar{T}$ dans $\bar{G}$. On fixe une paire de Borel \'epingl\'ee $\bar{{\cal E}}$ de $\bar{G}$ dont la paire de Borel sous-jacente soit $( \bar{B},\bar{T})$, o\`u $\bar{B}=B^*\cap \bar{G}$, et on munit $\bar{G}$ de l'unique action galoisienne $\sigma\mapsto \sigma_{\bar{G}}$ conservant $\bar{{\cal E}}$ et co\"{\i}ncidant sur  $\bar{T}=T^{*,\theta,0}$  avec l'action $\sigma\mapsto \omega_{\bar{G}}(\sigma)\sigma_{G^*}$. Il est clair que cette action galoisienne s'\'etend en une action sur le groupe $I_{\eta}=Z(G)^{\theta}\bar{G}$ et on munit ce groupe de cette action.  On introduit le groupe dual $\hat{\bar{G}}$, muni d'une paire de Borel \'epingl\'ee dont on note la paire sous-jacente $(\hat{\bar{B}},\hat{\bar{T}})$. On peut identifier $\hat{\bar{T}}$   \`a $\hat{T}/(1-\hat{\theta})(\hat{T})$, muni de l'action galoisienne  $\sigma\mapsto \omega_{\bar{G}}(\sigma)\sigma_{G^*}$. L'ensemble de racines $\Sigma^{\hat{\bar{G}}}(\hat{\bar{T}})$ est en bijection avec $\Sigma(\mu)$ et le sous-ensemble positif d\'etermin\'e par $\hat{\bar{B}}$ correspond au sous-ensemble $\Sigma_{+}(\mu)$.

 Consid\'erons une donn\'ee endoscopique de $(G,\tilde{G},{\bf a})$ de la forme ${\bf G}'=(G',{\cal G}',s\hat{\theta}) $, avec $s\in \hat{T}$. On la suppose elliptique. Soit $(\mu',\omega_{\bar{G}'})\in Stab(\tilde{G}'(F))$ tel que $(\mu',\omega_{\bar{G}'})\mapsto (\mu,\omega_{\bar{G}})$.    
  Posons $\hat{\bar{T}}_{ad}=\hat{\bar{T}}/Z(\hat{\bar{G}})$. C'est un quotient de $\hat{T}$.  Notons $\bar{s}$ l'image de $s$ dans ce groupe. Posons $\hat{\bar{H}}=Z_{\hat{\bar{G}}_{AD}}(\bar{s})^0$.  On munit ce groupe de la paire de Borel $(\hat{\bar{B}}_{ad}\cap \hat{\bar{H}},\hat{\bar{T}}_{ad})$. D'apr\`es la d\'efinition de 1.7, il y a une unique cocha\^{\i}ne $\sigma\mapsto \omega_{\bar{H}}(\sigma)$ de $\Gamma_{F}$ dans $W(\mu)$ de sorte que l'on ait l'\'egalit\'e
 $$\omega_{\bar{G}'}(\sigma)\omega_{G'}(\sigma)=\omega_{\bar{H}}(\sigma)\omega_{\bar{G}}(\sigma)$$
 pour tout $\sigma\in \Gamma_{F}$. Cette cocha\^{\i}ne  s'\'etend  en une cocha\^{\i}ne d\'efinie sur $W_{F}$. Puisque $W(\mu)=W^{\hat{\bar{G}}}=W^{\hat{\bar{G}}_{AD}}$, pour tout $w\in W_{F}$, on peut relever $\omega_{\bar{H}}(w)$ en un \'el\'ement $\bar{g}(w)\in \hat{\bar{G}}_{AD}$ qui normalise $\hat{\bar{T}}_{ad}$. On d\'efinit le sous-groupe $\bar{{\cal H}}\subset {^L(\bar{G}_{SC})}= \hat{\bar{G}}_{AD}\rtimes W_{F}$ engendr\'e par $\hat{\bar{H}}$ et les $(\bar{g}(w),w)$ pour $w\in W_{F}$. On a prouv\'e en [W1] 3.5 que ce groupe s'ins\'erait dans une extension scind\'ee
 $$1\to \hat{\bar{H}}\to \bar{{\cal H}}\to W_{F}\to 1.$$
 On peut donc munir $\hat{\bar{H}}$ d'une $L$-action galoisienne compatible avec cette extension. On introduit le groupe r\'eductif $\bar{H}$ d\'efini et quasi-d\'eploy\'e sur $F$ tel que $\hat{\bar{H}}$, muni de cette action galoisienne, soit le groupe dual de $\bar{H}$. On a prouv\'e en [W1] 3.5 que le triplet ${\bf \bar{H}}=(\bar{H},\bar{{\cal H}},\bar{s})$ \'etait une donn\'ee endoscopique de $\bar{G}_{SC}$. Il s'agit ici d'endoscopie non tordue, toute torsion et tout caract\`ere ont disparu. On fixe une paire de Borel $(B^{\bar{H}},T^{\bar{H}})$ de $\bar{H}$ d\'efinie sur $F$.
 
 {\bf Remarque.} Dans [W1],  la situation de d\'epart n'\'etait pas la m\^eme qu'ici, on partait d'un diagramme.   Mais on v\'erifie ais\'ement que les pr\'esentes hypoth\`eses sont suffisantes   pour assurer la validit\'e des propri\'et\'es \'enonc\'ees ci-dessus.  La m\^eme remarque vaut pour la suite. D'autre part, dans [W1],
 le corps de base \'etait local non-archim\'edien mais les constructions valent \'evidemment pour tout corps de base.  Enfin, on a d\'ej\`a repris cette construction dans la section 5 de [III], o\`u l'on a not\'e $\bar{{\bf G}}'$ la donn\'ee $\bar{{\bf H}}$. Il nous semble plus clair de revenir ici \`a la notation de [W1], la notation $\bar{G}'$ \'etant d\'ej\`a utilis\'ee.  
 
 \bigskip
 
 Fixons un \'el\'ement $\epsilon\in \tilde{G}'_{ss}(F)$ et une paire de Borel  $(B_{\epsilon},T_{\epsilon})$ v\'erifiant les conditions du (ii) du lemme 1.3 relativement \`a l'\'el\'ement $(\mu',\omega_{\bar{G}'})$. Le groupe $G'_{\epsilon}$ est quasi-d\'eploy\'e et est  muni de la paire de Borel $(B^*_{\flat},T_{\epsilon})$ d\'efinie dans ce lemme. Fixons une paire de Borel $(B^{G'},T^{G'})$ de $G'$ d\'efinie sur $F$. Par la construction de ce lemme, le tore $T_{\epsilon}$ s'identifie \`a $T^{G'}$, l'identification n'\'etant pas $\Gamma_{F}$-\'equivariante: elle identifie l'action galoisienne sur $T_{\epsilon}$ avec l'action $\sigma\mapsto \omega_{\bar{G}'}(\sigma)\sigma_{G'}$ sur $T^{G'}$. On a alors un diagramme
 $$(1) \qquad\begin{array}{ccccc} T^{\bar{H}}_{sc}&&&& T_{\epsilon,sc}\\ \downarrow&&&&\\    T^{\bar{H}}&&&&\downarrow \\ \parallel &&&& \\ \bar{T}_{sc}&&&& T_{\epsilon}\\ \downarrow&&&&  \\ \bar{T}&&&&\parallel\\ \downarrow&&&&\\ T^*&\to&T^*/(1-\theta^*)(T^*)&\simeq&T^{G'}\\ \end{array}$$
 (on a not\'e $T^{\bar{H}}_{sc}$, resp. $T_{\epsilon,sc}$,  $\bar{T}_{sc}$, les images r\'eciproques de $T^{\bar{H}}$ dans $\bar{H}_{SC}$, resp. de $T_{\epsilon}$ dans $G'_{\epsilon,SC}$, resp. de $\bar{T}$ dans $\bar{G}_{SC}$). Pour tout tore $S$, posons $X_{*,{\mathbb Q}}(S)=X_{*}(S)\otimes_{{\mathbb Z}}{\mathbb Q} $. Le lemme [W1] 3.6 affirme qu'il se d\'eduit de ce diagramme un isomorphisme
 $$X_{*,{\mathbb Q}}(T^{\bar{H}}_{sc}) \simeq X_{*,{\mathbb Q}}(T_{\epsilon,sc}) $$
 qui est \'equivariant pour les actions galoisiennes et gr\^ace auquel les deux groupes $\bar{H}_{SC}$ et $G'_{\epsilon,SC}$ forment une paire endoscopique non standard.
 
  Le lemme [W1] 3.6 nous dit que du diagramme (1) se d\'eduit un diagramme commutatif
 $$(2) \qquad \begin{array}{ccc}X_{*,{\mathbb Q}}(Z(G)^{\theta,0})&\to&X_{*,{\mathbb Q}}(Z(G')^0)\\ \downarrow&&\\ X_{*,{\mathbb Q}}(Z(\bar{G})^0)&&\downarrow\\ \downarrow&&\\  X_{*,{\mathbb Q}}(Z(\bar{G})^0)\oplus X_{*,{\mathbb Q}}(Z(\bar{H})^0)&\simeq& X_{*,{\mathbb Q}}(Z(G'_{\epsilon})^0)\\ \end{array}$$
 Il est form\'e d'applications injectives \'equivariantes pour les actions galoisiennes et la fl\`eche du bas est un isomorphisme.

 On a
 
 (3) supposons $(\mu',\omega_{\bar{G}'})$ elliptique; alors  la donn\'ee endoscopique $ \bar{{\bf H}}$ de $\bar{G}_{SC}$ est elliptique.
 
 Preuve.  L'hypoth\`ese implique que  $\epsilon$ est elliptique dans $\tilde{G}'(F)$.   Dans le diagramme (2), on prend les invariants par $\Gamma_{F}$. La fl\`eche verticale de droite devient un isomorphisme par ellipticit\'e de $\epsilon$. La fl\`eche horizontale du haut aussi par ellipticit\'e de ${\bf G}'$. Donc les deux fl\`eches de gauche deviennent aussi des isomorphismes. Cela entra\^{\i}ne $X_{*,{\mathbb Q}}(Z(\bar{H})^0)^{\Gamma_{F}}=\{0\}$ et l'assertion. $\square$

On a

(4) supposons $S(p_{\tilde{G}'}(\mu',\omega_{\bar{G}'}),\tilde{K}')\subset V$; alors   la donn\'ee $ \bar{{\bf H}}$ est non ramifi\'ee hors de $V$. 

La preuve est la m\^eme qu'en 1.7(4). Soit $v\not\in V$ et $\sigma$ dans le groupe d'inertie $I_{v}$. L'hypoth\`ese implique  $\omega_{\bar{G}'}(\sigma)=1$. On a $\omega_{G'}(\sigma)=1$ puisque ${\bf G}'$ n'est pas ramifi\'ee en $v$.   Donc $\omega_{\bar{G}}(\sigma)=\omega_{\bar{H}}(\sigma)^{-1}$. Cet \'el\'ement appartient \`a $ W^{\bar{G}}$ car $\omega_{\bar{H}}(\sigma)\in W^{\bar{G}}$ et il conserve l'ensemble des racines positives de $\bar{T}$ dans $\bar{G}$ car c'est le cas de $\omega_{\bar{G}}(\sigma)$.   Ces deux conditions entra\^{\i}nent qu'il est \'egal \`a $1$, c'est-\`a-dire $\omega_{\bar{H}}(\sigma)=1$. Dans le cas d'endoscopie non tordue, cela suffit \`a assurer que la donn\'ee $(\bar{H},\bar{{\cal H}},\bar{s})$ est non ramifi\'ee en $v$. $\square$

Posons quelques d\'efinitions. Soit $v$ une place de $F$. Soient $\bar{G}_{0}$ un groupe r\'eductif connexe d\'efini sur $F_{v}$ et $\psi_{v}:\bar{G}_{0}\to\bar{G}$ un torseur int\'erieur (pour l'action de $\Gamma_{F_{v}}$).  Soient $(B_{\flat},T_{\flat})$, $(B_{\bar{H}},T_{\bar{H}})$, $(B_{\natural},T_{\natural})$ et $(B_{\natural,0},T_{\natural,0})$ des paires de Borel respectivement de $G'_{\epsilon}$, $\bar{H}$, $\bar{G}$ et $\bar{G}_{0}$ d\'efinies sur $\bar{F}_{v}$. En conjuguant les trois premi\`eres paires en celles fix\'ees plus haut, on obtient des isomorphismes
$$(5) \qquad X_{*,{\mathbb Q}}(T_{\bar{H},sc})\simeq X_{*,{\mathbb Q}}(T_{\flat,sc})$$
et
$$(6) \qquad T_{\bar{H}}\simeq T_{\natural,sc}.$$
La donn\'ee de $\psi_{v}$ \'etablit aussi un isomorphisme
$$(7) \qquad T_{\natural}\simeq T_{\natural,0}.$$
Il se d\'eduit de (6) et (7) un isomorphisme
$$(8)\qquad T_{\bar{H}}\simeq T_{\natural,0,sc}$$
qui ne d\'epend pas de la paire $(B_{\natural},T_{\natural})$. 
On dit que $(B_{\flat},T_{\flat})$ et $(B_{\bar{H}},T_{\bar{H}})$ se correspondent si $T_{\flat}$ et $T_{\bar{H}}$ sont d\'efinis sur $F_{v}$ et que  (5) est \'equivariant pour l'action de $\Gamma_{F_{v}}$. On d\'efinit de m\^eme la notion de correspondance entre $(B_{\bar{H}},T_{\bar{H}})$ et $(B_{\natural},T_{\natural})$, resp. $(B_{\natural},T_{\natural})$ et $(B_{\natural,0},T_{\natural,0})$, resp.  $(B_{\bar{H}},T_{\bar{H}})$ et $(B_{\natural,0},T_{\natural,0})$, en rempla\c{c}ant l'isomorphisme (5) respectivement par (6), (7) et (8). 

Les m\^emes d\'efinitions peuvent \^etre pos\'ees sur le corps de base $F$, en rempla\c{c}ant simplement ci-dessus $F_{v}$ et $\bar{F}_{v}$ par $F$ et $\bar{F}$.
  
\bigskip

\subsection{La sous-somme attach\'ee \`a une donn\'ee endoscopique ${\bf H}$}
Consid\'erons la d\'efinition 5.1(6). Dans la somme en $(\mu',\omega_{\bar{G}'})$, on peut ajouter les conditions que $(\mu',\omega_{\bar{G}'})$ est  elliptique et que $S(p_{\tilde{G}'}(\mu',\omega_{\bar{G}'}),\tilde{K}')\subset V$. En effet, si elles ne sont pas v\'erifi\'ees, la distribution  $\underline{SA}^{{\bf G}'}(V,p_{\tilde{G}'}(\mu',\omega_{\bar{G}'}))$ est nulle. 
 Notons $E_{\hat{\bar{T}}_{ad},\star}(\bar{G}_{SC},V)$ l'ensemble des donn\'ees endoscopiques de $\bar{G}_{SC}$ de la forme ${\bf \bar{H}}=(\bar{H},\bar{{\cal H}},\bar{s})$, o\`u $\bar{s}\in \hat{\bar{T}}_{ad}$, qui sont elliptiques et non ramifi\'ees hors de $V$. A partir de  ${\bf G}'=(G',{\cal G}',s\hat{\theta})\in {\cal E}_{\hat{T}}(G,\omega,V)$ et de $(\mu',\omega_{\bar{G}'})\in Stab(\tilde{G}'(F))$  v\'erifiant les deux conditions ci-dessus et tel que $(\mu',\omega_{\bar{G}'})\mapsto (\mu,\omega_{\bar{G}})$, on a construit une donn\'ee endoscopique $ {\bf \bar{H}}$ de $\bar{G}_{SC}$, qui appartient \`a $E_{\hat{\bar{T}}_{ad},\star}(\bar{G}_{SC},V)$, cf. 5.2(3) et (4) . Pour ${\bf H}\in E_{\hat{\bar{T}}_{ad},\star}(\bar{G}_{SC},V)$, on note ${\cal J}({\bf H})$ la fibre de cette application. 
 On note $\underline{A}^{\tilde{G},{\cal E}}(V,{\bf H},\omega)$ la sous-somme de l'expression 5.1(6), o\`u on somme sur les triplets $({\bf G}',\mu',\omega_{\bar{G}'})\in {\cal J}({\bf H})$. Autrement dit
$$  \underline{A}^{\tilde{G},{\cal E}}(V,{\bf H},\omega)=\vert W^{G}(\mu)\vert^{-1} \vert Fix^G(\mu,\omega_{\bar{G}})\vert^{-1} \sum_{({\bf G}',\mu',\omega_{\bar{G}'})\in  {\cal J}({\bf H})}$$
 $$  i(\tilde{G},\tilde{G}',\mu',\omega_{\bar{G}'}) transfert(\underline{SA}^{{\bf G}'}(V,p_{\tilde{G}'}(\mu',\omega_{\bar{G}'}))).$$
 
 On a l'\'egalit\'e
 $$(1) \qquad \underline{A}^{\tilde{G},{\cal E}}(V,\mu,\omega_{\bar{G}},\omega)=\sum_{{\bf H}\in E_{\hat{\bar{T}}_{ad},\star}(\bar{G}_{SC},V)} \underline{A}^{\tilde{G},{\cal E}}(V,{\bf H},\omega).$$
 
 {\bf On fixe jusqu'en 8.1 une donn\'ee ${\bf H}\in E_{\hat{\bar{T}}_{ad},\star}(\bar{G}_{SC},V)$.}  On va \'etudier la distribution $ \underline{A}^{\tilde{G},{\cal E}}(V,{\bf H},\omega)$.

  \bigskip
 
  \subsection{Propri\'et\'es de relevance}
  Soit $v\in Val(F)$, pla\c{c}ons-nous sur le corps $F_{v}$. Consid\'erons l'ensemble $D_{v}$ des couples $(\eta_{v},r_{v})\in \tilde{G}(F_{v})\times G(\bar{F}_{v})$ tels que
  
  (1) $r_{v}\eta_{v}r_{v}^{-1}=\eta$;
  
  (2) en utilisant la paire de Borel $ad_{r_{v}^{-1}}(B^*,T^*)$ dans la construction de 1.2, on ait l'\'egalit\'e $(\mu_{\eta_{v}},\omega_{\eta_{v}})=(\mu,\omega_{\bar{G}_{v}})$.
  
  On peut traduire cette condition (2) de la fa\c{c}on suivante. La condition (1) implique que l'automorphisme $ad_{r_{v}}$ se restreint en un isomorphisme de $G_{\eta_{v}}$ sur $G_{\eta}=\bar{G}$. Le groupe $G_{\eta_{v}}$ \'etant muni de sa structure galoisienne naturelle et le groupe $\bar{G}$ \'etant muni de sa structure quasi-d\'eploy\'ee d\'efinie en 5.2, (2) \'equivaut \`a ce que $ad_{r_{v}}$ soit un torseur int\'erieur de $G_{\eta_{v}}$ sur $\bar{G}$ que l'on note $\psi_{r_{v}}$.  
  
  Le groupe $I_{\eta}=I_{\eta}(\bar{F}_{v})$ agit \`a gauche sur $D_{v}$ par $(x,(\eta_{v},r_{v}))\mapsto (\eta_{v},xr_{v})$. Le groupe $G(F_{v})$ agit \`a droite par $((\eta_{v},r_{v}),g)\mapsto
(g^{-1}\eta_{v}g_{v},r_{v}g)$.

 On se rappelle notre hypoth\`ese 5.1(3): l'image de ${\cal X}$ dans ${\bf Stab}(\tilde{G}(F_{v}))$ appartient \`a l'image de l'application $\chi^{\tilde{G}_{v}}$.  Montrons que
  
 (3) l'ensemble $D_{v}$ est non vide;
 
 \noindent fixons $ (\eta_{v},r_{v})\in D_{v}$ et rappelons que l'on a d\'efini l'ensemble ${\cal Y}_{\eta_{v}}$ en 1.2;  on a

 (4) l'ensemble $D_{v}$ co\"{\i}ncide avec l'ensemble des $(y^{-1}\eta_{v}y,r_{v}y)$ pour $y\in {\cal Y}_{\eta_{v}}$.
 
 Preuve. La proposition 1.2 vaut bien s\^ur aussi sur le corps de base $F_{v}$. L'hypoth\`ese 5.1(3) implique que l'on peut fixer $\eta_{v}\in \tilde{G}_{ss}(F_{v})$ et une paire de Borel $(B,T)$ conserv\'ee par $ad_{\eta_{v}}$ de sorte que le couple $(\mu_{\eta_{v}},\omega_{\eta_{v}})$ d\'eduit de ces donn\'ees ait m\^eme image que $(\mu,\omega_{\bar{G}_{v}})$ dans ${\bf Stab}(\tilde{G}(F_{v}))$. D'apr\`es 1.2(3), on peut modifier $(B,T)$ de sorte que  $(\mu_{\eta_{v}},\omega_{\eta_{v}})=(\mu,\omega_{\bar{G}_{v}})$. Fixons $r'_{v}\in G=G(\bar{F}_{v})$ tel que $(B,T)=ad_{{r'_{v}}^{-1}}(B^*,T^*)$. Par construction, l'\'egalit\'e $\mu_{\eta_{v}}=\mu$ \'equivaut \`a la relation $r'_{v}\eta_{v}{r'_{v}}^{-1}\in (1-\theta^*)(T^*)\eta$. En multipliant $r'_{v}$ \`a gauche par un \'el\'ement convenable de $T^*$, on obtient un \'el\'ement $r_{v}$ tel que (1)  soit v\'erifi\'ee. La propri\'et\'e (2) l'est aussi par construction. Alors $(\eta_{v},r_{v})\in D_{v}$, ce qui prouve (3). Soit $(\eta_{v},r_{v})\in D_{v}$. Pour $y\in {\cal Y}_{\eta_{v}}$, il est clair que l'\'el\'ement  $(y^{-1}\eta_{v}y,r_{v}y)$ v\'erifie (1). Il r\'esulte du (iii) de la proposition 1.2 que le couple  $(\omega_{y^{-1}\eta_{v}y},\omega_{y^{-1}\eta_{v}y})$ associ\'e \`a $y^{-1}\eta_{v}y$ et \`a la paire $ad_{y^{-1}r_{v}^{-1}}(B^*,T^*)$ est \'egal au couple $(\mu_{\eta_{v}},\omega_{\eta_{v}})$ associ\'e \`a $\eta_{v}$ et \`a la paire $ad_{r_{v}^{-1}}(B^*,T^*)$, donc \`a $(\mu,\omega_{\bar{G}_{v}})$. Donc  $(y^{-1}\eta_{v}y,r_{v}y)\in D_{v}$. Inversement, soit $(\eta'_{v},r'_{v})\in D_{v}$. La m\^eme proposition 1.2(iii) implique qu'il existe $y_{1}\in {\cal Y}_{\eta_{v}}$ de sorte que $\eta'_{v}=y_{1}^{-1}\eta_{v}y_{1}$ et $ad_{{r'_{v}}^{-1}}(B^*,T^*)=ad_{y_{1}^{-1}r_{v}^{-1}}(B^*,T^*)$. Cette derni\`ere \'egalit\'e implique $r'_{v}\in T^*r_{v}y_{1}$. En posant $T=ad_{r_{v}^{-1}}(T^*)$, cela \'equivaut \`a l'existence de $t\in T$ tel que $r'_{v}=r_{v}ty_{1}$. Les \'egalit\'es $r'_{v}\eta'_{v}{r'_{v}}^{-1}=\eta=r_{v}\eta_{v}r_{v}^{-1}$ et $\eta'_{v}=y_{1}^{-1}\eta_{v}y_{1}$ impliquent alors $t\eta_{v}t^{-1}=\eta_{v}$, donc $t\in T^{\theta}\subset I_{\eta_{v}}$. L'\'el\'ement $y=ty_{1}$ appartient \`a ${\cal Y}_{\eta_{v}}$ et on a l'\'egalit\'e $(\eta'_{v},r'_{v})=(y^{-1}\eta_{v}y,r_{v}y)$. Cela prouve (4). $\square$
  
  Soit  $(\eta_{v},r_{v})\in D_{v}$. La donn\'ee locale ${\bf H}_{v}$ \'etant une donn\'ee endoscopique de $\bar{G}_{SC}$ est aussi une donn\'ee endoscopique pour $G_{\eta_{v},SC}$, via le torseur $\psi_{r_{v}}$ introduit ci-dessus. On dira plus pr\'ecis\'ement que c'est une donn\'ee endoscopique pour $(G_{\eta_{v},SC},\psi_{r_{v}})$. On note $D_{v}^{rel}$ l'ensemble des $(\eta_{v},r_{v})\in D_{v}$ tels que ${\bf H}_{v}$ est relevante pour  $(G_{\eta_{v},SC},\psi_{r_{v}})$.  Pour deux \'el\'ements $(\eta_{v},r_{v})$ et $(\eta'_{v},r'_{v})$ de $D_{v}$ dans la m\^eme double classe modulo  les actions de $I_{\eta}$ \`a gauche et de  $G(F_{v})$ \`a droite, les couples $(G_{\eta_{v}},\psi_{r_{v}})$ et $(G_{\eta'_{v}},\psi_{r'_{v}})$ sont isomorphes. Il en r\'esulte que $D_{v}^{rel}$ est invariant par les actions de $I_{\eta}$  et $G(F_{v})$.

  On a d\'efini l'ensemble ${\cal J}({\bf H})$ en 5.3.  Introduisons  un ensemble ${\cal J}_{\star}({\bf H})$ a priori plus gros.  Notons $E_{\hat{T},\star}(\tilde{G},{\bf a},V)$ l'ensemble des donn\'ees endoscopiques de $(G,\tilde{G},{\bf a})$ qui sont de la forme ${\bf G}'=(G',{\cal G}',s\hat{\theta})$ avec $s\in \hat{T}$ et qui sont elliptiques et non ramifi\'ees hors de $V$. Fixons un ensemble de repr\'esentants ${\cal E}_{\hat{T},\star}(\tilde{G},{\bf a},V)$ des classes de $\hat{T}$-\'equivalence dans $E_{\hat{T},\star}(\tilde{G},{\bf a},V)$. La diff\'erence avec ${\cal E}_{\hat{T}}(\tilde{G},{\bf a},V)$ est que les donn\'ees ${\bf G}'$ ne sont pas suppos\'ees relevantes. Une partie de nos constructions   vaut aussi bien sans cette hypoth\`ese de relevance. On peut supposer ${\cal E}_{\hat{T}}(\tilde{G},{\bf a},V)\subset {\cal E}_{\hat{T},\star}(\tilde{G},{\bf a},V)$. 
 On note ${\cal J}_{\star}({\bf H})$ l'ensemble des triplets $({\bf G}',\mu',\omega_{\bar{G}'})$, o\`u  ${\bf G}'\in {\cal E}_{\hat{T},\star}(\tilde{G},{\bf a},V)$ et $(\mu',\omega_{\bar{G}'})\in {\bf Stab}(\tilde{G}'(F))$, tels que 
 
 - $(\mu',\omega_{\bar{G}'})\mapsto (\mu,\omega_{\bar{G}})$;
 
 - $(\mu',\omega_{\bar{G}'})$ est elliptique et $S(p_{\tilde{G}'}(\mu',\omega_{\bar{G}'}),\tilde{K}')\subset V$;
 
 - ${\bf H}$ est associ\'e \`a  $({\bf G}',\mu',\omega_{\bar{G}'})$ par la construction de 5.2.
 
 Soit $({\bf G}',\mu',\omega_{\bar{G}'})\in {\cal J}_{\star}({\bf H})$. On d\'efinit un \'el\'ement $\epsilon\in \tilde{G}'(F)$ comme en 5.2, dont on reprend les notations.  
 
   \ass{Lemme}{Supposons que la classe de conjugaison stable de $\epsilon$ dans $\tilde{G}'(F_{v})$ corresponde \`a une classe de conjugaison stable dans $\tilde{G}(F_{v})$. Alors $D_{v}^{rel}$ n'est pas vide.}
   
   Preuve. Par hypoth\`ese, il existe un diagramme 
 $(\epsilon,B',T',B,T,\eta_{v})$ joignant $\epsilon$ \`a un \'el\'ement semi-simple $\eta_{v}\in \tilde{G}(F_{v})$. Comme on l'a vu en [I] 1.10, on peut remplacer $B'$ par n'importe quel Borel contenant $T'$, quitte \`a changer $B$. On peut donc supposer que la paire $(B',T')$ est conjugu\'ee \`a la paire $(B_{\epsilon},T_{\epsilon})$ que l'on a fix\'ee en 5.2 par un \'el\'ement de $G'_{\epsilon}$. Comme on l'a vu dans la preuve de la proposition 1.2, la construction de ce paragraphe est insensible \`a une telle conjugaison. Donc, en utilisant  la paire $(B',T')$, on a l'\'egalit\'e $(\mu_{\epsilon},\omega_{\epsilon})=(\mu',\omega_{\bar{G}'_{v}})$. Fixons $r_{v}\in G(\bar{F}_{v})$ tel que $ad_{r_{v}}(B,T)=(B^*,T^*)$.  Dans la preuve du lemme 1.8(i), on a calcul\'e le couple $(\mu_{\eta_{v}},\omega_{\eta_{v}})$ d\'eduit  de la paire $(B,T)$. Il  est \'egal \`a l'image $(\mu,\omega_{\bar{G}_{v}})$ de $(\mu_{\epsilon},\omega_{\epsilon})$. Or $\mu_{\eta_{v}}$ est l'image de   $ad_{r_{v}}(\eta_{v})$  dans $(T^*/(1-\theta^*)(T^*))\times_{{\cal Z}(G)}{\cal Z}(\tilde{G})$. Donc $ad_{r_{v}}(\eta_{v})$ et $\eta$ ont m\^eme image dans cet ensemble. Quitte \`a multiplier $r_{v}$ \`a gauche par un \'el\'ement de $T^*$, on peut supposer $ad_{r_{v}}(\eta_{v})=\eta$. Alors le couple $(\eta_{v},r_{v})$ appartient  \`a $D_{v}$.  Il se d\'eduit de $(B,T)$ une paire de Borel $(B_{\natural},T_{\natural})=(B\cap G_{\eta_{v}},T\cap G_{\eta_{v}})$ de $G_{\eta_{v}}$, puis une paire de Borel $(B_{\natural,sc},T_{\natural,sc})$ de $G_{\eta_{v},SC}$. D'autre part, puisque $\bar{H}_{SC}$ et $G'_{\epsilon,SC}$ sont en situation d'endoscopie non standard, on peut fixer une paire de Borel $(B_{\bar{H}},T_{\bar{H}})$ de $\bar{H}$ (on rappelle que le corps de base est ici $F_{v}$) qui correspond \`a $(B'_{sc},T'_{sc})$. Puisque ${\bf H}_{v}$ est une donn\'ee endoscopique de $(G_{\eta_{v},SC},\psi_{r_{v}})$, les deux paires $(B_{\bar{H}},T_{\bar{H}})$ et $(B_{\natural,sc},T_{\natural,sc})$ d\'efinissent un isomorphisme de $T_{\natural,sc}$ sur $T_{\bar{H}}$. Les propri\'et\'es d'\'equivariance du diagramme 5.2(1) entra\^{\i}nent que cet isomorphisme est d\'efini sur $F_{v}$, autrement dit que les paires ci-dessus se correspondent.  Cela implique que ${\bf H}_{v}$ est relevante pour $(G_{\eta_{v},SC},\psi_{r_{v}})$. $\square$

 Inversement, soit $(\eta_{v},r_{v})\in D_{v}^{rel}$. On peut alors fixer des paires de Borel $(B_{\bar{H}},T_{\bar{H}})$ de $\bar{H}$ et $(B_{\natural},T_{\natural})$ de $G_{\eta_{v}}$ qui se correspondent via le torseur $\psi_{r_{v}}$. Les couples $ad_{r_{v}}(B_{\natural},T_{\natural})$ et $(\bar{B},\bar{T})$ sont des paires de Borel de $G_{\eta}$ sur $\bar{F}_{v}$. On peut fixer $x_{v}\in G_{\eta}$ de sorte que $ad_{x_{v}r_{v}}(B_{\natural},T_{\natural})=(\bar{B},\bar{T})$. Posons $(B,T)=ad_{x_{v}r_{v}}^{-1}(B^*,T^*)$.  On a $B\cap G_{\eta_{v}}=B_{\natural}$ et $T\cap G_{\eta_{v}}=T_{\natural}$. 
 D'autre part, puisque les groupes $G'_{\epsilon,SC}$ et $\bar{H}_{SC}$ sont en situation d'endoscopie non standard, on peut fixer une paire de Borel $(B_{\flat},T_{\flat})$ de $G'_{\epsilon}$ qui correspond \`a $(B_{\bar{H}},T_{\bar{H}})$. On peut fixer $u_{v}\in G'_{\epsilon}$ tel que $(B_{\flat},T_{\flat})=ad_{u_{v}}(B^*_{\flat},T_{\epsilon})$. Posons $(B',T')=ad_{u_{v}}(B_{\epsilon},T_{\epsilon})$. On a $B'\cap G'_{\epsilon}=B_{\flat}$ et $T'=T_{\flat}$.  De nouveau, les propri\'et\'es d'\'equivariance du diagramme 5.2(1) entra\^{\i}nent que
 
 (5) le sextuplet $(\epsilon,B',T',B,T,\eta_{v})$ est un diagramme.

  \bigskip
  
  \subsection{Les places hors de $V$}
  Soit $v$ une place de $F$, provisoirement quelconque. On a d\'efini en [I] 2.7 le groupe $G_{\sharp}=G/Z(G)^{\theta}$.  On a une suite de projections $G\to G_{\sharp}\to G_{AD}$. On note $\underline{G}_{\sharp}(F_{v})$ le sous-groupe des $g\in G(\bar{F}_{v})$ dont l'image dans $G_{\sharp}(\bar{F}_{v})$ appartient \`a $G_{\sharp}(F_{v})$. Pour $g\in \underline{G}_{\sharp}(F_{v})$, les automorphismes $ad_{g}$ de $G$ et de $\tilde{G}$ sont tous deux d\'efinis sur $F$. Pour $(\eta_{v},r_{v})\in D_{v}$, $\underline{G}_{\sharp}(F_{v})$ est inclus dans ${\cal Y}_{\eta_{v}}$. Il  r\'esulte de 5.4(4) que $\underline{G}_{\sharp}(F_{v})$ agit \`a droite sur $D_{v}$ par $((\eta_{v},r_{v}),g)\mapsto (g^{-1}\eta_{v}g,r_{v}g)$. Cette action \'etend celle de $G(F_{v})$.
  
  Supposons maintenant $v\not\in V$. Le groupe $K_{v}$ d\'etermine un sous-groupe hypersp\'ecial $K_{\sharp,v}$ de $G_{\sharp}(F_{v})$. On note $\underline{K}_{\sharp,v}$ le sous-groupe des \'el\'ements de $\underline{G}_{\sharp}(F_{v})$ dont l'image dans $G_{\sharp}(F_{v})$ appartient \`a $K_{\sharp,v}$.  Utilisons les notations de 1.5. On a
 
 (1) $\underline{K}_{\sharp,v}$ est inclus dans $Z(G)^{\theta}K_{v}^{nr}$; pour $g\in \underline{K}_{\sharp,v}$, l'automorphisme $ad_{g}$ de $\tilde{G}$ conserve $\tilde{K}_{v}$.
 
 Preuve.  La surjectivit\'e de l'application $G\to G_{\sharp}$ se prolonge en la surjectivit\'e de l'application analogue entre les sch\'emas en groupes associ\'es \`a $K_{v}$ et $K_{\sharp,v}$. On voit  ainsi qu'il y a une suite exacte
  $$1\to Z(G)^{\theta}\cap K_{v}^{nr}\to K_{v}^{nr}\to K_{\sharp,v}^{nr}\to 1.$$
  Tout \'el\'ement de $K_{\sharp,v}$ se rel\`eve donc en un \'el\'ement de $K_{v}^{nr}$. La premi\`ere assertion en r\'esulte.  La seconde a \'et\'e vue \`a la fin de la preuve de la proposition 2.1. $\square$ 
  
    On note $D_{v}^{nr}$ l'ensemble des $ (\eta_{v},r_{v})\in D_{v}$ tels qu'il existe $h\in G(F_{v})$ de sorte que  $\eta_{v}\in ad_{h}^{-1}(\tilde{K}_{v})$.  Puisque $v\not\in V$, on a $v\not\in S({\cal X},\tilde{K})$ d'apr\`es l'hypoth\`ese 5.1(1).   Le lemme 1.6 implique alors que 
    
    (2) $D_{v}^{nr}$ est non vide; pour $(\eta_{v},r_{v})\in D_{v}^{nr}$ et $h\in G(F_{v})$ tel que  $\eta_{v}\in ad_{h}^{-1}(\tilde{K}_{v})$, le groupe $G_{\eta_{v}}$ est non ramifi\'e et $G_{\eta_{v}}(F_{v})\cap ad_{h}^{-1}(K_{v})$ est un sous-groupe hypersp\'ecial de $G_{\eta_{v}}(F_{v})$. Plus pr\'ecis\'ement, en notant ${\cal K}_{\eta_{v}}$ le sch\'ema en groupes associ\'e \`a ce sous-groupe hypersp\'ecial, on a ${\cal K}_{\eta_{v}}(\mathfrak{o}_{E})=G_{\eta_{v}}(E)\cap ad_{h^{-1}}({\cal K}_{v}(\mathfrak{o}_{E}))$ pour toute extension non ramifi\'ee $E$ de $F_{v}$.

 Il est clair que $D_{v}^{nr}$ est invariant \`a gauche par $I_{\eta}$ et \`a droite par $G(F_{v})$. Remarquons que, puisque $ad_{g}$ est d\'efini sur $F_{v}$ pour tout $g\in \underline{G}_{\sharp}(F_{v})$, les ensembles $G(F_{v})\underline{K}_{\sharp,v}$ et $\underline{K}_{\sharp,v}G(F_{v})$ sont \'egaux et ce sont des groupes. 
 
 \ass{Lemme}{L'ensemble $D_{v}^{nr}$ est stable  par l'action \`a droite de $\underline{K}_{v,\sharp}$. Il forme une unique double classe modulo l'action \`a gauche de $I_{\eta}$ et \`a droite de $\underline{K}_{\sharp,v}G(F_{v})$.}
 
 Preuve.  Soient $(\eta_{v},r_{v})\in D_{v}^{nr}$ et $k\in \underline{K}_{\sharp,v}$. Posons $\eta'_{v}=k^{-1}\eta_{v}k$, $r'_{v}=r_{v}k$. Soit $h\in G(F_{v})$ tel que $\eta_{v}\in ad_{h}^{-1}(\tilde{K}_{v})$. Posons $h'=k^{-1}hk$. C'est un \'el\'ement de $G(F_{v})$. On a 
 $$\eta'_{v}\in ad_{k^{-1}h^{-1}}(\tilde{K}_{v})=ad_{(h')^{-1}k^{-1}}(\tilde{K}_{v})=ad_{(h')^{-1}}(\tilde{K}_{v}),$$
 la derni\`ere \'egalit\'e r\'esultant de (1). Donc $(\eta'_{v},r'_{v})\in D_{v}^{nr}$, ce qui prouve la premi\`ere assertion. Notons $D_{v}^{nr,0}$ le sous-ensemble des $(\eta_{v},r_{v})\in D_{v}$ tels que $\eta_{v}\in \tilde{K}_{v}$. Par d\'efinition, $D_{v}^{nr}$ est engendr\'e par $D_{v}^{nr,0}$ sous l'action \`a droite de $G(F_{v})$. La deuxi\`eme assertion du lemme r\'esulte de
 
 (3) $D_{v}^{nr,0}$ forme une unique double classe modulo l'action \`a gauche de $G_{\eta}$ et \`a droite de $\underline{K}_{\sharp,v}$.
 
 Fixons un \'el\'ement $(\eta_{v},r_{v})\in D_{v}^{nr,0}$. C'est loisible puisque $D_{v}^{nr}$ n'est pas vide.  Remarquons que  faire agir \`a gauche un \'el\'ement de $G_{\eta}$ sur $(\eta_{v},r_{v})$ revient \`a faire agir \`a droite un \'el\'ement de $G_{\eta_{v}}$. L'assertion \`a prouver revient \`a dire que tout \'el\'ement de $D_{v}^{nr,0}$ s'\'ecrit  $(y^{-1}\eta_{v}y,r_{v}y)$ pour un $y\in G_{\eta_{v}}\underline{K}_{\sharp,v}$. Soit $(\eta'_{v},r'_{v}) \in D_{v}^{nr,0}$. D'apr\`es 5.4(4), il existe $y\in {\cal Y}_{\eta_{v}}$ tel que  $(\eta'_{v},r'_{v}) =(y^{-1}\eta_{v}y,r_{v}y)$.  D'apr\`es le lemme [W1] 5.6(i), la relation $y^{-1}\eta_{v}y\in \tilde{K}_{v}$ entra\^{\i}ne 
 que $y\in G_{\eta_{v}}K_{v}^{nr}$. Ecrivons $y=xk$, avec $x\in G_{\eta_{v}}$ et $k\in K_{v}^{nr}$. On a encore $k\in {\cal Y}_{\eta_{v}}$. 
 Notons $Z(G)_{p'}$ le sous-groupe des \'el\'ements de $Z(G;\bar{F}_{v})$ d'ordre premier \`a $p$. C'est un sous-groupe de $K_{v}^{nr}$ ([W1] 5.5(1)). Le lemme 5.5 de [W1] implique que l'application
$$Z(G)_{p'}^{\theta}\to G_{\eta_{v}}\cap K_{v}^{nr}\backslash I_{\eta_{v}}\cap K_{v}^{nr} $$
est surjective.

Fixons un Frobenius $\phi\in \Gamma_{F_{v}}$.  Puisque  $k\in {\cal Y}_{\eta_{v}}$, on a $k\phi(k)^{-1}\in I_{\eta_{v}}$. Puisque $k\in K_{v}^{nr}$, on a aussi $k\phi(k)^{-1}\in K_{v}^{nr}$. D'apr\`es l'assertion ci-dessus, on peut \'ecrire $k\phi(k)^{-1}=g(\phi)z(\phi)$, avec $g(\phi)\in G_{\eta_{v}}\cap K_{v}^{nr}$ et $z(\phi)\in Z(G)_{p'}^{\theta}$.  Par les relations de cocycle habituelles, on prolonge $g(\phi)$ et $z(\phi)$ en des applications d\'efinies sur $\phi^{{\mathbb Z}}$. Par exemple, pour $n\in {\mathbb N}$, on pose 
$$g(\phi^n)=g(\phi)\phi(g(\phi))...\phi^{n-1}(g(\phi)).$$
Parce que $z(\phi)$ est d'ordre fini, on voit qu'il existe $N\geq1$ tel que $n\mapsto z(\phi^n)$ se factorise par ${\mathbb Z}/N{\mathbb Z}$.  Puisque $k\in K_{v}^{nr}$, il en est de m\^eme de l'application $n\mapsto k\phi^n(k)^{-1}$.  Puisque $g(\phi^n)=k\phi^n(k)^{-1}z(\phi^n)^{-1}$, il en est aussi de m\^eme de $n\mapsto g(\phi^n)$. Alors cette application d\'efinit un cocycle continu et  non ramifi\'e de $\Gamma_{F_{v}}$ dans $G_{\eta_{v}}\cap K_{v}^{nr}$. Un tel cocycle est un cobord. On peut donc fixer $x_{1}\in G_{\eta_{v}}\cap K_{v}^{nr}$ de sorte que $g(\phi)=x_{1}^{-1}\phi(x_{1})$. Posons $k_{1}=x_{1}k$. Alors $k_{1}\phi(k_{1})^{-1}=z(\phi)\in Z(G)^{\theta}$.  Donc $k_{1}\in \underline{G}_{\sharp}(F_{v})$. Puisque de plus $k_{1}\in K_{v}^{nr}$, cela implique $k_{1}\in \underline{K}_{\sharp,v}$. Alors $y=xx_{1}^{-1}k_{1}$ appartient \`a $G_{\eta_{v}}\underline{K}_{\sharp,v}$. Cela ach\`eve la preuve. $\square$

    On a
 
 (4) soit $ (\eta_{v},r_{v})\in D_{v}^{nr}$; alors ${\bf H}_{v}$ est relevante pour $(G_{\eta_{v},SC},\psi_{r_{v}})$.
 
 Il s'agit d'endoscopie non tordue. Puisque $G_{\eta_{v},SC}$ est quasi-d\'eploy\'e, tout se transf\`ere.
 
 \bigskip

\subsection{Une cons\'equence}

  \ass{Corollaire}{Supposons  qu'il existe $v\in V$ tel que $D_{v}^{rel}=\emptyset$. Alors $\underline{A}^{\tilde{G},{\cal E}}(V,{\bf H},\omega)=0$.}
 
 Preuve. Par d\'efinition, $\underline{A}^{\tilde{G},{\cal E}}(V,{\bf H},\omega)$ est une somme sur les \'el\'ements $({\bf G}',\mu',\omega_{\bar{G}'})\in {\cal J}({\bf H})$ des transferts des  distributions  $\underline{SA}^{{\bf G}'}(V,p_{\tilde{G}'}(\mu',\omega_{\bar{G}'}))$.  Fixons un tel triplet $({\bf G}',\mu',\omega_{\bar{G}'})$, auquel est associ\'e un \'el\'ement $\epsilon\in \tilde{G}'(F)$. La distribution $\underline{SA}^{{\bf G}'}(V,p_{\tilde{G}'}(\mu',\omega_{\bar{G}'}))$ est support\'ee par les \'el\'ements de $\tilde{G}'(F_{V})$ dont la partie semi-simple est stablement conjugu\'ee \`a $\epsilon$ (plus exactement, ces distributions vivent sur des donn\'ees auxiliaires, mais peu importe ici). Son transfert est nul s'il existe $v\in V$ tel que la classe de conjugaison stable de $\epsilon$ dans $\tilde{G}'(F_{v})$ ne correspond \`a aucune classe de conjugaison stable dans $\tilde{G}(F_{v})$. D'apr\`es le lemme 5.4, il est donc nul s'il existe $v\in V$ tel que $D_{v}^{rel}=\emptyset$.  $\square$
 
 {\bf Dor\'enavant, on suppose $D_{v}^{rel}\not=\emptyset$ pour tout $v\in V$.} Il r\'esulte de 5.5(3) et (4) qu'alors $D_{v}^{rel}\not=\emptyset$ pour tout $v\in Val(F)$. En cons\'equence
 
 (1) on a l'\'egalit\'e ${\cal J}_{\star}({\bf H})={\cal J}({\bf H})$.
 
 Preuve. Soit $({\bf G}',\mu',\omega_{\bar{G}'})\in {\cal J}_{\star}({\bf H})$. Le lemme 1.3 entra\^{\i}ne que $\tilde{G}'(F)$ est non-vide. Fixons un \'el\'ement $\epsilon$ comme dans ce lemme. Pour $v\in Val(F)$, soit $(\eta_{v},r_{v})\in D_{v}^{rel}$. On a construit en 5.4(5) un diagramme $(\epsilon,B',T',B,T,\eta_{v})$. Donc la donn\'ee locale ${\bf G}'_{v}$ est relevante. Par d\'efinition, ${\bf G}'$ est donc relevante et $({\bf G}',\mu',\omega_{\bar{G}'})$ appartient \`a $ {\cal J}_{\star}({\bf H})$. $\square$
 
 Dor\'enavant, pour $v\in Val(F)$, on notera $d_{v}$ un \'el\'ement de $D_{v}$. Quand on aura besoin de l'\'ecrire comme un couple $(\eta_{v},r_{v})$, on \'ecrira ce couple $(\eta[d_{v}],r[d_{v}])$. On supprimera le torseur $\psi_{r_{v}}$ de la notation: le groupe $G_{\eta[d_{v}]}$ sera toujours consid\'er\'e comme une forme int\'erieure de $\bar{G}$ via ce torseur. On utilisera des notations analogues dans diff\'erentes variantes de notre situation. Par exemple, on pose $D_{V}=\prod_{v\in V}D_{v}$ et $D_{V}^{rel}=\prod_{v\in V}D_{v}^{rel}$. On note $(\eta[d_{V}],r[d_{V}])$ le couple associ\'e \`a un \'el\'ement $d_{V}\in D_{V}$ etc...
 
 On notera $j$ un \'el\'ement de ${\cal J}({\bf H})$. Quand on a besoin de l'\'ecrire comme un triplet $({\bf G}',\mu',\omega_{\bar{G}'})$, on note ce triplet $({\bf G}'_{j},\mu'_{j},\omega_{\bar{G}_{j}'})$ et on affecte tous les objets relatifs \`a ce triplet d'un indice $j$. Par exemple, on fixe un \'el\'ement $\epsilon_{j}\in \tilde{G}'_{j}(F)$ v\'erifiant les conditions du lemme 1.3(ii). 
 
 On a
 
 (2) il existe    un sous-tore maximal $S_{\bar{H}}$ de $\bar{H}$, d\'efini sur $F$, tel que, pour toute place $v\in V$, le localis\'e $S_{\bar{H},v}$ soit elliptique dans $\bar{H}_{v}$ si $v$ est non-archim\'edienne et $S_{\bar{H},v}$ soit fondamental dans $\bar{H}_{v}$ si $v$ est archim\'edienne.
 
 Preuve. Pour toute place $v\in V$, on peut en tout cas fixer un sous-tore maximal $S_{v}$ de $\bar{H}_{v}$ d\'efini sur $F_{v}$ qui poss\`ede cette propri\'et\'e. On fixe un \'el\'ement $X_{v}\in \mathfrak{s}_{v}(F_{v})\cap \bar{\mathfrak{h}}_{reg}(F_{v})$. L'espace $\bar{\mathfrak{h}}(F)$ est dense dans $\bar{\mathfrak{h}}(F_{V})$. On peut donc fixer $X\in \bar{\mathfrak{h}} (F)$ tel que, pour toute place $v\in V$, $X$ soit arbitrairement proche de $X_{v}$.  Si $X$ est assez proche de $X_{v}$, sa classe de conjugaison par $\bar{H}(F_{v})$ coupe $\mathfrak{s}_{v}(F_{v})\cap \bar{\mathfrak{h}}_{reg}(F_{v})$. A fortiori, $X$ est semi-simple et r\'egulier. On note $ S_{\bar{H}}$ son commutant dans $\bar{H}$. C'est un sous-tore maximal d\'efini sur $F$. Pour tout $v\in V$, le localis\'e $S_{\bar{H},v}$ est conjugu\'e \`a $S_{v}$ donc poss\`ede la propri\'et\'e requise. $\square$

On fixe un tel tore, que l'on compl\`ete en une paire de Borel $(B_{\bar{H}},S_{\bar{H}})$ d\'efinie sur $\bar{F}$. Soit $j\in {\cal J}({\bf H})$. Puisque $\bar{H}_{SC}$ et $G'_{j,\epsilon_{j},SC}$ sont en situation d'endoscopie non standard, on peut fixer une paire de Borel  $(B_{\flat,j},S_{j})$ de $G'_{j,\epsilon_{j}}$ qui correspond \`a $(B_{\bar{H}},S_{\bar{H}})$. Puisque $(B_{\flat,j},S_{j})$ et $(B^*_{\flat,j},T_{\epsilon_{j}})$ sont deux paires de Borel de $G'_{j,\epsilon_{j}}$ d\'efinies sur $\bar{F}$, on peut fixer $u\in G'_{j,\epsilon_{j}}(\bar{F})$ tel que $(B_{\flat,j},S_{j})=ad_{u}(B^*_{\flat,j},T_{\epsilon_{j}})$. On pose $(B_{j},S_{j})=ad_{u}(B_{\epsilon_{j}},T_{\epsilon_{j}})$. C'est une paire de Borel de $G'_{j}$ d\'efinie sur $\bar{F}$ (avec $S_{j}$ d\'efini sur $F$) conserv\'ee par $ad_{\epsilon_{j}}$. On a $B_{j}\cap G'_{j,\epsilon_{j}}=B_{\flat,j}$.

Soit $v\in Val(F)$ et $d_{v}\in D_{v}^{rel}$. Si $v\not\in V$, supposons de plus $d_{v}\in D_{v}^{nr}$. Montrons que

(3) il existe une paire de Borel $(B_{\natural}[d_{v}],S_{\natural}[d_{v}])$ de $G_{\eta[d_{v}]}$ qui correspond \`a $(B_{\bar{H}},S_{\bar{H}})$.

Preuve. Le corps de base est ici $F_{v}$. Si $v\not\in V$, l'hypoth\`ese implique que $G_{\eta[d_{v}]}$ est quasi-d\'eploy\'e. L'assertion r\'esulte alors de [K1] corollaire 2.2. Supposons $v\in V$.  Puisque ${\bf H}_{v}$ est relevante pour $G_{\eta[d_{v}],SC}$, on peut en tout cas fixer des paires de Borel $(B'_{\bar{H}},S'_{\bar{H}})$ de $\bar{H}$ et  $(B'_{\natural,sc}[d_{v}],S'_{\natural,sc}[d_{v}])$ de $G_{\eta[d_{v}],SC}$ qui se correspondent. 
Puisque $\bar{G}_{SC}$ est quasi-d\'eploy\'e, on peut,  d'apr\`es le corollaire 2.2 de [K1],  compl\'eter ces deux paires par une paire $(\bar{B}'_{\natural,sc},\bar{S}'_{\natural,sc})$ de $\bar{G}_{SC}$ qui correspond \`a chacune de ces paires.   
On note $M'_{\bar{H}}$ le commutant de $A_{S'_{\bar{H}}}$ dans $\bar{H}$, $\bar{M}'_{sc}$ celui de $A_{\bar{S}'_{\natural,sc}}$ dans $\bar{G}_{SC}$ et $M'_{sc}[d_{v}]$ celui de $A_{S'_{\natural,sc}[d_{v}]}$ dans $G'_{\eta[d_{v}],SC}$. Les trois tores intervenant sont naturellement isomorphes et  les trois groupes $M'_{\bar{H}}$, $\bar{M}'_{sc}$ et $M_{sc}'[d_{v}]$ sont des Levi. On v\'erifie que $M_{sc}'[d_{v}]$ est une forme int\'erieure de $\bar{M}_{sc}'$ et que $M'_{\bar{H}}$ est un groupe endoscopique elliptique de $\bar{M}'_{sc}$, donc aussi de $M'_{sc}[d_{v}]$. Notons $M_{\bar{H}}$ le commutant de $A_{S_{\bar{H}}}$ dans $\bar{H}$. 
 Si $v$ est non archim\'edienne, $S_{\bar{H}}$ est elliptique, donc $M_{\bar{H}}=\bar{H}$, a fortiori $M_{\bar{H}}\supset M'_{\bar{H}}$. Si $v$ est archim\'edienne, $S_{\bar{H}}$ est fondamental. Cela implique que, quitte \`a effectuer une conjugaison, on peut supposer $M_{\bar{H}}\supset M'_{\bar{H}}$. Dans les deux cas, on a donc $A_{S_{\bar{H}}}\subset A_{S'_{\bar{H}}}$.  Le tore $A_{S_{\bar{H}}}$ s'identifie \`a des sous-tores de $A_{\bar{S}'_{\natural,sc}}$ et $A_{S'_{\natural,sc}[d_{v}]}$ dont on note $\bar{M}_{sc}$ et $M_{sc}[d_{v}]$ les commutants dans $\bar{G}_{SC}$ et $G_{\eta[d_{v}],SC}$. De nouveau, $M_{sc}[d_{v}]$ est une forme int\'erieure de $\bar{M}_{sc}$ et  $M_{\bar{H}}$ est un groupe endoscopique elliptique de $\bar{M}'_{sc}$, donc aussi de $M'_{sc}[d_{v}]$. Fixons un sous-groupe de  Borel $B_{0}$ de $M_{\bar{H}}$ contenant $S_{\bar{H}}$. 
Puisque $\bar{M}_{sc}$ est quasi-d\'eploy\'e, il r\'esulte du corollaire 2.2 de [K1] que l'on peut fixer une paire de Borel $(B_{1},\bar{S}_{sc})$ de $\bar{M}_{sc}$ qui correspond \`a $(B_{0},S_{\bar{H}})$. On a alors $A_{\bar{S}_{sc}}\simeq A_{S_{\bar{H}}}=A_{M_{\bar{H}}}\simeq A_{\bar{M}_{sc}}$. Cela implique que $\bar{S}_{sc}$ est elliptique dans $\bar{M}_{sc}$. D'apr\`es [K2] lemme 10.2, il existe donc une paire de Borel $(B_{2},S_{\natural,sc}[d_{v}])$ de $M_{sc}[d_{v}]$ qui correspond \`a $(B_{1},\bar{S}_{sc})$. Elle correspond aussi \`a $(B_{0},S_{\bar{H}})$. Fixons $u\in \bar{M}_{H}$ tel que $(B_{0},S_{\bar{H}})=ad_{u}((B'_{\bar{H}}\cap M_{\bar{H}}),S'_{\bar{H}})$ et posons $B_{\bar{H},0}=ad_{u}(B'_{\bar{H}})$. Fixons $x\in M_{sc}[d_{v}]$ tel que $(B_{2},S_{\natural,sc}[d_{v}])=ad_{x}((B'_{\natural,sc}[d_{v}]\cap M_{sc}[d_{v}]),S'_{\natural,sc}[d_{v}])$ et posons $B_{\natural,sc,2}[d_{v}]=ad_{x}(B'_{\natural,sc}[d_{v}])$. L'isomorphisme $S_{\bar{H}}\simeq S_{\natural,sc}[d_{v}]$ d\'eduit des paires de Borel $(B_{\bar{H},0},S_{\bar{H}})$ de $\bar{H}$ et $(B_{\natural,sc,2}[d_{v}],S_{\natural,sc}[d_{v}])$ de $G_{\eta[d_{v}],SC}$ est le m\^eme que celui d\'eduit des paires de Borel $(B_{0},S_{\bar{H}})$ de $M_{\bar{H}}$ et $(B_{2},S_{\natural,sc}[d_{v}])$ de $M_{sc}[d_{v}]$. Il est donc d\'efini sur $F_{v}$. Donc les paires de Borel $(B_{\bar{H},0},S_{\bar{H}})$ et  $(B_{\natural,sc,2}[d_{v}],S_{\natural,sc}[d_{v}])$ se correspondent. Comme on le sait, on peut remplacer le groupe $B_{\bar{H},0}$ par $B_{\bar{H}}$, quitte \`a remplacer $B_{\natural,sc,2}[d_{v}]$ par un sous-groupe de Borel $B_{\natural,sc}[d_{v}]$ convenable, cf. [I] 1.10. On obtient une paire $(B_{\natural,sc}[d_{v}],S_{\natural,sc}[d_{v}])$ de $G_{\eta[d_{v}],SC}$ qui correspond \`a $(B_{\bar{H}},S_{\bar{H}})$. Elle d\'etermine une paire $(B_{\natural}[d_{v}],S_{\natural}[d_{v}])$ de $G_{\eta[d_{v}]}$ qui correspond \`a $(B_{\bar{H}},S_{\bar{H}})$ au sens de 5.2. $\square$

Fixons une telle paire $(B_{\natural}[d_{v}],S_{\natural}[d_{v}])$. Comme en 5.4(5), on en d\'eduit   une paire de Borel de $G$, not\'ee ici $(B[d_{v}],S[d_{v}])$, qui est conserv\'ee par $ad_{\eta_{v}}$. 
La  preuve de 5.4(5) montre  que

(4) pour tout $j\in {\cal J}({\bf H})$, tout $v\in Val(F)$ et tout $d_{v}\in D_{v}^{rel}$ tel que $d_{v}\in D_{v}^{nr}$ si $v\not\in V$, le sextuplet  $(\epsilon_{j},B_{j},S_{j},B[d_{v}],S[d_{v}],\eta[d_{v}])$ est un diagramme sur $F_{v}$.

 \bigskip
 
 \bigskip

\subsection{Facteurs de transfert}
Soit $j\in {\cal J}({\bf H})$. On fixe des donn\'ees auxiliaires $G'_{j,1}$, $\tilde{G}'_{j,1}$, $C_{j,1}$, $\hat{\xi}_{j,1}$ pour ${\bf G}'_{j}$, non ramifi\'ees hors de $V$. On fixe un \'el\'ement $\epsilon_{j,1}\in \tilde{G}'_{j,1}(F)$ au-dessus de $\epsilon_{j}$. Pour $v\in Val(F)-V$, l'hypoth\`ese que $S(p_{\tilde{G}'_{j}}(\mu'_{j},\omega_{\bar{G}'_{j}}),\tilde{K}'_{j})\subset V$ et le lemme 1.6 assurent qu'il existe $\epsilon_{j,v}\in \tilde{K}'_{j,v}$ qui soit stablement conjugu\'e \`a $\epsilon_{j}$. Fixons $u_{v}\in G'_{j}(\bar{F}_{v})$ tel que $u_{v}^{-1}\epsilon_{j}u_{v}=\epsilon_{j,v}$ et $u_{v}\sigma(u_{v})^{-1}\in G'_{j,\epsilon_{j}}$ pour tout $\sigma\in \Gamma_{F_{v}}$. Posons $\epsilon_{j,1,v}=u_{v}^{-1}\epsilon_{j,1}u_{v}$. C'est un \'el\'ement au-dessus de $\epsilon_{j,v}$ et on v\'erifie qu'il appartient \`a $\tilde{G}'_{j,1}(F_{v})$. On pose $\tilde{K}'_{j,1,v}=K'_{j,1,v}\epsilon_{j,1,v}$. On peut choisir $x_{v}=1$ pour presque tout $v$. La famille $(\tilde{K}'_{j,1,v})_{v\not\in V}$ v\'erifie alors la condition de compatibilit\'e globale habituelle, c'est-\`a-dire que, pour $\delta_{1}\in \tilde{G}'_{j,1}(F)$, on a $\delta\in \tilde{K}'_{j,1,v}$ pour presque tout $v$. De ces choix d'espaces hypersp\'eciaux se d\'eduit un facteur de transfert $\Delta_{j,1,v}$ sur $\tilde{G}'_{j,1}(F_{v})\times \tilde{G}(F_{v})$ pour toute place $v\not\in V$. 
La relation 5.6(4) et le fait que la paire $(B_{j},S_{j})$ est d\'efinie sur $\bar{F}$ montrent que, non seulement  la donn\'ee ${\bf G}'_{j}$ est relevante, mais qu'elle v\'erifie l'hypoth\`ese {\bf Hyp} de [VI] 3.6 qui permet de d\'efinir un facteur de transfert global. De ce facteur global et des facteurs $\Delta_{j,1,v}$ pour $v\not\in V$ se d\'eduit comme dans cette r\'ef\'erence un facteur de transfert $\Delta_{j,1,V}$ sur $\tilde{G}'_{j,1}(F_{V})\times \tilde{G}(F_{V})$.  

Puisque $\bar{G}_{SC}$ est quasi-d\'eploy\'e,  la donn\'ee ${\bf H}$ est relevante pour  ce groupe. Puisqu'il est  aussi simplement connexe, on sait que l'on peut identifier $\bar{{\cal H}}$ au $L$-groupe de $\bar{H}$.  Autrement dit, on peut fixer des donn\'ees auxiliaires $\bar{H}_{1}$, $\bar{C}_{1}$ et $\hat{\bar{\xi}}_{1}$ telles que $\bar{H}_{1}=\bar{H}$ et $\bar{C}_{1}=\{1\}$.  On fixe de telles donn\'ees. Pour toute place $v$ et tout $d_{v}\in D_{v}^{rel}$,
  ces donn\'ees auxiliaires valent pour le groupe $G_{\eta[d_{v}],SC}$. Supposons $v\not\in V$ et $d_{v}\in D_{v}^{nr}$. Dans ce cas, ${\bf H}_{v}$ est une donn\'ee non ramifi\'ee pour $G_{\eta[d_{v}],SC}$. On dispose d'un facteur de transfert canonique pourvu que l'on fixe un sous-groupe hypersp\'ecial de $G_{\eta[d_{v}],SC}(F_{v})$ (dans la situation non tordue, la donn\'ee de ce sous-groupe suffit). Pour cela, on fixe un \'el\'ement $h_{v}\in G(F_{v})$ tel que $h_{v}^{-1}\eta[d_{v}]h_{v}\in \tilde{K}_{v}$. On choisit le sous-groupe image r\'eciproque dans $G_{\eta[d_{v}],SC}(F_{v})$ de $ad_{h_{v}}(K_{v})\cap G_{\eta[d_{v}]}(F_{v})$. Si $v\in V$, il n'y a pas de choix canonique de facteur de transfert, on en fixe un que l'on note $\Delta[d_{v}]$. Par le choix de nos donn\'ees auxiliaires, le facteur $\Delta[d_{v}]$ est quel que soit $v$ une fonction  sur  $\bar{H}(F_{v})\times G_{\eta[d_{v}],SC}(F_{v})$.

  Soient $v\in Val(F)$ et $d_{v}\in D_{v}^{rel}$. Soient  $\bar{Y}_{sc}\in \bar{\mathfrak{h}}_{SC}(F_{v})$, $Z_{2}\in \mathfrak{z}(\bar{H};F_{v})$ et $Z_{1}\in \mathfrak{z}(\bar{G};F_{v})$.  On suppose ces \'el\'ements en position g\'en\'erale.   Posons $\bar{Y}=\bar{Y}_{sc}+Z_{2}$. On suppose que la classe de conjugaison stable de $\bar{Y}$ se transf\`ere en une classe de conjugaison stable dans $\mathfrak{g}_{\eta[d_{v}],SC}(F_{v})$. On fixe un \'el\'ement $X[d_{v}]_{sc}$ dans cette classe. Puisqu'on a fix\'e un torseur int\'erieur entre $\bar{G}$ et $G_{\eta[d_{v}]}$, les centres de ces groupes s'identifient et on peut consid\'erer $Z_{1}$ comme un \'el\'ement de $\mathfrak{z}(G_{\eta[d_{v}]};F_{v})$. On pose $X[d_{v}]=Z_{1}+X[d_{v}]_{sc}\in \mathfrak{g}_{\eta[d_{v}]}(F_{v})$. Soit $j\in {\cal J}({\bf H})$. Puisque $\bar{H}_{SC}$ et $G'_{j,\epsilon_{j},SC}$ sont en situation d'endoscopie non standard, la classe de conjugaison stable de $\bar{Y}_{sc}$ se transf\`ere en une classe de conjugaison stable dans $\mathfrak{g}'_{j,\epsilon_{j},SC}(F_{v})$. On fixe un \'el\'ement $Y_{j,sc}$ dans cette classe. Par le diagramme 5.2(5), $ Z_{1}+Z_{2}$ s'identifie \`a un \'el\'ement $Z_{j}\in \mathfrak{z}(G'_{j,\epsilon_{j}};F_{v})$. On pose $Y_{j}=Y_{j,sc}+Z_{j}$. C'est un \'el\'ement de $\mathfrak{g}'_{j,\epsilon_{j}}(F_{v})$. Supposons les \'el\'ements de d\'epart assez proches de $0$ pour que  toutes les exponentielles qui suivent soient d\'efinies. Les classes de conjugaison stable de $exp(Y_{j})\epsilon_{j}$ et de $exp(X[d_{v}])\eta[d_{v}]$ se correspondent. Soit $Y_{j,1}\in \mathfrak{g}'_{j,1,\epsilon_{j,1}}(F_{v})$ au-dessus de $Y_{j}$. On dispose des deux facteurs de transfert $\Delta_{j,1,v}(exp(Y_{j,1})\epsilon_{j,1},exp(X[d_{v}])\eta[d_{v}])$ et $\Delta[d_{v}](exp(\bar{Y}),exp(X_{sc}[d_{v}])) $.
 
 \ass{Th\'eor\`eme}{(i) Si $v$ est finie, il existe une constante $\delta_{j}[d_{v}]$ telle que, pour des \'el\'ements comme ci-dessus tels que $Y_{j,1}$ soit assez proche de $0$, on ait l'\'egalit\'e
 $$\Delta_{j,1,v}(exp(Y_{j,1})\epsilon_{j,1},exp(X[d_{v}])\eta[d_{v}])=\delta_{j}[d_{v}]\Delta[d_{v}](exp(\bar{Y}),exp(X_{sc}[d_{v}])).$$
 
 (ii) Si $v$ est archim\'edienne,  il existe une constante $\delta_{j}[d_{v}]$ et un \'el\'ement $b_{j}\in X_{*}(Z(\hat{G}'_{j,1,\epsilon_{j,1}})^0)\otimes {\mathbb C}$ tels que, pour des \'el\'ements comme ci-dessus pour lesquels $Y_{j,1}$ est  assez proche de $0$, on ait l'\'egalit\'e
 $$\Delta_{j,1,v}(exp(Y_{j,1})\epsilon_{j,1},exp(X[d_{v}])\eta[d_{v}])=\delta_{j}[d_{v}]\Delta[d_{v}](exp(\bar{Y}),exp(X_{sc}[d_{v}]))exp(<b_{j},Y_{j,1}>).$$
 
 (iii) Si $v\not\in V$, on a $\delta_{j}[d_{v}]=\omega(h_{v})^{-1}$,  o\`u $h_{v}$ est l'\'el\'ement utilis\'e pour d\'efinir le sous-groupe hypersp\'ecial de $G_{\eta[d_{v}],SC}(F_{v})$.}

 Preuve.  L'assertion (i) est le th\'eor\`eme 3.9 de [W1]: les deux facteurs locaux co\"{\i}ncident \`a une constante pr\`es. L'assertion (ii) a \'et\'e vue en [V] 4.1. 
 Supposons  $v\not\in V$.  Dans le cas o\`u $h_{v}=1$, l'assertion (iii) est la proposition 5.9 de [W1]. On se ram\`ene \`a ce cas en introduisant le facteur $\Delta'_{j,1,v}$ normalis\'e \`a l'aide des espaces $ad_{h_{v}}(\tilde{K}_{v})$ et $\tilde{K}'_{j,1,v}$. On conserve le m\^eme facteur $\Delta[d_{v}]$. On a une nouvelle constante $\delta'_{j}[d_{v}]$ dont on vient de dire qu'elle valait $1$: quand on remplace $\tilde{K}_{v}$ par $ad_{h_{v}}(\tilde{K}_{v})$, on remplace en m\^eme temps $h_{v}$ par $1$. Mais on v\'erifie ais\'ement que $\Delta'_{j,1,v}=\omega(h_{v})\Delta_{j,1,v}$. L'assertion (iii) en r\'esulte.
$\square$

  \bigskip
 
 \subsection{D\'ebut du calcul de $\underline{A}^{\tilde{G},{\cal E}}(V,{\bf H},\omega)$}

 Soit $f\in C_{c}^{\infty}(\tilde{G}(F_{V}))$. Notre but est de calculer $I^{\tilde{G}}(\underline{A}^{\tilde{G},{\cal E}}(V,{\bf H},\omega),f)$. 
 
 Soit  $ ({\bf G}' ,\mu',\omega_{\bar{G}'})\in{\cal J}({\bf H})$. On se propose de calculer le terme
 $$ i(\tilde{G},\tilde{G}',\mu',\omega_{\bar{G}'})I^{\tilde{G}}(transfert(\underline{SA}^{{\bf G}'}(V, {\cal X}')),f)$$
 qui intervient dans la d\'efinition de 5.3, o\`u ${\cal X}'=p_{\tilde{G}'}(\mu',\omega_{\bar{G}'})$. On le r\'ecrit imm\'ediatement
  $$i(\tilde{G},\tilde{G}',\mu',\omega_{\bar{G}'})S^{{\bf G}'}(\underline{SA}^{{\bf G}'}(V, {\cal X}'),f^{{\bf G}'}).$$
 
 On utilise les d\'efinitions et notations des paragraphes pr\'ec\'edents associ\'ees \`a l'\'el\'ement $j=({\bf G}' ,\mu',\omega_{\bar{G}'})$, tout en supprimant cet indice $j$ pour simplifier la notation.  Le transfert $f^{{\bf G}'}$ s'identifie \`a un \'el\'ement $f_{1}\in SI_{\lambda_{1}}^{\infty}(\tilde{G}'_{1}(F_{V}))$. Alors $S^{{\bf G}'}(\underline{SA}^{{\bf G}'}(V, {\cal X}'),f^{{\bf G}'})$ est \'egal \`a $S_{\lambda_{1}}^{\tilde{G}'_{1}}(\underline{SA}^{\tilde{G}'_{1}}_{\lambda_{1}}(V,{\cal X}'),f_{1})$.  La distribution $\underline{SA}^{\tilde{G}'_{1}}_{\lambda_{1}}(V,{\cal X}')$ est d\'efinie par une formule similaire \`a  1.11(3).  On fixe une fonction $\phi\in C_{c}^{\infty}(\tilde{G}'_{1}(F_{V}))$ telle que
 $$f_{1}=\int_{C_{1}(F_{V})}\phi^{c}\lambda_{1}(c)\,dc.$$
 On a l'\'egalit\'e
 $$S^{{\bf G}'}(\underline{SA}^{{\bf G}'}(V, {\cal X}'),f^{{\bf G}'})=
 \tau'(C_{1})^{-1}\int_{C_{1}(F)\backslash C_{1}({\mathbb A}_{F})}\sum_{{\cal X}'_{1}\in Fib({\cal X}')}$$
 $$S^{\tilde{G}'_{1}}(\underline{SA}^{\tilde{G}'_{1}}(V,{\cal X}'_{1},c^V\tilde{K}_{1}^{_{'}V}),\phi^{c_{V}})\lambda_{1}(c)\,dc.$$
 Notons ${\cal X}'_{1}$ l'\'el\'ement de ${\bf Stab}(\tilde{G}'_{1}(F))$ param\'etrant la classe de conjugaison stable de $\epsilon_{1}$. On peut remplacer la somme sur $Fib({\cal X}')$ par la somme sur les $\xi{\cal X}'_{1}$ pour $\xi\in C_{1}(F)$, \`a condition de diviser par le nombre d'\'el\'ements du groupe $Fix({\cal X}'_{1})=\{\xi\in C_{1}(F); \xi{\cal X}'_{1}={\cal X}'_{1}\}$. Ensuite, par un calcul fait plusieurs fois, la somme en $\xi\in C_{1}(F)$ et l'int\'egrale en $c\in C_{1}(F)\backslash C_{1}({\mathbb A}_{F})$ se recomposent en une int\'egrale en $c\in C_{1}({\mathbb A}_{F})$. On obtient
 $$S^{{\bf G}'}(\underline{SA}^{{\bf G}'}(V, {\cal X}'),f^{{\bf G}'})=\tau'(C_{1})^{-1}\vert Fix({\cal X}'_{1})\vert ^{-1}\int_{C_{1}({\mathbb A}_{F})}S^{\tilde{G}'_{1}}(\underline{SA}^{\tilde{G}'_{1}}(V,{\cal X}'_{1},c^V\tilde{K}_{1}^{_{'}V}),\phi^{c_{V}})\lambda_{1}(c)\,dc.$$
 Soit $c\in C_{1}({\mathbb A}_{F})$. Par d\'efinition de $\underline{SA}^{\tilde{G}'_{1}}(V,{\cal X}'_{1},c^V\tilde{K}_{1}^{_{'}V})$, ce terme est nul sauf si ${\cal X}'_{1}$ est elliptique et $S({\cal X}'_{1},c\tilde{K}_{1}')\subset V$ (remarquons que cette inclusion ne d\'epend que des $c_{v}\tilde{K}'_{1,v}$ pour $v\not\in V$). Les conditions analogues ${\cal X}'$ elliptique et $S({\cal X}',\tilde{K}')\subset V$ sont satisfaites 
  d'apr\`es l'hypoth\`ese  sur $(\mu',\omega_{\bar{G}'})$. En revenant \`a leurs d\'efinitions, on voit que ces deux s\'eries d'hypoth\`eses sont \'equivalentes, \`a l'exception de la condition (nr3) de 1.6:   en une place $v\not\in V$,  cette condition pour ${\cal X}'$ et $\tilde{K}'_{v}$ n'implique pas la m\^eme condition pour ${\cal X}'_{1}$ et $c_{v}\tilde{K}_{1,v}'$. 
  Mais on a choisi les $\tilde{K}'_{1,v}$ de sorte que, pour $v\not\in V$, (nr3) soit v\'erifi\'ee pour $c_{v}=1$. Il r\'esulte alors de [VI] 2.5(13) qu'elle est satisfaite pour tout $v\not\in V$ si et seulement si $c^V
\in K_{C_{1}}^V$.  La fonction \`a int\'egrer est invariante par ce groupe et l'int\'egrale sur ce groupe dispara\^{\i}t. D'o\`u 
$$ S^{{\bf G}'}(\underline{SA}^{{\bf G}'}(V, {\cal X}'),f^{{\bf G}'})=\tau'(C_{1})^{-1}\vert Fix({\cal X}'_{1})\vert ^{-1}\int_{C_{1}(F_{V})} S^{\tilde{G}'_{1}}(\underline{SA}^{\tilde{G}'_{1}}(V,{\cal X}'_{1},\tilde{K}_{1}^{_{'}V}),\phi^{c_{V}})\lambda_{1}(c_{V})\,dc_{V}$$
$$=\tau'(C_{1})^{-1}\vert Fix({\cal X}'_{1})\vert ^{-1}I^{\tilde{G}'_{1}}(\underline{SA}^{\tilde{G}'_{1}}(V,{\cal X}'_{1},\tilde{K}_{1}^{_{'}V}),f_{1}).$$
Restreignons $f_{1}$ \`a un voisinage de la classe de conjugaison stable de $\epsilon_{1}$. On a d\'efini en [I] 4.8 un homomorphisme de descente d'Harish-Chandra  
$$desc_{\epsilon_{1}}^{st}:SI(\tilde{G}'_{1}(F_{V}))\to  SI(G'_{1,\epsilon_{1}}(F_{V})).$$
On pose $f_{\epsilon_{1}}=desc_{\epsilon_{1}}^{st}(f_{1})$.
Pour  d\'efinir correctement l'homomorphisme de descente, il faut en r\'ealit\'e remplacer l'espace de d\'epart, resp. d'arriv\'ee, par un espace de fonctions \`a support dans un voisinage convenable de $\epsilon_{1}$, resp. de l'\'el\'ement neutre.  Il est commode de le noter comme ci-dessus tout en se rappelant l'incorrection de cette notation. Concr\`etement, cela signifie que les int\'egrales orbitales stables de $f_{\epsilon_{1}}$ n'ont de sens qu'au voisinage de l'\'el\'ement neutre dans $G'_{1,\epsilon_{1}}(F_{V})$. Cela ne  g\^ene pas pour appliquer \`a cette fonction une distribution \`a support unipotent. 
 
  En appliquant la d\'efinition de 3.2, on obtient
$$ S^{{\bf G}'}(\underline{SA}^{{\bf G}'}(V, {\cal X}'),f^{{\bf G}'})=\tau'(C_{1})^{-1}\vert Fix({\cal X}'_{1})\vert ^{-1}\vert \Xi_{\epsilon_{1}}^{\Gamma_{F}}\vert ^{-1}\tau'(G'_{1})$$
$$\tau'(G'_{1,\epsilon_{1}})^{-1}S^{G'_{1,\epsilon_{1}}}(SA_{unip}^{G'_{1,\epsilon_{1}}}(V),f_{\epsilon_{1}}).$$
Introduisons la fonction $f_{\epsilon,sc}\in SI(G_{\epsilon,SC}(F_{V}))$ image de $f_{\epsilon_{1}}$ par l'homomorphisme $\iota_{G_{\epsilon,SC},G_{1,\epsilon_{1}}}$. De la m\^eme fa\c{c}on que ci-dessus, les int\'egrales orbitales stables de $f_{\epsilon,sc}$ n'ont de sens qu'au voisinage de l'\'el\'ement neutre dans $G_{\epsilon,SC}(F_{V})$. 
Puisque $G'_{\epsilon,SC}$ est simplement connexe, on a ${\cal A}_{G'_{\epsilon,SC}}=\{0\}$ donc $\tau'(G'_{\epsilon,SC})=\tau(G'_{\epsilon,SC})$. Ce dernier terme vaut $1$ d'apr\`es le th\'eor\`eme de Lai ([Lab1] th\'eor\`eme 1.2). En utilisant la proposition 4.6, on obtient
$$(1) \qquad  S^{{\bf G}'}(\underline{SA}^{{\bf G}'}(V, {\cal X}'),f^{{\bf G}'})=\tau'(C_{1})^{-1}\vert Fix({\cal X}'_{1})\vert ^{-1}\vert \Xi_{\epsilon_{1}}^{\Gamma_{F}}\vert ^{-1}\tau'(G'_{1})S^{G'_{\epsilon,SC}}(SA_{unip}^{G'_{\epsilon,SC}}(V),f_{\epsilon,sc}).$$
D'apr\`es le lemme 4.2(i), on a
 $$(2) \qquad \tau(C_{1})^{-1}\tau(G'_{1})=\tau(G').$$
Montrons que

(3) $\vert Fix({\cal X}'_{1})\vert \vert \Xi_{\epsilon_{1}}^{\Gamma_{F}}\vert =\vert \Xi_{\epsilon}^{\Gamma_{F}}\vert $.

On rappelle que $\Xi _{\epsilon}=Z_{G'}(\epsilon)/G'_{\epsilon}$. Pour $x\in Z_{G'}(\epsilon)$, l'\'el\'ement $x\epsilon_{1}x^{-1}$ est un \'el\'ement de la fibre de $G'_{1}$ au-dessus de $\epsilon$, donc de la forme $c\epsilon_{1}$ pour $c\in C_{1}$. Cela d\'efinit une application $Z_{G'}(\epsilon)\to C_{1}$. On voit que c'est un homomorphisme. Par ailleurs, le groupe $G'_{1,\epsilon_{1}}$ s'envoie surjectivement sur $G'_{\epsilon}$ et est l'image r\'eciproque de ce groupe dans $G'_{1}$. La suite suivante est donc exacte
$$1\to \Xi_{\epsilon_{1}}\to \Xi_{\epsilon}\to C_{1}.$$
D'o\`u aussi une suite exacte
$$1\to \Xi_{\epsilon_{1}}^{\Gamma_{F}}\to \Xi_{\epsilon}^{\Gamma_{F}}\to C_{1}(F).$$
Il suffit de montrer que l'image du second homomorphisme est exactement $Fix({\cal X}'_{1})$. Soit $x\in Z_{G'}(\epsilon)$ dont l'image dans $\Xi_{\epsilon}$ est fixe par $\Gamma_{F}$. Ecrivons $x\epsilon_{1}x^{-1}=c\epsilon_{1}$. La condition sur $x$ implique que $x\epsilon_{1}x^{-1}$ est stablement conjugu\'e \`a $\epsilon_{1}$. Donc $c\in Fix({\cal X}'_{1})$. Inversement, si $c\in Fix({\cal X}'_{1})$, il existe $u\in G'_{1}$ tel que $u^{-1}\epsilon_{1}u=c\epsilon_{1}$ et $u\sigma(u)^{-1} \in G'_{1,\epsilon_{1}}$ pour tout $\sigma\in \Gamma_{F}$. Soit $x$ l'image de $u$ dans $G'$. La premi\`ere condition sur $u$ entra\^{\i}ne que $x$ appartient \`a $Z_{G'}(\epsilon)$ et la seconde condition entra\^{\i}ne que l'image de $x$ dans $\Xi_{\epsilon}$ est fixe par $\Gamma_{F}$. Cela prouve  (3). 

Montrons que
$$(4) \qquad \vert \Xi_{\epsilon}^{\Gamma_{F}}\vert =\vert Fix^{G'}(\mu',\omega_{\bar{G}'})\vert $$
 
 On utilise la paire de Borel $(B_{\epsilon},T_{\epsilon})$ du lemme 1.3  pour construire $(\mu_{\epsilon},\omega_{\epsilon})$ comme en 1.2. Rappelons que ce couple est \'egal \`a $(\mu',\omega_{\bar{G}'})$ par d\'efinition de $\epsilon$ et $(B_{\epsilon},T_{\epsilon})$. A l'aide de la paire $(B_{\epsilon},T_{\epsilon})$, on identifie $W^{G'}$ au groupe de Weyl de $T_{\epsilon}$ dans $\mathfrak{g}'$. Il y a donc deux actions galoisiennes sur $W^{G'}$: l'action quasi-d\'eploy\'ee  not\'ee $\sigma\mapsto \sigma_{G^{_{'}*}}$ et celle provenant de l'action naturelle sur le normalisateur de $T_{\epsilon}$, que l'on note $\sigma\mapsto \sigma$. 
 Compl\`etons $(B_{\epsilon},T_{\epsilon})$  en une paire de Borel \'epingl\'ee ${\cal E}_{\epsilon}$.    L'action galoisienne $\sigma\mapsto \sigma_{G^{_{'}*}}=u_{{\cal E}_{\epsilon}}(\sigma)\circ \sigma$ conserve ${\cal E}_{\epsilon}$.  Parce que $T_{\epsilon}$ est d\'efini sur $F$, il en r\'esulte que $u_{{\cal E}_{\epsilon}}(\sigma)$ normalise $T_{\epsilon}$. Parce que $B_{\epsilon}\cap G'_{\epsilon}$ est d\'efini sur $F$, la cocha\^{\i}ne $u_{\epsilon}$ de la construction de 1.2 est \`a valeurs dans $T_{\epsilon,sc}$.  On voit alors  que, parce que $(\mu_{\epsilon},\omega_{\epsilon})=(\mu',\omega_{\bar{G}'})$, l'image dans $W^{G'}$ de $u_{{\cal E}_{\epsilon}}(\sigma)^{-1}$  est $\omega_{\bar{G}'}(\sigma)$. Il en r\'esulte que, pour $w\in W^{G'}$ et $\sigma\in \Gamma_{F}$, on a l'\'egalit\'e $\sigma(w)=\omega_{\bar{G}'}(\sigma)\sigma_{G^{_{'}*}}(w)\omega_{\bar{G}'}(\sigma)^{-1}$. 
 
 Notons $N$ l'ensemble des \'el\'ements $n\in G'$ tels que $ad_{n}$ conserve $\epsilon$ et la paire $(B_{\epsilon}\cap G'_{\epsilon},T_{\epsilon})$. Il est inclus dans $Z_{G'}(\epsilon)$ et on voit que cette inclusion se quotiente en un isomorphisme
 $$N/T_{\epsilon}\simeq \Xi_{\epsilon}.$$
 Cet isomorphisme identifie $\Xi_{\epsilon}^{\Gamma_{F}}$ au sous-groupe des \'el\'ements de $N/T_{\epsilon}$ fixes par l'action galoisienne naturelle.
 D'autre part, $N/T_{\epsilon}$ est un sous-ensemble de $W^{G'}$. On voit que c'est l'ensemble des $w\in W^{G'}$ tels que 
 
 - $w(\mu')=\mu'$ (cela traduit la condition  que $ad_{n}$ conserve $\epsilon$);
 
  - $w(\Sigma_{+}(\mu'))=\Sigma_{+}(\mu')$ (cela traduit la condition que $ad_{n}$ conserve la paire $(B_{\epsilon}\cap G'_{\epsilon},T_{\epsilon})$).
  
   D'apr\`es ce que l'on a dit ci-dessus,  l'element    
   $w$ est fixe par l'action galoisienne naturelle si et seulement si
  
 -  $w=\omega_{\bar{G}'}(\sigma)\sigma_{G^{_{'}*}}(w)\omega_{\bar{G}'}(\sigma)^{-1}$ pour tout $\sigma\in \Gamma_{F}$.
  
 Ces trois conditions  caract\'erisent le groupe $Stab^{G'}(\mu',\omega_{\bar{G}'})$, cf. 5.1. D'o\`u (4).
 
En utilisant (2), (3) et (4), la formule (1) se r\'ecrit
$$S^{{\bf G}'}(\underline{SA}^{{\bf G}'}(V, {\cal X}'),f^{{\bf G}'})= \tau'(G')\vert Fix^{G'}(\mu',\omega_{\bar{G}'})\vert ^{-1} S^{G'_{\epsilon,sc}}(SA_{unip}^{G'_{\epsilon,sc}}(V),f_{\epsilon,sc}).$$
Rappelons la formule 
$$i(\tilde{G},\tilde{G}',\mu',\omega_{\bar{G}'})= i(\tilde{G},\tilde{G}')\vert Out({\bf G}')\vert \vert W^{G'}(\mu')\vert \vert Fix^{G'}(\mu',\omega_{\bar{G}'})\vert$$
de 5.1(5) et la d\'efinition
$$i(\tilde{G},\tilde{G}')=\vert Out({\bf G}')\vert ^{-1} \vert det((1-\theta)_{\vert \mathfrak{A}_{G}/\mathfrak{A}_{\tilde{G}}})\vert ^{-1} \tau(G)\tau(G')^{-1} $$
$$\vert \pi_{0}((Z(\hat{G})/Z(\hat{G})\cap \hat{T}^{\hat{\theta},0})^{\Gamma_{F}})\vert ^{-1}\vert \pi_{0}(Z(\hat{G})^{\Gamma_{F},0}\cap \hat{T}^{\hat{\theta},0})\vert $$
de [VI] 5.1. 
Le groupe $W^{G'}(\mu')$ s'identifie \`a $W^{\bar{H}}$. On a l'\'egalit\'e $\tau'(G')=covol(\mathfrak{A}_{G',{\mathbb Z}})^{-1}\tau(G')$. On a suppos\'e que l'isomorphisme $\mathfrak{A}_{\tilde{G}}\simeq \mathfrak{A}_{G'}$ pr\'eservait les mesures. Il transforme le r\'eseau $\mathfrak{A}_{\tilde{G},{\mathbb Z}}=Hom(X^*(G)^{\Gamma_{F},\theta},{\mathbb Z})$ en le r\'eseau $\mathfrak{A}_{G'}$ car, dualement, $Z(\hat{G})^{\Gamma_{F},\hat{\theta},0}=Z(\hat{G}')^{\Gamma_{F}}$. Il en r\'esulte que $covol(\mathfrak{A}_{G',{\mathbb Z}})=covol(\mathfrak{A}_{\tilde{G},{\mathbb Z}})$, avec une d\'efinition \'evidente de ce dernier terme. Posons
$$C(\tilde{G})=\vert det((1-\theta)_{\vert \mathfrak{A}_{G}/\mathfrak{A}_{\tilde{G}}})\vert ^{-1} \tau(G)covol(\mathfrak{A}_{\tilde{G},{\mathbb Z}})^{-1}$$
$$\vert \pi_{0}((Z(\hat{G})/Z(\hat{G})\cap \hat{T}^{\hat{\theta},0})^{\Gamma_{F}})\vert ^{-1}\vert \pi_{0}(Z(\hat{G})^{\Gamma_{F},0}\cap \hat{T}^{\hat{\theta},0})\vert.$$
 On obtient alors l'\'egalit\'e
$$(4)\qquad  i(\tilde{G},\tilde{G}',\mu',\omega_{\bar{G}'})S^{{\bf G}'}(\underline{SA}^{{\bf G}'}(V, {\cal X}'),f^{{\bf G}'})= C(\tilde{G}) \vert W^{\bar{H}}\vert S^{G'_{\epsilon,sc}}(SA_{unip}^{G'_{\epsilon,sc}}(V),f_{\epsilon,sc}).$$

\bigskip
\subsection{Utilisation du th\'eor\`eme [VI] 5.6}
On poursuit notre calcul. Les deux groupes    $\bar{H}_{SC}$ et $G'_{\epsilon,SC}$ sont en situation d'endoscopie non standard. Plus pr\'ecis\'ement, notons $j_{*}:X_{*,{\mathbb Q}}(T^{\bar{H}}_{sc})\to X_{*,{\mathbb Q}}(T_{\epsilon,sc})$ l'isomorphisme d\'eduit du  diagramme 5.2(1). Alors $(\bar{H}_{SC},G'_{\epsilon,SC},j_{*})$ est un triplet endoscopique non standard. L'hypoth\`ese ${\cal X}\not\in {\bf Stab}_{excep}(\tilde{G}(F))$ implique que $N(\bar{H}_{SC},G'_{\epsilon,SC},j_{*})< dim(G_{SC})$, cf. [III] lemme 6.3. Nos hypoth\`eses de r\'ecurrence permettent d'appliquer le th\'eor\`eme [VI] 5.6.  On voit ais\'ement que $V$ v\'erifie les conditions de non-ramification impos\'ees dans cette r\'ef\'erence.  Il n'y a pas d'isomorphisme entre $SI(G'_{\epsilon,SC}(F_{V}))$ et $SI(\bar{H}_{SC}(F_{V}))$ mais il y en a  un par contre  entre $SI(\mathfrak{g}'_{\epsilon,SC}(F_{V}))$ et $SI(\bar{\mathfrak{h}}_{SC}(F_{V}))$. Via l'exponentielle, on d\'eduit de celui-ci un isomorphisme entre deux sous-espaces de $SI(G'_{\epsilon,SC}(F_{V}))$ et $SI(\bar{H}_{SC}(F_{V}))$, \`a savoir les sous-espaces de fonctions \`a support dans des voisinages convenables des \'el\'ements neutres. Comme pour les homomorphismes de descente, il est plus commode de consid\'erer cet isomorphisme comme une correspondance entre $SI(G'_{\epsilon,SC}(F_{V}))$ et $SI(\bar{H}_{SC}(F_{V}))$. 
 Introduisons la fonction $\bar{f}_{sc}\in SI(\bar{H}_{SC}(F_{V}))$  qui correspond ainsi  par endoscopie non standard \`a  $f_{\epsilon,sc}\in SI(G'_{\epsilon,SC}(F_{V}))$.  Ses int\'egrales orbitales stables n'ont de sens qu'au voisinage de l'\'el\'ement neutre de $\bar{H}_{SC}(F_{V})$ mais cela nous suffit.  Le th\'eor\`eme [VI] 5.6 transforme l'expression (4) du paragraphe pr\'ec\'edent en
$$(1) \qquad  i(\tilde{G},\tilde{G}',\mu',\omega_{\bar{G}'})S^{{\bf G}'}(\underline{SA}^{{\bf G}'}(V, {\cal X}'),f^{{\bf G}'})=C(\tilde{G})\vert W^{\bar{H}}\vert S^{\bar{H}_{SC}}(SA_{unip}^{\bar{H}_{SC}}(V),\bar{f}_{sc}).$$
Si $V$ \'etait r\'eduit \`a une seule place,  la fonction $\bar{f}_{sc}$ serait calcul\'ee par la formule [III] 5.2(6). Dans cette r\'ef\'erence, le corps de base \'etait non-archim\'edien. Comme on l'a dit en [V] 4.1, le m\^eme calcul vaut sur un corps de base archim\'edien. Dans ce cas, parce qu'on remonte nos fonctions au rev\^etement simplement connexe $\bar{H}_{SC}$, l'exponentielle  pertubatrice du (ii)  du th\'eor\`eme 5.7 dispara\^{\i}t. Le r\'esultat pour notre ensemble fini $V$ de places s'ensuit, en faisant le produit sur tous les $v\in V$. D\'ecrivons-le. Soit  $d_{V}=(d_{v})_{v\in V}\in D_{V}^{rel}$. Puisqu'on a fait dispara\^{\i}tre les espaces de mesures, on suppose implicitement fix\'ees des mesures sur les groupes  $G_{\eta[d_{V}]}(F_{V})$ et $G_{\eta[d_{V}],SC}(F_{V})$ (le choix fait en 4.1 des mesures de Tamagawa ne vaut pas ici puisque les groupes $G_{\eta[d_{V}]}$ et $G_{\eta[d_{V}],SC}$ ne sont pas d\'efinis sur $F$). Elles se d\'eduisent de mesures sur les alg\`ebres de Lie des groupes en question. On a l'isomorphisme
$$\mathfrak{g}_{\eta[d_{V}]}(F_{V})\simeq \mathfrak{z}(\bar{G};F_{V})\oplus \mathfrak{g}_{\eta[d_{V}],SC}(F_{V}).$$
Or $Z(\bar{G})^0$ est d\'efini sur $F$, on munit donc $\mathfrak{z}(\bar{G};F_{V})$ de la mesure de Tamagawa. On suppose que l'isomorphisme ci-dessus est compatible aux mesures. 
Posons $f[d_{V}]=desc_{\eta[d_{V}]}^{\tilde{G}}(f)$, cf. [I] 4.1 pour la d\'efinition de l'homomorphisme de descente $desc_{\eta[d_{V}]}^{\tilde{G}}$. C'est un \'el\'ement de $I(G_{\eta[d_{V}]}(F_{V}),\omega)$. Posons $f[d_{V}]_{sc}=\iota_{G_{\eta[d_{V}],SC},G_{\eta[d_{V}]}}(f[d_{V}])$. C'est un \'el\'ement de $I(G_{\eta[d_{V}],SC}(F_{V}))$. On a fix\'e le facteur de transfert $\Delta[d_{V}]$ en 5.7 (on supprime l'indice $j$ de cette r\'ef\'erence). Notons $\bar{f}[d_{V}]$ le transfert de $f[d_{V}]_{sc}$ \`a $\bar{H}(F_{V})$. C'est un \'el\'ement de $SI(\bar{H}(F_{V}))$. On note $\bar{f}[d_{V}]_{sc}$ son image par  $\iota_{\bar{H}_{SC},\bar{H}}$. C'est un \'el\'ement de $SI(\bar{H}_{SC}(F_{V}))$. Posons
$$c[d_{V}]=[I_{\eta[d_{V}]}(F_{V}):G_{\eta[d_{V}]}(F_{V})]^{-1}$$
et
$$\delta[d_{V}]=\prod_{v\in V}\delta[d_{v}],$$
avec la notation de 5.7 o\`u on supprime les indices $j$. 
Fixons un ensemble de repr\'esentants $\dot{D}_{V}^{rel}$ dans $D_{V}^{rel}$ de l'ensemble de doubles classes
$$I_{\eta}(\bar{F}_{V})\backslash D_{V}^{rel}/G(F_{V}).$$
La formule [III] 5.2(6) nous dit que
$$\bar{f}_{sc}=\sum_{d_{V}\in \dot{D}_{V}^{rel}}c[d_{V}]\delta[d_{V}]\bar{f}[d_{V}]_{sc}.$$
Plus exactement, les int\'egrales orbitales stables des deux membres co\"{\i}ncident dans un voisinage de l'\'el\'ement neutre de $\bar{H}_{SC}(F_{V})$. 

{\bf Remarque.} Cette formule, qui est issue de [W1] 3.11, n\'ecessite certaines compatibilit\'es dans nos choix de mesures. Pr\'ecis\'ement, les mesures sur $G'_{\epsilon}(F_{V})$ et $G'_{\epsilon,SC}(F_{V})$ se d\'eduisent l'une de l'autre par le choix d'une mesure sur $\mathfrak{z}(G'_{\epsilon};F_{V})$; les mesures sur $G_{\eta[d_{V}]}(F_{V})$ et $G_{\eta[d_{V}],SC}(F_{V})$ se d\'eduisent l'une de l'autre par le choix d'une mesure sur $\mathfrak{z}(G_{\eta[d_{V}]};F_{V})$; les mesures sur $\bar{H}(F_{V})$ et $\bar{H}_{SC}(F_{V})$ se d\'eduisent l'une de l'autre par le choix d'une mesure sur $\mathfrak{z}(\bar{H};F_{V})$. Alors les mesures sur $\mathfrak{z}(G'_{\epsilon};F_{V})$, $\mathfrak{z}(G_{\eta[d_{V}]};F_{V})$ et $\mathfrak{z}(\bar{H};F_{V})$ doivent \^etre compatibles avec l'isomorphisme
$$\mathfrak{z}(G'_{\epsilon};F_{V})\simeq \mathfrak{z}(G_{\eta[d_{V}]};F_{V})\oplus \mathfrak{z}(\bar{H};F_{V})$$
d\'eduit du diagramme 5.2(1). Cette compatibilit\'e est assur\'ee par le choix fait ci-dessus des mesures sur $G_{\eta[d_{V}]}(F_{V})$ et $G_{\eta[d_{V}],SC}(F_{V})$ et par le  choix des mesures de Tamagawa sur les autres groupes, cf. lemme 4.2(ii).

\bigskip
On a donc
$$ S^{\bar{H}_{SC}}(SA_{unip}^{\bar{H}_{SC}}(V),\bar{f}_{sc})=\sum_{d_{V}\in \dot{D}_{V}^{rel}}c[d_{V}]\delta[d_{V}] S^{\bar{H}_{SC}}(SA_{unip}^{\bar{H}_{SC}}(V),\bar{f}[d_{V}]_{sc}).$$
Appliquons la proposition 4.6. Elle se simplifie ici. En effet, puisque ${\bf H}$ est une donn\'ee endoscopique elliptique de $\bar{G}_{SC}$, on a $\mathfrak{A}_{\bar{H}}=\mathfrak{A}_{\bar{G}_{SC}}=\{1\}$, donc $\tau'(\bar{H})=\tau(\bar{H})$. On a de m\^eme $\mathfrak{A}_{\bar{H}_{SC}}=\{1\}$, donc $\tau'(\bar{H}_{SC})=\tau(\bar{H}_{SC})=1$ d'apr\`es le th\'eor\`eme de Lai ([Lab1] th\'eor\`eme 1.2).  On obtient que $\iota^*_{\bar{H}_{SC},\bar{H}}(SA_{unip}^{\bar{H}_{SC}}(V))=\tau(\bar{H})^{-1}SA_{unip}^{\bar{H}}(V)$. Pour tout $d_{V}\in \dot{D}_{V}^{rel}$, on a alors l'\'egalit\'e
$$S^{\bar{H}_{SC}}(SA_{unip}^{\bar{H}_{SC}}(V),\bar{f}[d_{V}]_{sc})=\tau(\bar{H})^{-1}S^{\bar{H}}(SA_{unip}^{\bar{H}}(V),\bar{f}[d_{V}]).$$

La formule (1) se r\'ecrit
$$ i(\tilde{G},\tilde{G}',\mu',\omega_{\bar{G}'})S^{{\bf G}'}(\underline{SA}^{{\bf G}'}(V, {\cal X}'),f^{{\bf G}'})=C(\tilde{G})\vert W^{\bar{H}}\vert \tau(\bar{H})^{-1}$$
$$\sum_{d_{V}\in \dot{D}_{V}^{rel}}c[d_{V}]\delta[d_{V}] S^{\bar{H}}(SA_{unip}^{\bar{H}}(V),\bar{f}[d_{V}]).$$
On a fix\'e le triplet $({\bf G}',\mu',\omega_{\bar{G}'})\in {\cal J}({\bf H})$ au d\'ebut du paragraphe pr\'ec\'edent. Faisons-le varier, en le notant $j$. Dans la formule ci-dessus, seul le terme $\delta[d_{V}]$ en d\'epend, on le note d\'esormais $\delta_{j}[d_{V}]$. En reprenant la d\'efinition de 5.3, la formule ci-dessus conduit \`a l'\'egalit\'e
  $$(2)\qquad I^{\tilde{G}}(\underline{A}^{\tilde{G},{\cal E}}(V,{\bf H},\omega),f)=\vert W^G(\mu)\vert ^{-1}\vert Fix^G(\mu,\omega_{\bar{G}})\vert^{-1} C(\tilde{G})\tau(\bar{H})^{-1}\vert W^{\bar{H}}\vert$$
 $$ \sum_{d_{V}\in \dot{D}_{V}^{rel}} c[d_{V}]S^{\bar{H}}(SA_{unip}^{\bar{H}}(V),\bar{f}[d_{V}])\sum_{j\in {\cal J}({\bf H})}\delta_{j}[d_{V}].$$
 
 \bigskip

\section{Calculs de facteurs de transfert}

\subsection{Rappels cohomologiques}
Rappelons quelques d\'efinitions usuelles. Soient $T_{1}$ et $T_{2}$ deux tores d\'efinis sur $F$ et $\varphi:T_{1}\to T_{2}$ un homomorphisme d\'efini sur $F$. On d\'efinit le groupe de cohomologie 
$$H^{1,0}(F;T_{1}\stackrel{\varphi}{\to}T_{2})=H^{1,0}(\Gamma_{F};T_{1}(\bar{F})\stackrel{\varphi}{\to}T_{2}(\bar{F})).$$
On d\'efinit le groupe
$H^{1,0}({\mathbb A}_{F};T_{1}\stackrel{\varphi}{\to}T_{2})$ comme la limite inductive sur les extensions galoisiennes finies $E$ de $F$ des groupes
$$H^{1,0}(Gal(E/F);T_{1}({\mathbb A}_{E})\stackrel{\varphi}{\to}T_{2}({\mathbb A}_{E})).$$
On peut aussi dire qu'en notant ${\mathbb A}_{\bar{F}}$ la limite inductive des ${\mathbb A}_{E}$,  c'est le groupe 
$$H^{1,0}(\Gamma_{F};T_{1}({\mathbb A}_{\bar{F}})\stackrel{\varphi}{\to}T_{2}({\mathbb A}_{\bar{F}})).$$  Pour toute place $v$, on d\'efinit le groupe
$$H^{1,0}(F_{v};T_{1}\stackrel{\varphi}{\to}T_{2})=H^{1,0}(\Gamma_{F_{v}};T_{1}(F_{v})\stackrel{\varphi}{\to}T_{2}(F_{v})).$$
Pour une place $v$ finie o\`u $T_{1}$ et $T_{2}$ sont non ramifi\'es, on d\'efinit 
$$H^{1,0}(\mathfrak{o}_{v};T_{1}\stackrel{\varphi}{\to}T_{2})=H^{1,0}(\Gamma_{v}^{nr};T_{1}(\mathfrak{o}_{v}^{nr})\stackrel{\varphi}{\to}T_{2}(\mathfrak{o}_{v}^{nr})),$$
o\`u on rappelle que $\Gamma_{v}^{nr}=Gal(F_{v}^{nr}/F_{v})$. 
Il  s'envoie injectivement dans  $H^{1,0}(F_{v};T_{1}\stackrel{\varphi}{\to}T_{2})$.
Le groupe $H^{1,0}({\mathbb A}_{F};T_{1}\stackrel{\varphi}{\to}T_{2})$ est isomorphe au produit restreint des $H^{1,0}(F_{v};T_{1}\stackrel{\varphi}{\to}T_{2})$, la restriction \'etant relative aux sous-groupes $H^{1,0}(\mathfrak{o}_{v};T_{1}\stackrel{\varphi}{\to}T_{2})$ d\'efinis pour presque tout $v$. On d\'efinit le groupe
$H^{1,0}({\mathbb A}_{F}/F;T_{1}\stackrel{\varphi}{\to}T_{2})$ comme la limite inductive comme ci-dessus des groupes 
$$H^{1,0}(Gal(E/F);T_{1}({\mathbb A}_{E})/T_{1}(E)\stackrel{\varphi}{\to}T_{2}({\mathbb A}_{E})/T_{2}(E)).$$
Ou encore comme
$$H^{1,0}(\Gamma_{F};T_{1}({\mathbb A}_{\bar{F}})/T_{1}(\bar{F})\stackrel{\varphi}{\to}T_{2}({\mathbb A}_{\bar{F}})/T_{2}(\bar{F})).$$
Dans l'appendice C de [KS], Kottwitz et Shelstad d\'efinissent une topologie sur ce groupe, qui en fait un groupe localement compact. Ils d\'efinissent un accouplement entre ce groupe et le groupe $H^{1,0}(W_{F};\hat{T}_{2}\stackrel{\hat{\varphi}}{\to}\hat{T}_{1})$. De cet accouplement  se d\'eduit un homomorphisme surjectif
 $$(1) \qquad H^{1,0}(W_{F};\hat{T}_{2}\stackrel{\hat{\varphi}}{\to}\hat{T}_{1})\to Hom_{cont}(H^{1,0}({\mathbb A}_{F}/F;T_{1}\stackrel{\varphi}{\to}T_{2}),{\mathbb C}^{\times}),$$
 o\`u, pour deux groupes topologiques $X$ et $Y$,  $Hom_{cont}(X,Y)$ d\'esigne le groupe des homomorphismes continus de $X$ dans $Y$. Il y a une suite exacte
 $$\hat{T}_{2}^{\Gamma_{F}}\to \hat{T}_{1}^{\Gamma_{F}}\to H^{1,0}(W_{F};\hat{T}_{2}\stackrel{\hat{\varphi}}{\to}\hat{T}_{1}).$$
 D'apr\`es le lemme C.2.C de [KS], on a
 
 (2) le noyau de (1) est l'image de $\hat{T}_{1}^{\Gamma_{F},0}$ par le second homomorphise de la suite ci-dessus. 

Consid\'erons maintenant

- deux autres tores $T'_{1}$ et $T'_{2}$ d\'efinis sur $F$ et un homomorphisme $\varphi':T'_{1}\to T'_{2}$ d\'efini sur $F$;

- deux groupes diagonalisables $Z_{1}$ et $Z_{2}$ d\'efinis sur $F$ et un homomorphisme $\psi:Z_{1}\to Z_{2}$ d\'efini sur $F$;

- des diagrammes commutatifs et \'equivariants pour les actions galoisiennes
$$\begin{array}{ccccccc}Z_{1}&\stackrel{\psi}{\to}&Z_{2}&\quad&Z_{1}&\stackrel{\psi}{\to}&Z_{2}\\ \downarrow&&\downarrow&&\downarrow&&\downarrow\\ T_{1}&\stackrel{\varphi}{\to}&T_{2}&&T'_{1}&\stackrel{\varphi'}{\to}&T'_{2}\\ \end{array}$$

On suppose que ces diagrammes sont des quasi-isomorphismes. Cela signifie qu'il s'en d\'eduit des isomorphismes entre groupes de cohomologie, c'est-\`a-dire entre

- le noyau $ker(\psi)$ et le noyau $ker(\varphi)$, resp. $ker(\varphi')$;

- le conoyau $coker(\psi)$ et le conoyau $coker(\varphi)$, resp. $coker(\varphi')$. 

 On a alors des homomorphismes naturels
 $$\begin{array}{ccccc}&&H^{1,0}(F;Z_{1}\stackrel{\psi}{\to}Z_{2})&&\\ &\swarrow&&\searrow&\\ 
H^{1,0}(F;T_{1}\stackrel{\varphi}{\to}T_{2})&&&&H^{1,0}(F;T'_{1}\stackrel{\varphi'}{\to}T'_{2})\\ \end{array}$$
Les deux fl\`eches descendantes sont des isomorphismes. On en d\'eduit un isomorphisme
$$H^{1,0}(F;T_{1}\stackrel{\varphi}{\to}T_{2})\simeq H^{1,0}(F;T'_{1}\stackrel{\varphi'}{\to}T'_{2}).$$
De m\^eme, pour toute place $v$, on a un isomorphisme
$$H^{1,0}(F_{v};T_{1}\stackrel{\varphi}{\to}T_{2})\simeq H^{1,0}(F_{v};T'_{1}\stackrel{\varphi'}{\to}T'_{2}).$$
On v\'erifie que, pour presque tout $v$, cet isomorphisme identifie $H^{1,0}(\mathfrak{o}_{v};T_{1}\stackrel{\varphi}{\to}T_{2})$ \`a $H^{1,0}(\mathfrak{o}_{v};T'_{1}\stackrel{\varphi'}{\to}T'_{2})$. On en d\'eduit un isomorphisme
$$H^{1,0}({\mathbb A}_{F};T_{1}\stackrel{\varphi}{\to}T_{2})\simeq H^{1,0}({\mathbb A}_{F};T'_{1}\stackrel{\varphi'}{\to}T'_{2}).$$
On a aussi un isomorphisme
$$H^{1,0}({\mathbb A}_{F}/F;T_{1}\stackrel{\varphi}{\to}T_{2})\simeq H^{1,0}({\mathbb A}_{F}/F;T'_{1}\stackrel{\varphi'}{\to}T'_{2}).$$
La preuve est plus d\'elicate mais routini\`ere et on la laisse au lecteur. Tous ces isomorphismes sont "fonctoriels" et compatibles aux suites exactes de cohomologie.

\bigskip 

\subsection{Groupes de cohomologie ab\'elienne}
On sait d\'efinir les groupes de cohomologie ab\'elienne d'un groupe r\'eductif connexe d\'efini sur $F$. Ce sont les groupes de cohomologie d'un complexe de tores. Consid\'erons l'exemple de $G$. Fixons un sous-tore maximal $T$ de $G$ d\'efini sur $F$. On peut prendre pour complexe $T_{sc}\to T$. Ainsi, on d\'efinit $H^0_{ab}(F;G)=H^{1,0}(F;T_{sc}\to T)$, $H^1_{ab}(F;G)=H^{2,1}(F;T_{sc}\to T)$. Les nombres $(i+1,i)$ en exposants indiquent que ces groupes classifient des cocycles qui sont des paires de cocha\^{\i}nes, la premi\`ere \'etant de degr\'e $i+1$ \`a valeurs dans $T_{sc}$, la seconde \'etant de degr\'e $i$ \`a valeurs dans $T$. On d\'efinit de m\^eme les groupes $H^{i}_{ab}(F_{v};G)$ pour $v\in Val(F)$, $H^{i}_{ab}({\mathbb A}_{F};G)$ et $H^{i}_{ab}({\mathbb A}_{F}/F;G)$. 
Le choix du tore $T$ n'importe pas. Plus g\'en\'eralement, introduisons $\underline{le}$ tore $T^*$ de $G$, muni de l'action quasi-d\'eploy\'ee et fixons un cocycle $\omega_{T'}:\Gamma_{F}\to W$. D\'efinissons le tore $T'$ comme \'etant \'egal \`a $T^*$, mais muni de l'action galoisienne $\sigma\mapsto \omega_{T'}(\sigma)\circ \sigma_{G^*}$.  On d\'efinit aussi $T'_{sc}$ comme \'etant \'egal \`a $T^*_{sc}$ muni de l'action pr\'ec\'edente.  Le tore $T'$ n'a pas de raison d'\^etre isomorphe \`a un sous-tore de $G$ mais les consid\'erations du paragraphe pr\'ec\'edent montrent que l'on peut aussi bien d\'efinir les groupes de cohomologie ab\'elienne de $G$ \`a l'aide du complexe $T'_{sc}\to T'$. En effet, les complexes $T_{sc}\to T$ et $T'_{sc}\to T'$ sont tous deux quasi-isomorphes au complexe $Z(G_{SC})\to Z(G)$. 

Soit $v\in Val(F)-V$. On peut choisir un sous-tore maximal $T_{v}$ de $G$ d\'efini   sur $F_{v}$ et non ramifi\'e. On d\'efinit alors
$$H^{i}_{ab}(\mathfrak{o}_{v};G)=H^{i+1,i}(\mathfrak{o}_{v};T_{sc,v}\to T_{v}) .$$
 Cela ne d\'epend pas du choix de $T_{v}$. Rappelons que $H^0_{ab}(F_{v};G)$ s'identifie \`a $G(F_{v})/\pi(G_{SC}(F_{v}))$. On a $H^1_{ab}(\mathfrak{o}_{v};G)=\{0\}$ tandis que $H^0_{ab}(\mathfrak{o}_{v};G)$ est l'image naturelle de $T_{v}(\mathfrak{o}_{v})$ dans $H^0_{ab}(F_{v};G)$ ([KS] lemme C.1.A).  L'assertion 1.5(2) peut se reformuler ainsi

(1) $H^0_{ab}(\mathfrak{o}_{v};G)$ est l'image naturelle de $K_{v}$ dans $H^0_{ab}(F_{v};G)$.

On notera plus simplement $G_{ab}(F_{v})=H^0_{ab}(F_{v};G)$ et $G_{ab}(\mathfrak{o}_{v})=H^0_{ab}(\mathfrak{o}_{v};G)$. 

La d\'efinition des groupes de cohomologie ab\'elienne s'\'etend aux groupes non connexes mais quasi-connexes, cf. [Lab2] 1.6. Il faut dans ce cas utiliser des complexes de tores de longueur $3$. Consid\'erons par exemple un \'el\'ement semi-simple $\gamma\in \tilde{G}(F)$, posons $I_{\gamma}=Z(G)^{\theta}G_{\gamma}$. Fixons un sous-tore maximal $T_{\natural}$ de $G_{\gamma}$, notons $T$ son commutant dans $G$ et $T_{\natural,sc}$ l'image r\'eciproque de $T_{\natural}$ dans $G_{\gamma,SC}$. Alors les groupes de cohomologie ab\'elienne de $I_{\gamma}$ sont d\'efinis \`a l'aide du complexe $T_{\natural,sc}\to T\stackrel{1-\theta}{\to}(1-\theta)(T)$. Par exemple, $H^1_{ab}(F;I_{\gamma})=H^{2,1,0}(F; T_{\natural,sc}\to T\stackrel{1-\theta}{\to}(1-\theta)(T))$. L'homomorphisme de complexes
$$\begin{array}{ccccc}Z(G_{\gamma,SC})&\to &Z(I_{\gamma})&\stackrel{1-\theta}{\to}&(1-\theta)(Z(G))\\ \downarrow&&\downarrow&&\downarrow\\ T_{\natural,sc}&\to& T&\stackrel{1-\theta}{\to}&(1-\theta)(T)\\ \end{array}$$
est un quasi-isomorphisme, c'est-\`a-dire qu'il s'en d\'eduit des
 isomorphismes entre  groupes de cohomologie.Les consid\'erations du paragraphe pr\'ec\'edent s'\'etendent aux complexes de longueur finie quelconque, avec les m\^emes cons\'equences.  Identifions $W^{G_{\gamma}}$ au groupe de Weyl de $T_{\natural}$ dans $G_{\gamma}$ et fixons un cocycle $\omega_{T'}:\Gamma_{F}\to W^{G_{\gamma}}$. D\'efinissons le tore $T'$ comme \'etant \'egal \`a $T$, muni de l'action galoisienne $\sigma\mapsto  \omega_{T'}\circ \sigma$. On d\'efinit de m\^eme les tores $T'_{\natural}$ et $T'_{\natural,sc}$. Alors
on peut aussi bien d\'efinir les groupes de cohomologie ab\'elienne de $I_{\gamma}$ \`a l'aide du complexe de tores $T'_{\natural,sc}\to T'\stackrel{1-\theta}{\to}(1-\theta)(T')$. 

\bigskip

\subsection{Un lemme de densit\'e}
Il y a un homomorphisme naturel de $H^0_{ab}(F;G)$ dans $H^0({\mathbb A}_{F};G)$ dont on note l'image $Im(H^0_{ab}(F;G))$. Cette image est discr\`ete pour la topologie naturelle de  $H^0({\mathbb A}_{F};G)$ ([KS] lemme C.3.A). D'autre part, pour toute place $v$, il y a un homomorphisme naturel $G(F_{v})\to G_{ab}(F_{v})=H^0_{ab}(F_{v};G)$. Il est continu et ouvert. Il est  surjectif si $v$ est finie. 
L'assertion 6.2(1) montre qu'il s'en d\'eduit un homomorphisme $G({\mathbb A}_{F})\to H^0_{ab}(F;G)$. Plus pr\'ecis\'ement, notons $V_{\infty}$ l'ensemble des places archim\'ediennes de $F$. En d\'efinissant de fa\c{c}on \'evidente le groupe $H^0_{ab}({\mathbb A}_{F}^{V_{\infty}};G)$, l'assertion 6.2(1) montre que l'homomorphisme $G({\mathbb A}_{F}^{V_{\infty}})\to H^0_{ab}({\mathbb A}_{F}^{V_{\infty}};G)$ est continu, ouvert et surjectif.

\ass{Lemme}{L'homomorphisme
$$G({\mathbb A}_{F})\to H^0_{ab}({\mathbb A}_{F};G)/Im(H^0_{ab}(F;G))$$
est continu, ouvert et surjectif.}

Preuve. Le fait qu'il soit continu et ouvert  r\'esulte de ce que l'on a dit ci-dessus. D'apr\`es le lemme C.5.A de [KS], la projection de $Im(H^0_{ab}(F;G))$ dans $\prod_{v\in V_{\infty}}H^0_{ab}(F_{v};G)$ est dense. Il revient au m\^eme de dire que l'image de $H^0_{ab}({\mathbb A}_{F}^{V_{\infty}};G)$ dans $H^0_{ab}({\mathbb A}_{F};G)/Im(H^0_{ab}(F;G))$ est dense, ou encore que l'image de $G({\mathbb A}_{F}^{V_{\infty}})$ dans ce quotient est dense. A fortiori, l'homomorphisme de l'\'enonc\'e est d'image dense. Cette image \'etant un sous-groupe ouvert, cela entra\^{\i}ne que cette image est le groupe $H^0_{ab}({\mathbb A}_{F};G)/Im(H^0_{ab}(F;G))$
tout entier. $\square$

\bigskip
 \subsection{Fibres de la descente}
 
 On conserve la donn\'ee ${\bf H}=(\bar{H},\bar{{\cal H}},\bar{s})\in E_{\hat{\bar{T}}_{ad},\star}(\bar{G},V)$ fix\'ee en 5.3, soumise \`a la condition $D_{v}^{rel}\not=\emptyset$ pour tout $v\in V$ pos\'ee en 5.6. 
    On veut  d\'ecrire  l'ensemble ${\cal J}({\bf H})$.   
    
   On a fix\'e en 5.2 des paires de Borel des groupes $\hat{\bar{G}}$ et $\hat{\bar{H}}$. On peut identifier le tore $\hat{\bar{T}}$ de la premi\`ere paire \`a $\hat{T}/(1-\hat{\theta})(\hat{T})$ et celui de la seconde \`a $\hat{\bar{T}}_{ad}=\hat{\bar{T}}/Z(\hat{\bar{G}})$. Les actions galoisiennes sur ces tores sont de la forme $\sigma\mapsto \sigma_{\bar{G}}=\omega_{\bar{G}}(\sigma)\sigma_{G^*}$ et $\sigma\mapsto \sigma_{\bar{H}}=\omega_{\bar{H}}(\sigma)\sigma_{\bar{G}}$, o\`u $\omega_{\bar{G}}$ est un cocycle \`a valeurs dans $W^{\theta}$ et $\omega_{\bar{H}}$ est un cocycle \`a valeurs dans $W^{\bar{G}}$. On a fix\'e en 5.6 une paire de Borel $(B_{\bar{H}},S_{\bar{H}})$ de $\bar{H}$. On peut identifier le tore dual $\hat{S}_{\bar{H}}$ \`a $\hat{\bar{T}}_{ad}$ 
  muni d'une action galoisienne $\sigma_{S}=\omega_{S,\bar{H}}(\sigma)\sigma_{\bar{H}}$, o\`u $\omega_{S,\bar{H}}$ est un cocycle de $\Gamma_{F}$ dans  $W^{\bar{H}}$. On pose
$$\omega_{S}(\sigma)=\omega_{S,\bar{H}}(\sigma)\omega_{\bar{H}}(\sigma)\omega_{\bar{G}}(\sigma).$$
On v\'erifie que $\omega_{S}:\Gamma_{F}\to W^{\theta}$ est un cocycle pour l'action quasi-d\'eploy\'ee $\sigma\mapsto \sigma_{G^*}$ de $\Gamma_{F}$ sur $W^{\theta}$.
Introduisons le tore $\hat{S}$ isomorphe \`a $\hat{T}$, muni de l'action galoisienne
$\sigma_{S}=\omega_{S}(\sigma)\sigma_{G^*}$. Alors $\hat{S}_{\bar{H}}$ s'identifie \`a $(\hat{S}/(1-\hat{\theta})(\hat{S}))/Z(\hat{\bar{G}})$. On fixe des $\chi$-data pour l'ensemble des racines du  tore $\hat{S}^{\hat{\theta},0}$ dans $\hat{G}^{\hat{\theta},0}$, d\'efinies sur $F$. On d\'efinit comme dans le cas local (cf. [I] 2.2) une cocha\^{\i}ne
$$\begin{array}{ccc}W_{F}&\to&\hat{G}^{\hat{\theta}}_{SC}\\ w&\mapsto& \hat{r}_{S}(w)\hat{n}_{G}(\omega_{S}(w))\\ \end{array}$$
D'apr\`es [LS] paragraphe 2.6, c'est un cocycle. Il prend ses valeurs dans le normalisateur de $\hat{T}$ dans $\hat{G}^{\hat{\theta}}_{SC}$. En le poussant en un cocycle \`a valeurs dans $W^{\theta}$, on obtient $\omega_{S}$ (relev\'e en un cocycle d\'efini sur $W_{F}$). 

Consid\'erons un \'el\'ement $({\bf G}',\mu',\omega_{\bar{G}'})\in {\cal J}(\bar{H})$. On \'ecrit ${\bf G}'=(G',{\cal G}',s\hat{\theta})$.    On se rappelle que $\omega_{\bar{G}'}(\sigma)\omega_{G'}(\sigma)=\omega_{\bar{H}}(\sigma)\omega_{\bar{G}}(\sigma)$ et que $W^{\bar{H}}=W^{G'}(\mu')$. On pose $\omega_{S,G'}(\sigma)=\omega_{S,\bar{H}}(\sigma)\omega_{\bar{G}'}(\sigma)$. On v\'erifie que $\omega_{S,G'}$ est un cocycle de $W_{F}$ dans $W^{G'}$ (muni de l'action galoisienne  provenant de $G'$). Le tore $\hat{S}^{\hat{\theta},0}$ s'identifie au sous-tore maximal $\hat{T}'=\hat{T}^{\hat{\theta},0}$ de $\hat{G}'$. Par cette identification, l'action $\sigma\mapsto \sigma_{S}$ co\"{\i}ncide avec $\sigma\mapsto \omega_{S,G'}(\sigma)\sigma_{G'}$.  Comme dans le cas local, des $\chi$-data que l'on a fix\'ees  se d\'eduisent de telles donn\'ees pour l'ensemble des racines du tore $\hat{S}^{\hat{\theta},0}$ dans le groupe $\hat{G}'$. On d\'efinit comme ci-dessus le cocycle
$$\begin{array}{ccc}W_{F}&\to&\hat{G}'_{SC}\\ w&\mapsto& \hat{r}_{S,G'}(w)\hat{n}_{G'}(\omega_{S,G'}(w)).\\ \end{array}$$
Pour $w\in W_{F}$, on fixe un \'el\'ement $g_{w}=(g(w),w)\in {\cal G}'$ tel que $ad_{g_{w}}$ co\"{\i}ncide avec $w_{G'}$ sur $\hat{G}'$. On pose
$$t_{S}(w)= \hat{r}_{S}(w)\hat{n}_{G}(\omega_{S}(w))g(w)^{-1}\hat{n}_{G'}(\omega_{S,G'}(w))^{-1}\hat{r}_{S,G'}(w)^{-1}.$$
C'est une cocha\^{\i}ne \`a valeurs dans $\hat{S}$. Ce n'est pas forc\'ement un cocycle, mais son image $\underline{t}_{S}:W_{F}\to \hat{S}/\hat{S}^{\hat{\theta},0}$ en est un, parce que $g(w)$ est bien d\'etermin\'e modulo $\hat{T}^{\hat{\theta},0}$. Pour la m\^eme raison, ce cocycle ne d\'epend pas du choix de $g_{w}$, ni de celui des \'epinglages n\'ecessaires pour d\'efinir les sections de Springer. Notons $\underline{s}$ l'image de $s$ dans  $\hat{S}_{ad}=\hat{S}/Z(\hat{G})\simeq \hat{T}_{ad}$. On v\'erifie que le couple $(\underline{t}_{S},\underline{s})$ est un \'el\'ement de $Z^{1,0}(W_{F};\hat{S}/\hat{S}^{\hat{\theta},0}\stackrel{1-\hat{\theta}}{\to}\hat{S}_{ad})$. Posons 
$$P=H^{1,0}(W_{F};\hat{S}/\hat{S}^{\hat{\theta},0}\stackrel{1-\hat{\theta}}{\to}\hat{S}_{ad})$$
et notons encore
   $(\underline{t}_{S},\underline{s})$ la classe dans  $P$ du cocycle pr\'ec\'edent.  On a ainsi d\'efini une application
$$ \begin{array}{cccc}{\bf p}:&{\cal J}({\bf H})&\to&  P\\ &j=({\bf G}',\mu',\omega_{\bar{G}'})&\mapsto&{\bf p}(j)= (\underline{t}_{S},\underline{s}).\\ \end{array}$$

L'ensemble ${\cal E}_{\hat{T}}(\tilde{G},\omega,V)$ est un ensemble de repr\'esentants de donn\'ees endoscopiques modulo $\hat{T}$-\'equivalence. Evidemment, l'application ci-dessus peut se d\'efinir sur toutes les donn\'ees et pas seulement sur un ensemble de repr\'esentants. Montrons qu'alors, elle se quotiente par cette $\hat{T}$-\'equivalence. En effet, rempla\c{c}ons la donn\'ee ${\bf G}'$ pr\'ec\'edente par une donn\'ee $\hat{T}$-\'equivalente. Cette nouvelle donn\'ee est de la forme $(G',x{\cal G}'x^{-1}, xs\hat{\theta}(x)^{-1}z)$, avec $x\in \hat{T}$ et $z\in Z(\hat{G})$. L'ensemble $Stab(\tilde{G}'(F))$ ne change pas et le couple $(\mu',\omega_{\bar{G}'})$ est encore un \'el\'ement de cet ensemble. Dans les constructions pr\'ec\'edentes, on peut remplacer le cocycle $w\mapsto \hat{r}_{S,G'}(w)\hat{n}_{G'}(\omega_{S,G'}(w))$ par  $w\mapsto x\hat{r}_{S,G'}(w)\hat{n}_{G'}(\omega_{S,G'}(w))x^{-1}$ et $g_{w}$ par $xg_{w}x^{-1}$, donc $g(w)$ par $xg(w)w_{G}(x)^{-1}$. Cela remplace $t_{S}(w)$ par
$$\hat{r}_{S}(w)\hat{n}_{G}(\omega_{S}(w))w_{G}(x)g(w)^{-1}\hat{n}_{G'}(\omega_{S,G'}(w))^{-1}\hat{r}_{S,G'}(w)^{-1}x^{-1}.$$
On a $\hat{r}_{S}(w)\hat{n}_{G}(\omega_{S}(w))\circ w_{G}=w_{S}$ sur $\hat{T}$, donc le terme  pr\'ec\'edent est $w_{S}(x)t_{S}(w)x^{-1}$. Evidemment, $\underline{s}$ est remplac\'e par $(1-\hat{\theta})(x)\underline{s}$. Mais le couple  form\'e du cocycle $w\mapsto w_{S}(x)\underline{t}_{S}(w)x^{-1}$ et de l'\'el\'ement $(1-\hat{\theta})(x)\underline{s}$ est cohomologue \`a $(\underline{t}_{S},\underline{s})$, ce qui d\'emontre l'assertion.

On a  des homomorphismes naturels
$$  P= H^{1,0}(W_{F};\hat{S}/\hat{S}^{\hat{\theta},0}\stackrel{1-\hat{\theta}}{\to}\hat{S}_{ad})\to H^0(W_{F};\hat{S}_{ad}/(1-\hat{\theta})(\hat{S}_{ad}))\to  H^0(W_{F};\hat{S}_{\bar{H}})=\hat{S}_{\bar{H}}^{\Gamma_{F}}.$$
On note ${\bf p}_{1}$ leur compos\'e.

  Le diagramme commutatif
$$\begin{array}{ccc}&&1\\ &&\downarrow\\ 1&\to&Z(\hat{G})\\ \downarrow&&\downarrow\\ \hat{S}/\hat{S}^{\hat{\theta},0}&\stackrel{1-\hat{\theta}}{\to}&\hat{S}\\ \downarrow&&\downarrow\\ \hat{S}/\hat{S}^{\hat{\theta},0}&\stackrel{1-\hat{\theta}}{\to}&\hat{S}_{ad}\\ \downarrow&&\downarrow\\ 1&&1\\ \end{array}$$
fournit un homomorphisme naturel
$$ {\bf p}'_{2}:P= H^{1,0}(W_{F};\hat{S}/\hat{S}^{\hat{\theta},0}\stackrel{1-\hat{\theta}}{\to}\hat{S}_{ad})\to H^1(W_{F};Z(\hat{G})),$$
puis
$$   {\bf p}_{2}:P\to H^1(W_{F};Z(\hat{G}))/ker^1(W_{F};Z(\hat{G})).$$

Soit $v\in Val(F)-V$, notons $I_{v}$ le groupe d'inertie de $\Gamma_{F_{v}}$. C'est aussi un sous-groupe de $W_{F_{v}}$. Remarquons que l'on n'a pas suppos\'e que le tore $S$ \'etait non ramifi\'e hors de $V$. On a un diagramme
$$\begin{array}{ccc}&&( \hat{S}_{sc}/(1-\hat{\theta})(\hat{S}_{sc}))^{I_{v}}\\ &&\parallel\\ &&H^{1,0}(I_{v};\hat{S}_{sc}/\hat{S}_{sc}^{\hat{\theta}}\stackrel{1-\hat{\theta}}{\to}\hat{S}_{sc})\\ && \downarrow \, \varphi_{v}\\ P=H^{1,0}(W_{F};\hat{S}/\hat{S}^{\hat{\theta},0}\stackrel{1-\hat{\theta}}{\to}\hat{S}_{ad})&\stackrel{res_{I_{v}}}{\to} &H^{1,0}(I_{v};\hat{S}/\hat{S}^{\hat{\theta},0}\stackrel{1-\hat{\theta}}{\to}\hat{S}_{ad}).\\  \end{array}$$
 On a not\'e comme toujours $\hat{S}_{sc}$ l'image r\'eciproque de $\hat{S}$ dans $\hat{G}_{SC}$. Son groupe de points fixes $\hat{S}_{sc}^{\hat{\theta}}$ est connexe et l'homomorphisme $1-\hat{\theta}$ est injectif sur $\hat{S}_{sc}/\hat{S}_{sc}^{\hat{\theta}}$. L'isomorphisme du haut en r\'esulte, par la suite exacte de [KS] p.119. L'homomorphisme $res_{I_{v}}$ est la restriction.

On  note 
 $ P({\bf H})$ l'ensemble des $p\in P$  tels que
 
 -   ${\bf p}_{1}(p)=\bar{s}$;
 
 -  ${\bf p}_{2}(p)={\bf a}$;
 
 -  pour tout $v\not\in V$,   $res_{I_{v}}(p)$ appartient \`a l'image de $\varphi_{v}$.

 \ass{Proposition}{L'application  ${\bf p}$ est injective. Son image est $P({\bf H})$.}
 
  Preuve. On commence par prouver l'injectivit\'e.  Consid\'erons deux \'el\'ements $({\bf G}'_{1}, \mu'_{1},\omega_{\bar{G}'_{1}})$ et  $({\bf G}'_{2}, \mu'_{2},\omega_{\bar{G}'_{2}})$ de ${\cal J}({\bf H})$ ayant m\^eme image par  ${\bf p}$. On affecte les termes attach\'es \`a chacune des donn\'ees d'un indice $1$ ou $2$. Les deux cocycles $(\underline{t}_{S,1},\underline{s}_{1})$ et $(\underline{t}_{S,2},\underline{s}_{2})$ sont cohomologues. On a prouv\'e que  remplacer la donn\'ee ${\bf G}'_{2}$ par une donn\'ee $\hat{T}$-\'equivalente rempla\c{c}ait le cocycle  $(\underline{t}_{S,2},\underline{s}_{2})$ par un cocycle cohomologue. En reprenant la preuve, on voit qu'\`a l'inverse, on peut remplacer ${\bf G}'_{2}$ par une donn\'ee $\hat{T}$-\'equivalente de sorte que les deux cocycles $(\underline{t}_{S,1},\underline{s}_{1})$ et $(\underline{t}_{S,2},\underline{s}_{2})$ soient \'egaux.
   Alors les images de $s_{1}$ et $s_{2}$ dans $\hat{T}_{ad}$ sont \'egales. A \'equivalence pr\`es, on peut supposer $s_{1}=s_{2}$. Alors $\hat{G}'_{1}=\hat{G}'_{2}$. Notons simplement $\hat{G}'$ ce groupe. Pour $\sigma\in \Gamma_{F}$, on a l'\'egalit\'e
 $$\omega_{\bar{G}'_{1}}(\sigma)\omega_{G'_{1}}(\sigma)\sigma_{G^*} =\omega_{\bar{G}'_{2}}(\sigma)\omega_{G'_{2}}(\sigma)\sigma_{G^*}$$
 parce que chacun des deux termes est \'egal \`a $\omega_{\bar{H}}(\sigma)\omega_{\bar{G}}(\sigma)\sigma_{G^*}$.
 Mais, pour $i=1,2$, $\omega_{G'_{i}}(\sigma)\sigma_{G^*}$ conserve l'ensemble des racines positives de $\hat{T}^{\hat{\theta},0}$ dans $\hat{G}'$, tandis que $\omega_{\bar{G}'_{i}}(\sigma)$ appartient \`a $W^{G'}$. Une d\'ecomposition en produits de  termes v\'erifiant ces propri\'et\'es est unique. D'o\`u les \'egalit\'es
 $$(1) \qquad \omega_{\bar{G}'_{1}}=\omega_{\bar{G}'_{2}} \text{ et }\omega_{G'_{1}}=\omega_{G'_{2}}.$$
 La seconde \'egalit\'e signifie que l'\'egalit\'e $\hat{G}'_{1}=\hat{G}'_{2}$ est compatible aux actions galoisiennes. Dualement, on peut supposer $G'_{1}=G'_{2}$. Pour $w\in W_{F}$, les termes qui interviennent dans la construction de $\underline{t}_{S,1}(w)$ et $\underline{t}_{S,2}(w)$ sont \'egaux, \`a l'exception peut-\^etre de $g_{1}(w)$ et $g_{2}(w)$. L'\'egalit\'e $\underline{t}_{S,1}(w)=\underline{t}_{S,2}(w)$ entra\^{\i}ne donc que $g_{2}(w)\in \hat{T}^{\hat{\theta},0}g_{1}(w)$ pour tout $w\in W_{F}$. Donc $g_{w,2} \in \hat{G}'g_{1,w}$. Pour $i=1,2$, ${\cal G}'_{i}$ est engendr\'e par $\hat{G}'$ et les $g_{w,i}$ pour $w\in W_{F}$. Donc ${\cal G}'_{1}={\cal G}'_{2}$.  Cela prouve que, quitte \`a remplacer ${\bf G}'_{2}$ par un \'el\'ement $\hat{T}$-\'equivalent, on a l'\'egalit\'e ${\bf G}'_{1}={\bf G}'_{2}$. Puisque les deux \'el\'ements de d\'epart appartiennent \`a un ensemble de repr\'esentants modulo $\hat{T}$-\'equivalence, ces deux \'el\'ements de d\'epart sont en fait \'egaux. Enfin, on a l'\'egalit\'e $\mu'_{1}=\mu'_{2}=\mu$.  D'o\`u, gr\^ace \`a (1), $({\bf G}'_{1}, \mu'_{1},\omega_{\bar{G}'_{1}})=({\bf G}'_{2}, \mu'_{2},\omega_{\bar{G}'_{2}})$. Cela prouve l'injectivit\'e de  ${\bf p}$.

Montrons maintenant que  ${\bf p}({\cal J}({\bf H}))\subset P({\bf H})$. Soit $j=({\bf G}',\mu',\omega_{\bar{G}'})\in {\cal J}({\bf H})$. On \'ecrit ${\bf G}'=(G',{\cal G}',s\hat{\theta})$.  On lui associe le cocycle $p=(\underline{t}_{S},\underline{s})$. Le fait que ${\bf p}_{1}(p)=\bar{s}$ est imm\'ediat.  

Concr\`etement, l'image $z={\bf p}'_{2}(p)$  se construit ainsi. Le cocycle $(\underline{t}_{S},\underline{s})$ se rel\`eve en la cocha\^{\i}ne $(\underline{t}_{S},s)$.  Le cocycle $z:W_{F}\to Z(\hat{G})$ est d\'efini par 
$$(1-\hat{\theta})(\underline{t}_{S}(w))=w_{S}(s)s^{-1}z(w).$$
Parce que $\hat{r}_{S}(w)\hat{n}_{G}(\omega_{S}(w))$ est fixe par $\hat{\theta}$, on a
$$ (1-\hat{\theta})(\underline{t}_{S}(w))=\hat{\theta}\left(\hat{r}_{S,G'}(w)\hat{n}_{G'}(\omega_{S,G'}(w))g(w)\right)g(w)^{-1}\hat{n}_{G'}(\omega_{S,G'}(w))^{-1}\hat{r}_{S,G'}(w)^{-1}$$
$$=s^{-1}s\hat{\theta}\left(\hat{r}_{S,G'}(w)\hat{n}_{G'}(\omega_{S,G'}(w))\right)s^{-1}s\hat{\theta}(g(w))g(w)^{-1}\hat{n}_{G'}(\omega_{S,G'}(w))^{-1}\hat{r}_{S,G'}(w)^{-1}.$$
Parce que $\hat{r}_{S,G'}(w)\hat{n}_{G'}(\omega_{S,G'}(w))\in \hat{G}'$, ce terme est fixe par $ad_{s}\circ\hat{\theta}$. D'autre part, on a une \'egalit\'e
$$s\hat{\theta}(g(w))=a(w)g(w)w_{G^*}(s),$$
o\`u $a$ est un cocycle \`a valeurs dans $Z(\hat{G})$, d'image ${\bf a}$ dans $H^1(W_{F};Z(\hat{G}))/ker^1(W_{F};Z(\hat{G}))$. On obtient
$$(1-\hat{\theta})(\underline{t}_{S}(w))=a(w)s^{-1}\hat{r}_{S,G'}(w)\hat{n}_{G'}(\omega_{S,G'}(w))g(w)w_{G^*}(s)g(w)^{-1}$$
$$\hat{n}_{G'}(\omega_{S,G'}(w))^{-1}\hat{r}_{S,G'}(w)^{-1}.$$
Mais 
$$ad_{\hat{n}_{G'}(\omega_{S,G'}(w))g(w)}\circ w_{G^*}=w_{S}.$$
D'o\`u
$$(1-\hat{\theta})(\underline{t}_{S}(w))=a(w)s^{-1}w_{S}(s).$$
Alors le cocycle $z$ est \'egal \`a $a$. Il s'ensuit que  ${\bf p}_{2}(p)={\bf a}$.

 Soit $v\in Val(F)-V$. 
  Pour $w\in I_{v}$, on peut prendre $g(w)=1$ puisque ${\bf G}'$ est non ramifi\'e en $v$. Alors les termes intervenant dans la d\'efinition de $t_{S}(w)$ appartiennent tous \`a $\hat{G}_{SC}$. On peut consid\'erer $t_{S}$ comme une cocha\^{\i}ne \`a valeurs dans $\hat{S}_{sc}$. C'est encore un cocycle. On en d\'eduit un cocycle encore not\'e $\underline{t}_{S}$ \`a valeurs dans $\hat{S}_{sc}/\hat{S}_{sc}^{\hat{\theta}}$. On fixe un \'el\'ement $s_{sc}\in \hat{S}_{sc}$ ayant m\^eme image que $s$ dans $\hat{S}_{ad}$. Le m\^eme calcul que ci-dessus montre que le couple $(\underline{t}_{S},s_{sc})$ est un cocycle et d\'efinit un \'el\'ement de $H^{1,0}(I_{v};\hat{S}_{sc}/\hat{S}_{sc}^{\hat{\theta}}\stackrel{1-\hat{\theta}}{\to}\hat{S}_{sc})$. Il est clair que $res_{I_{v}}(\underline{t}_{S},\underline{s})$ est l'image par $\varphi_{v}$ de ce cocycle. Cela ach\`eve de prouver que $p=(\underline{t}_{S},\underline{s})$ appartient \`a $P({\bf H})$.

  R\'eciproquement, soit $p\in P({\bf H})$. 
  On repr\'esente $p$ par un cocycle $(\underline{t}_{S},\underline{s})$. On rel\`eve $\underline{s}$ en un \'el\'ement $s\in \hat{T}$. On pose $\hat{G}'=Z_{\hat{G}}(s\hat{\theta})^0$. On munit ce groupe d'une paire de Borel \'epingl\'ee dont la paire sous-jacente soit $(\hat{B}\cap \hat{G}', \hat{T}^{\hat{\theta},0})$.
 
 Parce que $\hat{S}^{\hat{\theta},0}$ est connexe, l'homomorphisme
 $$H^1(W_{F}; \hat{S})\to H^1(W_{F};\hat{S}/\hat{S}^{\hat{\theta},0})$$
 est surjectif, cf. [Lan] p. 719 (1). On rel\`eve $\underline{t}_{S}$ en un cocycle $t_{S}$ \`a valeurs dans $\hat{S}$. Pour $w\in W_{F}$, posons
 $$u(w)=t_{S}(w)^{-1}\hat{r}_{S}(w)\hat{n}_{G}(\omega_{S}(w)),$$
 puis $u_{w}=(u(w),w)\in {^LG}$. On v\'erifie que $w\mapsto u_{w}$ est un homomorphisme de $W_{F}$ dans $^LG$.  Parce que ${\bf p}_{2}(p)={\bf a}$, il existe un cocycle $a:W_{F}\to Z(\hat{G})$, dont l'image dans $H^1( W_{F};Z(\hat{G}))/ker^1(F,Z(\hat{G}))$ est ${\bf a}$, tel que
 $$(2) \qquad (1-\hat{\theta})(t_{S}(w))=w_{S}(s)s^{-1}a(w)$$
 pour tout $w\in W_{F}$. Montrons que l'on a
 
 (3) $s\hat{\theta}(u(w))w_{G^*}(s)^{-1}=a(w)u(w)$ pour tout $w\in W_{F}$.
 
Parce que  $\hat{r}_{S}(w)\hat{n}_{G}(\omega_{S}(w))$ est fixe par $\hat{\theta}$, on a 
$$s\hat{\theta}(u(w))w_{G^*}(s)^{-1}=s\hat{\theta}(t_{S}(w))^{-1}\hat{r}_{S}(w)\hat{n}_{G}(\omega_{S}(w))w_{G^*}(s)^{-1}.$$
Parce que $w_{S}=ad_{\hat{r}_{S}(w)\hat{n}_{G}(\omega_{S}(w))}\circ w_{G^*}$, on obtient
$$s\hat{\theta}(u(w))w_{G^*}(s)^{-1}=s\hat{\theta}(t_{S}(w))^{-1}w_{S}(s)^{-1}\hat{r}_{S}(w)\hat{n}_{G}(\omega_{S}(w)).$$
En utilisant (2), on obtient
$$s\hat{\theta}(u(w))w_{G^*}(s)^{-1}=a(w)t_{S}(w)^{-1}\hat{r}_{S}(w)\hat{n}_{G}(\omega_{S}(w))=a(w)u(w).$$
Cela prouve (3). 

La relation (3) entra\^{\i}ne ais\'ement que $ad_{u_{w}}$ normalise $Z_{\hat{G}}(s\hat{\theta})$, donc aussi sa composante neutre $\hat{G}'$. Alors l'ensemble $\hat{G}'\{u_{w}; w\in W_{F}\}$ est un groupe. Notons-le ${\cal G}'$. On a la suite exacte
$$1\to \hat{G}'\to {\cal G}'\to W_{F}\to 1$$
qui est scind\'ee par l'homomorphisme $u$. Comme toujours, pour $w\in W_{F}$, on peut fixer un \'el\'ement $g_{w}\in \hat{G}'u_{w}$ tel que $ad_{g_{w}}$ conserve l'\'epinglage fix\'e de $\hat{G}'$. Alors $w\mapsto ad_{g_{w}}$ munit $\hat{G}'$ d'une action galoisienne pr\'eservant l'\'epinglage. On introduit le groupe r\'eductif connexe $G'$ sur $F$, quasi-d\'eploy\'e, tel que $\hat{G}'$, muni de son action galoisienne, soit le groupe dual de $G'$. La relation (3) montre que ${\bf G}'=(G',{\cal G}',s\hat{\theta})$ est une donn\'ee endoscopique pour $(G,\tilde{G},{\bf a})$. 

Montrons que

(4) cette donn\'ee endoscopique est non ramifi\'ee hors de $V$.

Soit $v\not\in V$. On sait  $res_{I_{v}}(p)$ appartient \`a l'image de $\varphi_{v}$.  Fixons un cocycle  $(\underline{t}'_{S,sc},s'_{sc})\in H^{1,0}(I_{v};\hat{S}_{sc}/\hat{S}_{sc}^{\hat{\theta}}\stackrel{1-\hat{\theta}}{\to}\hat{S}_{sc})$ tel que $res_{I_{v}}(p)=\varphi_{v}(\underline{t}'_{S,sc},s'_{sc})$. Cela signifie qu'il existe $x\in \hat{S}$ tel que l'on ait $s\in Z(\hat{G})xs'_{sc}$ et $\underline{t}_{S}(w)=\underline{t}'_{S,sc}(w)w_{S}(x)x^{-1}$ pour tout $w\in I_{v}$ (en identifiant un \'el\'ement de $G_{SC}$ \`a son image dans $G$). On \'ecrit $x=zx_{sc}$, avec $x_{sc}\in \hat{S}_{sc}$ et $z\in Z(\hat{G})$. On remplace $(\underline{t}'_{S,sc},s'_{sc})$ par le cocycle cohomologue $(\underline{t}_{S,sc},s_{sc})$ d\'efini par $\underline{t}_{S,sc}(w)=\underline{t}'_{S,sc}(w)w_{S}(x_{sc})x_{sc}^{-1}$ et $s_{sc}=s'_{sc}x_{sc}$. Les relations deviennent  $s\in Z(\hat{G})s_{sc}$ et $\underline{t}_{S}(w)=\underline{t}_{S,sc}(w)w_{S}(z)z^{-1}$ pour $\sigma\in I_{v}$. Mais toutes les actions galoisiennes co\"{\i}ncident sur $Z(\hat{G})$. Elles y sont non ramifi\'ees en $v$ puisque $v\not\in V_{ram}$. Donc $w_{S}(z)=z$ pour $w\in I_{v}$ et on a simplement $\underline{t}_{S}(w)=\underline{t}_{S,sc}(w)$ pour $w\in I_{v}$. Parce que $\hat{S}_{sc}^{\hat{\theta}}$ est connexe, le r\'esultat de [Lan] cit\'e ci-dessus permet de relever $\underline{t}_{S,sc}$ en un cocycle $t_{S,sc}$ \`a valeurs dans $\hat{S}_{sc}$. En se rappelant que le terme  $\hat{r}_{S}(w)\hat{n}_{G}(\omega_{S}(w))$ est naturellement un \'el\'ement de $\hat{G}_{SC}$, on d\'efinit $u_{sc}(w)=t_{S,sc}(w)^{-1}\hat{r}_{S}(w)\hat{n}_{G}(\omega_{S}(w))\in \hat{G}_{SC}$. Parce que $(t_{S,sc},s_{sc})$ est un cocycle, le m\^eme calcul qu'en (3) montre que
$$s_{sc}\hat{\theta}(u_{sc}(w))w_{G^*}(s_{sc})^{-1}=u_{sc}(w)$$
pour $w\in I_{v}$. On a en fait $w_{G^*}(s_{sc})=s_{sc}$ puisque l'action $w\mapsto w_{G^*}$ est non ramifi\'ee. Donc $u(w)$ appartient au groupe des points fixes de l'automorphisme $ad_{s_{sc}}\circ\hat{\theta}$ de $\hat{G}_{SC}$. Ce dernier groupe \'etant simplement connexe, le groupe de points fixes est connexe. Son image dans $\hat{G}$ est la composante neutre du groupe des points fixes de $ad_{s}\circ\hat{\theta}$, c'est-\`a-dire $\hat{G}'$. Il r\'esulte des d\'efinitions que $u(w)\in \hat{S}^{\hat{\theta},0}u_{sc}(w)$ pour $w\in I_{v}$. Donc $u(w)\in \hat{G}'$ pour $w\in I_{v}$. Alors $(1,w)=u(w)^{-1}u_{w}\in {\cal G}'$, ce qui est la condition pour que la donn\'ee ${\bf G}'$ soit non ramifi\'ee en $v$. Cela prouve (4).

Par construction, l'\'el\'ement $u(w)$ normalise $\hat{T}$. La relation (3) entra\^{\i}ne que son image dans $W$ est fixe par $\hat{\theta}$. Il en r\'esulte que l'\'el\'ement $g(w)$ d\'efini par $g_{w}=(g(w),w)$ a les m\^emes propri\'et\'es. On note $\omega_{G'}(w)$ son image dans $W^{\theta}$ et on note $\omega_{S,G'}(w)$ l'\'el\'ement de $W^{\theta}$ tel que l'image de $u(w)$ dans $W$ soit $\omega_{S,G'}(w) \omega_{G'}(w)$. Parce que $g_{w}\in \hat{G}'u_{w}$, on a en fait $\omega_{S,G'}(w)\in W^{G'}$. Il est clair que les applications $\omega_{G'}$ et $\omega_{S,G'}$, qui sont d\'efinies sur $W_{F}$, se factorisent en des applications continues sur $\Gamma_{F}$. La d\'efinition de $u_{w}$ entra\^{\i}ne l'\'egalit\'e
 $$\omega_{S,G'}(\sigma)\omega_{G'}(\sigma)=\omega_{S}(\sigma)=\omega_{S_{\bar{H}}}(\sigma)\omega_{\bar{H}}(\sigma)\omega_{\bar{G}}(\sigma)$$
 pour tout $\sigma\in \Gamma_{F}$. Comme en 1.7, l'\'el\'ement $\mu$ s'identifie \`a un \'el\'ement $\mu'\in T^{_{'}*}\times_{{\cal Z}(G')}{\cal Z}(\tilde{G}')$, o\`u $T^{_{'}*}$ est $\underline{le}$ tore de $G'$. Les calculs de syst\`emes de racines de ce paragraphe et l'hypoth\`ese ${\bf p}_{1}(p)=\bar{s}$ entra\^{\i}nent que le groupe $W^{\bar{H}}$ s'identifie au groupe $W^{G'}(\mu')$. Le terme $\omega_{S_{\bar{H}}}(\sigma)$ appartient \`a ce groupe. Posons 
 $$\omega_{\bar{G}'}(\sigma)=\omega_{S_{\bar{H}}}(\sigma)^{-1}\omega_{S,G'}(\sigma).$$
 On a $\omega_{\bar{G}'}(\sigma)\in W^{G'}$ et l'\'egalit\'e pr\'ec\'edente se r\'ecrit
 $$(5) \qquad \omega_{\bar{G}'}(\sigma)\omega_{G'}(\sigma)=\omega_{\bar{H}}(\sigma)\omega_{\bar{G}}(\sigma).$$
 D'o\`u aussi
 $$\omega_{\bar{G}'}(\sigma)\sigma_{G'}=\omega_{\bar{H}}(\sigma)\omega_{\bar{G}}(\sigma)\sigma_{G^*}.$$
 Le membre de droite fixe $\mu$ donc celui de gauche fixe $\mu'$. Le membre de droite conserve l'ensemble des racines positives de $\bar{H}$, donc celui de gauche conserve l'ensemble $\Sigma^{G'}_{+}(\mu')$. Donc $(\mu',\omega_{\bar{G}'})$ appartient \`a $Stab(\tilde{G}'(F))$. La relation (5) entra\^{\i}ne que cet \'el\'ement s'envoie sur $(\mu,\omega_{\bar{G}})$ par l'application de $Stab(\tilde{G}'(F))$ dans $Stab(\tilde{G}(F))$. D'apr\`es 1.7(3), $(\mu',\omega_{\bar{G}'})$ v\'erifie pour $v\not\in V$ les conditions  (nr1), (nr2) et (nr3) de 1.6, puisque $(\mu,\omega_{\bar{G}})$ les v\'erifie. La condition (nr4) r\'esulte de (5): pour $v\not\in V$ et $\sigma\in I_{v}$, on a $\omega_{G'}(\sigma)=1$ puisque ${\bf G}'$ est non ramifi\'e, $\omega_{\bar{H}}(\sigma)=1$ puisque ${\bf H}$ est non ramifi\'e, $\omega_{\bar{G}}(\sigma)=1$ puisque $v\not\in S({\cal X},\tilde{K})$ (on a impos\'e $S({\cal X},\tilde{K})\subset V$, cf. 5.1); d'o\`u $\omega_{\bar{G}'}(\sigma)=1$. Donc $S(p_{\tilde{G}'}(\mu',\omega_{\bar{G}'}),\tilde{K}')\subset V$.

En inversant la preuve de 5.2(3), on voit que les hypoth\`eses d'ellipticit\'e de ${\bf H}$ et de $(\mu,\omega_{\bar{G}})$ entra\^{\i}nent que ${\bf G}'$ est elliptique et que $(\mu',\omega_{\bar{G}'})$ l'est aussi.  En utilisant l'hypoth\`ese ${\bf p}_{1}(p)=\bar{s}$ et la relation (5), on voit  que les donn\'ees $({\bf G}',\mu',\omega_{\bar{G}'})$ s'envoient sur ${\bf H}$ par la construction du paragraphe 5.2. Evidemment, la donn\'ee ${\bf G}'$ n'a pas de raison d'appartenir \`a l'ensemble de repr\'esentants des classes de $\hat{T}$-\'equivalence que l'on a fix\'e, mais on peut la remplacer par l'\'el\'ement de cet ensemble qui lui est $\hat{T}$-\'equivalent. Tout cela d\'emontre que 
$({\bf G}',\mu',\omega_{\bar{G}'})$ appartient \`a l'ensemble ${\cal J}_{\star}({\bf H})$ introduit en 5.4. On rappelle que cet ensemble est d\'efini de fa\c{c}on analogue \`a ${\cal J}({\bf H})$, sauf que l'on supprime la condition que ${\bf G}'$ est relevante. Mais on a vu en 5.6(1) que, sous l'hypoth\`ese pos\'ee sur ${\bf H}$, cette condition de relevance \'etait automatique. Donc $({\bf G}',\mu',\omega_{\bar{G}'})\in {\cal J}({\bf H})$.  Il r\'esulte des constructions que l'image de cet \'el\'ement par ${\bf p}$ est l'\'el\'ement $p\in P({\bf H})$  dont on est parti. Cela ach\`eve la preuve. $\square$

\bigskip

\subsection{Dualit\'es}
On a d\'ej\`a introduit le groupe
$$P=H^{1,0}(W_{F};\hat{S}/\hat{S}^{\hat{\theta},0}\stackrel{1-\hat{\theta}}{\to}\hat{S}_{ad}).$$
 Introduisons le tore $S$ sur $F$   \'egal \`a $T^*$ muni de l'action galoisienne $\sigma\mapsto \omega_{S}(\sigma)\sigma_{G^*}$. Son dual est $\hat{S}$. Introduisons le groupe
$$Q=H^{1,0}({\mathbb A}_{F}/F;S_{sc}\stackrel{1-\theta}{\to}(1-\theta)(S)).$$
  Comme on l'a dit en 6.1, Kottwitz et Shelstad d\'efinissent  une topologie sur $Q$, pour laquelle ce groupe est localement compact, ainsi qu'un accouplement entre $P$ et $Q$. On a

(1)  l'accouplement entre  $P$ et $Q$ identifie $P$ au groupe $Hom_{cont}(Q,{\mathbb C}^{\times})$ des homomorphismes continus de $Q$ dans ${\mathbb C}^{\times}$.

D'apr\`es  6.1(2), il suffit de prouver que l'homomorphisme 
$$\hat{S}^{\Gamma_{F},0}\stackrel{1-\hat{\theta}}{\to}\hat{S}_{ad}^{\Gamma_{F},0}$$
est surjectif.   On a une projection $\hat{S}/(1-\hat{\theta})(\hat{S})\to \hat{S}_{\bar{H}}$ de noyau $Z(\hat{\bar{G}})$. Elle  est \'equivariante pour les actions galoisiennes. Le tore $S_{\bar{H}}$ est elliptique dans $\bar{H}_{v}$ pour toute place non-archim\'edienne $v\in V$. Cet ensemble de places est non vide puisque $V$ contient $V_{ram}$, lequel contient les places de caract\'eristique r\'esiduelle $2$, $3$ et $5$. Donc $S_{\bar{H}}$ est elliptique dans $\bar{H}$. De plus $\bar{H}$ est une donn\'ee elliptique pour $\bar{G}_{SC}$. Il en r\'esulte que $\hat{S}_{\bar{H}}^{\Gamma_{F},0}=\{1\}$. Donc $(\hat{S}/(1-\hat{\theta})(\hat{S}))^{\Gamma_{F},0}\subset Z(\hat{\bar{G}})^{\Gamma_{F},0}$. Mais $(\mu,\omega_{\bar{G}})$ est elliptique. Donc  $Z(\hat{\bar{G}})^{\Gamma_{F},0}$ est l'image naturelle de $Z(\hat{G})^{\Gamma_{F},0}$. Il en r\'esulte que  $(\hat{S}_{ad}/(1-\hat{\theta})(\hat{S}_{ad}))^{\Gamma_{F},0}=\{1\}$. C'est \'equivalent \`a la surjectivit\'e cherch\'ee. $\square$

  Pour  une place $v\in Val(F)$, on  pose
  $$P_{v}=H^{1,0}(W_{F_{v}};\hat{S}/\hat{S}^{\hat{\theta},0}\stackrel{1-\hat{\theta}}{\to}\hat{S}_{ad})$$
  et
  $$ Q_{v}=H^{1,0}(F_{v};S_{sc}\stackrel{1-\theta}{\to}(1-\theta)(S)).$$
 On note $res_{v}:P\to P_{v}$ l'homomorphisme de restriction et $\iota_{v}:Q_{v}\to Q$ l'homomorphisme naturel. On note
 $${\bf p}_{2,v}:P_{v}\to 
  H^1(W_{F_{v}};Z(\hat{G}))$$
  l'analogue local de l'homomorphisme ${\bf p}'_{2}$.
  
  On a rappel\'e en 5.5 le groupe $G_{\sharp}=G/Z(G)^{\theta}$. On a le diagramme commutatif
  $$\begin{array}{ccc}S_{sc}&\to& S\\ \downarrow&&\downarrow\\ S_{sc}&\to&S/Z(G)^{\theta}\\ \downarrow&&\qquad \downarrow 1-\theta\\ S_{sc}&\stackrel{1-\theta}{\to}&(1-\theta)(S)\\ \end{array}$$
  Comme on l'a dit en 6.2, le groupe $G_{ab}(F_{v}) $ peut se d\'efinir \`a l'aide du complexe $S_{sc}\to S$. De m\^eme, $G_{\sharp,ab}(F_{v}) $ peut se d\'efinir \`a l'aide du complexe $S_{sc}\to S/Z(G)^{\theta}$. On d\'eduit du diagramme ci-dessus un diagramme commutatif
  $$\begin{array}{ccccc}G_{ab}(F_{v})&&\stackrel{\boldsymbol{\zeta}_{v}}{\to}&& Q_{v}\\ &\searrow&&\qquad \nearrow \boldsymbol{\zeta}_{\sharp,v}&\\ &&G_{\sharp,ab}(F_{v})&&\\ \end{array}$$
  
   Enfin,   on note
 $${\bf q}_{1}:H^{1}({\mathbb A}_{F}/F; S_{\bar{H}})\to Q$$
 l'homomorphisme d\'eduit de l'homomorphisme  $ S_{\bar{H}}\to S_{sc}$ dual de $\hat{S}_{ad}\to \hat{S}_{\bar{H}}$.

 \ass{Proposition}{(i) Pour toute place $v$, le noyau de ${\bf p}_{2,v}\circ res_{v}$ est l'annulateur dans $P$ de   $\iota_{v}\circ \boldsymbol{\zeta}_{v}(G_{ab}(F_{v}))$.
 
 (ii) Pour toute place $v\not\in V$, le sous-groupe des $p\in P$ tels que $res_{I_{v}}$ appartienne \`a l'image de $\varphi_{v}$ est l'annulateur dans $P$ de   $\iota_{v}\circ  \boldsymbol{\zeta}_{\sharp,v}(G_{\sharp,ab}(\mathfrak{o}_{v})$.
 
 (iii) Le noyau de ${\bf p}_{1}$ est l'annulateur dans $P$ de l'image de ${\bf q}_{1}$.}
 
 Preuve. L'accouplement entre $H^{1}({\mathbb A}_{F}/F; S_{\bar{H}})$ et $\hat{S}_{\bar{H}}^{\Gamma_{F}}$  identifie le second groupe \`a celui des homomorphismes continus du premier dans ${\mathbb C}^{\times}$. Comme dans la preuve de (1), cela r\'esulte de [KS] lemme C.2.C et de l'ellipticit\'e de $S_{\bar{H}}$. Alors  ${\bf p}_{1}$ s'identifie \`a l'homomorphisme
 $$Hom_{cont}(Q,{\mathbb C}^{\times})\to Hom_{cont}(H^{1}({\mathbb A}_{F}/F; S_{\bar{H}}),{\mathbb C}^{\times})$$
 d\'eduit par dualit\'e de ${\bf q}_{1}$.
  L'assertion (iii) r\'esulte de la propri\'et\'e g\'en\'erale suivante: si $f:X\to Y$ est un homomorphisme continu entre groupes ab\'eliens localement compacts, le noyau  de  l'homomorphisme dual 
  $$f^D:Hom_{cont}(Y,{\mathbb C}^{\times})\to Hom_{cont}(X,{\mathbb C}^{\times})$$ 
  est l'annulateur de l'image de $f$.

 Soit $v$ une place de $F$. On a un diagramme
  $$\begin{array}{ccccc}P&\stackrel{res_{v}}{\to}&P_{v}&\stackrel{{\bf p}_{2,v}}{\to}&H^1(W_{F_{v}};Z(\hat{G}))\\ Q&\stackrel{\iota_{v}}{\leftarrow}&Q_{v}&\stackrel{\boldsymbol{\zeta}_{v}}{\leftarrow}& G_{ab}(F_{v})\\ \end{array}$$
  Il y a des accouplements entre $P$ et $Q$, entre $P_{v}$ et $Q_{v}$ et entre $H^1(W_{F_{v}};Z(\hat{G}))$ et $G_{ab}(F_{v})$. On a d\'ej\`a dit que, par le premier accouplement, $P$ s'identifiait \`a $Hom_{cont}(Q,{\mathbb C}^{\times})$. On sait aussi que, par le dernier,  $H^1(W_{F_{v}};Z(\hat{G}))$ s'identifie \`a $Hom_{cont}(G_{ab}(F_{v}),{\mathbb C}^{\times})$. On v\'erifie que ${\bf p}_{2,v}\circ res_{v}$ s'identifie \`a l'homomorphisme
  $$Hom_{cont}(G_{ab}(F_{v}),{\mathbb C}^{\times})\to Hom_{cont}(Q,{\mathbb C}^{\times})$$
  dual de $\iota_{v}\circ \boldsymbol{\zeta}_{v}$. L'assertion (ii) r\'esulte alors du m\^eme principe g\'en\'eral que ci-dessus.
  
  {\bf Remarque.} En g\'en\'eral, l'accouplement entre $P_{v}$ et $Q_{v}$ a un noyau dans $P_{v}$, \'egal \`a l'image naturelle de  l'homomorphisme $\hat{S}_{ad}^{\Gamma_{F_{v}},0}\to P_{v}$. C'est le conoyau de cet homomorphisme qui s'identifie \`a $Hom_{cont}(Q_{v},{\mathbb C}^{\times})$.
  \bigskip

 Soit $v\not\in V$. Notons $S_{\sharp}=S/Z(G)^{\theta}$. On a d\'ecrit le tore dual $\hat{S}_{\sharp}$ en [I] 2.7. On a une suite exacte
 $$1\to \hat{S}_{sc}/\hat{S}_{sc}^{\hat{\theta}}\stackrel{(\pi,1-\hat{\theta})}{\to} \hat{S}/\hat{S}^{\hat{\theta},0}\times \hat{S}_{sc}\to \hat{S}_{\sharp}\to 1.$$
 Il s'en d\'eduit un diagramme commutatif
 $$\begin{array}{ccc}\hat{S}_{sc}/\hat{S}_{sc}^{\hat{\theta}}&\stackrel{1-\hat{\theta}}{\to}&\hat{S}_{sc}\\ \downarrow&&\downarrow \\ \hat{S}/\hat{S}^{\hat{\theta},0}&\stackrel{1-\hat{\theta}}{\to}&\hat{S}_{ad}\\ \downarrow&&\downarrow\\ \hat{S}_{\sharp}&\to &\hat{S}_{ad}\\ \end{array}$$
 C'est un triangle distingu\'e  dans la cat\'egorie des complexes de tores. On obtient un diagramme
 $$\begin{array}{ccccc}&&&&H^{1,0}(I_{v};\hat{S}_{sc}/\hat{S}_{sc}^{\hat{\theta}}\stackrel{1-\hat{\theta}}{\to} \hat{S}_{sc})\\&&&&\downarrow \varphi_{v}\\ P&\stackrel{res_{v}}{\to}&P_{v}&\stackrel{Res_{I_{v}}}{\to}&H^{1,0}(I_{v};\hat{S}/\hat{S}^{\hat{\theta},0}\stackrel{1-\hat{\theta}}{\to}\hat{S}_{ad})\\ &&\downarrow {\bf p}_{2,\sharp,v}&&\downarrow \\&&  H^{1,0}(W_{F_{v}};\hat{S}_{\sharp}\to \hat{S}_{ad})&\stackrel{Res_{\sharp,I_{v}}}{\to} &  H^{1,0}(I_{v};\hat{S}_{\sharp}\to \hat{S}_{ad})\\ &&\parallel&&\parallel\\ &&H^1(W_{F_{v}};Z(\hat{G}_{\sharp}))&\stackrel{Res_{\sharp,I_{v}}}{\to}&H^1(I_{v};Z(\hat{G}_{\sharp}))\\ \end{array}$$
 Les deux premiers homomorphismes de la colonne de droite forment une suite exacte. D'autre part, on a l'\'egalit\'e $res_{I_{v}}=Res_{I_{v}}\circ res_{v}$ avec la notation ci-dessus. Donc, pour $p\in P$, la condition que  $res_{I_{v}}(p)$ appartienne \`a l'image de $\varphi_{v}$ est \'equivalente \`a ce que ${\bf p}_{2,\sharp,v}\circ res_{v}(p)$ appartienne au noyau de $Res_{\sharp,I_{v}}$. Le groupe $H^1(W_{F_{v}};Z(\hat{G}_{\sharp}))$ s'identifie \`a $Hom_{cont}(G_{\sharp,ab}(F_{v}),{\mathbb C}^{\times})$. D'apr\`es 1.5(4) et 6.2(1), 
 un \'el\'ement $\chi\in H^1(W_{F_{v}};Z(\hat{G}_{\sharp}))$ est annul\'e par $Res_{\sharp,I_{v}}$ si et seulement si $\chi$ annule  $G_{\sharp,ab}(\mathfrak{o}_{v})$. La condition que $res_{I_{v}}(p)$ appartienne \`a l'image de $\varphi_{v}$ \'equivaut donc \`a l'\'egalit\'e $<{\bf p}_{2,\sharp,v}\circ res_{v}(p),k>=1$ pour tout $k\in G_{\sharp,ab}(\mathfrak{o}_{v})$. Mais on a l'\'egalit\'e
  $<{\bf p}_{2,\sharp,v}\circ res_{v}(p),k> =<p,\iota_{v}\circ\boldsymbol{\zeta}_{\sharp,v}(k)>$. La condition ci-dessus \'equivaut donc \`a ce que $p$ 
    annule   $\iota_{v}\circ \boldsymbol{\zeta}_{\sharp,v}(G_{\sharp,ab}(\mathfrak{o}_{v}))$.    Cela prouve (ii). $\square$

 Supposons que $v$ soit une place hors de $V$ en laquelle $S$ est non ramifi\'ee.  On a alors
  
 (2) le groupe $\iota_{v}\circ \boldsymbol{\zeta}_{\sharp, v}(G_{\sharp,ab}(\mathfrak{o}_{v}))$ co\"{\i}ncide avec l'image naturelle dans $Q$ de $((1-\theta)(S))(\mathfrak{o}_{v})$. 
 
 Preuve. On a  une suite exacte
 $$1\to S^{\theta}/Z(G)^{\theta}\to S_{\sharp}\stackrel{1-\theta}{\to}(1-\theta)(S)\to 1.$$
 Il s'en d\'eduit une suite exacte
 $$1\to (S^{\theta}/Z(G)^{\theta})(\mathfrak{o}_{v}^{nr})\to S_{\sharp}(\mathfrak{o}_{v}^{nr})\stackrel{1-\theta}{\to}((1-\theta)(S))(\mathfrak{o}_{v}^{nr})\to 1.$$
   Mais le groupe $S^{\theta}/Z(G)^{\theta}$ est connexe en vertu de l'\'egalit\'e $S^{\theta}=S^{\theta,0}Z(G)^{\theta}$. En prenant les invariants par le groupe de Galois $\Gamma_{v}^{nr}$, le th\'eor\`eme de Lang implique la surjectivit\'e de l'homomorphisme
 $$S_{\sharp}(\mathfrak{o}_{v})\stackrel{1-\theta}{\to} ((1-\theta)(S))(\mathfrak{o}_{v}).$$
On a le diagramme commutatif
 $$\begin{array}{ccc}S_{\sharp}(\mathfrak{o}_{v})&\stackrel{1-\theta}{\to}&((1-\theta)(S))(\mathfrak{o}_{v})\\ \downarrow&&\downarrow\\ G_{\sharp,ab}(F_{v})&\to& Q_{v}.\\ \end{array}$$
D'apr\`es 1.5(2) et 6.2(1),    l'image de $S_{\sharp}(\mathfrak{o}_{v})$ dans $G_{\sharp,ab}(F_{v})$ co\"{\i}ncide avec  $G_{\sharp,ab}(\mathfrak{o}_{v})$.  Cela conclut.  $\square$

\bigskip

 \subsection{Description d'un annulateur}
 
 On note $P^0$ le sous-groupe des $p\in P$ tels que
 
 - ${\bf p}_{1}(p)=0$; 
 
 - pour toute place $v$, ${\bf p}_{2,v}\circ res_{v}(p)=0$;
 
 - pour toute place $v\not\in V$, $res_{I_{v}}(p)$ appartient \`a l'image de $\varphi_{v}$.
 
  Remarquons que l'ensemble $P({\bf H})\subset P$ introduit en 6.4 est soit vide, soit une unique classe modulo ce sous-groupe $P^0$. 
  
  On a d\'ej\`a d\'efini l'homomorphisme ${\bf q}_{1}:H^1({\mathbb A}_{F}/F;S_{\bar{H}})\to Q$. On pose $Q_{1}=H^1({\mathbb A}_{F}/F;S_{\bar{H}})$.
 
Pour toute place $v$, on a d\'efini des homomorphismes $\boldsymbol{\zeta}_{v}:G_{ab}(F_{v})\to Q_{v}$ et $\boldsymbol{\zeta}_{\sharp,v}:G_{\sharp,ab}(F_{v})\to Q_{v}$. Ils se globalisent en des homomorphismes $\boldsymbol{\zeta}:H^0_{ab}({\mathbb A}_{F};G)\to H^{1,0}({\mathbb A}_{F};S_{sc}\stackrel{1-\theta}{\to}(1-\theta)(S))$ et $\boldsymbol{\zeta}_{\sharp}:H^0_{ab}({\mathbb A}_{F};G_{\sharp})\to H^{1,0}({\mathbb A}_{F};S_{sc}\stackrel{1-\theta}{\to}(1-\theta)(S))$. En poussant le premier par l'application naturelle $H^{1,0}({\mathbb A}_{F};S_{sc}\stackrel{1-\theta}{\to}(1-\theta)(S))\to Q$, on obtient un homomorphisme $H^0_{ab}({\mathbb A}_{F};G)\to Q$. Il est clair qu'il annule l'image naturelle de $H^0_{ab}(F;G)$ dans l'espace de d\'epart. En notant $Im(H^0_{ab}(F;G))$ cette image et $Q_{2}=H^0_{ab}({\mathbb A}_{F};G)/Im(H^0_{ab}(F;G))$, on obtient un homomorphisme
${\bf q}_{2}:Q_{2}\to Q$. Notons $Q_{3}=\prod_{v\not\in V}H^0_{ab}(\mathfrak{o}_{v};G_{\sharp})$. C'est un sous-groupe de $H^0_{ab}({\mathbb A}_{F};G_{\sharp})$. En utilisant l'homomorphisme $\boldsymbol{\zeta}_{\sharp}$, on obtient de m\^eme un homomorphisme ${\bf q}_{3}:Q_{3}\to Q$. 

On note $Q_{0}$ le sous-groupe de $Q$ engendr\'e par les sous-groupes ${\bf q}_{i}(Q_{i})$ pour $i=1,2,3$.

 \ass{Lemme}{Le groupe $Q_{0}$ est un sous-groupe ouvert, ferm\'e et d'indice fini de $Q$. Le groupe $P^0$ est l'annulateur de $Q_{0}$ dans $P$. Le groupe $Q_{0}$ est l'annulateur de $P^0$ dans $Q$.}

Preuve. On a une suite d'homomorphismes
 $$(1) \qquad   ((1-\theta)(S))({\mathbb A}_{F})\to H^{1,0}({\mathbb A}_{F};S_{sc}\stackrel{1-\theta}{\to}(1-\theta)(S))\to Q.$$
 Pour $v\in Val_{F}$, le diagramme suivant est commutatif
 $$\begin{array}{ccccc}&&S(F_{v})&&\\ &\swarrow&&\searrow&\\ G(F_{v})&&&&(1-\theta)(S(F_{v}))\\ &&&&\downarrow\\ \downarrow&&&&((1-\theta)(S))(F_{v})\\&&&&\downarrow\\ G_{ab}(F_{v})&&\stackrel{\boldsymbol{\zeta}_{v}}{\to}&&Q_{v}\\ \end{array}$$
 Donc $Q_{0}$ contient l'image naturelle de $(1-\theta)(S(F_{v}))$. Ce groupe s'envoyant sur un sous-groupe ouvert d'indice fini de $((1-\theta)(S))(F_{v})$, $Q_{0}$ contient l'image naturelle d'un tel sous-groupe. L'assertion 6.5(2) montre qu'il contient aussi $((1-\theta)(S))(\mathfrak{o}_{v})$ pour presque toute place finie $v$. Donc $Q_{0}$ contient l'image par la suite d'homomorphismes (1) d'un sous-groupe ouvert de  $ ((1-\theta)(S))({\mathbb A}_{F})$. 
 D'apr\`es la d\'efinition de la topologie de $Q$ ([KS] page 147), $Q_{0}$ est donc un sous-groupe ouvert de $Q$.

 En [KS] page 151, Kottwitz et Shelstad d\'efinissent un homomorphisme
 $$(2) \qquad Q=H^{1,0}({\mathbb A}_{F}/F;S_{sc}\stackrel{1-\theta}{\to}(1-\theta)(S))\to coker(\mathfrak{A}_{S_{sc}}\stackrel{1-\theta}{\to} \mathfrak{A}_{(1-\theta)(S)}).$$
 Cet homomorphisme poss\`ede une section naturelle. On a l'\'egalit\'e $\mathfrak{A}_{S}=\mathfrak{A}_{S_{sc}}\times \mathfrak{A}_{G}$, d'o\`u
 $$\mathfrak{A}_{(1-\theta)(S)}=(1-\theta)(\mathfrak{A}_{S})=(1-\theta)(\mathfrak{A}_{S_{sc}})\times (1-\theta)(\mathfrak{A}_{G}).$$
 Le conoyau ci-dessus est donc isomorphe \`a $(1-\theta)(\mathfrak{A}_{G})$. En reprenant les d\'efinitions de [KS], on voit que l'image de ce groupe par la section de l'homomorphisme (2) co\"{\i}ncide avec l'image de 
 $$\mathfrak{A}_{G}\to \prod_{v\in Val_{\infty}(F)}G(F_{v})\stackrel{\prod_{v\in Val_{\infty}(F)} \iota_{v}\circ \boldsymbol{\zeta}_{v}}{\longrightarrow} Q.$$
 Notons $Q_{c}$ le noyau de (2). On obtient un isomorphisme
 $$Q\simeq Q_{c}\times (1-\theta)(\mathfrak{A}_{G}).$$
 C'est un hom\'eomorphisme et le groupe $Q_{c}$ est compact d'apr\`es [KS] lemme C.2.D. La description que l'on vient de donner de l'image de $(1-\theta)(\mathfrak{A}_{G})$ dans $Q$ montre que ce groupe est contenu dans $Q_{0}$. Donc $Q_{0}$ est le produit de $(1-\theta)(\mathfrak{A}_{G})$ et d'un sous-groupe ouvert de $Q_{c}$. Ce dernier \'etant compact, ce sous-groupe est aussi ferm\'e et d'indice fini. 
  D'o\`u la premi\`ere assertion de l'\'enonc\'e. 

 Par d\'efinition de $P^{0}$ et d'apr\`es la proposition 6.5, $P^{0}$ est l'annulateur du sous-groupe $Q'_{0}$ de $Q$ engendr\'e par les images des diff\'erents homomorphismes d\'ecrits dans cette proposition. C'est aussi l'annulateur de l'adh\'erence $Q''_{0}$ de ce sous-groupe. Tous les groupes d\'ecrits dans la proposition 6.5 sont inclus dans $Q_{0}$. Donc $Q'_{0}\subset Q_{0}$ et aussi $Q''_{0}\subset Q_{0}$ puisque $Q_{0}$ est ferm\'e. En sens inverse, $Q'_{0}$ contient ${\bf q}_{1}(Q_{1})$. Il contient $\iota_{v}\circ \boldsymbol{\zeta}_{v}(G_{ab}(F_{v}))$ pour tout $v$. Il est clair que ${\bf q}_{2}(Q_{2})$ est l'adh\'erence du groupe engendr\'e par ces sous-groupes quand $v$ parcourt $Val(F)$. Donc ${\bf q}_{2}(Q_{2})\subset Q''_{0}$. De m\^eme, $Q'_{0}$ contient  $\iota_{v}\circ \boldsymbol{\zeta}_{\sharp,v}(G_{ab}(\mathfrak{o}_{v}))$ pour tout $v\not\in V$. Le groupe ${\bf q}_{3}(Q_{3})$ est l'adh\'erence du groupe engendr\'e par ces sous-groupes quand $v$ parcourt $Val(F)-V$. Donc ${\bf q}_{3}(Q_{3})\subset Q''_{0}$. Cela d\'emontre que $Q''_{0}=Q_{0}$ donc que $P^0$ est l'annulateur de $Q_{0}$. 

Puisque $Q_{0}$ est un sous-groupe ouvert d'indice fini de $Q$, de la dualit\'e entre $P$ et $Q$ se d\'eduit une dualit\'e entre les groupes finis $P^0$ et $Q/Q_{0}$. Alors $Q_{0}$ est aussi l'annulateur de $P^0$ dans $Q$. $\square$

\bigskip

 \subsection{L'ensemble $D_{{\mathbb A}_{F}}$}
 
  Pour toute place $v\in Val(F)$, on a d\'efini l'ensemble $D_{v}$ en 5.4. La condition 5.1(3) signifie qu'il est non vide. On note $D_{{\mathbb A}_{F}}$ l'ensemble des familles $d=(d_{v})_{v\in Val(F)}$ telles que $d_{v}\in D_{v}$ pour tout $v$ et $\eta[d]\in \tilde{G}({\mathbb A}_{F})$, o\`u $\eta[d]=(\eta[d_{v}])_{v\in Val(F)}$. Soulignons qu'on n'impose aucune condition "globale" \`a la famille $r[d]=(r[d_{v}])_{v\in Val(F)}$.
  
On d\'efinit de m\^eme $D_{{\mathbb A}_{F}^V}$ en rempla\c{c}ant l'ensemble d'indices $Val(F)$ par $Val(F)-V$. On a l'\'egalit\'e $D_{{\mathbb A}_{F}}=D_{V}\times D_{{\mathbb A}_{F}^V}$. Pour un \'el\'ement  $d=(d_{v})_{v\in Val(F)}$,
  on a $\eta[d_{v}]\in \tilde{K}_{v}$ pour presque tout $v$, donc $ d_{v}\in D_{v}^{nr}$ pour presque tout $v$.  Inversement, on a dit en 
 5.5(2) que le lemme 1.6 impliquait que l'ensemble $D_{v}^{nr}$ \'etait non vide pour tout $v\not\in V$.  Puisque $D_{V}$ non vide d'apr\`es 5.1(3), l'ensemble $D_{{\mathbb A}_{F}}$ n'est pas vide lui non plus. On note $D_{{\mathbb A}_{F}^V}^{nr}$ le sous-ensemble des $ d=(d_{v})_{v\in Val(F)-V}\in D_{{\mathbb A}_{F}^V}$ tels que $ d_{v}\in D_{v}^{nr}$ pour tout $v\not\in V$. Il n'est pas vide lui non plus. 
 
 Soit $d=(d_{v})_{v\in Val(F)}\in D_{V}^{rel}\times D_{{\mathbb A}_{F}^V}^{nr}\subset D_{{\mathbb A}_{F}}$. On peut   fixer pour toute place $v\in Val(F)$ une paire de Borel $(B[d_{v}],S[d_{v}])$ v\'erifiant les conditions de 5.6. Rappelons celles-ci. Le tore $S[d_{v}]$ est d\'efini sur $F_{v}$. Le sous-groupe de Borel $B[d_{v}]$ est d\'efini sur $\bar{F}_{v}$.  La paire $(B[d_{v}],S[d_{v}])$ est conserv\'ee par $ad_{\eta[d_{v}]}$. Il existe $u\in G_{\eta[d_{v}]}$  tel que $ad_{r[d_{v}]u}(B[d_{v}],S[d_{v}])=(B^*,T^*)$ et que $ad_{r[d_{v}]u}$ se restreigne en un isomorphisme d\'efini sur $F_{v}$ de $S[d_{v}]$ sur $S$ (on rappelle que $S=T^*$ muni de l'action galoisienne $\sigma\mapsto \omega_{S}(\sigma)\sigma_{G^*}$). Nous allons imposer des conditions suppl\'ementaires "globales" \`a ces paires.

  On a fix\'e en 5.2 une paire de Borel \'epingl\'ee ${\cal E}^*$ de $G$, une paire de Borel \'epingl\'ee $\bar{{\cal E}}$ de $\bar{G}$  et des \'el\'ements $\nu\in T^*$ et $e\in Z(\tilde{G},{\cal E}^*)$.  On a not\'e $\theta^*=ad_{e}$. On fixe pour tout $\sigma\in \Gamma_{F}$ un \'el\'ement $u_{{\cal E}^*}(\sigma)\in G_{SC}(\bar{F})$ tel que $\sigma_{G^*}=ad_{u_{{\cal E}^*}(\sigma)}\circ \sigma$ conserve ${\cal E}^*$. On peut supposer que $\sigma\mapsto u_{{\cal E}^*}(\sigma)$ est continue et se factorise par un quotient fini de $\Gamma_{F}$. On peut aussi supposer $u_{{\cal E}^*}(1)=1$. D'autre part, les applications $\sigma\mapsto \omega_{\bar{G}}(\sigma)$ et  $\sigma\mapsto \omega_{S}(\sigma)$ sont  des cocycles de $\Gamma_{F}$ dans $W^{\theta^*}$ (muni de l'action quasi-d\'eploy\'ee). D'apr\`es [K1] corollaire 2.2, on peut fixer $x\in G_{SC}^{\theta^*}(\bar{F})$ tel que $x\sigma_{G^*}(x)^{-1}$ normalise $T^*$ et ait $\omega_{S}(\sigma)$ pour image dans $W$. On fixe une extension galoisienne finie $E$ de $F$ telle que
 
 - ${\cal E}^*$  et $\bar{{\cal E}}$ soient d\'efinies sur $E$ et $G$ soit d\'eploy\'e sur $E$;
 
 - $\nu\in T^*(E)$, $e\in Z(\tilde{G},{\cal E}^*;E)$, $x\in G_{SC}^{\theta^*}(E)$;
 
 - l'application $\sigma\mapsto u_{{\cal E}^*}(\sigma)$ se factorise par $Gal(E/F)$ et prend ses valeurs dans $G_{SC}(E)$:
 
 - l'application  $\sigma\mapsto \omega_{\bar{G}}(\sigma)$ se factorise par $Gal(E/F)$.
 
 Il en r\'esulte que toutes les actions galoisiennes co\"{\i}ncident sur $\Gamma_{E}$ et que tous les groupes qui interviennent sont d\'eploy\'es sur $E$. Utilisons les d\'efinitions de 1.5. Fixons un ensemble fini $V'$ de places de $F$, contenant $V$, de sorte que, pour toute place $v\not\in V'$ et toute place $w'$ de $E$ au-dessus de $v$, $E_{w'}/F_{v}$ soit non ramifi\'ee   et que les propri\'et\'es suivantes soient v\'erifi\'ees
 
 - $K_{w'}$ est le sous-groupe compact hypersp\'ecial issu de la paire de Borel \'epingl\'ee ${\cal E}^*$;
 
 - $e\in \tilde{K}_{w'}$,  $\nu\in T^*(\mathfrak{o}_{w'})$, $x\in K_{w'}$ et, pour tout $\sigma\in \Gamma_{F}$, $u_{{\cal E}^*}(\sigma)\in K_{w'}$.

  Soit $v\in Val(F)$. Rappelons que $v$ a \'et\'e prolong\'ee en une place $\bar{v}$ de $\bar{F}$, cf. [VI] 1.1. Le corps $\bar{F}_{v}$ a \'et\'e identifi\'e \`a la cl\^oture alg\'ebrique de $F_{v}$ dans le compl\'et\'e de $\bar{F}$ en $\bar{v}$. Par abus de notations, notons-le $\bar{F}_{\bar{v}}$.   Le groupe $\Gamma_{F_{v}}$ a \'et\'e identifi\'e au fixateur de $\bar{v}$ dans $\Gamma_{F}$. Notons-le plut\^ot $\Gamma_{\bar{v}}$. On notera sans plus de commentaire   $w$ la restriction de $\bar{v}$ \`a $E$. Soit $w'$ une autre place de $E$ divisant $v$.  
  On fixe une fois pour toutes un \'el\'ement $\tau\in  \Gamma_{F}$ telle que $\tau(w)=w'$ (avec $\tau=1$ dans le cas $w'=w$). Notons $\bar{v}'=\tau(\bar{v})$. De $\tau$ se d\'eduit un isomorphisme de $\bar{F}_{\bar{v}}$ sur $\bar{F}_{\bar{v}'}$. Pour toute vari\'et\'e alg\'ebrique $X$ d\'efini sur $F_{v}$, on a aussi un isomorphisme $\tau:X(\bar{F}_{\bar{v}})\to X(\bar{F}_{\bar{v}'})$.  
   Pour une paire  $(B[d_{v}],S[d_{v}])$ comme ci-dessus, le groupe $B[d_{v}]$ est pr\'ecis\'ement d\'efini sur $\bar{F}_{\bar{v}}$. En fait, le tore $S[d_{v}]$ est d\'eploy\'e sur $E_{w}$ donc tout sous-groupe de Borel contenant ce tore est d\'efini sur $E_{w}$. En particulier $B[d_{v}]$ est  d\'efini sur $E_{w}$. Notons
  $S[d]({\mathbb A}_{E})$ le produit restreint des $S[d_{v}](E_{w'})$ sur toutes les places $v\in Val(F)$ et les places $w'$ de $E$ divisant $v$. La restriction est relative aux sous-groupes $S[d_{v}](\mathfrak{o}_{w'})$ qui sont d\'efinis pour presque tous $v$ et $w'$.  Le groupe de Galois $Gal(E/F)$ agit naturellement sur $S[d]({\mathbb A}_{E})$.
  
  \ass{Proposition}{Soit  $d=(d_{v})_{v\in Val(F)}\in D_{V}^{rel}\times D_{{\mathbb A}_{F}^V}^{nr}\subset D_{{\mathbb A}_{F}}$.  On peut fixer 
  
  - pour toute place $v$ une paire de Borel  $(B[d_{v}],S[d_{v}])$ v\'erifiant les conditions de 5.6;
  
  - un \'el\'ement $g=(g_{w'})_{w'\in Val(E)}\in G_{SC}({\mathbb A}_{E})$;
  
  - un \'el\'ement $t=(t_{w'})_{w'\in Val(E)}\in ((1-\theta^*)(T^*))({\mathbb A}_{E})$;
  
  de sorte que les conditions suivantes soient v\'erifi\'ees:
  
  (i) si $v\not\in V'$ et $\eta[d_{v}]\in \tilde{K}_{v}$, alors, pour tout $w'$ divisant $v$, on a $S[d_{v}](\mathfrak{o}_{w'})\subset K_{w'}$, $g_{w'}\in K_{sc,w'}$ et $t_{w'}\in ((1-\theta^*)(T^*))(\mathfrak{o}_{w'})$;
  
  (ii) $ad_{g}(S[d]({\mathbb A}_{E}))=S({\mathbb A}_{E})$ et $ad_{g}$ se restreint en un isomorphisme de $S[d]({\mathbb A}_{E})$ sur $S({\mathbb A}_{E})$ qui est \'equivariant pour les actions galoisiennes;
  
  (iii) pour toute place $v\in Val(F)$, on a $ad_{g_{w}}(B[d_{v}],S[d_{v}])=(B^*,T^*)$ et $g_{w}\in T^*(\bar{F}_{\bar{v}})r[d_{v}]G_{\eta[d_{v}]}(\bar{F}_{\bar{v}})$;
  
  (iv) pour toute place $v\in Val(F)$ et toute place $w'$ de $E$ divisant $v$, on a l'\'egalit\'e
  $g_{w'}=x\tau_{G^*}(x)^{-1}u_{{\cal E}^*}(\tau)\tau(g_{w})$, o\`u $\tau\in \Gamma_{F}$ est l'\'el\'ement fix\'e tel que $\tau(w)=w'$;
  
  (v) $ad_{g}(\eta[d])=t\eta$.}
 
 {\bf Remarques.} (1) La condition (i) entra\^{\i}ne que $S[d]({\mathbb A}_{E})$ est contenu dans $G({\mathbb A}_{E})$. Cela donne un sens \`a la condition (ii).
 
 (2)   On peut choisir arbitrairement la paire $(B[d_{v}],S[d_{v}])$ pour un ensemble fini de places $v$, pourvu que ces paires satisfassent aux conditions de 5.6. 
 
 (3) Sous les hypoth\`eses de (i), la premi\`ere inclusion se g\'en\'eralise en $S[d_{v}](\mathfrak{o}_{v}^{nr})\subset K_{v}^{nr}$. En prenant les invariants par $\Gamma_{v}^{nr}$, on en d\'eduit $S[d_{v}](\mathfrak{o}_{v})\subset K_{v}$.
 
 \bigskip
 
 Preuve. Soit $v\in Val(F)$, fixons une paire de Borel $(B[d_{v}],S[d_{v}])$ v\'erifiant les conditions de 5.6. On va prouver l'existence de $g_{v}=(g_{w'})_{w'\vert  v}$ et $t_{v}=(t_{w'})_{w'\vert v}$ (o\`u $w'\vert  v$ signifie que $w'$ divise $v$) v\'erifiant les analogues des conditions (ii) \`a (v) o\`u l'on se restreint aux places de $E$ divisant $v$. Comme on l'a dit, on note $w$ la restriction de $\bar{v}$ \`a $E$. Les tores $S[d_{v}]$ et $T^*$ sont d\'eploy\'es sur $E_{w}$. Il en r\'esulte que les groupes de Borel $B[d_{v}]$ et $B^*$ sont d\'efinis sue $E_{w}$. Il existe donc $g_{w}\in G_{SC}(E_{w})$ tel que $ad_{g_{w}}(B[d_{v}],S[d_{v}])=(B^*,T^*)$. On fixe un tel \'el\'ement. 
 L'une des propri\'et\'es de la paire $(B[d_{v}],S[d_{v}])$ est qu'il  existe   $u\in G_{\eta[d_{v}]}(\bar{F}_{\bar{v}})$ tel que $ad_{r[d_{v}]u}(B[d_{v}],S[d_{v}])=(B^*,T^*)$. Alors $ad_{g_{w}u^{-1}r[d_{v}]^{-1}}$ conserve $(B^*,T^*)$. Cela implique 
  que  $g_{w}$ appartient \`a $T^*(\bar{F}_{\bar{v}})r[d_{v}]u$. La condition (iii) est donc satisfaite. Soit $w'$ une place de $E$ au-dessus de $v$. On note $\tau$ l'\'el\'ement fix\'e de $\Gamma_{F}$ tel que $\tau(w)=w'$. On d\'efinit $g_{w'}$ par l'\'egalit\'e de la condition (iv). Montrons que
 
    (4) $ad_{g_{w'}}( S[d_{v}])=T^*$;
   
   (5) pour $\sigma\in \Gamma_{\bar{v}'}$ et $s\in S[d_{v}](\bar{F}_{\bar{v}'})$, on a l'\'egalit\'e $ad_{g_{w'}}\circ \sigma (s)= \sigma_{S}\circ ad_{g_{w'}}(s)$, o\`u $\sigma_{S}=\omega_{S}(\sigma)\circ \sigma_{G^*}$;
   
   (6) pour $s\in S[d_{v}](\bar{F}_{\bar{v}})$, on a $ad_{g_{w'}}\circ \tau(s)=\tau_{S}\circ ad_{g_{w}}(s)$;
   
   (7) il existe $t_{w'}\in ((1-\theta^*)(T^*))(E_{w'})$ tel que $ad_{g_{w'}}(\eta[d_{v}])=t_{w'}\eta$.

   La preuve de (4), (5) et (6) est similaire \`a celle de [VI] 3.6(5) et (6). On la laisse au lecteur. Prouvons (7). Remarquons que, si l'on prouve l'existence de $t_{w'}\in (1-\theta^*)(T^*(\bar{F}_{\bar{v}'}))$ satisfaisant l'\'egalit\'e ci-dessus, on a n\'ecessairement $t_{w'}\in ((1-\theta^*)(T^*))(E_{w'})$ puisque les autres termes de cette \'egalit\'e sont d\'efinis sur $E_{w'}$. Pour $w'=w$, on sait qu'il existe $t_{0}\in T^*(\bar{F}_{\bar{v}})$ et $u\in G_{\eta[d_{v}]}(\bar{F}_{\bar{v}})$ tels que $g_{w}=t_{0}r[d_{v}]u$. On a alors 
   $$ad_{g_{w}}(\eta[d_{v}])=ad_{t_{0}}\circ ad_{r[d_{v}]}(\eta[d_{v}])=ad_{t_{0}}(\eta)=(1-\theta^*)(t_{0})\eta.$$
   D'o\`u l'assertion avec $t_{w}=(1-\theta^*)(t_{0})$. Pour une autre place $w'$, on a
   $$ad_{g_{w'}}(\eta[d_{v}])=ad_{x\tau_{G^*}(x)^{-1}u_{{\cal E}^*}(\tau)}\circ ad_{\tau(g_{w})}(\eta[d_{v}]).$$
   On a $\tau(\eta[d_{v}])=\eta[d_{v}]$ puisque $\eta[d_{v}]\in \tilde{G}(F_{v})$. Donc
   $$ad_{\tau(g_{w})}(\eta[d_{v}])=\tau\circ ad_{g_{w}}(\eta[d_{v}])=\tau(t_{w}\eta).$$
   L'application $\tau_{S}$ se prolonge en une application de $((1-\theta^*)(T^*))(\bar{F}_{\bar{v}})$ sur $((1-\theta^*)(T^*))(\bar{F}_{\bar{v}'})$ et on a
  $$ad_{x\tau_{G^*}(x)^{-1}u_{{\cal E}^*}(\tau)}\circ \tau(t_{w})=\tau_{S}(t_{w}).$$
  Il reste \`a prouver que 
   $$ad_{x\tau_{G^*}(x)^{-1}u_{{\cal E}^*}(\tau)}\circ \tau(\eta)\in (1-\theta^*)(T^*(\bar{F}))\eta.$$ 
  On a \'ecrit $\eta=\nu e$. Soit $z(\tau)\in Z(G)$ tel que $ad_{u_{{\cal E}^*}(\tau)}\circ \tau(e)=z(\tau)^{-1}e$. Parce que $x$ appartient \`a $G_{SC}^{\theta^*}$, on a 
   $$ad_{x\tau_{G^*}(x)^{-1}u_{{\cal E}^*}(\tau)}\circ \tau(\eta)=z(\tau)^{-1}ad_{x\tau_{G^*}(x)^{-1}u_{{\cal E}^*}(\tau)}\circ \tau(\nu)e=z(\tau)^{-1}\omega_{S}(\tau)\circ \tau_{G^*}(\nu)e.$$
   On peut d\'ecomposer $\omega_{S}(\tau) $ en $\omega_{S,\bar{G}}(\tau)\omega_{\bar{G}}(\tau)$. Les conditions impos\'ees en 1.1 \`a $\eta$ impliquent que $z(\tau)^{-1}\omega_{\bar{G}}(\tau)\circ \tau_{G^*}(\nu)$ appartient \`a $(1-\theta^*)(T^*(\bar{F}))\nu$. L'\'el\'ement $\omega_{S,\bar{G}}(\tau)$ appartient \`a $W^{\bar{G}}$ et tout \'el\'ement de ce groupe conserve l'ensemble $(1-\theta^*)(T^*(\bar{F}))\nu$. On obtient
  $$ad_{x\tau_{G^*}(x)^{-1}u_{{\cal E}^*}(\tau)}\circ \tau(\eta)\in (1-\theta^*)(T^*(\bar{F}))\nu e=(1-\theta^*)(T^*(\bar{F}))\eta.$$
  Cela prouve (7).
  
  Pour $w'$ divisant $v$, on a d\'efini le terme $g_{w'}$   et la relation (7) d\'efinit le terme $t_{w'}$. 
En posant $g_{v}=(g_{w'})_{w'\vert  v}$ et $t_{v}=(t_{w'})_{w'\vert v}$, la relation (7) entra\^{\i}ne la condition (v) restreinte aux places divisant $v$. Les relations (4), (5) et (6) entra\^{\i}nent la condition (ii) restreinte aux m\^emes places. 

Supposons maintenant que $v\not\in V'$ et que $\eta[d_{v}]$ appartient \`a $\tilde{K}_{v}$. On va prouver qu'en choisissant convenablement la paire $(B[d_{v}],S[d_{v}])$, on peut imposer la condition (i). 
   Puisque $v\not\in V$, on peut fixer une paire de Borel \'epingl\'ee ${\cal E}_{0}=(B_{0},T_{0},(E_{0,\alpha})_{\alpha\in \Delta})$ de $G$, d\'efinie sur $F_{v}$, dont est issu de groupe $K_{v}$. Puisque ${\cal E}^*$ et ${\cal E}_{0}$ sont toutes deux d\'efinies sur $E_{w}$, il existe $y_{ad}\in G_{AD}(E_{w})$ tel que $ad_{y_{ad}}({\cal E}_{0})={\cal E}^*$. L'automorphisme $ad_{y_{ad}}$ conserve le groupe $K_{w}$ puisque ce groupe est issu de chacune des deux paires. Donc $y_{ad}\in K_{ad,w}$. On sait que les deux applications
 $$\begin{array}{ccc}T^*_{ad}(\mathfrak{o}_{w})\times K_{sc,w}&\to& K_{ad,w}\\ (t,k)&\mapsto &t\pi(k)\\ \end{array}$$
 et
 $$T^*_{sc}(\mathfrak{o}_{v}^{nr})\to T^*_{ad}(\mathfrak{o}_{v}^{nr})$$
 sont surjectives (le corps $F_{v}^{nr}$ est ici un sous-corps de $\bar{F}_{\bar{v}}$). On peut donc fixer $t_{1}\in T^*_{sc}(\mathfrak{o}_{v}^{nr})$ et $y\in K_{sc,w}$ tels que $y_{ad}=\pi(t_{1} y)$. On a $t_{1} y\in K_{v}^{nr}$. 
 Posons $\eta_{1}=ad_{(t_{1} y)^{-1}}(\eta)$. On reprend maintenant la preuve du lemme 1.6. Puisque $\eta\in \tilde{K}_{w}$, on a $\eta_{1}\in \tilde{K}_{v}^{nr}$. De plus, $ad_{\eta_{1}}$ conserve $(B_{0},T_{0})$. Fixons $\nu_{0}\in T_{0}(\mathfrak{o}_{v}^{nr})$ et $e_{0}\in Z(\tilde{G},{\cal E}_{0})(F_{v}^{nr})$ tels que $\nu_{0}e_{0}\in \tilde{K}$, cf.  1.5.  On peut \'ecrire $\eta_{1}=\nu_{1}e_{0}$, avec $\nu_{1}\in T_{0}$. Puisque $\eta_{1}\in \tilde{K}_{v}^{nr}$, on a $\nu_{1}\in T_{0}(\mathfrak{o}_{v}^{nr})$.  Introduisons le cocycle $z:\Gamma_{F_{v}}\to Z(G)\cap T_{0}(\mathfrak{o}_{v}^{nr})$ tel que $\sigma(e_{0})=z(\sigma)^{-1}e_{0}$ et posons $\theta=ad_{e_{0}}$. La propri\'et\'e de d\'efinition de $\eta$ se transporte \`a $\eta_{1}$. C'est-\`a-dire qu'en identifiant $W$ au groupe de Weyl relatif \`a $T_{0}$, il existe une cocha\^{\i}ne $t:\Gamma_{F_{v}}\to (1-\theta)(T_{0}(\bar{F}_{\bar{v}}))$ telle que 
  $\omega_{\bar{G}}(\sigma)\sigma(\nu_{1})=z(\sigma)t(\sigma)\nu_{1}$ pour tout $\sigma\in \Gamma_{\bar{v}}$. 
  Puisque $\nu_{1}\in T_{0}(\mathfrak{o}_{v}^{nr})$ et que ce groupe est normalis\'e par $W$,
cette relation implique que $t(\sigma)\in ((1-\theta)( T_{0}))(\mathfrak{o}_{v}^{nr})$ et que $\sigma\mapsto t(\sigma)$ est un cocycle \`a valeurs dans ce groupe, si on munit celui-ci de l'action $\sigma\mapsto \sigma_{\bar{G}}=\omega_{\bar{G}}(\sigma)\circ \sigma$. Un tel cocycle est un cobord. De plus, l'hypoth\`ese $v\not\in V_{ram}$ implique que $1-\theta:T_{0}(\mathfrak{o}_{v}^{nr})\to ((1-\theta)( T_{0}))(\mathfrak{o}_{v}^{nr})$ est surjective. On peut donc fixer $t_{0}\in T_{0}(\mathfrak{o}_{v}^{nr})$ tel que $t(\sigma)=(1-\theta)(t_{0}\sigma_{\bar{G}}(t_{0})^{-1})$. Posons $\nu_{2}=\nu_{1}(1-\theta)(t_{0})$.  On a $\sigma_{\bar{G}}(\nu_{2})=z(\sigma)\nu_{2}$. On introduit le groupe $G_{SC}^{\theta}$ des points fixes de $\theta$ dans $G_{SC}$. De ${\cal E}_{0}$ se d\'eduit une paire de Borel \'epingl\'ee de ce groupe puis un sch\'ema en groupes ${\cal K}_{v}^1$. En appliquant 1.5(5), on construit $k\in {\cal K}_{v}^1(\mathfrak{o}_{w})$ tel que, pour tout $\sigma\in \Gamma_{F_{v}}$, $k^{-1}\sigma(k)$ normalise $T_{0}$ et ait $\omega_{\bar{G}}(\sigma)$ pour image dans $W^{\theta}$. On pose $\eta_{\star}=k\nu_{2}e_{0}k^{-1}$. Le m\^eme calcul qu'en 1.6  montre que $\eta_{\star}\in \tilde{K}_{v}$ et que la paire de Borel $(kB_{0}k^{-1}\cap G_{\eta_{\star}},kT_{0}k^{-1}\cap G_{\eta_{\star}})$ de $G_{\eta_{\star}}$ est d\'efinie sur $F_{v}$. En reprenant la preuve de [W1] lemme 5.4, on peut la compl\'eter  en une paire de Borel \'epingl\'ee d\'efinie sur $F_{v}$ de sorte que le sch\'ema en groupes ${\cal K}_{\star}$ issu de cette paire v\'erifie la condition ${\cal K}_{\star}(\mathfrak{o}_{v}^{nr})=K_{v}^{nr}\cap G_{\eta_{\star}}(\bar{F}_{\bar{v}})$. On note ${\cal K}_{\star,sc}$ le sch\'ema en groupes associ\'e dans le groupe $G_{\eta_{\star},SC}$. Pour $\sigma\in \Gamma_{F_{v}}$, on pose $\omega_{S,\bar{G}}(\sigma)=\omega_{S}(\sigma)\omega_{\bar{G}}(\sigma)^{-1}$. Par d\'efinition, c'est un \'el\'ement du groupe de Weyl $W^{\bar{G}}$, lequel s'identifie \`a $W^{G_{\eta_{\star}}}$. En appliquant de nouveau 1.5(5), on construit $h\in {\cal K}_{\star,sc}(\mathfrak{o}_{w})$ tel que, pour tout $\sigma\in \Gamma_{F_{v}}$, $h^{-1}\sigma(h)$ normalise $kT_{0}k^{-1}\cap G_{\eta_{\star}}$ et ait pour image $\omega_{S,\bar{G}}(\sigma)$ dans $W^{G_{\eta_{\star}}} $. Posons $(B_{\star},S_{\star})=ad_{hk}(B_{0},T_{0})$, $t_{\star}=t_{1} yt_{0}^{-1}y^{-1}$, $g_{\star}=yk^{-1}h^{-1}$ et $r_{\star}=t_{\star}g_{\star}$. On voit que $d_{\star}=(\eta_{\star},r_{\star})$ appartient \`a $ D_{v}$, que $\eta_{\star}$ appartient \`a $\tilde{K}_{v}$, que la paire  $(B_{\star},S_{\star})$ v\'erifie les conditions de 5.6 relatives \`a l'\'el\'ement $d_{\star}$,  que $t_{\star}\in T^*(\mathfrak{o}_{v}^{nr})$, que $g_{\star}\in K_{sc,w}$ et que $ad_{g_{\star}}$ envoie  la paire $(B_{\star},S_{\star})$ sur $(B^*,T^*)$.

  Revenons \`a  notre \'el\'ement quelconque $d_{v}=(\eta[d_{v}],r[d_{v}])\in D_{v}$ tel que $\eta[d_{v}]\in \tilde{K}_{v}$. D'apr\`es 5.5(3), on peut fixer $k_{\sharp}\in \underline{K}_{\sharp,v}$ et $u\in G_{\eta[d_{v}]}(\bar{F}_{\bar{v}})$ tels que $\eta_{\star}=k_{\sharp}^{-1}\eta[d_{v}] k_{\sharp}$ et $r_{\star}=r[d_{v}]uk_{\sharp}$. L'automorphisme $ad_{k_{\sharp}}$ envoie $G_{\eta_{\star}}$ sur $G_{\eta[d_{v}]}$ et est d\'efini sur $F_{v}$. On peut donc prendre pour paire $(B[d_{v}],S[d_{v}])$ la paire $ad_{k_{\sharp}}(B_{\star},S_{\star})$. Comme plus haut, l'application produit
 $$  K_{sc,w}\times T_{0}(\mathfrak{o}_{v}^{nr})\to K_{\sharp,w}$$
 est surjective. Mais $S_{\star}$ est conjugu\'e \`a $T_{0}$ par un \'el\'ement de $K_{sc,w}$. En conjuguant la propri\'et\'e ci-dessus, on obtient que  l'application produit
 $$  K_{sc,w}\times S_{\star}(\mathfrak{o}_{v}^{nr})\to K_{\sharp,w}$$
 est surjective. On peut donc \'ecrire $k_{\sharp}=zls$, avec $z\in Z(G)^{\theta}(\bar{F}_{\bar{v}})$, $s\in S_{\star}(\mathfrak{o}_{v}^{nr})$ et $l\in K_{sc,w}$. Posons $g_{w}=g_{\star}l^{-1}$ et $x=z^{-1}t_{\star}g_{\star}s^{-1}g_{\star}^{-1}$. On a $g_{w}\in K_{w,sc}$,  $ad_{g_{w}}(B[d_{v}],S[d_{v}])=(B^*,T^*)$
  et $g_{w}=x^{-1}r[d_{v}]u$. Puisque $x\in T^*(\bar{F}_{\bar{v}})$, la  condition (iii) est satisfaite. On a 
 $ad_{g_{w}}(\eta[d_{v}])= t_{w}\eta$, o\`u $t_{w}=(\theta^*-1)(x)$. Puisque $\eta[d_{v}]\in \tilde{K}_{v}$ et $g_{w}\in K_{sc,w}$, on a $ad_{g_{w}}(\eta[d_{v}])\in \tilde{K}_{w}$. On a aussi $\eta\in \tilde{K}_{w}$ par d\'efinition de $V'$. L'\'egalit\'e pr\'ec\'edente entra\^{\i}ne alors $t_{w}\in (1-\theta)(T^*(\bar{F}_{\bar{v}}))\cap K_{w}=((1-\theta)(T^*))(\mathfrak{o}_{w})$.
  Enfin, puisque $K_{w}$ est issu de ${\cal E}^*$, on a $T^*(\mathfrak{o}_{w})\subset K_{w}$. En conjuguant par $g_{w}^{-1}\in K_{sc,w}$, on en d\'eduit $S[d_{v}](\mathfrak{o}_{w})\subset K_{w}$. Cela satisfait les conditions (i), (iii) et (v) en la place $w$. Comme on l'a vu dans la premi\`ere partie de la preuve, la condition (iii) implique (ii).
  Pour  une autre place $w'$ de $E$ au-dessus de $v$, on construit $g_{w'}$ et $t_{w'}$ comme dans cette premi\`ere partie.
  Puisque le sch\'ema en groupes ${\cal K}_{v}$ est d\'efini sur $\mathfrak{o}_{v}$, on a $\tau(K_{sc,w})=K_{sc,w'}$. Les conditions impos\'ees \`a $V'$ entra\^{\i}nent que tous les termes de la formule (iv) appartiennent \`a $K_{sc,w'}$ donc $g_{w'}\in K_{sc,w'}$. Le m\^eme raisonnement que dans le cas de la place $w$ entra\^{\i}ne alors que $t_{w'}\in ((1-\theta)(T^*))(\mathfrak{o}_{w'})$ et que
  $S[d_{v}](\mathfrak{o}_{w'})\subset K_{w'}$.  Cela ach\`eve la preuve. $\square$

 \bigskip

 \subsection{L'ensemble $D_{F}$}
 
  On note $D_{F}$ l'ensemble des couples $d=(\eta[d],r[d])\in \tilde{G}(F)\times G(\bar{F})$ tels que
 
 - $r[d]\eta[d]r[d]^{-1}=\eta$;
 
 - en utilisant la paire de Borel $ad_{r[d]^{-1}}(B^*,T^*)$ dans la construction de 1.2, on ait l'\'egalit\'e $(\mu_{\eta[d]},\omega_{\eta[d]})=(\mu,\omega_{\bar{G}})$. 
 
 Soit $d\in D_{F}$. On a fix\'e en 5.2 une paire de Borel \'epingl\'ee $\bar{{\cal E}}$ de $\bar{G}$ d\'efinie sur $\bar{F}$. Posons $(B^*[d],T^*[d])=ad_{r[d]^{-1}}(B^*,T^*) $,  $(\bar{B}[d],\bar{T}[d])=(B^*[d]\cap G_{\eta[d]},T^*[d]\cap G_{\eta(d]})$ et $\bar{{\cal E}}[d]=ad_{r[d]^{-1}}(\bar{{\cal E}})$.  Alors $\bar{{\cal E}}[d]$ est une paire de Borel \'epingl\'ee de $G_{\eta[d]}$ d\'efinie sur $\bar{F}$ dont la paire de Borel sous-jacente est $(\bar{B}[d],\bar{T}[d])$. Pour tout $\sigma\in \Gamma_{F}$, fixons $\bar{u}[d](\sigma)\in G_{\eta[d],SC}(\bar{F})$ tel que $ad_{\bar{u}[d](\sigma)}\circ \sigma$ conserve $\bar{{\cal E}}[d]$.  On note $\sigma\mapsto \sigma_{G_{\eta[d]}^*}=ad_{\bar{u}[d](\sigma)}\circ \sigma$ l'action quasi-d\'eploy\'ee qui conserve $\bar{{\cal E}}[d]$. On suppose que $\sigma\mapsto \bar{u}[d](\sigma)$ est continue et que $\bar{u}[d](1)=1$. La deuxi\`eme condition ci-dessus signifie que $ad_{r[d]}$, qui envoie $T^*[d]$ sur $T^*$, entrelace l'action $\sigma\mapsto  \sigma_{G_{\eta[d]}^*}$ sur $T^*[d]$ avec l'action $\sigma\mapsto \omega_{\bar{G}}(\sigma)\circ \sigma_{G^*}$ sur $T^*$. Transportons par $ad_{r[d]^{-1}}$ le cocycle $\omega_{S,\bar{G}}$ en un cocycle \`a valeurs dans le groupe de Weyl de $G_{\eta[d]}$ relatif \`a $\bar{T}[d]$ (pour l'action quasi-d\'eploy\'ee). D'apr\`es [K1] corollaire 2.2, on peut fixer $\bar{x}[d]\in G_{\eta[d],SC}(\bar{F})$ tel que, pour tout $\sigma\in \Gamma_{F}$,  $\bar{x}[d] \sigma_{G_{\eta[d]}^*}(\bar{x}[d])^{-1}$ normalise $\bar{T}[d]$ et que son image dans le groupe de Weyl soit $\omega_{S,\bar{G}}(\sigma)$. Alors $ad_{r[d]}$ entrelace l'action $\sigma\mapsto ad_{\bar{x}[d] \sigma_{G_{\eta[d]}^*}(\bar{x}[d])^{-1}}\circ \sigma_{G_{\eta[d]}^*}$ sur $T^*[d]$ avec l'action $\sigma\mapsto \omega_{S}(\sigma)\circ \sigma_{G^*}$ sur $T^*$.  Fixons  une d\'ecomposition $r[d]=z[d]r[d]_{sc}$ avec $z\in Z(G;\bar{F})$ et $r[d]_{sc}\in G_{SC}(\bar{F})$.  Il est facile de traduire la condition pr\'ec\'edente par la propri\'et\'e 
  
 (1) pour tout $\sigma\in \Gamma_{F}$, il existe $t(\sigma)\in T^*_{sc}(\bar{F})$ tel que 
 $$x\sigma_{G^*}(x)^{-1}u_{{\cal E}^*}(\sigma)\sigma(r[d]_{sc})=t(\sigma)r[d]_{sc}\bar{x}[d]\sigma_{G_{\eta[d]}^*}(\bar{x}[d])^{-1}\bar{u}[d](\sigma).$$
 
   Pour $d\in D_{F}$ et $v\in Val(F)$, on note $d_{v}$ le m\^eme couple $(\eta[d],r[d])$, vu comme \'el\'ement de $\tilde{G}(F_{v})\times G(\bar{F}_{v})$. L'application $d\mapsto (d_{v})_{v\in Val(F)}$ est   une injection  $D_{F}\subset D_{{\mathbb A}_{F}}$. Soit $d\in D_{F}$. Supposons que l'image de $d$ dans $D_{{\mathbb A}_{F}}$ appartienne \`a $D_{V}^{rel}\times D_{{\mathbb A}_{F}^V}^{nr}$.    On reprend les constructions du paragraphe pr\'ec\'edent. On impose de plus au corps $E$ les conditions suivantes  
   
   - $z[d]\in Z(G;E)$, $r[d]_{sc}\in G_{SC}(E)$, $\bar{x}[d]\in G_{\eta[d],SC}(E)$ et $\bar{{\cal E}}$ est d\'efinie sur $E$;
   
   - l'application $\sigma\mapsto \bar{u}[d](\sigma)$ se factorise par $Gal(E/F)$ et prend ses valeurs dans $G_{\eta[d],SC}(E)$.
   
   Il en r\'esulte que l'application $\sigma\mapsto t(\sigma)$  se factorise par $Gal(E/F)$ et prend ses valeurs dans $T_{sc}^*(E)$. Construisons des paires de Borel et des \'el\'ements $g$ et $t$ v\'erifiant la proposition 6. 7. 
   
   \ass{Lemme}{     Sous ces hypoth\`eses, l'\'el\'ement $g$ appartient \`a $T^*_{sc}({\mathbb A}_{E})r[d]_{sc}G_{\eta[d],SC}({\mathbb A}_{E})$.}
   
   Preuve. 
   Soient $v$ une place de $F$ et $w'$ une place de $E$ au-dessus de $v$.     Montrons que
   
   (2) $g_{w'}\in T_{sc}^*(E_{w'})r[d]_{sc}G_{\eta[d],SC}(E_{w'})$.

   Supposons d'abord que $w'=w$ soit  la restriction de $\bar{v}$ \`a $E$.   Par construction, il existe $u_{1}\in G_{\eta[d]}(\bar{F}_{\bar{v}})$ de sorte que $ad_{u_{1}}(B[d_{v}],S[d_{v}])=(B^*[d],T^*[d])$. Les deux couples $(B[d_{v}]\cap G_{\eta[d]},S[d_{v}]\cap G_{\eta[d]})$ et $(B^*[d]\cap G_{\eta[d]},T^*[d]\cap G_{\eta[d]})$ sont des paires de Borel de $G_{\eta[d]}$ qui sont d\'efinies sur $E_{w}$. Il existe donc $u\in G_{\eta[d],SC}(E_{w})$ telle que la seconde soit l'image de la premi\`ere par $ad_{u}$. Puisque $u_{1}$ v\'erifie la m\^eme propri\'et\'e, on a $u_{1}\in (T^*[d]\cap G_{\eta[d]})u$. Alors on a aussi $ad_{u}(B[d_{v}],S[d_{v}])=(B^*[d],T^*[d])$. Donc $ad_{r[d]_{sc}u}(B[d_{v}],S[d_{v}])=(B^*,T^*)$. Puisque $g_{w}$ v\'erifie la m\^eme propri\'et\'e et que les deux \'el\'ements $r[d]_{sc}u$ et $g_{w}$ appartiennent \`a $G_{SC}(E_{w})$, ces deux \'el\'ements diff\`erent par multiplication \`a gauche par un \'el\'ement de $T^*_{sc}(E_{w})$. Cela prouve (2) pour la place $w$. 
   
   Consid\'erons maintenant une autre place $w'$ au-dessus de $v$.    Rappelons que
$g_{w'}=x\tau_{G^*}(x)^{-1}u_{{\cal E}^*}(\tau)\tau(g_{w})$. Ecrivons $g_{w}=t_{w}r[d]_{sc}u_{w}$, avec $t_{w}\in T^*_{sc}(E_{w})$ et $u_{w}\in G_{\eta[d],SC}(E_{w})$. Posons $t_{0,w'}=ad_{x\tau_{G^*}(x)^{-1}u_{{\cal E}^*}(\tau)}\circ \tau(t_{w})$.   Comme dans la preuve de 6.7(7), on a $t_{0,w'}=\tau_{S}(t_{w})\in T^*_{sc}(E_{w'})$. On a aussi $\tau(u_{w})\in G_{\eta[d],SC}(E_{w'})$. Puisque $$g_{w'}=t_{0,w'}x\tau_{G^*}(x)^{-1}u_{{\cal E}^*}(\tau)\tau(r[d]_{sc})\tau(u_{w}),$$ la relation (2) \`a prouver \'equivaut \`a
$$x\tau_{G^*}(x)^{-1}u_{{\cal E}^*}(\tau)\tau(r[d]_{sc})\in T^*_{sc}(E_{w'})r[d]_{sc}G_{\eta[d],SC}(E_{w'}).$$
Mais cette relation r\'esulte de (1) et des conditions impos\'ees \`a $E$. Cela prouve (2) en g\'en\'eral. 

Fixons un ensemble fini $V'$ de places de $E$ v\'erifiant les conditions du paragraphe pr\'ec\'edent ainsi que les conditions suivantes.   On impose   que, pour tout $v\not\in V'$ et pour toute place $w'$ de $E$ divisant $v$, on ait  d'abord
   
   - $z[d]\in K_{sc,w'}$ et $r[d]_{sc}\in K_{sc,w'}$.

   Il en r\'esulte que $\eta[d]\in \tilde{K}_{w'}\cap \tilde{G}(F_{v})=\tilde{K}_{v}$. On sait qu'alors le groupe $K_{v}[d]=K_{v}\cap G_{\eta[d]}(F_{v})$ est un sous-groupe compact hypersp\'ecial de $G_{\eta[d]}(F_{v})$ (cf.  lemme 1.6). Il s'en d\'eduit un tel sous-groupe $K_{sc,v}[d]$ de $G_{\eta[d],SC}(F_{v})$ puis un  tel sous-groupe $K_{sc,w'}[d]$ de $G_{\eta[d],SC}(E_{w'})$. On impose de plus que
   
   - $\bar{x}[d]\in K_{sc,w'}[d]$ et $\bar{u}[d](\sigma)\in K_{sc,w'}[d]$ pour tout $\sigma\in \Gamma_{F}$.
   
   Toutes ces conditions impliquent que l'on a aussi $t(\sigma)\in T^*_{sc}(\mathfrak{o}_{w'})$ pour tout $\sigma\in \Gamma_{F}$, o\`u $t(\sigma)$ est l'\'el\'ement figurant dans (1).

   Soient $v\not\in V'$ et $w'$ une place de $E$ divisant $v$. Puisque $g_{w'}$ conjugue $S[d_{v}]_{sc}(E_{w'})$ en $T^*_{sc}(E_{w'})$ (o\`u $S[d_{v}]_{sc}$ est l'image r\'eciproque de $S[d_{v}]$ dans $G_{SC}$), la relation (2) 
    \'equivaut \`a $g_{w'}\in r[d]_{sc}G_{\eta[d],SC}(E_{w'})S[d_{v}]_{sc}(E_{w'})$. Ecrivons  $g_{w'}=r[d]_{sc}us$, avec $u\in G_{\eta[d],SC}(E_{w'})$ et $s\in S[d_{v}]_{sc}(E_{w'})$. Cela entra\^{\i}ne $us=r[d]_{sc}^{-1}g_{w'}\in K_{sc,w'}$. Appliquons l'op\'erateur $\theta=ad_{\eta[d]}$. Puisque $\eta[d]\in \tilde{K}_{v}$, $\theta$ conserve $K_{sc,w'}$. Donc $u\theta(s)=\theta(us)\in K_{sc,w'}$. D'o\`u $(1-\theta)(s)=(u\theta(s))^{-1}us\in K_{sc,w'}$. Fixons une uniformisante $\varpi$ de $F_{v}$, qui est aussi une uniformisante de $E_{w'}$ puisque $E_{w'}/F_{v}$ est non ramifi\'ee. Puisque $S[d_{v}]_{sc}$ est d\'eploy\'e sur $E_{w'}$, l'\'el\'ement $s$ s'\'ecrit de fa\c{c}on unique $s_{0}x_{*}(\varpi)$ pour un \'el\'ement $s_{0}\in S[d_{v}]_{sc}(\mathfrak{o}_{w'})$ et un \'el\'ement $x_{*}\in X_{*}(S[d_{v}]_{sc})$. On a alors $(1-\theta)(s)=(1-\theta)(s_{0})((1-\theta)(x_{*}))(\varpi)$. Puisque $S[d_{v}]_{sc}(\mathfrak{o}_{w'})\subset K_{sc,w'}$ d'apr\`es le (i) de la proposition 6.7, cela entra\^{\i}ne $((1-\theta)(x_{*}))(\varpi)\in K_{sc,w'}$. Cette condition ne peut \^etre r\'ealis\'ee que si $(1-\theta)(x_{*})=0$. Alors $x_{*}\in X_{*}(S[d_{v}]_{sc}^{\theta,0})$ et $x_{*}(\varpi)\in S[d_{v}]_{sc}^{\theta,0}(E_{w'})$. Ce groupe est contenu dans $G_{SC,\eta[d]}(E_{w'})$.
   
   {\bf Remarque.} Le groupe $G_{SC,\eta[d]}$ est la composante neutre du commutant de $\eta[d]$ dans $G_{SC}$; on prend garde de le distinguer du groupe $G_{\eta[d],SC}$ qui est le rev\^etement simplement connexe  du groupe d\'eriv\'e de $G_{\eta[d]}$, ou encore le rev\^etement simplement connexe du groupe d\'eriv\'e de $G_{SC,\eta[d]}$. 
   
   On a $ux_{*}(\varpi)=(us)s_{0}^{-1}\in K_{sc,w'}$. Donc $ux_{*}(\varpi)$ appartient au groupe $ G_{SC,\eta[d]}(E_{w'})\cap K_{sc,w'}$, qui est un sous-groupe compact hypersp\'ecial de $G_{SC,\eta[d]}(E_{w'})$. D'apr\`es 1.5(2), il existe $s_{1}\in S[d_{v}]_{sc}(\mathfrak{o}_{w'})$ tel que $ux_{*}(\varpi)$ et $s_{1}$ aient m\^eme image dans $G_{SC,\eta[d],ab}(E_{w'})$. Puisque $u\in G_{\eta[d],SC}(E_{w'})$, son image dans $G_{SC,\eta[d],ab}(E_{w'})$ est l'identit\'e. Cela entra\^{\i}ne que l'image de $x_{*}(\varpi)s_{1}^{-1}$ est aussi l'identit\'e. Donc $x_{*}(\varpi)s_{1}^{-1}$ est l'image d'un \'el\'ement  $s_{2}\in G_{\eta[d],SC}(E_{w'})$. En notant $\bar{S}[d_{v}]_{sc}$ l'image r\'eciproque de $S[d_{v}]$ dans $G_{\eta[d],SC}$, l'\'el\'ement $s_{2}$ appartient    forc\'ement \`a $\bar{S}[d_{v}]_{sc}(E_{w'})$. Comme plus haut, on peut \'ecrire $s_{2}=\bar{s}_{0}\bar{x}_{*}(\varpi)$ avec $\bar{s}_{0}\in \bar{S}[d_{v}]_{sc}(\mathfrak{o}_{w'})$ et $\bar{x}_{*}\in X_{*}(\bar{S}[d_{v}]_{sc})$. L'unicit\'e de ces d\'ecompositions entra\^{\i}ne que $x_{*}=\bar{x}_{*}$. Cela entra\^{\i}ne que  $x_{*}(\varpi)$ appartient \`a $\bar{S}[d_{v}]_{sc}(E_{w'})$  donc aussi \`a $G_{\eta[d],SC}(E_{w'})$ (plus exactement, c'est l'image dans $G_{SC}(E_{w'})$ d'un \'el\'ement de ce groupe). Posons $\underline{u}=ux_{*}(\varpi)$.  Alors $\underline{u}$ appartient \`a $G_{\eta[d],SC}(E_{w'})$ et son image dans $G_{SC}(E_{w'})$ appartient \`a $K_{sc,w'}$. Donc $\underline{u}\in   K_{sc,w'}[d]$. On a $g_{w'}=r[d]_{sc}\underline{u}s_{0}\in r[d]_{sc}K_{sc,w'}[d]S[d_{v}]_{sc}(\mathfrak{o}_{w'})$.  Comme plus haut, cela implique par conjugaison que $g_{w'}\in T^*_{sc}(\mathfrak{o}_{w'})r[d]_{sc}K_{sc,w'}[d]$. Cette relation est v\'erifi\'ee pour tout $v\not\in V'$ et tout $w'$ divisant $v$. Jointe \`a (2), elle entra\^{\i}ne le lemme. $\square$

 \bigskip

   \subsection{Un r\'esultat d'annulation}
 
  Soit $d_{V}=(d_{v})_{v\in V}\in D_{V}^{rel}$. On note $D_F[d_V]$ l'ensemble des $d\in D_{F}$ tels que
  
  - la projection de $d$ dans $D_{{\mathbb A}_{F}^V}$ appartient \`a $D_{{\mathbb A}_{F}^V}^{nr}$;
  
  - la projection de $d$ dans $D_{V}$ appartient \`a $I_{\eta}d_{V}G(F_{V})$.

  Pour $j\in {\cal J}({\bf H})$, on a d\'efini la constante $\delta_{j}[d_{V}]$ en 5.9.  
   
 \ass{Proposition}{Soit $d_{V}\in D_{V}^{rel}$. Si $D_F[d_V]=\emptyset$, alors   $$\sum_{j\in {\cal J}({\bf H})}\delta_{j}[d_{V}]=0.$$} 
 
 Preuve. Ce r\'esultat est trivial si ${\cal J}({\bf H})$ est vide. On suppose cet ensemble non vide. Alors l'application ${\bf p}$ de 6.4 l'identifie \`a $P({\bf H})$, qui est une unique classe modulo le sous-groupe $P^0\subset P$ de 6.6. Fixons $p\in P({\bf H})$. Soit $p^0\in P^0$. Posons $j={\bf p}^{-1}(p)$ et $j'={\bf p}^{-1}(p^0p)$. On va calculer le rapport $\delta_{j'}[d_{V}]\delta_{j}[d_{V}]^{-1}$.

 On dispose des paires de Borel $(B_{j},S_{j})$ de $G'_{j}$ et $(B_{j'},S_{j'})$ de $G'_{j'}$ dont les groupes de Borel sont d\'efinis sur $\bar{F}$ et les tores sont d\'efinis sur $F$. Par construction, on a des isomorphismes 
 $$X_{*,{\mathbb Q}}(S_{j})\simeq X_{*,{\mathbb Q}}(S_{\bar{H}})\oplus X_{*,{\mathbb Q}}(Z(\bar{G})^0)\simeq X_{*,{\mathbb Q}}(S_{j'})$$
 qui sont \'equivariants pour les actions galoisiennes. En fait, l'isomorphisme compos\'e provient d'un isomorphisme $S_{j}\simeq S_{j'}$ sur $\bar{F}$, qui est le compos\'e de $S_{j}\simeq T^*/(1-\theta^*)(T^*)\simeq S_{j'}$.  Puisque l'isomorphisme $X_{*,{\mathbb Q}}(S_{j})\simeq X_{*,{\mathbb Q}}(S_{j'})$ qui s'en d\'eduit fonctoriellement est \'equivariant pour les actions galoisiennes, l'isomorphisme $S_{j}\simeq S_{j'}$ est lui-m\^eme d\'efini sur $F$.
 
  On compl\`ete la donn\'ee $d_{V}$ en fixant un \'el\'ement $d^V=(d_{v})_{v\not\in V}\in D_{{\mathbb A}_{F}^V}^{nr}$ et en posant $d=d_{V}\times d^V=(d_{v})_{v\in Val(F)}$. On applique \`a cet \'el\'ement la proposition 6.7. On en d\'eduit des paires de Borel $(B[d_{v}],S[d_{v}])$ pour toute place $v\in Val(F)$ et des \'el\'ements $g\in G_{SC}({\mathbb A}_{E})$ et $t\in ((1-\theta^*)(T^*))({\mathbb A}_{E})$.  Pour toute place $v$, 
  les sextuplets 
 $$(\epsilon_{j},B_{j},S_{j},B[d_{v}],S[d_{v}],\eta[d_{v}])$$
 et
 $$(\epsilon_{j'},B_{j'},S_{j'},B[d_{v}],S[d_{v}],\eta[d_{v}])$$
  sont des diagrammes. 
 
 Pour $v\in V$, on fixe des \'el\'ements $\bar{Y}_{sc,v}\in \mathfrak{s}_{\bar{H}} (F_{v})$,  $Z_{1}\in \mathfrak{z}(\bar{G};F_{v})$ et $Z_{2}\in \mathfrak{z}(\bar{H};F_{v})$. On les suppose en position g\'en\'erale et proches de $0$.  On construit comme en 5.7 un \'el\'ement $X[d_{v}]\in \mathfrak{g}_{\eta[d_{v}]}(F_{v})$ dont on peut supposer qu'il appartient \`a $\mathfrak{s}[d_{v}]^{\theta}(F_{v})$. 
  On construit les \'el\'ements $Y_{j,v}$ et $Y_{j',v}$ comme en 5.7.   On peut supposer que $Y_{j,v}$ appartient \`a $ \mathfrak{s}_{j}(F_{v})$, plus pr\'ecis\'ement que $Y_{j,v}$ est l'image de $X[d_{v}]$ par l'homomorphisme $\mathfrak{s}[d_{v}](F_{v})\to  \mathfrak{s}_{j}(F_{v})$ provenant du premier diagramme ci-dessus. On peut supposer que $Y_{j',v} $ v\'erifie des propri\'et\'es analogues. Alors $Y_{j,v}$ et $Y_{j',v}$
  se correspondent via l'isomorphisme ci-dessus entre $S_{j}$ et $S_{j'}$.  On en fixe des rel\`evements $Y_{j,1,v}\in \mathfrak{g}'_{j,1,\epsilon_{j,1}}(F_{v})$ de $Y_{j,v}$ et $Y_{j',1,v}\in \mathfrak{g}'_{j',1,\epsilon_{j',1}}(F_{v})$ de $Y_{j',v}$, que l'on suppose proches de $0$ . On pose $x[d_{v}]=exp(X[d_{v}])$, $y_{j,v}=exp(Y_{j,v})$, $y_{j',v}=exp(Y_{j',v})$, $y_{j,1,v}=exp(Y_{j,1,v})$, $y_{j',1,v}=exp(Y_{j',1,v})$. Les sextuplets
$$(y_{j,v}\epsilon_{j},B_{j},S_{j},B[d_{v}],S[d_{v}],x[d_{v}]\eta[d_{v}])$$
 et
 $$(y_{j',v}\epsilon_{j'},B_{j'},S_{j'},B[d_{v}],S[d_{v}],x[d_{v}]\eta[d_{v}])$$
 sont des diagrammes. 
 
 Pour $v\not\in V$, on fixe des \'el\'ements $x[d_{v}]\in S[d_{v}]^{\theta,0}(F_{v})$, $y_{j,v}\in S_{j}(F_{v})$ et $y_{j',v}\in S_{j'}(F_{v})$ de sorte que les sextuplets ci-dessus soient encore des diagrammes. On impose que $x[d_{v}]\eta[d_{v}]$  est fortement r\'egulier. Pour presque tout $v$, $S[d_{v}]$ est non ramifi\'e et poss\`ede une structure naturelle sur $\mathfrak{o}_{v}$. On impose $x[d_{v}]\in S[d_{v}]^{\theta,0}(\mathfrak{o}_{v})$ pour presque tout $v$. D'apr\`es le (i) de la proposition 6.7, cela entra\^{\i}ne que $x[d_{v}]\eta[d_{v}]\in \tilde{K}_{v}$ pour presque tout $v$. Cela entra\^{\i}ne aussi que $y_{j,v}\in \tilde{K}'_{j,v}$ et $y_{j',v}\in \tilde{K}'_{j',v}$ pour presque tout $v$.  Enfin, on impose que $x[d_{v}]\eta[d_{v}]$ v\'erifie la condition [VI] 3.6(2) pour presque tout $v$. C'est loisible car, d'apr\`es la m\^eme preuve que celle de [VI] 3.6(15), on peut choisir pour presque tout $v$ un \'el\'ement $x[d_{v}]\in S[d_{v}]^{\theta,0}(\mathfrak{o}_{v})$ tel que cette condition soit v\'erifi\'ee. 
 On rel\`eve $y_{j,v}$ et $y_{j',v}$ en des \'el\'ements $y_{j,1,v}$ et $y_{j',1,v}$.  On suppose comme il est loisible que  $y_{j,1,v}\in \tilde{K}'_{j,1,v}$ et $y_{j',1,v}\in \tilde{K}'_{j',1,v}$ pour presque tout $v$. 
 
On supprime les indices $v$ pour noter les produits sur toutes les places de $F$. Par exemple $x[d]=\prod_{v\in Val(F)}x[d_{v}]$.  On remplace ces indices $v$ par $V$ pour noter les produits sur $v\in V$, par exemple $x[d]_{V}=\prod_{v\in V}x[d_{v}]$.  
 Par d\'efinition $\delta_{j'}[d_{V}]\delta_{j}[d_{V}]^{-1}$ est la limite quand les termes $Y_{j,1,V}$ et $Y_{j',1,V}$ tendent vers $0$ de
 $$\Delta_{j',1,V}(y_{j',1}\epsilon_{j',1},x[d_{V}]\eta[d_{V}])\Delta_{j,1,V}(y_{j,1}\epsilon_{j,1},x[d_{V}]\eta[d_{V}])
^{-1}.$$
Par d\'efinition, ce rapport est le produit du quotient des facteurs globaux de [VI] 3.6
$$(1) \qquad \Delta_{j',1}(y_{j',1}\epsilon_{j',1},x[d]\eta[d])\Delta_{j,1}(y_{j}\epsilon_{j,1},x[d]\eta[d])
^{-1}$$
et du produit pour toute place $v\not\in V$ des quotients des facteurs normalis\'es
$$(2) \qquad \Delta_{j',1,v}(y_{j',1}\epsilon_{j',1},x[d_{v}]\eta[d_{v}])^{-1}\Delta_{j,1,v}(y_{j,1}\epsilon_{j,1},x[d_{v}]\eta[d_{v}]).$$

On a

(3) quitte \`a remplacer en un nombre fini de places hors de $V$ les termes $x[d_{v}]$, $y_{j,1,v}$ et $y_{j',1,v}$ par d'autres termes assez proches des \'el\'ements neutres, le quotient (2) vaut $1$ pour tout $v\not\in V$.

Comme on l'a vu en [VI] 3.6, les deux termes de ce rapport valent $1$ en presque toute place, disons pour $v\not\in V''$, o\`u $V''$ est un ensemble fini de places contenant $V$. Pour $v\in V''-V$, rempla\c{c}ons les termes $x[d_{v}]$, $y_{j,1,v}$ et $y_{j',1,v}$ par d'autres termes assez proches des \'el\'ements neutres. Le th\'eor\`eme 5.7(i) dit que le rapport est alors \'egal \`a $\delta_{j}[d_{v}]\delta_{j'}[d_{v}]^{-1}$, ces termes \'etant d\'efinis comme dans ce paragraphe. Le (iii) du m\^eme th\'eor\`eme dit que ce rapport vaut $1$. D'o\`u (3).

On doit calculer le quotient (1), ou plus exactement sa limite quand les composantes $y_{j,1,v}$ et $y_{j',1,v}$ tendent vers $1$ pour tout $v\in V$ (dans la suite, on dira simplement "sa limite"). Comme ci-dessus, on peut s'autoriser \`a remplacer en un nombre fini de places hors de $V$ les termes $x[d_{v}]$, $y_{j,1,v}$ et $y_{j',1,v}$ par d'autres termes assez proches des \'el\'ements neutres. On se reporte \`a la d\'efinition de [VI] 3.6. L'automorphisme $ad_{g}$ envoie $S[d]({\mathbb A}_{E})$ sur $S({\mathbb A}_{E})$ de fa\c{c}on \'equivariante pour l'action de $Gal(E/F)$.  En particulier, il envoie $x[d]$ sur un \'el\'ement de $S({\mathbb A}_{F})$. que l'on note $x_{S}[d]=\prod_{v\in Val(F)}x_{S}[d_{v}]$.  D'autre part, on a  $ad_{g}(\eta[d])=t\eta$. D'o\`u $ad_{g}(x[d]\eta[d])=x_{S}[d]t\eta$.

Posons
$$\Delta_{j',j,II}[d]=\Delta_{j',II}(y_{j'}\epsilon_{j'},x[d]\eta[d])\Delta_{j,II}(y_{j}\epsilon_{j'},x[d]\eta[d])^{-1}.$$
Montrons que

(4)  quitte \`a remplacer en un nombre fini de places hors de $V$ les termes $x[d_{v}]$, $y_{j,1,v}$ et $y_{j',1,v}$ par d'autres termes assez proches des \'el\'ements neutres, on a
$$lim \,\Delta_{j',j,II}[d] =1.$$

   On note $\Sigma(S)_{res}$ l'ensemble des restrictions \`a $S^{\theta^*,0}=T^{*,\theta^*,0}$ des racines de $S$ dans $G$ et $\Sigma(S)_{res,ind}$ le sous-ensemble des \'el\'ements indivisibles.  Les termes  intervenant dans $\Delta_{j',j,II}[d]$ sont produits de termes index\'es par les orbites de $\Gamma_{F}$ dans l'ensemble $\Sigma(S)_{res,ind}$.  

Consid\'erons une orbite galoisienne d'un \'el\'ement $\alpha_{res}\in\Sigma(S)_{res,ind}$ qui se rel\`eve en une racine $\alpha\in \Sigma(S)$ de type $1$. Une telle racine contribue au num\'erateur comme au d\'enominateur de $\Delta_{j',j,II}[d]$ soit par $1$, soit par
$$(5) \qquad\chi_{\alpha_{res}}(\frac{(N\alpha)(x_{S}[d]\nu)-1}{a_{\alpha_{res}}})$$
(le terme $ t$ dispara\^{\i}t quand on applique $N\alpha$).
Si $(N\alpha)(\nu)\not=1$, on peut fixer un ensemble $V''$ de places contenant $V$ et assez grand pour que, pour $v\not\in V''$, on ait
$$\chi_{\alpha_{res},v}(\frac{(N\alpha)(x_{S}[d_{v}]\nu)-1}{a_{\alpha_{res}}})=\chi_{\alpha_{res},v}(\frac{(N\alpha)(\nu)-1}{a_{\alpha_{res}}})=1.$$
Pour $v\in V''-V$, quitte \`a remplacer nos donn\'ees $x[d_{v}]$ etc... par des termes assez proches des \'el\'ements neutres, on peut supposer que l'on a
$$\chi_{\alpha_{res},v}(\frac{(N\alpha)(x_{S}[d_{v}]\nu)-1}{a_{\alpha_{res}}})=\chi_{\alpha_{res},v}(\frac{(N\alpha)(\nu)-1}{a_{\alpha_{res}}}).$$
Pour $v\in V$, on a en tout cas
$$lim\,\chi_{\alpha_{res},v}(\frac{(N\alpha)(x_{S}[d_{v}]\nu)-1}{a_{\alpha_{res}}})=\chi_{\alpha_{res},v}(\frac{(N\alpha)(\nu)-1}{a_{\alpha_{res}}}).$$
Ainsi, la limite de l'expression (5) est $ \chi_{\alpha_{res}}(\frac{(N\alpha)(\nu)-1}{a_{\alpha_{res}}})$, qui vaut $1$  puisque les \'el\'ements qui interviennent appartiennent \`a $F_{\alpha}$ et que le caract\`ere $\chi_{\alpha_{res}}$ est  automorphe. Donc la racine $\alpha_{res}$ ne contribue pas   \`a $lim\,\,\Delta_{j',j,II}[d]$ si $(N\alpha)(\nu)\not=1$.
Supposons $(N\alpha)(\nu)=1$. Dans ce cas, $\alpha_{res}$ est une racine de $\bar{G}$. Mais alors, par construction, on a les \'egalit\'es $(N\hat{\alpha})(s_{j'})=(N\hat{\alpha})(s_{j})= (N\hat{\alpha})(\bar{s})$. Les conditions $(N\hat{\alpha})(s_{j'})=1$ et $(N\hat{\alpha})(s_{j})=1$ sont \'equivalentes. Donc $\alpha_{res}$ contribue au num\'erateur  de $\Delta_{j',j,II}[d]$ si et seulement si elle contribue au d\'enominateur. Quand elles contribuent, leur contributions sont toutes deux \'egales.  
Donc $\alpha_{res}$ ne contribue pas au quotient.  

Consid\'erons maintenant une orbite galoisienne d'un \'el\'ement $\alpha_{res}\in\Sigma(S)_{res,ind}$ qui se rel\`eve en  une racine $\alpha\in \Sigma(S)$ de type $2$. Une telle racine contribue au num\'erateur comme au d\'enominateur de $\Delta_{j',j,II}[d]$ soit par $1$, soit par $\chi_{\alpha_{res}}(\frac{(N\alpha)(x_{S}[d]\nu)^2-1}{a_{\alpha_{res}}})$, soit par $\chi_{\alpha_{res}}((N\alpha)(x_{S}[d]\nu)+1)$.  Si $(N\alpha)(\nu)\not=\pm 1$, le m\^eme raisonnement que ci-dessus montre que $\alpha_{res}$ ne contribue pas \`a $lim\,\,\Delta_{j',j,II}[d]$.  Supposons $(N\alpha)(\nu)=1$. Le m\^eme raisonnement montre que la troisi\`eme contribution possible est triviale et que l'on peut remplacer la seconde par $\chi_{\alpha_{res}}((N\alpha)(x_{S}[d]\nu)-1)$. La condition $(N\alpha)(\nu)=1$ signifie que $\alpha_{res}$ est une racine de $\bar{G}$ et  $2N\hat{\alpha}$ est une racine de $\hat{\bar{G}}$. On a alors les \'egalit\'es $(N\hat{\alpha})(s_{j'})^2=(N\hat{\alpha})(s_{j})^2= (N\hat{\alpha})(\bar{s})^2$. Les conditions $(N\hat{\alpha})(s_{j'})^2=1$ et $(N\hat{\alpha})^2(s_{j})=1$ sont \'equivalentes et de nouveau, $\alpha_{res}$ contribue de la m\^eme fa\c{c}on au num\'erateur et au d\'enominateur de $\Delta_{j',j,II}[d]$. Donc $\alpha_{res}$ ne contribue pas au quotient.  Supposons enfin $(N\alpha)(\nu)=-1$.  Le m\^eme raisonnement que ci-dessus montre que les deux contributions non triviales possibles sont toutes deux \'egales \`a $\chi_{\alpha_{res}}((N\alpha)(x_{S}[d]\nu)+1)$. Ce terme intervient au num\'erateur si et seulement si $(N\hat{\alpha})(s_{j'})\not=1$ et au d\'enominateur si  et seulement si $(N\hat{\alpha})(s_{j})\not=1 $. La condition $(N\alpha)(\nu)=-1$ signifie que $2\alpha_{res}$ est une racine de $\bar{G}$ et  $N\hat{\alpha}$ est une racine de $\hat{\bar{G}}$.  On a de nouveau $(N\hat{\alpha})(s_{j'})=(N\hat{\alpha})(s_{j})= (N\hat{\alpha})(\bar{s})$. Les conditions 
$(N\hat{\alpha})(s_{j'})\not=1$ et $(N\hat{\alpha})(s_{j})\not=1$ sont \'equivalentes et, comme pr\'ec\'edemment, $\alpha_{res}$ ne contribue pas \`a $\Delta_{j',j,II}[d]$. Cela prouve (4). 

La limite du quotient (1) est donc la m\^eme que celle du quotient
$$(6) \qquad \Delta_{j',imp}(y_{j',1}\epsilon_{j',1},x[d]\eta[d])\Delta_{j',imp}(y_{j}\epsilon_{j,1},x[d]\eta[d])^{-1}.$$

On utilise les d\'efinitions de [VI] 3.6.  On y remplace les $T$ par des $S$ et on ajoute des indices $j$ ou $j'$. Le d\'enominateur de (6) est l'inverse d'un produit $<h,\hat{h}>$ dans
 $$H^{1,0}({\mathbb A}_{F}/F;S_{sc}\stackrel{1-\theta}{\to}{\cal S}_{j,1})\times H^{1,0}(W_{F};\hat{{\cal S}}_{j,1}\stackrel{1-\hat{\theta}}{\to}\hat{S}_{ad}).$$
 L'\'el\'ement $h\in H^{1,0}({\mathbb A}_{F}/F;S_{sc}\stackrel{1-\theta}{\to}{\cal S}_{j,1})$ intervenant est un couple form\'e d'un cocycle $V_{S}$ \`a valeurs dans $S_{sc}({\mathbb A}_{E})/S_{sc}(E)$ et d'un \'el\'ement de ${\cal S}_{j,1}({\mathbb A}_{E})/{\cal S}_{j,1}(E)$. Rappelons la d\'efinition de $V_{S}$.   Avec les notations de 6.7, on a
 $$V_{S}(\sigma)=x\sigma_{G^*}(x)^{-1}u_{{\cal E}^*}(\sigma)\sigma(g)g^{-1}$$
 pour tout $\sigma\in \Gamma_{F}$. Comme on le voit, ce terme ne d\'epend pas de $j$. On note $e_{j}\in {\cal Z}(\tilde{G}'_{j};E)$ l'image naturelle de $e$ et on rel\`eve cet \'el\'ement en un \'el\'ement $e_{j,1}\in {\cal Z}(\tilde{G}'_{j,1};E)$. On \'ecrit $\epsilon_{j,1}=\mu_{j,1}e_{j,1}$. 
 L'\'el\'ement de ${\cal S}_{j,1}({\mathbb A}_{E})/{\cal S}_{j,1}(E)$ intervenant est l'image du couple $(x_{S}[d]t\nu,y_{j,1}\mu_{j,1})\in S({\mathbb A}_{E})\times S_{j,1}({\mathbb A}_{E})$ dans ${\cal S}_{j,1}({\mathbb A}_{E})/{\cal S}_{j,1}(E)$. Un calcul d\'ej\`a fait de nombreuses fois montre que la contribution de $(x_{S}[d],y_{j,1})$ est la valeur en ce point d'un   caract\`ere automorphe de ${\cal S}_{j,1}({\mathbb A}_{F})$. On peut fixer un ensemble fini de places $V''$ contenant $V$ tel que le caract\`ere vaille $1$ aux composantes hors de $V''$ de ce point, puis supposer les  composantes dans $V''-V$ assez proches de l'\'el\'ement neutre pour que le caract\`ere vaille aussi $1$ sur ces composantes.  Modulo ces modifications,  la limite de la valeur du caract\`ere quand les composantes dans $V$ tendent vers l'\'el\'ement neutre  vaut $1$. Donc, pour calculer notre limite, on peut remplacer le couple $(x_{S}[d]t\nu,y_{j',1}\mu_{j,1})$ par $(t\nu,\mu_{j,1})$. Mais $\nu$ et $\mu_{j,1}$ sont des \'el\'ements de $S(E)$ et $S_{j,1}(E)$. Ils disparaissent dans  ${\cal S}_{j,1}({\mathbb A}_{E})/{\cal S}_{j,1}(E)$ . Le terme $t$ est l'image de ce m\^eme terme consid\'er\'e comme un \'el\'ement de $((1-\theta)(S))({\mathbb A}_{E})/((1-\theta)(S))(E)$. On v\'erifie que le couple $(V_{S},t)$ d\'efinit un cocycle dans $Z^{1,0}({\mathbb A}_{F}/F;S_{sc}\stackrel{1-\theta}{\to}(1-\theta)(S))$. Notre \'el\'ement $h$ est l'image de ce cocycle par l'homomorphisme naturel
 $$Q=H^{1,0}({\mathbb A}_{F}/F;S_{sc}\stackrel{1-\theta}{\to}(1-\theta)(S))\to H^{1,0}({\mathbb A}_{F}/F;S_{sc}\stackrel{1-\theta}{\to}{\cal S}_{j,1}).$$
 Le tore dual de $(1-\theta)(S)$ est $\hat{S}/\hat{S}^{\hat{\theta},0}$.
 Par compatibilit\'e des produits, on a l'\'egalit\'e $<h,\hat{h}>=<(V_{S}, t),\hat{h}'>$, o\`u $\hat{h'}$ est l'image de $\hat{h}$ par l'homomorphisme dual
$$H^{1,0}(W_{F};\hat{{\cal S}}_{j,1}\stackrel{1-\hat{\theta}}{\to}\hat{S}_{ad})\to P= H^{1,0}(W_{F}; \hat{S}/\hat{S}^{\hat{\theta},0}\stackrel{1-\hat{\theta}}{\to}\hat{S}_{ad}).$$
Il suffit de comparer les d\'efinitions pour constater que $\hat{h}'$ n'est autre que  ${\bf p}(j)=p$. La limite du d\'enominateur de (6) est donc \'egale \`a $<(V_{S}, t),p>^{-1}$. Un m\^eme calcul s'applique au num\'erateur: sa limite est $<(V_{S}, t),p^0p>^{-1}$. On obtient que la limite de (6) est
$<(V_{S}, t),p^0>^{-1}$. Cela ach\`eve le calcul de notre rapport
 $$\delta_{j'}[d_{V}]\delta_{j}[d_{V}]^{-1}=<(V_{S},t),p^0>^{-1}.$$
 
 La somme de l'\'enonc\'e de la proposition est donc \'egale \`a
  $$\delta_{j}[d_{V}]\sum_{p^0\in P^0}<(V_{S}, t),p^0>^{-1}.$$
Cette somme est nulle si $(V_{S}, t)$ n'appartient pas \`a l'annulateur de $P^0$, c'est-\`a-dire \`a $Q_{0}$ d'apr\`es le lemme 6.6.  Pour \'etablir la proposition, il suffit donc de prouver que, si $(V_{S},t)$  appartient \`a $Q_{0}$, alors $D_F[d_V]\not=\emptyset$.  

On va d'abord prouver

(7) supposons que $(V_{S},t)$ appartienne \`a $Q_{0}$; alors il existe un \'el\'ement $d'$ v\'erifiant les m\^emes hypoth\`eses que $d$ et tel que son couple associ\'e $(V'_{S},t')$ appartienne \`a ${\bf q}_{1}(Q_{1})$.  

D'apr\`es le lemme 6.3, l'application $G({\mathbb A}_{F})\to Q_{2}$ est surjective. D'apr\`es 6.2(1), l'application naturelle $\prod_{v\not\in V}K_{\sharp,v}\to Q_{3}$ est surjective. L'hypoth\`ese signifie qu'il existe $\beta\in Q_{1}$, $h=(h_{v})_{v\in Val(F)}\in G({\mathbb A}_{F})$ et $k_{\sharp}=(k_{\sharp,v})_{v\not\in V}\in \prod_{v\not\in V}K_{\sharp,v}$ de sorte que $(V_{S},t)={\bf q}_{1}(\beta){\bf q}_{2}(h){\bf q}_{3}(k_{\sharp})$, o\`u on identifie $h$ et $k_{\sharp}$ \`a leurs images dans $Q_{2}$ et $Q_{3}$. Pour presque tout $v$, $h_{v}$ appartient \`a $K_{v}$. Alors l'image de $h_{v}$ par ${\bf q}_{2}$ est la m\^eme que l'image par ${\bf q}_{3}$ de l'image naturelle de $h_{v}$ dans $K_{\sharp,v}$. Quitte \`a modifier $k_{\sharp,v}$, on peut donc supposer $h_{v}=1$. On peut donc fixer un ensemble fini $V''$ de places de $F$, contenant $V$, tel que $h_{v}=1$ pour $v\not\in V''$. On peut imposer de plus 
  que $S$ est non ramifi\'e hors de $V''$.   Pour $v\not\in V''$, l'image de $K_{\sharp,v}$ dans $G_{\sharp,ab}(F_{v})$ co\"{\i}ncide avec celle de $(S[d_{v}]/Z(G)^{\theta})(\mathfrak{o}_{v})$ (cf. 1.5(2)). On peut donc supposer que $k_{\sharp,v}$ est un \'el\'ement  de $(S[d_{v}]/Z(G)^{\theta})(\mathfrak{o}_{v})$. Pour tout $v$, on rel\`eve $k_{\sharp,v}$ en un \'el\'ement $\underline{k}_{\sharp,v}\in \underline{K}_{\sharp,v}$ (on rappelle que ce groupe est l'image r\'eciproque de $K_{\sharp,v}$ dans $G(\bar{F}_{v})$). Pour unifier l'\'ecriture, on pose $\underline{k}_{\sharp,v}=1$ pour $v\in V$. Pour tout $v$, on
   d\'efinit l'\'el\'ement $d'_{v}=(\eta[d'_{v}],r[d'_{v}])$ par $\eta[d'_{v}]=ad_{(h_{v}\underline{k}_{\sharp,v})^{-1}}(\eta[d_{v}])$, $r[d'_{v}]=r[d_{v}]h_{v}\underline{k}_{\sharp,v}$. La famille $d'=(d'_{v})_{v\in Val(F)}$ appartient encore \`a $D_{{\mathbb A}_{F}}$ car, pour $v\not\in V''$,  $\eta[d'_{v}]=(1-\theta)(\underline{k}_{\sharp, v})\eta[d_{v}]\in ((1-\theta)(S[d_{v}]))(\mathfrak{o}_{v})\eta[d_{v}]$ et cet ensemble est contenu dans $ K_{v}\eta[d_{v}]$ pour  presque tout $v$ d'apr\`es le (i) de la proposition 6.7. Pour $v\not\in V$, $d_{v}$ appartient \`a $D_{v}^{nr}$ et le lemme 5.5 entra\^{\i}ne que $d'_{v}\in D_{v}^{nr}$. 
     La projection de $d'$  dans $D_{V}$ appartient \`a $d_{V}G(F_{V})$ puisque $d'_{v}=d_{v}h_{v}$ pour $v\in V$.   Pour tout $v\in Val(F)$, posons $(B[d'_{v}],S[d'_{v}])=ad_{(h_{v}\underline{k}_{\sharp,v})^{-1}}(B[d_{v}],S[d_{v}])$. Puisque $ad_{(h_{v}\underline{k}_{\sharp,v})^{-1}}$ est d\'efini sur $F_{v}$, cette paire v\'erifie pour $d'_{v}$ les m\^emes conditions que $(B[d_{v}],S[d_{v}])$ pour $d_{v}$. De plus, on a $(B[d'_{v}],S[d'_{v}])=(B[d_{v}],S[d_{v}])$ pour $v\not\in V''$ puisqu'alors  $h_{v}=1$ et $\underline{k}_{\sharp,v}\in S[d_{v}](\bar{F}_{v})$.  Soit $v$ une place de $F$, notons $w$ la restriction de $\bar{v}$ \`a $E$. Les deux paires $(B[d_{v}],S[d_{v}])$ et $(B[d'_{v}],S[d'_{v}])$ \'etant d\'eploy\'ees sur $E_{w}$, on peut fixer ${\bf g}_{w}\in G_{SC}(E_{w})$ tel que $ad_{{\bf g}_{w}}(B[d'_{v}],S[d'_{v}])=(B[d_{v}],S[d_{v}])$. Les deux paires \'etant \'egales pour $v\not\in V''$, on suppose ${\bf g}_{w}=1$ dans ce cas. Puisque $h_{v}\underline{k}_{\sharp,v}$ v\'erifie la m\^eme condition que ${\bf g}_{w}$, c'est-\`a-dire $ad_{h_{v}\underline{k}_{\sharp,v}}(B[d'_{v}],S[d'_{v}])=(B[d_{v}],S[d_{v}])$, il existe $c_{w}\in S[d_{v}]$ tel que $h_{v}\underline{k}_{\sharp,v}=c_{w}{\bf g}_{w}$. On a $ad_{{\bf g}_{w}}(\eta[d'_{v}])={\bf t}_{w}\eta[d_{v}]$, o\`u ${\bf t}_{w}=(\theta-1)(c_{w})$. Cette relation entra\^{\i}ne que  ${\bf t}_{w}\in ((1-\theta)(S[d_{v}]))(E_{w})$.  De plus, on a $c_{w}=\underline{k}_{\sharp,v}$ pour $v\not\in V''$, donc ${\bf t}_{w}\in ((1-\theta)(S[d_{v}]))(\mathfrak{o}_{v})$ dans ce cas. Pour une autre place $w'$ de $E$ au-dessus de $v$, on a fix\'e $\tau\in \Gamma_{F}$ tel que $\tau(w)=w'$, on pose ${\bf g}_{w'}=\tau({\bf g}_{w})$ et ${\bf t}_{w'}=\tau({\bf t}_{w})$. On pose ${\bf g}=({\bf g}_{w'})_{w'\in Val(E)}$ et ${\bf t}=({\bf t}_{w'})_{w'\in Val(E)}$.  On voit alors que les paires $(B[d'_{v}],S[d'_{v}])$ pour $v\in Val(F)$ et les \'el\'ements $g'=g{\bf g}$ et $t'=t\,ad_{g}({\bf t})$ satisfont les conditions de la proposition 6.7 pour l'\'el\'ement $d'$. 

On peut  reprendre nos constructions pour l'\'el\'ement $d'$ et ces donn\'ees auxiliaires. On affecte d'un $'$ les objets obtenus.
  On calcule
$$V'_{S}(\sigma)=V_{S}(\sigma)ad_{g}(\sigma({\bf g}){\bf g}^{-1})$$
pour tout $\sigma\in \Gamma_{F}$. Le couple $(V'_{S},t')$ est donc le produit de $(V_{S},t)$ et du cocycle $(ad_{g}(\chi),ad_{g}({\bf t}))$, o\`u $\chi(\sigma)=\sigma({\bf g}){\bf g}^{-1}$. Pour $v\in Val(F)$, on d\'efinit le cocycle localis\'e $(ad_{g_{w}}(\chi_{v}),ad_{g_{w}}({\bf t}_{w}))$, o\`u $\chi_{v}(\sigma)=\sigma({\bf g}_{w}){\bf g}_{w}^{-1}$ pour $\sigma\in \Gamma_{F_{v}}$. 
Alors $(ad_{g}(\chi),ad_{g}({\bf t}))$ est l'image naturelle dans $Q$ de  l'\'el\'ement de $H^{1,0}({\mathbb A}_{F};S_{sc}\stackrel{1-\theta}{\to}(1-\theta)(S))$ dont la composante en une place $v$ est  $(ad_{g_{w}}(\chi_{v}),ad_{g_{w}}({\bf t}_{w}))$. Ce dernier est un \'el\'ement de $H^{1,0}(F_{v};S_{sc}\stackrel{1-\theta}{\to}(1-\theta)(S))$. D'autre part, il r\'esulte des d\'efinitions que ${\bf q}_{2}(h_{v}){\bf q}_{3}(k_{\sharp,v})$ est l'image de $(h_{v},k_{\sharp,v})$ par la suite d'applications suivantes:

- on envoie $(h_{v},k_{\sharp,v})$ sur le produit des images de $h_{v}$ et $k_{\sharp,v}$ dans $G_{\sharp,ab}(F_{v})$;

- on envoie $G_{\sharp,ab}(F_{v})$ dans  $H^{1,0}(F_{v};S_{sc}\stackrel{1-\theta}{\to}(1-\theta)(S))$ par la suite d'applicatons
$$G_{\sharp,ab}(F_{v})\simeq H^{1,0}(F_{v}; S_{sc}[d_{v}]\to S[d_{v}]/Z(G)^{\theta,0})\to H^{1,0}(F_{v};S_{sc}[d_{v}]\stackrel{1-\theta}{\to}(1-\theta)(S[d_{v}]))$$
$$\stackrel{ad_{g_{w}}}{\simeq}H^{1,0}(F_{v};S_{sc}\stackrel{1-\theta}{\to}(1-\theta)(S));$$

- on envoie $H^{1,0}(F_{v};S_{sc}\stackrel{1-\theta}{\to}(1-\theta)(S))$ dans $Q$ par la m\^eme application que plus haut.

En se rappelant les \'egalit\'es $h_{v}\underline{k}_{\sharp,v}=c_{w}{\bf g}_{w}$ et  ${\bf t}_{w}=(\theta-1)(c_{w})$ et les d\'efinitions, on voit que l'image de $(h_{v},k_{\sharp,v})$ dans $H^{1,0}(F_{v}; S_{sc}[d_{v}]\to S[d_{v}]/Z(G)^{\theta,0})$ est $(\chi_{v}^{-1},c_{w})$, donc son image dans $H^{1,0}(F_{v};S_{sc}\stackrel{1-\theta}{\to}(1-\theta)(S))$ est $(ad_{g_{w}}(\chi_{v}^{-1}),ad_{g_{w}}({\bf t}_{w}^{-1}))$. On obtient que ${\bf q}_{2}(h_{v}){\bf q}_{3}(k_{\sharp,v})$ est l'inverse de l'image de $(ad_{g_{w}}(\chi_{v}),ad_{g_{w}}({\bf t}_{w}))$ dans $Q$. Donc
 $(V'_{S},t')$ est \'egal au produit de $(V_{S},t)$ et de  ${\bf q}_{2}(h)^{-1}{\bf q}_{3}(k_{\sharp})^{-1}$.
D'o\`u $(V'_{S},t')={\bf q}_{1}(\beta)$. Cela prouve (7).

Il reste \`a prouver

(8) supposons que $(V_{S},t)$ appartienne \`a ${\bf q}_{1}(Q_{1})$; alors $D_{F}[d_{V}]\not=\emptyset$.

Supposons $(V_{S},t)={\bf q}_{1}(\beta)$, o\`u $\beta\in Q_{1}=H^1({\mathbb A}_{F}/F;S_{\bar{H}})$. 
On rel\`eve $\beta$ en une cocha\^{\i}ne encore not\'ee $\beta$  \`a valeurs dans $S_{\bar{H}}({\mathbb A}_{\bar{F}})$, que l'on pousse en une cocha\^{\i}ne \`a valeurs dans $S_{sc}({\mathbb A}_{\bar{F}})$.  Il est maintenant plus commode d'identifier $S$ \`a $T^*$, muni de l'action galoisienne $\sigma\mapsto \sigma_{S}$. Notons que le tore  $S_{\bar{H}}$ s'identifie \`a l'image r\'eciproque  $\bar{T}_{sc}$ dans $\bar{G}_{SC}$ du sous-tore maximal $\bar{T}=T^{\theta^*,0}$ de $\bar{G}$. L'\'egalit\'e $(V_{S}, t)=q(\beta)$ signifie qu'il existe un \'el\'ement $t_{sc}\in T^*_{sc}({\mathbb A}_{\bar{F}})$ de sorte que l'on ait les relations
$$ (9) \qquad  V_{S}(\sigma)\sigma_{S}(t_{sc})^{-1}t_{sc}\in \beta(\sigma) T^*_{sc}(\bar{F})$$
pour tout $\sigma\in \Gamma_{F}$ et
$$tt_{sc}^{-1}\theta(t_{sc})\in ((1-\theta)(T^*))(\bar{F}).$$
 Fixons un tel \'el\'ement. On a $((1-\theta)(T^*))(\bar{F})=(1-\theta)(T^*(\bar{F}))$. Il existe donc $t'\in T^*(\bar{F})$ tel que $t^{-1}(1-\theta)(t_{sc})=(1-\theta)(t')$. On \'ecrit $t'=t'_{sc}z'$, avec $t'_{sc}\in T^*_{sc}(\bar{F})$ et $z'\in Z(G;\bar{F})$. On peut remplacer $t_{sc}$ par $t_{sc}(t'_{sc})^{-1}$ sans perturber la relation (9). On peut donc supposer
qu'il existe $z'\in Z(G;\bar{F})$ tel que $t^{-1}(1-\theta)(t_{sc})=(1-\theta)(z')$. On peut remplacer $g$ par $t_{sc}^{-1}g$. Cela ne perturbe pas les conditions impos\'ees \`a cet \'el\'ement. Evidemment, le nouvel \'el\'ement obtenu n'a pas de raison d'appartenir \`a $G_{SC}({\mathbb A}_{E})$. Mais on peut \`a ce point oublier cette propri\'et\'e et consid\'erer $g$ comme un \'el\'ement de $G_{SC}({\mathbb A}_{\bar{F}})$ (ou bien, on \'etend $E$ de sorte que $t_{sc}$ appartienne \`a $T_{sc}^*({\mathbb A}_{E})$). En modifiant ainsi $g$, le cocycle $V_{S}$ est remplac\'e par $\sigma\mapsto  V_{S}(\sigma)\sigma_{S}(t_{sc})^{-1}t_{sc}$ et l'\'el\'ement $t$ est remplac\'e par $t(\theta-1)(t_{sc})$. Donc, pour ce nouveau choix d'\'el\'ement $g$, on a
$$(10) \qquad V_{S}(\sigma)\in \beta(\sigma)T^*_{sc}(\bar{F})$$
pour tout $\sigma\in \Gamma_{F}$ et $t^{-1}=(1-\theta)(z')$, avec $z'\in Z(G;\bar{F})$. 

  Pour $\sigma\in \Gamma_{F}$, on pose
$$A(\sigma)=\beta(\sigma) V_{S}(\sigma)^{-1}\,\,\text{ et }\,\, X(\sigma)=A(\sigma)x\sigma_{G^*}(x)^{-1}=\beta(\sigma)g\sigma(g)^{-1}u_{{\cal E}^*}(\sigma)^{-1}.$$
D'apr\`es (10), $A$ est une cocha\^{\i}ne \`a valeurs dans $T^*_{sc}(\bar{F})$, donc $X$ prend ses valeurs dans $G_{sc}(\bar{F})$.   On calcule
 $$ ad_{X(\sigma)u_{{\cal E}^*}(\sigma)}\circ\sigma (\eta)=ad_{\beta(\sigma)g\sigma(g)^{-1}}\circ \sigma(\eta)=ad_{\beta(\sigma) g}(\sigma(g)^{-1}\sigma(\eta)\sigma(g))=ad_{\beta(\sigma)g} \circ\sigma(g^{-1}\eta g).$$
Par d\'efinition de $g$ et $t$, on a $g^{-1}t\eta g= \eta[d]$. Puisque $t=(\theta-1)(z')$ est central, cela implique $g^{-1}\eta g=t^{-1}\eta[d]$.
 L'\'el\'ement $\eta[d]$ est fixe par $\Gamma_{F}$. D'o\`u 
 $$ad_{\beta(\sigma)} \circ\sigma(g^{-1}\eta g)=ad_{\beta(\sigma)g}(\sigma(t^{-1})\eta[d])=\sigma(t)^{-1}ad_{\beta(\sigma)}\circ ad_{g}(\eta[d])$$
 $$=\sigma(t^{-1})ad_{\beta(\sigma)}(t\eta)=t\sigma(t)^{-1}ad_{\beta(\sigma)}(\eta)$$
 toujours parce que $t$ est central. Mais $\beta(\sigma)\in S_{\bar{H}}({\mathbb A}_{\bar{F}})$, a fortiori $\beta(\sigma)\in \bar{G}_{SC}({\mathbb A}_{\bar{F}})$ donc $ad_{\beta(\sigma)}$ fixe $\eta$. On obtient finalement
 $$ad_{X(\sigma)u_{{\cal E}^*}(\sigma)}\circ \sigma(\eta)=t\sigma(t)^{-1}\eta.$$ 
Puisque $t$ est central, on a l'\'egalit\'e $G_{ t\sigma(t)^{-1}\eta}=G_{\eta}=\bar{G}$. L'\'egalit\'e pr\'ec\'edente implique que l'automorphisme $ad_{X(\sigma)}\circ \sigma_{G^*}$ de $G(\bar{F})$ conserve $\bar{G}(\bar{F})$. Notons $\psi(\sigma)$ sa restriction \`a $\bar{G}$. On a introduit l'action quasi-d\'eploy\'ee $\sigma\mapsto \sigma_{\bar{G}}$ sur $\bar{G}$ relative \`a la paire de Borel \'epingl\'ee fix\'ee en 5.2. Rappelons  que l'on a
$$ad_{X(\sigma)}\circ \sigma_{G^*}=ad_{A(\sigma)}ad_{x\sigma_{G^*}(x)^{-1}}\circ \sigma_{G^*}.$$
Par d\'efinition de $x$,  $x\sigma_{G^*}(x)^{-1} $ normalise $T^*$ et son image dans le groupe $W$ est $\omega_{S}(\sigma)$. De plus, $A(\sigma)\in T^*_{sc}(\bar{F})$.   On en d\'eduit que $\psi(\sigma)$ conserve $\bar{T}=T^{\theta^*,0}$ et que sa restriction \`a $\bar{T}$ coincide avec $\omega_{S,\bar{G}}(\sigma)\sigma_{\bar{G}}$, o\`u $\omega_{S,\bar{G}}(\sigma)=\omega_{S}(\sigma)\omega_{\bar{G}}(\sigma)^{-1}$. On a $\omega_{S,\bar{G}}(\sigma)\in W^{\bar{G}}$ et le corollaire 2.2 de [K1] permet de fixer un \'el\'ement $\bar{x}\in \bar{G}_{SC}(\bar{F})$ tel que $\bar{x}\sigma_{\bar{G}}(\bar{x})^{-1}$ normalise $\bar{T}$ et que son image dans $W^{\bar{G}}$ soit $\omega_{S,\bar{G}}(\sigma)$. Alors les automorphismes $\psi(\sigma)$ et $ad_{\bar{x}\sigma_{\bar{G}}(\bar{x})^{-1}}\circ \sigma_{\bar{G}}$ de $\bar{G}(\bar{F})$ co\"{\i}ncident sur $\bar{T}$. Il existe donc $\bar{t}_{sc}(\sigma)\in \bar{T}_{sc}(\bar{F})$ tel que 
$$ad_{\bar{t}_{sc}(\sigma)}\circ\psi(\sigma)= ad_{\bar{x}\sigma_{\bar{G}}(\bar{x})^{-1}}\circ \sigma_{\bar{G}}.$$
L'\'el\'ement $\beta(\sigma)$ est un rel\`evement dans $\bar{T}_{sc}({\mathbb A}_{\bar{F}})$ de notre cocycle initial \`a valeurs dans $\bar{T}_{sc}({\mathbb A}_{\bar{F}})/\bar{T}_{sc}(\bar{F})$. 
On peut le multiplier par l'\'el\'ement   $\bar{t}_{sc}(\sigma)\in \bar{T}_{sc}(\bar{F})$. On voit que cela ne perturbe aucune des relations ci-dessus  mais cela remplace $\psi(\sigma)$ par $ad_{\bar{t}_{sc}(\sigma)}\circ\psi(\sigma)$. Apr\`es ce remplacement, on obtient simplement 
$$(11) \qquad \psi(\sigma)=ad_{\bar{x}\sigma_{\bar{G}}(\bar{x})^{-1}}\circ \sigma_{\bar{G}}.$$
En  utilisant les d\'efinitions, on en d\'eduit
$$(12)\qquad \sigma_{\bar{G}}=ad_{\sigma_{\bar{G}}(\bar{x})\bar{x}^{-1}}\circ  \psi(\sigma)=ad_{\sigma_{\bar{G}}(\bar{x})\bar{x}^{-1}}\circ  ad_{X(\sigma)}\circ \sigma_{G^*}$$
$$=ad_{\sigma_{\bar{G}}(\bar{x})\bar{x}^{-1}}\circ ad_{X(\sigma)u_{{\cal E}^*}(\sigma)}\circ \sigma=
ad_{\sigma_{\bar{G}}(\bar{x})\bar{x}^{-1}}\circ ad_{\beta(\sigma)g\sigma(g)^{-1}}\circ \sigma .$$
Il  r\'esulte de (11) que   $\psi$ est un homomorphisme de $\Gamma_{F}$ dans le groupe des automorphismes de $\bar{G}$ (il s'agit ici d'automorphismes de groupes abstraits). En revenant \`a la d\'efinition de $\psi$, on obtient que, pour $\sigma_{1},\sigma_{2}\in \Gamma_{F}$, les automorphismes
$$ad_{X(\sigma)_{1}}\circ\sigma_{1,G^*}\circ ad_{X(\sigma_{2})}\circ \sigma_{2,G^*}$$
et
$$ad_{X(\sigma_{1}\sigma_{2})}\circ (\sigma_{1}\sigma_{2})_{G^*}$$
co\"{\i}ncident sur $\bar{G}$. Cela \'equivaut \`a dire que $ X(\sigma_{1})\sigma_{1,G^*}(X(\sigma_{2}))X(\sigma_{1}\sigma_{2})^{-1}$ commute \`a $\bar{G}$. En utilisant la d\'efinition de la cocha\^{\i}ne $X$ et en se rappelant que $ad_{x\sigma_{G^*}(x)^{-1}}\circ\sigma_{G^*}=\sigma_{S}$ sur $T^*$, on calcule
$$ X(\sigma_{1})\sigma_{1,G^*}(X(\sigma_{2}))X(\sigma_{1}\sigma_{2})^{-1}=A(\sigma_{1})\sigma_{1,S}(A(\sigma_{2}))A(\sigma_{1}\sigma_{2})^{-1}=\partial_{S}(A)(\sigma_{1},\sigma_{2}).$$
On a not\'e $\partial_{S}$ la diff\'erentielle calcul\'ee pour l'action $\sigma\mapsto \sigma_{S}$. Donc $\partial_{S}(A)$ prend ses valeurs dans le commutant de $\bar{G}$. Par d\'efinition, on a
$$\partial_{S}(A)=\partial_{S}(\beta)\partial_{S}(V_{S}^{-1}).$$
On calcule ais\'ement $\partial_{S}(V_{S}^{-1})=\partial(u_{{\cal E}^*})^{-1}$, o\`u $\partial$ est la diff\'erentielle calcul\'ee pour l'action naturelle. On sait que $\partial(u_{{\cal E}^*})$ est \`a valeurs dans $Z(G_{SC})$, a fortiori dans le commutant de $\bar{G}$. Donc $\partial_{S}(\beta)$ prend aussi ses valeurs dans ce commutant. Puisque $\beta$ devient un cocycle quand on le pousse dans $\bar{T}_{sc}({\mathbb A}_{\bar{F}})/\bar{T}_{sc}(\bar{F})$, $\partial_{S}(\beta)$ est \`a valeurs dans $\bar{T}_{sc}(\bar{F})$. L'intersection de ce groupe avec le commutant de $\bar{G}$ est $Z(\bar{G}_{SC};\bar{F})$. Donc $\partial_{S}(\beta)$ est \`a valeurs dans ce dernier groupe. 

D\'efinissons une cocha\^{\i}ne $\bar{X}:\Gamma_{F}\to \bar{G}_{SC}({\mathbb A}_{\bar{F}})$ par
$\bar{X}(\sigma)=\beta(\sigma)^{-1}\bar{x}\sigma_{\bar{G}}(\bar{x})^{-1}$. Un calcul simple montre que $\partial_{\bar{G}}(\bar{X})=\partial_{S}(\beta)^{-1}$, o\`u $\partial_{\bar{G}}$ est la diff\'erentielle pour l'action $\sigma\mapsto \sigma_{\bar{G}}$. Donc $\bar{X}$ se pousse en un cocycle de $\Gamma_{F}$ dans $\bar{G}_{SC}({\mathbb A}_{\bar{F}})/Z(\bar{G}_{SC};\bar{F})$. Parce que $\bar{G}_{SC}$ est simplement connexe, le th\'eor\`eme 2.2 de [K2] dit que l'application naturelle
$$H^1(\Gamma_{F};\bar{G}_{SC}(\bar{F})/Z(G_{SC};\bar{F}))\to H^1(\Gamma_{F};\bar{G}_{SC}({\mathbb A}_{\bar{F}})/Z(\bar{G}_{SC};\bar{F}))$$
est surjective. Ainsi, on peut fixer $\bar{a}\in \bar{G}_{SC}({\mathbb A}_{\bar{F}})$ tel que 
$$ \bar{a}^{-1}\bar{X}(\sigma)\sigma_{\bar{G}}(\bar{a})\in \bar{G}_{SC}(\bar{F})$$
 pour tout $\sigma\in \Gamma_{F}$. 
 En utilisant la d\'efinition de $\bar{X}(\sigma)$ et la relation (12), cela \'equivaut \`a
 $$\bar{a}^{-1}\beta(\sigma)^{-1}\bar{x}\sigma_{\bar{G}}(\bar{x})^{-1}ad_{\sigma_{\bar{G}}(\bar{x})\bar{x}^{-1}\beta(\sigma)g\sigma(g)^{-1}}(\sigma(\bar{a}))\in \bar{G}_{SC}(\bar{F}),$$
 ou encore
 $$\bar{a}^{-1}g\sigma(\bar{a}^{-1}g)^{-1}\sigma(g)g^{-1}\beta(\sigma)^{-1}\bar{x}\sigma_{\bar{G}}(\bar{x})^{-1}\in \bar{G}_{SC}(\bar{F}).$$
 Puisque $\bar{x}\in \bar{G}_{SC}(\bar{F})$ et que $ \beta(\sigma)g\sigma(g)^{-1}=X(\sigma)u_{{\cal E}^*}(\sigma) $, cela \'equivaut aussi \`a
 $$(13) \qquad \bar{a}^{-1}g\sigma(\bar{a}^{-1}g)^{-1}\in \bar{G}_{SC}(\bar{F})X(\sigma)u_{{\cal E}^*}(\sigma).$$
A fortiori $\bar{a}^{-1}g\sigma(\bar{a}^{-1}g)^{-1}\in G_{SC}(\bar{F})$.  
  L'application $\sigma\mapsto \bar{a}^{-1}g\sigma(g^{-1}\bar{a})$ est un cocycle \`a valeurs dans $G_{SC}(\bar{F})$ qui est \'evidemment trivial quand on le pousse  dans $G_{SC}({\mathbb A}_{\bar{F}})$. D'apr\`es le th\'eor\`eme [Lab2] 1.6.9, il est trivial. On peut donc fixer $\dot{g}_{sc}\in G_{SC}(\bar{F})$ et $h\in G_{SC}({\mathbb A}_{F})$ tel que $g=\bar{a}\dot{g}_{sc}h$. Posons $\dot{g}=z'\dot{g}_{sc}$ et $\dot{d}=(\dot{g}^{-1}\eta\dot{g},\dot{g})$. Montrons que
  
 (14) $\dot{d}\in  D_F[d_V]$.
 
  Puisque $\dot{g}= z'\bar{a}^{-1}gh^{-1}$ et que $\bar{a}\in \bar{G}$, on a 
  $$\dot{g}^{-1}\eta\dot{g}=ad_{hg^{-1}}((\theta-1)(z')\eta)=ad_{hg^{-1}}(t\eta)=ad_{h}(\eta[d]).$$
  Donc $\dot{g}^{-1}\eta\dot{g}\in \tilde{G}({\mathbb A}_{F})$.  Puisque c'est aussi un \'el\'ement de $\tilde{G}(\bar{F})$, c'est un \'el\'ement de $\tilde{G}(F)$. Notons $\dot{\psi}$ la restriction de $ad_{\dot{g}}$ \`a $G_{\dot{g}^{-1}\eta\dot{g}}$, qui est un isomorphisme de ce groupe sur $\bar{G}$. Pour $\sigma\in \Gamma_{F}$, calculons
  $$\sigma(\dot{\psi})\circ\dot{\psi}^{-1}=\sigma_{\bar{G}}\circ \dot{\psi}\circ \sigma^{-1}\circ \dot{\psi}^{-1}.$$
  En utilisant (12), c'est la restriction \`a $G_{\dot{g}^{-1}\eta\dot{g}}$ de
  $$ad_{\sigma_{\bar{G}}(\bar{x})\bar{x}^{-1}}\circ ad_{X(\sigma)u_{{\cal E}^*}(\sigma)}\circ \sigma\circ ad_{\dot{g}}\circ\sigma^{-1}\circ ad_{\dot{g}^{-1}}=ad_{\sigma_{\bar{G}}(\bar{x})\bar{x}^{-1}}\circ ad_{X(\sigma)u_{{\cal E}^*}(\sigma)}\circ ad_{\sigma(\dot{g})\dot{g}^{-1}}.$$
  D'apr\`es les d\'efinitions, on a 
  $$ad_{\sigma(\dot{g})\dot{g}^{-1}}=ad_{\sigma(\dot{g}_{sc})\dot{g}_{sc}^{-1}}=ad_{\sigma(\bar{a}^{-1}g)g^{-1}\bar{a}}.$$
  Posons $\bar{u}(\sigma)=X(\sigma)u_{{\cal E}^*}(\sigma)\sigma(\bar{a}^{-1}g)g^{-1}\bar{a}$. On obtient que $\sigma(\dot{\psi})\circ\dot{\psi}^{-1}$ est la restriction \`a $\bar{G}$ de $ad_{\sigma_{\bar{G}}(\bar{x})\bar{x}^{-1}\bar{u}(\sigma)}$. Mais $\bar{x}\in \bar{G}_{SC}(\bar{F})$ et (13) montre que l'on a aussi $\bar{u}(\sigma)\in \bar{G}_{SC}(\bar{F})$.  Donc $\dot{\psi}$ est un torseur int\'erieur de $G_{\dot{g}^{-1}\eta\dot{g}}$ sur $\bar{G}$. Cela prouve que $\dot{d}\in D_{F}$. On peut \'etendre  le corps $E$ de sorte que tous les \'el\'ements intervenant appartiennent \`a $G({\mathbb A}_{E})$. Soit $v$ une place de $F$, notons $w$ la restriction de $\bar{v}$ \`a $E$. On se rappelle que $g_{w}\in T^*(\bar{F}_{\bar{v}})r[d_{v}]G_{\eta[d_{v}]}(\bar{F}_{\bar{v}})$. Ecrivons conform\'ement $g_{w}=t_{w}r[d_{v}]u_{w}$. Les \'egalit\'es $ad_{g_{w}}(\eta[d_{v}])=t\eta=(\theta^*-1)(z')\eta$ et $ad_{r[d_{v}]}(\eta[d_{v}])=\eta$ impliquent que $t_{w}\in (z')^{-1}T^{\theta^*} \subset (z')^{-1}I_{\eta}$. Il r\'esulte des d\'efinitions que 
  $$\dot{g}=\bar{a}_{w}^{-1}z'g_{w}h_{v}^{-1}=\bar{a}_{w}^{-1}z't_{w}r[d_{v}]u_{w}h_{v}^{-1}={\bf u}_{w}r[d_{v}]h_{v}^{-1},$$
  o\`u
  ${\bf u}_{w}=\bar{a}_{w}^{-1}z't_{w} ad_{r[d_{v}]}(u_{w})$. L'\'el\'ement ${\bf u}_{w}$ est un produit d'\'el\'ements de $I_{\eta}(\bar{F}_{\bar{v}})$. L'\'el\'ement $h_{v}$ appartient \`a $G(F_{v})$. Donc $\dot{g}\in I_{\eta}(\bar{F}_{\bar{v}})r[d_{v}]G(F_{v})$.  Les propri\'et\'es de $d_{v}$ se propagent donc \`a $\dot{d}$, c'est-\`a-dire que la projection de $\dot{d}$ dans $D_{{\mathbb A}_{F}^V}$ appartient \`a $D_{{\mathbb A}_{F}^V}^{nr}$ et sa projection dans $D_{V}$ appartient \`a $I_{\eta}d_{V}G(F_{V})$. Cela prouve (14).

  Evidemment, (14) d\'emontre (8), ce qui 
    ach\`eve la preuve de la proposition. $\square$

  {\bf Remarque.} On peut prouver la r\'eciproque, \`a savoir que, si $D_F[d_V]\not=\emptyset$,  alors $(V_{S},t)$ appartient \`a $Q_{0}$. Dans ce cas, l'application $j\mapsto \delta_{j}[d_{V}]$ est constante sur ${\cal J}({\bf H})$. Nous n'utiliserons pas ce r\'esultat, en le rempla\c{c}ant par la proposition du paragraphe suivant.
 
 \bigskip
 
 \subsection{Comparaison de deux facteurs de transfert}
Soit $d\in D_{F}$. On suppose

(1)  la projection de $d$ dans $D_{{\mathbb A}_{F}^V}$ appartient \`a $D_{{\mathbb A}_{F}^V}^{nr}$.

Notons $d_{V}$ la projection de $d$ dans $D_{V}$. On suppose

(2) $d_{V}\in D_{V}^{rel}$.

  L'\'el\'ement $\eta[d]$ appartient \`a $\tilde{G}(F)$, donc le groupe $G_{\eta[d]}$ est d\'efini sur $F$. Le torseur $\psi_{r[d]}$ est d\'efini sur $\bar{F}$.  La construction de [VI] 3.6 s'applique \`a la donn\'ee endoscopique $\bar{H}$ de $G_{\eta[d],SC}$ et fournit un facteur $\Delta[d_{V}]$ canonique, pourvu que l'on ait choisi en toute place $v\not\in V$ un sous-groupe compact hypersp\'ecial de $G_{\eta[d],SC}(F_{v})$. Pour cela, comme en 5.7, on fixe en toute place $v\not\in V$ un \'el\'ement $h_{v}\in G(F_{v})$ tel que $ad_{h_{v}^{-1}}(\eta[d])\in \tilde{K}_{v}$. C'est loisible d'apr\`es (1). On peut supposer que $h_{v}=1$ pour presque tout $v$. On choisit le sous-groupe image r\'eciproque dans $G_{\eta[d],SC}(F_{v})$ de $ad_{h_{v}}(K_{v})\cap G_{\eta[d]}(F_{v})$. On pose $h=(h_{v})_{v\not\in V}$ et on note plus pr\'ecis\'ement $\Delta[d_{V},h]$ le facteur de transfert canonique attach\'e \`a ce choix de compacts. Utilisons ce facteur dans la d\'efinition des constantes  de 5.7. On note $\delta_{j}[d_{V},h]$ le produit de ces constantes sur les places $v\in V$.

\ass{Proposition}{Pour tout $j\in {\cal J}({\bf H})$, on a l'\'egalit\'e $\delta_{j}[d_{V},h]=\omega(h)$.}

Preuve. On reprend les constructions de la preuve pr\'ec\'edente, en utilisant les notations des paragraphes 6.7 et 6.8. Pour toute place $v$, on note $\bar{S}[d_{v}]_{sc}$ l'image r\'eciproque de $S[d_{v}]^{\theta,0}$  dans $G_{\eta[d],SC}$. Pour $v\in V$, on pose $\bar{y}_{v}=exp(\bar{Y}_{sc,v})$ et $x_{sc}[d_{v}]=exp(X_{sc}[d_{v}])$ (cf. 5.7 pour ces notations; $d_{v}$ est l'image de $d$ dans $D_{v}$). Pour $v\not\in V$,  on a fix\'e dans la preuve pr\'ec\'edente un \'el\'ement $x[d_{v}]\in S[d_{v}]^{\theta,0}(F_{v})$ v\'erifiant diverses conditions. On v\'erifie facilement qu'on peut lui imposer de plus que son image dans $G_{\eta[d],AD}$ est l'image d'un \'el\'ement $x_{sc}[d_{v}]\in \bar{S}[d_{v}]_{sc}(F_{v})$. On fixe un tel \'el\'ement et on note $\bar{y}_{v}\in  S_{\bar{H}}(F_{v})$ l'\'el\'ement qui lui correspond par l'isomorphisme entre $ \bar{S}[d_{v}]_{sc}$ et $ S_{\bar{H}}$ d\'etermin\'e par les paires de Borel $(B_{\bar{H}},S_{\bar{H}})$ de $\bar{H}$ et $( \bar{B}[d_{v}]_{sc},\bar{S}[d_{v}]_{sc})$ de $G_{\eta[d],SC}$, o\`u $\bar{B}[d_{v}]_{sc}$ est l'image r\'eciproque dans ce groupe de $B[d_{v}]\cap G_{\eta[d]}$. Le terme $\delta_{j}[d_{V},h]$ est la limite quand  le terme $Y_{j,1,V}$ tend vers $0$ de 
$$\Delta_{j,1,V}(y_{j,1}\epsilon_{j,1},x[d_{V}]\eta)\Delta[d_{V},h](\bar{y}_{V},x_{sc}[d_{V}])^{-1}.$$
Comme dans le paragraphe pr\'ec\'edent, on note simplement $lim$ cette notion de limite.
Ce rapport est le produit des rapports des facteurs globaux et du produit sur les places $v\not\in V$ des rapports de facteurs normalis\'es
$$\Delta_{j,1,v}(y_{j,1}\epsilon_{j,1},x[d_{V}]\eta)^{-1}\Delta[d_{v},h](\bar{y}_{V},x_{sc}[d_{V}]).$$

Montrons que

(3) quitte \`a remplacer en un nombre fini de places hors de $V$ les termes $x[d_{v}]$, $x_{sc}[d_{v}]$ et $y_{j,1,v}$ par des \'el\'ements assez proches des \'el\'ements neutres, le  terme ci-dessus vaut $\omega(h_{v})$ pour tout $v\not\in V$.

On peut fixer un ensemble fini $V''$ de places contenant $V$ de sorte que le rapport vaille $1$ et $\omega(h_{v})=1$ pour tout $v\not\in V''$. Pour $v\in V''-V$, on remplace nos termes par des termes assez proches de $1$. Les assertions (i) et (iii) du th\'eor\`eme 5.8 disent qu'alors le rapport vaut $\omega(h_{v})$. D'o\`u (3).

On est ramen\'e \`a calculer le rapport des facteurs globaux. Il se d\'ecompose encore en produit du rapport des facteurs $\Delta_{II}$ et de celui des facteurs $\Delta_{imp}$. Montrons que

(4) quitte \`a remplacer en un nombre fini de places hors de $V$ les termes $x[d_{v}]$, $x_{sc}[d_{v}]$ et $y_{j,1,v}$ par des \'el\'ements assez proches des \'el\'ements neutres, on a
$$lim\, \prod_{v\in Val(F)}\Delta_{II,v}(y_{j}\epsilon_{j},x[d_{v}]\eta)\Delta_{II}[d_{v},h](\bar{y}_{v},x_{sc}[d_{v}])^{-1}=1.$$

Pour toute place $v$, les relations (1) et (2) de [W1] 10.3 nous disent que le rapport intervenant ci-dessus a une limite quand les donn\'ees $x[d_{v}]$ etc... tendent vers les \'el\'ements neutres et elles calculent cette limite. Notons-la $\ell_{v}$. En se reportant aux d\'efinitions de [W1] 9.3 et 10.3, on voit que, pour toute  \'el\'ement  $\alpha_{res}\in\Sigma(S)_{res,ind}$, il existe un \'el\'ement $c_{\alpha_{res}}\in F_{\alpha_{res}}^{\times}$ de sorte que

- pour $\sigma\in \Gamma_{F}$, $c_{\sigma(\alpha_{res})}=\sigma(c_{\alpha_{res}})$;

- pour tout $v$, $\ell_{v}=\prod_{\alpha_{res}\in\Sigma(S)_{res,ind}/\Gamma_{F_{v}}}\chi_{\alpha_{res},w}(c_{\alpha_{res}})$, o\`u $w$ est ici la restriction de $\bar{v}$ \`a  $F_{\alpha_{res}}$  et $\Sigma(S)_{res,ind}/\Gamma_{F_{v}}$ est le quotient de $\Sigma(S)_{res,ind}$ par l'action de $\Gamma_{F_{v}}$.

Les propri\'et\'es des $\chi$-data (qui sont bien s\^ur d\'efinies sur $F$) et des $c_{\alpha_{res}}$ permettent de r\'ecrire
$$\ell_{v}=\prod_{\alpha_{res}\in\Sigma(S)_{res,ind}/\Gamma_{F}}\prod_{w\vert v}\chi_{\alpha_{res},w}(c_{\alpha_{res}}),$$
o\`u le produit en $w$ est sur les places de $F_{\alpha_{res}}$ au-dessus de $v$. 

On peut fixer un ensenble fini $V''$ de places contenant $V$ de sorte que, pour $v\not\in V''$,

- $\Delta_{II,v}(y_{j}\epsilon_{j},x[d_{v}]\eta)\Delta_{II}[d_{v},h](\bar{y}_{v},x_{sc}[d_{v}])^{-1}=1$;

- $\ell_{v}=1$. 

Pour $v\in V''-V$, on remplace les donn\'ees $x[d_{v}]$ etc... par des \'el\'ements assez proches des \'el\'ements neutres pour que le rapport des facteurs $\Delta_{II,v}$ vaille $\ell_{v}$. Alors la limite de l'assertion (4) est $\prod_{v\in Val(F)}\ell_{v}$.  Mais ce produit vaut
$$\prod_{\Sigma(S)_{res,ind}/\Gamma_{F}}\prod_{w\in Val(F_{\alpha_{res}})}\chi_{\alpha_{res},w}(c_{\alpha_{res}}).$$
Cela vaut $1$ par la formule du produit. D'o\`u (4).

Il reste \`a calculer la limite du rapport des facteurs globaux
$$\Delta_{j,imp}(y_{j,1}\epsilon_{j,1},x[d]\eta)\Delta_{imp}[d](\bar{y},x_{sc}[d])^{-1}.$$
On a montr\'e dans la preuve de la proposition 6.9 que le premier terme avait une limite, \'egale \`a $<(V_{S}, t),{\bf p}(j)>^{-1}$ avec les notations de ce paragraphe. La m\^eme preuve montre que le second terme a aussi une limite, \'egale \`a $<V_{\bar{S}},\bar{s}>^{-1}$, o\`u $V_{\bar{S}}:\Gamma_{F}\to  S_{\bar{H}}({\mathbb A}_{\bar{F}})/S_{\bar{H}}(\bar{F})$ est un cocycle analogue \`a $V_{S}$. Il reste \`a prouver l'\'egalit\'e
$$ (5) \qquad <(V_{S}, t),j>=<V_{\bar{S}},\bar{s}>.$$
 Reprenons la construction de $(V_{S},t)$ du paragraphe pr\'ec\'edent. Ici, l'\'el\'ement $d$ appartient \`a $ D_{F}$  et on peut appliquer le lemme 6.8.  On peut donc  \'ecrire $g=t_{1}r[d]_{sc}u$, avec $t_{1}\in T^*_{sc}({\mathbb A}_{E})$ et $u\in G_{\eta[d],SC}({\mathbb A}_{E})$. Cela entra\^{\i}ne $g=z[d]^{-1}t_{1}r[d]u$. L'\'el\'ement $t\in ((1-\theta^*)(T^*))({\mathbb A}_{E})$ a \'et\'e d\'efini par $ad_{g}(\eta[d])=t\eta$. La relation pr\'ec\'edente entra\^{\i}ne $t=(1-\theta^*)(t_{1}z[d]^{-1})$. Pour $\sigma\in \Gamma_{F}$, on a
 $$V_{S}(\sigma)=x\sigma_{G^*}(x)^{-1}u_{{\cal E}^*}(\sigma)\sigma(g)g^{-1}$$
 $$=\sigma_{S}(t_{1})x\sigma_{G^*}(x)^{-1}u_{{\cal E}^*}(\sigma)\sigma(r[d]_{sc})r[d]_{sc}^{-1}ad_{r[d]_{sc}}(\sigma(u)u^{-1})t_{1}^{-1}= \sigma_{S}(t_{1})t_{1}^{-1}V'_{S}(\sigma),$$
 o\`u
 $$V'_{S}(\sigma)=x\sigma_{G^*}(x)^{-1}u_{{\cal E}^*}(\sigma)\sigma(r[d]_{sc})r[d]_{sc}^{-1}ad_{r[d]}(\sigma(u)u^{-1}).$$
 Le couple form\'e du cocycle $\sigma\mapsto \sigma_{S}(t_{1})t_{1}^{-1}$ et de l'\'el\'ement $(1-\theta^*)(t_{1})$ est un cobord. On peut le supprimer.  Le terme $(1-\theta^*)(z[d])^{-1}$ appartient \`a $(1-\theta^*)(S)(E)$ et dispara\^{\i}t dans notre groupe de cohomologie $Q$. On obtient que $(V_{S},t)$ (ou plus exactement son image dans $Q$ est l'image par l'homomorphisme naturel
 $$H^{1,0}({\mathbb A}_{\bar{F}}/\bar{F}; S_{sc})\to Q= H^{1,0}({\mathbb A}_{\bar{F}}/\bar{F}; S_{sc}\stackrel{1-\theta}{\to}(1-\theta)(S))$$
 du cocycle $V'_{S}$.

 Pour  calculer $V_{\bar{S}}$,  on utilise  les objets d\'efinis en 6.8.  Posons $S[d]=\prod_{v\in Val(F)}S[d_{v}]$ et $\bar{S}[d]=\prod_{v\in Val(F)}(S[d_{v}]\cap G_{\eta[d]})$. Notons aussi $\bar{S}^*[d]$ le tore $\bar{T}[d]$ muni de l'action $\sigma\mapsto \sigma_{S}$ transport\'ee par $ad_{r[d]^{-1}}$. Notons $\bar{S}[d]_{sc}$ et $\bar{S}^*[d]_{sc}$ les images r\'eciproques de $\bar{S}[d]$ et $\bar{S}^*[d]$ dans $G_{\eta[d],SC}$. Puisque $ad_{g}$ se restreint en un isomorphisme d\'efini sur $F$ de $S[d]$ sur $S$, $ad_{u}$ se restreint en un isomorphisme d\'efini sur $F$ de $\bar{S}[d]$ sur $\bar{S}^*[d]$. La d\'efinition de [VI] 3.6 fournit un cocycle $V'_{\bar{S}}$ \`a valeurs dans $\bar{S}^*[d]_{sc}({\mathbb A}_{E})/\bar{S}^*[d]_{sc}(E)$ d\'efini par
 $$V'_{\bar{S}}(\sigma)=\bar{x}[d]\sigma_{G_{\eta[d]^*}}(\bar{x}[d])^{-1}\bar{u}[d](\sigma)\sigma(u)u^{-1}.$$
 Le cocycle $V_{\bar{S}}$ est l'image de $V'_{\bar{S}}$ par l'isomorphisme $ad_{r[d]}:\bar{S}^*[d]_{sc}\to \bar{S}_{sc}=S_{\bar{H}}$. On voit sur les formules ci-dessus que $V'_{S}(\sigma)$ est le produit de l'image naturelle de $V_{\bar{S}}(\sigma)$ et de termes appartenant \`a $S_{sc}(E)$. Il en r\'esulte que $(V_{S},t)$ est l'image de $V_{\bar{S}}$ par l'homomorphisme
 $${\bf q}_{1}:H^1({\mathbb A}_{\bar{F}}/\bar{F};S_{\bar{H}})\to Q= H^{1,0}({\mathbb A}_{\bar{F}}/\bar{F}; S_{sc}\stackrel{1-\theta}{\to}(1-\theta)(S)).$$
 Puisque l'application
 $${\bf p}_{1}:P=H^{1,0}(W_{F}; \hat{S}/\hat{S}^{\hat{\theta},0}\stackrel{1-\hat{\theta}}{\to}\hat{S}_{ad})\to \hat{S}_{\bar{H}}^{\Gamma_{F}}$$
 est duale de ${\bf q}_{1}$, cela entra\^{\i}ne que
 $$<(V_{S},t),{\bf p}(j)>=<V_{\bar{S}},{\bf p}_{1}\circ{\bf p}(j)>.$$
 Mais ${\bf p}_{1}\circ{\bf p}(j)=\bar{s}$ d'apr\`es la proposition 6.4. Cela prouve (5) et la proposition. $\square$

\bigskip

\section{Le cas o\`u $D_{F}[d_{V}]$ est non vide}

\bigskip

\subsection{Une proposition de nullit\'e}
Dans cette section, on fixe un \'el\'ement $d_{V}\in D_{V}^{rel}$. 
On a

(1) l'ensemble de classes $I_{\eta}(\bar{F})\backslash D_F[d_V]/G(F)$ est fini.

Preuve. Quand $d$ parcourt $ D_F[d_V]$, les $\eta[d]$ sont des \'el\'ements de $\tilde{G}(F)$ qui appartiennent \`a la m\^eme classe de conjugaison stable. Pour $v\not\in V$, leurs classes de conjugaison par $G(F_{v})$ coupent $\tilde{K}_{v}$. D'apr\`es le lemme [VI] 2.1, ils sont contenus dans un nombre fini de classes de conjugaison par $G(F)$. Fixons un sous-ensemble fini $X_{0}\subset D_F[d_V]$ tel que, pour tout $d\in D_F[d_V]$, il existe $d_{0}\in X_{0}$ tel que $\eta[d]$ et $\eta[d_{0}]$ soient conjugu\'es par un \'el\'ement de $G(F)$. L'ensemble $ I_{\eta}(\bar{F})\backslash Z_{G}(\eta;\bar{F})$ est fini. Fixons-en un ensemble de repr\'esentants $Z_{0}$. Pour $d\in D_F[d_V]$ et $z\in Z_{0}$, l'\'el\'ement  $zd$ peut appartenir ou non \`a $D_F[d_V]$. Notons $X_{1}$ l'ensemble des $ zd_{0}$ qui appartiennent \`a $D_F[d_V]$, pour $d_{0}\in X_{0}$ et $z\in Z_{0}$. Soit $d\in D_F[d_V]$. On peut choisir $d_{0}\in X_{0}$ et $x\in G(F)$ tels que $\eta[d]=x\eta[d_{0}]x^{-1}$. Alors les deux \'el\'ements $r[d_{0}]$ et $r[d]x$ conjuguent $\eta[d_{0}]$ en $\eta$. Ils diff\`erent donc par multiplication  \`a gauche par un \'el\'ement de  $Z_{G}(\eta;\bar{F})$, que l'on peut \'ecrire $u^{-1}z$, avec $z\in Z_{0}$ et $u\in I_{\eta}(\bar{F})$. On a alors $ zd_{0}= udx$. L'\'el\'ement de droite appartient \`a $D_F[d_V]$ donc aussi celui de gauche. Ce dernier appartient donc \`a $X_{1}$. L'\'egalit\'e pr\'ec\'edente montre que cet \'el\'ement a m\^eme image que $d$ dans notre ensemble de doubles classes. Cela montre que tout \'el\'ement de cet ensemble de doubles classes est l'image d'un \'el\'ement de $X_{1}$, lequel est fini. $\square$

 Fixons un ensemble de repr\'esentants $\dot{D}_{F}[d_{V}]$ de l'ensemble de doubles classes 
 
\noindent $I_{\eta}(\bar{F})\backslash D_F[d_V]/G(F)$. Consid\'erons un ensemble fini $V'$ de places de $F$ contenant $V$ et tel que

(2) pour tout $d\in \dot{D}_{F}[d_{V}]$ et tout $v\not\in V'$, on a $\eta[d]\in \tilde{K}_{v}$. 

Soit $d\in \dot{D}_{F}[d_{V}]$. On fixe un ensemble de repr\'esentants ${\cal U}[V',d]$ de l'ensemble de  classes
$$G_{\eta[d]}(F_{V'}^V)/Z(G_{\eta[d]};F_{V'}^V)G_{\eta[d],SC}(F_{V'}^V).$$
On suppose que $1\in {\cal U}[V',d]$.
Pour $v\in V'-V$, on fixe $h_{v}[d]\in G(F_{v})$ tel que $ad_{h_{v}[d]^{-1}}(\eta[d])\in \tilde{K}_{v}$. On pose $h[d]=(h_{v}[d])_{v\in V'-V}$. Pour $u=(u_{v})_{v\in V'-V}\in {\cal U}[V',d]$, munissons $G_{\eta[d],SC}({\mathbb A}_{F}^V)$ du sous-groupe compact $K_{sc}^V[d,u]$ d\'efini ainsi: pour tout $v\not\in V$, $K_{sc,v}[d,u]$ est l'image r\'eciproque dans $G_{\eta[d],SC}(F_{v})$ de $K_{v}\cap G_{\eta[d]}(F_{v})$ si $v\not\in V'$ et de $ad_{u_{v}h_{v}[d]}(K_{v})\cap G_{\eta[d]}(F_{v})$ si $v\in V'-V$. Ce groupe permet de d\'efinir un facteur de transfert canonique sur
$$\bar{H}(F_{V})\times G_{\eta[d],SC}(F_{V})$$
 que nous notons  $\Delta_{V}[d,u]$. 
 
En 5.9, on a d\'efini une fonction $\bar{f}[d_{V}]\in SI(\bar{H}(F_{V}))$.  Pour $d\in \dot{D}_{F}[d_{V}]$, la "composante en $V$" de $d$ n'est pas forc\'ement $d_{V}$, cette composante appartient seulement \`a $I_{\eta}(\bar{F}_{V})d_{V}G(F_{V})$. On peut n\'eanmoins appliquer la d\'efinition de 5.9 en rempla\c{c}ant $d_{V}$ par cette composante en $V$. Dans cette d\'efinition,  il intervient un facteur de transfert. On peut prendre $\Delta_{V}[d,u]$ pour un \'el\'ement $u\in {\cal U}[V',d]$. On note alors $\bar{f}[d,u]$ la fonction obtenue. Posons
$$(3) \qquad \bar{\varphi}[V',d_{V}]=  \sum_{d\in \dot{D}_{F}[d_{V}]}\vert {\cal U}[V',d]\vert ^{-1}   \sum_{u\in {\cal U}[V',d]}\omega(uh[d]) \bar{f}[d,u].$$
Cette expression d\'epend de plusieurs donn\'ees auxiliaires, en particulier de l'ensemble de places $V'$ puisque  les ensembles ${\cal U}[V',d]$ en d\'ependent.
Les paragraphes 7.3 \`a 7.9 sont consacr\'es \`a prouver la proposition suivante.

\ass{Proposition}{Supposons ${\cal J}({\bf H})=\emptyset$. Alors, si $V'$ est assez grand, on a $\bar{\varphi}[V',d_{V}]=0$.}

\bigskip

\subsection{Premier calcul d'une expression intervenant en 5.9}
Consid\'erons l'expression 

(1) $  \bar{f}[d_{V}]\sum_{j\in {\cal J}({\bf H})}\delta_{j}[d_{V}]$

\noindent qui appara\^{\i}t dans l'expression 5.9(2).

\ass{Corollaire}{L'expression (1) est nulle si $D_{F}[d_{V}]=\emptyset$. Supposons $D_{F}[d_{V}]$ non vide. Alors cette expression (1) est \'egale \`a $\vert P^0\vert \vert \dot{D}_{F}[d_{V}]\vert ^{-1} \bar{\varphi}[V',d_{V}]$ pourvu que $V'$ soit assez grand. }

Preuve. La premi\`ere assertion r\'esulte de la proposition 6.9. Supposons $D_{F}[d_{V}]$ non vide. Si ${\cal J}({\bf H})$ est vide, l'expression (1) est nulle et $\bar{\varphi}[V',d_{V}]$ aussi d'apr\`es la proposition 7.1. Supposons ${\cal J}({\bf H})$ non vide. L'expression  (1) ne d\'epend que de la double classe $I_{\eta}(\bar{F}_{V})d_{V}G(F_{V})$. On peut remplacer ce terme $d_{V}$ par la composante en $V$ d'un \'el\'ement $d\in \dot{D}_{F}[d_{V}]$.  On peut utiliser dans les d\'efinitions le facteur de transfert $\Delta_{V}[d,u]$ pour un \'el\'ement $u\in {\cal U}[V',d]$. La fonction $\bar{f}[d_{V}]$ est alors remplac\'ee par $\bar{f}[d,u]$.  La proposition 6.10 calcule les facteurs $\delta_{j}[d_{V}]$ pour nos choix: ils sont \'egaux \`a $\omega(uh[d])$, o\`u $u$ et $h[d]$ sont compl\'et\'es en des \'el\'ements de $G({\mathbb A}_{F})$ de composantes $1$ hors de $V'-V$. 
L'expression (1) devient
  $$\vert {\cal J}({\bf H})\vert  \omega(uh[d]) \bar{f}[d,u].$$
  On peut moyenner en $d$ et $u$ et on obtient l'expression
  $$ \vert {\cal J}({\bf H})\vert \vert \dot{D}_{F}[d_{V}]\vert ^{-1}  \bar{\varphi}[V',d_{V}].$$ 
D'apr\`es 6.4 et 6.6, on a $\vert {\cal J}({\bf H})\vert =\vert P^0\vert $ puisque ${\cal J}({\bf H})\not=\emptyset$. D'o\`u la deuxi\`eme assertion du corollaire. $\square$.

\bigskip

 \subsection{Mise en place de la situation}
 
La conclusion de la proposition 7.1 est triviale si $D_{F}[d_{V}]$ est vide. {\bf On suppose  jusqu'en 7.15  que $D_{F}[d_{V}]\not=\emptyset$}. 
On fixe un \'el\'ement de cet ensemble, que l'on note $d_{\star}=(\eta_{\star},r_{\star})$. Pour simplifier les notations, on pose simplement
 $$I_{\star}=I_{\eta_{\star}},\,\, \bar{G}_{\star}=G_{\eta_{\star}}.$$
 Ces groupes sont d\'efinis sur $F$ et $ad_{r_{\star}}$ se restreint en un torseur int\'erieur de $\bar{G}_{\star}$ sur $\bar{G}$. Fixons un sous-tore maximal $\bar{T}_{\star}$ de $\bar{G}_{\star}$ d\'efini sur $F$. Notons $T_{\star}$ son commutant dans $G$. On peut fixer $\bar{r}_{\star}\in \bar{G}_{\star}(\bar{F})$ tel que $ad_{\bar{r}_{\star}}(\bar{T}_{\star})=ad_{r_{\star}^{-1}}(T^*)\cap \bar{G}_{\star}$.  On a $ad_{r_{\star}\bar{r}_{\star}}(\bar{T}_{\star})=\bar{T}$ ($=\bar{G}\cap T^*$) et  $ad_{r_{\star}\bar{r}_{\star}}(T_{\star})=T^*$. Il existe  un cocycle $\omega_{T_{\star}}:\Gamma_{F}\to W^{\bar{G}}$ tel que $ad_{r_{\star}\bar{r}_{\star}}$ entrelace l'action galoisienne naturelle sur $\bar{T}_{\star}$ avec l'action $\sigma\mapsto \omega_{T_{\star}}(\sigma)\circ \sigma_{\bar{G}}$ sur $\bar{T}$. Puisque $\sigma_{\bar{G}}=\omega_{\bar{G}}(\sigma)\circ \sigma_{G^*}$, $ad_{r_{\star}\bar{r}_{\star}}$ entrelace l'action naturelle sur $T_{\star}$ avec l'action $\sigma\mapsto \omega_{T_{\star}}(\sigma)\omega_{\bar{G}}(\sigma)\circ \sigma_{G^*}$ sur $T^*$. On a d\'efini le tore $S$ comme \'etant $T^*$ muni de l'action galoisienne $\sigma\mapsto \omega_{S}(\sigma)\circ \sigma_{G^*}$ et on a $\omega_{S}(\sigma)=\omega_{S,\bar{G}}(\sigma)\omega_{\bar{G}}(\sigma)$, o\`u $\omega_{S,\bar{G}}(\sigma)\in W^{\bar{G}}$. Par l'isomorphisme $ad_{r_{\star}\bar{r}_{\star}}$, on identifie $W^{\bar{G}}$ au groupe de Weyl $W^{\bar{G}_{\star}}$ de $\bar{G}_{\star}$ relatif \`a $\bar{T}_{\star}$. On d\'efinit $\omega_{S,T_{\star}}(\sigma)=\omega_{S,\bar{G}}(\sigma)\omega_{T_{\star}}(\sigma)^{-1}$. C'est un \'el\'ement de $W^{\bar{G}_{\star}}$ et  $S$ s'identifie au tore $T_{\star}$ muni de l'action $\sigma\mapsto  \omega_{S,T_{\star}}(\sigma)\circ \sigma$ (ce dernier $\sigma$ \'etant l'action naturelle sur $T_{\star}$). Dans cette section, on identifie ainsi $S$ \`a $T_{\star}$. On note $\bar{S}=S\cap \bar{G}_{\star}$, c'est-\`a-dire $\bar{S}=\bar{T}_{\star}$ muni de l'action $\sigma\mapsto  \omega_{S,T_{\star}}(\sigma)\circ \sigma$. On note $\bar{S}_{sc}$ l'image r\'eciproque de $\bar{S}$ dans $\bar{G}_{\star,SC}$. Ce tore s'identifie \`a notre pr\'ec\'edent tore $S_{\bar{H}}$. On note aussi $S_{sc}$ l'image r\'eciproque de $S$ dans $G_{SC}$. 
 
 En appliquant les constructions de 6.7 \`a notre \'el\'ement $d_{\star}=(\eta_{\star},r_{\star})$, on fixe pour toute place $v\in Val(F)$ une paire de Borel $(B[d_{\star,v}],S[d_{\star,v}])$  v\'erifiant toutes les propri\'et\'es de ce paragraphe. On note simplement $S_{\star,v}=S[d_{\star,v}]$. Pour une extension galoisienne finie $E$ de $F$, on   note $S_{ \star}({\mathbb A}_{E})$ le produit restreint des $S_{\star,v}(E_{w'})$ sur les places $v$ de $F$ et sur les places $w'$ de $E$ divisant $v$. La restriction est relative aux sous-groupes $S_{\star,v}(\mathfrak{o}_{w'})$ qui se d\'efinissent naturellement en presque tout couple $(v,w')$ de places. On note $S_{\star}({\mathbb A}_{\bar{F}})$ la limite inductive de $S_{ \star}({\mathbb A}_{E})$ quand $E$ parcourt les extensions galoisiennes finies de $F$. D'apr\`es le lemme 6.8, on peut fixer $u_{\star}\in \bar{G}_{\star}({\mathbb A}_{\bar{F}})$ tel que $ad_{r_{\star}u_{\star}}(S_{\star}({\mathbb A}_{\bar{F}}))=T^*({\mathbb A}_{\bar{F}})$ et que $ad_{r_{\star}u_{\star}}$ entrelace l'action galoisienne naturelle de $\Gamma_{F}$ sur $S_{\star}({\mathbb A}_{\bar{F}})$ avec l'action $\sigma\mapsto \omega_{S}(\sigma)\circ \sigma_{G^*}$ sur $T^*({\mathbb A}_{\bar{F}})$. Il en r\'esulte que $ad_{\bar{r}_{\star}^{-1}u_{\star}}$ se restreint en un isomorphisme  \'equivariant pour les actions de $\Gamma_{F}$ de $S_{\star}({\mathbb A}_{\bar{F}})$ sur $S({\mathbb A}_{\bar{F}})$. De m\^eme, pour tout $v$, cet isomorphisme se restreint en un isomorphisme \'equivariant pour les actions de $\Gamma_{F_{v}}$ de $S_{\star,v}(\bar{F}_{v})$ sur $S(F_{v})$. Des tores $S_{\star,v}$ se d\'eduisent comme ci-dessus des tores $\bar{S}_{\star,v}$, $\bar{S}_{\star,v,sc}$ et $S_{\star,v,sc}$ et des isomorphismes analogues relient ces tores \`a $\bar{S}$, $\bar{S}_{sc}$, $S_{sc}$. 
 
 Ainsi, quand on consid\`ere les points sur $\bar{F}_{v}$, resp. ${\mathbb A}_{\bar{F}}$, on peut identifier  $S$  \`a $S_{\star}$. La diff\'erence essentielle est que $S$ \`a une structure sur $\bar{F}$, donc le groupe $S(\bar{F})$ est bien d\'efini tandis que $S_{\star}(\bar{F})$ n'a pas de sens. Toutefois, le groupe $Z(I_{\star};{\mathbb A}_{\bar{F}})$  se plonge naturellement dans $S({\mathbb A}_{\bar{F}})$ comme dans $S_{\star}({\mathbb A}_{\bar{F}})$. Puisque $\bar{r}_{\star}$ et $u_{\star}$ appartiennent \`a $\bar{G}_{\star}$, l'isomorphisme $S({\mathbb A}_{\bar{F}})\simeq S_{\star}({\mathbb A}_{\bar{F}})$ se restreint en l'identit\'e de $Z(I_{\star};{\mathbb A}_{\bar{F}})$. Puisque $Z(I_{\star})$ est d\'efini sur $F$, on obtient le diagramme commutatif
 $$\begin{array}{ccccc}&&S(\bar{F})&\to&S({\mathbb A}_{\bar{F}})\\ &\nearrow&&\nearrow&\\ Z(I_{\star};\bar{F})&\to&Z(I_{\star};{\mathbb A}_{\bar{F}})&&\downarrow\\ &&&\searrow&\\ &&&&S_{\star}({\mathbb A}_{\bar{F}})\\ \end{array}$$
 Une propri\'et\'e analogue vaut pour  $\bar{S}_{sc}$ et $\bar{S}_{\star,sc}$, le groupe $Z(I_{\star})$ \'etant remplac\'e par $Z(\bar{G}_{\star,SC})$. 
 
 En vertu de ce que l'on a dit en 6.2, les groupes de cohomologie ab\'elienne de $I_{\star}$, \`a coefficients dans $F$, $F_{v}$, ${\mathbb A}_{F}$ ou ${\mathbb A}_{F}/F$, peuvent se calculer \`a l'aide du complexe $\bar{S}_{sc}\to S\stackrel{1-\theta}{\to}(1-\theta)(S)$. Sur $F_{v}$ ou ${\mathbb A}_{F}$, on peut aussi bien utiliser le complexe $\bar{S}_{\star,sc}\to S_{\star}\stackrel{1--\theta}{\to}(1-\theta)(S_{\star})$.
 
 On pose ${\cal Y}_{\star}={\cal Y}_{\eta_{\star}}$. Rappelons que c'est l'ensemble des $y\in G(\bar{F})$ tels que $y\sigma(y)^{-1}\in I_{\star}$. Pour $y\in {\cal Y}_{\star}$, posons $\eta[y]=ad_{y^{-1}}(\eta_{\star})$. L'application $y\mapsto ( \eta[y],r_{\star}y)$  est une bijection de ${\cal Y}_{\star}$ sur  $D_{F}$ (la preuve est la m\^eme qu'en 5.4(4)). Gr\^ace au lemme 5.5, l'inverse de cette bijection identifie le sous-ensemble $D_F[d_V]\subset D_{F}$ au sous-ensemble des $y\in {\cal Y}_{\star}$ tels que

- pour tout $v\in V$, $y\in I_{\star}(\bar{F}_{v})G(F_{v})$;

- pour tout $v\not\in V$, $y\in I_{\star}(\bar{F}_{v})\underline{K}_{\sharp,v}G(F_{v})$.

On note ${\cal Y}_{\star}[d_{V}]$ cet ensemble. L'inverse de la bijection ci-dessus identifie $\dot{D}_{F}[d_{V}]$ \`a un ensemble de repr\'esentants de l'ensemble de doubles classes $I_{\star}(\bar{F})\backslash {\cal Y}_{\star}[d_{V}]/G(F)$. On note cet ensemble $\dot{{\cal Y}}_{\star}[d_{V}]$.

   On a

(1) pour tout $y\in {\cal Y}_{\star}[d_V]$, il existe $y'\in I_{\star}(\bar{F})y$ tel que $y'\sigma(y')^{-1}$ appartienne au centre de $I_{\star}$ pour tout $\sigma\in \Gamma_{F}$.

 Preuve. Pour $y\in  {\cal Y}_{\star}[d_V]$, notons $\chi_{y}$ le cocycle $\sigma\mapsto y\sigma(y)^{-1}$ de $\Gamma_{F} $ dans $I_{\star}$ et $\chi_{y,ad}$ son image naturelle \`a valeurs dans $\bar{G}_{\star,AD}$. L'automorphisme $ad_{y}$ se restreint en un torseur int\'erieur de $G_{\eta[y]}$ sur $\bar{G}_{\star}$. La classe d'isomorphisme de $G_{\eta[y]}$ est d\'etermin\'e par l'\'el\'ement $\chi_{y,ad}\in H^1(\Gamma_{F};\bar{G}_{\star,AD})$. Or, par d\'efinition de $ {\cal Y}_{\star}[d_V]$, ce cocycle est localement trivial en toute place $v\in Val(F)$. Parce que le groupe $\bar{G}_{\star,AD}$ est adjoint, l'application
$$H^1(F;\bar{G}_{\star,AD})\to \oplus_{v\in Val(F)}H^1(F_{v};\bar{G}_{\star,AD})$$
est injective ([S] corollaire 5.4). Donc $\chi_{y,ad}$ est trivial. Quitte \`a multiplier $y$ \`a gauche par un \'el\'ement de $\bar{G}_{\star}(\bar{F})\subset I_{\star}(\bar{F})$, on peut donc supposer que $ \chi_{y,ad}(\sigma)=1$ pour tout $\sigma\in \Gamma_{F}$. Cela entra\^{\i}ne que $\chi_{y}(\sigma)$ appartient au centre de $I_{\star}$. $\square$

Remarquons que la propri\'et\'e de $y'$ entra\^{\i}ne que $ad_{y'}$ se restreint en un isomorphisme d\'efini sur $F$ de $I_{\eta[y]}$ sur $I_{\star}$.  

 Par la bijection de   ${\cal Y}_{\star}[d_{V}]$ sur $D_{F}[d_{V}]$, la multiplication \`a gauche par $I_{\star}(\bar{F})$ correspond \`a la multiplication \`a gauche par $I_{\eta}(\bar{F})$. On voit qu'une telle multiplication ne modifie pas les termes intervenant dans la d\'efinition de la fonction $\bar{\varphi}[V',d_{V}]$. En cons\'equence, on peut supposer que tous les \'el\'ements de l'ensemble de repr\'esentants $\dot{\cal Y}_{\star}[d_{V}]$ satisfont la conclusion de (1). On ne perd rien non plus \`a supposer que $d_{\star}$ appartient \`a l'ensemble $\dot{D}_{F}[d_{V}]$. Il revient au m\^eme de supposer que $1$ appartient \`a $\dot{\cal Y}_{\star}[d_{V}]$. Pour $d\in \dot{D}_{F}[d_{V}]$, on a d\'efini en 7.1 des termes $h[d]$, ${\cal U}[V',d]$ et, pour $u\in {\cal U}[V',d]$, des termes $K_{sc,v}[d,u]$, $\Delta_{V}[d,u]$,  $\bar{f}[d,u]$. Si $d$ correspond \`a $y\in \dot{\cal Y}_{\star}[d_{V}]$, on note aussi ces termes $h[y]$, ${\cal U}[V',y]$, $K_{sc,v}[y,u]$, $\Delta_{V}[y,u]$ et $\bar{f}[y,u]$. On fera une exception pour le terme $d_{\star}$ correspondant \`a $1\in \dot{\cal Y}_{\star}[d_{V}]$. Dans ce cas, on notera  $h_{\star}=h[1]$  et $ \bar{f}_{\star}[u]=\bar{f}[1,u]$.

 \bigskip
  
  \subsection{Une premi\`ere propri\'et\'e de nullit\'e}
  Rappelons (cf. [I] 2.7) que la donn\'ee ${\bf H}$ de $\bar{G}_{SC}$ d\'efinit un caract\`ere automorphe de $\bar{G}_{AD}({\mathbb A}_{F})$, que nous notons ici $\omega_{{\bf H}}$. Il se transporte en un caract\`ere automorphe  de $\bar{G}_{\star,AD}({\mathbb A}_{F})$. On a un homomorphisme naturel $I_{\star}({\mathbb A}_{F})\to \bar{G}_{\star,AD}({\mathbb A}_{F})$ gr\^ace auquel le caract\`ere $\omega_{{\bf H}}$ devient un caract\`ere de $I_{\star}({\mathbb A}_{F})$. De m\^eme, pour $y\in \dot{{\cal Y}}_{\star}(d_{V}]$, on a un caract\`ere de $I_{\eta[y]}({\mathbb A}_{F})$ que l'on note encore $\omega_{{\bf H}}$ et qui n'est autre que le compos\'e du caract\`ere de $I_{\star}({\mathbb A}_{F})$ par l'isomorphisme $ad_{y}$ entre ces deux groupes. On a aussi le caract\`ere $\omega$ qui se restreint \`a chacun des groupes $I_{\eta[y]}({\mathbb A}_{F})$ et $I_{\star}({\mathbb A}_{F})$. Remarquons que 
  
  (1)  l'isomorphisme $ad_{y}$ conserve le caract\`ere $\omega$.
  
  En effet, soient $v\in Val(F)$ et $x\in I_{\eta[y]}(F_{v})$. Le commutateur $ad_{y}(x)x^{-1}$ est l'image naturelle de l'\'el\'ement $ad_{y_{sc}}(x_{sc})x_{sc}^{-1}$ de $G_{SC}$, o\`u $y_{sc}$ et $x_{sc}$ sont des \'el\'ements de $G_{SC}$ ayant m\^eme image que $y$ et $x$ dans $G_{AD}$. Pour $\sigma\in \Gamma_{F_{v}}$, la condition que $y\sigma(y)^{-1}$ commute \`a $I_{\eta}$, cf. 7.3(1), implique que $\sigma(y_{sc})^{-1}y_{sc}$ commute \`a $x_{sc}$.  Puisque de plus $\sigma(x_{sc})\in Z(G_{SC})x_{sc}$, on v\'erifie que $ad_{y_{sc}}(x_{sc})x_{sc}^{-1}$ est fixe par $\sigma$, donc $ad_{y_{sc}}(x_{sc})x_{sc}^{-1}\in G_{SC}(F_{v})$. Le caract\`ere $\omega$ est trivial sur ce groupe. Donc $\omega(ad_{y}(x)x^{-1})=1$, ce qui prouve l'assertion (1). $\square$
  
 Pour $v\not\in V$, posons $K_{\star,v}=ad_{h_{\star,v}}(K_{v})\cap \bar{G}_{\star}(F_{v})$. C'est un sous-groupe compact hypersp\'ecial de  de $\bar{G}_{\star}(F_{v})$.   Imposons \`a l'ensemble de places $V'$ de 7.1 la condition
  
  (2) l'ensemble $\bar{G}_{\star}(F)\bar{G}_{\star}(F_{V'})\prod_{v\not\in V'}K_{\star,v}$ est dense dans $\bar{G}_{\star}({\mathbb A}_{F})$.
  
  \ass{Lemme}{Si $\omega_{{\bf H}}$ et $\omega$ ne co\"{\i}ncident pas sur $I_{\star}({\mathbb A}_{F})$, alors $\bar{\varphi}[V',d_{V}]=0$. }
  
  Preuve. Soient $y\in \dot{{\cal Y}}_{\star}[d_{V}]$ et $u\in {\cal U}[V',y]$. Rappelons la construction de la fonction $\bar{f}[y,u]$. On descend la fonction $f$ en une fonction, disons $f[y]$, sur $G_{\eta[y]}(F_{V})$. On d\'efinit la fonction $f[y]_{sc}$ sur $G_{\eta[y],SC}(F_{V})$ par $f[y]_{sc}=\iota_{G_{\eta[y],SC},G_{\eta[y]}}(f[y])$. Alors $\bar{f}[y,u]$ est le transfert de $f[y]_{sc}$  \`a $\bar{H}(F_{V})$ pour le facteur de transfert $\Delta_{V}[y,u]$. Notons que cette construction ne d\'epend que de l'image de $f$ dans $I(\tilde{G}(F_{V}),\omega)$. Le groupe $Z_{G}(\eta[y];F_{V})$ agit naturellement sur $I(G_{\eta[y],SC}(F_{V}))$.     Parce que l'on part d'un \'el\'ement de $I(\tilde{G}(F_{V}),\omega)$, la fonction $f[y]_{sc}$ appartient au sous-espace isotypique de $I(G_{\eta[y],SC}(F_{V}))$ sur lequel $Z_{G}(\eta[y];F_{V})$ agit par le caract\`ere $\omega$. Le groupe $G_{\eta[y],AD}(F_{V})$ agit aussi sur   $G_{\eta[y],SC}(F_{V})$ et sur $I(G_{\eta[y],SC}(F_{V}))$. Par cette action, le facteur de transfert $\Delta_{V}[y,u]$ se transforme par le caract\`ere $\omega_{{\bf H}}^{-1}$. Donc le transfert se factorise par la projection sur le sous-espace isotypique de $I(G_{\eta[y],SC}(F_{V}))$ sur lequel $G_{\eta[y],AD}(F_{V})$ agit par le caract\`ere $\omega_{{\bf H}}$. Le groupe $I_{\eta[y]}(F_{V})$ s'envoie \`a la fois dans $Z_{G}(\eta[y];F_{V})$ et dans $G_{\eta[y],AD}(F_{V})$ et les deux actions co\"{\i}ncident sur ce groupe. Cela entra\^{\i}ne que le transfert $\bar{f}[y,u]$ de $f[y]_{sc}$ est nul si les deux caract\`eres ne co\"{\i}ncident pas sur $I_{\eta[y]}(F_{V})$. D'apr\`es (1), cette condition de co\"{\i}ncidence \'equivaut \`a la m\^eme condition portant sur les caract\`eres de  $I_{\star}(F_{V})$. On obtient que $\bar{\varphi}[V',d_{V}]=0$ si les deux caract\`eres $\omega_{{\bf H}}$ et $\omega$ de  $I_{\star}(F_{V})$ sont distincts. 
  
  Soient de nouveau $y\in \dot{{\cal Y}}_{\star}[d_{V}]$ et $u\in {\cal U}[V',y]$. Rappelons que l'on a suppos\'e $1\in {\cal U}[V',y]$. Par construction, le rapport $\Delta_{V}[y,u]/\Delta_{V}[y,1]$ est \'egal au produit sur les places $v\in V'-V$ des rapports de facteurs locaux normalis\'es $\Delta_{v}[y,1]/\Delta_{v}[y,u_{v}]$. Le facteur $\Delta_{v}[y,1]$ est relatif au sous-groupe compact $K_{sc,v}[y,1] $ et le facteur $\Delta_{v}[y,u]$ est relatif au sous-groupe $K_{sc,v}[y,u]=ad_{u_{v}}(K_{sc,v}[y,1])$. Pour $\bar{y}\in \bar{H}(F_{v})$ et $x\in G_{\eta[y],SC}(F_{v})$, on a par simple transport de structure l'\'egalit\'e $\Delta_{v}[y,u](\bar{y},ad_{u_{v}}(x))=\Delta_{v}[y,1](\bar{y},x)$. Mais le premier terme est \'egal \`a $\omega_{{\bf H}}(u_{v})^{-1}\Delta_{v}[y,u](\bar{y},x)$. On en d\'eduit 
 $\Delta_{v}[y,1]/\Delta_{v}[y,u_{v}]=\omega_{{\bf H}}(u_{v})^{-1}$ puis $\Delta_{V}[y,u]/\Delta_{V}[y,1]=\omega_{{\bf H}}(u)^{-1}$. Donc $\bar{f}[y,u]=\omega_{{\bf H}}(u)^{-1}\bar{f}[y,u]$. Dans la d\'efinition 7.1(3), la somme en $u$ se r\'ecrit donc
 $$\omega(h[y])\bar{f}[y,1]\sum_{u\in {\cal U}[V',y]}\omega(u)\omega_{{\bf H}}(u)^{-1}.$$
 Remarquons que les deux caract\`eres $\omega$ et $\omega_{{\bf H}}$ sont  \'evidemment triviaux sur  $G_{\eta[y],SC}(F_{V'}^V)$ et le sont aussi sur $Z(G_{\eta[y]};F_{V'}^V)$: pour $\omega_{{\bf H}}$, cela resulte de sa construction comme image r\'eciproque d'un caract\`ere de $G_{\eta[y],AD}(F_{V'}^V)$; pour $\omega$, cela r\'esulte de l'hypoth\`ese 5.1(4).   La somme en $u$ ci-dessus est donc \'egale \`a la somme des valeurs d'un caract\`ere du groupe quotient
 $$G_{\eta[y]}(F_{V'}^V)/  Z(G_{\eta[y]};F_{V'}^V)G_{\eta[y],SC}(F_{V'}^V).$$
 Elle est nulle si ce caract\`ere n'est pas trivial, autrement dit si $\omega$ et $\omega_{{\bf H}}$ ne co\"{\i}ncident pas sur $G_{\eta[y]}(F_{V'}^V)$. De nouveau, cela \'equivaut \`a ce que les caract\`eres correspondant de $\bar{G}_{\star}(F_{V'}^V)$ sont distincts. On obtient que $\bar{\varphi}[V',d_{V}]=0$ si les deux caract\`eres $\omega_{{\bf H}}$ et $\omega$ de  $\bar{G}_{\star}(F_{V'}^V)$ sont distincts. 
 
 Soit $v\not\in V$.  Comme on l'a expliqu\'e dans la preuve du lemme 1.6, les conditions impos\'ees \`a $V$ nous autorisent \`a appliquer les r\'esultats de [W1]. Le lemme 5.6(ii) de cette r\'ef\'erence implique
 
 (3) $I_{\star}(F_{v})=\bar{G}_{\star}(F_{v})(ad_{h_{\star,v}}(K_{v})\cap I_{\star}(F_{v}))$.
 
Du groupe $K_{\star,v}$ se d\'eduit un sous-groupe compact hypersp\'ecial $K_{\star,ad,v}$ de $\bar{G}_{\star,AD}(F_{v})$. Montrons que

(4) l'image de $ad_{h_{\star}}(K_{v})\cap I(F_{v})$ dans $\bar{G}_{\star,AD}(F_{v})$ est contenue dans $K_{\star,ad,v}$.

On a dit dans la preuve du lemme 5.5 que l'application
$$Z(G)_{p'}^{\theta}\to   ad_{h_{\star,v}}(K^{nr}_{v})\cap \bar{G}_{\star}\backslash ad_{h_{\star,v}}(K^{nr}_{v})\cap  I_{\star}$$
\'etait surjective.   Il en r\'esulte que l'image de  $ ad_{h_{\star,v}}(K_{v}^{nr})\cap I_{\star}(F_{v})$ dans $\bar{G}_{\star,AD}(F_{v})$ est contenue dans celle de $ ad_{h_{\star,v}}(K_{v}^{nr})\cap \bar{G}_{\star}$, c'est-\`a-dire celle de $K_{\star,v}^{nr}$. Celle-ci est contenue dans $K_{\star,ad,v}^{nr}$. L'intersection de ce groupe avec $\bar{G}_{\star,AD}(F_{v})$ est \'egale \`a $K_{\star,ad,v}$, d'o\`u (4).

    Puisque $V\supset V_{ram}$, $\omega$ est trivial sur $ ad_{h_{\star,v}}(K_{v})\cap I_{\star}(F_{v})$. La donn\'ee ${\bf H}$ est non ramifi\'ee hors de $V$. Donc le caract\`ere $\omega_{{\bf H}}$ est trivial sur $K_{\star,ad,v}$.  D'apr\`es (4), le caract\`ere de $I_{\star}(F_{v})$ qui s'en d\'eduit est trivial sur $ad_{h_{\star,v}}(K_{v})\cap I_{\star}(F_{v})$. 
    
    Achevons la preuve du lemme. Supposons $\bar{\varphi}[V',d_{V}]\not=0$. On a vu que cela impliquait que $\omega$ et $\omega_{{\bf H}}$ co\"{\i}ncidaient sur $I_{\star}(F_{V})$ et sur $\bar{G}_{\star}(F_{V'}^V)$. Ces caract\`eres co\"{\i}ncident aussi sur $K_{\star,v}$ pour $v\not\in V$. Puisque ces caract\`eres sont automorphes, la condition (2) implique qu'ils co\"{\i}ncident sur tout $\bar{G}_{\star}({\mathbb A}_{F})$. Pour $v\not\in V$, la propri\'et\'e (3) et le fait qu'ils sont triviaux sur $ad_{h_{\star,v}}(K_{v})\cap I_{\star}(F_{v})$ impliquent alors qu'ils co\"{\i}ncident sur $I_{\star}(F_{v})$. Donc ils co\"{\i}ncident sur $I_{\star}({\mathbb A}_{F})$, ce qui prouve le lemme. $\square$

 \bigskip
 
 \subsection{Description de l'ensemble $\dot{{\cal Y}}_{\star}[d_{V}]$}
 Les ensembles $H^1_{ab}(F,G)$ et $ H^1({\mathbb A}_{F},G)$ s'envoient naturellement dans $H^1_{ab}({\mathbb A}_{F};G)$. On note
  $H^1_{ab}(F,G)\times_{H^1_{ab}({\mathbb A}_{F};G)}H^1({\mathbb A}_{F},G)$ leur produit fibr\'e  au-dessus de   $H^1_{ab}({\mathbb A}_{F};G)$. On a un diagramme commutatif
  $$(1) \qquad \begin{array}{ccc}H^1(F;I_{\star})&\to&H^1_{ab}(F,I_{\star})\times_{H^1_{ab}({\mathbb A}_{F};I_{\star})}H^1({\mathbb A}_{F},I_{\star}) \\ \downarrow&&\downarrow\\ H^1(F;G) &\to &H^1_{ab}(F,G)\times_{H^1_{ab}({\mathbb A}_{F};G)}H^1({\mathbb A}_{F},G)\\ \end{array}$$
  Soit $v$ une place finie de $F$. Alors l'application naturelle $H^1(F_{v};I_{\star})\to H^1_{ab}(F_{v};I_{\star})$ est bijective ([Lab2] proposition 1.6.7). Supposons $v\not\in V$. Alors le groupe $I_{\star}$ est non ramifi\'e et on d\'efinit le sous-groupe $H_{ab}^1(\mathfrak{o}_{v};I_{\star})\subset H^1_{ab}(F_{v};I_{\star})$.  Rappelons que $H^1_{ab}({\mathbb A}_{F};I_{\star})$ s'identifie au produit restreint des $H^1_{ab}(F_{v};I_{\star})$ sur toutes les places $v\in Val(F)$, la restriction \'etant relative aux sous-groupes $H^1_{ab}(\mathfrak{o}_{v};I_{\star})$ d\'efinis presque partout. On note $H^1_{ab}(\mathfrak{o}^V;I_{\star})$ le produit des $H_{ab}^1(\mathfrak{o}_{v};I_{\star})$ sur toutes les places $v\not\in V$.
  
  Notons ${\mathbb U}$ le sous-ensemble des $u\in H^1(F;I_{\star})$ qui v\'erifient les trois conditions
  
  (2) l'image de $u$ dans $H^1_{ab}(F;G)$ par l'application issue du diagramme (1) est nulle;
  
  (3) pour $v\in V$, l'image de $u$ dans $H^1(F_{v};I_{\star})$ est nulle;
  
  (4) l'image de $u$ dans $H^1_{ab}({\mathbb A}_{F}^V;I_{\star})$ appartient \`a $H^1_{ab}(\mathfrak{o}^V;I_{\star})$. 
  
  On a une application naturelle $H^1(F;Z(I_{\star}))\to H^1(F;I_{\star})$. Montrons que
  
  (5) son noyau $Ker$  est un sous-groupe et l'application se quotiente en une injection $H^1(F;Z(I_{\star}))/Ker\to H^1(F;I_{\star})$. 
  
  Preuve. Soient $z_{1}$ et $z_{2}$ deux \'el\'ements de $Ker$, que l'on rel\`eve en des cocycles. On  peut choisir $i_{1}$ et $i_{2}$ dans $I_{\star}$ tels que, pour $j=1,2$ et pour tout $\sigma\in \Gamma_{F}$, on ait $z_{j}(\sigma)=i_{j}\sigma(i_{j})^{-1}$. Puisque $z_{2}(\sigma)$ est central dans $I_{\star}$, on a 
  $$z_{1}(\sigma)z_{2}(\sigma)=i_{1}\sigma(i_{1})^{-1}z_{2}(\sigma)=i_{1}z_{2}(\sigma)\sigma(i_{1})^{-1}=i_{1}i_{2}\sigma(i_{1}i_{2})^{-1}.$$
  Donc $z_{1}z_{2}\in Ker$. On a aussi
  $$z_{1}(\sigma)^{-1}=i_{1}^{-1}z_{1}(\sigma)^{-1}i_{1}=i_{1}^{-1}\sigma(i_{1})i_{1}^{-1}i_{1}=i_{1}^{-1}\sigma(i_{1}).$$
  Donc $z_{1}^{-1}\in Ker$ et $Ker$ est bien un sous-groupe. Soient $z_{1}$ et $z_{2}$ deux \'el\'ements de $H^1(F;Z(I_{\star}))$, que l'on rel\`eve en des cocycles. Leurs images dans $H^1(F;I_{\star})$ co\"{\i}ncident si et seulement s'il existe $i\in I_{\star}$ de sorte que, pour tout $\sigma\in \Gamma_{F}$, on ait $z_{1}(\sigma)=iz_{2}(\sigma)\sigma(i)^{-1}$. Parce que $z_{2}$ est \`a valeurs centrales, cela \'equivaut \`a $z_{1}(\sigma)z_{2}(\sigma)^{-1}=i\sigma(i)^{-1}$. Mais l'existence de $i$ v\'erifiant cette condition \'equivaut \`a $z_{1}z_{2}^{-1}\in Ker$. D'o\`u la seconde assertion de (5). $\square$
  
  Notons $H^1_{Z}(F;I_{\star})$ l'image de $H^1(F;Z(I_{\star}))$ dans $H^1(F;I_{\star})$. Gr\^ace \`a (5), cet ensemble est naturellement un groupe. On d\'efinit une application ${\cal Y}_{\star}\to H^1(F;I_{\star})$ qui, \`a $y\in {\cal Y}_{\star}[d_{V}]$, associe la classe du cocycle $\sigma\mapsto y\sigma(y)^{-1}$.

\ass{Lemme}{Cette application se restreint en une bijection de $\dot{{\cal Y}}_{\star}[d_{V}]$ sur ${\mathbb U}$. L'ensemble ${\mathbb U}$ est un sous-groupe de $H_{Z}^1(F;I_{\star})$. }

Preuve.    Il est imm\'ediat que notre application  se quotiente en une bijection de 
$$(6) \qquad I_{\star}(\bar{F})\backslash {\cal Y}_{\star}/G(F)$$
 dans le noyau de l'application
$$(7)\qquad H^1(F;I_{\star})\to H^1(F;G).$$
 Rappelons que $\dot{{\cal Y}}_{\star}[d_{V}]$ est un ensemble de repr\'esentants de l'image de ${\cal Y}_{\star}[d_{V}]\subset {\cal Y}_{\star}$ dans l'ensemble de doubles classes (6).  L'application  se restreint donc en  une injection de  $\dot{{\cal Y}}_{\star}[d_{V}]$ dans le noyau de (7). D'apr\`es la d\'efinition de 7.3, son image est form\'ee des \'el\'ements $u$ de ce noyau qui v\'erifient la condition
(3) et 

(8) pour tout $v\not\in V$, il existe $k\in \underline{K}_{\sharp,v}$ tel que l'image de $u$ dans $H^1(F_{v};I_{\star})$ est cohomologue au cocycle  $\sigma\mapsto k\sigma(k)^{-1}$.

Rappelons que, par d\'efinition de $\underline{K}_{\sharp,v}$, ce dernier cocycle prend ses valeurs dans $Z(G)^{\theta}\subset Z(I_{\star})$. Montrons que 

(9) les conditions (4) et (8) sont \'equivalentes. 

Soit $v\not\in V$. Fixons un sous-tore maximal $T_{\natural}$ de $\bar{G}_{\star}$ d\'efini sur $F_{v}$ et non ramifi\'e. Notons $T$ son commutant dans $G$ et $T_{\natural,sc}$ l'image r\'eciproque de $T_{\natural}$ dans $\bar{G}_{\star,SC}$. Par d\'efinition
$$H^1_{ab}(\mathfrak{o}_{v};I)=H^{2,1,0}(\mathfrak{o}_{v};T_{\natural,sc}\to T\stackrel{1-\theta}{\to}(1-\theta)(T)).$$
Le diagramme
$$\begin{array}{ccccc}T_{\natural,sc}&\to& T&\stackrel{1-\theta}{\to}&(1-\theta)(T)\\ \downarrow&&\downarrow&&\downarrow\\ T_{\natural,sc}&\to &T/Z(G)^{\theta}&\stackrel{1-\theta}{\to}&(1-\theta)(T)\\ \downarrow&&\downarrow&& \\ T&\to &T/Z(G)^{\theta}&&\\ \end{array}$$
est un triangle exact dans la cat\'egorie des complexes de tores. On en d\'eduit une suite exacte
$$(10) \qquad H^{1,0}(\mathfrak{o}_{v}; T\to T/Z(G)^{\theta})\to H^{2,1,0}(\mathfrak{o}_{v};T_{\natural,sc}\to T\stackrel{1-\theta}{\to}(1-\theta)(T))$$
$$\to H^{2,1,0}(\mathfrak{o}_{v};T_{\natural,sc}\to T/Z(G)^{\theta}\stackrel{1-\theta}{\to}(1-\theta)(T)).$$
D'autre part, le diagramme
$$\begin{array}{ccccc}T_{\natural,sc}&\to& T^{\theta}/Z(G)^{\theta}&&\\ \downarrow&&\downarrow&&\\ T_{\natural,sc}&\to &T/Z(G)^{\theta}&\stackrel{1-\theta}{\to}&(1-\theta)(T)\\ \end{array}$$
induit un isomorphisme de cohomologie
$$H^{2,1}(\mathfrak{o}_{v}; T_{\natural,sc}\to T^{\theta}/Z(G)^{\theta})\simeq H^{2,1,0}(\mathfrak{o}_{v};T_{\natural,sc}\to T/Z(G)^{\theta}\stackrel{1-\theta}{\to}(1-\theta)(T)).$$
Or le premier groupe est nul ([KS] lemme C.1.A). On en d\'eduit que le premier homomorphisme de la suite (10) est surjectif. Cet homomorphisme se r\'ecrit
$$H^0_{ab}(\mathfrak{o}_{v};G_{\sharp})\to H_{ab}^1(\mathfrak{o}_{v};I_{\star})$$
o\`u on rappelle que $G_{\sharp}=G/Z(G)^{\theta}$. 
L'ensemble de d\'epart  est un sous-groupe de $H^0_{ab}(F_{v};G_{\sharp})$, qui est un quotient de $G_{\sharp}(F_{v})$. L'assertion 1.5(2) \'equivaut \`a dire que $H^0_{ab}(\mathfrak{o}_{v};G_{\sharp})$ est l'image de $K_{\sharp,v}$ dans $H^0_{ab}(F_{v};G_{\sharp})$. Puisque $\underline{K}_{\sharp,v}$ s'envoie surjectivement sur $K_{\sharp,v}$, on obtient une application surjective
$$(11) \qquad \underline{K}_{\sharp,v}\to  H_{ab}^1(\mathfrak{o}_{v};I_{\star}).$$
On a alors deux fa\c{c}ons d'envoyer $\underline{K}_{\sharp,v}$ dans $H^1(F_{v},I_{\star})$. D'abord celle de la relation (8): \`a $k\in \underline{K}_{\sharp,v}$ on associe le cocycle $\sigma\mapsto k\sigma(k)^{-1}$. On peut aussi envoyer $k$ en un \'el\'ement  de $H^1_{ab}(\mathfrak{o}_{v},I_{\star})$ par l'application pr\'ec\'edente, puis on  rel\`eve celui-ci en un \'el\'ement de $H^1(F_{v},I_{\star})$ en utilisant la bijectivit\'e de l'application $H^1(F_{v},I_{\star})\to H^1_{ab}(F_{v},I_{\star})$ ([Lab2] prop. 1.6.7). En inspectant les d\'efinitions, on s'aper\c{c}oit que les deux applications obtenues co\"{\i}ncident. Puisque l'application (11) est surjective,  la condition (8) \'equivaut donc \`a ce que l'image de $u$ dans $H^1_{ab}(F_{v};I_{\star})$ appartienne \`a  $H^1_{ab}(\mathfrak{o}_{v};I_{\star})$, ce qui est la condition (4). Cela prouve (9). 

Montrons que, pour $u\in H^1(F;I_{\star})$ v\'erifiant les conditions (3) et (4), on a

(12) l'image de $u$ dans $H^1(F;G)$ est nulle si et seulement si son image  dans $H^1_{ab}(F;G)$ est nulle.

La premi\`ere condition implique la seconde. Inversement,  notre \'el\'ement $u$ v\'erifie (8) d'apr\`es (9). Cela entra\^{\i}ne que, pour $v\not\in V$, son image dans $H^1(F_{v};G)$ est nulle. En ajoutant (3), l'image de $u$ dans $H^1({\mathbb A}_{F};G)$ est nulle. Si de plus l'image de $u$ dans $H^1_{ab}(F;G)$ est nulle, son image dans le terme sud-est du diagramme (1) est nulle. Or l'application du bas de ce diagramme est bijective ([Lab2] th\'eor\`eme 1.6.10). Donc l'image de $u$ dans $H^1(F;G)$ est nulle. Cela prouve (12).

On a vu que l'image de $\dot{{\cal Y}}_{\star}[d_{V}]$ dans $H^1(F;I_{\star})$ est l'ensemble des $u$ v\'erifiant les conditions (3) et (8) et dont l'image dans $H^1(F;G)$ est nulle. Gr\^ace \`a (9) et (12), c'est exactement l'ensemble ${\mathbb U}$. 

D'apr\`es 7.3(1), l'image de $\dot{{\cal Y}}_{\star}[d_{V}]$ dans $H^1(F;I_{\star})$ est contenue dans $H^1_{Z}(F;I_{\star})$. Il reste \`a prouver que cette image ${\mathbb U}$ est un sous-groupe.  On v\'erifie facilement que les applications compos\'ees
  $$H^1(F;Z(I_{\star}))\to H^1(F;I_{\star})\to H^1_{ab}(F;G)$$
  et
  $$H^1(F;Z(I_{\star}))\to H^1(F;I_{\star})\to H^1_{ab}({\mathbb A}_{F}^V;I_{\star})$$
  sont des homomorphismes de groupes. On en d\'eduit que l'ensemble des \'el\'ements de $H^1_{Z}(F;I_{\star})$ qui v\'erifient les conditions (2) et (4) est un sous-groupe de $H^1_{Z}(F;I_{\star})$. D'autre part,  pour $v\in V$, on a le diagramme commutatif
  $$\begin{array}{ccc}H^1(F;Z(I_{\star}))&\to& H^1(F;I_{\star})\\ \downarrow&&\downarrow\\ H^1(F_{v};Z(I_{\star}))&\to&H^1(F_{v};I_{\star})\\ \end{array}$$
  Un \'el\'ement de $H^1_{Z}(F;I_{\star})$ v\'erifie (3) si et seulement si c'est l'image d'un \'el\'ement de $H^1(F;Z(I_{\star}))$ dont l'image dans $H^1(F_{v};Z(I_{\star}))$ pout $v\in V$ appartient au noyau de l'application du bas. La m\^eme preuve qu'en (5) montre que ce noyau est un groupe. Donc l'ensemble des \'el\'ements de $H^1_{Z}(F;I_{\star})$ qui v\'erifient (3) est l'image d'un sous-groupe de $H^1(F;Z(I_{\star}))$. Cela conclut. $\square$ 
  
 \bigskip
 \subsection{D\'efinition d'un homomorphisme ${\bf q}_{\infty}$}
 Consid\'erons les groupes 
 $$Q_{1}=H^1({\mathbb A}_{F}/F;\bar{S}_{sc}) ,\,\,Q_{2}=H^0_{ab}({\mathbb A}_{F};G)/Im(H^0_{ab}(F;G)),$$
 $$Q_{3}=H^0_{ab}(\mathfrak{o}^V;G_{\sharp})=\prod_{v\not\in V}H^0_{ab}(\mathfrak{o}_{v};G_{\sharp})$$
 d\'efinis en 6.6.
  Posons
 $$Q_{\times}=Q_{1}\times Q_{2}\times Q_{3}.$$
 Posons
 $$Q_{1,2}=I_{\star}({\mathbb A}_{F}),\,\, Q_{1,3}=H^0_{ab}(\mathfrak{o}^V;I_{\star}/Z(G)^{\theta}),\,\, Q_{2,3}=H_{ab}^0(\mathfrak{o}^V;G).$$
Le groupe $Q_{1,2}$ s'envoie naturellement dans $G({\mathbb A}_{F})$ lequel s'envoie dans $Q_{2}$.  Il s'envoie aussi dans $\bar{G}_{\star,AD}({\mathbb A}_{F})$ puis dans $H^0_{ab}({\mathbb A}_{F};\bar{G}_{\star,AD})$. Ce dernier groupe  s'identifie \`a $H^{1,0}({\mathbb A}_{F};\bar{S}_{sc}\to S^{\theta}/Z(I_{\star}))$, lequel s'envoie dans $H^1({\mathbb A}_{F};\bar{S}_{sc})$ puis dans $Q_{1}$. Puisque $I_{\star}/Z(G)^{\theta}$ est un sous-groupe de $G_{\sharp}$, $Q_{1,3}$ s'envoie naturellement dans $Q_{3}$.  Le groupe $Q_{1,3}$ est un sous-groupe de $H^0_{ab}({\mathbb A}_{F};I_{\star}/Z(G)^{\theta})$ que l'on peut identifier \`a $H^{1,0}({\mathbb A}_{F};\bar{S}_{sc}\to S^{\theta}/Z(G)^{\theta})$. Ce dernier s'envoie naturellement dans $H^1({\mathbb A}_{F};\bar{S}_{sc})$ puis dans $Q_{1}$. Par composition, on obtient un homomorphisme $Q_{1,3}\to Q_{1}$.   Le groupe $Q_{2,3}$ s'envoie naturellement dans $Q_{2}$ et $Q_{3}$. Toutes ces applications sont des homomorphismes. Ainsi, pour $1\leq j<j'\leq 3$, on a des homomorphismes
 $$\begin{array}{ccc}&&Q_{j}\\ &\nearrow &\\Q_{j,j'}&&\\ &\searrow &\\ &&Q_{j'}\\ \end{array}$$
 On en d\'eduit un homomorphisme
 ${\bf q}_{j,j'}:Q_{j,j'}\to Q_{\times}$ dont la composante dans $Q_{j}$ est l'homomorphisme pr\'ec\'edent, la composante dans $Q_{j'}$ est l'oppos\'e de l'homomorphisme pr\'ec\'edent et la derni\`ere composante est triviale. On note $Q_{\infty}$ le quotient de $Q_{\times}$ par le groupe engendr\'e par les images des homomorphismes ${\bf q}_{j,j'}$. 
  
  Soit $u:\Gamma_{F}\to Z(I_{\star};\bar{F})$ un cocycle dont l'image dans $H^1(F;I_{\star})$ appartient \`a ${\mathbb U}$. Pour $v\in V$, l'image de $u$ dans $H^1(F_{v};I_{\star})$ est nulle. On peut fixer $i_{v}\in I_{\star}$ de sorte que $u(\sigma)=i_{v}\sigma(i_{v})^{-1}$ pour tout $\sigma\in \Gamma_{F_{v}}$. Pour $v\not\in V$, $u$ v\'erifie la condition 7.5(8). On peut fixer $k_{v}\in \underline{K}_{\sharp,v}$ et $i_{v}\in I_{\star}$ de sorte que $u(\sigma)=i_{v}k_{v}\sigma(i_{v}k_{v})^{-1}$ pour tout $\sigma\in \Gamma_{F_{v}}$. Montrons que
  
  (1) on peut supposer $k_{v}\in \underline{K}_{\sharp,v}\cap K_{v}^{nr}$ et $i_{v}\in I_{\star}(F_{v}^{nr})\cap K_{v}^{nr}$ pour presque tout $v$.
  
  Preuve.    D'apr\`es le lemme 7.5, on peut fixer $y\in \dot{{\cal Y}}_{\star}[d_{V}]$ tel que $u$ soit cohomologue au cocycle $\sigma\mapsto y\sigma(y)^{-1}$. Quitte \`a multiplier $y$ \`a gauche par un \'el\'ement de $I_{\star}(\bar{F})$, on obtient un \'el\'ement $y\in G(\bar{F})$ tel que $u(\sigma)=y\sigma(y)^{-1}$ pour tout $\sigma\in \Gamma_{F}$.  On fixe un ensemble fini $V''$ de places de $F$, contenant $V$, tel que, pour $v\not\in V''$,

  - $y\in K_{v}^{nr}$;
  
  - $\eta_{\star}\in \tilde{K}_{v}$.
  
  Soit $v\not\in V''$.   D'apr\`es 5.5(1), on a $\underline{K}_{\sharp,v}=Z(G)^{\theta}(\underline{K}_{\sharp,v}\cap K_{v}^{nr})$. On peut donc \'ecrire $k_{v}=z_{v}k'_{v}$ avec $z_{v}\in Z(G)^{\theta}$ et $k'_{v}\in \underline{K}_{\sharp,v}\cap K_{v}^{nr}$. Puisque $Z(G)^{\theta}\subset I_{\star}$,  on peut remplacer le couple $(i_{v},k_{v})$ par $(i_{v}z_{v},k'_{v})$.  Apr\`es ce remplacement, on a $k_{v}\in \underline{K}_{\sharp,v}\cap K_{v}^{nr}$. L'\'egalit\'e $i_{v}k_{v}\sigma(i_{v}k_{v})^{-1}=y\sigma(y)^{-1}$ pour tout $\sigma\in \Gamma_{F_{v}}$ implique qu'il existe $g_{1}\in G(F_{v})$ tel que $y=i_{v}k_{v}g_{1}$. Posons $g=ad_{k_{v}}(g_{1})$. Puisque $\underline{K}_{\sharp,v}$ normalise $G(F_{v})$, on a encore $g\in G(F_{v})$.  On a $y=i_{v}gk_{v}$. D'o\`u 
  $$ad_{g^{-1}}(\eta_{\star})=ad_{k_{v}y^{-1}i_{v}}(\eta_{\star})=ad_{k_{v}y^{-1}}(\eta_{\star})\in ad_{k_{v}y^{-1}}(\tilde{K}_{v}^{nr})=\tilde{K}_{v}^{nr}.$$
   Puisque $ad_{g^{-1}}(\eta_{\star})\in \tilde{G}(F_{v})$, cela entra\^{\i}ne $ad_{g^{-1}}(\eta_{\star})\in \tilde{K}_{v}$.  Comme on l'a expliqu\'e dans la preuve du lemme 1.6, on peut appliquer les r\'esultats de [W1] 5.6. Le (ii) du lemme de cette r\'ef\'erence implique que $g\in \bar{G}_{\star}(F_{v})K_{v}$. Ecrivons $g=xk'$, avec $x\in \bar{G}_{\star}(F_{v})$ et $k'\in K_{v}$. On peut remplacer notre couple $(i_{v},k_{v})$ par $(i_{v}x,k'k_{v})$. Notons simplement $(i_{v},k_{v})$ ce nouveau couple. Il  v\'erifie les m\^emes propri\'et\'es que l'ancien mais v\'erifie de plus $y=i_{v}k_{v}$. Donc $i_{v}k_{v}\in K_{v}^{nr}$.  Puisque $k_{v}\in K_{v}^{nr}$, cela entra\^{\i}ne $i_{v}\in I_{\star}(F_{v}^{nr})\cap K_{v}^{nr}$. $\square$
  
  On suppose v\'erifi\'ee la condition (1).  Soit $v\in Val(F)$. Si $v\in V$, on a $i_{v}\sigma(i_{v})^{-1}=u(\sigma)$ donc  $i_{v}\sigma(i_{v})^{-1}\in Z(I_{\star})$ pour tout $\sigma\in \Gamma_{F_{v}}$. La m\^eme propri\'et\'e vaut si $v\not\in V$ car alors 
   $k_{v}\sigma(k_{v})^{-1}\in Z(G)^{\theta}\subset Z(I_{\star})$. Donc l'image $i_{v,ad}$ de $i_{v}$ dans $\bar{G}_{\star,AD}$ appartient \`a $\bar{G}_{\star,AD}(F_{v})$. On a des applications naturelles
   $$\bar{G}_{\star,AD}(F_{v})\to H^1(F_{v};Z(\bar{G}_{\star,SC}))\to H^1(F_{v};\bar{S}_{sc}).$$
   Pour presque tout $v$, le groupe $K_{v}\cap \bar{G}_{\star}(F_{v})$ est un sous-groupe compact hypersp\'ecial de $\bar{G}_{\star}(F_{v})$. Il d\'etermine un tel sous-groupe $\bar{K}_{\star,v,ad}$ de $\bar{G}_{\star,AD}(F_{v})$. La condition (1) entra\^{\i}ne que $i_{v,ad}$ appartient \`a ce sous-groupe pour presque tout $v$. Si de plus, $\bar{S}_{sc}$ est non ramifi\'e en $v$, on voit que l'image de $i_{v,ad}$ par l'application pr\'ec\'edente appartient au sous-groupe $H^1(\mathfrak{o}_{v};\bar{S}_{sc})$, qui est nul. On a ainsi construit pour tout $v$ un \'el\'ement de $H^1(F_{v};\bar{S}_{sc})$ qui est nul pour presque tout $v$. La collection de ces termes est un \'el\'ement de $H^1({\mathbb A}_{F};\bar{S}_{sc})$. On l'envoie dans $Q_{1}$ par l'application naturelle. Notons $u_{1}$ l'\'el\'ement de $Q_{1}$ obtenu. 
      
   Fixons comme dans la preuve de (1)  un \'el\'ement $y\in G(\bar{F})$ tel que $u(\sigma)=y\sigma(y)^{-1}$ pour tout $\sigma\in \Gamma_{F}$.  Fixons une d\'ecomposition $y=z\pi(y_{sc})$, avec $z\in Z(G;\bar{F})$ et $y_{sc}\in G_{SC}(\bar{F})$. Pour unifier les notations, posons $k_{v}=1$ pour $v\in V$.  Pour tout $v\in Val(F)$, fixons une d\'ecomposition $i_{v}k_{v}=z_{v}\pi(x_{sc,v})$, avec $z_{v}\in Z(G;\bar{F}_{v})$ et $x_{sc,v}\in G_{SC}(\bar{F}_{v})$. Pour $\sigma\in \Gamma_{F_{v}}$, les termes $y_{sc}\sigma(y_{sc})^{-1}$ et $x_{sc,v}\sigma(x_{sc,v})^{-1}$ ont m\^eme image dans $G_{AD}$: c'est l'image de $u(\sigma)$. Celle-ci appartient \`a l'image dans $G_{AD}$ de $Z(I_{\star})$, laquelle est contenu dans celle de $S_{\star,v}$. Il en r\'esulte que les termes ci-dessus appartiennent \`a $S_{\star, sc,v}(\bar{F}_{v})$ et que leur rapport appartient \`a $Z(G_{SC};\bar{F}_{v})$. Posons $\chi_{v}(\sigma)=\sigma(x_{sc,v})x_{sc,v}^{-1}y_{sc}\sigma(y_{sc})^{-1}$. L'\'egalit\'e 
   $$i_{v}k_{v}\sigma(i_{v}k_{v})^{-1}=u(\sigma)=y\sigma(y)^{-1}$$
   entra\^{\i}ne que le couple $(\chi_{v},zz_{v}^{-1})$ est un cocycle de $\Gamma_{F_{v}}$ dans le complexe $Z(G_{SC})\to Z(G)$, que l'on pousse en un cocycle \`a valeurs dans le complexe $S_{sc}\to S$. On obtient ainsi un \'el\'ement de $H^0(F_{v};S_{sc}\to S)=H^0_{ab}(F_{v};G)$. Il est clair que cet \'el\'ement ne d\'epend pas des d\'ecompositions choisies de $y$ et $i_{v}k_{v}$. Pour presque tout $v$, on a $y_{sc}\in K_{sc,v}^{nr}$ et $z\in K_{v}^{nr}$. L'hypoth\`ese (1) permet de choisir des \'el\'ements $x_{sc,v}\in K_{sc,v}$ et $z_{v}\in K_{v}^{nr}$. On v\'erifie alors que l'\'el\'ement de $H^0_{ab}(F_{v};G)$ que l'on vient de construire appartient \`a $H^0_{ab}(\mathfrak{o}_{v};G)$. La collection de ces \'el\'ements appartient \`a $H^0_{ab}({\mathbb A}_{F};G)$. On pousse cet \'el\'ement en un \'el\'ement de $Q_{2}$ que l'on note $u_{2}$. On a choisi pour le construire un \'el\'ement $y$. Mais on ne peut modifier $y$ qu'en le multipliant \`a droite par un \'el\'ement $g\in G(F)$. On v\'erifie qu'une telle multiplication multiplie l'\'el\'ement de $H^0_{ab}({\mathbb A}_{F};G)$ que l'on a construit par l'image de $g^{-1}$ dans ce groupe par la suite d'applications
   $$G(F)\to H^0_{ab}(F;G)\to H^0_{ab}({\mathbb A}_{F};G).$$
   Donc l'image $u_{2}$ dans $Q_{2}$ est inchang\'ee.

   Pour tout $v\not\in V$, $k_{v}$ a une image naturelle dans $H^0_{ab}(\mathfrak{o}_{v};G_{\sharp})$, cf. 6.2(1). On note $u_{3}$ le produit de ces \'el\'ements dans $Q_{3}$.

  Notons ${\bf q}_{\infty}(u)$ l'image dans $Q_{\infty}$ du triplet $(u_{1},u_{2},u_{3})\in Q_{\times}$.
  
  \ass{Lemme}{(i) L'\'el\'ement ${\bf q}_{\infty}(u)$ ne d\'epend pas des choix effectu\'es.
  
  (ii) L'application $u\mapsto {\bf q}_{\infty}(u)$ se quotiente en un homomorphisme de ${\mathbb U}$ dans $Q_{\infty}$.}
  
  {\bf Notation.} On notera encore ${\bf q}_{\infty}$ cet homomorphisme de ${\mathbb U}$ dans $Q_{\infty}$.
  \bigskip
  
  Preuve. Le cocycle $u$ \'etant fix\'e, les choix sont ceux des \'el\'ements $i_{v}$ et $k_{v}$. Pour une place $v\in Val(F)$,  consid\'erons d'autres choix $i'_{v}$, $k'_{v}$. Notons $u'_{1}$ etc... les analogues de $u_{1}$ etc... construits \`a l'aide de ces nouveaux choix.   Supposons d'abord $v\in V$. On a alors $k_{v}=k'_{v}=1$ et  $u'_{3}=u_{3}$. La relation $i_{v}\sigma(i_{v})^{-1}=u(\sigma)=i'_{v}\sigma(i'_{v})^{-1}$ implique qu'il existe $g\in G(F_{v})$ tel que $i'_{v}=i_{v}g$. Puisque $i_{v}$ et $i'_{v}$ appartiennent \`a $I_{\star}$, on a $g\in I_{\star}(F_{v})$. Cet \'el\'ement $g$ s'envoie en un \'el\'ement  de $Q_{1,2}$. On voit que le couple $(u'_{1},u'_{2})$ est le produit de $(u_{1},u_{2})$ et de l'image par ${\bf q}_{1,2}$ de cet \'el\'ement de $Q_{1,2}$. Donc les images dans $Q_{\infty}$ de $(u_{1},u_{2},u_{3})$ et de $(u'_{1},u'_{2},u'_{3})$ sont \'egales. Supposons maintenant $v\not\in V$. 
 On a fix\'e un \'el\'ement $h_{\star,v}\in G(F_{v})$ tel que $ad_{h_{\star, v}^{-1}}(\eta_{\star})\in \tilde{K}_{v}$.   Remarquons que, puisque $k_{v}\sigma(k_{v})^{-1}$ est central pour tout $\sigma\in \Gamma_{F_{v}}$, on a l'\'egalit\'e
 $$ h_{\star,v}k_{v} \sigma(h_{\star,v}k_{v})^{-1}=k_{v}\sigma(k_{v})^{-1}.$$
 On a donc $u(\sigma)=i_{v}h_{\star,v}k_{v}\sigma(i_{v}h_{\star,v}k_{v})^{-1}$ pour tout $\sigma\in \Gamma_{F_{v}}$. On a une relation analogue avec $i'_{v}$ et $k'_{v}$. Cela entra\^{\i}ne
  qu'il existe $g_{1}\in G(F_{v})$ tel que $i'_{v}h_{\star,v}k'_{v}=i_{v}h_{\star,v}k_{v}g_{1}$. Posons $g=ad_{k_{v}}(g_{1})$. On a encore $g\in G(F_{v})$ et on a $i'_{v}h_{\star,v}k'_{v}=i_{v}h_{\star,v}gk_{v}$. On a alors
  $$ad_{g^{-1}}\circ ad_{h_{\star,v}^{-1}}(\eta_{\star})=ad_{k_{v}(k'_{v})^{-1}}\circ ad_{h_{\star,v}^{-1}}\circ ad_{(i'_{v})^{-1}i_{v}}(\eta_{\star})=ad_{k_{v}(k'_{v})^{-1}}\circ ad_{h_{\star,v}^{-1}}(\eta_{\star})$$
  $$\in ad_{k_{v}(k'_{v})^{-1}}(\tilde{K}_{v})=\tilde{K}_{v}.$$
  D'apr\`es le lemme 5.6(ii) de [W1], cela implique $g\in G_{ad_{h_{\star,v}^{-1}}(\eta_{\star})}(F_{v})K_{v}=ad_{h_{\star,v}^{-1}}(\bar{G}_{\star}(F_{v}))K_{v}$. Ecrivons $g=ad_{h_{\star,v}^{-1}}(a)b$, avec $a\in \bar{G}_{\star}(F_{v})$ et $b\in K_{v}$. Posons $i''_{v}=i_{v}a$ et $k''_{v}=bk_{v}$. On voit que le couple $(i''_{v},k''_{v})$ est encore un choix possible, donnant naissance \`a des termes $u''_{1}$ etc... On a de plus $i'_{v}h_{\star,v}k'_{v}=i''_{v}h_{\star,v}k''_{v}$. L'\'el\'ement $a$  a une image naturelle   $\underline{a}\in Q_{1,2}$. L'\'el\'ement $b$  a une image naturelle $\underline{b}\in Q_{2,3}$. On v\'erifie que $(u''_{1},u''_{2},u''_{3})$ est le produit de $(u_{1},u_{2},u_{3})$ et de ${\bf q}_{1,2}(\underline{a}){\bf q}_{2,3}(\underline{b})^{-1}$. Donc les images dans $Q_{\infty}$ de   $(u''_{1},u''_{2},u''_{3})$ et de $(u_{1},u_{2},u_{3})$  sont \'egales. Pour simplifier les notations, on peut supposer maintenant $i''_{v}=i_{v}$ et $k''_{v}=k_{v}$. On a alors l'\'egalit\'e $i'_{v}h_{\star,v}k'_{v}=i_{v}h_{\star,v}k_{v}$. Posons  $j=i_{v}^{-1}i'_{v}=h_{\star,v}k_{v}(k'_{v})^{-1}h_{\star,v}^{-1}$. Alors $j\in I_{\star}(\bar{F}_{v})\cap ad_{h_{\star,v}}( \underline{K}_{\sharp,v})$. D'apr\`es 5.5(1), cette intersection est \'egale au produit de $Z(G)^{\theta}$ et de 
 $$ I_{\star}(\bar{F}_{v})\cap ad_{h_{\star,v}}( \underline{K}_{\sharp,v}\cap K_{v}^{nr}).$$
 Le groupe $\bar{G}_{\star}(F_{v})\cap ad_{h_{\star,v}}(K_{v})$ est un sous-groupe compact hypersp\'ecial de $\bar{G}_{\star}(F_{v})$, qui donne naissance \`a un tel sous-groupe $\bar{K}_{\star, ad,v}$ de $\bar{G}_{\star,AD}(F_{v})$. La propri\'et\'e pr\'ec\'edente entra\^{\i}ne que l'image $j_{ad}$ de $j$ dans $\bar{G}_{\star,AD}(F_{v})$ appartient \`a $\bar{K}_{\star,ad,v}$. Elle d\'efinit donc un \'el\'ement $\underline{j}$ de $Q_{1,3}$. On voit que $u'_{2}=u_{2}$ tandis que $(u'_{1},u'_{3})$ est le produit de $(u_{1},u_{3})$ et de ${\bf q}_{1,3}(\underline{j})$. De nouveau, les images dans $Q_{\infty}$ de $(u_{1},u_{2},u_{3})$ et de $(u'_{1},u'_{2},u'_{3})$ sont \'egales.  Cela prouve le (i) de l'\'enonc\'e. 
  
  Consid\'erons deux cocycles $u$ et $u'$ de $\Gamma_{F}$ dans $Z(I_{\star};\bar{F})$  qui ont m\^eme image dans $H^1(F;I_{\star})$, cette image appartenant \`a ${\mathbb U}$.   Alors on peut fixer $i\in I_{\star}(\bar{F})$ tel que $u'(\sigma)=iu(\sigma)\sigma(i)^{-1}$ pour tout $\sigma\in \Gamma_{F}$. Cette relation implique que $i\sigma(i)^{-1}\in Z(I_{\star};\bar{F})$. Des donn\'ees $i_{v}$ et $k_{v}$ \'etant fix\'ees pour tout $v$ pour le cocycle $u$, on peut choisir pour $u'$ les donn\'ees $i'_{v}=ii_{v}$ et $k'_{v}=k_{v}$. On note $u_{1}$ etc... les termes associ\'es \`a $u$ et aux donn\'ees $i_{v}$ et $k_{v}$ et $u'_{1}$ etc... ceux associ\'es \`a $u'$ et aux donn\'ees $i'_{v}$ et $k'_{v}$. On a trivialement $u'_{3}=u_{3}$. La relation $i\sigma(i)^{-1}\in Z(I_{\star};\bar{F})$ pour tout $\sigma\in \Gamma_{F}$ implique que l'image $i_{ad}$ de $i$ dans $\bar{G}_{\star,AD}$ appartient \`a $\bar{G}_{\star,AD}(F)$. On voit que $u'_{1}$ est le produit de $u_{1}$ et de l'image de $i_{ad}$ par la suite d'applications naturelles
$$\bar{G}_{\star,AD}(F)\to  H^1(F;Z(\bar{G}_{\star,SC}))\to H^1(F;\bar{S}_{sc})\to H^1({\mathbb A}_{F}/F;\bar{S}_{sc})=Q_{1}.$$
Or la derni\`ere application ci-dessus est nulle, donc $u'_{1}=u_{1}$. Dans la construction de $u_{2}$, on a choisi un \'el\'ement $y\in G(\bar{F})$ tel que $u(\sigma)=y\sigma(y)^{-1}$ pour tout $\sigma\in \Gamma_{F}$. On peut choisir pour $u'$ l'\'el\'ement $y'=iy$. On v\'erifie alors que $u'_{2}=u_{2}$. Donc ${\bf q}_{\infty}(u')={\bf q}_{\infty}(u)$.  Cela prouve que l'application ${\bf q}_{\infty}$ se quotiente en une application de ${\mathbb U}$ dans $Q_{\infty}$. 
  
  Consid\'erons deux cocycles $u$ et $u'$ de $\Gamma_{F}$ dans $Z(I_{\star};\bar{F})$  dont les images dans $H^1(F;I_{\star})$ appartiennent \`a ${\mathbb U}$. Posons $u''=uu'$. Choisissons pour toute place $v$ des donn\'ees $i_{v}$ et $k_{v}$ pour $u$ et des donn\'ees $i'_{v}$ et $k'_{v}$ pour $u'$. Pour tout $v$ et tout $\sigma\in \Gamma_{F_{v}}$, on a
  $$u''(\sigma)=u(\sigma)u'(\sigma)=i_{v}k_{v}\sigma(i_{v}k_{v})^{-1}u'(\sigma)=k_{v}\sigma(k_{v})^{-1}i_{v}\sigma(i_{v})^{-1}u'(\sigma),$$
  parce que $k_{v}\sigma(k_{v})^{-1}\in Z(G)$. Puis 
   $$u''(\sigma)=k_{v}\sigma(k_{v})^{-1}i_{v}u'(\sigma)\sigma(i_{v})^{-1}$$
   parce que $u'(\sigma)\in Z(I_{\star})$. Puis
   $$u''(\sigma)=k_{v}\sigma(k_{v})^{-1}i_{v}i'_{v}k'_{v}\sigma(i'_{v}k'_{v})^{-1}\sigma(i_{v})^{-1}= k_{v}\sigma(k_{v})^{-1}i_{v}i'_{v}k'_{v}\sigma(i_{v}i'_{v}k'_{v})^{-1}$$
   $$=i_{v}i'_{v}k'_{v}k_{v}\sigma(k_{v})^{-1}\sigma(i_{v}i'_{v}k'_{v})^{-1}$$
   toujours parce que $k_{v}\sigma(k_{v})^{-1}\in Z(G)$. D'o\`u
   $$u''(\sigma)=i_{v}i'_{v}k'_{v}k_{v}\sigma(i_{v}i'_{v}k'_{v}k_{v})^{-1}.$$
   Pour $u''$, on peut donc choisir pour donn\'ees $i''_{v}=i_{v}i'_{v}$ et $k''_{v}=k'_{v}k_{v}$. On note $u_{1}$ etc..., $u'_{1}$ etc..., $u''_{1}$ etc... les termes construits avec ces diff\'erentes donn\'ees.  Il est imm\'ediat que $u''_{1}=u_{1}u'_{1}$ et $u''_{3}=u_{3}u'_{3}$. Pour construire le terme $u_{2}$, on doit choisir un \'el\'ement $y\in G(\bar{F})$ tel que $u(\sigma)=y\sigma(y)^{-1}$ pour tout $\sigma\in \Gamma_{F}$ et  des d\'ecompositions   $y=z\pi(y_{sc})$ et $i_{v}k_{v}=z_{v}\pi(x_{sc,v})$. Pour cette derni\`ere, on peut choisir des d\'ecompositions  $i_{v}=a_{v}\pi(i_{sc,v})$ et $k_{v}=b_{v}\pi(k_{sc,v})$ avec $a_{v},b_{v}\in Z(G)$ et $i_{sc,v},k_{sc,v}\in G_{SC}$. On pose alors $z_{v}=a_{v}b_{v}$, $x_{sc,v}=i_{sc,v}k_{sc,v}$. On choisit des termes analogues pour $u'$, que l'on affecte d'un $'$. Pour $u''$, on choisit un terme $y''$ et une d\'ecomposition $y''=z''\pi(y''_{sc})$. On peut choisir $i''_{sc,v}=i_{sc,v}i'_{sc,v}$, $a''_{v}=a_{v}a'_{v}$, $k''_{sc,v}=k'_{sc,v}k_{sc,v}$, $b''_{v}=b'_{v}b_{v}$. Pour $\sigma\in \Gamma_{F_{v}}$, posons
   $$\chi_{v}(\sigma)=\sigma(x_{sc,v})x_{sc,v}^{-1}\sigma(x'_{sc,v})(x'_{sc,v})^{-1}x''_{sc,v}\sigma(x''_{sc,v})^{-1}.$$
   Parce que $k_{v}\sigma(k_{v})^{-1}\in Z(G)$, on a $k_{sc,v}\sigma(k_{sc,v})^{-1}\in Z(G_{SC})$. Notons $I_{\star,sc}$ l'image r\'eciproque dans $G_{SC}$ de l'image de $I_{\star}$ dans $G_{AD}$. On a $i_{sc,v}\in I_{\star,sc}$. Parce que $i_{v}\sigma(i_{v})^{-1}\in Z(I_{\star})$, on a $i_{sc,v}\sigma(i_{sc,v})^{-1}\in Z(I_{\star,sc})$. De m\^emes propri\'et\'es valent pour $k'_{sc,v}$, $k''_{sc,v}$, $i'_{sc,v}$ et $i''_{sc,v}$. On calcule alors
   $$(2) \qquad \chi_{v}(\sigma)=\sigma(k_{sc,v})k_{sc,v}^{-1}\sigma(k'_{sc,v})(k'_{sc,v})^{-1}k''_{sc,v}\sigma(k''_{sc,v})^{-1}$$
   $$\sigma(i_{sc,v})i_{sc,v}^{-1}\sigma(i'_{sc,v})(i'_{sc,v})^{-1}i''_{sc,v}\sigma(i''_{sc,v})^{-1}.$$
   On a
$$ k''_{sc,v}\sigma(k''_{sc,v})^{-1}= k'_{sc,v}k_{sc,v}\sigma(k_{sc,v}) ^{-1}  \sigma(k'_{sc,v})^{-1}=k'_{sc,v}\sigma(k'_{sc,v}) ^{-1}k_{sc,v}\sigma(k_{sc,v}) ^{-1}$$
et les six premiers termes de (2) disparaissent. On a aussi
$$\sigma(i_{sc,v})i_{sc,v}^{-1}\sigma(i'_{sc,v})(i'_{sc,v})^{-1}=\sigma(i_{sc,v})\sigma(i'_{sc,v})(i'_{sc,v})^{-1}i_{sc,v}^{-1}=\sigma(i''_{sc,v})(i''_{sc,v})^{-1}$$
et les six derniers termes de (2) disparaissent. D'o\`u $\chi_{v}(\sigma)=1$. 
 En appliquant les d\'efinitions, on  voit alors que  
    $u_{2}u'_{2}(u''_{2})^{-1}$ est l'image dans $Q_{2}$ du cocycle $(\xi,zz'(z'')^{-1})\in H^0_{ab}( {\mathbb A}_{F};G)$, o\`u $\xi$ est  d\'efini par 
    $$\xi(\sigma)=y_{sc}\sigma(y_{sc})^{-1}y'_{sc}\sigma(y'_{sc})^{-1}\sigma(y''_{sc})(y''_{sc})^{-1}.$$
    Ce cocycle est l'image d'un \'el\'ement de $H^0_{ab}(F;G)$. Or $H^0_{ab}(F;G)$ s'envoie sur $0$ dans $Q_{2}$. D'o\`u $u''_{2}=u_{2}u'_{2}$. Cela prouve que l'application ${\bf q}_{\infty}$, quotient\'ee en une application d\'efinie sur ${\mathbb U}$,  est un homomorphisme. $\square$
  
  \bigskip
  
  \subsection{L'image de l'homomorphisme ${\bf q}_{\infty}$}
  Rappelons que
  $$Q=H^{1,0}({\mathbb A}_{F}/F;S_{sc}\stackrel{1-\theta}{\to}(1-\theta)(S)).$$
  Pour $j=1,2,3$, on a d\'efini en 6.6 un homomorphisme ${\bf q}_{j}:Q_{j}\to Q$. On en d\'eduit 
  un homomorphisme produit
 $$(1)\qquad  Q_{\times}=Q_{1}\times Q_{2}\times Q_{3}\to Q.$$
 Montrons que
 
 (2) pour $1\leq j<j'\leq3$, le compos\'e  de l'homomorphisme (1) et de $ {\bf q}_{j,j'}$   est nul.
 
 Prouvons-le pour $(j,j')=(1,2)$, la d\'emonstration \'etant similaire pour les autres couples. Il s'agit de prouver que le diagramme suivant
$$\begin{array}{ccccc}&&Q_{1,2}&&\\ &\swarrow&&\searrow&\\ Q_{1}&&&&Q_{2}\\ &{\bf q}_{1}\searrow\quad&&\quad\swarrow{\bf q}_{2}&\\ &&Q&&\\ \end{array}$$
 est commutatif, o\`u les fl\`eches du haut sont celles d\'efinies en 7.6.

  Le compos\'e des homomorphismes de gauche est compos\'e de
 $$Q_{1,2}=I_{\star}({\mathbb A}_{F})\to H^0_{ab}({\mathbb A}_{F};I_{\star})\to H^0_{ab}({\mathbb A}_{F};\bar{G}_{\star,AD})=H^{1,0}({\mathbb A}_{F};\bar{S}_{sc}\to S^{\theta}/Z(I_{\star}))$$
 $$\to H^1({\mathbb A}_{F};\bar{S}_{sc})\to Q_{1}=H^1({\mathbb A}_{F}/F;\bar{S}_{sc})\to H^{1,0}({\mathbb A}_{F}/F;S_{sc}\stackrel{1-\theta}{\to}(1-\theta)(S)).$$
 On peut remplacer les deux derniers homomorphismes par
 $$H^1({\mathbb A}_{F};\bar{S}_{sc})\to H^{1,0}({\mathbb A}_{F};S_{sc}\stackrel{1-\theta}{\to}(1-\theta)(S))\to
 H^{1,0}({\mathbb A}_{F}/F;S_{sc}\stackrel{1-\theta}{\to}(1-\theta)(S)).$$
Le compos\'e des homomorphismes de droite du diagramme est compos\'e de 
 $$Q_{1,2}=I_{\star}({\mathbb A}_{F})\to H^0_{ab}({\mathbb A}_{F};I_{\star})\to H^0_{ab}({\mathbb A}_{F};G)\to Q_{2}=H^0_{ab}({\mathbb A}_{F};G)/Im(H^0_{ab}(F;G)$$
 $$\to H^0_{ab}({\mathbb A}_{F}/F;G)=H^{1,0}({\mathbb A}_{F}/F;S_{sc}\to S)\to H^{1,0}({\mathbb A}_{F}/F;S_{sc}\stackrel{1-\theta}{\to}(1-\theta)(S)).$$
 Il est \'egal au compos\'e de
  $$Q_{1,2}=I_{\star}({\mathbb A}_{F})\to H^0_{ab}({\mathbb A}_{F};I_{\star})\to H^0_{ab}({\mathbb A}_{F};G)=H^{1,0}({\mathbb A}_{F};S_{sc}\to S)$$
  $$\to  H^{1,0}({\mathbb A}_{F};S_{sc}\stackrel{1-\theta}{\to}(1-\theta)(S))\to
 H^{1,0}({\mathbb A}_{F}/F;S_{sc}\stackrel{1-\theta}{\to}(1-\theta)(S)).$$
 Ainsi, nos deux homomorphismes se factorisent par des homomorphismes
 $$H^0_{ab}({\mathbb A}_{F};I_{\star})\to H^{1,0}({\mathbb A}_{F};S_{sc}\stackrel{1-\theta}{\to}(1-\theta)(S)),$$
 que l'on d\'ecompose en produits d'homomorphismes locaux 
 $$H^0_{ab}(F_{v};I_{\star})\to H^{1,0}(F_{v};S_{sc}\stackrel{1-\theta}{\to}(1-\theta)(S)).$$
 On fixe $v$ et il suffit de montrer que ces homomorphismes locaux sont \'egaux. 
 Rappelons que $H^0_{ab}(F_{v};I_{\star})=H^{1,0}(F_{v};\bar{S}_{sc}\to S^{\theta})$. On v\'erifie que le premier  homomorphisme est d\'eduit du compos\'e des homomorphismes suivants de complexes de tores
 $$\begin{array}{ccc}\bar{S}_{sc}&\to &S^{\theta}\\ \downarrow &&\downarrow\\ \bar{S}_{sc}&\to&0\\ \downarrow&&\downarrow\\ S_{sc}&\stackrel{1-\theta}{\to}&(1-\theta)(S)\\ \end{array}$$
 Le second est d\'eduit du compos\'e des homomorphismes suivants de complexes de tores
$$\begin{array}{ccc}\bar{S}_{sc}&\to &S^{\theta}\\ \downarrow&&\downarrow\\  S_{sc}&\to&S\\ \downarrow&&\,\,\qquad \downarrow 1-\theta\\ S_{sc}&\stackrel{1-\theta}{\to}&(1-\theta)(S)\\ \end{array}$$
  Mais les deux homomorphismes compos\'es
   $$\begin{array}{ccc}\bar{S}_{sc}&\to &S^{\theta}\\ \downarrow&& \downarrow \\ S_{sc}&\stackrel{1-\theta}{\to}&(1-\theta)(S)\\ \end{array}$$
  sont \'egaux. Cela prouve l'\'egalit\'e de nos homomorphismes. D'o\`u (2). 
 
 Gr\^ace \`a (2), l'homomorphisme (1)  se quotiente en un homomorphisme
 ${\bf q}_{0}:Q_{\infty}\to Q$. Son image est \'egale au groupe $Q_{0}$ d\'efini en 6.6. 
 
 \ass{Lemme}{L'image de l'homomorphisme ${\bf q}_{\infty}$ est \'egale au noyau de ${\bf q}_{0}$. }
 
 Preuve. Montrons d'abord que l'image de ${\bf q}_{\infty}$ est contenue dans le noyau de ${\bf q}_{0}$. On introduit les notations suivantes pour les homomorphismes naturels
 $$\begin{array}{ccccc}\bar{G}_{\star,SC}&&\stackrel{\bar{\pi}_{\star}}{\to}&&G\\ &\bar{\pi}_{\star,sc}\searrow\,\,&&\,\,\nearrow\pi&\\&&G_{SC}&&\\ \end{array}$$
 Soit $u:\Gamma_{F}\to Z(I_{\star})$ un cocycle dont l'image dans $H^1(F;I_{\star})$ appartient \`a ${\mathbb U}$. On fixe $y\in G(\bar{F})$ tel que $u(\sigma)=y\sigma(y)^{-1}$ pour tout $\sigma\in \Gamma_{F}$. On fixe $y_{sc}\in G_{SC}(\bar{F})$ et $z\in Z(G;\bar{F})$ tels que $y=z\pi(y_{sc})$. Pour toute place $v$, on fixe des \'el\'ements $i_{v}$ et $k_{v}$ comme en 7.6. On fixe 
 
 - $z_{v}\in Z(G;\bar{F}_{v})$ et $x_{sc,v}\in G_{SC}(\bar{F}_{v})$ tels que $i_{v}k_{v}=z_{v}\pi(x_{sc,v})$;
 
 - $\zeta_{v}\in Z(I_{\star};\bar{F}_{v})$ et $\bar{i}_{sc,v}\in \bar{G}_{\star,SC}(\bar{F}_{v})$ tels que $i_{v}=\zeta_{v}\bar{\pi}_{\star}(\bar{i}_{sc,v})$;
 
 - $b_{v}\in Z(G;\bar{F}_{v} )$ et $k_{sc,v}\in G_{SC}(\bar{F}_{v})$ tels que $k_{v}=b_{v}\pi(k_{sc,v})$ (on prend $b_{v}=1$ et $k_{sc,v}=1$ pour $v\in V$).
 
 Reprenons les constructions des termes $u_{1}$, $u_{2}$ et $u_{3}$ de 7.6. De l'isomorphisme $ad_{\bar{r}_{\star}^{-1}u_{\star}}:\bar{S}_{\star}\to \bar{S}$ se d\'eduisent des isomorphismes
$$H^1({\mathbb A}_{F};\bar{S}_{\star,sc})\to H^1({\mathbb A}_{F};\bar{S}_{sc}),$$
$$ H^{1,0}({\mathbb A}_{F};S_{\star,sc}\to S_{\star})\to H^{1,0}({\mathbb A}_{F};S_{sc}\to S),$$
 $$H^{1,0}({\mathbb A}_{F};S_{\star,sc}\to  S_{\star}/Z(G)^0)\to H^{1,0}({\mathbb A}_{F};S_{sc}\to  S/Z(G)^0).$$
 Les espaces d'arriv\'ee s'envoient ensuite respectivement dans $Q_{1}$, $Q_{2}$ et $Q_{3}$. Les espaces de d\'epart s'envoient donc eux-aussi dans ces groupes.
   
 Le terme $u_{1}$ est l'image   par l'application ainsi d\'efinie de  l'\'el\'ement de $H^1({\mathbb A}_{F};\bar{S}_{\star,sc})$ dont la composante en $v$ est le cocycle $\sigma\mapsto \bar{i}_{sc,v}\sigma(\bar{i}_{sc,v})^{-1}$.   Le terme $u_{2}$ est l'image de l'\'el\'ement de $H^{1,0}({\mathbb A}_{F};S_{\star,sc}\to S_{\star})$ dont la composante en $v$ est le couple form\'e du cocycle $\sigma\mapsto y_{sc}\sigma(y_{sc})^{-1}\sigma(x_{sc,v})x_{sc,v}^{-1}$ et de l'\'element $zz_{v}^{-1}$ de $S_{\star}$.
  Le terme $u_{3}$ est l'image de  l'\'el\'ement de $H^{1,0}({\mathbb A}_{F};S_{\star,sc}\to  S_{\star}/Z(G)^0)$ dont la composante en $v$ est le couple form\'e du cocycle $\sigma\mapsto k_{sc,v}\sigma(k_{sc,v})^{-1}$ et de l'image de $b_{v}$ dans $S_{\star}/Z(G)^0$. 
 
  On envoie $u_{1}$, $u_{2}$ et $u_{3}$ dans $Q$ et on fait le produit. On obtient l'image naturelle d'un \'el\'ement de $H^{1,0}({\mathbb A}_{F};S_{\star,sc}\stackrel{1-\theta}{\to}(1-\theta)(S_{\star}))$ dont la composante en $v$ est $(\chi_{v},(1-\theta)(zz_{v}^{-1}b_{v}))$, o\`u $\chi_{v}$ est d\'efini par
 $$\chi_{v}(\sigma)=\bar{\pi}_{\star,sc}(\bar{i}_{sc,v}\sigma(\bar{i}_{sc,v})^{-1})k_{sc,v}\sigma(k_{sc,v})^{-1} y_{sc}\sigma(y_{sc})^{-1}\sigma(x_{sc,v})x_{sc,v}^{-1}.$$
 L'\'el\'ement $\sigma(x_{sc,v})x_{sc,v}^{-1}$ appartient \`a $S_{\star,sc,v}(\bar{F}_{v})$ et on ne change pas $\chi_{v}(\sigma)$ en le conjuguant par cet \'el\'ement. On obtient
 $$\chi_{v}(\sigma)=X_{v}(\sigma)y_{sc}\sigma(y_{sc})^{-1},$$
 o\`u
 $$X_{v}(\sigma)=\sigma(x_{sc,v})x_{sc,v}^{-1}\bar{\pi}_{\star,sc}(\bar{i}_{sc,v}\sigma(\bar{i}_{sc,v})^{-1})k_{sc,v}\sigma(k_{sc,v})^{-1} .$$
 Parce que $k_{v}\sigma(k_{v})^{-1}\in Z(G)$, on a $k_{sc,v}\sigma(k_{sc,v})^{-1} \in Z(G_{SC})$ et on peut r\'ecrire
 $$X_{v}(\sigma)=\sigma(x_{sc,v})x_{sc,v}^{-1}\bar{\pi}_{\star,sc}(\bar{i}_{sc,v})k_{sc,v}\sigma(k_{sc,v})^{-1}\sigma(\bar{\pi}_{\star,sc}(\bar{i}_{sc,v}))^{-1}$$
 $$=\sigma(x_{sc,v})x_{sc,v}^{-1} \tau_{v}^{-1}x_{sc,v}\sigma(x_{sc,v})^{-1}\sigma(\tau_{v}),$$
 o\`u
  $\tau_{v}=x_{sc,v}k_{sc,v}^{-1}\bar{\pi}_{\star,sc}(\bar{i}_{sc,v})^{-1}$. Il r\'esulte des d\'efinitions que $\pi(\tau_{v})=z_{v}^{-1}b_{v}\zeta_{v}$. Ce dernier terme appartient \`a $Z(G;\bar{F}_{v})Z(I_{\star};\bar{F}_{v})\subset S_{\star,v}(\bar{F}_{v})$. Donc $\tau_{v}\in S_{\star,sc,v}(\bar{F}_{v})$.  En particulier, il commute \`a $x_{sc,v}\sigma(x_{sc,v})^{-1}$ et on obtient simplement $X_{v}(\sigma)=\tau_{v}^{-1}\sigma(\tau_{v})$. L'\'el\'ement $\sigma(\tau_{v})$ commute aussi \`a $\chi_{v}(\sigma)$ et on obtient
  $$\chi_{v}(\sigma)=\sigma(\tau_{v})^{-1}\chi_{v}(\sigma)\sigma(\tau_{v})=\sigma(\tau_{v})^{-1}X_{v}(\sigma)y_{sc}\sigma(y_{sc})^{-1}\sigma(\tau_{v})=
 \tau_{v}^{-1}y_{sc}\sigma(y_{sc})^{-1}\sigma(\tau_{v}).$$
 Le couple $(\chi_{v},(1-\theta)(zz_{v}^{-1}b_{v}))$ est cohomologue au cocycle $(\chi'_{v},(1-\theta)(zz_{v}^{-1}b_{v}\pi(\tau_{v}^{-1})))$, o\`u $\chi'_{v}(\sigma)=y_{sc}\sigma(y_{sc})^{-1}$. On calcule $z_{v}^{-1}b_{v}\pi(\tau_{v}^{-1})=\zeta_{v}^{-1}$. Or cet \'el\'ement appartient \`a $Z(I_{\star})$ donc est annul\'e par $1-\theta$. Notre cocycle est donc cohomologue au cocycle $(\chi'_{v},(1-\theta)(z))$. Notons que $\chi'_{v}$ prend ses valeurs dans $Z(I_{\star, sc})$ et que $z\in Z(G)$. L'automorphisme $ad_{\bar{r}_{\star}^{-1}u_{\star}}$ est l'identit\'e sur ces groupes. Notre cocycle se transporte donc en un cocycle d\'efini par les m\^emes formules, \`a valeurs cette fois dans le complexe $S_{sc}(\bar{F}_{v})\stackrel{1-\theta}{\to}((1-\theta)(S))(\bar{F}_{v})$.  
 Il devient alors   la composante en $v$ d'un cocycle de $\Gamma_{F}$ dans le complexe $S_{sc}(\bar{F})\stackrel{1-\theta}{\to}((1-\theta)(S))(\bar{F})$. Donc son image dans $Q=H^{1,0}({\mathbb A}_{F}/F;S_{sc}\stackrel{1-\theta}{\to}(1-\theta)(S))$ est nulle. Cela prouve que l'image de ${\bf q}_{\infty}$ est contenue dans le noyau de ${\bf q}_{0}$. 
 
 D\'emontrons la r\'eciproque.  Consid\'erons des \'el\'ements $q_{j}\in Q_{j}$, pour $j=1,2,3$, tels que ${\bf q}_{0}(q_{1},q_{2},q_{3})=0$. On rel\`eve $q_{1}$ en une cocha\^{\i}ne $\dot{\beta}:\Gamma_{F}\to \bar{S}_{sc}({\mathbb A}_{\bar{F}})$ telle que $\partial \dot{\beta}$ prend ses valeurs dans $\bar{S}_{sc}(\bar{F})$. On a not\'e $\partial$ la diff\'erentielle.  D'apr\`es le lemme 6.3, l'application  $G({\mathbb A}_{F})\to H^0_{ab}({\mathbb A}_{F};G)/Im(H^0_{ab}(F;G))=Q_{2}$ est surjective.  On rel\`eve  $q_{2}$ en un \'el\'ement $g=(g_{v})_{v\in Val(F)}$ de $G({\mathbb A}_{F})$. 
 On sait que l'application naturelle $\underline{K}^V\to H^0_{ab}(\mathfrak{o}^V;G_{\sharp})=Q_{3}$ est surjective. On  rel\`eve $q_{3}$ en un \'el\'ement $k^V=\prod_{v\not\in V}k_{v}\in \underline{K}^V$.  Pour unifier les notations, on pose $k_{v}=1$ pour $v\in V$.  Fixons un ensemble fini $V''$ de places de $F$ contenant $V$ et tel que, pour $v\not\in V''$, le tore $S_{\star,v}$ soit non ramifi\'e et $g_{v}$ appartienne \`a $K_{v}$. Pour $v\not\in V''$, on sait que $K_{v}$ et $S_{\star,v}(\mathfrak{o}_{v})$ ont m\^eme image dans $H^0_{ab}(F_{v};G)$, cf. 1.5(2). Puisque seule compte l'image de $g_{v}$ dans ce groupe, on peut supposer $g_{v}\in S_{\star,v}(\mathfrak{o}_{v})$ pour $v\not\in V$. Pour la m\^eme raison, on peut supposer $k_{v}\in S_{\star,v}(\mathfrak{o}_{v}^{nr})$ pour $v\not\in V''$. Pour $v\in V'$, on choisit des \'el\'ements $g_{sc,v}, k_{sc,v}\in G_{SC}(\bar{F}_{v})$ et $a_{v},b_{v}\in Z(G;\bar{F}_{v})$ de sorte que $g_{v}=a_{v}\pi(g_{sc,v})$, $k_{v}=b_{v}\pi(k_{sc,v})$. On suppose $b_{v}=1$ et $k_{sc,v}=1$ pour $v\in V$. Pour $v\not\in V''$, on pose $g_{sc,v}=1$, $k_{sc,v}=1$ et $a_{v}=g_{v}$, $b_{v}=k_{v}$. Pour tout $v$, $a_{v}$ et $b_{v}$ sont des \'el\'ements de $S_{\star,v}(\bar{F}_{v})$.  
 Pour tout $v$ et tout $\sigma\in \Gamma_{F_{v}}$, on pose  $\gamma_{v}(\sigma)=g_{sc,v}\sigma(g_{sc,v})^{-1}$ et $\kappa_{v}(\sigma)=k_{sc,v}\sigma(k_{sc,v})^{-1}$. Ce sont des cocycles \`a valeurs dans $Z(G_{SC};\bar{F}_{v})$ et ils sont triviaux si $v\not\in V''$. Le couple $(\gamma_{v},a_{v})$ est l'image naturelle de $g_{v}$ dans $H^{1,0}(F_{v};S_{\star,sc,v}\to S_{\star,v})\simeq H^0_{ab}(F_{v};G)$. On note $\dot{a}_{v}$ l'image de $a_{v}$ dans $S(\bar{F}_{v})$ par l'isomorphisme $ad_{\bar{r}_{\star}^{-1}u_{\star}}$. Alors $(\gamma_{v},\dot{a}_{v})$ est l'image naturelle de $g_{v}$ dans $H^{1,0}(F_{v};S_{sc}\to S)\simeq H^0_{ab}(F_{v};G)$. Avec des notations similaires, $(\kappa_{v},\dot{b}_{v})$ est l'image naturelle de $k_{v}$ dans $H^{1,0}(F_{v};S_{sc}\to S/Z(G)^{\theta})\simeq H^0_{ab}(F_{v};G_{\sharp})$. Nos cocycles s'\'etendent en des cocycles ad\'eliques $\gamma$ et $\kappa$ et nos termes $a_{v}$, $\dot{a}_{v}$ etc... se regroupent en des termes ad\'eliques $a$, $\dot{a}$ etc...
 
 La condition ${\bf q}_{0}(q_{1},q_{2},q_{3})=0$ signifie que la cocha\^{\i}ne  $(\bar{\pi}_{\star, sc}(\dot{\beta})\kappa\gamma, (1-\theta)(\dot{b}\dot{a}))$, \`a valeurs dans le complexe $S_{sc}({\mathbb A}_{\bar{F}})\stackrel{1-\theta}{\to}(1-\theta)(S({\mathbb A}_{\bar{F}}))$, est cohomologue \`a une cocha\^{\i}ne \`a valeurs dans le complexe $S_{sc}(\bar{F})\stackrel{1-\theta}{\to}(1-\theta)(S(\bar{F}))$. On peut donc fixer $\dot{x}\in S_{sc}({\mathbb A}_{\bar{F}})$ de sorte que
 
   $(1-\theta)(\dot{b}\dot{a}\pi(\dot{x})^{-1})\in (1-\theta)(S(\bar{F}))$
  
  \noindent et que, en posant $\delta(\sigma)=\dot{x}\sigma(\dot{x})^{-1}\bar{\pi}_{\star,sc}(\dot{\beta}(\sigma))\kappa(\sigma)\gamma(\sigma)$, on ait
 
 (4)  pour tout $\sigma\in \Gamma_{F}$, $ \delta(\sigma)\in S_{sc}(\bar{F})$.

On peut \'ecrire $(1-\theta)(\dot{b}\dot{a}\pi(\dot{x})^{-1})=(1-\theta)(z'\pi(\dot{x}'))$, avec $z'\in Z(G;\bar{F})$ et $\dot{x}'\in S_{sc}(\bar{F})$. On peut remplacer  $\dot{x}$ par $\dot{x}\dot{x}'$. Cela ne perturbe pas la relation (4). Mais la  relation pr\'ec\'edente devient

 $(1-\theta)(\dot{b}\dot{a}\pi(\dot{x})^{-1})=(1-\theta)(z')$, avec $z'\in Z(G;\bar{F})$. 

Parce que $\dot{\beta}$ prend ses valeurs dans $\bar{S}_{sc}$, $\bar{\pi}_{\star,sc}(\dot{\beta})$ prend ses valeurs dans $S_{sc}^{\theta}$. De plus, les couples $(\gamma,(1-\theta)(\dot{a}))$ et $(\kappa,(1-\theta)(\dot{b}))$ sont des cocycles. Cela permet de calculer 

(5) $(1-\theta)\circ\pi(\delta(\sigma))=(1-\theta)(\sigma(z')(z')^{-1})$ pour tout $\sigma\in \Gamma_{F}$. 

A fortiori $(1-\theta)\circ\pi(\delta(\sigma))\in Z(G;\bar{F})$. Cela entra\^{\i}ne $\delta(\sigma)\in S_{sc}^{\theta}(\bar{F})Z(G_{SC};\bar{F})$. 
Rappellons que l'on note  $I_{\star,sc}$ l'image r\'eciproque dans $G_{SC}$ de l'image de $I_{\star}$ dans $G_{AD}$. Le groupe $S_{sc}^{\theta}(\bar{F})Z(G_{SC};\bar{F})$ est le produit de $Z(I_{\star,sc};\bar{F})$ et de $\bar{\pi}_{\star,sc}(\bar{S}_{sc}(\bar{F}))$. Seule compte pour nous l'image de $\dot{\beta}$ dans $H^1({\mathbb A}_{F}/F;\bar{S}_{sc})$. On peut modifier $\dot{\beta}$ par une cocha\^{\i}ne \`a valeurs dans $\bar{S}_{sc}(\bar{F})$. Par une telle modification, on peut donc supposer

(6) $\delta(\sigma)\in Z(I_{\star,sc};\bar{F})$ pour tout $\sigma\in \Gamma_{F}$. 

Transportons par $ad_{u_{\star}^{-1}\bar{r}_{\star}}$ la cocha\^{\i}ne $\dot{\beta}$ en une cocha\^{\i}ne $\beta$ \`a valeurs dans $\bar{S}_{\star,sc}({\mathbb A}_{\bar{F}})$. Transportons de m\^eme $\dot{x}$ en un \'el\'ement $x\in S_{\star,sc}({\mathbb A}_{\bar{F}})$. Puisque l'isomorphisme $ad_{u_{\star}^{-1}\bar{r}_{\star}}$ est \'equivariant pour les actions galoisiennes et est l'identit\'e sur $Z(I_{\star,sc})$, nos relations se conservent. C'est-\`a-dire que l'on a

 (7)  $\delta(\sigma)=x\sigma(x)^{-1}\bar{\pi}_{\star, sc}(\beta(\sigma))\kappa(\sigma)\gamma(\sigma)$ pour tout $\sigma\in \Gamma_{F}$;
 
 (8)  $(1-\theta)(ba\pi(x)^{-1})=(1-\theta)(z')$.

D'apr\`es (6),  $\partial \delta$ prend  ses valeurs dans $Z(I_{\star,sc};\bar{F})$. Mais $\partial \delta=\bar{\pi}_{\star,sc}(\partial \beta)$. Donc $\partial \beta$ prend ses valeurs dans $Z(\bar{G}_{\star,SC};\bar{F})$.  Alors $\beta$ se pousse en un cocycle \`a valeurs dans $\bar{G}_{\star,SC}({\mathbb A}_{\bar{F}})/Z(\bar{G}_{\star,SC};\bar{F})$. Parce que $\bar{G}_{\star,SC}$ est simplement connexe, le th\'eor\`eme 2.2 de [K2] dit que l'application
$$H^1(\Gamma_{F};\bar{G}_{\star,SC}(\bar{F})/Z(\bar{G}_{\star,SC} ;\bar{F}))\to H^1(\Gamma_{F};\bar{G}_{\star,SC}({\mathbb A}_{\bar{F}})/Z(\bar{G}_{\star,SC};\bar{F}))$$
est surjective. On peut donc fixer $\bar{m}\in \bar{G}_{\star,SC}({\mathbb A}_{\bar{F}})$ tel que, en posant $\bar{\alpha}(\sigma)=\bar{m}\beta(\sigma)^{-1}\sigma(\bar{m})^{-1}$, on ait $\bar{\alpha}(\sigma)\in \bar{G}_{\star,SC}(\bar{F})$ pour tout $\sigma\in \Gamma_{F}$. Posons $m=\bar{\pi}_{\star,sc}(\bar{m})$ et $\alpha(\sigma)=\bar{\pi}_{\star,sc}(\bar{\alpha}(\sigma))$. En utilisant (7), on obtient
$$\alpha(\sigma)=m\delta(\sigma)^{-1}x\sigma(x)^{-1}\kappa(\sigma)\gamma(\sigma)\sigma(m)^{-1}.$$
Mais $\delta(\sigma)\in Z(I_{\star,sc})$ commute \`a $m$ et \`a $\alpha(\sigma)$. D'o\`u
$$(9) \qquad \alpha(\sigma)\delta(\sigma)=mx\sigma(x)^{-1}\kappa(\sigma)\gamma(\sigma)\sigma(m)^{-1}.$$
La cocha\^{\i}ne $\alpha\delta$ prend ses valeurs dans $G_{SC}(\bar{F})$. C'est un cocycle car le terme de droite ci-dessus en est un. Montrons que

(10) ce cocycle est localement trivial. 

  On peut fixer une extension finie $E$ de $F$ telle que tous nos \'el\'ements et toutes nos cocha\^{\i}nes prennent leurs valeurs dans ${\mathbb A}_{E}$. Soit $v\in Val(F)$. Comme d'habitude, on note $\bar{v}$ le prolongement fix\'e de $v$ \`a $\bar{F}$ et $w$ sa restriction \`a $E$. Les termes $x$ et $m$ ont des composantes locales $x_{w}$ et $m_{w}$. Pour simplifier, on les note $x_{v}$ et $m_{v}$. 
  Ces notations seront utilis\'ees dans la suite de la preuve. Pour $\sigma\in \Gamma_{F_{v}}$, on a
 $\gamma_{v}(\sigma)=g_{sc,v}\sigma(g_{sc,v})^{-1}$ et $\kappa_{v}(\sigma)=k_{sc,v}\sigma(k_{sc,v})^{-1}$, ces deux \'el\'ements appartenant \`a $Z(G_{SC})$ et valant $1$ si $v\not\in V''$. On en d\'eduit que la composante dans $E_{w}$ de $\alpha(\sigma)\delta(\sigma)$ est \'egale \`a
 $m_{v}x_{v}k _{sc,v}g_{sc,v}\sigma(m_{v}x_{v}k_{sc,v}g_{sc,v})^{-1}$. Ce cocycle est un cobord, d'o\`u
 l'assertion (10).
 
  Parce que $G_{SC}$ est simplement connexe, l'application
$$H^1(F;G_{SC})\to H^1({\mathbb A}_{F};G_{SC})$$
est injective ([Lab2] th\'eor\`eme 1.6.9). Il en r\'esulte que l'image de $\alpha\delta$ dans $H^1(F;G_{SC})$ est triviale.  On peut donc fixer $Y_{sc}\in G_{SC}$ tel que $\alpha(\sigma)\delta(\sigma)=Y_{sc}\sigma(Y_{sc})^{-1}$ pour tout $\sigma\in \Gamma_{F}$. Le calcul de locale trivialit\'e que l'on vient de faire implique que, pour toute place $v$, il existe $h_{sc,v}\in G_{SC}(F_{v})$ tel que
$$(11)\qquad Y_{sc}=m_{v}x_{v}k_{sc,v}g_{sc,v}h_{sc,v}.$$
Posons $Y=z'\pi(Y_{sc})$. Montrons que

(12) $Y$ appartient \`a ${\cal Y}_{\star}[d_{V}]$. 
 
 Pour $\sigma\in \Gamma_{F}$, on a $Y\sigma(Y)^{-1}=\pi(\alpha(\sigma)\delta(\sigma))z'\sigma(z')^{-1}$. On a $\pi(\alpha(\sigma))\in \bar{G}_{\star}(\bar{F})\subset I_{\star}(\bar{F})$. On a $\delta(\sigma)\in Z(I_{\star,sc};\bar{F})$ donc $\pi(\delta(\sigma))\in Z(I_{\star};\bar{F})Z(G;\bar{F})$ et aussi $\pi(\delta(\sigma))z'\sigma(z')^{-1}\in Z(I_{\star};\bar{F})Z(G;\bar{F})$. La relation (5) entra\^{\i}ne que cet \'el\'ement est invariant par $\theta$. Donc il appartient \`a $Z(I_{\star};\bar{F})Z(G;\bar{F})^{\theta}$ qui est inclus dans $Z(I_{\star};\bar{F})$. Donc $Y\sigma(Y)^{-1}\in I_{\star}(\bar{F})$. Cela prouve que $Y\in {\cal Y}_{\star}$. Soit $v\in Val(F)$.  La relation (8) entra\^{\i}ne qu'il existe $\xi_{v}\in S_{\star,v}(\bar{F}_{v})^{\theta}$ tel que 
 
 (13)  $z'=\xi_{v}\pi(x_{v})^{-1}b_{v} a_{v}$. 
 
 Utilisons (11). Puisque $z'$ est central, on a
 $$Y=\pi(Y_{sc})z'=\pi(m_{v}x_{v})z'\pi(k_{sc,v} g_{sc,v}h_{sc,v})=\pi(m_{v})\xi_{v}b_{v} a_{v}\pi(k_{sc,v}g_{sc,v}h_{sc,v}).$$
 Par construction, on a $a_{v}\in Z(G;\bar{F}_{v})$ si $v\in V''$ et $k_{sc,v}=1$ si $v\not\in V''$. En tout cas, $a_{v}$ et $k_{sc,v}$ commutent. L'\'egalit\'e pr\'ec\'edente se r\'ecrit
 $$(14) \qquad Y=\pi(m_{v})\xi_{v}b_{v}\pi(k_{sc,v}) a_{v}\pi( g_{sc,v})\pi(h_{sc,v})=\pi(m_{v})\xi_{v}k_{v}g_{v}\pi(h_{sc,v}) .$$
 Cette \'egalit\'e  d\'ecompose $Y$ en le produit d'un \'el\'ement de $I_{\star}(\bar{F}_{v})$, \`a savoir $\pi(m_{v})\xi_{v}$, d'un \'el\'ement $k_{v}$ \'egal \`a $1$ si $v\in V$ et qui appartient \`a $\underline{K}_{v}$ si $v\not\in V$, et d'un \'el\'ement de $G(F_{v})$, \`a savoir $g_{v}\pi(h_{sc,v})$. C'est exactement la condition pour que $Y$ appartienne \`a ${\cal Y}_{\star}[d_{V}]$. Cela prouve (12).  
 
 Fixons $j\in I(\bar{F})$ et $g'\in G(F)$ tels que l'\'el\'ement $y=jYg'$ appartienne \`a  notre ensemble de repr\'esentants $\dot{{\cal Y}}_{\star}[d_{V}]$. Pour $\sigma\in \Gamma_{F}$, posons $u(\sigma)=y\sigma(y)^{-1}$. Alors $u$ est un cocycle de $\Gamma_{F}$ dans $Z(I_{\star})$ qui appartient \`a ${\mathbb U}$ (lemme 7.5). Pour achever la preuve du lemme, il suffit de prouver que l'image de $(q_{1},q_{2},q_{3})$ dans $Q_{\infty}$ est \'egale \`a ${\bf q}_{\infty}(u)$. Reprenons la construction de ${\bf q}_{\infty}(u)$ de 7.6. On dispose d\'ej\`a de l'\'el\'ement $y\in G(\bar{F})$ tel que $u(\sigma)=y\sigma(y)^{-1}$ pour tout $\sigma\in \Gamma_{F}$. On doit fixer pour toute place $v$ des \'el\'ements $i_{v}\in I(\bar{F}_{v})$ et $k'_{v}$ avec $k'_{v}=1$ si $v\in V$ et $k'_{v}\in \underline{K}_{\sharp,v}$ si $v\not\in V$, de sorte que $u(\sigma)=i_{v}k'_{v}\sigma(i_{v}k'_{v})^{-1}$ pour tout $\sigma\in \Gamma_{F_{v}}$. L'\'egalit\'e (14) montre que l'on peut choisir
 $i_{v}=j\pi(m_{v})\xi_{v}$ et $k'_{v}=k_{v}$. Ces \'el\'ements v\'erifient la relation 7.6(1) pour presque tout $v$. C'est clair pour $k_{v}$ puisque $k_{v}\in S_{\star,v}(\mathfrak{o}_{v}^{nr})$ pour $v\not\in V''$.  Gr\^ace \`a (13), on a $i_{v}=jz'\pi(m_{v}x_{v})b_{v}^{-1}a_{v}^{-1}$. Les termes $b_{v}$ et $a_{v}$ appartiennent \`a $K_{v}^{nr}$ pour $v\not\in V'$. Les termes $j$ et $z'$ sont d\'efinis sur $\bar{F}$, donc appartiennent \`a $K_{v}^{nr}$ pour presque tout $v$. Les termes $m_{v}$ et $x_{v}$ sont les composantes en $w$  de termes ad\'eliques donc v\'erifient la m\^eme condition. On peut donc utiliser ces termes $i_{v}$ et $k_{v}$ dans la construction de 7.6.
 
  L'\'el\'ement $u_{3}$ construit dans ce paragraphe est trivialement \'egal \`a $q_{3}$. 
 
 Avant de calculer $u_{1}$, on a besoin d'un r\'esultat pr\'eliminaire. On fixe une d\'ecomposition $j=\bar{\pi}_{\star}(j_{sc})\zeta_{j}$ avec $j_{sc}\in \bar{G}_{\star,SC}(\bar{F})$ et $\zeta_{j}\in Z(I_{\star};\bar{F})$. Introduisons le cocycle $\psi:\Gamma_{F}\to \bar{G}_{\star,SC}({\mathbb A}_{\bar{F}})$ d\'efini par $\psi(\sigma)=(j_{sc}\bar{m})^{-1}\sigma(j_{sc}\bar{m})$. Montrons que
 
 (15) $\psi$ prend ses valeurs dans $\bar{S}_{\star,sc}({\mathbb A}_{\bar{F}})$; on a $\psi(\sigma)\beta(\sigma)\in Z(\bar{G}_{\star,SC};\bar{F})$ pour tout $\sigma\in \Gamma_{F}$. 
 
 On a vu dans la preuve de (12)  que $Y\sigma(Y)^{-1}\in \pi(\alpha(\sigma))Z(I_{\star};\bar{F})$. On a aussi $jY\sigma(jY)^{-1}=y\sigma(y)^{-1}\in Z(I_{\star};\bar{F})$ par d\'efinition de $y$. Il en r\'esulte que $j\pi(\alpha(\sigma))\sigma(j)^{-1}\in Z(I_{\star};\bar{F}) $. Il en r\'esulte que $j_{sc}\bar{\alpha}(\sigma)\sigma(j_{sc})^{-1}\in Z(\bar{G}_{\star,SC};\bar{F})$. 
 En rempla\c{c}ant $\bar{\alpha}(\sigma)$ par sa valeur $\bar{m}\beta(\sigma)^{-1}\sigma(\bar{m})^{-1}$ et par inversion, on obtient $\sigma(j_{sc}\bar{m})\beta(\sigma)(j_{sc}\bar{m})^{-1}\in Z(\bar{G}_{\star,SC};\bar{F})$.  On peut aussi bien conjuguer cette relation par $(j_{sc}\bar{m})^{-1}$ et on obtient la seconde assertion de (15). La premi\`ere en r\'esulte imm\'ediatement.

 L'\'el\'ement $u_{1}$ est l'image dans $Q_{1}$ d'un cocycle ad\'elique. Calculons sa composante en une place $v$. On note $i_{v,ad}$ l'image de $i_{v}$ dans $\bar{G}_{\star,AD}$. Elle appartient \`a $\bar{G}_{\star,AD}(F_{v})$. Alors $u_{1,v}$ est   l'image de $i_{v,ad}$ par l'application
 $$\bar{G}_{\star,AD}(F_{v})\to H^1(\Gamma_{F_{v}};Z(\bar{G}_{\star,SC}))\to H^1(F_{v};\bar{S}_{sc}).$$
 On fixe une d\'ecomposition $\xi_{v}=\bar{\pi}_{\star}(\xi_{sc,v})\zeta_{v}$ avec $\xi_{sc,v}\in \bar{S}_{\star, sc,v}(\bar{F}_{v})$ et $\zeta_{v}\in Z(I_{\star};\bar{F}_{v})$. On a alors l'\'egalit\'e $i_{v}=\bar{\pi}_{\star}(j_{sc}\bar{m}_{w}\xi_{sc,v})\zeta_{j}\zeta_{v}$.   Pour $\sigma\in \Gamma_{\bar{v}}$, on a
  $$u_{1,v}(\sigma)= j_{sc}\bar{m}_{w}\xi_{sc,v}\sigma(j_{sc}\bar{m}_{w}\xi_{sc,v})^{-1}.$$
 Ce terme appartient \`a  $Z(\bar{G}_{\star ,SC};\bar{F}_{v})$
 On peut aussi bien  conjuguer le terme ci-dessus par $\sigma (\xi_{sc,v})\bar{m}_{w}^{-1}j_{sc}^{-1}$ et on obtient
 $$u_{1,v}(\sigma)=\xi_{sc,v}\sigma(j_{sc}\bar{m}_{w})^{-1}j_{sc}\bar{m}_{w}\sigma(\xi_{sc,v})^{-1},$$
 autrement dit
 $$u_{1,v}(\sigma)=\xi_{sc,v}\psi_{v}(\sigma)^{-1}\sigma(\xi_{sc,v})^{-1}.$$
 Tous les termes appartiennent \`a $S_{\star, sc,v}(\bar{F}_{v})$ et leur produit appartient \`a $Z(\bar{G}_{\star,SC};\bar{F}_{v})$. On peut aussi bien conjuguer chaque terme par $\bar{r}_{\star}^{-1}u_{\star}$. D'o\`u
 $$u_{1,v}(\sigma)=\dot{\xi}_{sc,v}\dot{ \psi}_{v}(\sigma)^{-1} \sigma(\dot{\xi}_{sc,v})^{-1},$$
 o\`u $\dot{\xi}_{sc,v}$ et $\dot{\psi}_{v}(\sigma)$ sont les images de $\xi_{sc,v}$ et $\psi_{v}(\sigma)$ dans $\bar{S}_{sc}(\bar{F}_{v})$. Cela montre que les cocycles $u_{1,v}$ et $\dot{\psi}_{v}^{-1}$ sont cohomologues. Il en r\'esulte que l'image de $u_{1}$ dans $Q_{1}$ est la m\^eme que celle du cocycle $\dot{\psi}^{-1}$. En conjuguant la propri\'et\'e (14) par $\bar{r}_{\star}^{-1}u_{\star}$, on voit que $\dot{\psi}^{-1}$ et $\dot{\beta}$ diff\`erent par une cocha\^{\i}ne \`a valeurs dans  $Z(\bar{G}_{\star,SC};\bar{F})$. Leurs images dans $Q_{1}=H^1({\mathbb A}_{F}/F;\bar{S}_{sc})$ sont donc \'egales.  Puisque l'image de $\dot{\beta}$ est $q_{1}$, cela d\'emontre que 
   $u_{1}=q_{1}$. 
 
 Pour calculer $u_{2}$, on doit fixer une d\'ecomposition $y=\pi(y_{sc})z$ avec $y_{sc}\in G_{SC}(\bar{F})$ et $z\in Z(G;\bar{F})$. Pour tout $v$, on fixe une d\'ecomposition $i_{v}k_{v}=\pi(x_{sc,v})z_{v}$ avec $x_{sc,v}\in G_{SC}(\bar{F}_{v})$ et $z_{v}\in Z(G;\bar{F}_{v})$. Alors $u_{2}$ est l'image du cocycle ad\'elique \`a valeurs dans le complexe $S_{sc}\to S$ dont la composante en $v$ est la suivante. Le premier terme est $\sigma\mapsto \sigma(x_{sc,v})x_{sc,v}^{-1}y_{sc}\sigma(y_{sc})^{-1}$ et le second est $zz_{v}^{-1}$. Notons que le premier terme est \`a valeurs centrales, il est \'egal \`a son conjugu\'e   $\sigma\mapsto x_{sc,v}^{-1}y_{sc}\sigma(x_{sc,v}^{-1}y_{sc})$. Or $y=jYg'=j\pi(m_{v})\xi_{v}k_{v}g_{v}\pi(h_{sc,v})g'=i_{v}k_{v}g_{v}\pi(h_{sc,v})g'$. Il en r\'esulte que $g_{v}\pi(h_{sc,v})g'=\pi(x_{sc,v}^{-1}y_{sc})zz_{v}^{-1}$. Donc le terme $u_{2}$ appara\^{\i}t comme l'image naturelle dans $Q_{2}$ de l'\'el\'ement $g\pi(h_{sc})g'\in G({\mathbb A}_{F})$. L'\'el\'ement $\pi(h_{sc})$ s'envoie sur $0$ dans $H^0_{ab}({\mathbb A}_{F};G)$. L'\'el\'ement $g'$ s'envoie sur un \'el\'ement de $H^0_{ab}(F;G)$, qui devient nul dans $Q_{2}$. Donc $u_{2}$ est l'image naturelle de $g$, c'est-\`a-dire $q_{2}$. Cela ach\`eve la d\'emonstration. $\square$
 
 \bigskip
 
 \subsection{Un caract\`ere de $Q_{\infty}$}
 Dans ce paragraphe, on suppose
 
 (1) $\omega$ et $\omega_{{\bf H}}$ co\"{\i}ncident sur $I_{\star}({\mathbb A}_{F})$.
 
 Puisque $\bar{S}_{sc}\simeq S_{\bar{H}}$ est un sous-tore de $\bar{H}$, l'\'el\'ement $\bar{s}\in Z(\hat{\bar{H}})^{\Gamma_{F}}$ est aussi un \'el\'ement de $\hat{\bar{S}}_{ad}^{\Gamma_{F}}$. On a un produit sur
 $$H^1({\mathbb A}_{F}/F;\bar{S}_{sc})\times \hat{\bar{S}}_{ad}^{\Gamma_{F}}.$$
 Donc $\bar{s}$ d\'efinit un caract\`ere   du premier groupe, lequel n'est autre que  $Q_{1}$. On le note $\omega_{1}$.   Relevons l'\'el\'ement ${\bf a}\in H^1(W_{F};Z(\hat{G}))/ker^1(W_{F};Z(\hat{G}))$ en un cocycle encore not\'e ${\bf a}$ \`a valeurs dans $Z(\hat{G})$. Il se pousse en un \'el\'ement de $H^{1,0}(W_{F};\hat{S}\to \hat{S}_{ad})$, qui, par les dualit\'es usuelles, d\'efinit un caract\`ere de $H^{1,0}({\mathbb A}_{F}/F;S_{sc}\to S)$. Via le plongement de $Q_{2}$ dans ce groupe, on r\'ecup\`ere un caract\`ere $\omega_{2}$ de $Q_{2}$.   On v\'erifie que le compos\'e de $\omega_{2}$ avec l'application naturelle  $G({\mathbb A}_{F})\to Q_{2}$ n'est autre que $\omega$. Puisque qu'on a d\'ej\`a dit que cette application etait surjective, cela fournit une autre d\'efinition de $\omega_{2}$. On note $\omega_{3}$ le caract\`ere trivial de $Q_{3}$. Notons $\omega_{\times}$ le caract\`ere de $Q_{\times}$ dont la composante sur $Q_{j}$ est $\omega_{j}$ pour tout $j=1,2,3$. Montrons que
 
 (2) pour $1\leq j<j'\leq3$, $\omega_{\times}$ est trivial sur l'image de ${\bf q}_{j,j'}$.
 
Preuve.  Pour $j=2$ et $j'=3$, cela r\'esulte de la non-ramification de $\omega$ hors de $V$, qui implique que $\omega$ est trivial sur $K_{v}$ pour $v\not\in V$. Pour $j=1$, cela va r\'esulter de la propri\'et\'e 
 
 (3) le compos\'e de $\omega_{1}$ et de l'homomorphisme $H^0_{ab}({\mathbb A}_{F};\bar{G}_{\star,AD})\to Q_{1}$ est le caract\`ere $\omega_{{\bf H}}$,
 
 \noindent que l'on prouvera ci-dessous.  
  L'homomorphisme $Q_{1,2}\to Q_{1}$   se factorisant par un homomorphisme naturel $Q_{1,2}\to H^0_{ab}({\mathbb A}_{F};\bar{G}_{\star,AD})$, l'assertion (2) pour $j=1$ et $j'=2$ r\'esulte, gr\^ace \`a (3), de l'hypoth\`ese que $\omega$ co\"{\i}ncide sur $I_{\star}({\mathbb A}_{F})$ avec le caract\`ere de ce groupe d\'eduit de $\omega_{{\bf H}}$.  L'homomorphisme $Q_{1,3}\to Q_{1}$   se factorisant par un homomorphisme naturel $Q_{1,3}\to H^0_{ab}({\mathbb A}_{F};\bar{G}_{\star,AD})$, l'assertion (2) pour $j=1$ et $j'=3$ r\'esulte, gr\^ace \`a (3), de la non-ramification de $\omega_{{\bf H}}$ hors de $V$ qui implique que $\omega_{{\bf H}}$ est trivial sur $K_{\star,ad,v}$ pour $v\not\in V$. Cela prouve (2). $\square$
  
  Preuve de (3). On reprend la construction de $\omega_{{\bf H}}$ donn\'ee en [I] 2.7. Elle se simplifie puisqu'il n'y a pas ici de torsion. On fixe un rel\`evement $\bar{s}_{sc}\in \hat{\bar{T}}_{sc}$ de $\bar{s}\in \hat{\bar{T}}_{ad}$. Pour tout $w\in W_{F}$, on fixe $(\bar{g}(w),w)\in \hat{ {\cal H}}$ tel que l'action galoisienne sur $\hat{\bar{H}}$ soit $w\mapsto w_{\bar{H}}=ad_{\bar{g}_{w}}\circ w_{\bar{G}}$. On fixe un rel\`evement $\bar{g}_{sc}(w)\in \hat{\bar{G}}_{SC}$ de $\bar{g}(w)$. On d\'efinit le cocycle $a_{sc}(w)=\bar{s}_{sc}\bar{g}_{sc}(w)w_{\bar{G}}(\bar{s}_{sc})^{-1}\bar{g}_{sc}(w)^{-1}$. Il est \`a valeurs dans $Z(\hat{\bar{G}}_{SC})$ et $\omega_{\bar{H}}$ est le caract\`ere de $\bar{G}_{AD}({\mathbb A}_{F})$ associ\'e \`a ce cocycle.  Remarquons que l'on a aussi $a_{sc}(w)=\bar{s}_{sc}w_{\bar{H}}(\bar{s}_{sc})^{-1}$. On peut identifier $\hat{\bar{S}}_{sc}$ au tore $\hat{\bar{T}}_{sc}$ muni d'une action galoisienne $\sigma\mapsto \sigma_{S}$. Celle-ci est de la forme $\sigma_{S}=\omega_{S,\bar{H}}(\sigma)\circ \sigma_{\bar{H}}$, o\`u $\omega_{S,\bar{H}}$ est un cocycle \`a valeurs dans $W^{\bar{H}}$. Un \'el\'ement $\omega_{S,\bar{H}}(\sigma)$ se repr\'esente comme l'action adjointe d'un \'el\'ement du centralisateur de $\bar{s}_{sc}$ dans $\hat{\bar{G}}_{SC}$. Il en r\'esulte que 
  $$w_{S}(\bar{s}_{sc})=\omega_{S,\bar{H}}(w)(w_{\bar{H}}(\bar{s}_{sc}))=\omega_{S,\bar{H}}(w)(a_{sc}(w)^{-1}\bar{s}_{sc})=a_{sc}(w)^{-1}\bar{s}_{sc}.$$
  Donc $a_{sc}(w)=\bar{s}_{sc}w_{S}(\bar{s}_{sc})^{-1}$.  Ce cocycle est l'image naturelle de l'\'el\'ement de $H^{1,0}(W_{F};\hat{\bar{S}}_{sc}\to \hat{\bar{S}}_{ad})$ dont la premi\`ere composante est $a_{sc}$ et la seconde est triviale. Mais ce cocycle est cohomologue \`a celui dont la premi\`ere composante est triviale et la seconde est $\bar{s}$.  Celui-ci d\'efinit le caract\`ere $\omega_{1}$ de $Q_{1}$. L'assertion (3) r\'esulte de cela par dualit\'e. $\square$ 
 
 Gr\^ace \`a (2), le caract\`ere $\omega_{\times}$ de $Q_{\times}$ se quotiente en un caract\`ere  not\'e $\omega_{\infty}$ de $Q_{\infty}$. 
 
 \ass{Lemme}{Si $\omega_{\infty}$ n'est pas trivial sur l'image de l'homomorphisme ${\bf q}_{\infty}$, alors $\bar{\varphi}[V',d_{V}]=0$.}
 
 Preuve. La fonction $\bar{\varphi}[V',d_{V}]$ est d\'efinie par la formule 7.1(3).  Conform\'ement aux bijections de 7.3, on la r\'ecrit en rempla\c{c}ant l'ensemble de sommation $\dot{D}_{F}[d_{V}]$ et ses \'el\'ements $d$ par $\dot{\cal Y}_{\star}[d_{V}]$ et ses \'el\'ements $y$.   L'hypoth\`ese (1) pos\'ee ci-dessus permet de simplifier la formule 7.1(3). La preuve du lemme 7.4 montre en effet que, sous cette hypoth\`ese (1), on a l'\'egalit\'e $\omega(uh[y])\bar{f}[y,u]=\omega(h[y])\bar{f}[y,1]$ pour tout $y\in \dot{{\cal Y}}_{\star}[d_{V}]$ et tout $u\in {\cal U}[V',y]$. Posons simplement $\bar{f}[y]=\bar{f}[y,1]$. On a donc
 $$(4) \qquad \bar{\varphi}[V',d_{V}]= \sum_{y\in \dot{{\cal Y}}_{\star}[d_{V}]}\omega(h[y])\bar{f}[y].$$
Soit $y\in \dot{{\cal Y}}_{\star}[d_{V}]$, notons $u_{y}:\Gamma_{F}\to Z(I_{\star};\bar{F})$ le cocycle d\'efini par $u_{y}(\sigma)=y\sigma(y)^{-1}$.  On note encore $u_{y}$ son image dans $H^1(F;I_{\star})$. C'est un \'el\'ement de ${\mathbb U}$ (lemme 7.5). On va prouver
$$(5) \qquad \omega(h[y])\bar{f}[y]=\omega_{\infty}\circ {\bf q}_{\infty}(u_{y})^{-1}\omega(h_{\star})\bar{f}_{\star},$$
o\`u on rappelle que  $h_{\star}=h[1]$ et $\bar{f}_{\star}=\bar{f}[1]$. 
En admettant ce r\'esultat, et en utilisant le lemme 7.5, la formule (4) se r\'ecrit
$$\bar{\varphi}[V',d_{V}]= \omega(h_{\star})\bar{f}_{\star}\sum_{u\in {\mathbb U}}\omega_{\infty}\circ {\bf q}_{\infty}(u)^{-1}.$$
Le lemme en r\'esulte.

Il s'agit donc de prouver (5). On consid\`ere $y$  comme ci-dessus. Soit ${\bf x}\in \bar{H}(F_{V})$ un \'el\'ement en position g\'en\'erale et proche de $1$. Par d\'efinition,
$$S^{\bar{H}}({\bf x},\bar{f}[y])=\sum_{x}\Delta_{V}[y]({\bf x},x)I^{G_{\eta[y],SC}}(x,f[y]_{sc}).$$
Ici, $x$ parcourt, modulo conjugaison par $G_{\eta[y],SC}(F_{V})$, l'ensemble des \'el\'ements de ce groupe qui correspondent \`a ${\bf x}$.  Le facteur $\Delta_{V}[y]$ est le facteur de transfert canonique associ\'e aux choix de sous-groupes compacts hypersp\'eciaux $K_{sc,v}[y]$ de $G_{\eta[y],SC}(F_{v})$ pour $v\not\in V$. On rappelle que $K_{sc,v}[y]$ est l'image r\'eciproque de $K_{v}[y]=ad_{h_{v}[y]}(K_{v})\cap G_{\eta[y]}(F_{v})$.  Enfin, on a rappel\'e dans la preuve du lemme 7.4 la construction de la fonction $f[y]_{sc}$.  
 D'apr\`es 7.3(1), $ad_{y}$ est un  isomorphisme d\'efini sur $F$ de $G_{\eta[y]}$ sur $\bar{G}_{\star}=G_{\eta_{\star}}$. Il se rel\`eve en un tel isomorphisme de $G_{\eta[y],SC}$ sur $\bar{G}_{\star,SC}$ que l'on note encore $ad_{y}$.  Cet isomorphisme envoie l'ensemble des \'el\'ements de $G_{\eta[y],SC}(F_{V})$ qui correspondent \`a ${\bf x}$ sur l'ensemble  des \'el\'ements de $\bar{G}_{\star,SC}(F_{V})$ qui correspondent \`a ${\bf x}$. Notons  ${\cal X}({\bf x})$ ce dernier ensemble. D\'efinissons presque partout sur $\bar{H}(F_{V})\times \bar{G}_{\star,SC}(F_{V})$ une fonction $\Delta'_{V}[y]$ par $\Delta'_{V}[y]({\bf x}',x')=\Delta_{V}[y]({\bf x}',ad_{y^{-1}}(x'))$. C'est le facteur de transfert   associ\'e aux sous-groupes compacts hypersp\'eciaux $ad_{y}(K_{sc,v}[y])$ de $\bar{G}_{\star,SC}(F_{v})$ pour $v\not\in V$. 
  Par simple transport de structure, on obtient
$$(6) \qquad S^{\bar{H}}({\bf x},\bar{f}[y])=\sum_{x\in {\cal X}({\bf x})}\Delta'_{V}[y]({\bf x},x)I^{\bar{G}_{\star,SC}}(x,f[y]_{sc}\circ ad_{y^{-1}}).$$

Pour toute place $v$, on fixe des \'el\'ements $i_{v}$ et $k_{v}$ associ\'es \`a $u_{y}$ comme en 7.6. En particulier, on a $u_{y}(\sigma)=i_{v}k_{v}\sigma(i_{v}k_{v})^{-1}$ pour tout $\sigma\in \Gamma_{F_{v}}$. Cela implique qu'il existe un unique $g_{v}\in G(F_{v})$ tel que $y=i_{v}k_{v}g_{v}$. Pour $v\not\in V$, les \'el\'ements $h_{v}[y]$ et $h_{\star,v}=h_{v}[1]$ de $G(F_{v})$ ont \'et\'e fix\'es tels que $ad_{h_{\star,v}^{-1}}(\eta_{\star})\in \tilde{K}_{v}$ et $ad_{h_{v}[y]^{-1}}(\eta[y])\in \tilde{K}_{v}$.  On a 
$$ad_{h_{v}[y]^{-1}}(\eta[y])\in \tilde{K}_{v}\iff ad_{h_{v}[y]^{-1}y^{-1}}(\eta_{\star})\in \tilde{K}_{v} \iff ad_{h_{v}[y]^{-1}g_{v}^{-1}k_{v}^{-1}i_{v}^{-1}}(\eta_{\star})\in \tilde{K}_{v}.$$
Puisque $ad_{i_{v}^{-1}}(\eta_{\star})=\eta_{\star}$, ces relations sont encore \'equivalentes \`a 
$ ad_{h_{v}[y]^{-1}g_{v}^{-1}k_{v}^{-1}}(\eta_{\star})\in \tilde{K}_{v}$. Posons $g'_{v}=ad_{k_{v}}(g_{v}h_{v}[y])$. On sait que $g'_{v}\in G(F_{v})$. La relation pr\'ec\'edente \'equivaut \`a $ad_{k_{v}^{-1}(g'_{v})^{-1}}(\eta_{\star})\in \tilde{K}_{v}$, ou encore $ad_{(g'_{v})^{-1}}(\eta_{\star})\in \tilde{K}_{v}$ puisque $ad_{k_{v}}$ conserve cet ensemble. On a encore l'\'equivalence
$$ad_{(g'_{v})^{-1}}(\eta_{\star})\in \tilde{K}_{v}\iff ad_{(g'_{v})^{-1}h_{\star,v}}(ad_{h_{\star,v}^{-1}}(\eta_{\star}))\in \tilde{K}_{v}.$$
Ainsi qu'on l'a d\'ej\`a fait plusieurs fois, on peut appliquer \`a $ad_{h_{\star,v}^{-1}}(\eta_{\star})$ le lemme 5.6(ii) de [W1]. Il implique $h_{\star,v}^{-1}g'_{v}\in G_{ad_{h_{\star,v}^{-1}}(\eta_{\star})}(F_{v})K_{v}$. On peut donc fixer $m_{v}\in \bar{G}_{\star}(F_{v})$ et $k'_{v}\in K_{v}$ de sorte que $h_{\star,v}^{-1}g'_{v}=ad_{h_{\star,v}^{-1}}(m_{v})k'_{v}$. Cela \'equivaut \`a $k_{v}g_{v}h_{v}[y]=m_{v}h_{\star,v}k'_{v}k_{v}$. Le groupe $ad_{y}(K_{sc,v}[y])$ est l'image r\'eciproque dans $\bar{G}_{\star,SC}(F_{v})$ de $ad_{yh_{v}[y]}(K_{v})\cap \bar{G}_{\star}(F_{v})$. On a 
$$ad_{yh_{v}[y]}(K_{v})=ad_{i_{v}k_{v}g_{v}h_{v}[y]}(K_{v})=ad_{i_{v}m_{v}h_{\star,v}k'_{v}k_{v}}(K_{v})=ad_{i_{v}m_{v}h_{\star,v}}(K_{v}),$$
puisque $ad_{k'_{v}k_{v}}$ conserve $K_{v}$. Notons $K_{\star,v}=ad_{h_{\star,v}}(K_{v})\cap \bar{G}_{\star}(F_{v})$ et $K_{\star,sc,v}$ son image r\'eciproque dans $\bar{G}_{\star,SC}(F_{v})$. On obtient l'\'egalit\'e $ad_{yh_{v}[y]}(K_{v})\cap \bar{G}_{\star}(F_{v})=ad_{i_{v}m_{v}}(K_{\star,v})$. L'automorphisme $ad_{i_{v}m_{v}}$ de $\bar{G}_{\star}$ se rel\`eve en un automorphisme de $\bar{G}_{\star,SC}$ not\'e de la m\^eme fa\c{c}on. D'o\`u l'\'egalit\'e 
$$(7) \qquad ad_{y}(K_{sc,v}[y])=ad_{i_{v}m_{v}}(K_{\star,sc,v}).$$
On note $\Delta_{\star,V}=\Delta_{V}[1]$. Rappelons que c'est le facteur de transfert sur  $\bar{H}(F_{V})\times \bar{G}_{\star,SC}(F_{V})$  associ\'e aux choix de compacts $K_{\star,sc,v}$ pour $v\not\in V$. A l'aide de la formule (7), le m\^eme calcul que dans la preuve du lemme  7.4 conduit \`a l'\'egalit\'e
$$\Delta'_{V}[y]=\Delta_{\star,V}\prod_{v\not\in V}\omega_{\bar{H}}(i_{v,ad}m_{v,ad})^{-1},$$
o\`u $i_{v,ad}$ et $m_{v,ad}$ sont les images de $i_{v}$ et $m_{v}$ dans $\bar{G}_{F,AD}(F_{v})$. Soit  $v\not\in V$. D'apr\`es l'hypoth\`ese (1), on a $\omega_{\bar{H}}(m_{v,ad})=\omega(m_{v})$. Par d\'efinition, on a
$$m_{v}=ad_{k_{v}}(g_{v}h_{v}[y])(k'_{v})^{-1}h_{\star,v}^{-1}.$$
Le caract\`ere $\omega$ est non ramifi\'e en $v$, donc trivial sur $K_{v}$. On a d\'ej\`a remarqu\'e qu'il \'etait invariant par conjugaison par un \'el\'ement de $G_{\sharp}(F_{v})$. Il en r\'esulte que $\omega(m_{v})=\omega(g_{v}h_{v}[y]h_{\star,v}^{-1})$. On rappelle que $h[y]=(h_{v}[y])_{v\not\in V}$ et que de m\^eme $h_{\star}=(h_{\star,v})_{v\not\in V}$.  La formule  plus haut se r\'ecrit
$$ \Delta'_{V}[y]=\omega(h_{\star})\omega(h[y])^{-1}\Delta_{\star,V}\prod_{v\not\in V}\omega_{\bar{H}}(i_{v,ad})^{-1}\omega(g_{v})^{-1}.$$
La formule (6) se r\'ecrit
$$(8) \qquad S^{\bar{H}}({\bf x},\bar{f}[y])=\omega(h_{\star})\omega(h[y])^{-1}\left(\prod_{v\not\in V}\omega_{\bar{H}}(i_{v,ad})^{-1}\omega(g_{v})^{-1}\right)$$
$$\sum_{x\in {\cal X}({\bf x})}\Delta_{\star,V}({\bf x},x)I^{\bar{G}_{\star,SC}}(x,f[y]_{sc}\circ ad_{y^{-1}}).$$
 La fonction $f[y]_{sc}\circ ad_{y^{-1}}$ est l'image par $\iota_{\bar{G}_{\star,SC},\bar{G}_{\star}}$ de $f[y]\circ ad_{y^{-1}}$. Pour un \'el\'ement $x\in \bar{G}_{\star}(F_{V})$ en position g\'en\'erale et proche de $1$, on a 
 $$I^{\bar{G}_{\star}}(x,\omega,f[y]\circ ad_{y^{-1}})=I^{G_{\eta[y]}}(ad_{y^{-1}}(x),\omega,f[y])$$
 $$=
 I^{\tilde{G}}(ad_{y^{-1}}(x)\eta[y],\omega,f)=I^{\tilde{G}}(ad_{y^{-1}}(x\eta_{\star}),\omega,f).$$
 Le terme $ad_{y^{-1}}(x\eta_{\star})$  est un \'el\'ement de $\tilde{G}(F_{V})$. On rappelle que $k_{v}=1$ pour $v\in V$. On a donc pr\'ecis\'ement $ad_{y^{-1}}(x\eta_{\star})=ad_{g_{V}^{-1}i_{V}^{-1}}(x\eta_{\star})$ o\`u, par exemple, $g_{V}=(g_{v})_{v\in V}$. On a 
 $$I^{\tilde{G}}(ad_{g_{V}^{-1}i_{V}^{-1}}(x\eta_{\star}),\omega,f)=\omega(g_{V})^{-1}I^{\tilde{G}}(ad_{i_{V}^{-1}}(x\eta_{\star}),\omega,f)$$
 $$=\omega(g_{V})^{-1}I^{\tilde{G}}(ad_{i_{V}^{-1}}(x)\eta_{\star},\omega,f)=\omega(g_{V}^{-1})I^{\bar{G}_{\star}}(ad_{i_{V}^{-1}}(x),\omega,f_{\star}),$$
 o\`u on rappelle que $f_{\star}=desc_{\eta_{\star}}^{\tilde{G}}(f)$. L'automorphisme $ad_{i_{V}}$ de $\bar{G}_{\star}$ est d\'efini sur $F_{V}$ et conserve le caract\`ere $\omega_{\bar{H}}$ de $\bar{G}_{\star}(F_{V})$, lequel est \'egal \`a $\omega$. Il en r\'esulte que  
  $$ I^{\bar{G}_{\star}}(ad_{i_{V}^{-1}}(x),\omega,f_{\star})= I^{\bar{G}_{\star}}(x,\omega,f_{\star}\circ ad_{i_{V}^{-1}}).$$
  En rassemblant ces calculs, on obtient l'\'egalit\'e 
  $f[y]\circ ad_{y^{-1}}=\omega(g_{V}^{-1})f_{\star}\circ ad_{i_{V}^{-1}}$.
   D'o\`u aussi $f[y]_{sc}\circ ad_{y^{-1}}=\omega(g_{V}^{-1})f_{\star,sc}\circ ad_{i_{V,ad}^{-1}}$ o\`u $ i_{V,ad}$ est  l'image de $i_{V}$ dans $\bar{G}_{\star,AD}(F_{V})$. Pour $x\in {\cal X}({\bf x})$, on a
$$I^{\bar{G}_{\star,SC}}(x,f[y]_{sc}\circ ad_{y^{-1}})=\omega(g_{V}^{-1})I^{\bar{G}_{\star,SC}}(x,f_{\star,sc}\circ ad_{i_{V,ad}^{-1}})=\omega(g_{V}^{-1})I^{\bar{G}_{\star,SC}}(ad_{i_{V,ad}^{-1}}(x),f_{\star,sc}).$$
L'automorphisme $ad_{i_{V,ad}}$ conserve ${\cal X}({\bf x})$. Par changement de variables, l'\'egalit\'e (8) se r\'ecrit
$$S^{\bar{H}}({\bf x},\bar{f}[y])=\omega(h_{\star})\omega(h[y])^{-1}\omega(g)^{-1}\sum_{x\in {\cal X}({\bf x})}\Delta_{\star,V}({\bf x},ad_{i_{V,ad}}(x))I^{\bar{G}_{\star,SC}}(x,f_{\star,sc})\prod_{v\not\in V}\omega_{\bar{H}}(i_{v,ad})^{-1},$$
o\`u $g=(g_{v})_{v\in Val(F)}$. Mais on a l'\'egalit\'e $\Delta_{\star,V}({\bf x},ad_{i_{V,ad}}(x))=\omega_{\bar{H}}(i_{V,ad}^{-1})\Delta_{\star,V}({\bf x},x)$, cf. [I] 2.7. D'o\`u
$$S^{\bar{H}}({\bf x},\bar{f}[y])=\omega(h_{\star})\omega(h[y])^{-1}\omega(g)^{-1}\omega_{\bar{H}}(i_{ad})^{-1}\sum_{x\in {\cal X}({\bf x})}\Delta_{\star,V}({\bf x}, x)I^{\bar{G}_{\star,SC}}(x,f_{\star,sc}),$$
o\`u $i=(i_{v})_{v\in Val(F)}$. La somme du membre de droite n'est autre que $S^{\bar{H}}({\bf x},\bar{f}_{\star})$. D'o\`u l'\'egalit\'e
$$\bar{f}[y]=\omega(h_{\star})\omega(h[y])^{-1}\omega(g)^{-1}\omega_{\bar{H}}(i_{ad})^{-1}\bar{f}_{\star}.$$
En reprenant les d\'efinitions, on voit que $\omega(g)\omega_{\bar{H}}(i_{ad})=\omega_{\infty}\circ {\bf q}_{\infty}(u_{y})$. L'\'egalit\'e (5) s'en d\'eduit, ce qui ach\`eve la d\'emonstration. $\square$

\bigskip

\subsection{ Preuve de la proposition 7.1}
 On suppose $D_{F}[d_{V}]$  non vide car sinon la proposition 7.1 est triviale, ainsi qu'on l'a dit en 7.3.  On utilise les constructions de 7.3. On suppose que l'ensemble de places $V'$ satisfait les conditions 7.1(2) et 7.4(2). Supposons  $\bar{\varphi}[V',d_{V}]\not=0$. D'apr\`es le lemme 7.4, $\omega$ et $\omega_{\bar{H}}$ co\"{\i}ncident sur $I_{\star}({\mathbb A}_{F})$, c'est-\`a-dire que l'hypoth\`ese (1) de 7.8 est v\'erifi\'ee. Le lemme de ce paragraphe implique que $\omega_{\infty}$ est trivial sur l'image de ${\bf q}_{\infty}$. Puisque cette image est le noyau de ${\bf q}_{0}$, cf. lemme 7.7, $\omega_{\infty}$ se quotiente en un caract\`ere de l'image de ${\bf q}_{0}$. Cette image est le groupe $Q_{0}$ de 6.6. C'est un sous-groupe ouvert, ferm\'e et d'indice fini de $Q$, cf. lemme 6.6. On v\'erifie ais\'ement que le caract\`ere ainsi d\'efini de $Q_{0}$ est continu. Il se prolonge donc en un caract\`ere $\omega_{Q}$ de $Q$. D'apr\`es 6.5(1), ce caract\`ere s'identifie \`a un \'el\'ement $p\in P$. Par construction, on a

- le compos\'e de $\omega_{Q}$ et de l'homomorphisme ${\bf q}_{1}:Q_{1}\to Q$ est le caract\`ere de $Q_{1}$ associ\'e \`a l'\'el\'ement $\bar{s}\in \hat{\bar{S}}_{ad}^{\Gamma_{F}}$;

- le compos\'e de $\omega_{Q}$ et de l'homomorphisme $G_{ab}({\mathbb A}_{F})\to Q_{2}\stackrel{{\bf q}_{2}}{\to}Q$ est le caract\`ere $\omega$;

- le compos\'e de $\omega_{Q}$ et de l'homomorphisme ${\bf q}_{3}:Q_{3}\to Q$ est trivial.

  Les m\^emes calculs qu'en 6.5 montrent que ces propri\'et\'es sont respectivement \'equivalentes \`a

- ${\bf p}_{1}(p)=\bar{s}$;

- ${\bf p}_{2}(p)={\bf a}$;

- pour tout $v\not\in V$, $res_{I_{v}}(p)$ appartient \`a l'image de $\varphi_{v}$.

Autrement dit, $p\in P({\bf H})$. Donc cet ensemble n'est pas vide. Alors ${\cal J}({\bf H})$ ne l'est pas non plus d'apr\`es la proposition 6.4. Cela prouve la proposition 7.1. $\square$

\bigskip

\subsection{Calcul d'une constante}
On a fix\'e une fois pour toutes une mesure sur $\mathfrak{A}_{\tilde{G}}$. En 5.8, on a d\'efini le r\'eseau $\mathfrak{A}_{\tilde{G},{\mathbb Z}}=Hom(X^*(G)^{\Gamma_{F},\theta},{\mathbb Z})\subset \mathfrak{A}_{\tilde{G}}$ et on a not\'e $covol(\mathfrak{A}_{\tilde{G},{\mathbb Z}})$ le volume du quotient $\mathfrak{A}_{\tilde{G}}/\mathfrak{A}_{\tilde{G},{\mathbb Z}}$. La donn\'ee ${\cal X}$ est elliptique,  cf. 5.1(2). Pour tout $d\in D_{F}[d_{V}]$, l'\'el\'ement $\eta[d]$ est donc elliptique (cf. fin de 1.2). Il en r\'esulte que $\mathfrak{A}_{G_{\eta[d]}}=\mathfrak{A}_{\tilde{G}}$. On d\'efinit dans cet espace le r\'eseau $\mathfrak{A}_{G_{\eta[d]},{\mathbb Z}}=Hom(X^*(G_{\eta[d]})^{\Gamma_{F}},{\mathbb Z})$ et on note 
$covol(\mathfrak{A}_{G_{\eta[d]},{\mathbb Z}})$ le volume du quotient $\mathfrak{A}_{\tilde{G}}/\mathfrak{A}_{G_{\eta[d]},{\mathbb Z}}$. 

\ass{Proposition}{Supposons $D_{F}[d_{V}]$ non vide. Alors, pour tout $d\in D_{F}[d_{V}]$,  on a l'\'egalit\'e
$$\vert P^0\vert \vert \dot{D}_{F}[d_{V}]\vert ^{-1} = C(\tilde{G})^{-1}\tau(G_{\eta[d]})[ I_{\eta[d]}(F):G_{\eta[d]}(F)]^{-1}[ I_{\eta[d]}(F_{V}):G_{\eta[d]}(F_{V})]  covol(\mathfrak{A}_{G_{\eta[d]},{\mathbb Z}})^{-1}.$$}
On renvoie \`a 5.8 pour la d\'efinition de $C(\tilde{G})$. La preuve occupe les paragraphes 7.11 \`a 7.15. 

\bigskip

\subsection{Calcul de $\vert P^0\vert  $}
Labesse d\'efinit des groupes de cohomologie ab\'elienne de $I_{\star}\backslash G$, cf. [Lab1] 3.3. Consid\'erons  comme en 7.3 un sous-tore maximal $\bar{T}_{\star}$ de $I_{\star}$ d\'efini sur $F$, notons $T_{\star}$ son commutant dans $G$ et introduisons les images r\'eciproques $\bar{T}_{\star,sc}$ de $\bar{T}_{\star}$ dans $\bar{G}_{\star,SC}$ et $T_{\star,sc}$ de $T_{\star}$ dans $G_{SC}$. On a un complexe de tores
$$\bar{T}_{\star,sc}\to T_{\star}\times T_{\star,sc}\to (1-\theta)(T_{\star})\times T_{\star}.$$
En notant $\bar{\pi}_{\star}:\bar{G}_{\star,SC}\to G$, $\bar{\pi}_{\star,sc}:\bar{G}_{\star,SC}\to G_{SC}$ et $\pi:G_{SC}\to G$ les homomorphismes naturels, les fl\`eches sont
$$x\mapsto (\bar{\pi}_{\star}(x),-\bar{\pi}_{\star,sc}(x))$$
pour la premi\`ere et
$$(y,z)\mapsto ((1-\theta)(y),y+\pi(z))$$
pour la seconde. 
  On pose $$H^0_{ab}(F;I_{\star}\backslash G)=H^{2,1,0}(F;\bar{T}_{\star,sc}\to T_{\star}\times T_{\star,sc}\to (1-\theta)(T_{\star})\times  T_{\star}).$$
On d\'efinit de m\^eme $H^0_{ab}({\mathbb A}_{F};I_{\star}\backslash G)$, $H^0_{ab}({\mathbb A}_{F}/F;I_{\star}\backslash G)$ et $H^0_{ab}(F_{v};I_{\star}\backslash G)$ pour $v\in Val(F)$.   Pour $v\not\in V$, on peut choisir $T_{\star}$ non ramifi\'e en $v$ et on pose 
$$H^0_{ab}(\mathfrak{o}_{v};I_{\star}\backslash G)=H^{2,1,0}(\mathfrak{o}_{v};\bar{T}_{\star,sc}\to T_{\star}\times T_{\star,sc}\to (1-\theta)(T_{\star})\times T_{\star}).$$
On pose aussi $H^0_{ab}(\mathfrak{o}^V;I_{\star}\backslash G)=\prod_{v\not\in V}H^0_{ab}(\mathfrak{o}_{v};I_{\star}\backslash G)$.
Des notations analogues seront utilis\'ees dans la suite.
Ces d\'efinitions ne d\'ependent pas du choix de $\bar{T}_{\star}$, \`a isomorphismes canoniques pr\`es. 

 Les isomorphismes
$$\begin{array}{ccc}T_{\star}\times T_{\star,sc}&\to &T_{\star}\times T_{\star,sc}\\ (y,z)&\mapsto&( y+\pi(z),-z)\\ \end{array}$$
et
$$\begin{array}{ccc}(1-\theta)(T_{\star})\times T_{\star}&\to &T_{\star}\times (1-\theta)(T_{\star})\\ (y',z')&\mapsto&(z',y'-(1-\theta)(z'))\\ \end{array}$$
fournissent un isomorphisme entre le complexe (1) et la somme des deux complexes
$$(2) \qquad \bar{T}_{\star,sc}\stackrel{\bar{\pi}_{\star,sc}}{\to} T_{\star,sc}\stackrel{1-\theta}{\to} (1-\theta)(T_{\star}) $$
et $T_{\star}\stackrel{id}{\to }T_{\star}$. Puisque le deuxi\`eme complexe est cohomologiquement trivial, la cohomologie de $I_{\star}\backslash G$ peut se d\'efinir \`a l'aide du complexe (2). 
Comme en 6.2, on voit que l'on peut remplacer le complexe (2) par
$$(3) \qquad \bar{S}_{sc}\stackrel{\bar{\pi}_{\star,sc}}{\to} S_{sc}\stackrel{1-\theta}{\to} (1-\theta)(S).$$
En effet, les deux complexes sont quasi-isomorphes au complexe
$$Z(\bar{G}_{\star,SC})\stackrel{\bar{\pi}_{\star,sc}}{\to} Z(I_{\star,sc})\stackrel{1-\theta}{\to} (1-\theta)(Z(G)).$$

On a une suite exacte
$$H^0_{ab}(F;G)\to H^0_{ab}(F;I_{\star}\backslash G)\to H^1_{ab}(F;I_{\star})\to H^1_{ab}(F;G).$$
Labesse d\'efinit
$\mathfrak{E}(I_{\star},G;F)$ comme le conoyau de la premi\`ere fl\`eche, ou encore le noyau de la troisi\`eme. Cf. [Lab1] 3.3. On d\'efinit de m\^eme $\mathfrak{E}(I_{\star},G;{\mathbb A}_{F})$, $\mathfrak{E}(I_{\star},G;F_{v})$ pour $v\in Val(F)$ et $\mathfrak{E}(I_{\star},G;\mathfrak{o}_{v})$ pour $v\not\in V$. 

{\bf Attention.} Labesse d\'efinit par contre $\mathfrak{E}(I_{\star},G;{\mathbb A}_{F}/F)$ comme le conoyau de l'homomorphisme compos\'e
$$H^0_{ab}({\mathbb A}_{F};G)\to H^0_{ab}({\mathbb A}_{F}/F,G)\to H^0_{ab}({\mathbb A}_{F}/F;I_{\star}\backslash G).$$

\bigskip
  On a un diagramme commutatif
$$\begin{array}{ccccccc}H^0_{ab}(\mathfrak{o}^V;G)&\to& H^0_{ab}(\mathfrak{o}^V;I_{\star}\backslash G)&\to& \mathfrak{E}(I_{\star},G;\mathfrak{o}^V)&\to& 1\\ \downarrow&&\downarrow&&\downarrow&& \\ H^0_{ab}({\mathbb A}_{F};G)&\to& H^0_{ab}({\mathbb A}_{F};I_{\star}\backslash G)&\to& \mathfrak{E}(I_{\star},G;{\mathbb A}_{F})&\to& 1\\ \downarrow&&\downarrow&&\downarrow&& \\ H^0_{ab}({\mathbb A}_{F};G)&\to& H^0_{ab}({\mathbb A}_{F}/F;I_{\star}\backslash G)&\to& \mathfrak{E}(I_{\star},G;{\mathbb A}_{F}/F)&\to& 1\\ \end{array}$$
dont les suites horizontales sont exactes. Notons $e^V:\mathfrak{E}(I_{\star},G;\mathfrak{o}^V)\to \mathfrak{E}(I_{\star},G;{\mathbb A}_{F}/F)$ le compos\'e des deux derni\`eres fl\`eches verticales.

\ass{Lemme}{On a l'\'egalit\'e $\vert P^0\vert =\vert \mathfrak{E}(I_{\star},G;{\mathbb A}_{F}/F)/Im(e^V)\vert $.}

Preuve. Comme on l'a dit ci-dessus, on peut utiliser le complexe (3) pour calculer la cohomologie ab\'elienne de $I_{\star}\backslash G$. Il s'en d\'eduit une suite exacte de cohomologie
$$H^1({\mathbb A}_{F}/F;\bar{S}_{sc})\to H^{1,0}({\mathbb A}_{F}/F;S_{sc}\stackrel{1-\theta}{\to}(1-\theta)(S))\stackrel{\iota}{\to} H^0_{ab}({\mathbb A}_{F}/F;I_{\star}\backslash G)\to H^2({\mathbb A}_{F}/F;\bar{S}_{sc}).$$
On se rappelle que $\bar{S}_{sc}\simeq \bar{S}_{\bar{H}}$ est un sous-tore elliptique de $\bar{H}$ et que ${\bf H}$ est une donn\'ee endoscopique elliptique de $\bar{G}_{\star,SC}$. Il en r\'esulte que $\bar{S}_{sc}$ est elliptique, donc $H^2({\mathbb A}_{F}/F;\bar{S}_{sc})=0$ d'apr\`es les isomorphismes de Tate-Nakayama ( [K3] 3.4.2.1). Avec les notations des paragraphes pr\'ec\'edents, la suite exacte ci-dessus se r\'ecrit 
$$Q_{1}\stackrel{{\bf q}_{1}}{\to }Q\stackrel{\iota}{\to} H^0_{ab}({\mathbb A}_{F}/F;I_{\star}\backslash G)\to 0$$
On voit facilement que l'homomorphisme naturel
$$(4) \qquad H^0_{ab}({\mathbb A}_{F};G)\to H^0_{ab}({\mathbb A}_{F}/F;I_{\star}\backslash G)$$
est le compos\'e de l'homomorphisme naturel du groupe de d\'epart dans $Q_{2}$ et de $\iota\circ {\bf q}_{2}$.  Notons $Q'$ le quotient de $Q$  par le sous-groupe engendr\'e par ${\bf q}_{1}(Q_{1})$ et ${\bf q}_{2}(Q_{2})$. On obtient   que $\iota$ se quotiente en un isomorphisme entre  $Q'$   et le conoyau de (4), c'est-\`a-dire $\mathfrak{E}(I_{\star},G;{\mathbb A}_{F}/F)$. Montrons que

(5)  cet isomorphisme envoie l'image dans $Q'$ de ${\bf q}_{3}(Q_{3})$ sur $Im(e^V)$.

 Du diagramme
$$\begin{array}{ccccc}&&S_{sc}&\to&S/Z(G)^{\theta}\\&&\parallel&&\qquad \downarrow 1-\theta\\
 \bar{S}_{sc}&\stackrel{\bar{\pi}_{\star, sc}}{\to} &S_{sc}&\stackrel{1-\theta}{\to}& (1-\theta)(S)\\ \end{array}$$
 de complexes de tores se d\'eduit un homomorphisme
 $$H^0_{ab}({\mathbb A}_{F};G_{\sharp})\to H^0_{ab}({\mathbb A}_{F};I_{\star}\backslash G).$$
 Il se restreint en un homomorphisme
 $$(6)\qquad H^0_{ab}(\mathfrak{o}^V;G_{\sharp})\to H^0_{ab}(\mathfrak{o}^V;I_{\star}\backslash G).$$
 On obtient un diagramme
 $$\begin{array}{ccc}H^0_{ab}(\mathfrak{o}^V;G_{\sharp})&\to& H^0_{ab}(\mathfrak{o}^V;I_{\star}\backslash G)\\ \downarrow&&\downarrow\\ H^0_{ab}({\mathbb A}_{F};G_{\sharp})&\to &H^0_{ab}({\mathbb A}_{F};I_{\star}\backslash G)\\   \downarrow&&\downarrow\\  H^{1,0}({\mathbb A}_{F};S_{sc}\stackrel{1-\theta}{\to}(1-\theta)(S))&\to& H^0_{ab}({\mathbb A}_{F};I_{\star}\backslash G)\\ \downarrow&&\downarrow\\ 
 H^{1,0}({\mathbb A}_{F}/F;S_{sc}\stackrel{1-\theta}{\to}(1-\theta)(S))&\to& H^0_{ab}({\mathbb A}_{F}/F;I_{\star}\backslash G)\\ \downarrow&&\downarrow\\
Q'&\stackrel{\iota}{\to}&\mathfrak{E}(I_{\star},G;{\mathbb A}_{F}/F)\\ \end{array}$$
Il est clair que ce diagramme est commutatif. Par d\'efinition, l'image de ${\bf q}_{3}(Q_{3})$ dans $Q'$ est l'image de la suite verticale de gauche tandis que $Im(e^V)$ est l'image de celle de droite. Pour prouver (5), il suffit donc de prouver que l'homomorphisme (6) est surjectif. On peut aussi bien fixer $v\not\in V$ et d\'emontrer que l'analogue local de (6) est surjectif. Comme on l'a dit, on peut remplacer le tore $S$ et ses divers avatars $S_{sc}$ etc... par un tore $T_{\star}$ et ses avatars $T_{\star,sc}$ etc... comme au d\'ebut du paragraphe et on peut supposer que ce tore est non ramifi\'e en $v$.  L'analogue local de (6) se factorise en
$$(7) \qquad H^0_{ab}(\mathfrak{o}_{v};G_{\sharp})=H^{1,0}(\mathfrak{o}_{v};T_{\star,sc}\to T_{\star}/Z(G)^0)\to H^{1,0}(\mathfrak{o}_{v};T_{\star,sc}\stackrel{1-\theta}{\to}(1-\theta)(T_{\star}))$$
$$\to H^{1,0}(\mathfrak{o}_{v};\bar{T}_{\star,sc}\to T_{\star,sc}\stackrel{1-\theta}{\to}(1-\theta)(T_{\star}))=H^0_{ab}(\mathfrak{o}_{v};I_{\star}\backslash G).$$
La deuxi\`eme fl\`eche est surjective car son conoyau s'envoie injectivement dans $H^2(\mathfrak{o}_{v};\bar{T}_{\star,sc})$ qui est nul ([KS] lemme C.1.A). L'homomorphisme entre complexes de tores donnant naissance \`a la premi\`ere fl\`eche se compl\`ete en une suite exacte de complexes de tores
$$\begin{array}{ccc}&&1\\ &&\downarrow\\ 1&\to&T_{\star}^{\theta}/Z(G)^{\theta}\\ \downarrow&&\downarrow\\ T_{\star,sc}&\to& T_{\star}/Z(G)^0\\ \downarrow&&\qquad \downarrow 1-\theta\\ T_{\star,sc}&\stackrel{1-\theta}{\to}&(1-\theta)(T_{\star}))\\ \downarrow&&\downarrow\\ 1&&1\\ \end{array}$$
Le conoyau de la premi\`ere fl\`eche de (7) s'envoie donc injectivement dans $H^1(\mathfrak{o}_{v};T_{\star}^{\theta}/Z(G)^{\theta})$ qui est nul car $T_{\star}^{\theta}/Z(G)^{\theta}$ est connexe ([KS] lemme C.1.A). Cette fl\`eche est donc surjective. La surjectivit\'e des deux fl\`eches de (7) entra\^{\i}ne l'assertion cherch\'ee. Cela prouve (5). 

Cette assertion prouve que le nombre d'\'el\'ements de $\mathfrak{E}(I_{\star},G;{\mathbb A}_{F}/F)/Im(e^V)$ est \'egal \`a celui du quotient de $Q$ par le sous-groupe engendr\'e par les ${\bf q}_{j}(Q_{j})$ pour $j=1,2,3$. Ce sous-groupe n'est autre que $Q_{0}$. Mais le lemme 6.6 dit que $P^0$ et $Q/Q_{0}$ sont des groupes (finis) duaux. D'o\`u
 l'\'egalit\'e
$$\vert P^0\vert =\vert Q/Q_{0}\vert .$$
 Le lemme en r\'esulte. $\square$
 
 \bigskip
 
 \subsection{Un premier calcul de $\vert P^0\vert \vert {\mathbb U}\vert ^{-1}$}
 On d\'efinit usuellement le nombre de Tamagawa $\tau(G)$ de $G$, qui est calcul\'e par la formule rappel\'ee en 3.2. Labesse \'etend la d\'efinition aux groupes quasi-connexes (Lab1] 1.2). On dispose donc du nombre de Tamagawa $\tau(I_{\star})$ de $I_{\star}$. Labesse d\'efinit aussi le groupe $D(I_{\star},G)$ comme le conoyau de l'homomorphisme
 $$H^1_{ab}({\mathbb A}_{F}/F;I_{\star})\to H^1_{ab}({\mathbb A}_{F}/F;G).$$
 C'est un groupe fini dont on note $d(I_{\star},G)$ le nombre d'\'el\'ements. 
 
 \ass{Lemme}{On a l'\'egalit\'e 
 $$\vert P^0\vert \vert {\mathbb U}\vert ^{-1}=\tau(I_{\star})\tau(G)^{-1}d(I_{\star},G)\vert (I_{\star}/\bar{G}_{\star})(F_{V})\vert .$$}

Preuve. Soit $v\in Val(F)$.   On peut  calculer le groupe $H^1_{ab}(F_{v};I_{\star})$  \`a l'aide du complexe 
$$\bar{S}_{sc}\to S\stackrel{1-\theta}{\to} (1-\theta)(S).$$
Celui-ci est \'equivalent au complexe
$$\bar{S}_{sc}\to S^{\theta}$$
Ici, le groupe $S^{\theta}$ n'est plus un tore mais c'est un groupe diagonalisable. En utilisant ce complexe, on a un homomorphisme naturel
$$H^1_{ab}(F_{v};I_{\star})\to H^1_{ab}(F_{v};S^{\theta}/S^{\theta,0}).$$
Ces deux groupes sont finis. 
Labesse les munit en [Lab1] 2.3 de mesures. On s'aper\c{c}oit en utilisant sa d\'efinition que la mesure d'un point est la m\^eme dans chacun des groupes (deuxi\`eme suite exacte de la page 417 de loc. cit.). Dans le deuxi\`eme groupe, la mesure d'un point est $\vert (S^{\theta}/S^{\theta,0})(F_{v})\vert ^{-1}$. L'homomorphisme naturel  $S^{\theta}/S^{\theta,0}\to I_{\star}/\bar{G}_{\star}$ est un isomorphisme. La mesure d'un point dans $H^1_{ab}(F_{v};I_{\star})$ est donc $\vert (I_{\star}/\bar{G}_{\star})(F_{v})\vert ^{-1}$. Par d\'efinition, $\mathfrak{E}(I_{\star},G;F_{v})$ est un sous-groupe de $H^1(F_{v};I_{\star})$. On le munit de la mesure induite. La mesure d'un point est donc la m\^eme que pr\'ec\'edemment.  

Supposons $v\not\in V$. L'homomorphisme $\mathfrak{E}(I_{\star},G;\mathfrak{o}_{v})\to \mathfrak{E}(I_{\star},G;F_{v})$ est injectif. En effet, les deux groupes  se plongent respectivement dans $H^1_{ab}(\mathfrak{o}_{v};I_{\star})$ et $H^1_{ab}(F_{v};I_{\star})$ et l'homomorphisme $H^1_{ab}(\mathfrak{o}_{v};I_{\star}) \to H^1_{ab}(F_{v};I_{\star})$ est injectif. Identifions $\mathfrak{E}(I_{\star},G;\mathfrak{o}_{v})$ \`a un sous-groupe de $\mathfrak{E}(I_{\star},G;F_{v})$. Montrons que

(1) $mes(\mathfrak{E}(I_{\star},G;\mathfrak{o}_{v}))=1$.

D'apr\`es la d\'efinition de la mesure, cela \'equivaut \`a

$$(2) \qquad \vert \mathfrak{E}(I_{\star},G;\mathfrak{o}_{v})\vert=\vert (I_{\star}/\bar{G}_{\star})(F_{v})\vert .$$

On utilise ici la d\'efinition de $\mathfrak{E}(I_{\star},G;\mathfrak{o}_{v})$ comme conoyau de
$$H^0_{ab}(\mathfrak{o}_{v};G)\to H^0_{ab}(\mathfrak{o}_{v};I_{\star}\backslash G).$$
On fixe un  tore $\bar{T}_{\star}$ comme en 7.11, non ramifi\'e en $v$. L'homomorphisme ci-dessus devient
$$H^{1,0}_{ab}(\mathfrak{o}_{v};T_{\star,sc}\to T_{\star})\to H^{2,1,0}_{ab}(\mathfrak{o}_{v};\bar{T}_{\star,sc}\to T_{\star,sc}\stackrel{1-\theta}{\to}(1-\theta)(T_{\star})).$$
On a une suite de cohomologie
$$H^1_{ab}(\mathfrak{o}_{v};\bar{T}_{\star,sc})\to H^{1,0}_{ab}(\mathfrak{o}_{v};T_{\star,sc}\stackrel{1-\theta}{\to}(1-\theta)(T_{\star}))$$
$$\to H^{2,1,0}_{ab}(\mathfrak{o}_{v};\bar{T}_{\star,sc}\to T_{\star,sc}\stackrel{1-\theta}{\to}(1-\theta)(T_{\star}))\to H^2_{ab}(\mathfrak{o}_{v};\bar{T}_{\star,sc}).$$
Les deux groupes extr\^emes sont nuls. Donc la fl\`eche centrale est un isomorphisme et $\mathfrak{E}(I_{\star},G;\mathfrak{o}_{v})$ devient le conoyau de l'homomorphisme
$$H^{1,0}_{ab}(\mathfrak{o}_{v};T_{\star,sc}\to T_{\star})\to H^{1,0}_{ab}(\mathfrak{o}_{v}; T_{\star,sc}\stackrel{1-\theta}{\to}(1-\theta)(T_{\star})).$$
Ces deux groupes se calculent comme respectivement $T_{\star}(\mathfrak{o}_{v})/\pi(T_{\star,sc}(\mathfrak{o}_{v}))$ et 

\noindent $((1-\theta)(T_{\star}))(\mathfrak{o}_{v})/(1-\theta)\circ \pi(T_{\star,sc}(\mathfrak{o}_{v}))$, cf. [KS] lemme C.1.A. Donc $\mathfrak{E}(I_{\star},G;\mathfrak{o}_{v})$ s'identifie au conoyau de
$$T_{\star}(\mathfrak{o}_{v})\stackrel{1-\theta}{\to}((1-\theta)(T_{\star}))(\mathfrak{o}_{v}).$$
De la suite exacte
$$1\to T_{\star}^{\theta}(\mathfrak{o}_{v}^{nr})\to T_{\star}(\mathfrak{o}_{v}^{nr})\to  ((1-\theta)(T_{\star}))(\mathfrak{o}_{v}^{nr})\to 1$$
(o\`u $T_{\star}^{\theta}(\mathfrak{o}_{v}^{nr})=T_{\star}^{\theta}\cap T_{\star}(\mathfrak{o}_{v}^{nr})$) et de la nullit\'e des $H^1$ pour les groupes connexes r\'esulte que le conoyau ci-dessus s'identifie \`a $H^1(\Gamma_{v}^{nr};T_{\star}^{\theta}(\mathfrak{o}_{v}^{nr}))$. On note ce groupe $H^1(\mathfrak{o}_{v};T_{\star}^{\theta})$. On a une suite exacte
$$H^1(\mathfrak{o}_{v};T_{\star}^{\theta,0})\to H^1(\mathfrak{o}_{v};T_{\star}^{\theta})\to H^1(\mathfrak{o}_{v};T_{\star}^{\theta}/T_{\star}^{\theta,0})\to H^2(\mathfrak{o}_{v};T_{\star}^{\theta,0}),$$
avec une d\'efinition \'evidente du troisi\`eme groupe. Les deux groupes extr\^emes sont nuls. Donc la fl\`eche centrale est un isomorphisme. Le groupe $T_{\star}^{\theta}(\mathfrak{o}_{v}^{nr})/T_{\star}^{\theta,0}(\mathfrak{o}_{v}^{nr})$ est fini. Le lemme 5.5 de [W1] implique qu'il est isomorphe \`a $T_{\star}^{\theta}(\bar{F}_{v})/T_{\star}^{\theta,0}(\bar{F}_{v})$, lequel est isomorphe \`a $I_{\star}(\bar{F}_{v})/\bar{G}_{\star}(\bar{F}_{v})$. Il r\'esulte de ces deux faits que les groupes suivants on m\^eme nombre d'\'el\'ements: $H^1(\mathfrak{o}_{v};T_{\star}^{\theta}/T_{\star}^{\theta,0})$, $H^0(\mathfrak{o}_{v};T_{\star}^{\theta}/T_{\star}^{\theta,0})$, $H^0(F_{v};T_{\star}^{\theta}/T_{\star}^{\theta,0})$, $H^0(F_{v};I_{\star}/\bar{G}_{\star})$.  
En rassemblant ces \'egalit\'es, on obtient (2), d'o\`u (1). 

On a une suite exacte  
$$(3) \qquad 1\to B_{ab}(I_{\star},G)\to \mathfrak{E}(I_{\star},G;F)\to \mathfrak{E}(I_{\star},G;{\mathbb A}_{F})\to \mathfrak{E}(I_{\star},G;{\mathbb A}_{F}/F)$$
(premi\`ere suite de la page 427 de [Lab1]). 
On n'aura pas besoin de conna\^{\i}tre le groupe $B_{ab}(I_{\star},G)$, disons seulement qu'il est fini. Le conoyau de la  derni\`ere fl\`eche est fini. Les deux derniers groupes sont munis de topologies. La derni\`ere fl\`eche envoie $ \mathfrak{E}(I_{\star},G;{\mathbb A}_{F})$ sur un sous-groupe ouvert de $ \mathfrak{E}(I_{\star},G;{\mathbb A}_{F}/F)$ et ce dernier groupe est compact. Le produit sur $v\in Val(F)$ des mesures locales d\'efinies ci-dessus  donne une mesure sur $\mathfrak{E}(I_{\star},G;{\mathbb A}_{F})$. 
 On munit les deux premiers groupes de la suite (1) de la mesure de comptage et le dernier de la mesure compatible avec cette suite et les mesures d\'efinies sur les autres groupes. En utilisant cette mesure, on peut r\'ecrire le lemme 7.11 sous la forme
 $$(4) \qquad \vert P^0\vert =mes(\mathfrak{E}(I_{\star},G;{\mathbb A}_{F}/F))/mes(Im(e^V)).$$
Rappelons que $e^V$ est le compos\'e de la suite
$$\mathfrak{E}(I_{\star},G;\mathfrak{o}^V)\to \mathfrak{E}(I_{\star},G;{\mathbb A}_{F})\to \mathfrak{E}(I_{\star},G;{\mathbb A}_{F}/F)$$
On a vu ci-dessus que la premi\`ere fl\`eche \'etait injective. On s'en sert pour identifier $\mathfrak{E}(I_{\star},G;\mathfrak{o}^V)$ \`a un sous-groupe de $\mathfrak{E}(I_{\star},G;{\mathbb A}_{F})$. Notons ${\mathbb U}_{ab}$ l'image r\'eciproque de $\mathfrak{E}(I_{\star},G;\mathfrak{o}^V)$ dans $\mathfrak{E}(I_{\star},G;F)$ par l'homomorphisme de la suite (3). Il r\'esulte des d\'efinitions des mesures que
$$(5) \qquad mes(Im(e^V))=mes(\mathfrak{E}(I_{\star},G;\mathfrak{o}^V))\vert {\mathbb U}_{ab}\vert ^{-1}\vert B_{ab}(I_{\star},G)\vert .$$
Vu comme sous-groupe de $\mathfrak{E}(I_{\star},G;{\mathbb A}_{F})$,  $\mathfrak{E}(I_{\star},G;\mathfrak{o}^V)$ est le produit sur toutes les places $v\in Val(F)$ du sous-groupe $\{0\}$ si $v\in V$, du sous-groupe $\mathfrak{E}(I_{\star},G;\mathfrak{o}_{v})$ de $\mathfrak{E}(I_{\star},G;F_{v})$ si $v\not\in V$. Sa mesure est le produit des mesures de ces sous-groupes. Puisque la mesure d'un point dans $\mathfrak{E}(I_{\star},G;F_{v})$ est $\vert (I_{\star}/\bar{G}_{\star})(F_{v})\vert ^{-1}$, (1) entra\^{\i}ne
$$(6) \qquad mes(\mathfrak{E}(I_{\star},G;\mathfrak{o}^V))=(\prod_{v\in V}\vert (I_{\star}/\bar{G}_{\star})(F_{v})\vert ^{-1}).$$ 
Notons $ker^1(F;G)$ le noyau de l'application $H^1(F;G)\to H^1({\mathbb A}_{F};G)$. On d\'efinit de m\^eme $ker^1(F;I_{\star})$. On a une application naturelle $ker^1(F;I_{\star})\to ker^1(F;G)$. On note $B(I_{\star},G)$ son noyau (c'est-\`a-dire, puisqu'il ne s'agit pas de groupes, l'ensemble des \'el\'ements de $ker^1(F;I_{\star})$ qui deviennent triviaux dans $H^1(F;G)$). Rappelons que ${\mathbb U}$ est un sous-ensemble de $H^1(F;I_{\star})$ et que ${\mathbb U}_{ab}$ est un sous-ensemble de $\mathfrak{E}(I_{\star},G;F)$, lequel est un sous-ensemble de $H^1_{ab}(F;I_{\star})$. Montrons que

(7) l'application naturelle $H^1(F;I_{\star})\to H^1_{ab}(F;I_{\star})$ se restreint en une application surjective ${\mathbb U}\to {\mathbb U}_{ab}$ dont toutes les fibres ont m\^eme nombre d'\'el\'ements;

(8) l'image r\'eciproque de $B_{ab}(I_{\star},G)$ par cette application est $B(I_{\star},G)$.

Par d\'efinition, ${\mathbb U}_{ab}$ est l'ensemble des $u_{ab}\in H^1_{ab}(F;I_{\star})$ tels que

(9) l'image de $u_{ab}$ dans $H^1_{ab}(F;G)$ est nulle;

(10) l'image de $u_{ab}$ dans $H^1_{ab}({\mathbb A}_{F}^V;I_{\star})$ appartient \`a $H^1_{ab}(\mathfrak{o}^V;I_{\star})$;

(11) l'image de $u_{ab}$ dans $H^1_{ab}(F_{V};I_{\star})$ est nulle. 

 On se rappelle l'application
$$H^1(F;I_{\star})\to H^{1}_{ab}(F;I_{\star})\times_{H^1_{ab}({\mathbb A}_{F};I_{\star})}H^1({\mathbb A}_{F};I_{\star})$$
de 7.5(1). Pour $v\not\in V$, l'application $H^1(F_{v};I_{\star})\to H^1_{ab}(F_{v};I_{\star})$ est bijective ([Lab2] proposition 1.6.7). L'application ci-dessus s'identifie donc \`a
$$(12) \qquad H^1(F;I_{\star})\to H^{1}_{ab}(F;I_{\star})\times_{H^1_{ab}(F_{V};I_{\star})}H^1(F_{V};I_{\star}).$$
Notons ${\bf 1}$ l'\'el\'ement trivial de $H^1(F_{V};I_{\star})$. La condition (11) \'equivaut \`a ce que $(u_{ab},{\bf 1})$ appartienne au produit fibr\'e ci-dessus. En particulier,  ${\mathbb U}_{ab}\times \{{\bf 1}\}$ est un sous-ensemble de ce produit fibr\'e. 
Notons $u\mapsto (u_{ab},{\bf u})$ l'application (12). Un \'el\'ement $u\in H^1(F;I_{\star})$ appartient \`a ${\mathbb U}$ si et seulement s'il v\'erifie les conditions (2), (3) et (4) de 7.5. La condition (3) \'equivaut \`a l'\'egalit\'e
  ${\bf u}={\bf 1}$. Les conditions (2) et (4) \'equivalent aux conditions (9) et (10) ci-dessus. Jointes au fait que $(u_{ab},{\bf u})$ appartient au produit fibr\'e de droite de la relation (12), cela \'equivaut comme on vient de le dire  \`a la relation $u_{ab}\in {\mathbb U}_{ab}$. Ainsi ${\mathbb U}$ est exactement l'image r\'eciproque de ${\mathbb U}_{ab}\times {\bf 1}$ par l'application (12). Puisque cette application est surjective ([Lab2] th\'eor\`eme 1.6.10), cela d\'emontre les premi\`eres assertions de (7). Cela d\'emontre aussi que les fibres de l'application ${\mathbb U}\to {\mathbb U}_{ab}$ sont des fibres de l'application (12). Notons $K(I_{\star})$ le noyau de cette application. Soit $u\in H^1(F;I_{\star})$. On voit que la fibre de (12) au-dessus de l'image de $u$ est isomorphe \`a $K(I_{\star,u})$, o\`u $I_{\star, u}$ est la forme int\'erieure de $I_{\star}$ associ\'ee au cocycle $u_{ad}$. Il n'est pas difficile de prouver en g\'en\'eral  que $K(I_{\star,u})$ a m\^eme nombre d'\'el\'ements que $K(I_{\star})$. Dans notre cas, on s'int\'eresse aux \'el\'ements $u\in {\mathbb U}$ et l'assertion est triviale puisque, d'apr\`es 7.3(1) et le lemme 7.5, on a $I_{\star,u}\simeq I_{\star}$ pour tout $u\in {\mathbb U}$. Cela ach\`eve de prouver (7). 
  
  D'apr\`es (3), $B_{ab}(I_{\star},G)$ est l'ensemble des \'el\'ements de ${\mathbb U}_{ab}$ d'image nulle dans $H^1_{ab}({\mathbb A}_{F};I_{\star})$. Donc son image r\'eciproque dans ${\mathbb U}$ est ${\mathbb U}\cap ker^1(F;I_{\star})$. On a vu en 7.5(12) que l'image dans $H^1(F;G)$ d'un \'el\'ement de ${\mathbb U}$ \'etait nulle. Donc ${\mathbb U} \cap ker^1(F;I_{\star})\subset B(I_{\star},G)$. Inversement, un \'el\'ement de $B(I_{\star},G)$ est d'image nulle dans $H^1({\mathbb A}_{F};I_{\star})$, a fortiori d'image nulle dans $H^1_{ab}({\mathbb A}_{F};I_{\star})$. Il v\'erifie donc les conditions (3) et (4) de 7.5. Il est aussi d'image nulle dans $H^1(F;G)$. Il v\'erifie donc aussi la condition (2) de 7.5, d'apr\`es 7.5(12). Donc il appartient \`a ${\mathbb U}$. Il appartient aussi \`a $ker^1(F;I_{\star})$, donc $B(I_{\star},G)\subset {\mathbb U}\cap ker^1(F;I_{\star})$. Finalement, ces deux ensembles sont \'egaux, ce qui ach\`eve la preuve de (8). 

  On d\'eduit de (7) et (8) l'\'egalit\'e
  $$(13) \qquad \vert {\mathbb U}_{ab}\vert ^{-1}\vert B_{ab}(I_{\star},G)\vert =\vert {\mathbb U}\vert ^{-1}\vert B(I_{\star},G)\vert.$$
  D'apr\`es la proposition 3.7 de [Lab1], on a
  $$(14)\qquad mes(\mathfrak{E}(I_{\star},G;{\mathbb A}_{F}/F)=\tau(I_{\star})\vert B(I_{\star},G)\vert d(I_{\star},G)\tau(G)^{-1}.$$
  En mettant bout-\`a-bout les \'egalit\'es (4), (5), (6), (13) et (14), on obtient le lemme. $\square$

\bigskip

\subsection{Comparaison de deux mesures de Tamagawa}
Pour tout groupe alg\'ebrique $H$ d\'efini sur $F$, on note $X^*(H)$  le groupe de caract\`eres alg\'ebriques de $ H$. On a un homomorphisme naturel de restriction $X^*(I_{\star})\to X^*(\bar{G}_{\star})$ dont les noyau et conoyau sont finis. On a donc aussi un homomorphisme
$$(1) \qquad X^*(I_{\star})^{\Gamma_{F}}\to X^*(\bar{G}_{\star})^{\Gamma_{F}}$$
qui a les m\^emes propri\'et\'es. Notons $x(I_{\star})$ le quotient du nombre d'\'el\'ements de son noyau par le nombre d'\'el\'ements de son conoyau.

\ass{Lemme}{On a l'\'egalit\'e
$$\tau(I_{\star})=\tau(\bar{G}_{\star})x(I_{\star}) [I_{\star}(F):\bar{G}_{\star}(F)] ^{-1}[ I_{\star}(F_{V}):\bar{G}_{\star}(F_{V})]\vert (I_{\star}/\bar{G}_{\star})(F_{V})\vert ^{-1}.$$}

Preuve. 
On a rappel\'e en 4.1 la d\'efinition de la mesure de Tamagawa sur $\bar{G}_{\star}({\mathbb A}_{F})$. Elle est de la forme $dx=\ell_{\bar{G}_{\star}}^{-1}\otimes_{v\in Val(F)}dx_{v}$, o\`u $\ell_{\bar{G}_{\star}}$ est le terme principal du d\'eveloppement en $s=1$ d'une  certaine fonction $L$ et o\`u, pour toute place $v$, $dx_{v}$ est une certaine mesure sur $\bar{G}_{\star}(F_{v})$. Une d\'efinition analogue  vaut pour le groupe quasi-connexe $I_{\star}$. La mesure de Tamagawa sur $I_{\star}({\mathbb A}_{F})$ est de la forme $di=\ell_{I_{\star}}^{-1}\otimes_{v\in Val(F)}di_{v}$, avec des notations \'evidentes. En inspectant la d\'efinition de [Lab1] 1.2, on s'aper\c{c}oit que

- on a $\ell_{I_{\star}}=\ell_{\bar{G}_{\star}}$;

- pour tout $v\in Val(F)$, la restriction de $di_{v}$ \`a l'ouvert $\bar{G}_{\star}(F_{v})$ de $I_{\star}(F_{v})$ est \'egale \`a $\vert (I_{\star}/\bar{G}_{\star})(F_{v})\vert ^{-1} dx_{v}$.

Pour d\'efinir les nombres de Tamagawa, on a aussi besoin de mesures sur le groupe $\mathfrak{A}_{\bar{G}_{\star}}$. Ce groupe s'identifie naturellement  \`a $Hom(X^*(\bar{G}_{\star})^{\Gamma_{F}},{\mathbb R})$, ou encore \`a $Hom(X^*(I_{\star})^{\Gamma_{F}},{\mathbb R})$. On d\'efinit les r\'eseaux $\mathfrak{A}_{\bar{G}_{\star},{\mathbb Z}}=Hom(X^*(\bar{G}_{\star})^{\Gamma_{F}},{\mathbb Z})$ et $\mathfrak{A}_{I_{\star},{\mathbb Z}}=Hom(X^*(I_{\star})^{\Gamma_{F}},{\mathbb Z})$. On  note $da_{\bar{G}_{\star}}$ la mesure de Haar sur $\mathfrak{A}_{\bar{G}_{\star}}$ telle que $\mathfrak{A}_{\bar{G}_{\star},{\mathbb Z}}$ soit de covolume $1$. Selon [Lab1] 1.2, on note $da_{I_{\star}}$ la mesure de Haar sur $\mathfrak{A}_{\bar{G}_{\star}}$ telle que $\mathfrak{A}_{I_{\star},{\mathbb Z}}$ soit de covolume $\vert X^*(I_{\star}/\bar{G}_{\star})^{\Gamma_{F}}\vert ^{-1}$. Alors $\tau(\bar{G}_{\star})$ est la mesure de $\bar{G}_{\star}(F)\backslash \bar{G}_{\star}({\mathbb A}_{F})/\mathfrak{A}_{\bar{G}_{\star}}$, $\bar{G}_{\star}({\mathbb A}_{F})$ et $\mathfrak{A}_{\bar{G}_{\star}}$ \'etant munis des mesures $dx$ et $da_{\bar{G}_{\star}}$. Et $\tau(I_{\star})$ est la mesure de $I_{\star}(F)\backslash I_{\star}({\mathbb A}_{F})/ \mathfrak{A}_{\bar{G}_{\star}}$, $I_{\star}({\mathbb A}_{F})$ et $\mathfrak{A}_{\bar{G}_{\star}}$ \'etant munis des mesures $di$ et $da_{I_{\star}}$. 
Le r\'eseau $\mathfrak{A}_{\bar{G}_{\star},{\mathbb Z}}$ est contenu dans $\mathfrak{A}_{I_{\star},{\mathbb Z}}$ et l'indice de ce sous-groupe est \'egal au nombre d'\'el\'ements du conoyau de l'homomorphisme (1). Le groupe $ X^*(I_{\star}/\bar{G}_{\star})^{\Gamma_{F}}$ est le noyau  cet homomorphisme. Il en r\'esulte que $da_{I_{\star}}=x(I_{\star})^{-1}da_{\bar{G}_{\star}}$. On peut donc utiliser  l'unique mesure $da_{\bar{G}_{\star}}$ sur $\mathfrak{A}_{\bar{G}_{\star}}$ et d\'efinir $\tau(I_{\star})$ comme le produit de $x(I_{\star})$ et de la mesure   de $I_{\star}(F)\backslash I_{\star}({\mathbb A}_{F})/ \mathfrak{A}_{\bar{G}_{\star}}$, $I_{\star}({\mathbb A}_{F})$ et $\mathfrak{A}_{\bar{G}_{\star}}$ \'etant munis des mesures $di$ et $da_{\bar{G}_{\star}}$.

 Pour $v\not\in V$, on a d\'efini le groupe hypersp\'ecial $K_{\star,v}=ad_{h_{\star,v}}(K_{v})\cap \bar{G}_{\star}(F_{v})$. Posons $K_{I_{\star},v}=ad_{h_{\star,v}}(K_{v})\cap I_{\star}(F_{v})$. 
   Montrons que

(2) pour toute place $v\not\in V$, l'application $K_{I_{\star},v}\to (I_{\star}/\bar{G}_{\star})(F_{v})$ est surjective; les groupes  $ (I_{\star}/\bar{G}_{\star})(F_{v})$, $K_{I_{\star},v}/K_{\star,v}$ et $I_{\star}(F_{v})/\bar{G}_{\star}(F_{v})$ sont isomorphes. 

On peut fixer une paire de Borel d\'efinie sur $F_{v}$  de $\bar{G}_{\star}$ dont est issu  le sous-groupe hypersp\'ecial $K_{\star,v}$. On note son tore $\bar{T}_{0}$. Notons $T_{0}$ son commutant dans $G$. On a alors $I_{\star}/\bar{G}_{\star}\simeq T_{0}^{\theta}/T_{0}^{\theta,0}$. On a d\'ej\`a dit que, pour un tel tore, on a $T_{0}^{\theta}(\bar{F}_{v})/T_{0}^{\theta,0}(\bar{F}_{v})=T_{0}^{\theta}(\mathfrak{o}_{v}^{nr})/T_{0}^{\theta,0}(\mathfrak{o}_{v}^{nr})$. Parce que $H^1(\mathfrak{o}_{v};T_{0}^{\theta,0})=0$, l'homomorphisme 
$$T_{0}^{\theta}(\mathfrak{o}_{v})\to H^0(\mathfrak{o}_{v};T_{0}^{\theta}/T_{0}^{\theta,0})=(T_{0}^{\theta}/T_{0}^{\theta,0})(F_{v})=(I_{\star}/\bar{G}_{\star})(F_{v})$$
est surjectif. Par ailleurs, le lemme 5.5 de [W1] entra\^{\i}ne que $T_{0}^{\theta}(\mathfrak{o}_{v}^{nr})$ est engendr\'e par $T_{0}^{\theta,0}(\mathfrak{o}_{v}^{nr})$ et par  les \'el\'ements de $Z(G)^{\theta}$ d'ordre fini premier \`a la caract\'eristique r\'esiduelle $p$. Ces deux groupes \'etant contenus dans $K_{v}^{nr}$, on a $T_{0}^{\theta}(\mathfrak{o}_{v}^{nr})\subset K_{v}^{nr}$, donc 
$T_{0}^{\theta}(\mathfrak{o}_{v})\subset K_{I_{\star},v}$. La premi\`ere assertion de (2) en r\'esulte. La seconde en est une cons\'equence imm\'ediate. 

Fixons un ensemble fini $V''$ de places de $F$, contenant $V$ et tel que l'on ait l'\'egalit\'e
$$\bar{G}_{\star}({\mathbb A}_{F})=\bar{G}_{\star}(F)(\bar{G}_{\star}(F_{V''})\times K_{\star}^{V''}),$$
o\`u $K_{\star}^{V''}=\prod_{v\not\in V''}K_{\star,v}$. Il r\'esulte de (2) que l'on a 
$$I_{\star}({\mathbb A}_{F})=I_{\star}(F_{V''})\times \bar{G}_{\star}({\mathbb A}_{F}^{V''})K_{I_{\star}}^{V''}=\bar{G}_{\star}({\mathbb A}_{F})(I_{\star}(F_{V''})\times K_{I_{\star}}^{V''}).$$
L'\'egalit\'e pr\'ec\'edente entra\^{\i}ne
$$(3) \qquad I_{\star}({\mathbb A}_{F})=\bar{G}_{\star}(F)(I_{\star}(F_{V''})\times K_{I_{\star}}^{V''}),$$
a fortiori
$$I_{\star}(\mathfrak{A}_{F})=I_{\star}(F)(I_{\star}(F_{V''})\times K_{I_{\star}}^{V''}).$$
Notons $\Xi_{\bar{G}_{\star}}$ la projection dans $\bar{G}_{\star}(F_{V''})$ de $\bar{G}_{\star}(F)\cap (\bar{G}_{\star}(F_{V'''})\times K_{\star}^{V''})$. D\'efinissons de m\^eme $\Xi_{I_{\star}}$. Alors on a les \'egalit\'es
$$(4) \qquad \tau(\bar{G}_{\star})=mes(K_{\star}^{V''})mes(\Xi_{\bar{G}_{\star}}\backslash \bar{G}_{\star}(F_{V''})/\mathfrak{A}_{\bar{G}_{\star}}),$$
$$(5) \qquad \tau(I_{\star})=x(I_{\star})mes(K_{I_{\star}}^{V''})mes(\Xi_{I_{\star}}\backslash I_{\star}(F_{V''})/\mathfrak{A}_{\bar{G}_{\star}}).$$
Pour tout $v\not\in V''$, la mesure $di_{v}$ se restreint \`a $\bar{G}_{\star}(F_{v})$ en la mesure $dx_{v}$ multipli\'ee par $\vert (I_{\star}/\bar{G}_{\star})(F_{v})\vert ^{-1}$. La relation (2) entra\^{\i}ne que $mes(K_{I_{\star},v})=mes(K_{\star,v})$. Donc 

(6) $mes(K_{\star}^{V''})=mes(K_{I_{\star}}^{V''})$.

 Fixons un ensemble de repr\'esentants ${\cal U}$ du quotient $\bar{G}_{\star}(F_{V''})\backslash I_{\star}(F_{V''})$. Alors 
$\Xi_{\bar{G}_{\star}}\backslash I_{\star}(F_{V''})/\mathfrak{A}_{\bar{G}_{\star}}$ est r\'eunion disjointe des sous-ensembles $\Xi_{\bar{G}_{\star}}\backslash \bar{G}_{\star}(F_{V''})u/\mathfrak{A}_{\bar{G}_{\star}}$ pour $u\in {\cal U}$. La comparaison des mesures locales montre que chacun de ces sous-ensembles a pour mesure 
$$\vert (I_{\star}/\bar{G}_{\star})(F_{V''})\vert ^{-1} mes(\Xi_{\bar{G}_{\star}}\backslash \bar{G}_{\star}(F_{V''})/\mathfrak{A}_{\bar{G}_{\star}}).$$
Puisque $\vert {\cal U}\vert =[ I_{\star}(F_{V''}):\bar{G}_{\star}(F_{V''})] $, on obtient
$$mes (\Xi_{\bar{G}_{\star}}\backslash I_{\star}(F_{V''})/\mathfrak{A}_{\bar{G}_{\star}})=\vert (I_{\star}/\bar{G}_{\star})(F_{V''})\vert ^{-1}[ I_{\star}(F_{V''}):\bar{G}_{\star}(F_{V''})] mes(\Xi_{\bar{G}_{\star}}\backslash \bar{G}_{\star}(F_{V''})/\mathfrak{A}_{\bar{G}_{\star}}),$$
ou encore
$$mes(\Xi_{I_{\star}}\backslash I_{\star}(F_{V''})/\mathfrak{A}_{\bar{G}_{\star}})=[\Xi_{I_{\star}}:\Xi_{\bar{G}_{\star}}]^{-1}\vert (I_{\star}/\bar{G}_{\star})(F_{V''})\vert ^{-1}$$
$$[I_{\star}(F_{V''}):\bar{G}_{\star}(F_{V''})] mes(\Xi_{\bar{G}_{\star}}\backslash \bar{G}_{\star}(F_{V'''})/\mathfrak{A}_{\bar{G}_{\star}}).$$
La relation (2) entra\^{\i}ne que
$$\vert (I_{\star}/\bar{G}_{\star})(F_{V''})\vert ^{-1}[ I_{\star}(F_{V''}):\bar{G}_{\star}(F_{V''})]=\vert (I_{\star}/\bar{G}_{\star})(F_{V})\vert ^{-1}[ I_{\star}(F_{V}):\bar{G}_{\star}(F_{V})].$$
D'autre part, (3) entra\^{\i}ne que $I_{\star}(F)=\bar{G}_{\star}(F)\Xi_{I_{\star}}$. Donc l'application naturelle $\Xi_{I_{\star}}/\Xi_{\bar{G}_{\star}}\to I_{\star}(F)/\bar{G}_{\star}(F)$ est bijective. La relation ci-dessus se r\'ecrit
$$(7) \qquad mes(\Xi_{I_{\star}}\backslash I_{\star}(F_{V''})/\mathfrak{A}_{\bar{G}_{\star}})=[ I_{\star}(F):\bar{G}_{\star}(F)] ^{-1}\vert (I_{\star}/\bar{G}_{\star})(F_{V})\vert ^{-1}[ I_{\star}(F_{V}):\bar{G}_{\star}(F_{V})] $$
$$mes(\Xi_{\bar{G}_{\star}}\backslash \bar{G}_{\star}(F_{V'''})/\mathfrak{A}_{\bar{G}_{\star}}).$$
Le lemme r\'esulte de (4), (5), (6) et (7). $\square$

\bigskip

\subsection{Calcul de $d(I_{\star},G)$}
En appliquant la d\'efinition de 7.10 \`a l'\'el\'ement $\eta_{\star}$, on d\'efinit le covolume  $covol(\mathfrak{A}_{\bar{G}_{\star},{\mathbb Z}})$.
 D'autre part, on a d\'efini en 5.8 une constante not\'ee $C(\tilde{G})$. 

\ass{Lemme}{On a l'\'egalit\'e
$$d(I_{\star},G)=x(I_{\star})^{-1}\tau(G)C(\tilde{G})^{-1} covol(\mathfrak{A}_{\bar{G}_{\star},{\mathbb Z}})^{-1}.$$}

Preuve. Par d\'efinition, $D(I_{\star},G)$ est le   conoyau de l'homomorphisme
$$H^1_{ab}({\mathbb A}_{F}/F;I_{\star})\to H^1_{ab}({\mathbb A}_{F}/F;G),$$
c'est-\`a-dire de l'homomorphisme
$$H^{2,1,0}({\mathbb A}_{F}/F;\bar{S}_{sc}\to S\stackrel{1-\theta}{\to}(1-\theta)(S))\to H^{2,1}({\mathbb A}_{F}/F;S_{sc}\to S).$$
On a une suite exacte
$$H^{1,0}({\mathbb A}_{F}/F;S\stackrel{1-\theta}{\to}(1-\theta)(S))\to H^{2,1,0}({\mathbb A}_{F}/F;\bar{S}_{sc}\to S\stackrel{1-\theta}{\to}(1-\theta)(S))\to H^2({\mathbb A}_{F}/F;\bar{S}_{sc}).$$
Comme on l'a dit en 7.11, le dernier groupe est nul car $\bar{S}_{sc}$ est un tore elliptique. Donc $D(I_{\star},G)$ est aussi le conoyau de l'homomorphisme compos\'e
$$(1) \qquad H^{1,0}({\mathbb A}_{F}/F;S\stackrel{1-\theta}{\to}(1-\theta)(S))\to H^{2,1}({\mathbb A}_{F}/F;S_{sc}\to S).$$
Dans la cat\'egorie des complexes de tores, le triangle
$$\begin{array}{ccc}S_{sc}&\to&S\\ \parallel&&\qquad\downarrow 1-\theta\\ S_{sc}&\stackrel{1-\theta}{\to}&(1-\theta)(S)\\ \,\, \downarrow \pi&&\parallel\\ S&\stackrel{1-\theta}{\to}&(1-\theta)(S)\\ \end{array}$$
est distingu\'e. Il donne naissance \`a une suite exacte
$$H^{1,0}({\mathbb A}_{F}/F;S\stackrel{1-\theta}{\to}(1-\theta)(S))\to H^{2,1}({\mathbb A}_{F}/F;S_{sc}\to S)$$
$$\to H^{2,1}({\mathbb A}_{F}/F;S_{sc}\stackrel{1-\theta}{\to}(1-\theta)(S))\to H^{2,1}({\mathbb A}_{F}/F;S_{sc}\stackrel{1-\theta}{\to}(1-\theta)(S)).$$
En inspectant les d\'efinitions, on voit que le premier homomorphisme ci-dessus est \'egal \`a celui de (1). Donc $D(I_{\star},G)$ est aussi le noyau du dernier homomorphisme ci-dessus. Le lemme C.2.A de [KS] 
entra\^{\i}ne par dualit\'e que ce noyau a m\^eme nombre d'\'el\'ements que le conoyau de l'homomorphisme dual
$$(2) \qquad H^{1,0}(\Gamma_{F};X^*((1-\theta)(S))\to X^*(S))\to H^{1,0}(\Gamma_{F};X^*((1-\theta)(S))\to X^*(S_{sc})).$$
La fl\`eche $X^*((1-\theta)(S))\to X^*(S)$ envoie un caract\`ere $x^*$ de $(1-\theta)(S)$ sur le caract\`ere $x^*\circ (1-\theta)$ de $S$. La fl\`eche $X^*((1-\theta)(S))\to X^*(S_{sc})$ envoie un caract\`ere $x^*$ de $(1-\theta)(S)$ sur le caract\`ere $x^*\circ (1-\theta)\circ \pi$ de $S_{sc}$. Montrons que

(3) le groupe $H^{1,0}(\Gamma_{F};X^*((1-\theta)(S))\to X^*(S_{sc}))$ est isomorphe \`a $\pi_{0}((Z(\hat{G})/Z(\hat{G})\cap \hat{T}^{\hat{\theta},0})^{\Gamma_{F}})$.

On a    $X^*((1-\theta)(S))=X_{*}(\hat{S}/\hat{S}^{\hat{\theta},0})$ et $X^*(S_{sc})=X_{*}(\hat{S}_{ad})$.  Pour simplifier, notons ces groupes $X_{1}$ et $X_{2}$. Notons $\varphi:X_{1}\to X_{2}$ l'homomorphisme d\'efini ci-dessus ainsi que les homomorphismes qui s'en d\'eduisent fonctoriellement. On a la suite exacte de complexes de $\Gamma_{F}$-modules 
$$\begin{array}{ccc}1&&1\\ \downarrow&&\downarrow\\ X_{1}&\stackrel{\varphi}{\to}& X_{2}\\  \downarrow&&\downarrow\\ X_{1}\otimes_{{\mathbb Z}}{\mathbb Q}&\stackrel{\varphi}{\to}&X_{2}\otimes_{{\mathbb Z}}{\mathbb Q}\\ 
\downarrow&&\downarrow\\ X_{1}\otimes_{{\mathbb Z}}{\mathbb Q}/{\mathbb Z}&\stackrel{\varphi}{\to} &X_{2}\otimes_{{\mathbb Z}}{\mathbb Q}/{\mathbb Z}\\ \downarrow&&\downarrow\\  1&&1\\ \end{array}$$
D'o\`u une suite exacte
$$(4) \qquad H^{0,.}(\Gamma_{F};X_{1}\otimes_{{\mathbb Z}}{\mathbb Q}\to X_{2}\otimes_{{\mathbb Z}}{\mathbb Q})\to H^{0,.}(\Gamma_{F};X_{1}\otimes_{{\mathbb Z}}{\mathbb Q}/{\mathbb Z}\to X_{2}\otimes_{{\mathbb Z}}{\mathbb Q}/{\mathbb Z})$$
$$\to H^{1,0}(\Gamma_{F};X_{1}\to X_{2})\to H^{1,0}(\Gamma_{F};X_{1}\otimes_{{\mathbb Z}}{\mathbb Q}\to X_{2}\otimes_{{\mathbb Z}}{\mathbb Q}).$$
On a not\'e $H^{0,.}$ le groupe not\'e simplement $H^0$ par Kottwitz et Shelstad. On a aussi une suite exacte
$$H^0(\Gamma_{F};X_{1}\otimes_{{\mathbb Z}}{\mathbb Q})\to H^0(\Gamma_{F};X_{2}\otimes_{{\mathbb Z}}{\mathbb Q})\to H^{1,0}(\Gamma_{F};X_{1}\otimes_{{\mathbb Z}}{\mathbb Q}\to X_{2}\otimes_{{\mathbb Z}}{\mathbb Q})\to H^1(\Gamma_{F};X_{1}\otimes_{{\mathbb Z}}{\mathbb Q}).$$
Le dernier groupe est limite inductive de groupes $H^1(Gal(E/F);X_{1}\otimes_{{\mathbb Z}}{\mathbb Q})$ sur les extensions galoisiennes finies $E$ de $F$ telles que $\Gamma_{E}$ agisse trivialement sur $X_{1}$. Or ces groupes sont nuls puisque $Gal(E/F)$ est fini et $X_{1}\otimes_{{\mathbb Z}}{\mathbb Q}$ est divisible. La suite ci-dessus se r\'ecrit
$$(5) \qquad X_{1}^{\Gamma_{F}}\otimes_{{\mathbb Z}}{\mathbb Q}\to X_{2}^{\Gamma_{F}}\otimes_{{\mathbb Z}}{\mathbb Q}\to H^{1,0}(\Gamma_{F};X_{1}\otimes_{{\mathbb Z}}{\mathbb Q}\to X_{2}\otimes_{{\mathbb Z}}{\mathbb Q})\to 0.$$
  Il r\'esulte des d\'efinitions que
$$ X_{2}\otimes_{{\mathbb Z}}{\mathbb Q}=\varphi(X_{1}\otimes_{{\mathbb Z}}{\mathbb Q})\oplus X_{2}^{\hat{\theta}}\otimes_{{\mathbb Z}}{\mathbb Q}.$$
D'o\`u
$$(6) \qquad X_{2}^{\Gamma_{F}}\otimes_{{\mathbb Z}}{\mathbb Q}=\varphi(X_{1}^{\Gamma_{F}}\otimes_{{\mathbb Z}}{\mathbb Q})\oplus X_{2}^{\Gamma_{F},\hat{\theta}}\otimes_{{\mathbb Z}}{\mathbb Q}.$$
Le tore $\bar{S}_{sc} $ est elliptique  donc $X_{*}(\bar{S}_{sc})^{\Gamma_{F}}=0$. Il en r\'esulte que $X_{*}(\bar{S})^{\Gamma_{F}}=X_{*}(Z(\bar{G}_{\star})^0)^{\Gamma_{F}}$. On a d\'ej\`a remarqu\'e que $\eta_{\star}$ \'etait elliptique dans $\tilde{G}(F)$. Le groupe ci-dessus est donc $X_{*}(Z(G)^{\theta,0})^{\Gamma_{F}}$. D'autre part $X_{*}(\bar{S})=X_{*}(S)^{\theta}$. On obtient  l'\'egalit\'e $X_{*}(S)^{\Gamma_{F},\theta}=X_{*}(Z(G)^{\theta,0})^{\Gamma_{F}}$. Cela entra\^{\i}ne $X_{*}(S_{sc})^{\Gamma_{F},\theta}=0$. Cela \'equivaut \`a  $X_{2}^{\Gamma_{F},\hat{\theta}}=0$. Il r\'esulte alors de (6) que la premi\`ere application de (5) est surjective. En cons\'equence, $H^{1,0}(\Gamma_{F};X_{1}\otimes_{{\mathbb Z}}{\mathbb Q}\to X_{2}\otimes_{{\mathbb Z}}{\mathbb Q})=0$. 

Il r\'esulte des d\'efinitions que les deux premiers groupes de la suite (4) sont respectivement les noyaux des homomorphismes
$$X_{1}^{\Gamma_{F}}\otimes_{{\mathbb Z}}{\mathbb Q}\to X_{2}^{\Gamma_{F}}\otimes_{{\mathbb Z}}{\mathbb Q}$$
et 
$$(X_{1}\otimes_{{\mathbb Z}}{\mathbb Q}/{\mathbb Z})^{\Gamma_{F}}\to (X_{2}\otimes_{{\mathbb Z}}{\mathbb Q}/{\mathbb Z})^{\Gamma_{F}}.$$
On peut identifier ${\mathbb Q}/{\mathbb Z}$ au groupe des racines de l'unit\'e dans ${\mathbb C}^{\times}$. Alors l'application $x_{*}\otimes \zeta\mapsto x_{*}(\zeta)$ identifie  $X_{1}\otimes_{{\mathbb Z}}{\mathbb Q}/{\mathbb Z}$ au sous-groupe de torsion $(\hat{S}/\hat{S}^{\hat{\theta},0})_{tors}\subset \hat{S}/\hat{S}^{\hat{\theta},0}$. Le sous-groupe $(X_{1}\otimes_{{\mathbb Z}}{\mathbb Q}/{\mathbb Z})^{\Gamma_{F}}$ s'identifie \`a $(\hat{S}/\hat{S}^{\hat{\theta},0})_{tors}^{\Gamma_{F}}$. Notons $K$ le noyau de l'homomorphisme
$$\hat{S}/\hat{S}^{\hat{\theta},0}\to \hat{S}_{ad}.$$
Alors le deuxi\`eme groupe de la suite (4) s'identifie \`a $K_{tors}^{\Gamma_{F}}$. On voit que l'image du premier groupe est $(K^{\Gamma_{F},0})_{tors}$. A ce point, on d\'eduit de (4) une suite exacte
$$(K^{\Gamma_{F},0})_{tors}\to K_{tors}^{\Gamma_{F}}\to H^{1,0}(\Gamma_{F};X_{1}\to X_{2})\to 0$$
On v\'erifie facilement que l'homomorphisme naturel
$$(K^{\Gamma_{F},0})_{tors}\backslash K_{tors}^{\Gamma_{F}}\to \pi_{0}(K^{\Gamma_{F}})$$
est bijectif. D'autre part, le noyau de l'application $t\mapsto (1-\theta)(t_{ad})$ de $\hat{S}$ dans $\hat{S}_{ad}$ est $Z(\hat{G})\hat{S}^{\hat{\theta},0}$. Il en r\'esulte que $K\simeq Z(\hat{G})/Z(\hat{G})\cap \hat{S}^{\hat{\theta},0}$. Le tore $\hat{S}$ est \'egal \`a $\hat{T}$ muni de l'action galoisienne $\sigma\mapsto \omega_{S}(\sigma)\circ \sigma_{G^*}$. Sur $Z(\hat{G})$, cette action co\"{\i}ncide avec $\sigma\mapsto \sigma_{G^*}$. Donc $Z(\hat{G})/Z(\hat{G})\cap \hat{S}^{\hat{\theta},0}=Z(\hat{G})/Z(\hat{G})\cap \hat{T}^{\hat{\theta},0}$ (en tant que groupes munis d'actions galoisiennes). On en d\'eduit 
$$(K^{\Gamma_{F},0})_{tors}\backslash K_{tors}^{\Gamma_{F}}\simeq \pi_{0}(K^{\Gamma_{F}})\simeq \pi_{0}((Z(\hat{G})/Z(\hat{G})\cap \hat{T}^{\hat{\theta},0})^{\Gamma_{F}})$$
puis (3).

L'application (2) s'ins\`ere dans un diagramme de suites exactes de cohomologie
$$\begin{array}{ccc}H^0(\Gamma_{F};X^*((1-\theta)(S)))&=&H^0(\Gamma_{F};X^*((1-\theta)(S)))\\ \downarrow&&\downarrow\\ H^0(\Gamma_{F};X^*(S))&\to&H^0(\Gamma_{F};X^*(S_{sc}))\\ \downarrow&&\downarrow\\ H^{1,0}(\Gamma_{F};X^*((1-\theta)(S))\to X^*(S))&\to& H^{1,0}(\Gamma_{F};X^*((1-\theta)(S))\to X^*(S_{sc}))\\ \downarrow&&\downarrow\\ H^1(\Gamma_{F};X^*((1-\theta)(S)))&=&H^1(\Gamma_{F};X^*((1-\theta)(S)))\\ \end{array}$$
Il en r\'esulte formellement que le noyau de (2) est l'image dans $H^{1,0}(\Gamma_{F};X^*((1-\theta)(S))\to X^*(S))$ du noyau de
$$H^0(\Gamma_{F};X^*(S))\to H^0(\Gamma_{F};X^*(S_{sc})).$$
Mais le noyau de $X^*(S)\to X^*(S_{sc})$ n'est autre que $X^*(G)$. Donc le noyau de (2) est l'image naturelle dans $H^{1,0}(\Gamma_{F};X^*((1-\theta)(S))\to X^*(S))$ de $X^*(G)^{\Gamma_{F}}$. D'autre part, on a la suite exacte
$$1\to S^{\theta}\to S\stackrel{1-\theta}{\to}(1-\theta)(S)\to 1$$
Le groupe $S^{\theta}$ n'est pas connexe mais est diagonalisable. Il en r\'esulte une suite exacte
$$0\to X^*((1-\theta)(S))\to X^*(S)\to X^*(S^{\theta})\to 0$$
Il en r\'esulte formellement l'\'egalit\'e
$$H^{1,0}(\Gamma_{F};X^*((1-\theta)(S))\to X^*(S))=X^*(S^{\theta})^{\Gamma_{F}}.$$
L'homomorphisme de $X^*(G)^{\Gamma_{F}}$ dans ce groupe n'est autre que l'application de restriction. En notant $Im(X^*(G)^{\Gamma_{F}})$ son image, on obtient que l'image de l'application (2) est isomorphe \`a $X^*(S^{\theta})^{\Gamma_{F}}/Im(X^*(G)^{\Gamma_{F}})$. On a dit que $d(I_{\star},G)$ \'etait \'egal au nombre d'\'el\'ements du conoyau de l'homomorphisme (2). Gr\^ace \`a (3) et au r\'esultat pr\'ec\'edent, on obtient
$$ d(I_{\star},G)=\vert \pi_{0}((Z(\hat{G})/Z(\hat{G})\cap \hat{T}^{\hat{\theta},0})^{\Gamma_{F}})\vert [ X^*(S^{\theta})^{\Gamma_{F}}:Im(X^*(G)^{\Gamma_{F}})] ^{-1}.$$
Notons $Im(X^*(G)^{\Gamma_{F},\theta})$ l'image dans $X^*(S^{\theta})^{\Gamma_{F}}$ du sous-groupe $X^*(G)^{\Gamma_{F},\theta}\subset X^*(G)^{\Gamma_{F}}$. Kottwitz et Shelstad ont calcul\'e
$$[ Im(X^*(G)^{\Gamma_{F}}):Im(X^*(G)^{\Gamma_{F},\theta})]=\vert det((1-\theta)_{\vert \mathfrak{A}_{G}/\mathfrak{A}_{\tilde{G}}})\vert \vert \pi_{0}(\hat{T}^{\hat{\theta},0}\cap Z(\hat{G})^{\Gamma_{F},0})\vert ^{-1},$$
cf. [KS] calcul de l'expression 6.4.14. 
L'expression ci-dessus se transforme en
$$ d(I_{\star},G)=\vert \pi_{0}((Z(\hat{G})/Z(\hat{G})\cap \hat{T}^{\hat{\theta},0})^{\Gamma_{F}})\vert \vert det((1-\theta)_{\vert \mathfrak{A}_{G}/\mathfrak{A}_{\tilde{G}}})\vert \vert \pi_{0}(\hat{T}^{\hat{\theta},0}\cap Z(\hat{G})^{\Gamma_{F},0})\vert ^{-1}$$
$$[ X^*(S^{\theta})^{\Gamma_{F}}:Im(X^*(G)^{\Gamma_{F},\theta})] ^{-1}.$$
En se rappelant la d\'efinition de $C(\tilde{G})$ donn\'ee en 5.8, on peut r\'ecrire
$$d(I_{\star},G)=C(\tilde{G})^{-1}\tau(G)covol(\mathfrak{A}_{\tilde{G},{\mathbb Z}})^{-1}[X^*(S^{\theta})^{\Gamma_{F}}:Im(X^*(G)^{\Gamma_{F},\theta})]^{-1}.$$
Pour obtenir le lemme, il reste \`a prouver l'\'egalit\'e
$$(7) \qquad x(I_{\star})covol(\mathfrak{A}_{\bar{G}_{\star},{\mathbb Z}})=[ X^*(S^{\theta})^{\Gamma_{F}}:Im(X^*(G)^{\Gamma_{F},\theta})] covol(\mathfrak{A}_{\tilde{G},{\mathbb Z}}) .$$
Le groupe $I_{\star}$ est isomorphe au quotient de $Z(I_{\star})\times \bar{G}_{\star,SC}$ par le sous-groupe form\'e des $(\bar{\pi}_{\star}(z)^{-1},z)$ pour $z\in Z(\bar{G}_{\star,SC})$. Il en r\'esulte que $X^*(I_{\star})$ est le groupe des caract\`eres du groupe $Z(I_{\star})/\bar{\pi}_{\star}(Z(\bar{G}_{\star,SC}))$. L'homomorphisme naturel de ce groupe dans $S^{\theta}/\bar{\pi}_{\star}(\bar{S}_{sc})$ est un isomorphisme. Donc $X^*(I_{\star})=X^*(S^{\theta}/\bar{\pi}_{\star}(\bar{S}_{sc}))$. On obtient une suite exacte
$$0\to X^*(I_{\star})\to X^*(S^{\theta})\to X^*(\bar{S}_{sc}).$$
D'o\`u aussi une suite exacte
$$0\to X^*(I_{\star})^{\Gamma_{F}}\to X^*(S^{\theta})^{\Gamma_{F}}\to X^*(\bar{S}_{sc})^{\Gamma_{F}}.$$
Mais $\bar{S}_{sc}$ est un tore elliptique donc le dernier groupe est nul. Finalement 
$$(8) \qquad X^*(I_{\star})^{\Gamma_{F}}= X^*(S^{\theta})^{\Gamma_{F}}.$$
On a un diagramme commutatif d'homomorphismes naturels
$$(9) \qquad \begin{array}{ccccccc}&&&&&&X^*(G)^{\Gamma_{F},\theta}\\ &&&&&\swarrow&\quad\downarrow \phi\\ 0&\to&X^*(I_{\star}/\bar{G}_{\star})^{\Gamma_{F}}&\to&X^*(I_{\star})^{\Gamma_{F}}&\to&X^*(\bar{G}_{\star})^{\Gamma_{F}}\\ \end{array}$$
La suite horizontale est exacte. La fl\`eche verticale $\phi$ s'inscrit dans un diagramme commutatif
$$\begin{array}{ccc}X^*(G)^{\Gamma_{F},\theta}&\to&X^*(Z(G)^{0})^{\Gamma_{F},\theta}\\ \quad\downarrow \phi&&\downarrow\\ X^*(\bar{G}_{\star})^{\Gamma_{F}}&\to&X^*(Z(\bar{G}_{\star})^0)^{\Gamma_{F}}\\ \end{array}$$
Les fl\`eches horizontales sont injectives. La fl\`eche de droite est injective: cela r\'esulte de l'ellipticit\'e de $\eta_{\star}$. Donc la fl\`eche de gauche est aussi injective. En revenant au diagramme (9), on en d\'eduit un diagramme commutatif \`a fl\`eches injectives
$$\begin{array}{ccccc}&&X^*(G)^{\Gamma_{F},\theta}&&\\ &\phi_{1}\swarrow&&\searrow\phi&\\ X^*(I_{\star})^{\Gamma_{F}}/X^*(I_{\star}/\bar{G}_{\star})^{\Gamma_{F}}&&\stackrel{\phi_{2}}{\to}&&X^*(\bar{G}_{\star})^{\Gamma_{F}}\\ \end{array}$$
Par d\'efinition, on a
$$x(I_{\star})=\vert X^*(I_{\star}/\bar{G}_{\star})^{\Gamma_{F}}\vert  \vert coker(\phi_{2})\vert ^{-1}.$$
L'injectivit\'e de $\phi_{1}$ et la relation (8) entra\^{\i}nent que
$$[X^*(S^{\theta})^{\Gamma_{F}}:Im(X^*(G)^{\Gamma_{F},\theta})] =\vert X^*(I_{\star}/\bar{G}_{\star})^{\Gamma_{F}}\vert \vert coker(\phi_{1})\vert .$$
Il en r\'esulte que
$$x(I_{\star})[X^*(S^{\theta})^{\Gamma_{F}}:Im(X^*(G)^{\Gamma_{F},\theta})]^{-1}=\vert coker(\phi_{2})\vert ^{-1}\vert coker(\phi_{1})\vert^{-1}=\vert coker(\phi)\vert ^{-1}.$$
D'apr\`es les d\'efinitions, on a aussi
$$\vert coker(\phi)\vert =\vert \mathfrak{A}_{\tilde{G},{\mathbb Z}}/\mathfrak{A}_{\bar{G}_{\star},{\mathbb Z}}\vert =covol(\mathfrak{A}_{\bar{G}_{\star},{\mathbb Z}})covol(\mathfrak{A}_{\tilde{G},{\mathbb Z}})^{-1}.$$
L'\'egalit\'e (7) r\'esulte des deux \'egalit\'es ci-dessus. Cela ach\`eve la d\'emonstration. $\square$

\bigskip

\subsection{ Preuve de la proposition 7.10}
Les lemmes 7.12, 7.13 et 7.14 conduisent \`a 
  l'\'egalit\'e
$$\vert P^0\vert \vert {\mathbb U}\vert ^{-1}=C(\tilde{G})\tau(\bar{G}_{\star})[I_{\star}(F):\bar{G}_{\star}(F)] ^{-1}[I_{\star}(F_{V}):\bar{G}_{\star}(F_{V})] covol(\mathfrak{A}_{\bar{G}_{\star}})^{-1}.$$
On a aussi $\vert {\mathbb U}\vert =\vert \dot{D}_{F}[d_{V}]\vert $ d'apr\`es le lemme 7.5. Enfin, l'assertion 7.3(1) entra\^{\i}ne que pour tout $d\in D_{F}[d_{V}]$, il y a des isomorphismes compatibles et d\'efinis sur $F$ de $I_{\star}$ sur $I_{\eta[d]}$ et de $\bar{G}_{\star}$ sur $G_{\eta[d]}$. On peut donc remplacer dans le deuxi\`eme membre de l'\'egalit\'e ci-dessus les termes $I_{\star}$ et $\bar{G}_{\star}$ par $I_{\eta[d]}$ et $G_{\eta[d]}$. La proposition 7.10 en r\'esulte. 

\bigskip

\subsection{ Calcul final}
On l\`eve l'hypoth\`ese $D_{F}[d_{V}]\not=\emptyset$ pos\'ee en 7.3.
Consid\'erons l'expression

(1) $  \bar{f}[d_{V}]\sum_{j\in {\cal J}({\bf H})}\delta_{j}[d_{V}]$

\noindent qui appara\^{\i}t dans la formule 5.9(2). 

\ass{Corollaire}{Pourvu que l'ensemble de places $V'$ soit assez grand, l'expression (1) est \'egale \`a
$$C(\tilde{G})^{-1} \sum_{d\in \dot{D}_{F}[d_{V}]}\tau(G_{\eta[d]})[ I_{\eta[d]}(F):G_{\eta[d]}(F)] ^{-1}[ I_{\eta[d]}(F_{V}):G_{\eta[d]}(F_{V})]  covol(\mathfrak{A}_{G_{\eta[d]}})^{-1}$$
$$\vert {\cal U}[V',d]\vert ^{-1}   \sum_{u\in {\cal U}[V',d]}\omega(uh[d]) \bar{f}[d,u].$$}

Cela r\'esulte du corollaire 7.2 et de la proposition 7.10.

\bigskip

\section{Preuve du th\'eor\`eme 3.3}

\bigskip

\subsection{Suite du calcul de la section 5}
Utilisons le corollaire 7.16 pour transformer la formule 5.9(2). Il y a deux simplifications. La constante $C(\tilde{G})$ appara\^{\i}t deux fois et ses deux occurences se compensent. Le terme $c[d_{V}]$ d\'efini en 5.9 est \'egal \`a $[I_{\eta[d]}(F_{V}):G_{\eta[d]}(F_{V})]^{-1}$ pour tout $d\in \dot{D}_{F}[d_{V}]$ et ce terme compense l'un de ceux intervenant dans le corollaire 7.16. D'autre part, il intervient dans ce corollaire un ensemble fini de places $V'$ qui doit \^etre assez grand. Cette notion d\'epend a priori des \'el\'ements fix\'es dans la section 7, \`a savoir $d_{V}$ et ${\bf H}$.  Mais ces donn\'ees $d_{V}$ et ${\bf H}$ parcourent des ensembles finis. On peut donc fixer un ensemble $V'$ assez grand pour que le corollaire 7.16 soit valable pour tous $d_{V}$ et ${\bf H}$.  
On obtient alors
$$I^{\tilde{G}}(\underline{A}^{\tilde{G},{\cal E}}(V,{\bf H},\omega),f)=\vert W^G(\mu)\vert^{-1} \vert Fix^G(\mu,\omega_{\bar{G}})\vert ^{-1} \tau(\bar{H})^{-1}\vert W^{\bar{H}}\vert $$
$$\sum_{d_{V}\in \dot{D}_{V}^{rel}}\sum_{d\in \dot{D}_{F}[d_{V}]}\tau(G_{\eta[d]})[ I_{\eta[d]}(F):G_{\eta[d]}(F)] ^{-1}covol(\mathfrak{A}_{G_{\eta[d]},{\mathbb Z}})^{-1}$$
$$\vert {\cal U}[V',d]\vert ^{-1}\sum_{u\in {\cal U}[V',d]}\omega(uh[d])S^{\bar{H}}(SA_{unip}^{\bar{H}}(V),\bar{f}[d,u]).$$
 Notons $D_{F}^{V-nr}$ l'ensemble des $d\in D_{F}$ tels que la projection de $d$ dans $D_{{\mathbb A}_{F}^V}$ appartienne \`a $D_{{\mathbb A}_{F}^V}^{nr}$. L'ensemble $D_{F}[d_{V}]$ a \'et\'e d\'efini en 6.9. Dans ce paragraphe, on supposait $d_{V}$ relevante mais c'\'etait inutile pour cette d\'efinition. On peut alors d\'efinir $\dot{D}_{F}[d_{V}]$ comme en 7.1 sans cette hypoth\`ese de relevance. Notons $\dot{D}_{F}^{V-nr}$ la r\'eunion des  $\dot{D}_{F}[d_{V}]$ pour $d_{V}\in \dot{D}_{V}$.    C'est un ensemble de repr\'esentants de l'ensemble de doubles classes
$$I_{\eta}(\bar{F})\backslash D_{F}^{V-nr}/G(F).$$
La m\^eme preuve qu'en 7.1(1) montre que c'est un ensemble fini. 
Gr\^ace \`a 5.5(4), la double somme ci-dessus en $d_{V}\in \dot{D}_{V}^{rel}$ et $d\in \dot{D}_{F}[d_{V}]$ se simplifie en une somme sur les $d\in \dot{D}_{F}^{V-nr}$ tels que ${\bf H}$ soit relevante pour $G_{\eta[d],SC}$. D'o\`u
$$I^{\tilde{G}}(\underline{A}^{\tilde{G},{\cal E}}(V,{\bf H},\omega),f)=\vert W^G(\mu)\vert^{-1} \vert Fix^G(\mu,\omega_{\bar{G}})\vert ^{-1} \tau(\bar{H})^{-1}\vert W^{\bar{H}}\vert $$
$$\sum_{d\in \dot{D}_{F}^{V-nr},{\bf H} \text{ relevante pour }G_{\eta[d],SC}}\tau(G_{\eta[d]})[ I_{\eta[d]}(F):G_{\eta[d]}(F)] ^{-1}covol(\mathfrak{A}_{G_{\eta[d]},{\mathbb Z}})^{-1}$$
$$\vert {\cal U}[V',d]^{-1}\sum_{u\in {\cal U}[V',d]}\omega(uh[d])S^{\bar{H}}(SA_{unip}^{\bar{H}}(V),\bar{f}[d,u]).$$
On peut maintenant lever la premi\`ere hypoth\`ese faite sur ${\bf H}$, \`a savoir que $D_{V}^{rel}$ est non vide. Si elle n'est pas v\'erifi\'ee, le membre de droite est nul car la somme en $d$ est vide. Le membre de gauche est nul lui-aussi d'apr\`es le corollaire 5.6. L'\'egalit\'e est donc encore v\'erifi\'ee. 
Rappelons que la fonction $\bar{f}[d,u]$ qui intervient ici est le transfert \`a $\bar{H}(F_{V})$ d'une fonction $f[d]_{sc}$ sur $G_{\eta[d],SC}(F_{V})$, le facteur de transfert \'etant d\'etermin\'e par $u$. Il convient d'introduire la donn\'ee ${\bf H}$ dans la notation et de noter plut\^ot cette fonction $\bar{f}^{{\bf H}}[d,u]$.   
Utilisons maintenant l'\'egalit\'e 5.3(1). On obtient
$$(1) \qquad I^{\tilde{G}}(\underline{A}^{\tilde{G},{\cal E}}(V,\mu,\omega_{\bar{G}},\omega),f)=\vert W^G(\mu)\vert ^{-1}\vert Fix^G(\mu,\omega_{\bar{G}})\vert ^{-1} $$
$$\sum_{d\in \dot{D}_{F}^{V-nr} }\tau(G_{\eta[d]})[ I_{\eta[d]}(F):G_{\eta[d]}(F)] ^{-1}covol(\mathfrak{A}_{G_{\eta[d]},{\mathbb Z}})^{-1}\vert {\cal U}[V',d]\vert ^{-1}\sum_{u\in {\cal U}[V',d]}\omega(uh[d])X[d,u],$$
o\`u
$$ X[d,u]=\sum_{{\bf H}\in E_{\hat{T}_{ad},\star}(\bar{G}_{SC},V),{\bf H}\text{ relevante pour }G_{\eta[d],SC}}  \tau(\bar{H})^{-1}\vert W^{\bar{H}}\vert S^{\bar{H}}(SA_{unip}^{\bar{H}}(V),\bar{f}^{{\bf H}}[d,u]).$$

\bigskip

\subsection{Elimination de la somme en ${\bf H}$}
Fixons $d\in \dot{D}_{F}^{V-nr}$ et $u\in {\cal U}[V',d]$.  
Notons $E_{\hat{\bar{T}}_{ad}}(G_{\eta[d],SC},V)$ l'ensemble de sommation qui intervient dans la d\'efinition de $X[d,u]$ , c'est-\`a-dire l'ensemble des ${\bf H}\in E_{\hat{\bar{T}}_{ad},\star}(\bar{G}_{SC},V)$ qui sont relevantes pour $G_{\eta[d],SC}$. Fixons comme toujours un ensemble ${\cal E}(G_{\eta[d],SC},V)$ de repr\'esentants des classes d'\'equivalence de donn\'ees endoscopiques de $G_{\eta[d],SC}$ qui sont elliptiques, non ramifi\'ees hors de $V$ et relevantes.  Il y a une application naturelle $E_{\hat{\bar{T}}_{ad}}(G_{\eta[d],SC},V)\to {\cal E}(G_{\eta[d],SC},V)$ qui, \`a ${\bf H}$ dans l'ensemble de d\'epart, associe le repr\'esentant de sa classe d'\'equivalence. Cette application est surjective: toute donn\'ee endoscopique de $G_{\eta[d],SC}$ est \'equivalente \`a une donn\'ee $(\bar{H},\bar{{\cal H}},\bar{s})$ telle que $\bar{s}$ appartienne \`a $\hat{\bar{T}}_{ad}$. Les nombres d'\'el\'ements des fibres de cette application se calculent comme en 5.1. On obtient que la fibre contenant un \'el\'ement ${\bf H}\in E_{\hat{\bar{T}}_{ad}}(G_{\eta[d],SC},V)$ a pour nombre d'\'el\'ements $\vert W^{\bar{G}}\vert \vert W^{\bar{H}}\vert ^{-1}\vert Out({\bf H})\vert ^{-1}$.  On obtient
$$X[d,u]=\vert W^{\bar{G}}\vert\sum_{{\bf H}\in {\cal E}(G_{\eta[d],SC},V)}\tau(\bar{H})^{-1}\vert Out({\bf H})\vert ^{-1}S^{\bar{H}}(SA_{unip}^{\bar{H}}(V),\bar{f}^{{\bf H}}[d,u]).$$
En utilisant la d\'efinition [VI] 5.1, on voit que $\tau(\bar{H})^{-1}\vert Out({\bf H})\vert ^{-1}=\tau(G_{\eta[d],SC})^{-1}i(G_{\eta[d],SC},\bar{H})$.  Ici $\tau(G_{\eta[d],SC})=1$ puisque $G_{\eta[d],SC}$ est simplement connexe ([Lab1] th\'eor\`eme 1.2). D'autre part, 
$$S^{\bar{H}}(SA_{unip}^{\bar{H}}(V),\bar{f}^{{\bf H}}[d,u])=I^{G_{\eta[d],SC}}(transfert_{u}(SA_{unip}^{\bar{H}}(V)),f[d]_{sc}),$$
o\`u $transfert_{u}$ d\'esigne le transfert relatif au facteur de transfert $\Delta_{V}[d,u]$ d\'efini en 7.1. On obtient
$$X[d,u]=\vert W^{\bar{G}}\vert  I^{G_{\eta[d],SC}}(B[d,u],f[d]_{sc}),$$
o\`u
$$B[d,u]=\sum_{{\bf H}\in {\cal E}(G_{\eta[d],SC},V)}i(G_{\eta[d],SC},\bar{H})transfert_{u}(SA_{unip}^{\bar{H}}(V)).$$
En se reportant \`a la d\'efinition de 1.10, on voit que $B[d,u]$ n'est autre que la distribution $A^{G_{\eta[d],SC},{\cal E}}_{unip}(V)$. Plus exactement, on se rappelle que cette distribution d\'epend d'un choix de sous-groupes compacts hypersp\'eciaux des groupes $G_{\eta[d],SC}(F_{v})$ pour $v\not\in V$. Puisqu'on utilise le facteur de transfert $\Delta_{V}[d,u]$, ces groupes sont les $K_{sc,v}[d,u]$ d\'efinis en 7.1. Notons $K_{v}[d]$ le sous-groupe compact hypersp\'ecial de $G_{\eta[d]}(F_{v})$ d\'efini par $K_{v}[d]=ad_{h_{v}[d]}(K_{v})\cap G_{\eta[d]}(F_{v})$. Avec la notation de 4.3, on a $\prod_{v\not\in V}K_{sc,v}[d,u]={^{u}K}[d]_{sc}^V$. On a donc pr\'ecis\'ement $B[d,u]=  
 A^{G_{\eta[d],SC},{\cal E}}_{unip}(V, {^{u}K}[d]_{sc}^V)$. On est ici dans une situation sans torsion et on peut appliquer le th\'eor\`eme d'Arthur 3.1. Donc $B[d,u]=A^{G_{\eta[d],SC}}_{unip}(V, {^{u}K}[d]_{sc}^V)$ et
 $$ X[d,u]=\vert W^{\bar{G}}\vert  I^{G_{\eta[d],SC}}(A^{G_{\eta[d],SC}}_{unip}(V, {^{u}K}[d]_{sc}^V),f[d]_{sc}).$$

  \bigskip
  
  \subsection{Elimination des rev\^etements simplement connexes}
 Fixons $d\in \dot{D}_{F}^{V-nr}$. Gr\^ace \`a la derni\`ere formule ci-dessus, on a
  $$\vert {\cal U}[V',d]\vert ^{-1}\sum_{u\in {\cal U}[V',d]}\omega(uh[d])X[d,u]=\vert W^{\bar{G}}\vert  \omega(h[d])I^{G_{\eta[d],SC}}(A^{G_{\eta[d],SC}}_{unip;G_{\eta[d]},\omega,V'}(V),f[d]_{sc}),$$
  o\`u $A^{G_{\eta[d],SC}}_{unip;G_{\eta[d]},\omega,V'}(V)$ est d\'efini en 4.3. Puisque $f[d]_{sc}=\iota_{G_{\eta[d],SC},G_{\eta[d]}}(f[d])$, on a 
  $$I^{G_{\eta[d],SC}}(A^{G_{\eta[d],SC}}_{unip;G_{\eta[d]},\omega,V'}(V),f[d]_{sc})=I^{G_{\eta[d]}}(\iota^*_{G_{\eta[d],SC},G_{\eta[d]}}(A^{G_{\eta[d],SC}}_{unip;G_{\eta[d]},\omega,V'}(V)),f[d]).$$
  On applique la proposition 4.3. L'ensemble de places $V'$ doit \^etre assez grand, cette notion d\'ependant de $d$. Puisqu'on a d\'ej\`a dit que l'ensemble $\dot{D}_{F}^{V-nr}$ \'etait fini, on peut supposer que l'ensemble $V'$ fix\'e en 8.1 satisfait cette condition pour tout $d$. La proposition 4.3 nous dit que 
  $$\iota^*_{G_{\eta[d],SC},G_{\eta[d]}}(A^{G_{\eta[d],SC}}_{unip;G_{\eta[d]},\omega,V'}(V))=\tau'(G_{\eta[d],SC})\tau'(G_{\eta[d]})^{-1}A^{G_{\eta[d]}}_{unip}(V,\omega,K[d]^V).$$
  Puisque $G_{\eta[d],SC}$ est simplement connexe, on a $\mathfrak{A}_{G_{\eta[d],SC}}=\{0\}$ et $\tau(G_{\eta[d],SC})=1$ ([Lab1] th\'eor\`eme 1.2). Donc aussi $\tau'(G_{\eta[d],SC})=1$. On a aussi par d\'efinition 
  $$\tau'(G_{\eta[d]})=\tau(G_{\eta[d]})covol(\mathfrak{A}_{G_{\eta[d]},{\mathbb Z}})^{-1}.$$
   D'o\`u
$$\vert {\cal U}[V',d]\vert ^{-1}\sum_{u\in {\cal U}[V',d]}\omega(uh[d])X[d,u]=\vert W^{\bar{G}}\vert  \omega(h[d])  \tau(G_{\eta[d]})^{-1}covol(\mathfrak{A}_{G_{\eta[d]},{\mathbb Z}})$$
$$I^{G_{\eta[d]}}(A^{G_{\eta[d]}}_{unip}(V,\omega,K[d]^V),f[d]).$$
Reportons cette valeur dans la formule 8.1(1). Comme on l'a dit en 1.1, le groupe $W^G(\mu)$ s'identifie \`a $W^{\bar{G}}$. On obtient
 $$(1) \qquad I^{\tilde{G}}(\underline{A}^{\tilde{G},{\cal E}}(V,\mu,\omega_{\bar{G}},\omega),f)=  \vert Fix^G(\mu,\omega_{\bar{G}})\vert ^{-1}$$
 $$\sum_{d\in \dot{D}_{F}^{V-nr} } [ I_{\eta[d]}(F):G_{\eta[d]}(F)] ^{-1}\omega(h[d])I^{G_{\eta[d]}}(A^{G_{\eta[d]}}_{unip}(V,\omega,K[d]^V),f[d]).$$
 
 \bigskip
 
 \subsection{Fin de la preuve}
 Rappelons qu'il nous suffit de prouver la proposition 5.1, c'est-\`a-dire l'\'egalit\'e
 $$(1)\qquad I^{\tilde{G}}(\underline{A}^{\tilde{G},{\cal E}}(V,\mu,\omega_{\bar{G}},\omega),f)=I^{\tilde{G}}(A^{\tilde{G}}(V,{\cal X},\omega),f).$$
 Supposons que ${\cal X}$ n'appartienne pas \`a l'image de l'application $\chi^{\tilde{G}}$, c'est-\`a-dire ne corresponde \`a aucune classe de conjugaison stable dans $\tilde{G}_{ss}(F)$. Alors l'ensemble $D_{F}$ est vide, a fortiori $\dot{D}_{F}^{V-nr}$ aussi et la formule  8.3(1) montre que $I^{\tilde{G}}(\underline{A}^{\tilde{G},{\cal E}}(V,\mu,\omega_{\bar{G}},\omega),f)=0$. D'apr\`es la d\'efinition 1.9, on a aussi $I^{\tilde{G}}(A^{\tilde{G}}(V,{\cal X},\omega),f)=0$. L'\'egalit\'e (1) ci-dessus en r\'esulte.
 
 Nous supposons maintenant que ${\cal X}$ est l'image par $\chi^{\tilde{G}}$ d'une classe de conjugaison stable ${\cal O}\subset \tilde{G}(F)$. Cette classe est unique d'apr\`es la proposition 1.2 et c'est la classe de conjugaison stable de $\eta[d]$ pour tout $d\in D_{F}$. Notons $Cl({\cal O})$ l'ensemble des classes de conjugaison par $G(F)$ contenues dans ${\cal O}$. Il y a une application naturelle $D_{F}\to Cl({\cal O})$ qui, \`a $d\in D_{F}$, associe la classe de conjugaison de $\eta[d]$. Elle est surjective et se quotiente en une application de $I_{\eta}(\bar{F})\backslash D_{F}/G(F)$ dans $Cl({\cal O})$. 
 Etendons l'ensemble $\dot{D}_{F}^{V-nr}$ en un ensemble de repr\'esentants $\dot{D}_{F}$ de l'ensemble de doubles classes $I_{\eta}(\bar{F})\backslash D_{F}/G(F)$.  L'application ci-dessus devient une application $cl:\dot{D}_{F}\to Cl({\cal O})$. Soit $d\in \dot{D}_{F}$, posons ${\cal C}=cl(d)$. Utilisons la d\'efinition [VI] 2.3 ainsi que les hypoth\`eses (1), (2) et (4) pos\'ees en 5.1. Consid\'erons la propri\'et\'e
 
 (2) pour tout $v\not\in V$, la classe de conjugaison par $G(F_{v})$ engendr\'ee par ${\cal C}$ coupe $\tilde{K}_{v}$. 
 
 Elle \'equivaut \`a ce que, pour tout $v\not\in V$, l'image de $d$ dans $D_{v}$ appartienne \`a $D_{v}^{nr}$, cf. 5.5. Ou encore \`a
  $d\in \dot{D}_{F}^{V-nr}$.
 Si (2) n'est pas v\'erifi\'ee, on a $A^{\tilde{G}}(V,{\cal C},\omega)=0$. Si (2) est v\'erifi\'ee, la d\'efinition [VI] 2.3(7) conduit \`a l'\'egalit\'e
 $$I^{\tilde{G}}(A^{\tilde{G}}(V,{\cal C},\omega),f)=[Z_{G}(\eta[d];F):G_{\eta[d]}(F)] ^{-1}\omega(h[d])I^{G_{\eta[d]}}(A^{G_{\eta[d]}}_{unip}(V,\omega,K[d]^V),f[d]).$$ 
 Il r\'esulte alors de 8.3(1)  que
  $$(3) \qquad I^{\tilde{G}}(\underline{A}^{\tilde{G},{\cal E}}(V,\mu,\omega_{\bar{G}},\omega),f)=  \vert Fix^G(\mu,\omega_{\bar{G}})\vert ^{-1}$$
 $$\sum_{d\in \dot{D}_{F}}[Z_{G}(\eta[d];F):I_{\eta[d]}(F)] I^{\tilde{G}}(A^{\tilde{G}}(V,cl(d),\omega),f).$$
 L'ensemble de sommation peut maintenant \^etre infini mais seuls un nombre fini de termes sont non nuls.
 
 On va prouver
 
 (4) pour tout $d\in \dot{D}_{F}$, la fibre de l'application $cl$ au-dessus de $cl(d)$ a pour nombre d'\'el\'ements $\vert Fix^G(\mu,\omega_{\bar{G}})\vert [Z_{G}(\eta[d];F):I_{\eta[d]}(F)]^{-1}$.
 
 L'\'el\'ement $d\in \dot{D}_{F}$ est fix\'e. On rappelle que l'on note ${\cal Y}_{\eta[d]}$ l'ensemble des $y\in G(\bar{F})$ tels que $y\sigma(y)^{-1}\in I_{\eta[d]}(\bar{F})$ pour tout $\sigma\in \Gamma_{F}$. L'application $y\mapsto (y^{-1}\eta[d]y,r[d]y)$ est une bijection de ${\cal Y}_{\eta[d]}$ sur $D_{F}$ (la preuve est analogue \`a celle de 5.4(4)). L'inverse de cette bijection identifie $\dot{D}_{F}$ \`a un ensemble de repr\'esentants de l'ensemble de doubles classes
 $$I_{\eta[d]}(\bar{F})\backslash {\cal Y}_{\eta[d]}/G(F).$$
  L'\'el\'ement $y^{-1}\eta[d]y$ est conjugu\'e \`a $\eta[d]$ par un \'el\'ement de $G(F)$ si et seulement si $y\in Z_{G}(\eta[d];\bar{F})G(F)$. Ou encore, puisque $y$ est suppos\'e appartenir \`a $ {\cal Y}_{\eta[d]}$, si et seulement si $y\in (Z_{G}(\eta[d];\bar{F})\cap {\cal Y}_{\eta[d]})G(F)$. Donc la fibre de $cl$ au-dessus de $cl(d)$ est en bijection avec l'ensemble de doubles classes
 $$(5) \qquad I_{\eta[d]}(\bar{F})\backslash (Z_{G}(\eta[d];\bar{F})\cap {\cal Y}_{\eta[d]})G(F)/G(F).$$
 
 Posons $\Xi[d]=I_{\eta[d]}(\bar{F})\backslash Z_{G}(\eta[d];\bar{F})$. C'est un groupe fini qui est muni d'une action naturelle de  $\Gamma_{F}$. Il r\'esulte des d\'efinitions que  
  $$\Xi[d]^{\Gamma_{F}}=I_{\eta[d]}(\bar{F})\backslash (Z_{G}(\eta[d];\bar{F})\cap {\cal Y}_{\eta[d]}).$$
 Ce groupe contient contient le sous-groupe $\Xi_{F}[d]=I_{\eta[d]}(F)\backslash Z_{G}(\eta[d];F)$. On voit que l'application naturelle de $ \Xi[d]^{\Gamma_{F}}$ dans l'ensemble (5) se quotiente en une bijection de $\Xi[d]^{\Gamma_{F}}/\Xi_{F}[d]$ sur cet ensemble. On en d\'eduit que le nombre d'\'el\'ements de la fibre de $cl$ au-dessus de $cl(d)$ est \'egal \`a
 $$(6)\qquad \vert  \Xi[d]^{\Gamma_{F}} \vert [Z_{G}(\eta[d];F):I_{\eta[d]}(F)]^{-1}.$$
 
  On va identifier l'ensemble $\Xi[d]$ et son action galoisienne. Munissons le groupe $W^{\theta^*}$ de l'action galoisienne $\sigma\mapsto \sigma_{\bar{G}}$ d\'efinie par $\sigma_{\bar{G}}(w)=\omega_{\bar{G}}(\sigma)\sigma_{G^*}(w)\omega_{\bar{G}}(\sigma)^{-1}$.    Le groupe $W^{\theta^*}$ agit naturellement sur $(T^*/(1-\theta^*)(T^*))\times_{{\cal Z}(G)}{\cal Z}(\tilde{G})$, cf. 1.1.  Notons $Fix^G(\mu)$ le groupe des $w\in W^{\theta^*}$ tels que $w\mu=\mu$ et $w(\Sigma_{+}(\mu))=\Sigma_{+}(\mu)$, cf. 1.1 pour les notations.  Pour tout $\sigma\in \Gamma_{F}$, $\omega_{\bar{G}}(\sigma)\circ \sigma_{G^*}$ fixe $\mu$ et conserve $\Sigma_{+}(\mu)$. On en d\'eduit que l'action  galoisienne $\sigma\mapsto \sigma_{\bar{G}}$ sur $W^{\theta^*}$ conserve $Fix^G(\mu)$. D'apr\`es la d\'efinition de 5.1, $Fix^G(\mu,\omega_{\bar{G}})$ n'est autre que l'ensemble de points fixes par cette action dans $Fix^G(\mu)$. On va montrer
 
 (7) il existe un isomorphisme de $\Xi[d]$ sur $Fix^G(\mu)$ qui entrelace l'action galoisienne naturelle sur $\Xi[d]$ avec l'action $\sigma\mapsto \sigma_{\bar{G}}$ sur $Fix^G(\mu)$. 
 
 En admettant cela, on a
$$ \vert \Xi[d]^{\Gamma_{F}}\vert =\vert Fix^G(\mu,\omega_{\bar{G}})\vert $$
et (4) r\'esulte alors de la formule (6). Prouvons (7).  Posons $\Xi=I_{\eta}\backslash Z_{G}(\eta)=I_{\eta}(\bar{F})\backslash Z_{G}(\eta;\bar{F})$. Soit $\xi\in \Xi$ et relevons $\xi$ en un \'el\'ement $x\in Z_{G}(\eta)$. Quitte \`a multiplier $x$ \`a gauche par un \'el\'ement de $I_{\eta}$, on peut supposer que $ad_{x}$ conserve $T^{*,\theta^*,0}$ et $\bar{B}=B^*\cap \bar{G}$.   Il en r\'esulte que $ad_{x}$ conserve $T^*$ donc $x$ s'envoie sur un \'el\'ement $w\in W$. Cet \'el\'ement est invariant par $\theta^*$.   L'\'el\'ement $x$ est  uniquement d\'etermin\'e modulo multiplication \`a gauche par un \'el\'ement de $T^{*,\theta^*}$. Donc $w$ est bien d\'etermin\'e. Le fait que $x$ appartient \`a $Z_{G}(\eta)$ entra\^{\i}ne que $w$ fixe $\mu$. Parce que $ad_{x}$ conserve $\bar{B}=B^*\cap \bar{G}$, $w$ conserve $\Sigma_{+}(\mu)$. Autrement dit, $w\in Fix^G(\mu)$. Cela d\'efinit une application $\xi\mapsto w$ de $\Xi$ dans $Fix^G(\mu)$. Il est imm\'ediat que c'est un homomorphisme de groupes. Son noyau est l'ensemble des $\xi\in \Xi$ qui se rel\`event en un \'el\'ement $x\in T^*\cap Z_{G}(\eta)$. Mais ce groupe est \'egal \`a $T^{*,\theta^*}$ qui est contenu dans $I_{\eta}$. Donc le noyau est nul. Montrons que l'application $\xi\mapsto w$ est surjective. Soit $w\in Fix^G(\mu)$. Ecrivons $\eta=\nu e$ avec $\nu\in T^* $ et $e\in Z(\tilde{G};{\cal E}^*)$. On  rel\`eve $w$ en un \'el\'ement $x\in G_{e}$ qui normalise $T^*$. Parce que $w$ fixe $\mu$, il existe $t\in T^*$ tel que $ad_{x}(\nu)=(1-\theta^*)(t)\nu$. Ou encore, puisque $ad_{x}$ fixe $e$, $ad_{x}(\eta)=(1-\theta^*)(t)\eta$.  L'\'el\'ement $t^{-1}x$ rel\`eve encore $w$ et appartient \`a $Z_{G}(\eta)$.  En notant $\xi$ l'image de $t^{-1}x$ dans $\Xi$, on voit que $\xi$ s'envoie sur $w$ par l'application pr\'ec\'edente. Cela prouve que l'application $\xi\mapsto w$ est un isomorphisme de $\Xi$ sur $Fix^G(\mu)$. Posons $r=r[d]$. Alors $ad_{r}$ est un isomorphisme de $\Xi[d]$ sur $\Xi$. par composition, on obtient un isomorphisme de $\Xi[d]$ sur $Fix^G(\mu)$. Pour \'etudier les actions galoisiennes, on doit se rappeler que, par d\'efinition de $D_{F}$, $\omega_{\bar{G}}$ co\"{\i}ncide avec le cocycle $\omega_{\eta[d]}$ calcul\'e comme en 1.2 \`a l'aide de la paire de Borel $(ad_{r^{-1}}(B^*),ad_{r^{-1}}(T^*))$. En reprenant les d\'efinitions de 1.2, on voit que cela se traduit par la propri\'et\'e suivante:

(8) pour tout $\sigma\in \Gamma_{F}$, il existe $\bar{g}(\sigma)\in \bar{G}$ tel que $\bar{g}(\sigma)r\sigma(r)^{-1}u_{{\cal E}^*}(\sigma)^{-1}$ normalise $T^*$ et s'envoie sur l'\'el\'ement $\omega_{\bar{G}}(\sigma)\in W$.

 Soient $\sigma\in \Gamma_{F}$ et $\xi[d]\in \Xi[d]$. L'\'el\'ement $\xi[d]$ s'envoie par $ad_{r}$ sur un \'el\'ement $\xi\in \Xi$, que l'on rel\`eve comme ci-dessus en un \'el\'ement $x\in Z_{G}(\eta)$, qui s'envoie sur un \'el\'ement $w\in Fix^G(\mu)$. L'\'el\'ement $\sigma(\xi[d])$ s'envoie par $ad_{r}$ sur un \'el\'ement $\xi'\in \Xi$, que l'on rel\`eve en un \'el\'ement $x'\in Z_{G}(\eta)$, qui s'envoie sur un \'el\'ement $w'\in Fix^G(\mu)$. Relevons $\xi[d]$ en l'\'el\'ement $ad_{r^{-1}}(x)$. Alors $\xi'$ est l'image dans $\Xi$ de l'\'el\'ement $ad_{r}\circ \sigma\circ ad_{r^{-1}}(x)=ad_{r\sigma(r)^{-1}}\circ \sigma(x)$. La conjugaison par $\bar{G}$ conserve $Z_{G}(\eta)$ et se quotiente en l'identit\'e de $\Xi$. Donc $\xi'$ est aussi l'image dans $\Xi$ de l'\'el\'ement $x''=ad_{\bar{g}(\sigma)r\sigma(r)^{-1}}\circ \sigma(x)$. On a
 $$(9) \qquad x''=ad_{\bar{g}(\sigma)r\sigma(r)^{-1}u_{{\cal E}^*}(\sigma)^{-1}}\circ ad_{u_{{\cal E}^*}(\sigma)}\circ \sigma(x).$$
 L'\'el\'ement $x$ normalise $T^*$. Les automorphismes $ad_{\bar{g}(\sigma)r\sigma(r)^{-1}u_{{\cal E}^*}(\sigma)^{-1}}$ et $ad_{u_{{\cal E}^*}(\sigma)}\circ \sigma$ aussi. Il en r\'esulte que $x''$ normalise $T^*$. Puisque $x'$ et $x''$ rel\`event tous deux $\xi'$, $x''$ se d\'eduit de $x$ par multiplication \`a gauche par un \'el\'ement du normalisateur de $T^*$ dans $I_{\eta}$. L'image $w''$ de $x''$ dans $W$ se d\'eduit de l'image $w'$ de $x'$ par multiplication \`a gauche par un \'el\'ement de $W(\mu)$. D'apr\`es (8) et (9), on a $w''=\omega_{\bar{G}}(\sigma)\circ \sigma_{G^*}(w)=\sigma_{\bar{G}}(w)$. Cela implique que $w''\in Fix^G(\mu)$. Puisque $w'$ et $w''$ sont deux \'el\'ements de cet ensemble qui se d\'eduisent l'un de l'autre par multiplication \`a gauche par un \'el\'ement de $W(\mu)$, ils sont \'egaux (parce que les \'el\'ements de $Fix^G(\mu)$ conservent $\Sigma_{+}(\mu)$). On obtient $w'=w''=\sigma_{\bar{G}}(w)$. Mais cela prouve que l'isomorphisme de $\Xi[d]$ sur $Fix^G(\mu)$ entrelace l'action galoisienne naturelle sur $\Xi[d]$ avec l'action $\sigma\mapsto \sigma_{\bar{G}}$ sur $Fix^G(\mu)$. Cela ach\`eve la preuve de (7) et de (4). 

 En utilisant (4), l'\'egalit\'e (3) se r\'ecrit
 $$I^{\tilde{G}}(\underline{A}^{\tilde{G},{\cal E}}(V,\mu,\omega_{\bar{G}},\omega),f)=\sum_{{\cal C}\in Cl({\cal O})}I^{\tilde{G}}(A^{\tilde{G}}(V,{\cal C},\omega),f).$$
 D'apr\`es la d\'efinition de 1.9, le membre de droite ci-dessus est \'egal \`a celui de l'\'egalit\'e (1). Cela prouve cette \'egalit\'e et cela ach\`eve la preuve du th\'eor\`eme 3.3.
 
 \bigskip
 
 \section{Preuve du th\'eor\`eme [VI] 5.6}
 
 \bigskip
 
 \subsection{Rappel de l'\'enonc\'e du th\'eor\`eme}
 On rappelle bri\`evement l'\'enonc\'e du th\'eor\`eme [VI] 5.6, en renvoyant \`a cette r\'ef\'erence pour plus de d\'etails. On consid\`ere un triplet endoscopique non standard $(G_{1},G_{2},j_{*})$ d\'efini sur $F$. Pour $i=1,2$, on fixe une paire de Borel $(B_{i},T_{i})$ de $G_{i}$ d\'efinie sur $F$. On note $\Sigma_{i}$ l'ensemble des racines de $T_{i}$ dans $\mathfrak{g}_{i}$. Le terme $j_{*}$ est un isomorphisme de $X_{*}(T_{1})\otimes_{{\mathbb Z}}{\mathbb Q}$ sur  $X_{*}(T_{2})\otimes_{{\mathbb Z}}{\mathbb Q}$. Il y a une bijection $\tau:\Sigma_{2}\to \Sigma_{1}$ et une fonction $b:\Sigma_{2}\to {\mathbb Q}_{>0}$ de sorte que, pour tout $\alpha_{2}\in \Sigma_{2}$, on ait l'\'egalit\'e $j_{*}(\check{\alpha}_{1})=b(\alpha_{2})\check{\alpha}_{2}$, o\`u $\alpha_{1}=\tau(\alpha_{2})$ et $\check{\alpha}_{i}$ est la coracine associ\'ee \`a $\alpha_{i}$.
 
 Pour toute place $v$ de $F$, on a un isomorphisme
 $$(1) \qquad D^{st}_{g\acute{e}om}(\mathfrak{g}_{1}(F_{v}))\otimes Mes(G_{1}(F_{v}))^*\simeq D^{st}_{g\acute{e}om}(\mathfrak{g}_{2}(F_{v}))\otimes Mes(G_{2}(F_{v}))^*$$
 qui se restreint en un isomorphisme analogue pour les distributions \`a support nilpotent.
 Par l'exponentielle,  il s'en d\'eduit un isomorphisme
  $$(2) \qquad D^{st}_{unip}(G_{1}(F_{v}))\otimes Mes(G_{1}(F_{v}))^*\simeq D^{st}_{unip}(G_{2}(F_{v}))\otimes Mes(G_{2}(F_{v}))^*.$$
 On fixe un ensemble fini de places $V$ de sorte que
 
 - $V$ contient les places archim\'ediennes de $F$;
 
 - $G_{1}$ et $G_{2}$ sont non ramifi\'es hors de $V$;
 
 - pour $v\not\in V$, notons $p$ la caract\'eristique r\'esiduelle de $F_{v}$ et $e_{v}=[F_{v}:{\mathbb Q}_{p}]$; alors $p>e_{v}N(G_{i})+1$ pour $i=1,2$, o\`u $N(G_{i})$ est l'entier  d\'efini en [W1] 4.3;
 
 - les valeurs de la fonction $b$ sont des unit\'es hors de $V$.

 \ass{Th\'eor\`eme}{Sous ces hypoth\`eses, les distributions $SA_{unip}^{G_{1}}(V)$ et $SA_{unip}^{G_{2}}(V)$ se correspondent par le produit tensoriel sur les $v\in V$ des isomorphismes (2)  ci-dessus.}

 Pour simplifier, on  note encore $j_{*}$ toute application d\'eduite de l'application $j_{*}$ primitive. Par exemple, on note $j_{*}$ les isomorphismes (1) et (2). Pour $i=1,2$, on pose $M_{i,0}=T_{i}$. On a un isomorphisme $j_{*}:\mathfrak{A}_{M_{1,0}}\simeq \mathfrak{A}_{M_{2,0}}$. On fixe des mesures sur ces espaces qui se correspondent par cet isomorphisme. 
 On se d\'ebarrasse des espaces de mesures  comme dans les paragraphes pr\'ec\'edents. C'est-\`a-dire que, en une place $v\not\in V$, on  utilise les mesures canoniques. Pour $i=1,2$, on munit  $G_{i}({\mathbb A})$ de la mesure de Tamagawa, cf. 4.1. De ces choix se d\'eduit une mesure sur $G_{i}(F_{V})$. Pour des Levi $M_{1}\in {\cal L}(M_{1,0})$ et $M_{2}\in {\cal L}(M_{2,0})$ qui se correspondent, on fait des choix analogues. 
 \bigskip
 
 \subsection{Le lemme fondamental pond\'er\'e non standard}
 Dans ce paragraphe et le suivant, on fixe une place $v\not\in V$ et  {\bf on consid\`ere que le corps de base est $F_{v}$}. Pour $i=1,2$, soit $M_{i}\in {\cal L}(M_{i,0})$. On a d\'efini en [II] 4.2 une fonction sur $D^{st}_{g\acute{e}om}(M_{i}(F_{v}))$. Dans la situation g\'en\'erale de cette r\'ef\'erence, elle \'etait not\'ee $s_{\tilde{M}_{i}}^{\tilde{G}_{i}}(. ,\tilde{K}_{i})$. Elle ne d\'ependait  que de la classe de conjugaison par $G_{i,AD}(F_{v})$ du sous-espace hypersp\'ecial $\tilde{K}_{i}$. Ici, la situation n'est pas tordue, on prend  $\tilde{K}_{i}=K_{i}$  et la classe de conjugaison par  $G_{i,AD}(F_{v})$ de ce groupe est uniquement d\'etermin\'ee. On peut donc noter simplement $s_{M_{i}}^{G_{i}}$ notre fonction. 
 
 Supposons que $M_{1}$ et $M_{2}$ se correspondent. On a d\'efini en [III] 6.5 une constante $c_{M_{1},M_{2}}^{G_{1},G_{2}}\in {\mathbb C}^{\times}$.

 \ass{Proposition}{Soit $\boldsymbol{\delta}_{1}\in D^{st}_{g\acute{e}om}(M_{1}(F_{v}))$. Supposons que son  support  est  assez voisin de l'origine. Alors on a l'\'egalit\'e
 $$s_{M_{1}}^{G_{1}}(\boldsymbol{\delta}_{1})= c_{M_{1},M_{2}}^{G_{1},G_{2}}s_{M_{2}}^{G_{2}}(j_{*}(\boldsymbol{\delta}_{1})).$$} 
 
 Preuve. Pour $i=1,2$, il y a une analogue $s_{\mathfrak{m}_{i}}^{\mathfrak{g}_{i}}$ de la fonction $s_{M_{i}}^{G_{i}}$, d\'efinie sur $D^{st}_{g\acute{e}om}(\mathfrak{m}_{i}(F_{v}))$. 
  Pour $d_{1}\in D^{st}_{g\acute{e}om}(\mathfrak{m}_{1}(F_{v}))$ \`a support r\'egulier dans $\mathfrak{g}_{1}(F_{v})$, on a l'\'egalit\'e
$$s_{\mathfrak{m}_{1}}^{\mathfrak{g}_{1}}(d_{1})=  c_{M_{1},M_{2}}^{G_{1},G_{2}}s_{\mathfrak{m}_{2}}^{\mathfrak{g}_{2}}(j_{*}(d_{1})).$$
Cette \'egalit\'e est la conjecture 3.7 de [W2]. Une preuve est annonc\'ee par Chaudouard et Laumon ([CL]).

    Pour $d_{i}\in D^{st}_{g\acute{e}om}(\mathfrak{m}_{i}(F_{v}))$ \`a support assez voisin de l'origine, on peut d\'efinir $exp(d_{i})\in D^{st}_{g\acute{e}om}(G_{i}(F_{v}))$. Il r\'esulte des d\'efinitions que 
 l'on  a l'\'egalit\'e $s_{M_{i}}^{G_{i}}(exp(d_{i}))=s_{\mathfrak{m}_{i}}^{\mathfrak{g}_{i}}(d_{i})$. Pour $\boldsymbol{\delta}_{1}\in D^{st}_{g\acute{e}om}(M_{1}(F_{v}))$ \`a support   r\'egulier dans $G_{1}(F_{v})$ et assez voisin de l'origine, il existe $d_{1}\in  D^{st}_{g\acute{e}om}(\mathfrak{m}_{1}(F_{v}))$ \`a support  r\'egulier dans $\mathfrak{g}_{1}(F_{v})$ et assez voisin de l'origine de sorte que $\boldsymbol{\delta}_{1}=exp(d_{1})$. Les consid\'erations ci-dessus entra\^{\i}nent l'\'egalit\'e de l'\'enonc\'e pour un tel $\boldsymbol{\delta}_{1}$.
 
  Il reste \`a lever l'hypoth\`ese que le support de $\boldsymbol{\delta}_{1}$ est r\'egulier dans $G_{1}(F_{v})$. Cela se fait en deux temps comme dans la section 4 de [II]. Consid\'erons d'abord une classe de conjugaison stable semi-simple ${\cal O}\subset M_{1}(F_{v})$ qui est $G_{1}$-\'equisinguli\`ere et assez proche de l'origine. Soit $\boldsymbol{\delta}_{1}\in D_{g\acute{e}om}^{st}({\cal O})$, c'est-\`a-dire que son support est form\'e d'\'el\'ements dont la partie semi-simple appartient \`a ${\cal O}$. D'apr\`es le lemme [II] 2.2, on peut fixer $\boldsymbol{\delta}'_{1}\in D_{g\acute{e}om}^{st}(M_{1}(F_{v}))$, \`a support r\'egulier dans $G_{1}(F_{v})$ aussi voisin qu'on le veut de ${\cal O}$, de sorte que $\boldsymbol{\delta}_{1}=g_{M_{1}}^{M_{1}}(\boldsymbol{\delta}'_{1})$. D'apr\`es la proposition [III] 6.7, on a l'\'egalit\'e $j_{*}(\boldsymbol{\delta}_{1})=g_{M_{2}}^{M_{2}}\circ j_{*}(\boldsymbol{\delta}'_{1})$. D'apr\`es la relation [II] 4.5(2), on a les \'egalit\'es
 $$s_{M_{1}}^{G_{1}}(\boldsymbol{\delta}'_{1})=s_{M_{1}}^{G_{1}}(g_{M_{1}}^{M_{1}}(\boldsymbol{\delta}'_{1}))=s_{M_{1}}^{G_{1}}(\boldsymbol{\delta}_{1}),$$
 $$s_{M_{2}}^{G_{2}}(j_{*}(\boldsymbol{\delta}'_{1}))=s_{M_{2}}^{G_{2}}(g_{M_{2}}^{M_{2}}\circ j_{*}(\boldsymbol{\delta}'_{1}))=s_{M_{2}}^{G_{2}}(j_{*}(\boldsymbol{\delta}_{1})).$$
 L'\'egalit\'e de l'\'enonc\'e d\'ej\`a prouv\'ee pour $\boldsymbol{\delta}'_{1}$ entra\^{\i}ne alors l'\'egalit\'e analogue pour $\boldsymbol{\delta}_{1}$.
 
Soit maintenant $\boldsymbol{\delta}_{1}$ soumis \`a la seule restriction que son support est assez voisin de l'origine. Posons $\boldsymbol{\delta}_{2}=j_{*}(\boldsymbol{\delta}_{1})$. Pour $i=1,2$, on introduit une variable $a_{i}\in A_{M_{i}}(F_{v})$ en position g\'en\'erale et proche de $1$. D'apr\`es [II] 4.7, le germe en $1$ de la  fonction $a_{i}\mapsto s_{M_{i}}^{G_{i}}(a_{i}\boldsymbol{\delta}_{i})$ est \'equivalent \`a un \'el\'ement de  l'espace $U_{M_{i}}^{G_{i}}$ dont le terme constant est \'egal \`a  $s_{M_{i}}^{G_{i}}(\boldsymbol{\delta}_{i})$ (cf. [II] 4.6 pour les d\'efinitions).   Consid\'erons l'application qui, \`a une fonction $\varphi$ de $a_{2}$, associe la fonction $a_{1}\mapsto \varphi(j_{*}(a_{1}))$. Les consid\'erations de [III] 6.5 et l'hypoth\`ese $v\not\in V$ entra\^{\i}nent que cette application respecte l'\'equivalence. Elle envoie $U_{M_{2}}^{G_{2}}$ sur $U_{M_{1}}^{G_{1}}$ et sa restriction \`a $U_{M_{2}}^{G_{2}}$ commute \`a l'application "terme constant". Il en r\'esulte que le germe en $1$ de la fonction $a_{1}\mapsto s_{M_{2}}^{G_{2}}(j_{*}(a_{1})\boldsymbol{\delta}_{2})$ est \'equivalent \`a un \'el\'ement de  l'espace $U_{M_{1}}^{G_{1}}$ dont le terme constant est \'egal \`a  $s_{M_{2}}^{G_{2}}(\boldsymbol{\delta}_{2})$.  Puisque $a_{1}\boldsymbol{\delta}_{1}$ est \`a support $G_{1}$-\'equisingulier, on a d\'emontr\'e l'\'egalit\'e
$$s_{M_{1}}^{G_{1}}(a_{1}\boldsymbol{\delta}_{1})= c_{M_{1},M_{2}}^{G_{1},G_{2}}s_{M_{2}}^{G_{2}}(j_{*}(a_{1})\boldsymbol{\delta}_{2}).$$
Le germe de cette fonction en $1$ est \'equivalent \`a un \'el\'ement de $U_{M_{1}}^{G_{1}}$ et on  a deux fa\c{c}ons de calculer son terme constant: la premi\`ere donne $s_{M_{1}}^{G_{1}}(\boldsymbol{\delta}_{1})$; la seconde donne $c_{M_{1},M_{2}}^{G_{1},G_{2}}s_{M_{2}}^{G_{2}}(\boldsymbol{\delta}_{2})$. D'o\`u l'\'egalit\'e de ces deux termes, ce qui ach\`eve la d\'emonstration. $\square$

\bigskip

\subsection{Extension aux Levi}
 Pour $i=1,2$, soient $M_{i}$ et $L_{i}$ deux Levi de $G_{i}$ tels que $M_{i,0}\subset M_{i}\subset L_{i}$. Les d\'efinitions de [II] 4.2 se simplifient comme dans le paragraphe pr\'ec\'edent: on a une fonction $s_{M_{i}}^{L_{i}}$ sur 
 $D^{st}_{g\acute{e}om}(M_{i}(F_{v}))$. On suppose que $M_{1}$ et $M_{2}$, resp. $L_{1}$ et $L_{2}$,  se correspondent.

 On note comme toujours $L_{i,SC}$ le rev\^etement simplement connexe de $L_{i}$ et on note $M_{i,sc}$, resp. $T_{i,sc}$, l'image r\'eciproque de $M_{i}$, resp. $T_{i}$, dans $L_{i,SC}$. De $j_{*}$ se d\'eduit un isomomorphisme $j_{*,sc}:X_{*}(T_{1,sc})\otimes_{{\mathbb Z}}{\mathbb Q}\to X_{*}(T_{2,sc})\otimes_{{\mathbb Z}}{\mathbb Q}$. Le triplet $(L_{1,SC},L_{2,SC},j_{*,sc})$ est encore endoscopique non standard. 
 
  On a d\'efini en [III] 3.5 un homomorphisme surjectif $\iota^*_{M_{i,sc},M_{i}}:D^{st}_{unip}(M_{i,sc}(F_{v}))\to D^{st}_{unip}(M_{i}(F_{v}))$.

 \ass{Lemme}{(i) Pour $i=1,2$ et pour $\boldsymbol{\delta}_{i,sc}\in D^{st}_{unip}(M_{i,sc}(F_{v}))$, on a l'\'egalit\'e
 $$s_{M_{i}}^{L_{i}}(\iota^*_{M_{i,sc},M_{i}}(\boldsymbol{\delta}_{i,sc}))=s_{M_{i,sc}}^{L_{i,SC}}(\boldsymbol{\delta}_{i,sc}).$$
 
 (ii) Pour $\boldsymbol{\delta}_{1}\in D^{st}_{unip}(M_{1}(F_{v}))$, on a l'\'egalit\'e
 $$s_{M_{1}}^{L_{1}}(\boldsymbol{\delta}_{1})=c_{M_{1,sc},M_{2,sc}}^{L_{1,SC},L_{2,SC}}s_{M_{2}}^{L_{2}}(j_{*}(\boldsymbol{\delta}_{1})).$$}
 
 Preuve de (i). Le temps de cette preuve, l'indice $i=1,2$ est fix\'e et on le supprime pour simplifier.  Rappelons que les fonctions $s_{M}^L$ et $s_{M_{sc}}^{L_{SC}}$ sont d\'eduites de fonctions primitives $r_{M}^L(.,K^L)$ et $r_{M_{sc}}^{L_{SC}}(.,K^L_{sc})$ qui calculent les int\'egrales orbitales pond\'er\'ees non-invariantes des fonctions caract\'eristiques  des compacts $K^L$ et $K^L_{sc}$. Fixons un ensemble de repr\'esentants ${\cal U}$ de l'ensemble de  doubles classes $Z(L)^0(F_{v})\backslash M(F_{v})/\pi(M_{sc}(F_{v}))$. Soit $\boldsymbol{\gamma}_{sc}\in D_{unip}(M_{sc}(F_{v}))$. Posons $\boldsymbol{\gamma}'_{sc}=\vert {\cal U}\vert ^{-1}\sum_{u\in {\cal U}}ad_{u}(\boldsymbol{\gamma}_{sc})$. Il r\'esulte du lemme [III] 3.3 que l'on a l'\'egalit\'e
 $$r_{M}^L(\iota^*_{M_{sc},M}(\boldsymbol{\gamma}_{sc}),K^L)=r_{M_{sc}}^{L_{SC}}(\boldsymbol{\gamma}'_{sc},K^L_{sc}).$$
 
 {\bf Remarque.} Dans le lemme cit\'e apparaissaient des espaces de mesures que l'on a ici fait dispara\^{\i}tre. En utilisant la preuve du lemme 4.2 du pr\'esent article, on   montre que cette disparition est l\'egitime compte tenu de nos  choix de mesures canoniques.
  
  \bigskip
  
 Une distribution stable est invariante par l'action du groupe adjoint. Pour $\boldsymbol{\delta}_{sc}\in D^{st}_{unip}(M_{sc}(F_{v}))$, on a donc $\boldsymbol{\delta}'_{sc}=\boldsymbol{\delta}_{sc}$ et l'\'egalit\'e se simplifie en
 $$r_{M}^L(\iota^*_{M_{sc},M}(\boldsymbol{\delta}_{sc}),K^L)=r_{M_{sc}}^{L_{SC}}(\boldsymbol{\delta}_{sc},K^L_{sc}).$$
 Cette \'egalit\'e se propage alors aux fonctions $s_{M}^L$ et $s_{M_{sc}}^{L_{SC}}$ par la m\^eme preuve qu'en [III] 3.6. Cela prouve le (i) de l'\'enonc\'e.
 
 Preuve de (ii). Le triplet $(L_{1,SC},L_{2,SC},j_{*,sc})$ \'etant endoscopique non standard, la proposition 8.2 fournit  l'\'egalit\'e
 $$s_{M_{1,sc}}^{L_{1,SC}}( \boldsymbol{\delta}_{1,sc})=c_{M_{1,sc},M_{2,sc}}^{L_{1,SC},L_{2,SC}}s_{M_{2,sc}}^{L_{2,SC}}( j_{*,sc}(\boldsymbol{\delta}_{1,sc}))$$
 pour tout $\boldsymbol{\delta}_{1,sc}\in D^{st}_{unip}(M_{1,sc}(F_{v}))$. En utilisant (i), on obtient
$$s_{M_{1}}^{L_{1}}(\iota^*_{M_{1,sc},M_{1}}(\boldsymbol{\delta}_{1,sc}))= c_{M_{1,sc},M_{2,sc}}^{L_{1,SC},L_{2,SC}}s_{M_{2}}^{L_{2}}(\iota^*_{M_{2,sc},M_{2}}\circ j_{*,sc}(\boldsymbol{\delta}_{1,sc})).$$
Evidemment, $\iota^*_{M_{2,sc},M_{2}}\circ j_{*,sc}=j_{*}\circ \iota^*_{M_{1,sc},M_{1}}$. Puisque tout \'el\'ement $\boldsymbol{\delta}_{1}\in D^{st}_{unip}(M_{1}(F_{v}))$ peut s'\'ecrire sous la forme $\boldsymbol{\delta}_{1}=\iota^*_{M_{1,sc},M_{1}}(\boldsymbol{\delta}_{1,sc})$, on en d\'eduit le (ii) de l'\'enonc\'e. $\square$
 
 \bigskip
 
 \subsection{Globalisation}
 
On revient \`a notre corps de base $F$.  Soit $U$ un ensemble fini de places de $F$ disjoint de $V$ et non vide. Pour $i=1,2$, soit $M_{i}\in {\cal L}(M_{i,0})$. Comme en 9.2, les constructions du paragraphe 2.2 se simplifient dans notre situation et donnent naissance \`a une fonction que l'on note $s_{M_{i},U}^{G_{i}}$ sur $D_{g\acute{e}om}^{st}(M_{i}(F_{U}))$. Supposons que $M_{1}$ et $M_{2}$ se correspondent. On d\'efinit une constante $c_{M_{1},M_{2}}^{G_{1},G_{2}}$ de la m\^eme fa\c{c}on qu'en [III] 6.5, en y rempla\c{c}ant le corps de base local de cette r\'ef\'erence par notre corps $F$. Rappelons la d\'efinition. Soit $n\geq1$ un entier tel que $nj_{*}$ envoie $X_{*}(T_{1})$ dans $X_{*}(T_{2})$. L'homomorphisme dual de $nj_{*}$ envoie $X^*(T_{2})$ dans $X^*(T_{1})$. Il s'identifie \`a un homomorphisme de $X_{*}(\hat{T}_{2})$ dans $X_{*}(\hat{T}_{1})$, qui d\'efinit un homomorphisme de $\hat{T}_{2}$ dans $\hat{T}_{1}$. Celui-ci se restreint en un homomorphisme
$$\hat{j}_{n}:Z(\hat{M}_{2})^{\Gamma_{F}}\to Z(\hat{M}_{1})^{\Gamma_{F}}$$
qui est surjectif et de noyau fini. On pose
$$c_{M_{1},M_{2}}^{G_{1},G_{2}}=n^{-a_{M_{2}}}\vert ker(\hat{j}_{n})\vert ,$$
o\`u $a_{M_{2}}=dim(A_{M_{2}})=dim(A_{M_{1}})$. Cela ne d\'epend pas du choix de $n$. 
\ass{Proposition}{Soit $\boldsymbol{\delta}_{1}\in D_{unip}^{st}(M_{1}(F_{U}))$. On a l'\'egalit\'e
$$s_{M_{1},U}^{G_{1}}(\boldsymbol{\delta}_{1})=c_{M_{1},M_{2}}^{G_{1},G_{2}}s_{M_{2},U}^{G_{2}}(j_{*}(\boldsymbol{\delta}_{1})).$$}

Preuve. On peut supposer $\boldsymbol{\delta}_{1}=\otimes_{v\in U}\boldsymbol{\delta}_{1,v}$. On a \'ecrit en 2.2  la formule de d\'ecomposition
$$(1) \qquad s_{M_{1},U}^{G_{1}}(\boldsymbol{\delta}_{1})=\sum_{L_{1}^U\in {\cal L}(M_{1,U})} e^{G_{1}}_{M_{1,U}}(M_{1},L_{1}^U)\prod_{v\in U}s_{M_{1,v}}^{L_{1}^v}(\boldsymbol{\delta}_{1,v}),$$
Pour tout $L_{1}^U\in {\cal L}(M_{1,U})$ et tout $v\in U$, notons $L_{2}^v$ l'\'el\'ement de ${\cal L}(M_{2,v})$ qui correspond \`a $L_{1}^v$. Le lemme 9.3(ii) entra\^{\i}ne l'\'egalit\'e
$$s_{M_{1,v}}^{L_{1}^v}(\boldsymbol{\delta}_{1,v})=c_{M_{1,v,sc},M_{2,v,sc}}^{L_{1,SC}^v,L_{2,SC}^v}s_{M_{2,v}}^{L_{2}^v}(j_{*}(\boldsymbol{\delta}_{1,v})).$$
Posons $L_{2}^U=(L_{2}^v)_{v\in U}$.  Cette famille est un \'el\'ement de ${\cal L}(M_{2,U})$. On prouvera ci-dessous l'\'egalit\'e
$$(2) \qquad e^{G_{1}}_{M_{1,U}}(M_{1},L_{1}^U)\prod_{v\in U}c_{M_{1,v,sc},M_{2,v,sc}}^{L_{1,SC}^v,L_{2,SC}^v}=c_{M_{1},M_{2}}^{G_{1},G_{2}}e^{G_{2}}_{M_{2,U}}(M_{2},L_{2}^U).$$
Admettons-la. La formule (1) se transforme en
$$ s_{M_{1},U}^{G_{1}}(\boldsymbol{\delta}_{1})=c_{M_{1},M_{2}}^{G_{1},G_{2}}\sum_{L_{1}^U\in {\cal L}(M_{1,U})}e^{G_{2}}_{M_{2,U}}(M_{2},L_{2}^U)\prod_{v\in U}s_{M_{2,v}}^{L_{2}^v}(j_{*}(\boldsymbol{\delta}_{1,v})).$$
L'application $L_{1}^{U}\mapsto L_{2}^{U}$ est une bijection de ${\cal L}(M_{1,U})$ sur ${\cal L}(M_{2,U})$. On peut remplacer la somme en $L_{1}^{U}$ ci-dessus par une somme en $L_{2}^{U}\in {\cal L}(M_{2,U})$. Gr\^ace \`a la formule similaire \`a (1), le membre de droite ci-dessus est \'egal \`a
$$c_{M_{1},M_{2}}^{G_{1},G_{2}}s_{M_{2},U}^{G_{2}}(j_{*}(\boldsymbol{\delta}_{1})),$$
ce qui d\'emontre l'\'egalit\'e de l'\'enonc\'e.

Il reste \`a prouver (2). D'apr\`es la d\'efinition de [VI] 4.2, on a pour $i=1,2$ l'\'egalit\'e 
$$ e^{G_{i}}_{M_{i,U}}(M_{i},L_{i}^U)=d^{G_{i}}_{M_{i,U}}(M_{i},L_{i}^U)k^{G_{i}}_{M_{i,U}}(M_{i},L_{i}^U)^{-1}.$$
Le terme $d^{G_{i}}_{M_{i,U}}(M_{i},L_{i}^U)$ est le rapport entre deux mesures sur l'espace ${\cal A}_{M_{i,U}}^{G_{i}}$, cf. [VI] 1.4. Il r\'esulte de nos d\'efinitions que les espaces et mesures ne d\'ependent pas de l'indice $i$. Donc 
$$d^{G_{1}}_{M_{1,U}}(M_{1},L_{1}^U)=d^{G_{2}}_{M_{2,U}}(M_{2},L_{2}^U).$$
Il suffit donc de prouver l'\'egalit\'e
$$(3) \qquad k^{G_{2}}_{M_{2,U}}(M_{2},L_{2}^U)\prod_{v\in U}c_{M_{1,v,sc},M_{2,v,sc}}^{L_{1,SC}^v,L_{2,SC}^v}=c_{M_{1},M_{2}}^{G_{1},G_{2}}k^{G_{1}}_{M_{1,U}}(M_{1},L_{1}^U).$$
Fixons un entier $n$ comme au d\'ebut du paragraphe. L'application $\hat{j}_{n}$ a des avatars locaux $\hat{j}_{n,v}$. Consid\'erons le diagramme
$$\begin{array}{ccc}Z(\hat{M}_{2})^{\Gamma_{F}}&\stackrel{\hat{j}_{n}}{\to}&Z(\hat{M}_{1})^{\Gamma_{F}}\\ \downarrow&&\downarrow\\ \prod_{v\in V}Z(\hat{M_{2,v}})^{\Gamma_{F_{v}}}/Z(\hat{L_{2}^v})^{\Gamma_{F_{v}}}&\stackrel{\prod_{v\in U}\hat{j}_{n,v}}{\to}&\prod_{v\in V}Z(\hat{M_{1,v}})^{\Gamma_{F_{v}}}/Z(\hat{L_{1}^v})^{\Gamma_{F_{v}}}\\ \end{array}$$
Les fl\`eches verticales sont les homomorphismes naturels. Le diagramme est commutatif. Tous les homomorphismes sont surjectifs,  de noyaux finis. Il r\'esulte des d\'efinitions que

- $k^{G_{2}}_{M_{2,U}}(M_{2},L_{2}^U)$ est le nombre d'\'el\'ements du noyau de la fl\`eche verticale de gauche;

-  $k^{G_{1}}_{M_{1,U}}(M_{1},L_{1}^U)$ est le nombre d'\'el\'ements du noyau de la fl\`eche verticale de droite;

- $c_{M_{1},M_{2}}^{G_{1},G_{2}}$ est le nombre d'\'el\'ements du noyau de la fl\`eche horizontale du haut, multipli\'e par $n^{-a_{M_{2}}}$;

- $\prod_{v\in U}c_{M_{1,v,sc},M_{2,v,sc}}^{L_{1,SC}^v,L_{2,SC}^v}$ est le nombre d'\'el\'ements de la fl\`eche horizontale de droite, multipli\'e par $\prod_{v\in U}n^{-a_{M_{2,v}}+a_{L_{2}^v}}$.

En utilisant le chemin sud-ouest du diagramme, on voit que le membre de gauche de (3) est le nombre d'\'el\'ements du noyau de la fl\`eche compos\'ee, multipli\'e par $\prod_{v\in U}n^{-a_{M_{2,v}}+a_{L_{2}^v}}$. En utilisant le chemin nord-est, on voit que le membre de droite de (3) est le nombre d'\'el\'ements du m\^eme noyau, multipli\'e par $n^{-a_{M_{2}}}$. Parce que $L_{2}^{U}$ appartient \`a ${\cal L}(M_{2,U})$, on v\'erifie l'\'egalit\'e
$$a_{M_{2}}=\sum_{v\in U}a_{M_{2,v}}-a_{L_{2}^v}.$$
L'\'egalit\'e (3) en r\'esulte, ce qui ach\`eve la d\'emonstration. $\square$

 \bigskip
 
 \subsection{G\'en\'eralisation du th\'eor\`eme 9.1}
  Soient $M_{1}\in {\cal L}(M_{1,0})$ et $M_{2}\in {\cal L}(M_{2,0})$ deux Levi qui se correspondent.
  
  \ass{Proposition}{Supposons $M_{1}\not=G_{1}$. Alors on a l'\'egalit\'e
  $$c_{M_{1},M_{2}}^{G_{1},G_{2}}j_{*}(SA^{M_{1}}_{unip}(V))=SA^{M_{2}}_{unip}(V).$$}
  
  Preuve. Comme en 9.3, de $j_{*}$ se d\'eduit un isomorphisme $j_{*,sc}$ tel que le triplet $(M_{1,SC},M_{2,SC},j_{*,sc})$ soit endoscopique non standard. Puisqu'on suppose $M_{1}\not=G_{1}$, nos hypoth\`eses de r\'ecurrence habituelles nous permettent d'appliquer le th\'eor\`eme 9.1 \`a ce triplet. Donc
  $$j_{*,sc}(SA^{M_{1,SC}}_{unip}(V))=SA^{M_{2,SC}}_{unip}(V).$$
  Pour $i=1,2$, on a l'\'egalit\'e
  $$\tau'(M_{i})^{-1}SA^{M_{i}}_{unip}(V)= \tau'(M_{i,SC})^{-1} \iota^*_{M_{i,SC},M_{i}}(SA^{M_{i,SC}}_{unip}(V))$$
 d'apr\`es la proposition 4.6.  Puisque $\iota^*_{M_{2,SC},M_{2}}\circ j_{*,sc}=j_{*}\circ \iota^*_{M_{1,SC},M_{1}}$, on obtient l'\'egalit\'e
 $$\frac{\tau'(M_{2})\tau'(M_{1,SC})}{\tau'(M_{1})\tau'(M_{2,SC})}j_{*}(SA^{M_{1}}_{unip}(V))=SA^{M_{2}}_{unip}(V).$$
 Il suffit de prouver l'\'egalit\'e
 $$(1) \qquad \frac{\tau'(M_{2})\tau'(M_{1,SC})}{\tau'(M_{1})\tau'(M_{2,SC})}=c_{M_{1},M_{2}}^{G_{1},G_{2}}.$$
 Pour $i=1,2$, $M_{i,SC}$ est simplement connexe. On a d\'ej\`a dit plusieurs fois que cela impliquait $\tau'(M_{i,SC})=1$. Par d\'efinition (cf. 4.1), on a
 $$\tau'(M_{i})=\tau(M_{i})covol(\mathfrak{A}_{M_{i},{\mathbb Z}})^{-1}=\vert \pi_{0}(Z(\hat{M}_{i})^{\Gamma_{F}})\vert  \vert ker^1(F;Z(\hat{M}_{i}))\vert ^{-1}covol(\mathfrak{A}_{M_{i},{\mathbb Z}})^{-1}.$$
 Le groupe $\hat{M}_{i}$ est un Levi du groupe $\hat{G}_{i}$ qui est adjoint. Il en r\'esulte que $Z(\hat{M}_{i})^{\Gamma_{F}}$ est connexe. D'autre part, on a
 $ker^1(F;Z(\hat{M}_{i}))\simeq ker^1(F;Z(\hat{G}_{i}))$ d'apr\`es le lemme [VI] 6.1. Le nombre d'\'el\'ements de ce groupe est \'egal \`a $\tau(G_{i})\vert \pi_{0}(Z(\hat{G}_{i})^{\Gamma_{F}})\vert ^{-1}$. Or ce nombre vaut $1$ parce que $G_{i}$ est simplement connexe. 
 D'o\`u $\tau'(M_{i})=covol(\mathfrak{A}_{M_{i},{\mathbb Z}})^{-1}$.  Le membre de gauche de (1) est \'egal \`a 
 $$covol(\mathfrak{A}_{M_{1},{\mathbb Z}})covol(\mathfrak{A}_{M_{2},{\mathbb Z}})^{-1}.$$
 Rappelons que l'on a un isomorphisme $j_{*}=\mathfrak{A}_{M_{1}}\simeq \mathfrak{A}_{M_{2}}$ compatible aux mesures sur ces espaces. Donc
$$covol(\mathfrak{A}_{M_{1},{\mathbb Z}})=vol(\mathfrak{A}_{M_{1}}/\mathfrak{A}_{M_{1},{\mathbb Z}})=vol(\mathfrak{A}_{M_{2}}/j_{*}(\mathfrak{A}_{M_{1},{\mathbb Z}})).$$
On note $covol(j_{*}(\mathfrak{A}_{M_{1},{\mathbb Z}}))$ ce dernier volume.  Les r\'eseaux $j_{*}(\mathfrak{A}_{M_{1},{\mathbb Z}})$ et $\mathfrak{A}_{M_{2},{\mathbb Z}}$ sont commensurables. 
 Choisissons un entier $n\geq1$ tel que $nj_{*}( \mathfrak{A}_{M_{1},{\mathbb Z}})\subset \mathfrak{A}_{M_{2},{\mathbb Z}}$. Alors
$$covol(j_{*}(\mathfrak{A}_{M_{1},{\mathbb Z}}))=n^{-a_{M_{2}}}covol(nj_{*}(\mathfrak{A}_{M_{1},{\mathbb Z}}))  =n^{-a_{M_{2}}}[\mathfrak{A}_{M_{2},{\mathbb Z}}:nj_{*}(\mathfrak{A}_{M_{1},{\mathbb Z}})]covol(\mathfrak{A}_{M_{2},{\mathbb Z}}).$$
On en d\'eduit que le membre de gauche de (1) est \'egal \`a
$$(2) \qquad n^{-a_{M_{2}}}[\mathfrak{A}_{M_{2},{\mathbb Z}}:nj_{*}(\mathfrak{A}_{M_{1},{\mathbb Z}}].$$

Modifions l'hypoth\`ese sur $n$ en supposant que $nj_{*}(X_{*}(T_{1}))\subset X_{*}(T_{2})$. Comme on l'a dit en 9.4, on a alors un homomorphisme $\hat{j}_{n}:Z(\hat{M}_{2})^{\Gamma_{F}}\to Z(\hat{M}_{1})^{\Gamma_{F}}$. Puisqu'il s'agit de tores complexes, on a aussi un homomorphisme $\hat{j}_{*,n}:X_{*}(Z(\hat{M}_{2}))^{\Gamma_{F}}\to X_{*}(Z(\hat{M}_{1}))^{\Gamma_{F}}$. Le nombre d'\'el\'ements du noyau de $\hat{j}_{n}$ est \'egal \`a celui du conoyau de $\hat{j}_{*,n}$. La d\'efinition de 9.4 conduit donc \`a l'\'egalit\'e
$$c_{M_{1},M_{2}}^{G_{1},G_{2}}=n^{-a_{M_{2}}}\vert coker(\hat{j}_{*,n})\vert .$$
Pour $i=1,2$, les groupes $X_{*}(Z(\hat{M}_{i}))$ et $X^*(M_{i})$ sont isomorphes. L'homomorphisme $\hat{j}_{*,n}$ s'identifie \`a un homomorphisme $X^*(M_{2})^{\Gamma_{F}}\to X^*(M_{1})^{\Gamma_{F}}$ dont on d\'eduit par dualit\'e un homomorphisme $\mathfrak{A}_{M_{1},{\mathbb Z}}\to \mathfrak{A}_{M_{2},{\mathbb Z}}$. On voit que ce dernier n'est autre que $nj_{*}$. Il en r\'esulte d'une part que $n$ satisfait aussi l'hypoth\`ese pos\'ee avant l'\'egalit\'e (2), d'autre part que le conoyau de $\hat{j}_{*,n}$ a m\^eme nombre d'\'el\'ements que celui  du conoyau de l'homomorphisme $nj_{*}:\mathfrak{A}_{M_{1},{\mathbb Z}}\to \mathfrak{A}_{M_{2},{\mathbb Z}}$. 
Donc $c_{M_{1},M_{2}}^{G_{1},G_{2}}$ est lui-aussi \'egal au terme (2), ce qui prouve (1) et la proposition. $\square$

\subsection{Extension de l'ensemble fini de places}
\ass{Lemme}{Supposons qu'il existe un ensemble fini $S$ de places de $F$ contenant $V$ tel que le th\'eor\`eme 9.1 soit v\'erifi\'e pour cet ensemble $S$. Alors il l'est pour l'ensemble $V$.}

Preuve. On peut supposer $S\not=V$. On pose $U=S-V$. Soit $M_{1}\in {\cal L}(M_{1,0})$. Comme en 2.3(3), on peut \'ecrire
$$(1) \qquad SA^{M_{1}}_{unip}(S)=\sum_{\ell=1,...,n_{M_{1}}}Sk_{\ell,U}^{M_{1}}\otimes SA_{\ell,V}^{M_{1}},$$
avec des $Sk_{\ell,U}^{M_{1}}\in D^{st}_{unip}(M_{1}(F_{U}))$ et des $SA_{\ell,V}^{M_{1}}\in D^{st}_{unip}(M_{1}(F_{V}))$. En adaptant les notations, la proposition 2.3(ii) implique
$$(2) \qquad SA^{G_{1}}(V)=\sum_{M_{1}\in {\cal L}(M_{1,0})}\vert W^{M_{1}}\vert \vert W^{G_{1}}\vert^{-1}\sum_{\ell=1,...,n_{M_{1}}}s_{M_{1},U}^{G_{1}}(Sk_{\ell,U}^{M_{1}})(SA_{\ell,V}^{M_{1}})^{G_{1}}.$$ 
Pour tout $M_{1}\in {\cal L}(M_{1,0})$, notons $M_{2}\in {\cal L}(M_{2,0})$ le Levi correspondant. Si $M_{1}\not=G_{1}$, la proposition 9.5 appliqu\'ee \`a l'ensemble $S$ dit que
 $$c_{M_{1},M_{2}}^{G_{1},G_{2}}j_{*}(SA^{M_{1}}_{unip}(S))=SA^{M_{2}}_{unip}(S).$$
 Si $M_{1}=G_{1}$, la constante $c_{G_{1},G_{2}}^{G_{1},G_{2}}$ vaut $1$ et l'\'egalit\'e ci-dessus reste valable d'apr\`es l'hypoth\`ese de l'\'enonc\'e. On d\'eduit alors de (1) l'\'egalit\'e
 $$ SA^{M_{2}}_{unip}(S)=\sum_{\ell=1,...,n_{M_{1}}}Sk_{\ell,U}^{M_{2}}\otimes SA_{\ell,V}^{M_{2}},$$
 o\`u
 $$(3) \qquad Sk_{\ell,U}^{M_{2}}=c_{M_{1},M_{2}}^{G_{1},G_{2}}j_{*}(Sk_{\ell,U}^{M_{1}})$$
 et
 $$(4) \qquad SA_{\ell,V}^{M_{2}}=j_{*}(SA_{\ell,V}^{M_{1}}).$$
  De (4) se  d\'eduit l'\'egalit\'e
$$(SA_{\ell,V}^{M_{2}})^{G_{2}}=j_{*}((SA_{\ell,V}^{M_{1}})^{G_{1}}).$$
De (3) et de la proposition 9.4 se d\'eduit  l'\'egalit\'e
$$s_{M_{2},U}^{G_{1}}(Sk_{\ell,U}^{M_{2}})=s_{M_{1},U}^{G_{1}}(Sk_{\ell,U}^{M_{1}}).$$
Le terme $SA^{G_{2}}(V)$ est calcul\'e par une \'egalit\'e similaire \`a (2).  En utilisant les \'egalit\'es ci-dessus, on voit que le terme de droite de cette \'egalit\'e est \'egal \`a l'image par $j_{*}$ de celui de (2). D'o\`u l'\'egalit\'e $SA^{G_{2}}(V)=j_{*}(SA^{G_{1}}(V))$. $\square$

\bigskip

\subsection{Preuve du th\'eor\`eme 9.1}
En [III] 6.3, on a attach\'e \`a notre triplet endoscopique non standard un triplet particulier $(G,\tilde{G},\omega)$. En fait $\omega=1$ et on le supprime des notations. On consid\`ere ce triplet et on fixe un \'el\'ement ${\cal X}\in {\bf Stab}_{excep}(\tilde{G}(F))$, cf. 3.3. Le lemme 9.6 nous autorise \`a agrandir l'ensemble de places $V$. On peut donc supposer que $V$ contient $S({\cal X},\tilde{K})$. On reprend maintenant la d\'emonstration des sections 5 \`a 8 qui calcule la distribution $\underline{A}^{\tilde{G},{\cal E}}(V,{\cal X})$. On a une premi\`ere simplification car l'ensemble $Fib({\cal X})$ de 5.1 est r\'eduit \`a un \'el\'ement. En effet, comme on l'a dit en 3.3, ${\cal X}$ correspond \`a un \'el\'ement de ${\cal Z}(\tilde{G})^{\Gamma_{F}}$. Notons $\mu$ son image naturelle dans $(T^*/(1-\theta^*)(T^*))\times {\cal Z}(\tilde{G})$.  L'unique \'el\'ement de $Fib({\cal X})$ est $(\mu, 1)$. L'assertion du th\'eor\`eme 3.3 est donc \'equivalente \`a celle de la proposition 5.1. Comme on l'a not\'e en [III] 7.7, le groupe $\bar{G}$ associ\'e \`a $(\mu,1)$ comme en 1.1 est isomorphe au groupe $G_{1}$ de notre triplet endoscopique non standard. En particulier, il est simplement connexe. 

La d\'emonstration des sections 5 \`a 8 vaut jusqu'au point o\`u on avait utilis\'e le th\'eor\`eme [VI] 5.6, c'est-\`a-dire jusqu'en 5.9.   Au d\'ebut de ce paragraphe, on a une donn\'ee ${\bf H}\in E_{\hat{T}_{ad},\star}(\bar{G}_{SC},V)$ et un triplet $({\bf G}',\mu',\omega_{\bar{G}'})\in {\cal J}({\bf H})$. Il s'en d\'eduit un triplet endoscopique non standard $(\bar{H}_{SC},G'_{\epsilon,SC},j_{*})$. Si $N (\bar{H}_{SC},G'_{\epsilon,SC},j_{*})< dim(G_{SC})$, nos hypoth\`eses de r\'ecurrence nous permettent d'appliquer le th\'eor\`eme [VI] 5.6 \`a ce triplet et  on a encore l'\'egalit\'e 5.9(1), c'est-\`a-dire
$$(1) \qquad i(\tilde{G},\tilde{G}',\mu',\omega_{\bar{G}'})S^{{\bf G}'}(\underline{SA}^{{\bf G}'}(V,{\cal X}'),f^{{\bf G}'})=C(\tilde{G})\vert W^{\bar{H}}\vert S^{\bar{H}_{SC}}(SA_{unip}^{\bar{H}_{SC}}(V),\bar{f}_{sc}).$$

 Le lemme [III] 6.3 entra\^{\i}ne que l'on a en tout cas  $N (\bar{H}_{SC},G'_{\epsilon,SC},j_{*})\leq dim(G_{SC})$. Reste le cas o\`u cette in\'egalit\'e est une \'egalit\'e. Dans ce cas, le lemme cit\'e entra\^{\i}ne que la donn\'ee ${\bf G}'=(G',{\cal G}',\tilde{s})$ est \'equivalente \`a la donn\'ee maximale de $(G,\tilde{G})$. Puisque ${\bf G}'$ appartient \`a l'ensemble ${\cal E}_{\hat{T}}(\tilde{G},V)$ d\'efini en 5.1, on voit facilement qu'il n'y a qu'une telle donn\'ee: on peut supposer $\tilde{s}=\hat{\theta}$ et ${\cal G}'=\hat{G}^{\hat{\theta}}\rtimes W_{F}$. La donn\'ee $(\mu',\omega_{\bar{G}'})$ est elle-aussi unique: $\mu'$ est l'image naturelle de $\mu$ dans ${\cal Z}(\tilde{G}')$ et $\omega_{\bar{G}'}=1$. Ces unicit\'es impliquent celle de ${\bf H}$: c'est la donn\'ee principale de $\bar{G}$, c'est-\`a-dire ${\bf H}=(\bar{G},{^L\bar{G}},1)$. Supposons ces conditions  v\'erifi\'ees. Alors le triplet $(\bar{H}_{SC},G'_{\epsilon,SC},j_{*})$ est \'egal \`a notre triplet de d\'epart $(G_{1},G_{2},j_{*})$. Notons $f_{1}=\bar{f}_{sc}$ et $f_{2}=f_{\epsilon,sc}$ avec les notations de 5.8 et 5.9. On a $f_{1}\in  SI(G_{1}(F_{V}))$, $f_{2}\in SI(G_{2}(F_{V}))$  et les fonctions $f_{1}\circ exp$ et $f_{2}\circ exp$ d\'efinies au voisinage de $0$ dans les alg\`ebres de Lie se correspondent par endoscopie non standard. On ne conna\^{\i}t pas l'\'egalit\'e (1) mais on peut en tout cas \'ecrire 
$$ i(\tilde{G},\tilde{G}',\mu',\omega_{\bar{G}'})S^{{\bf G}'}(\underline{SA}^{{\bf G}'}(V,{\cal X}'),f^{{\bf G}'})=C(\tilde{G})\vert W^{\bar{H}}\vert (S^{\bar{H}_{SC}}\left(SA_{unip}^{\bar{H}_{SC}}(V),\bar{f}_{sc})+X(f)\right),$$
o\`u
$$(2) \qquad X(f)=S^{G_{2}}(SA^{G_{2}}_{unip}(V),f_{2})-S^{G_{1}}(SA^{G_{1}}_{unip}(V),f_{1}).$$

Pour toutes les donn\'ees ${\bf H}$ sauf la donn\'ee maximale de $\bar{G}$, le calcul de 5.9 est donc valable. Pour la donn\'ee maximale, la formule 5.9(2) doit \^etre corrig\'ee: on ajoute au membre de droite le terme $X(f)$ multipli\'e par une constante. Il est facile de calculer celle-ci: c'est $C(\tilde{G})$. Le calcul se poursuit et on obtient finalement non pas l'\'egalit\'e 8.4(1), mais l'\'egalit\'e
$$I^{\tilde{G}}(\underline{A}^{\tilde{G},{\cal E}}(V,\mu,\omega_{\bar{G}}),f)=I^{\tilde{G}}(A^{\tilde{G}}(V,{\cal X}),f)+C(\tilde{G})X(f).$$
Comme on l'a dit ci-dessus, le membre de gauche est \'egal \`a $I^{\tilde{G}}(\underline{A}^{\tilde{G},{\cal E}}(V,{\cal X}),f)$.  D'apr\`es nos hypoth\`eses de r\'ecurrence, le th\'eor\`eme [VI] 5.4 est connu pour $(G,\tilde{G})$.  Comme on l'a vu en 3.6, l'\'egalit\'e du th\'eor\`eme 3.3 est donc v\'erifi\'ee. Cela entra\^{\i}ne   l'\'egalit\'e
$$I^{\tilde{G}}(\underline{A}^{\tilde{G},{\cal E}}(V,\mu,\omega_{\bar{G}}),f)=I^{\tilde{G}}(A^{\tilde{G}}(V,{\cal X}),f).$$
Il en r\'esulte que $X(f)=0$. Le m\^eme raisonnement qu'en [III] 7.7 montre que, pour tout $\varphi\in SI(G_{1}(F_{V}))$, il existe $f$ tel que $f_{1}$ co\"{\i}ncide avec $\varphi$ au voisinage de l'unit\'e. L'\'egalit\'e $X(f)=0$ pour tout $f$ entra\^{\i}ne donc l'\'egalit\'e voulue $j_{*}(SA^{G_{1}}_{unip}(V))=SA^{G_{2}}_{unip}(V)$. Cela prouve le th\'eor\`eme 9.1. $\square$

{\bf Bibliographie}

\bigskip

[A1] J. Arthur: {\it A stable trace formula I. General expansions}, Journal of the Inst. of Math. Jussieu 1 (2002), p. 175-277

[A2] -----------: {\it A stable trace formula II. Global descent}, Inventiones Math. 143 (2001), p. 157-220

[CL] P.-H. Chaudouard, G. Laumon: {\it Le lemme fondamental pond\'er\'e I: constructions g\'eom\'etriques}, Compositio Math. 146 (2010), p. 1416-1506

[K1] R. Kottwitz: {\it Rational conjugacy classes in reductive groups}, Duke Math. J. 49 (1982), p. 785-806

[K2] --------------: {\it Stable trace formula: elliptic singular terms}, Math. Ann. 275 (1986), p. 365-399

[K3] --------------: {\it Stable trace formula:  cuspidal tempered terms}, Duke Math. J. 51 (1984), p. 611-650

[KS] -------------, D. Shelstad: {\it Foundations of twisted endoscopy}, Ast\'erisque 255 (1999)

[Lab1] J.-P. Labesse: {\it Nombres de Tamagawa des groupes r\'eductifs quasi-connexes}, manuscripta math. 104 (2001), p. 407-430

[Lab2] -----------------: {\it Cohomologie, stabilisation et changement de base}, Ast\'erisque 257 (1999)

[Lab3] -----------------: {\it Stable twisted trace formula: elliptic terms}, Journal of the Inst. of Math. Jussieu 3 (2004), p. 473-530

[Lan] R. P. Langlands: {Stable conjugacy: definitions and lemmas}, Can. J. Math. 31 (1979), p. 700-725

[LS] -------------------, D. Shelstad: {\it On the definition of transfer factors}, Math. Ann. 278 (1987), p. 219-271

[Oe] J. Oesterl\'e: {\it Nombres de Tamagawa et groupes unipotents en caract\'eristique $p$}, Inventiones Math. 78 (1984), p. 13-88

[S] J.-L. Sansuc: {\it Groupe de Brauer et arithm\'etique des groupes alg\'ebriques lin\'eaires sur un corps de nombres}, Journal f\"ur die r. und ang. Math. 327 (1981), p. 12-80

[W1] J.-L. Waldspurger: {\it L'endoscopie tordue n'est pas si tordue}, Memoirs AMS 908 (2008)

[W2] ------------------------: {\it A propos du lemme fondamental pond\'er\'e tordu}, Math. Annalen 343 (2009), p. 103-174

[I], [II], [III], [V] ----------------------:{ \it  Stabilisation de la formule des traces tordue I: endoscopie tordue sur un corps local}, {\it II: int\'egrales orbitales et endoscopie sur un corps local non-arhcim\'edien; d\'efinitions et \'enonc\'es des r\'esultats}, {\it III: int\'egrales orbitales et endoscopie sur un corps local non-arhcim\'edien; r\'eductions et preuves}, {\it V: int\'egrales orbitales et endoscopie sur le corps r\'eel}, pr\'epublications 2014

[VI] C. Moeglin, J.-L. Waldspurger: {\it Stabilisation de la formule des traces tordue VI: la partie g\'eom\'etrique de cette formule}, pr\'epublication 2014

\bigskip

e-mail:   jean-loup.waldspurger@imj-prg.fr

 \end{document}